\documentclass[12pt]{article}
\author{Constantin N. Beli}
\title{Representations of quadratic lattices over dyadic local fields}
\usepackage{amsfonts,amsmath,amssymb}

\def\a{\alpha} \def\b{\beta} \def\c{\gamma} 
\def\D{\Delta}  \def\e{\varepsilon} 
  \def\h{\frac} 
\def\j{\infty}   
\def\m{\lim}    
\def\p{\partial}   \def\s{\sigma}
   
 \def\te{\theta}  
  \def\z{\longrightarrow}
\def\zz{\Longleftrightarrow}\def\({\overline} \def\){\underline}
\def\<{\cdot} \def\go{\mathfrak}
\def\>{~~~~~~~}
\def\be{\begin{equation}}
\def\ee{\end{equation}}

\def\sb{\subset} \def\sp{\supset} \def\sbq{\subseteq}
\def\spq{\supseteq} \def\ti{\times} \def\od{{\rm ord}\,} \def\oo{{\cal
O}} \def\pp{\perp} \def\ss{{\go s}} \def\nn{{\go n}} \def\ff{\dot{F}}
\def\ooo{{\oo^\ti}} \def\mo{{\rm mod}~}  \def\hh{{\rm H}}
\def\aa{A(0,0)} \def\ab{A(2,2\rho )} \def\fs{\ff^2}  \def\os{\ooo\fs}
\def\p{\go p} \def\*{\sharp} 
\def\rr{{\cal R}} \def\WW{{\go W}} \def\0{} 
\def\1{^{-1}}  \def\dd{{\go d}} \def\aaa{{\cal
A}} \def\bbb{{\cal B}} \def\[{\prec} \def\]{\succ}
\def\bm{\left(\begin{array}} \def\em{\end{array}\right)}
\def\ev{\equiv} \def\ap{\cong}  
 \def\N{{\rm N}} 
\def\la{\langle} \def\ra{\rangle} \def\rep{{\rightarrow\!\!\!\! -}}
\def\notrep{{\not\rightarrow\!\!\!\! -}}
 \def\m2{~(\mo 2)} \def\no{\noindent}
 \def\btm{\begin{thm}}
\def\etm{\end{tm}}
 \def\blem{\begin{lem}}
\def\elem{\end{lem}}

\newtheorem{theorem}{Theorem}[section]
\newtheorem{proposition}[theorem]{Proposition}
\newtheorem{lemma}[theorem]{Lemma}
\newtheorem{definition}{Definition}
\newtheorem{corollary}[theorem]{Corollary}

\newtheorem{bof}[theorem]{}
\newtheorem{teorema}{Theorem}

\def\qed{\mbox{$\Box$}\vspace{\baselineskip}}
\def\pf{$Proof.$} 
\def\bco{\begin{corollary}} \def\eco{\end{corollary}}
\def\bpr{\begin{proposition}} \def\epr{\end{proposition}}
\def\bdf{\begin{definition}} \def\edf{\end{definition}}
\def\btm{\begin{theorem}} \def\etm{\end{theorem}}
\def\blm{\begin{lemma}} \def\elm{\end{lemma}}
\def\bff{\begin{bof}\rm} \def\eff{\end{bof}}
\def\btr{\begin{teorema}} \def\etr{\end{teorema}}

\def\de{\newcommand} \de\tm[1]{{\no\bf Theorem~#1}} 
 \def\mb{\mathbb} 
\def\RR{{\mb R}} \def\QQ{{\mb Q}}  \def\ZZ{{\mb Z}}
\def\NN{{\mb N}} \def\AA{{\mathsf a}}  \def\BB{{\mathsf b}} 
\def\FF{{\go f}}  \def\GG{{\go g}} \def\ww{{\go w}}
\de\lm[1]{{\no\bf Lemma~#1}}
\de\df[1]{{\no\bf Definition~#1}} \de\co[1]{{\no\bf Corollary~#1}}
\de\tp[1]{\te (#1 )} \de\ts[1]{\te (O^-(#1 ))} \de\ty[1]{\te
(O(#1 ))} \de\tx[1]{\te (#1 )} \de\up[1]{(1+\p^{#1} )\fs}
 \de\upn[2]{(1+\p^{#1})\fs\cap\N (#2 )} \de\xt[2]{\te (#1 /#2 )}
 \de\ups[1]{((1+\p^{#1})\fs )^*} \de\upo[1]{(1+\p^{#1} )\ooo^2}
\de\upon[2]{(1+\p^{#1})\ooo^2\cap\N (#2 )}
\de\lr[1]{\longrightarrow^{\!\!\!\!\!\!\!\! #1}}
\de\lf[1]{\longleftarrow^{\!\!\!\!\!\!\!\! #1}}
\DeclareMathOperator\ord{ord}

\begin{document}
\maketitle

\section*{Introduction}

An important part in the arithmetic theory of quadratic forms is
understanding the local theory. It is well known that the local theory
of the quadratic forms is significantly more difficult when the base
field is dyadic, i.e. when it is a finite extension of $\QQ_2$. The
2-adic case (unramified finite extensions of $\QQ_2$) is somewhat
easier than the general dyadic case. Typically problems are first
proved over non-dyadics, then over 2-adics and finally over general
dyadics. 

We now give a brief list of some of the most important local problems
in the local theory of the quadratic forms, together with their
history and current status. 

1. Classification. Deciding if two quadratic lattices $M,N$ are
isometric, i.e. if $M\ap N$. 

2. Representation. Deciding if a lattice $M$ represents a lattice
$N$, i.e. if $N\rep M$. 

3. Calculation of the integral spinor norm group $\te (O^+(M))$.

4. Calculation of the relative integral spinor norm group $\te
(X(M/N))$ when $N\rep M$.

5. Primitive representation. Deciding if a lattice $M$ represents
primitively a lattice $N$, i.e. if $N\rep^*M$. 

6. Calculation of the primitive relative integral spinor group $\te
(X^*(M/N))$ when $N\sbq^*M$.

Problem 1. was solved completely. The dyadic case was done by O'Meara
[OM, 93:28]. Problem 2., which is the subject of this paper, was
solved by O'Meara [OM1] in the 2-adic case and by Riehlm [R] when the
bigger lattice $M$ is modular. Problem 3. was solved by Earnest and
Hsia [EH1] in the 2-adic case and by the author [B1] in the general
dyadic case, by using BONGs, a new way of describing quadratic
lattices. Recently Lv and Xu [LX] gave an alternate solution in terms
of the traditional Jordan decompositions. Problem 4. was first
considered by Hsia, Shao and Xu in [HSX], where this group was
introduced. The authors proved that $\theta (X(M/N))$ is a group and
they calculated it in the non-dyadic case. Later Shao calculated this
group in the 2-adic case in his PhD thesis [OSU1], a result that is
extremely complicated. The problem was solved in the general case by
the author in his PhD thesis [OSU2]. Very little is known about
problems 5. and 6.. Problem 5. was solved by James [J] in the
non-diadic case when the bigger lattice is modular. 

There are several global results which are incomplete due to lack of
knowledge of the local theory in the dyadic case. Often these results
are stated under the condition that the quadratic forms involved
``behave well'' at the dysdic primes or that all dyadic primes of the
underlying field are 2-adic, i.e. unramified. 

One example of such incomplete results is the so-called aritmetic
Springer theorem. Let $F\sbq E$ be two number fields with $[E:F]$ odd
and let $\oo_F,\oo_E$ be their rings of integers. Given two quadratic
lattices $M,N$ over $\oo_F$ and $\tilde M=\oo_E\otimes_{\oo_F}M$,
$\tilde N=\oo_E\otimes_{\oo_F}N$, we want to know if $\tilde
M\ap\tilde N$ implies $M\ap N$ and, more generaly, if $\tilde
N\rep\tilde M$ implies $N\rep M$. In the case when $M$ is undefined
these two problems can be tackled by using spinor genera and they
depend on some local results. If $\go P$ is a prime of $E$ lying over
a prime $\p$ of $F$ and $\N_{E_{\go P}/F_\p}:F_\p\to E_{\go P}$ is the
norm map then one has to prove that $\N_{E_{\go P}/F_\p}(\theta
(O^+(\tilde M_{\go P}))\sbq\theta (O^+(M_\p))$ and if $[E_{\go
P}:F_\p]$ is odd then $M_p\ap N_\p$ iff $\tilde M_{\go P}\ap\tilde
N_{\go P}$. For the representation problem one has to prove that if
$N_\p\sbq M_\p$ then $\N_{E_{\go P}/F_\p}(\theta (X(\tilde M_{\go
P}/\tilde N_{\go  P}))\sbq\theta (X(M_\p/N_\p ))$ and if $[E_{\go
P}:F_\p]$ is odd then $N_p\rep M_\p$ iff $\tilde N_{\go P}\rep\tilde
M_{\go P}$. The proof of these four statements depend on knowing how
to solve problems 1, 2, 3 and 4 above. When Earnest and Hsia first
considered this problem in [EH2] only problem 1. was completely
solved so their result depended on some restrictions on the
genus at dyadic primes. The most general result is due to Xu. In
[X] he proved the norm principle $\N_{E_{\go P}/F_\p}(\theta
(X(\tilde M_{\go P}/\tilde N_{\go P}))\sbq\theta (X(M_\p/N_\p ))$,
without calculating $\theta (X(\tilde M_{\go P}/\tilde N_{\go P})$ and
$\theta (X(M_\p/N_\p ))$, by using some reduction formulas and
induction on rank. He also proved some necessary conditions for
$\tilde N_{\go P}\rep\tilde M_{\go P}$, which in turn imply some
conditions on $M_\p,N_\p$. When $F_\p$ is 2-adic these conditions
imply $N_\p\rep M_\p$ by O'Meara's 2-adic representatin theorem from
[OM1]. So he obtained the representation Springer theorem when $2$ is
unramified in $F$. The norm principle can be proved also directly if
we use the formulas for the relative spinor norm group from [OSU2]. And
by using the main result of this paper, Theorem 2.1, one can prove the
equivalence between $\tilde N_{\go P}\rep\tilde M_{\go P}$ and
$N_\p\rep M_\p$ when $[E_{\go P}:F_\p ]$ is odd for arbitrary dyadic
fields. Hence the arithmetic Springer theorem for spinor genera holds
in the most general setting.

The first mathematician to solve nontrivial problems regarding
quadratic forms over dyadic fields was O'Meara. In his renown Theorem
93:28 he solved the difficult problem of classifying lattices over
dyadic fields. Later, in 1958 [OM1], he tackled the problem of
representation of lattices over local fields. He successfully
solved this problem in the case when the base fields is either
non-dyadic or 2-adic (i.e. a nonramified extension of $\QQ_2$). As
one can see, while the result in the non-dyadic case involves only
one condition, the 2-adic case is much more complicated and it
requires ten conditions. (There are five conditions involved in
the fact that the smaller lattice is of ``lower type'' and there
are the five conditions I-V of the main theorem.) In 1962 C. Riehm
[R] solved the same problem over arbitrary dyadic fields, when the
larger lattice is modular. Again the result is quite complicated
and difficult to apply. Since then there was no progress towards a
solution in the general dyadic case except for some neccesary
conditions obtained by Xu in [X, \S4]. Although thy are not stated in
terms of BONGs, Xu's results have correspondences in our paper. His
[X, Proposition 4.2] follows from condition (i) of Theorem 2.1 (parts
(1) and (2)(ii) and (iii)) and and from Lemma 6.5 (parts (2)(i), (3)
and (4)). And his [X, Proposition 4.8] is related to our Lemma
2.19. (See also the remark following Lemma 2.19.)

In this paper we completely solve the problem of representation
for lattices over general dyadic fields. Our result, Theorem 2.1,
is given in terms of bases of norm generators, BONGs for short. A
simple algorithm fo finding a good BONG is provided in [B2, \S7]. 
However our main theorem can be translated in terms of the traditional
Jordan splittings, as stated in Lemma 3.10. The result from Lemma 3.10
is given in terms of the so called approximations (Definitions 9 and
10), which can be obtained in terms of Jordan splittings as stated in
Lemmas 3.2 and 3.7. (See also [B2, \S5].)

The usefulness of using BONGs rather then Jordan splittings was
already seen in [B1], where we computed the spinor norm group $\theta
(O^+(L))$ for lattices $L$ over arbitrary dyadic fields. It seems that
by using of BONGs, even though the proofs are still quite long (this
paper is a good example), the results are expressed in a very compact
form. See e.g. Theorem 2.1 of this paper compared with the results
from [OM1] and [R] or the formulas for the relative spinor norm groups
in in the dyadic case from [OSU2], compared to those in the 2-adic
case from [OSU1]. Also the use of BONGs makes it easyer to recognise
patterns and make educated guesses on general results from partcular
cases in small dimensions. E.g. the main result of this paper, Theorem
2.1, was guessed after the case of rank 4 was understood. 
\vskip 3mm

In [B3] we translated O'Meara's classification theorem for lattices
over dyadic fields in terms of BONGs. That was a first step towards
the more difficult problem of representation. Therefore in this paper
we will use terminology and results from [B3]. We will also use
notations, definitions and results from [B1]. We now give a summary of
the main results from [B1] and [B3] we use here. An overview of the
results from [B1] can also be found in the introduction of [B3].

All the quadratic spaces and lattices in this paper will be over a
dyadic field $F$. We denote by $\oo$ the ring of integers, $\p$
the prime ideal, $\ooo :=\oo\setminus\p$ the group of units,
$e:=\ord 2$ and $\pi$ is a fixed prime element. For $a\in\fs$ we
denote its quadratic defect by $\dd (a)$ and we take $\D =1-4\rho$
a fixed unit with $\dd (\D)=4\oo$.

If $x_1,\ldots,x_n$ is an orthogonal set of vectors with $Q(x_i)=a_i$
we say that $V\ap [a_1,\ldots,a_n]$ relative to $x_1,\ldots,x_n$ if
$V=Fx_1\pp\cdots\pp Fx_n$ and we say that $L\ap\la a_1,\ldots,a_n\ra$
relative to $x_1,\ldots,x_n$ if $L=\oo x_1\pp\cdots\pp\oo x_n$.

We denote by $d:\ff/\fs\z\NN\cup\{\j\}$ the order of the ``relative
quadratic defect'' $d(a):=\ord a\1\dd (a)$. If $a=\pi^R\e$ then
$d(a)=0$ if $R$ is odd and $d(a)=d(\e )=\ord\dd (\e )$ if $R$ is
even. Thus $d(\ff )=\{ 0,1,3,\ldots,2e-1,2e,\j\}$. $d$ satisfies the
domination principle $d(ab)\geq\min\{ d(a),d(b)\}$. Also for
$\a\in\RR\cap\{\j\}$ we denote $\up\a :=\{ a\in\ff\mid d(a)\geq\a\}$
and $\upo\a :=\{ a\in\oo\mid d(a)\geq\a\}$. 

We denote by $(\cdot,\cdot )_\p :\ff/\fs\times\ff/\fs\z\{\pm 1\}$
the Hilbert symbol. If $a\in\ff$ we denote $\N (a)=\{ b\in\ff|
(a,b)_\p =1\} =\N (F(\sqrt a)/F)$. If $b\in\ff$ and $d(a)+d(b)>2e$
then $(a,b)_\p =1$. This result doesn't hold anymore if
$d(a)+d(b)\leq 2e$. In fact if $a\notin\fs$ then there is
$b\in\ff$ with $d(b)=2e-d(a)$ s.t. $(a,b)_\p =-1$. Thus
$\up\a\sbq\N (a)$ iff $\a +d(a)>2e$. For further properties of $d$
and $(\cdot,\cdot )_\p$ see also $\S 8$, Lemmas 8.1 and 8.2.

An element $x$ of a lattice $L$ is called norm generator of $L$ if
$\nn L=Q(x)\oo$. A set $x_1,\ldots,x_n$ is called a basis of norm
generators (BONG) of $L$ if $x_1$ is a norm generator of $L$ and
$x_2,\ldots,x_n$ is a BONG of $pr_{x_1^\pp}L$. A BONG uniquely determine
a lattice so if $x_1,\ldots,x_n$ is a BONG for $L$ we will write $L=\[
x_1,\ldots,x_n\]$. If moreover $Q(x_i)=a_i$ we say that $L\ap\[
a_1,\ldots,a_n\]$ relative to the BONG $x_1,\ldots,x_n$. If $L\ap\[
a_1,\ldots,a_n\]$ then $\det L=a_1\cdots a_n$. 

If $L$ is binary with $\nn L=\a\oo$ we denote by $a(L)=\det
L\a^{-2}$ and by $R(L)=\ord vol L-2\ord\nn L=\ord a(L)$.
$a(L)\in\ff/\ooo^2$ is an invariant of $L$ and it determines the
class of $L$ up to scaling. If $L\ap\[\a,\b\]$ then $a(L)=\h\b\a$.

The set of all possible values of $a(L)$ where $L$ is an arbitrary binary
lattice is denoted by $\aaa =\aaa_F$. If $a\in\ff$ and $\ord a=R$ then
$a\in\aaa$ iff $R+2e\geq 0$ and $R+d(-a)\geq 0$. Some further properties of
$\aaa$ can be found in [B1] (see the comments to [B1, Lemma 3.5]):

If $a(L)=a=\pi^R\e$ with $d(-a)=d$ then:

$L$ is nonmodular, proper modular or improper modular iff $R>0$, $R=0$
resp. $R<0$.

If $R$ is odd then $R>0$.

The inequality $R+2e\geq 0$ becomes equality iff $a\in -\h 14\ooo^2$ or
$a\in -\h\D 4\ooo^2$. We have $a(L)=-\h 14$ resp. $a(L)=-\h\D 4$ when
$L\ap\pi^r\aa$ resp. $\pi^r\ab$ for some integer $r$.

The inequality $R+d\geq 0$ becomes equality iff $a\in -\h\D 4\ooo^2$.
\vskip 3mm

If $L\ap\[ a_1,\ldots,a_n\]$ relative to some BONG $x_1,\ldots,x_n$ and $\ord
a_i=R_i$ we say that the BONG is good if $R_i\leq R_{i+2}$.

A set $x_1,\ldots,x_n$ of orthogonal vectors with $Q(x_i)=a_i$ and
$\ord a_i=R_i$ is a good BONG for some lattice iff $R_i\leq
R_{i+2}$ for all $1\leq i\leq n-2$ and $a_{i+1}/a_i\in\aaa$ for
all $1\leq i\leq n-1$. The condition that $a(\[ a_i,a_{i+1}\]
)=a_{i+1}/a_i\in\aaa$ is equivalent to $R_{i+1}-R_i+2e\geq 0$ and
$R_{i+1}-R_i+d(-a_ia_{i+1})\geq 0$. If $R_{i+1}-R_i=-2e$ then
$a_{i+1}/a_i\in -\h 14\ooo^2$ or $-\h\D 4\ooo^2$ so $\[
a_i,a_{i+1}\]$ is $\ap\pi^r\aa$ or $\pi^r\ab$ for some $r$. In
fact since $\nn\[ a_i,a_{i+1}\] =a_i\oo$ we have $\[
a_i,a_{i+1}\]\ap\h 12a_i\aa$ or $\h 12a_i\ab$. Also if
$R_{i+1}-R_i$ is odd then $R_{i+1}-R_i>0$. (See the properties of
$\aaa$ above.)

Good BONGs enjoy some properties similar to those of the orthogonal
bases. If $L\ap\[ a_1,\ldots,a_n\]$ then relative to some good BONG
$x_1,\ldots,x_n$ and $\ord a_i=R_i$ then $L^\*\ap\[ a_1\1,\ldots,a_n\1\]$
relative to the good BONG $x_n^\*,\ldots,x_1^\*$ where $x_i^\* =Q(x)\1
x_i$. Also if for some $1\leq i\leq j\leq n$ we have $\[
x_i,\ldots,x_j\]\ap\[ b_i,\ldots,b_j\]$ relative to some other good BONG
$y_i,\ldots,y_j$ then $L\ap\[
a_1,\ldots,a_{i-1}b_i,\ldots,b_j,a_{i+1},\ldots,a_n\]$ relative to the good
BONG $x_1,\ldots,x_{i-1}y_i,\ldots,y_j,x_{i+1},\ldots,x_n$.

The most important notion introduced in [B3] are the invariants $\a_i$'s. If
$L\ap\[ a_1,\ldots,a_n\]$ relative to some good BONG and $\ord a_i=R_i$ then for
any $1\leq i\leq n-1$ we define $\a_i$ as the minimum of the set 
\begin{multline*}
\{ (R_{i+1}-R_i)/2+e\}\cup\{
R_{i+1}-R_j+d(-a_ja_{j+1})\mid 1\leq j\leq i\}\\
\cup\{ R_{j+1}-R_i+d(-a_ja_{j+1})\}\mid i\leq j\leq n-1\}.
\end{multline*}

Now $R_i=\ord a_i$ and $\a_i$ defined above are independent of the
choice of the good BONG so we denote them by $R_i(L)$ and $\a_i(L)$.

Some properties of $\a_i$:

The sequence $(R_i+\a_i)$ is increasing and the sequence
$(-R_{i+1}+\a_i)$ is decreasing.

$\a_i\in ([0,2e]\cap\ZZ )\cup ([2e,\j )\cap\h 12\ZZ )$.

$\a_i\geq 0$ with equality iff $R_{i+1}-R_i=-2e$.

$\a_i$ is $<2e$, $=2e$ or $>2e$ iff $R_{i+1}-R_i$ is so.

If $R_{i+1}-R_i\geq 2e$ then $\a_i=(R_{i+1}-R_i)/2+e$.

If $R_{i+1}-R_i\leq 2e$ then $\a_i\leq R_{i+1}-R_i$ with equality iff
$R_{i+1}-R_i$ is equal to $2e$ or it is odd.

In the following cases $\a_i$ is uniquely determined by $R_{i+1}-R_i$:
If $R_{i+1}-R_i$ is equal to $-2e,2-2e,2e-2$ or it is $\geq 2e$ then
$\a_i=(R_{i+1}-R_i)/2+e$. If $R_{i+1}-R_i$ is odd then $\a_i=\min\{
(R_{i+1}-R_i)/2+e,R_{i+1}-R_i\}$.
\vskip 3mm

The invariants $R_i$ and $\a_i$ of a lattice $L$ can be expressed in
terms of Jordan splittings. If $L=L_1\pp\cdots\pp L_t$ is a Jordan
splitting, $L_{(k)}:=L_1\pp\cdots\pp L_k$, $n_k=\dim L_{(k)}$,
$\ss_k=\ss L_k$, $\GG_k=\GG L^{\ss_k}$, $\AA_k$ is a norm generator
for $\GG_k$, $\FF_k$ are the invariants introduced in [OM] $\S$93E,
page 264, $r_k=\ord\ss_k$ and $u_k=\ord\AA_k$ then $R_i$'s are in
one-to-one correspondence with $t$, $\dim L_k$, $\ss L_k$ and
$\AA_k\oo$, and $\a_i$'s are in one-to-one correspondence with
$\ww_k$'s and $\FF_k$'s. 

We now state the main result of [B3], Theorem 3.1:
\vskip 3mm
Let $L\ap\[ a_1,\ldots,a_n\]$ and $K\ap\[ b_1,\ldots,b_n\]$ relative
to some good BONGs with $R_i(L)=R_i$, $R_i(K)=S_i$, $\a_i(L)=\a_i$ and
$\a_i(K)=\b_i$. If $FL\ap FK$ then $L\ap K$ iff:

(i) $R_i=S_i$ for $1\leq i\leq n$.

(ii) $\a_i=\b_i$ for $1\leq i\leq n-1$.

(iii) $d(a_1\cdots a_ib_1\cdots b_i)\geq\a_i$ for $1\leq i\leq n-1$.

(iv) $[b_1,\ldots,b_{i-1}]\rep [a_1,\ldots,a_i]$ for any $2\leq i\leq
n-1$ s.t. $\a_{i-1}+\a_i>2e$.
\vskip 3mm

If $K=K_1\pp\cdots\pp K_{t'}$, $\ss K_k=\ss'_k$, $\GG'_k=\GG
K^{\ss'_k}$, $\ww'_k=\ww K^{\ss'_k}$, $\FF'_k$ are the $\FF_k$
invariants corresponding to $K$ and $\BB_k$ is a norm generator
for $K^{\ss'_k}$ then: (i) is equivalent to $t=t'$, $\dim L_k=\dim
K_k$, $\ss_k=\ss'_k$ and $\ord\AA_k=\ord\BB_k$; assuming that (i)
holds then (ii) is equivalent to $\ww_k=\ww'_k$ and
$\FF_k=\FF'_k$; assuming that (i) and (ii) hold then (iii) is
equivalent to $\AA_k\ap\BB_k\mo\ww_k$ for $1\leq k\leq t$ and the
condition (i) of [OM] 93:28; if (i)-(iii) hold then (iv) is
equivalent to conditions (ii) and (iii) of [OM] 93:28.
\vskip 5mm

\section{Preliminary results}

In this section we give some definitions and results which will be
used throughout the paper. For the purpose of understanding the main
theorem one only needs the definition of $d[\e a_{1,i}b_{1,j}]$. So the
reader who is interested only in the main result only needs to read
Defintion 1 and skip the rest of the section. 

\subsection{The invariants $d[\e a_{1,i}b_{1,j}]$ and $d[\e a_{i,j}]$}

In order to shorten our formulas we will use the following notations.
\vskip 4mm
{\bf Notation} If $x_1,x_2,\ldots$ is a sequence in $\ff$ (or
$\ff/\fs$) then we denote $x_{i,j}:= x_ix_{i+1}\cdots x_j$ for any
$i\leq j$. Also by convention $x_{i,i-1}=1$.

\bdf Let $M\ap\[ a_1,\ldots,a_m\]$ and $N\ap\[ b_1,\ldots,b_n\]$ with
$\a_i(M)=\a_i$ and $\a_i(N)=\b_i$.

If $0\leq i\leq m$, $0\leq j\leq n$ and $\e\in\ff$ let 
$$d[\e a_{1,i}b_{1,j}]:=\min\{ d(\e a_{1,i}b_{1,j}),\a_i,\b_j\}.$$

If $0\leq i-1\leq j\leq m$ and $\e\in\ff$ let
$$d[\e a_{i,j}]:=\min\{ d(\e a_{i,j}),a_{i-1},\a_j\}.$$

(If $i\in\{ 0,m\}$ or $j\in\{ 0,n\}$ then we ignore $\a_i$ resp. $\b_j$
from the definition of $d[\e a_{1,i}b_{1,j}]$. Also if $i-1$ or $j$
belongs to $\{ 0,m\}$ then we ignore $\a_{i-1}$ resp $\a_j$ from the
definition of $d[\e a_{i,j}]$.)
\edf

Note that in $\ff/\fs$ we have $a_{i,j}=a_{1,i-1}a_{1,j}$ so $d(\e
a_{i,j})=d(\e a_{1,i-1}a_{1,j})$. Hence  $d[\e
a_{i,j}]:=d[\e a_{1,i-1}a_{1,j}]$, i.e. the invariants $d[\e a_{i,j}]$
are particular cases of the invariants $d[\e a_{1,i}b_{1,,j}]$ when
the two lattices coincide.

In general we will take $\e =\pm 1$. More precisely we will take $d[\e
a_{1.i}b_{1,j}]$ with $i\geq i$, $i\ev j\m2$ and $\e =(-1)^{(i-j)/2}$
and $d[\e a_{i,j}]$ with $j\ev i-1\m2$ and $\e =(-1)^{(j-i+1)/2}$.

\bff {\bf Remark} Since $d[-a_{i,i+1}]=\min\{
d(-a_{i,i+1}),\a_{i-1},\a_{i+1}\}$ [B3, Corollary 2.5(i)] writes as 
$$\a_i=\min\{ (R_{i+1}-R_i)/2+e,R_{i+1}-R_i+d[-a_{i,i+1}]\}.$$ 
In particular, $\a_i\leq R_{i+1}-R_i+d[-a_{i,i+1}]$ so
$d[-a_{i,i+1}]\geq R_i-R_{i+1}+\a_i$.

Also note that if $1\leq j\leq i$ then $R_{i+1}-R_j+d[-a_{j,j+1}]\geq
R_{i+1}-R_{j+1}+\a_j\geq\a_i$ and if $i\leq j\leq n-1$ then
$R_{j+1}-R_i+d[-a_{j,j+1}]\geq R_j-R_i+\a_j\geq\a_i$. Hence the
inequalities $R_{i+1}-R_j+d(-a_{j,j+1})\geq\a_i$ for $j\leq i$ and
$R_{j+1}-R_i+d(-a_{j,j+1})\geq\a_i$ for $j\geq i$ which follow from
the definition of $\a_i$ also hold when we replace $d(-a_{j,j+1})$ by
$d[-a_{j,j+1}]$.
\eff

\blm Definition 1 is independent of the choice of the BONGs.
\elm

\pf Let $M\ap [a_1',\ldots,a_m']$ relative to another good BONG. We
have $d(a_{1,i}a'_{1,i})\geq\a_i$ by [B3, Theorem 1] so $d(\e
a'_{1,i}b_{1,j})\geq\min\{ d(\e
a_{1,i}b_{1,j}),d(a_{1,i}a'_{1,j})\}\geq\min\{ d(\e
a_{1,i}b_{1,j}),\a_i\}$. Thus $d[\e a'_{1,i}b_{1,j}]=\min\{ d(\e
a'_{1,i}b_{1,j}),\a_i,\b_j\}\geq\min\{ d(\e
a_{1,i}b_{1,j}),\a_i,\b_j\} =d[\e a_{1,i}b_{1,j}]$. Similarly
$d[\e a_{1,i}b_{1,j}]\geq d[\e a'_{1,i}b_{1,j}]$ so $d[\e
a_{1,i}b_{1,j}]=d[\e a'_{1,i}b_{1,j}]$. Similarly we prove that
$d[\e a_{1,i}b_{1,j}]$ is independent of the choice of the BONG
$b_1,\ldots,b_n$. \qed

\bff Note that the numbers $d[\e a_{1,i}b_{1,j}]$ satisfy similar
domination principle as $d(\e a_{1,i}b_{1,j})$. Namely If $K\ap\[
c_1,\ldots,c_k\]$ is a third lattice with $\a_i(K)=\c_i$ then
$d[\e\e'a_{1,i}c_{1,l}]\geq\\ \min\{ d[\e
a_{1,i}b_{1,j}],d[\e'b_{1,j}c_{1,l}]\}$. Indeed by the domination
principle [B1, 1.1] we have $d(\e\e'a_{1,i}c_{1,l})\geq\min\{ d(\e
a_{1,i}b_{1,j}),d(\e'b_{1,j}c_{1,l})\}$, which implies that\\
$d[\e\e'a_{1,i}c_{1,l}]= \min\{
d(\e\e'a_{1,i}c_{1,l}),\a_i,\b_j\}\geq\min\{ d(\e
a_{1,i}b_{1,j}),d(\e'b_{1,j}c_{1,l}),\a_i,\b_j,\c_l\} =\\ \min\{ d[\e
a_{1,i}b_{1,j}],d[\e'b_{1,j}c_{1,l}]\}$. Similarly, if some of the
lattices $M,N,K$ are the same, we get inequalities involving also
terms of the form $d[\e a_{i,j}]$.
\eff

The following lemma is trivial but it will be very useful.

\blm Let $x,y,z\in\RR$ and $a,b\in\ff$. In each of the following three
cases

(a) $x<y$ and $x+d(a)\geq\min\{ y+d(b),z\}$

(b) $x\leq y$ and $x+d(a)>\min\{ y+d(b),z\}$

(c) $y+d(a)\geq z$

\no we have $\min\{ y+d(b),z\} =\min\{ y+d(ab),z\}$.
\elm
\pf Suppose $x<y$ and $x+d(a)\geq\{ y+d(b),z\}$. Denote $A=\min\{
y+d(b),z\}$. If $y+d(b)\leq z$ then $A=y+d(b)$. Thus $x+d(a)\geq
y+d(b)$ and we also have $x<y$ so $d(a)>d(b)$. By the domination
principle we have $d(b)=d(ab)$ so $\min\{ y+d(b),z\} =\min\{
y+d(ab),z\}$. If $y+d(b)>z$ then $A=z$. So $y+d(a)>x+d(a)\geq z$ and
$y+d(b)>z$. It follows that $y+d(ab)>z$ so $\min\{ y+d(ab),z\}
=z=\min\{ y+d(b),z\}$. The case when $x\leq y$ and $x+d(a)>\min\{
y+d(b),z\}$ is similar.

If $y+d(a)\geq z$ then $\min\{ y+d(b),z\} =\min\{ y+d(b),z,y+d(a)\}$
and  $\min\{ y+d(ab),z\} =\min\{ y+d(ab),z,y+d(a)\}$. But by
domination principle $\min\{ d(b),d(a)\} =\min\{ d(ab),d(a)\}$ so
$\min\{ y+d(b),z\} =\min\{ y+d(ab),z\}$ \qed

{\bf Remark} In the view of 1.3 if $a,b,ab$ are of the form $\e
a_{1,i}b_{1,j}$ or $\e a_{i,j}$ the result above also holds with
$d(a),d(b),d(ab)$ replaced by $d[a],d[b],d[ab]$ because these three
also satisfy the domination principle.

\subsection{A representation lemma for quadratic spaces}

\blm In all cases bellow exactly an even number of the four
statements are true.

(i) $[b_1,\ldots,b_{i-1}]\rep [a_1,\ldots,a_i]$, $[c_1,\ldots,c_{i-1}]\rep
[b_1,\ldots,b_i]$, \\
$[c_1,\ldots,c_{i-1}]\rep [a_1,\ldots,a_i]$ and
$(a_{1,i}b_{1,i},b_{1,i-1}c_{1,i-1})_\p =1$

(ii) $[b_1,\ldots,b_i]\rep [a_1,\ldots,a_{i+1}]$, $[c_1,\ldots,c_{i-1}]\rep
[b_1,\ldots,b_i]$, \\
$[c_1,\ldots,c_{i-1}]\rep [a_1,\ldots,a_i]$ and
$(a_{1,i}b_{1,i},-a_{1,i+1}c_{1,i-1})_\p =1$

(iii) If $[b_1,\ldots,b_{i-1}]\rep [a_1,\ldots,a_i]$,
$[c_1,\ldots,c_{i-2}]\rep 
[b_1,\ldots,b_{i-1}]$, \\
$[c_1,\ldots,c_{i-1}]\rep [a_1,\ldots,a_i]$ and
$(b_{1,i-1}c_{1,i-1},-a_{1,i}c_{1,i-2})_\p =1$.

In particular, if two sentences are true then the remaining two
are equivalent; and if three statements are true then so is the
fourth. \elm \pf Note first that if $ab=cd$ in $\ff/\fs$ then
$[a,b]=[c,d]$ iff $a\rep [c,d]$ i.e. iff $a^2\rep [ac,ad]$ i.e.
iff $(ac,ad)_\p =1$. Note that $(ac,ad)_\p =(ac,-ac)_\p (ac,ad)_\p
=(ac,-cd)_\p$.

The idea of the proof is the fact that, given $V_1,V_2,V_3,V_4$
quadratic spaces with the same dimension and determinant, exactly
an even number of the statements $V_1\ap V_2$, $V_2\ap V_3$,
$V_3\ap V_4$ and $V_4\ap V_1$ are false and so an even number of
them are true. This happens because there are at most 2 classes of
quadratic spaces of given dimension and determinant. For some good
choices of $V_1,V_2,V_3,V_4$ the four statements above are
equivalent to the four statements from the cases (i)-(iii) of our
lemma.

(i) We take $V_1=[a_1,\ldots,a_i,b_i]$,
$V_2=[b_1,\ldots,b_i,a_{1,i}b_{1,i-1}]$,\\
$V_3=[c_1,\ldots,c_{i-1},b_{1,i}c_{1,i-1},a_{1,i}b_{1,i-1}]$,
$V_4=[c_1,\ldots,c_{i-1},a_{1,i}c_{1,i-1},b_i]$.

We have $[b_1,\ldots,b_{i-1}]\rep [a_1,\ldots,a_i]$ iff $[a_1,\ldots,a_i]\ap
[b_1,\ldots,b_{i-1},a_{1,i}b_{1,i-1}]$ i.e. iff $V_1\ap
V_2$. $[c_1,\ldots,c_{i-1}]\rep [b_1,\ldots,b_i]$ iff $[b_1,\ldots,b_i]\ap
[c_1,\ldots,c_{i-1},b_{1,i}c_{1,i-1}]$ i.e. iff $V_2\ap
V_3$. $[c_1,\ldots,c_{i-1}]\rep [a_1,\ldots,a_i]$ iff  $[a_1,\ldots,a_i]\ap
[c_1,\ldots,c_{i-1},a_{1,i}c_{1,i-1}]$ i.e. iff $V_1\ap V_4$. Finally,
$V_3\ap V_4$ iff $[a_{1,i}c_{1,i-1},b_i]\ap
[b_{1,i}c_{1,i-1},a_{1,i}b_{1,i-1}]$. But
this is equivalent to $1=(b_ib_{1,i}c_{1,i-1},b_ia_{1,i}b_{1,i-1})_\p
=(b_{1,i-1}c_{1,i-1},a_{1,i}b_{1,i})_\p$.

(ii) We take $V_1=[a_1,\ldots,a_{i+1}]$,
$V_2=[b_1,\ldots,b_i,a_{1,i+1}b_{1,i}]$,\\
$V_3=[c_1,\ldots,c_{i-1},b_{1,i}c_{1,i-1},a_{1,i+1}b_{1,i}]$,
$V_4=[c_1,\ldots,c_{i-1},a_{1,i}c_{1,i-1},a_{i+1}]$.

We have $[b_1,\ldots,b_i]\rep [a_1,\ldots,a_{i+1}]$ iff $V_1\ap V_2$,
$[c_1,\ldots,c_{i-1}]\rep [b_1,\ldots,b_i]$ iff $V_2\ap V_3$ and
$[c_1,\ldots,c_{i-1}]\rep [a_1,\ldots,a_i]$ iff $V_1\ap V_4$.
Also $V_3\ap V_4$ is equivalent to
$[b_{1,i}c_{1,i-1},a_{1,i+1}b_{1,i}]\ap [a_{1,i}c_{1,i-1},a_{i+1}]$
i.e. to $1=(a_{1,i+1}b_{1,i}a_{i+1},-a_{1,i}c_{1,i-1}a_{i+1})_\p
=(a_{1,i}b_{1,i},-a_{1,i+1}b_{1,i-1})_\p$.

(iii) We take $V_1=[a_1,\ldots,a_i]$,
$V_2=[b_1,\ldots,b_{i-1},a_{1,i}b_{1,i-1}]$,\\
$V_3=[c_1,\ldots,c_{i-2},b_{1,i-1}c_{1,i-2},a_{1,i}b_{1,i-1}]$,
$V_4=[c_1,\ldots,c_{i-1},a_{1,i}c_{1,i-1}]$.

We have $[b_1,\ldots,b_{i-1}]\rep [a_1,\ldots,a_i]$ iff $V_1\ap V_2$,
$[c_1,\ldots,c_{i-2}]\rep [b_1,\ldots,b_{i-1}]$ iff $V_2\ap V_3$ and
$[c_1,\ldots,c_{i-1}]\rep [a_1,\ldots,a_i]$ iff $V_1\ap V_4$. Also $V_3\ap
V_4$ is equivalent to $[b_{1,i-1}c_{1,i-2},a_{1,i}b_{1,i-1}]\ap
[c_{i-1},a_{1,i}c_{1,i-1}]$ i.e. to
$1=(b_{1,i-1}c_{1,i-2}c_{i-1},-c_{i-1}a_{1,i}c_{1,i-1})_\p
=(b_{1,i-1}c_{1,i-1},-a_{1,i}c_{1,i-2})_\p $. \qed

\subsection{The ordered set $(\bbb,\leq )$}

\bdf If $n\in\NN$ let $\bbb_n=\{ (x_1,\ldots,x_n)\in\RR^n~|~x_i\leq
x_{i+2}~\forall~1\leq i\leq n-2\}$ and let $\bbb =\cup_{n\geq
0}\bbb_n$. On $\bbb$ we introduce the following relation: We say that
$(x_1,\ldots,x_m)\leq (y_1,\ldots,y_n)$ if $m\geq n$ and for any $1\leq
i\leq n$ we have either $x_i\leq y_i$ or $1<i<m$ and $x_i+x_{i+1}\leq
y_{i-1}+y_i$.
\edf

Note that the condition $x_i\leq x_{i+2}$ is equivalent to
$x_i+x_{i+1}\leq x_{i+1}+x_{i+2}$ i.e. to the fact that the sequence
$(x_i+x_{i+1})$ is increasing.

\blm Suppose that $x=(x_1,\ldots,x_m),y=(y_1,\ldots,y_n)\in\bbb$ and $x\leq
y$. Then:

(i) $x_i+x_{i+1}\leq y_i+y_{i+1}~\forall~1\leq i\leq n-1$.

(ii) If $x_i+x_{i+1}>y_{i-1}+y_i$ then $x_{i+1}>y_{i-1}$

(iii) If $x_{i+1}\geq y_{i-1}$ then $x_i\leq y_i$
\elm
\pf (i) If $x_i\leq y_i$ and $x_{i+1}\leq y_{i+1}$ then
$x_i+x_{i+1}\leq y_i+y_{i+1}$. If $x_i>y_i$ then $x_i+x_{i+1}\leq
y_{i-1}+y_i\leq y_i+y_{i+1}$. If $x_{i+1}>y_{i+1}$ then $x_i+x_{i+1}\leq
x_{i+1}+x_{i+2}\leq y_i+y_{i+1}$.

(ii) We have $x_i\leq y_i$ since otherwise $x_i+x_{i+1}\leq
y_{i-1}+y_i$. This, together with $x_i+x_{i+1}>y_{i-1}+y_i$, implies
$x_{i+1}>y_{i-1}$.

(iii) Suppose the contrary. Then $x_i+x_{i+1}\leq y_{i-1}+y_i$. This,
together with $x_{i+1}\geq y_{i-1}$, implies $x_i\leq y_i$ \qed

\blm $\leq$ is an order relation on $\bbb$.
\elm
\pf Take $x,y\in\bbb$, $x=(x_1,\ldots,x_m)$, $y=(y_1,\ldots,y_n)$. If $x\leq
y$ and $y\leq x$ then $m\geq n$ and $n\leq n$ so $m=n$. Since $x\leq
y$ and $y\leq x$ we have $x_1\leq y_1$ and $y_1\leq x_1$ so
$x_1=y_1$. Suppose that $x\neq y$. Let $i$ be the smallest index
s.t. $x_i\neq y_i$. Suppose $x_i>y_i$. Since $x\leq y$ we have by
Lemma 1.6(iii), $x_{i+1}<y_{i-1}=x_{i-1}$. Contradiction. Similarly if
$y_i<x_i$. Thus $x=y$.

Take now $z=(z_1,\ldots,z_k)\in\bbb$. Suppose that $x\leq y$ and $y\leq
z$. This implies $m\geq n$ and $n\geq k$ so $m\geq k$. Let $1\leq
i\leq k$. If $x_i\leq y_i$ and $y_i\leq z_i$ then $x_i\leq z_i$. If
$x_i>y_i$ then $1<i<m$ and $x_i+x_{i+1}\leq y_{i-1}+y_i$. Also
$y_{i-1}+y_i\leq z_{i-1}+z_i$ by Lemma 1.6(i) and so $x_i+x_{i+1}\leq
z_{i-1}+z_i$. If $y_i>z_i$ then $1<i<n\leq m$ and $y_i+y_{i+1}\leq
z_{i-1}+z_i$. Also $x_i+x_{i+1}\leq y_i+y_{i+1}$ by Lemma 1.6(i) and
so $x_i+x_{i+1}\leq z_{i-1}+z_i$. Thus $x\leq z$. \qed

\bdf For any $\kappa\in\RR$ we define $\bbb_n(\kappa )=\{
(x_1,\ldots,x_n)\in\bbb_n\mid x_i\leq x_{i+1}+\kappa~\forall~1\leq
i\leq n-1\}$ and $\bbb (\kappa )=\cup_{n\geq 0}\bbb_n(\kappa )$.
\edf

\blm If $x=(x_1,\ldots,x_m),y=(y_1,\ldots,y_n)\in\bbb (\kappa )$ and
$x\leq y$ then:

(i) If $x_{i+1}-x_i\geq \kappa$ or $y_i-y_{i-1}\geq \kappa$ then
$x_i\leq y_i$. 

(ii) If $x_{i+1}-y_i\geq \kappa$ then $x_i\leq y_i$ and, if $i\leq
n-1$, $x_{i+1}\leq y_{i+1}$.
\elm
\pf (i) Suppose that $x_i>y_i$. By Lemma 1.6(iii) we have
$x_{i+1}<y_{i-1}$. If $x_{i+1}-x_i\geq \kappa$ then also $y_{i-1}\leq
y_i+\kappa$ and so $x_i+\kappa\leq x_{i+1}<y_{i-1}\leq y_i+\kappa$ so
$x_i\leq y_i$. If $y_i-y_{i-1}\geq \kappa$ then also $x_i\leq
x_{i+1}+\kappa$ and so $x_i-\kappa\leq x_{i+1}<y_{i-1}\leq y_i-\kappa$
so $x_i\leq y_i$. 

(ii) If $x_i>y_i$ then $y_i+\kappa\leq x_{i+1}<y_{i-1}\leq
y_i+\kappa$. Contradiction. If $x_{i+1}>y_{i+1}$ then
$x_{i+1}-\kappa\leq x_{i+2}<y_i\leq
x_{i+1}-\kappa$. Contradiction. \qed 

\bff If $x=(x_1,\ldots,x_n)\in\RR^n$ we denote $x^\*
=(-x_n,\ldots,-x_1)$. Obviously $(x^\* )^\* =x$.
One can easily see that $x\in\bbb_n$ or $x\in\bbb_n(\kappa )$ iff
$x^\*\in\bbb_n$ resp. $x^\*\in\bbb_n(\kappa )$.

Also if $x,y\in\bbb_n$ then $x\leq y$ iff $y^\*\leq x^\*$.
\eff

\bff If $L\ap\[ a_1,\ldots,a_n\]$ and $R_i=R_i(L)=\ord a_i$ we denote
$\rr (L)=(R_1,\ldots,R_n)$.

Note that $\rr (L)\in\bbb$ because $R_i\leq R_{i+2}$ by the definition
of the good BONGs. Moreover $R_{i+1}-R_i=\ord a_{i+1}/a_i\geq -2e$ (we
have $a_{i+1}/a_i\in\aaa$). Thus $\rr (L)\in\bbb (2e)$.

Also $L^\*\ap\[ a_n\1,\ldots,a_1\1\]$ so $\rr (L^\*
)=(-R_n,\ldots,-R_1)=\rr (L)^\*$.
\eff

\bff If $R_i(L)=R_i$ and $\a_i(L)=\a_i$ we denote
$$\WW (L)=(R_1+\a_1,R_2-\a_1,\ldots,R_{n-1}+\a_{n-1},R_n-\a_{n-1}).$$
By [B3, Lemma 2.2], the sequence $(R_i+\a_i)$ is increasing and
$(-R_{i+1}+\a_i)$ is decreasing so $(R_{i+1}-\a_i)$ is
increasing. This implies that $\WW (L)\in\bbb$.

Since $R_i(L^\* )=-R_{n+1-i}$ and $\a_i(L^\* )=\a_{n-i}$ we have $\WW
(L^\*
)=(-R_n+\a_{n-1},-R_{n-1}-\a_{n-1},\ldots,-R_2+\a_1,-R_1-\a_1)=\WW
(L)^\*$. 
\eff

Throughout this paper we make the following conventions.

{\bf Convention 1} Unless otherwise specified, all BONGs are assumed
to be good. (The only time we have a bad BONG is in the Lemma 9.6.)

{\bf Convention 2} If $L$ is a lattice of rank $n$ and $R_i=R_i(L)$
then we make the assumption that $R_i\ll 0$ if $i<0$ and $R_i\gg 0$ if
$i>n$. 

{\bf Convention 3} Typically when we speak of one lattice we denote it
by $L$ with $L\ap\[ a_1.\ldots,a_n\]$, $R_i(L)=R_i$ and
$\a_i(L)=\a_i$. 

If we speak of two lattices we denote them by $M,N$. When we have a
third lattice we denote it by $K$.  

We write $M\ap\[ a_1,\ldots,a_m\]$, $N\ap\[ b_1,\ldots,b_n\]$ and
$K=\[ c_1,\ldots,c_k\]$, $R_i(M)=R_i$, $R_i(N)=S_i$, $R_i(K)=T_i$,
$\a_i(M)=\a_i$, $\a_i(N)=\b_i$ and $\a_i(K)=\c_i$. 

We will also denote $A_i(M,N)=A_i$, $A_i(N,K)=B_i$ and $A_i(M,K)=C_i$,
where $A_i(M,N)$ is the invariant from Definition 4, \S2.

\section{Main theorem}

In this section we state our main result, the representation Theorem
2.1 and we also give some consequences and equivalent conditions for
the conditions (i)-(iv) of Theorem 2.1.

The main results we use in this section are:

The fact that $(R_i+\a_i)$ increases and $(-R_{i+1}+\a_i)$ decreses,
which can be written as
$$\a_i\geq R_j-R_i+\a_j\text{ if }i\geq j\text{ and }\a_i\geq
R_{i+1}-R_{j+1}+\a_j\text{ if }i\leq j.$$

The inequality $d[-a_{i,i+1}]\geq R_i-R_{i+1}+\a_i$ from Remark 1.1. 

The domination principle and Lemma 1.4 for the $d[\e a_{1,i}b_{1,j}]$
and $d[\e a_{i,j}]$ invariants. 

The consequences of Lemmas 1.6 and 1.8 following from Theorem 2.1(i)
(see 2.2). 

The representation Lemma 1.5. 
\vskip 3mm

Throughout this section $M\ap\[ a_1,\ldots,a_m\]$ and $N\ap\[
b_1,\ldots,b_n\]$ relative to good BONGs, $R_i=R_i(M)=\ord a_i$,
$S_i=R_i(M)=\ord b_i$, $\a_i=\a_i(M)$, $\b_i=\a_i (N)$.

\bdf For any $1\leq i\leq\min\{m-1,n\}$ we define
$$A_i:=\min\{ (R_{i+1}-S_i)/2+e, R_{i+1}-S_i+d[-a_{1,i+1}b_{1,i-1}],
R_{i+1}+R_{i+2}-S_{i-1}-S_i+d[a_{1,i+2}b_{1,i-2}]\}.$$
If $n\leq m-2$ we define
$$S_{n+1}+A_{n+1}:=\min\{ R_{n+2}+d[-a_{1,n+2}b_{1,n}],
R_{n+2}+R_{n+3}-S_n+d[a_{1,n+3}b_{1,n-1}]\}.$$

(The terms that don't make sense are ignored. I.e. if $i=1$ or $m-1$
then we ignore $R_{i+1}+R_{i+2}-S_{i-1}-S_i+d[a_{1,i+2}b_{1,i-2}]$
from the definition of $A_i$; and if $n=m-2$ then we ignore
$R_{n+2}+R_{n+3}-S_n+d[a_{1,n+3}b_{1,n-1}]$ from the definition of
$S_{n+1}+A_{n+1}$.)

In the view of Lemma 1.2 $A_i$ doesn't depend on BONGs so we can write
$A_i=A_i(M,N)$. Same for $S_{n+1}+A_{n+1}$.
\edf

The definition of $S_{n+1}+A_{n+1}$ when $n\leq m-2$ is justified by
Convention 1. When we take $i=n+1$ in the definition of $A_i$ and add
$S_{n+1}$ we get $S_{n+1}+A_{n+1}=\min\{ (R_{n+2}+S_{n+1})/2+e,
R_{n+2}+d[-a_{1,n+2}b_{1,n}],
R_{n+2}+R_{n+3}-S_n+d[a_{1,n+3}b_{1,n-1}]\}$. But
$(R_{n+2}+S_{n+1})/2+e\gg 0$, as $S_{n+1}\gg 0$, so it can be
removed. 

\btm Assume that $FN\rep FM$. Then $N\rep M$ iff:

(i) For any $1\leq i\leq n$ we have either $R_i\leq S_i$ or $1<i<m$
and $R_i+R_{i+1}\leq S_{i-1}+S_i$.

(ii) For any $1\leq i\leq\min\{ m-1,n\}$ we have
$d[a_{1,i}b_{1,i}]\geq A_i$.

(iii) For any $1<i\leq\min\{ m-1,n+1\}$ s.t. $R_{i+1}>S_{i-1}$
and $A_{i-1}+A_i>2e+R_i-S_i$ we have $[b_1,\ldots,b_{i-1}]\rep
[a_1,\ldots,a_i]$. 

(iv) If $1<i\leq\min\{ m-2,n+1\}$ s.t. $S_i\geq R_{i+2}>S_{i-1}+2e\geq
R_{i+1}+2e$ then $[b_1,\ldots,b_{i-1}]\rep [a_1,\ldots,a_{i+1}]$. (If
$i=n+1$ we ignore the condition $S_i\geq R_{i+2}$.)
\etm

Note that if $n\leq m-2$ and $i=n+1$ then condition from (iii) can be
written as $A_n+S_{n+1}+A_{n+1}>2e+R_{n+1}$ so it makes sense since
$S_{n+1}+A_{n+1}$ is defined, although $S_{n+1}$ and $A_{n+1}$ are
not. 

Also note that $R_{i+1}-S_i+d[-a_{1,i+1}b_{1,i-1}]\geq A_i$ and
$R_i-S_{i-1}+d[-a_{1,i}b_{1,i-2}]\geq A_{i-1}$, so a necessary
condition for the inequality from (iii) to hold is that
$(R_i-S_{i-1}+d[-a_{1,i}b_{1,i-2}])+(R_{i+1}-S_i+d[-a_{1,i+1}b_{1,i-1}])>
2e+R_i-S_i$, i.e. that
$d[-a_{1,i}b_{1,i-2}]+d[-a_{1,i+1}b_{1,i-1}]>2e+S_{i-1}-R_{i+1}$. In
Lemma 2.16 we will show that if $R_{i+1}>S_{i-1}$ and the lattices
$M,N$ satisfy the conditions (i) and (ii) of Theorem 2.1, then this
condition is also sufficient. Therefore the condition (iii) of Theorem
2.1 may be replaced by:
\medskip

{\it (iii') For any $1<i\leq\min\{ m-1,n+1\}$ s.t. $R_{i+1}>S_{i-1}$
and $d[-a_{1,i}b_{1,i-2}]+d[-a_{1,i+1}b_{1,i-1}]>2e+S_{i-1}-R_{i+1}$
we have $[b_1,\ldots,b_{i-1}]\rep [a_1,\ldots,a_i]$.} 
\medskip

Before starting the proof we make some rematks regarding conditions
(i)-(iv) of the theorem.

\bff Condition (i) simply means that $\rr (M)\leq\rr (N)$ in $\bbb$ so
we can use Lemmas 1.6 and 1.8 (with $\kappa =2e$). So if (i) is
satisfied then $R_i+R_{i+1}\leq S_i+S_{i+1}$ if
$R_i+R_{i+1}>S_{i-1}+S_i$ then $R_{i+1}>S_{i-1}$ and if
$R_{i+1}>S_{i-1}$ then $R_i\leq S_i$. Also if $R_{i+1}-R_i\geq 2e$ or
$S_i-S_{i-1}\geq 2e$ then $R_i\leq S_i$ and if $R_{i+1}-S_i\geq 2e$
then $R_i\leq S_i$ and $R_{i+1}\leq S_{i+1}$. 
\eff

\subsection*{Dualization}

\bff Let's see how the properties (i)-(iv) behave under
dualization when $m=n$. We have $N^\*\ap\[ a_1^\*,\ldots,a_n^\*\]$ and
$M^\*\ap\[ b_1^\*,\ldots,b_n^\*\]$ with $a_i^\* =b_{n+1-i}\1$ and $b_i^\*
=a_{n+1-i}\1$. 

We have $R_i^\* :=R_i(N^\* )=-S_{n+1-i}$, $S_i^\* :=R_i(M^\*
)=-R_{n+1-i}$, $\a_i^\* :=\a_i(N^\* )=\b_{n-i}$ and $\b_i^\*
:=\a_i(M^\* )=\a_{n-i}$. 

Since $a_{1,n}b_{1,n}\in\fs$ we have $\e a_{1,i}^\* b_{1,j}^\* =\e
b_{n+1-i,n}a_{n+1-j,n}=\e a_{1,n-j}b_{1,n-i}$ (in $\ff/\fs$). Since
also $\a_i^\* =\b_{n-i}$ and $\b_i^\* =\a_{n-i}$ we get $d[\e
a_{1,i}^\* b_{1,j}^\* ]=d[\e a_{1,n-j}b_{1,n-i}]$. Similarly,
$d[\e a_{i,j}^\* ]=d[\e b_{n-j+1,n-i+1}]$ and $d[\e b_{i,j}^\* ]=d[\e
a_{n-j+1,n-i+1}]$. 

As a consequence one may see that $A_i^\* :=A_i(N^\*,M^\* )=A_{n-i}$
for $1\leq i\leq n-1$.

Also condition $[b_1^\*,\ldots,b_j^\* ]\rep [a_1^\*,\ldots,a_i^\* ]$
means $[a_{n+1-j},\ldots,a_n]\rep [b_{n+1-i},\ldots,b_n]$ which is
equivalent to $[b_1,\ldots,b_{n-i}]\rep [a_1,\ldots,a_{n-j}]$. (We have
$[a_1,\ldots,a_n]\ap [b_1,\ldots,b_n]$.)
\eff

\bff
In many cases, given a formula or a statement for $M,N$ at some index
$i$, we will consider the same formula or statement for $N^\*,M^\*$ at
index $n+1-i$. When we do so we have:

$R_{i+k},S_{i+k}$ become $-S_{i-k},-R_{i-k}$.

$d[\e a_{i+k,i+l}]$, $d[\e b_{i+k,i+l}]$ become $d[\e b_{i-l,i-k}]$,
$d[\e a_{i-l,i-k}]$.

$\a_{i+k},\b_{i+k},A_{i+k}$ become $\b_{i-k-1},\a_{i-k-1},A_{i-k-1}$

$d[\e a_{1,i+k}b_{1,i+l}]$ becomes $d[\e a_{1,i-l-1}b_{1,i-k-1}]$.

$[b_1,\ldots,b_{i+k}]\rep [a_1,\ldots,a_{i+l}]$ becomes
$[b_1,\ldots,b_{i-k-1}]\rep [a_1,\ldots,a_{i-l-1}]$.

Sometimes instead of $n+1-i$ we take $n-i$. In this case all indices
in the dual formulas and statements described above are increased by
$1$. 

Note that even if $m>n$ we can still use the duality argument. The
reader may repeat the same reasoning for the dual statement according
to the rules above.
\eff

\bff
One can easily see that conditions (i)-(iv) for $M,N$ are
equivalent to those for $N^\* ,M^\*$. (Conditions (i) and (iii)
for $N^\* ,M^\*$ at an index $i$ correspond to the similar
conditions for $M,N$ at $n+1-i$. For (ii) and (iv) they correspond
to the similar conditions at $n-i$.) \eff

\subsection*{The invariants $A'_i$, $\( A_i$, $\( A'_i$ and condition
2.1(ii)} 

\bdf We define $\a'_i:=R_{i+1}-R_i+d[-a_{i,i+1}]$ and
$$A'_i:=\min\{ R_{i+1}-S_i+d[-a_{1,i+1}b_{1,i-1}],
R_{i+1}+R_{i+2}-S_{i-1}-S_i+d[a_{1,i+2}b_{1,i-2}]\}.$$

We have $A_i=\min\{ (R_{i+1}-S_i)/2+e,A'_i\}$ and, in the view of 1.1,
$\a_i=\min\{ (R_{i+1}-R_i)/2+e,\a'_i\}$.
\edf

\bff The following numbers are $\geq A'_i$: 

$R_{i+1}-S_i+\a_{i+1}$ and $R_{i+1}-S_i+\b_{i-1}$

$R_{i+1}+R_{i+2}-S_{i-1}-S_i+\a_{i+2}$ and
$R_{i+1}+R_{i+2}-S_{i-1}-S_i+\b_{i-2}$

$(R_{i+1}+R_{i+2})/2-S_i+e$ and $R_{i+1}-(S_{i-1}+S_i)/2+e$. 

In particular, they are $\geq A_i$. (The numbers that don't make sense
are ignored.)

Indeed $\a_{i+1},\b_{i-1}\geq d[-a_{1,i+1}b_{1,i-1}]$ and
$\a_{i+2},\b_{i-2}\geq d[a_{1,i+2}b_{1,i-2}]$ so the first four
numbers are $\geq A'_i$. Also
$(R_{i+1}+R_{i+2})/2-S_i+e=R_{i+1}-S_i+(R_{i+2}-R_{i+1})/2+e\geq
R_{i+1}-S_i+\a_{i+1}\geq A'_i$ and
$R_{i+1}-(S_{i-1}+S_i)/2+e=R_{i+1}-S_i+(S_i-S_{i-1})/2+e\geq
R_{i+1}-S_i+\b_{i-1}\geq A'_i$.

If $M,N$ satisfy the 2.1(i) then also
$(R_{i+1}-R_i)/2+e,(S_{i+1}-S_i)/2+e\geq A_i$. 

Indeed, we have either $R_i\leq S_i$ so $(R_{i+1}-R_i)/2+e\geq
(R_{i+1}-S_i)/2+e\geq A_i$ or $R_i+R_{i+1}\leq S_{i-1}+S_i$ so
$(R_{i+1}-R_i)/2+e\geq R_{i+1}-(S_{i-1}+S_i)/2+e\geq A_i$. Similarly
we have $R_{i+1}\leq S_{i+1}$ so $(S_{i+1}-S_i)/2+e\geq
(R_{i+1}-S_i)/2+e\geq A_i$ or $R_{i+1}+R_{i+2}\leq S_i+S_{i+1}$ so
$(S_{i+1}-S_i)/2+e\geq (R_{i+1}+R_{i+2})/2-S_i+e\geq A_i$.

Also note that $\a_i,\b_i\geq d[a_{1,i}b_{1,i}]$, i.e. if 2.1(ii) is
satisfied at $i$ then $\a_i,\b_i\geq A_i$.
\eff

\blm In the term $R_{i+1}+R_{i+2}-S_{i-1}-S_i+d[a_{1,i+2}b_{1,i-2}]$
from the definition of $A_i$ and $A'_i$ $d[a_{1,i+2}b_{1,i-2}]$ can be
replaced as follows: 

(i) If $R_{i+1}\geq S_{i-1}$ it can be replaced by
$d[-a_{1,i}b_{1,i-2}]$. 

(ii) If $R_{i+2}\geq S_i$ it can be replaced by
$d[-a_{1,i+2}b_{1,i}]$. 

(iii) If both $R_{i+1}\geq S_{i-1}$ and $R_{i+2}\geq S_i$ it can be
replaced by $d[a_{1,i}b_{1,i}]$.
\elm
\pf Since $A_i=\min\{ (R_{i+1}-S_i)/2+e,A'_i\}$, it is enough to prove
the statements (i)-(iii) for $A'_i$. 

(i) If $R_{i+1}\geq S_{i-1}$ then
$R_{i+1}+R_{i+2}-S_{i-1}-S_i+d[-a_{1+1,i+2}]\geq
R_{i+2}-S_i+d[-a_{i+1,i+2}]\geq R_{i+1}-S_i+\a_{i+1}\geq
R_{i+1}-S_i+d[-a_{1,i+1}b_{1,i-1}]$ so by Lemma 1.4(c) we have $\min\{
R_{i+1}-S_i+d[-a_{1,i+1}b_{1,i-1}],
R_{i+1}+R_{i+2}-S_{i-1}-S_i+d[-a_{1,i}b_{1,i-2}]\} =\min\{
R_{i+1}-S_i+d[-a_{1,i+1}b_{1,i-1}],
R_{i+1}+R_{i+2}-S_{i-1}-S_i+d[a_{1,i+2}b_{1,i-2}]\} =A'_i$

(ii) is similar but this time $R_{i+2}\geq S_i$ implies
$R_{i+1}+R_{i+2}-S_{i-1}-S_i+d[-b_{i-1,i}]\geq
R_{i+1}-S_{i-1}+d[-b_{i-1,i}]\geq R_{i+1}-S_i+\b_{i-1}\geq
R_{i+1}-S_i+d[-a_{1,i+1}b_{1,i-1}]$.

(iii) Since both $R_{i+1}\geq S_{i-1}$ and $R_{i+2}\geq S_i$ hold we
can use procedures from both (i) and (ii). So we replace
$d[a_{1,i+2}b_{1,i-2}]$ first by $d[-a_{1,i}b_{1,i-2}]$ and then by
$d[a_{1,i}b_{1,i}]$. \qed

\bdf Let $1\leq i\leq\min\{ m-1,n\}$. If $R_{i+1}+R_{i+2}>S_{i-1}+S_i$
we define
$$\( A'_i=\min\{ R_{i+1}-S_i+d[-a_{1,i+1}b_{1,i-1}],
2R_{i+1}-S_{i-1}-S_i+\a_{i+1}, R_{i+1}+R_{i+2}-2S_i+\b_{i-1}\}$$
and $\( A_i=\min\{ (R_{i+1}-S_i)/2+e,\( A'_i\}$.

(If $i\in\{ 1,m-1\}$ then $2R_{i+1}-S_{i-1}-S_i+\a_{i+1}$ and
$R_{i+1}+R_{i+2}-2S_i+\b_{i-1}$ don't make sense so they are
ignored. So in this case we have $\(
A'_i=A'_i=R_{i+1}-S_i+d[-a_{1,i+1}b_{1,i-1}]$ and $\( A_i=A_i=\min\{
(R_{i+1}-S_i)/2+e, R_{i+1}-S_i+d[-a_{1,i+1}b_{1,i-1}]$.)
\edf

\bff{\bf Remarks} By Convention 2 if $i\in\{ 1,m-1\}$ then
$R_{i+1}+R_{i+2}>S_{i-1}+S_i$. 

If $R_{i+1}\geq S_{i-1}$ then $2R_{i+1}-S_{i-1}-S_i+\a_{i+1}\geq
R_{i+1}-S_i+\a_{i+1}\geq R_{i+1}-S_i+d[-a_{1,i+1}b_{1,i-1}]$ so it can
be ignored in the formulas for $\( A_i$ and $\( A'_i$. Similarly if
$R_{i+2}\geq S_i$ then $R_{i+1}+R_{i+2}-2S_i+\b_{i-1}\geq
R_{i+1}-S_i+\b_{i-1}\geq R_{i+1}-S_i+d[-a_{i+1}b_{1,i-1}]$ so it can
be ignored. 

Since $R_{i+1}+R_{i+2}>S_{i-1}+S_i$ we always have $R_{i+1}>S_{i-1}$
or $R_{i+2}>S_i$ so at least one of $2R_{i+1}-S_{i-1}-S_i+\a_{i+1}$
and $R_{i+1}+R_{i+2}-2S_i+\b_{i-1}$ can be removed.
\eff

\blm Suppose that $R_{i+1}+R_{i+2}>S_{i-1}+S_i$. Then if
$d[a_{1,i}b_{1,i}]\geq A_i$ or $d[a_{1,i}b_{1,i}]\geq\( A_i$ we have
$A_i=\( A_i$ and $A'_i=\( A'_i$. 

In particular, $d[a_{1,i}b_{1,i}]\geq A_i$ iff
$d[a_{1,i}b_{1,i}]\geq\( A_i$. 
\elm
\pf We use Lemma 1.4(a) with $x=0$, $y=R_{i+1}+R_{i+2}-S_{i-1}-S_i>0$ and
$z=\min\{ (R_{i+1}-S_i)/2+e,
R_{i+1}-S_i+d[-a_{1,i+1}b_{1,i-1}]\}$. Note that
$(R_{i+1}+R_{i+2})/2-S_i+e\geq R_{i+1}-S_i+\a_{i+1}\geq z$ and
$R_{i+1}-(S_{i-1}+S_i)/2+e\geq R_{i+1}-S_i+\b_{i-1}\geq z$ (see the proof
of 2.6) We have three cases:

a. $R_{i+1}\geq S_{i-1}$ and $R_{i+2}\geq S_i$. We have $\( A_i=z$ by
Remark 2.8 and $A_i=\min\{ y+d[a_{1,i}b_{1,i}], z\}$ by Lemma
2.7(iii). If $d[a_{1,i}b_{1,i}]\geq A_i$ then $y+d[a_{1,i}b_{1,i}]>
A_i=\min\{ y+d[a_{1,i}b_{1,i}],z\}$ so $A_i=z=\( A_i$. Same happens if
$d[a_{1,i}b_{1,i}]\geq\( A_i$ since $\( A_i\geq A_i$. 

b. $R_{i+1}<S_{i-1}$. This implies  $R_{i+2}>S_i$. So $A_i=\min\{
y+d[-a_{1,i+2}b_{1,i}],z\}$ by Lemma 2.7(ii) and $\( A_i=\min\{
z,2R_{i+1}-S_{i-1}-S_i+\a_{i+1}\}$ by Remark 2.8. We have
$\a_{i+1}=\min\{ (R_{i+2}-R_{i+1})/2+e,
R_{i+2}-R_{i+1}+d[-a_{i+1,i+2}]\}$ so $\( A_i=\min\{ z,
2R_{i+1}-S_{i-1}-S_i+(R_{i+2}-R_{i+1})/2+e,
R_{i+1}+R_{i+2}-S_{i-1}-S_i+d[-a_{i+1,i+2}]\} =\min\{
y+d[-a_{i+1,i+2}],z,R_{i+1}+(R_{i+1}+R_{i+2})/2-S_{i-1}-S_i+e\}$. But
$R_{i+1}+(R_{i+1}+R_{i+2})/2-S_{i-1}-S_i+e>
R_{i+1}+(S_{i-1}+S_i)/2-S_{i-1}-S_i+e=R_{i+1}-(S_{i-1}+S_i)/2+e\geq z$
so it can be ignored. Thus $\( A_i=\min\{ y+d[-a_{i+1,i+2}],z\}$. Since
$y>0=x$, by Lemma 1.4(a), if $d[a_{1,i}b_{1,i}]\geq A_i=\min\{
y+d[-a_{1,i+2}b_{1,i}],z\}$ then $A_i=\min\{ y+d[-a_{i+1,i+2}],z\}
=\( A_i$. Similarly if $d[a_{1,i}b_{1,i}]\geq\( A_i$.

c. $R_{i+2}<S_i$. Similar with case b.. (If $m=n$ it follows from
b. by duality at index $n-i$. See 2.4.) 

So $A_i=\( A_i$. If $A_i=\( A_i\neq (R_{i+1}-S_i)/2+e$
then $A_i=A'_i$ and $\( A_i=\( A'_i$ so $A'_i=\( A'_i$. If
$A_i=(R_{i+1}-S_i)/2+e$ then $(R_{i+1}+R_{i+2})/2-S_i+e,
R_{i+1}-(S_{i-1}+S_i)/2+e\geq A'_i\geq A_i=(R_{i+1}-S_i)/2+e$ so
$R_{i+1}\geq S_{i-1}$ and $R_{i+2}\geq S_i$. By Remark 2.8 $\(
A'_i=R_{i+1}-S_i+d[-a_{1,i+1}b_{1,i-1}]$ and by Lemma 2.7(iii)
$A'_i=\min\{ R_{i+1}-S_i+d[-a_{1,i+1}b_{1,i-1}],
R_{i+1}+R_{i+2}-S_{i-1}-S_i+d[a_{1,i}b_{1,i}]\}$. But
$d[a_{1,i}b_{1,i}]\geq A_i=(R_{i+1}-S_i)/2+e\geq
(S_{i-1}-S_i)/2+e$ and $R_{i+1}+R_{i+2}-S_{i-1}-S_i>0$ so
$R_{i+1}+R_{i+2}-S_{i-1}-S_i+d[a_{1,i}b_{1,i}]>
(R_{i+1}+R_{i+2}-S_{i-1}-S_i)/2+(S_{i-1}-S_i)/2+e=
(R_{i+1}+R_{i+2})/2-S_i+e\geq A'_i$ (see 2.6). So
$A'_i=R_{i+1}-S_i+d[-a_{1,i+1}b_{1,i-1}]=\( A'_i$. \qed

\bco Suppose that $M,N$ satisfy the condition (i) of the main
theorem. If $R_{i+1}-S_i>2e$ then $d[a_{1,i}b_{1,i}]\geq A_i$ iff
$a_{1,i}b_{1,i}\in\fs$.
\eco
\pf If $2\leq i\leq m-2$ then $R_{i+2}\geq R_{i+1}-2e\geq S_i$ and
$R_{i+1}>S_i+2e\geq S_{i-1}$. This implies that
$R_{i+1}+R_{i+2}>S_{i-1}+S_i$ and by Remark 2.8 $\( A_i=\min\{
(R_{i+1}-S_i)/2+e, R_{i+1}-S_i+d[-a_{1,i+1}b_{1,i-1}]\}$. Same happens
if $i\in\{ 1,m-1\}$, by definition.

But $R_{i+1}-S_i>2e$ so $R_{i+1}-S_i+d[-a_{1,i+1}b_{1,i-1}]\geq
R_{i+1}-S_i>(R_{i+1}-S_i)/2+e$ so $\( A_i=(R_{i+1}-S_i)/2+e$. By
Lemma 2.9 we have $d[a_{1,i}b_{1,i}]\geq A_i$ iff
$d[a_{1,i}b_{1,i}]\geq\( A_i=(R_{i+1}-S_i)/2+e$. Note that
$R_{i+1}-S_i>2e$ implies by 2.2 that $R_i\leq S_i$ and
$R_{i+1}\leq S_{i+1}$. Hence both $R_{i+1}-R_i$ and $S_{i+1}-S_i$ are
$\geq R_{i+1}-S_i>2e$. By [B3, Lemma 2.7(ii)]
$\a_i=(R_{i+1}-R_i)/2+e\geq (R_{i+1}-S_i)/2+e$ and
$\b_i=(S_{i+1}-S_i)/2+e\geq (R_{i+1}-S_i)/2+e$. Since
$d[a_{1,i}b_{1,i}]=\min\{ d(a_{1,i}b_{1,i}),\a_i,\b_i\}$ we have
$d[a_{1,i}b_{1,i}]\geq (R_{i+1}-S_i)/2+e$ iff $d(a_{1,i}b_{1,i})\geq
(R_{i+1}-S_i)/2+e$. But   $(R_{i+1}-S_i)/2+e>2e$ so this is equivalent
to $a_{1,i}b_{1,i}\in\fs$. \qed 

\blm If $R_{i+1}+R_{i+2}\leq S_{i-1}+S_i$ then:

(i) If $R_{i+1}\leq S_{i-1}$ then $d[a_{1,i}b_{1,i}]\geq A_i$ iff
$R_{i+1}-R_{i+2}+\a_{i+1}\geq A_i$.

(ii) If $R_{i+2}\leq S_i$ then $d[a_{1,i}b_{1,i}]\geq A_i$ iff
$S_{i-1}-S_i+\b_{i-1}\geq A_i$.
\elm
\pf Since $R_{i+1}+R_{i+2}\leq S_{i-1}+S_i$ we have $d[a_{1,i+2}b_{1,i-2}]\geq
R_{i+1}+R_{i+2}-S_{i-1}-S_i+d[a_{1,i+2}b_{1,i-2}]\geq A_i$.

(i) We have $d[-b_{i-1,i}]\geq S_{i-1}-S_i+\b_{i-1}\geq
R_{i+1}-S_i+\b_{i-1}\geq A_i$ which, together with $d[a_{1,i+2}b_{1,i-2}]\geq
A_i$, implies $d[-a_{1,i+2}b_{1,i}]\geq A_i$. Thus
$d[a_{1,i}b_{1,i}]\geq A_i$ iff $d[-a_{i+1,i+2}]\geq A_i$. On the
other hand $\a_{i+1}=\min\{ (R_{i+2}-R_{i+1})/2+e,
R_{i+2}-R_{i+1}+d[-a_{i+1,i+2}]\}$ so $R_{i+1}-R_{i+2}+\a_{i+1}=\min\{
(R_{i+1}-R_{i+2})/2+e, d[-a_{i+1,i+2}]\}$. But $R_{i+1}+R_{i+2}\leq
S_{i-1}+S_i$ implies $(R_{i+1}-R_{i+2})/2+e\geq
R_{i+1}-(S_{i-1}+S_i)/2+e\geq A_i$ so $d[-a_{i+1,i+2}]\geq A_i$ iff
$R_{i+1}-R_{i+2}+\a_{i+1}\geq A_i$. 

(ii) is similar. If $m=n$ it follows by duality at index $n-i$. \qed

\subsection*{Essential indices}

\bdf An index $1\leq i\leq\min\{ m,n+1\}$ is called essential if
$R_{i+1}>S_{i-1}$ and  $R_{i+1}+R_{i+2}>S_{i-2}+S_{i-1}$. (By
Convention 2 the conditions that don't make sense are ignored. Thus
$1$ is always essential and if $R_3>S_1$ then $2$ is essential. If
$n\geq m-1$ then $m$ is essential. If $n\geq m-2$ and $R_m>S_{m-2}$
then $m-1$ is essential.)
\edf

{\bf Remark} If $m=n$ then $i$ is essential for $M,N$ iff $n+1-i$ is
essential for $N^\*,N^\*$.

\blm If both $i$ and $i+1$ are not essential then condition 2.1(ii) is
vacuous at index $i$. 
\elm
\pf If $R_{i+1}+R_{i+2}>S_{i-1}+S_i$ then we have both
$R_{i+1}+R_{i+2}>S_{i-2}+S_{i-1}$ and $R_{i+2}+R_{i+3}>S_{i-1}+S_i$. We also
have $R_{i+1}>S_{i-1}$ or $R_{i+2}>S_i$. Therefore either $i$ or $i+1$ is
essential. Thus we may assume that $R_{i+1}+R_{i+2}\leq S_{i-1}+S_i$. This
implies $R_{i+1}\leq S_{i-1}$ or $R_{i+2}\leq S_i$.

If $R_{i+1}\leq S_{i-1}$ then by Lemma 2.11(i) condition 2.1(ii) is
equivalent to $R_{i+1}-R_{i+2}+\a_{i+1}\geq A_i$. Since $i+1$ is not
essential we have $R_{i+2}\leq S_i$ or $R_{i+2}+R_{i+3}\leq
S_{i-1}+S_i$. In the first case $R_{i+1}-R_{i+2}+\a_{i+1}\geq
R_{i+1}-S_i+\a_{i+1}\geq A_i$ and in the second case
$R_{i+1}-R_{i+2}+\a_{i+1}\geq R_{i+1}-R_{i+3}+\a_{i+2}\geq
R_{i+1}+R_{i+2}-S_{i-2}-S_{i-1}+\a_{i+2}\geq A_i$.

If $R_{i+2}\leq S_i$ then by Lemma 2.11(ii) condition 2.1(ii) is
equivalent to $S_{i-1}-S_i+\b_{i-1}\geq A_i$. Since $i$ is not
essential we have $R_{i+1}\leq S_{i-1}$ or $R_{i+1}+R_{i+2}\leq
S_{i-2}+S_{i-1}$. In the first case $S_{i-1}-S_i+\b_{i-1}\geq
R_{i+1}-S_i+\b_{i-1}\geq A_i$ and in the second case
$S_{i-1}-S_i+\b_{i-1}\geq S_{i-2}-S_i+\b_{i-2}\geq
R_{i+1}+R_{i+2}-S_{i-1}-S_i+\b_{i-2}\geq A_i$. \qed

\blm If $i$ is not essential then condition 2.1(iii) is vacuous at
index $i$.
\elm
\pf Suppose that $R_{i+1}>S_{i-1}$ and $A_{i-1}+A_i>2e+R_i-S_i$. By
2.6 $(R_{i+1}+R_{i+2})/2-S_i+e+ R_i-(S_{i-2}+S_{i-1})/2+e\geq
A_i+A_{i-1}>2e+R_i-S_i$ so $R_{i+1}+R_{i+2}>S_{i-2}+S_{i-1}$. Thus $i$
is essential. \qed

\subsection*{Condition 2.1(iii)}

\blm Assume that $M,N$ satisfy 2.1(i) and (ii) and $A_i\neq A'_i$ for
some $1\leq i\leq\min\{ n,m-1\}$. Then the conditions
$R_{j+1}>S_{j-1}$ and $A_{j-1}+A_j>2e+R_j-S_j$ from (iii) are
fullfiled for $j=i,i+1$. (If $i=1$ then only for $j=i+1=2$
and if $i=m-1$ then only for $j=i=m-1$.)
\elm
\pf We only prove our statement for $j=i$. The case $j=i+1$ is
similar. If $m=n$ it follows by duality at index $n-i$. 

Condition $A_i\neq A'_i$ means $A_i=(R_{i+1}-S_i)/2+e<A'_i$. This
implies $(R_{i+1}+R_{i+2})/2-S_i+e\geq A'_i>(R_{i+1}-S_i)/2+e$ so
$R_{i+2}>S_i$ and $R_{i+1}-(S_{i-1}+S_i)/2+e\geq
A'_i>(R_{i+1}-S_i)/2+e$ so $R_{i+1}>S_{i-1}$. Also
$R_{i+1}-S_i+d[-a_{1,i+1}b_{1,i-1}]\geq A'_i>(R_{i+1}-S_i)/2+e$ so
$d[-a_{1,i+1}b_{1,i-1}]>e-(R_{i+1}-S_i)/2$. It follows that
$\a_{i+1},\b_{i-1}>e-(R_{i+1}-S_i)/2$. 

$R_{i+1}>S_{i-1}$ implies $R_i\leq S_i$ and, by Lemma 2.7(ii)
$A_{i-1}=\min\{ (R_i-S_{i-1})/2+e, R_i-S_{i-1}+d[-a_{1,i}b_{1,i-2}],
R_i+R_{i+1}-S_{i-2}-S_{i-1}+d[-a_{1,i+1}b_{1,i-1}]\}$. 

Since $R_{i+1}>S_{i-1}$ and $R_i\leq S_i$ we have
$(R_{i+1}-S_i)/2+e+(R_i-S_{i-1})/2+e>2e+R_i-S_i$.

We have $d[-b_{i-1,i}]\geq S_{i-1}-S_i+\b_{i-1}\geq
S_{i-1}-S_i+e-(R_{i+1}-S_i)/2$ and $d[a_{1,i}b_{1,i}]\geq
A_i=(R_{i+1}-S_i)/2+e>S_{i-1}-S_i+e-(R_{i+1}-S_i)/2$. (Because
$R_{i+1}>S_{i-1}$. By domination principle
$d[-a_{1,i}b_{1,i-2}]>S_{i-1}-S_i+e-(R_{i+1}-S_i)/2$ and so
$(R_{i+1}-S_i)/2+e+R_i-S_{i-1}+d[-a_{1,i}b_{1,i-2}]>2e+R_i-S_i$. 

Finally, $d[-a_{1,i+1}b_{1,i-1}]>e-(R_{i+1}-S_i)/2$ so
$(R_{i+1}-S_i)/2+e+R_i+R_{i+1}-S_{i-2}-S_{i-1}+d[-a_{1,i+1}b_{1,i-1}]>
2e+R_i+R_{i+1}-S_{i-2}-S_{i-1}>2e+R_i-S_i$. (We have $R_{i+1}>S_{i-1}$
and $S_i\geq S_{i-2}$.)

In consequence
$A_{i-1}+A_i=(R_{i+1}-S_i)/2+e+A_{i-1}>2e+R_i-S_i$. \qed 

\bff If $M,N$ satisfy 2.1(i) and (ii) then $A_i\neq A'_i$ iff
$d[-a_{1,i+1}b_{1,i-1}]>e-(R_{i+1}-S_i)/2$. The necessity is
obvious: If $A_i\neq A'_i$ then
$R_{i+1}-S_i+d[-a_{1,i+1}b_{1,i-1}]\geq A'_i>(R_{i+1}-S_i)/2+e$ so
$d[-a_{1,i+1}b_{1,i-1}]>e-(R_{i+1}-S_i)/2$.

Conversely, assume $d[-a_{1,i+1}b_{1,i-1}]>e-(R_{i+1}-S_i)/2$. It
follows that
$R_{i+1}-S_i+d[-a_{1,i+1}b_{1,i-1}]>(R_{i+1}-S_i)/2+e$. Also 
$(R_{i+2}-R_{i+1})/2+e\geq\a_{i+1}>e-(R_{i+1}-S_i)/2$ and
$(S_i-S_{i-1})/2+e\geq\b_{i-1}>e-(R_{i+1}-S_i)/2$. So
$R_{i+2}>S_i$ and $R_{i+1}>S_{i-1}$. Thus $\( A'_i$ is defined and we
have $\( A'_i=R_{i+1}-S_i+d[-a_{1,i+1}b_{1,i-1}]$. By Lemma 2.9
$A'_i=\( A'_i=R_{i+1}-S_i+d[-a_{1,i+1}b_{1,i-1}]>(R_{i+1}-S_i)/2+e$
and so $A_i=(R_{i+1}-S_i)/2+e<A'_i$. \eff

\blm If If $M,N$ satisfy 2.1(i) and (ii) and $R_{i+1}>S_{i-1}$, then
$A_{i-1}+A_i>2e+R_i-S_i$ is equivalent to
$d[-a_{1,i}b_{1,i-2}]+d[-a_{1,i+1}b_{1,i-1}]>2e+S_{i-1}-R_{i+1}$.
\elm 
\pf First note that $A_{i-1}+A_i>2e+R_i-S_i$ iff
$A'_{i-1}+A'_i>2e+R_i-S_i$. Indeed the ``only if'' implication follows
from $A'_{i-1}\geq A_{i-1}$ and $A'_i\geq A_i$. Conversely,
$A'_{i-1}+A'_i>2e+R_i-S_i$ implies $A_{i-1}+A_i>2e+R_i-S_i$ if
$A'_{i-1}=A_{i-1}$ and $A'_i=A_i$. Otherwise either of
$A'_{i-1}\neq A_{i-1}$ and $A'_i\neq A_i$ implies
$A_{i-1}+A_i>2e+R_i-S_i$ by Lemma 2.14.

So we have to prove that $A'_{i-1}+A'_i>2e+R_i-S_i$ is equivalent
to
$d[-a_{1,i}b_{1,i-2}]+d[-a_{1,i+1}b_{1,i-1}]>2e+S_{i-1}-R_{i+1}$.

Since $R_{i+1}>S_{i-1}$, by Lemma 2.7(i) and (ii), we have
$A'_{i-1}=\min\{ R_i-S_{i-1}+d[-a_{1,i}b_{1,i-2}],
R_i+R_{i+1}-S_{i-2}-S_{i-1}+d[-a_{1,i+1}b_{1,i-1}]\}$ and
$A'_i=\min\{ R_{i+1}-S_i+d[-a_{1,i+1}b_{1,i-1}],
R_{i+1}+R_{i+2}-S_{i-1}-S_i+d[-a_{1,i}b_{1,i-2}]\}$. It follows that
$R_i-S_{i-1}+d[-a_{1,i}b_{1,i-2}]\geq A'_{i-1}$, with equality iff
$R_{i+1}-S_{i-2}+d[-a_{1,i+1}b_{1,i-1}]\geq d[-a_{1,i}b_{1,i-2}]$, and
$R_{i+1}-S_i+d[-a_{1,i+1}b_{1,i-1}]\geq A'_i$, with equality iff
$R_{i+2}-S_{i-1}+d[-a_{1,i}b_{1,i-2}]\geq d[-a_{1,i+1}b_{1,i-1}]$.

If $A'_{i-1}+A'_i>2e+R_i-S_i$, then
$(R_i-S_{i-1}+d[-a_{1,i}b_{1,i-2}])+(R_{i+1}-S_i+d[-a_{1,i+1}b_{1,i-1}])\geq
A'_{i-1}+A'_i>2e+R_i-S_i$, so
$d[-a_{1,i}b_{1,i-2}]+d[-a_{1,i+1}b_{1,i-1}]>2e+S_{i-1}-R_{i+1}$.

Conversely, assume that
$d[-a_{1,i}b_{1,i-2}]+d[-a_{1,i+1}b_{1,i-1}]>2e+S_{i-1}-R_{i+1}$. We
have $d[-a_{1,i}b_{1,i-2}]\leq\beta_{i-2}\leq (S_{i-1}-S_{i-2})/2+e$
and
$d[-a_{1,i+1}b_{1,i-1}]\leq\alpha_{i+1}\leq (R_{i+2}-R_{i+1})/2+e$. 

Then $d[-a_{1,i}b_{1,i-2}]+d[-a_{1,i+1}b_{1,i-1}]>2e+S_{i-1}-R_{i+1}
=S_{i-2}-R_{i+1}+2((S_{i-1}-S_{i-2})/2+e)\geq
S_{i-2}-R_{i+1}+2d[-a_{1,i}b_{1,i-2}]$, so
$R_{i+1}-S_{i-2}+d[-a_{1,i+1}b_{1,i-1}]>d[-a_{1,i}b_{1,i-2}]$, which
implies that $A'_{i-1}=R_i-S_{i-1}+d[-a_{1,i}b_{1,i-2}]$.

Similarly,
$d[-a_{1,i}b_{1,i-2}]+d[-a_{1,i+1}b_{1,i-1}]>2e+S_{i-1}-R_{i+1}
=S_{i-1}-R_{i+2}+2((R_{i+2}-R_{i+1})/2+e)\geq
S_{i-1}-R_{i+2}+2d[-a_{1,i+1}b_{1,i-1}]$, so
$R_{i+2}-S_{i-1}+d[-a_{1,i}b_{1,i-2}]>d[-a_{1,i+1}b_{1,i-1}]$, which
implies that $A'_i=R_{i+1}-S_i+d[-a_{1,i+1}b_{1,i-1}]$.

It follows that $A'_{i-1}+A'_i
=(R_i-S_{i-1}+d[-a_{1,i}b_{1,i-2}])+(R_{i+1}-S_i+d[-a_{1,i+1}b_{1,i-1}])
>2e+R_i-S_i$. \qed

In practice, we will use the following weaker version:

\bco Suppose that $M,N$ satisfy 2.1(i) and (ii) and
$R_{i+1}>S_{i-1}$. If $d[-a_{1,i+1}b_{1,i-1}]$ or
$d[-a_{1,i}b_{1,i-2}]$ is $>2e+S_{i-1}-R_{i+1}$ then
$A_{i-1}+A_i>2e+R_i-S_i$.
\eco

\blm Suppose that $M,N$ satisfy 2.1(i) and (ii), $R_{i+1}>S_{i-1}$ and
$A_{i-1}+A_i>2e+R_i-S_i$. We have:

(i) $\a_i+A_{i-1}>2e$ or $\a_i+d[-a_{1,i+1}b_{1,i-1}]>2e$.

(ii) $\b_{i-1}+A_i>2e$ or $\b_{i-1}+d[-a_{1,i}b_{1,i-2}]>2e$.
\elm
\pf (i) Note that $R_{i+1}-S_i+d[-a_{1,i+1}b_{1,i-1}]\geq A_i$ so
$d[-a_{1,i+1}b_{1,i-1}]\geq S_i-R_{i+1}+A_i$. We have $\a_i =\min\{
(R_{i+1}-R_i)/2+e,R_{i+1}-R_i+d[-a_{i,i+1}]\}$. If
$\a_i=(R_{i+1}-R_i)/2+e$ then
$(\a_i+A_{i-1})+(\a_i+d[-a_{1,i+1}b_{1,i-1}])\geq
2\a_i+A_{i-1}+S_i-R_{i+1}+A_i=S_i-R_i+2e+A_{i-1}+A_i>4e$ which implies
that $\a_i+A_{i-1}$ or $\a_i+d[-a_{1,i+1}b_{1,i-1}]$ is $>2e$ so we
are done.

If $\a_i=R_{i+1}-R_i+d[-a_{i,i+1}]$ then we have either $\a_i\geq
R_{i+1}-R_i+d[-a_{1,i+1}b_{1,i-1}]\geq S_i-R_i+A_i$ or $\a_i\geq
R_{i+1}-R_i+d[a_{1,i-1}b_{1,i-1}]\geq R_{i+1}-R_i+A_{i-1}$. In the
first case we have $\a_i+A_{i-1}\geq S_i-R_i+A_i+A_{i-1}>2e$ so we are
done. In the second case $\a_i+d[-a_{1,i+1}b_{1,i-1}]\geq
R_{i+1}-R_i+A_{i-1}+S_i-R_{i+1}+A_i>2e$ and again we are done.

(ii) is similar. If $m=n$ it follows from (i) by duality at index
$n+1-i$. \qed 

\subsection*{A generalization of 2.1(iv)}

\blm If $M,N$ satisfy 2.1(i)-(iv) then $[b_1,\ldots,b_j]\rep
[a_1,\ldots,a_{l-1}]$ whenever $R_l-S_j>2e$.
\elm
\pf Note that if $l=i+2$, $j=i-1$ we obtain a stronger statement than
in (iv). In fact, except for the case from (iv), our lemma will
follow from (i)-(iii).

First note that $R_{l+1}+2e\geq R_l>S_j+2e\geq S_{j-1}$ so
$R_l+R_{l+1}>S_{j-1}+S_j\geq R_{j-1}+R_j$. It implies that $j\leq
l$. (This obviously happens also if $l=m$ or $j=1$, when $R_{l+1}$
resp. $S_{j-1}$ is not defined.) We cannot have $j=l$ since it
would imply $R_j>S_j+2e$ so $R_l+R_{l+1}=R_j+R_{j+1}\leq
S_{j-1}+S_j$. Hence $j<l$. We have the following cases:

$l-j=1$. Assume $R_{i+1}-S_i>2e$. We want to prove that
$[b_1,\ldots,b_i]\ap [a_1,\ldots,a_i]$. From Corollary 2.10 we have
$a_{1,i}b_{1,i}\in\fs$. Since $d[-a_{1,i+1}b_{1,i-1}]\geq
0>e-(R_{i+1}-S_i)/2$ we have $A_i\neq A'_i$ by 2.15. By Lemma 2.14 we
get $R_{i+1}>S_{i-1}$ and $A_{i-1}+A_i>2e+R_i-S_i$. By 2.1(iii) we
have $[b_1,\ldots,b_{i-1}]\rep [a_1,\ldots,a_i]$. Together with
$a_{1,i}b_{1,i}\in\fs$, this implies $[b_1,\ldots,b_i]\ap
[a_1,\ldots,a_i]$. 

$l-j=2$. Assume $R_{i+1}-S_{i-1}>2e$. We want to prove that
$[b_1,\ldots,b_{i-1}]\rep [a_1,\ldots,a_i]$. We have $R_{i+1}>S_{i-1}$
and $d[-a_{1,i+1}b_{1,i-1}]\geq 0>2e+S_{i-1}-R_{i+1}$. So
$A_{i-1}+A_i>2e+R_i-S_i$ by Corollary 2.17 and we can use 2.1(iii).

$l-j=3$. Assume $R_{i+2}-S_{i-1}>2e$. We want to prove that
$[b_1,\ldots,b_{i-1}]\rep [a_1,\ldots,a_{i+1}]$. We may assume that
$R_{i+1}>S_{i-1}$ or $R_{i+2}>S_i$ since otherwise we can use
condition (iv) of the main theorem. We also may assume
$a_{1,i+1}b_{1,i-1}\in -\fs$ since otherwise the statement is
trivial. Hence $d(-a_{1,i+1}b_{1,i-1})=\j$ which implies
$d[-a_{1,i+1}b_{1,i-1}]=\min\{ \a_{i+1},\b_{i-1}\}$. But
$\a_{i+1}\geq R_{i+2}-R_{i+1}$ if $R_{i+2}-R_{i+1}\leq 2e$ and
$\a_{i+1}>2e$ otherwise. Similarly for $\b_{i-1}$. Therefore
$d[-a_{1,i+1}b_{1,i-1}]\geq\min\{
R_{i+2}-R_{i+1},S_i-S_{i-1},2e\}$. It is enough to prove that
$[b_1,\ldots,b_{i-1}]\rep [a_1,\ldots,a_i]$ or $[b_1,\ldots,b_i]\rep
[a_1,\ldots,a_{i+1}]$. To do this we need to prove that
$R_{j+1}>S_{j-1}$ and $A_{j-1}+A_j>2e+R_j-S_j$ hold for $j=i$ or
$i+1$. We have two cases:

1. $R_{i+2}-R_{i+1}\leq S_i-S_{i-1}$. If $R_{i+1}\leq S_{i-1}$ we get
by adding $R_{i+2}\leq S_i$ which contradicts our assumption. Thus we
may assume that $R_{i+1}>S_{i-1}$. We want to prove that
$A_{i-1}+A_i>2e+R_i-S_i$. To do this we use Corollary 2.17. Since
$d[-a_{1,i+1}b_{1,i-1}]\geq\min\{ R_{i+2}-R_{i+1},S_i-S_{i-1},2e\}
=\min\{ R_{i+2}-R_{i+1},2e\}$ it is enough to prove that $\min\{
R_{i+2}-R_{i+1},2e\} >2e+S_{i-1}-R_{i+1}$. Since $R_{i+1}>S_{i-1}$ we
have $2e>2e+S_{i-1}-R_{i+1}$ and since $R_{i+2}-S_{i-1}>2e$ we have
$R_{i+2}-R_{i+1}>2e+S_{i-1}-R_{i+1}$ so we are done.

2. $R_{i+2}-R_{i+1}\geq S_i-S_{i-1}$. If $R_{i+2}\leq S_i$ we get by
subtracting $R_{i+1}\leq S_{i-1}$ which contradicts our
assumption. Thus we may assume that $R_{i+2}>S_i$. We want to prove
that $A_i+A_{i+1}>2e+R_{i+1}-S_i$. Again we use Corollary 2.17. Since
$d[-a_{1,i+1}b_{1,i-1}]\geq\min\{ R_{i+2}-R_{i+1},S_i-S_{i-1},2e\}
=\min\{ S_i-S_{i-1},2e\}$ it is enough to prove that $\min\{
S_i-S_{i-1},2e\} >2e+S_i-R_{i+2}$. Since $R_{i+2}>S_i$ we have
$2e>2e+S_i-R_{i+2}$ and since $R_{i+2}-S_{i-1}>2e$ we have
$S_i-S_{i-1}>2e+S_i-R_{i+2}$ so we are done. \qed

{\bf Remark} One can also prove that if $R_l-S_j>2e$ then also $\[
a_1,\ldots,a_{l-1}\]$ and $\[ b_1,\ldots,b_j\]$ satisfy the conditions
(i)-(iv) of Theorem 2.1 so if Theorem 2.1 is true then $\[
b_1,\ldots,b_j\]\rep\[ a_1,\ldots,a_{l-1}\]$.

\subsection*{Reduction to the case $m=n$}

\bdf We say that $N\leq M$ if $FN\rep FM$ and $M,N$ satisfy the
conditions (i)-(iv) of the Theorem 2.1.
\edf

At this point the condition $N\leq M$ depends on the choice
of BONGs so we will mean $N\leq M$ relative to some BONGs. The
independence of BONGs will be proved in $\S 3$.

Of course the aim of this paper is to prove that $N\rep M$ iff $N\leq
M$.

\blm Suppose that $m>n$ and $K$ is a third lattice and $s\in\ZZ$. If
$s\gg 0$ and $FM\ap FN\pp FK$ then $N\leq M$ iff $N\pp\p^sK\leq M$.
\elm
\pf Since $s\gg 0$ we have $R_i(\p^s K)=R_i(K)+2s\gg 0$. Let $\p^s
K\ap\[ b_{n+1},\ldots,b_m\]$ and $S_i=\od b_i$ for $n+1\leq i\leq m$. We
have $S_i\gg 0$ for $n+1\leq i\leq m$. In particular $S_n\leq
S_{n+1}$, $S_{n-1}\leq S_{n+1}$ and, if $m\geq n+2$, $S_n\leq
S_{n+2}$. This implies, by [B1, Corollary 4.4(v)], that $\(
N:=N\pp\p^sK\ap\[ b_1,\ldots,b_m\]$ relative to a good BONG. So we have
to prove that for $s$ large we have $N\leq M$ iff $\( N\leq M$.

Condition 2.1(i) for $i\leq n$ is the same in both cases. If $i>n$
then $S_i\gg 0$ so $R_i\leq S_i$.

If $1\leq i\leq n-1$ then, by the definition of $\b_i=\a_i(N)$ and
$\a_i(\( N)$, we have $\a_i(\( N)=\min (\{\b_i\}\cup\{
S_{j+1}-S_i+d(-b_{j,j+1})\mid n\leq j<m\} )$. But $S_{j+1}$ is
large  for $n\leq j$ so $S_{j+1}-S_i+d(-b_{j,j+1})\geq
S_{j+1}-S_i\geq\b_i$. Thus $\a_i(\( N)=\b_i$ for $1\leq i\leq n-1$. We
keep the notation $\b_i$ for $\a_i(\( N)$ also when $n\leq i\leq
m-1$. In particular, $S_n\ll S_{n+1}$ so $\b_n=(S_{n+1}-S_n)/2+e\gg
0$. As a consequence the value of $d[\e a_{1,i}b_{1,j}]$ coresponding
to $M,\( N$ is the same as the one for $M,N$ when $j\leq n-1$. When
$j=n$ we have  $d[\e a_{1,i}b_{1,n}](M,\( N)=\min\{ d[\e
a_{1,i}b_{1,n}],\b_n\}$. Since $\b_n\gg 0$ we usually have $d[\e
a_{1,i}b_{1,n}](M,\( N)=d[\e a_{1,i}b_{1,n}]$. The only exception is
when $i\in\{ 0,m\}$ and $\e a_{1,i}b_{1,n}\in\fs$. In this case $d[\e
a_{1,i}b_{1,n}]=d(\e a_{1,i}b_{1,n})=\j$ while $d(\e
a_{1,i}b_{1,n})(M,\( N)=\b_n=(S_{n+1}-S_n)/2+e\gg 0$.

As a consequence, $A_i(M,\( N)=A_i(M,N)$ for $1\leq i\leq n$. So
condition (ii) is the same for both pairs of lattices. Same
for the condition (iii) when $i\leq n$. When $i=n+1$ we usually have
$(S_{n+1}+A_{n+1})(M,\( N)=(S_{n+1}+A_{n+1})(M,N)$ so the condition
2.1(iii) is the same. (See Definition 4 and the following remark). The
only exception is when $m=n+2$ and $-a_{1,m}b_{1,n}$ is a square. In
this case we have $(S_{n+1}+A_{n+1})(M,N)
=R_{n+2}+d[-a_{1,n+2}b_{1,n}]=R_m+d[-a_{1,m}b_{1,n}]=\j$, while
$A_{n+1}(M,\( N)=\min\{ (R_{n+2}-S_{n+1})/2+e,
R_{n+2}-S_{n+1}+(S_{n+1}-S_n)/2+e\}$ so $(S_{n+1}+A_{n+1})(M,\(
N)=\min\{ (R_{n+2}+S_{n+1})/2+e, R_{n+1}+(S_{n+1}-S_{n-1})/2+e\}\gg
0$. Thus condition $A_n+(S_{n+1}+A_{n+1})>2e+R_{n+1}$ is satisfied in
both cases. Condition 2.1(ii) for $i>n$ is true because $S_i\gg 0$ so
$A_i\leq (R_{i+1}-S_i)/2+e<0\leq d[a_{1,i}b_{1,i}]$. Condition
2.1(iii) for $i>n+1$ is vacuous since $S_{i-1}\gg 0$ so
$R_{i+1}<S_{i-1}$. 

Finally, 2.1(iv) is the same for both $M,N$ and $M,\( N$ when $i\leq
n+1$. (At $i=n+1$ note that the condition $S_i\geq R_{i+2}$ is ignored
for the lattices $M,N$ while for $M,\( N$ it is always satisfied  since
$S_{n+1}\gg 0$.) When $i>n+2$ the condition 2.1(iv) for $M,\( N$ is
vacuous since we have $S_{i-1}\gg 0$ and so $R_{i+1}<S_{i-1}+2e$. \qed

\blm The main theorem is true if it is true when $m=n$.
\elm
\pf Suppose that $m>n$ and $FM=FN\pp V$ where $V$ is a quadratic space
with $\dim V=m-n$. Let $K$ be a lattice on $V$. If $N\sbq M$ then
$N\pp (M\cap V)\sbq M$. Since both $M\cap V$ and $K$ are lattices on
$V$ we have $\p^sK\sbq M\cap V$ for $s$ large enough so $N\pp\p^s
K\sbq M$. Thus $N\rep M$ iff $N\pp\p^sK\rep M$ for $s$ large
enough. On the other hand by Lemma 2.20 we have $N\leq M$ iff
$N\pp\p^sK\leq M$ for $s$ large enough. Therefore the main theorem is
true for $M,N$ iff it is true for $M,N\pp\p^sK$. \qed

From now on we restrict to the case $m=n$. In particular this allows
the use of duality which can shorten the proof.

\subsection*{Idea for the proof of necessity}
The proof of necessity consists of two steps:

(I) Proof of transitivity of $\leq$.

(II) Proof that $N\sbq M$ implies $N\leq M$ in the case when
$[M:N]\spq\p$. 

Indeed, if $N\rep M$ then we may assume that $N\sbq M$. If
$[M:N]\sb\p$ then we can go from $M$ to $N$ through a set of minimal
inclusions $M=M_0\sp M_1\sp\ldots\sp M_k=N$ with $[M_{i-1}:M_i]=\p$. By
(II) we have $M_{i-1}\leq M_i$ for $1\leq i\leq k$ and by transitivity
of $\leq$ we get $N\leq M$. 

In order to prove the necessity in the case $[M:N]=\p$ we need to
use Jordan decompositions so we need to translate the main theorem
in terms of Jordan decompositions. This will be done in the next
section.

\section{Approximations. Main theorem in terms of Jordan
decompositions} 

Let $L\ap [a_1,\ldots,a_n]$ be relative to a good BONG and let
$L=L_1\pp\cdots\pp L_t$ be a Jordan decomposition.

\bdf If $1\leq i\leq n-1$ then an element $X_i\in\ff/\fs$ is called an
approximation for $a_{1,i}$ if $d(a_{1,i}X_i)\geq\a_i$. We denote this
by $X_i\sim a_{1,i}$. For $i=0$ or $n$ we say that $X_0\sim a_{1,0}$
if $X_0=a_{1,0}=1$ and we say that $X_n\sim a_{1,n}$ if
$X_n=a_{1,n}=\det FL$ (in $\ff/\fs$).
\edf

\bdf Let $1\leq i\leq n-1$ and let $V_i$ be a quadratic space of
dimension $i$ s.t. $\det V_i\sim a_{1,i}$. We say that $V_i$ is an
approximation to the left for $[a_1,\ldots,a_i]$ if $\a_{i-1}+\a_i>2e$ or
$i=1$ implies $[a_1,\ldots,a_{i-1}]\rep V_i$; an approximation to the
right if $\a_i+\a_{i+1}>2e$ or $i=n-1$ implies $V_i\rep
[a_1,\ldots,a_{i+1}]$; an approximation if it is an approximation both to
the left and right. We denote this by $V_i\sim_l$(resp. $\sim_r$ or
$\sim$)$[a_1,\ldots,a_i]$. 
\edf

{\bf Remark} If $i=1$ the $[a_1,\ldots,a_{i-1}]$ is $0$ (the quadratic
space of dimension zero) and the condition $[a_1,\ldots,a_{i-1}]\rep V_i$
from the definition of $\sim_l$ means $0\rep V_1$ and it is always
satisfied. If $i=n-1$ then $[a_1,\ldots,a_{i+1}]\ap FL$ so the condition
$V_i\rep [a_1,\ldots,a_{i+1}]$ from the definition of $\sim_r$ means
$V_{n-1}\rep FL$.
\vskip 0.3cm

We want now to find approximates for $a_{1,i}$ and $[a_1,\ldots,a_i]$ in
terms of the Jordan decomposition $L=L_1\pp\cdots\pp L_t$. Like in [B3] we
denote $r_k=\ord\ss L_k$, we take $\AA_k$ a norm generator for
$L^{\ss L_k}$ and define $u_k=\ord\AA_k=\ord\nn L^{\ss
L_k}$. Also denote $L_{(k)}=L_1\pp\cdots\pp L_k$, $n_k=\dim L_{(k)}$ and by
$\FF_k$ the invariants defined in [OM] 93:24 for $1\leq k\leq t-1$. 
\vskip 4mm

{\bf Notation} If $V,W$ are two quadratic spaces and $W\rep V$ then we
denote by $V\top W$ a quadratic space $U$ s.t. $V\ap W\pp U$.

\bff {\bf Remark} Note that if $\a_{i-1}+\a_i\leq 2e$
(resp. $\a_i+\a_{i+1}\leq 2e$ or both) Then $V_i\sim_l$(resp. $\sim_r$ or
$\sim$)$[a_1,\ldots,a_i]$ whenever $\dim V_i=i$ and $\det V_i\sim
a_{1,i}$. (E.g. we can take $V_i=[1,\ldots,1,X_i]$, where $X_i$ is an
arbitrary approximation for $a_{1,i}$.) By [B3, Lemma 3.3] the condition
$\a_{i-1}+\a_i>2e$ implies $R_{i+1}>R_{i-1}$ so $i=n_k$ or $i=n_k+1$
for some $k$. (Otherwise $n_k+1<i<n_{k+1}$ for some $k$ so
$R_{i-1}=R_{i+1}$.) Similarly $\a_i+\a_{i+1}>2e$ implies $R_{i+2}>R_i$
so $i=n_k-1$ or $i=n_k$ for some $k$. Hence in order to have
$\a_{i-1}+\a_i>2e$ or $\a_i+\a_{i+1}>2e$ then $i$ must be of the form
$n_k-1,n_k$ or $n_k+1$.
\eff

\blm Let $n_{k-1}\leq i<n_k$  with $1\leq k\leq t$.

(i) If $i\ev n_{k-1}\m2$ then $X_i=(-1)^{(i-n_{k-1})/2}\det
FL_{(k-1)}$ approximates $a_{1,i}$.

(ii) If $i\ev n_{k-1}+1\m2$ then $X_i=\AA_k\det FL_{(k-1)}$
approximates $a_{1,i}$.
\elm
\pf Suppose that $L\ap [a_1,\ldots,a_n]$ relative to the good BONG
$x_1,\ldots,x_n$. Let $L=L^1\pp\cdots\pp L^m$ be a maximal norm
splitting with all binary components improper s.t. $x_1,\ldots,x_n$ is
obtained by putting together the BONGs of $L^1,\ldots,L^m$. Let
$L=K_1\pp\cdots\pp K_t$ be the Jordan composition obtained by putting
together $L^j$'s with same scale. We have $K_{(k)}\ap
[a_1,\ldots,a_{n_k}]$ relative to $x_1,\ldots,x_{n_k}$. 

Take first $i=n_{k-1}$. We have $a_{1,i}=\det FK_{(k-1)}$ and
$X_i=\det FL_{(k-1)}$. By [OM] 92:28 we have $\det L_{(k-1)}/\det
K_{(k-1)}\ap 1\mod\FF_{k-1}$, which is equivalent to
$d(X_ia_{1,i})=d(\det FL_{(k-1)}\det FK_{(k-1)})\geq\ord\FF_{k-1}$. By
[B3, Lemma 2.16(ii)] we have either $\ord\FF_{k-1}=\a_i$ or both
$\ord\FF_{k-1}$ and $\a_i$ are $>2e$. Thus
$d(X_ia_{1,i})\geq\ord\FF_{k-1}$ implies that $d(X_ia_{1,i})\geq\a_i$
so $X_i$ approximates $a_{1,i}$. The same happens if $k=1$, when
$L_{(k-1)}=L_{(0)}=0$ and so $X_0=\det 0=1$. 

Suppose now that $i=n_{k-1}+1<n_k$. By [B3, Corollary 2.17(i)] we have
$\a_i=\ord\AA_k\1\ww_k$. By [B3, Lemma 2.13(iii)] $\BB_k:=a_i$ is a 
norm generator for $L^{\ss L_k}$. Hence $\AA_k\ap\BB_k\mod\ww_k$,
which is equivalent to $d(\AA_k\BB_k)\geq\ord\AA_k\1\ww_k=\a_i$. Also
$i-1=n_{k-1}$ so by the previous case $\det FL_{(k-1)}$ approximates
$a_{1,i-1}$ so $d(\det FL_{(k-1)}a_{1,i-1})\geq\a_{i-1}\geq
R_i-R_{i+1}+\a_i\geq\a_i$. (We have $R_i=u_k\geq 2r_k-u_k=R_{i+1}$.)
This, together with $d(\AA_ka_i)=d(\AA_k\BB_k)\geq\a_i$, implies
$d(\AA_k\det L_{(k-1)}a_{1,i})\geq\a_i$ so $X_i=\AA_k\det L_{(k-1)}$
approximates $a_{1,i}$. Note that the inequality $d(\det
FL_{(k-1)}a_{1,i-1})\geq\a_i$ also holds if $k=1$, when $i=1$, since
$\det FL_{(0)}a_{1,0}=1$.

For the general case we claim that if $n_{k-1}\leq i$ and $i+2<n_k$
and $X_i$ approximates $a_{1,i}$ then $-X_i$ approximates
$a_{1,i+2}$. Indeed, $X_i\sim a_{1,i}$ means
$d(X_ia_{1,i})\geq\a_i$. Since $n_{k-1}+1\leq i+1,i+3\leq n_k$ we have
$R_{i+1}=R_{i+3}$ ([B3, Lemma 2.13(i)]), which implies
$R_{i+1}+\a_{i+1}=R_{i+2}+\a_{i+2}$ ([B3, Corollary 2.3(i)]). So
$d(X_ia_{1,i})\geq\a_i\geq R_{i+1}-R_{i+3}+\a_{i+2}=\a_{i+2}$ and
$d(-a_{i+1,i+2})\geq R_{i+1}-R_{i+2}+\a_{i+1}=\a_{i+2}$. By domination
principle we have $d(-X_ia_{1,i+2})\geq\a_i=\a_{i+2}$. So $-X_i$
approximates $a_{1,i+2}$. Note that the inequality
$d(X_ia_{1,i})\geq\a_{i+2}$ also holds if $i=0$ since
$d(X_0a_{1,0})=d(1)=\j$.

Since $\det FL_{(k-1)}$ approximates $a_{1,n_{k-1}}$ we can use the
result above and we get by induction that if $n_{k-1}\leq
i=n_{k-1}+2l<n_k$ then $(-1)^l\det FL_{(k-1)}=(-1)^{(i-n_{k-1})/2}\det
FL_{(k-1)}$ approximates $a_{1,i}$. Similarly $\pm\AA_k\det
FL_{(k-1)}$ approximates $a_{1,i}$ if $n_{k-1}\leq
i=n_{k-1}+1+2l<n_k$. But since we have the liberty of choosing $\AA_k$
or $-\AA_k$ as the norm generator for $L^{\ss L_k}$ the $\pm$ sign can be
dropped. \qed

We have $L^\*\ap [a_1^\*,\ldots,a_n^\* ]$ with $a_i^\* =a_{n+1-i}\1$
so $a_{1,i}^\*  =a_{n+1-i,n}\1 =a_{1,n-i}(\det L)\1$. Also $\a_i^\*
=\a_{n-i}$. It follows that if $X_{n-i}$ is an approximation for
$a_{1,n-i}$ then $X_i^\* :=X_{n-i}\det FL$ is an approximation for
$a_{1,i}^\*$. (If $1\leq i\leq n-1$ then $d(a_{1,i}^\* X_{n-i}\det
FL)=d(a_{1,n-i}X_{n-i})\geq\a_{n-i}=\a_i^\*$. If $i=0$ or $i=n$ we just
note that $X_0=1$ and $X_n=\det FL$ so $X_0^\* =1$ and $X_n^\* =\det
FL=\det FL^\*$.)

Also  we have the Jordan decomposition $L^\* =L_t^\*\pp\cdots\pp L_1^\*$
so $L_{(k)}^\* =L_t^\*\pp\cdots\pp L_{t+1-k}^\* =(L_{(t+1-k)}^* )^\*$ and
$n_k^\* :=\dim L_{(k)}^\* =n-n_{t-k}$. Also $\AA_i^\*
=\pi^{-2r_{t+i-1}}\AA_{t+1-i}$.

\bco Let $n_{k-1}<i\leq n_k$ with $1\leq k\leq t$.

(i) If $i\ev n_k\m2$ then $X_i=(-1)^{(n_k-i)/2}\det FL_{(k)}$
approximates $a_{1,i}$.

(ii) If $i\ev n_k+1\m2$ then $X_i=\AA_k\det FL_{(k)}$ approximates
$a_{1,i}$.
\eco
\pf Let $l=t+1-k$ and $j=n-i$. Then $n_{k-1}<i\leq n_k$ implies
$n-n_k\leq n-i<n-n_{k-1}$ i.e. $n_{t-k}^\*\leq n-i<n_{t-k+1}^\*$. If
$i\ev n_k+1\m2$ then $n-i\ev n-n_k+1=n_{t-k}^\* +1\m2$. By Lemma 3.2
we can take $X_{n-i}^\* =\AA_{t+1-k}^\*\det FL_{(t-k)}^\*$. But in
$\ff/\fs$ we have $\AA_{t+1-k}^\* =\AA_k$ and $\det FL_{(t-k)}^\*
=\det FL_{(k+1)}^*=\det FL_{(k)}\det FL$ so $X_{n-i}^\* =\AA_k\det
FL_{(k)}\det FL$. Since $X_{n-i}^\*$ is an approximation for
$a_{1,n-i}^\*$, $X_i=X_{n-i}^\*\det FL=\AA_k\det FL_{(k)}$ is an
approximation for $a_{1,i}$. Similarly for $i\ev n_k\m2$ we take
$X_{n-i}^\* =(-1)^{((n-i)-(n-n_k))/2}\det FL_{(t-k)}^\*
=(-1)^{(n_k-i)/2}\det FL_{(k)}\det FL$ so $X_i=X_{n-i}^\*\det
FL=(-1)^{(n_k-i)/2}\det FL_{(k)}$. \qed

\blm Let $1\leq i\leq n-1$. We have:

(i) If $i=n_k$ and $V_i=FL_{(k)}$ then $V_i\sim [a_1,\ldots,a_i]$.

(ii) If $i=n_k+1<n_{k+1}$ and $V_i=FL_{(k)}\pp [\AA_{k+1}]$ then
$V_i\sim_l[a_1,\ldots,a_i]$.

(iii) If $i=n_k-1>n_{k-1}$ and $V_i=FL_{(k)}\top [\AA_k]$ then
$V_i\sim_r[a_1,\ldots,a_i]$.
\elm
\pf We use the same splitting $L=K_1\pp\cdots\pp K_t$ from the proof of
Lemma 3.2. If $i=n_t$ then $[a_1,\ldots,a_i]\ap FK_{(k)}$. Also by
[B3, Lemma 2.13(iii)] we can choose $\BB_k=\pm a_{n_{k-1}+1}$ or
$\pm\pi^{2u_k-2r_k}a_{n_k}$ as a norm generator for $L^{\ss L_k}$. Thus if
$i=n_k+1$ then $[a_1,\ldots,a_i]=[a_1,\ldots,a_{i-1}]\pp [a_i]\ap
FK_{(k)}\pp [\BB_{k+1}]$, where we took $\BB_{k+1}=a_{n_k+1}=a_i$ and
if $i=n_k-1$ then $[a_1,\ldots,a_i]=[a_1,\ldots,a_{i+1}]\top [a_{i+1}]\ap
FK_{(k)}\top [\BB_k]$, where we took
$\BB_k=\pi^{2u_k-2r_k}a_{n_k}=a_{i+1}$ (in $\ff/\fs$).

(i) We have $V_i=FL_{(k)}$ so $\dim V_i=n_k=i$ and $\det V_i=\det
FL_{(k)}$, which is an approximation for $a_{1,i}$ by Lemma 3.2(i).

Suppose that $\a_{i-1}+\a_i>2e$. Suppose first that $\dim
L_{k-1}>1$ i.e. that $i=n_k>n_{k-1}+1$. By [B3, Lemma 3.7(i)]
$\a_{i-1}+\a_i>2e$ implies $\FF_k\sb 4\AA_k\ww_k\1$ so by O'Meara's
93:28(iv) $FK_{(k)}\rep FL_{(k)}\pp [\BB_k]$.  But this implies
$FK_{(k)}\top [\BB_k]\rep
FL_{(k)}$, i.e. $[a_1,\ldots,a_{i-1}]\rep FL_{(k)}=V_i$. (We have
$i-1=n_k-1$.) If $\dim L_{k-1}=1$ i.e. if $i=n_k=n_{k-1}+1$ then by
[B3, Lemma 3.7(iii)] $\a_{i-1}+\a_i>2e$ implies that
$\FF_k\sb 4\AA_k\ww_k\1$ or $\FF_{k-1}\sb 4\AA_k\ww_k\1$. In the first
case the proof goes as before. In the second case we have by O'Meara's
93:28(iii) $FK_{(k-1)}\rep FL_{(k-1)}\pp [\BB_k]$. But we can choose
$\BB_k$ s.t. $L_k\ap\la\BB_k\ra$ so we have $FL_{(k-1)}\pp [\BB_k]\ap
FL_{(k)}=V_i$. Thus $[a_1,\ldots,a_{i-1}]\ap FK_{(k-1)}\rep V_i$. (We
have $i-1=n_{k-1}$.) Note that this also holds when $i=1$, i.e. when
$k=1$ and $L_1$ is unary, (see the remark following Definition 10).

Suppose now that $\a_i+\a_{i+1}>2e$. Suppose first that $\dim
L_{k+1}>1$ i.e. that $i+1=n_k+1<n_{k+1}$. By [B3, Lemma 3.7(ii)]
$\a_i+\a_{i+1}>2e$ implies $\FF_k\sb 4\AA_{k+1}\ww_{k+1}\1$ so by
O'Meara's 93:28(iii) we get $V_i=FL_{(k)}\rep FK_{(k)}\pp
[\BB_{k+1}]\ap [a_1,\ldots,a_{i+1}]$. (We have $i+1=n_k+1$.) If $\dim
L_{k+1}=1$, i.e. if $i+1=n_k+1=n_{k+1}$ then, by [B3, Lemma 3.7(iii)],
$\a_i+\a_{i+1}>2e$ implies that $\FF_k\sb 4\AA_{k+1}\ww_{k+1}\1$
or $\FF_{k+1}\sb 4\AA_{k+1}\ww_{k+1}\1$. The first case was treated
above. In the second case we use O'Meara's 93:28(iv) and we get
$FL_{(k+1)}\rep FK_{(k+1)}\pp [\BB_{k+1}]$ or $FL_{(k+1)}\top
[\BB_{k+1}]\rep FK_{(k+1)}\ap [a_1,\ldots,a_{i+1}]$. (We have
$i+1=n_{k+1}$.) But  we can choose $\BB_{k+1}$
s.t. $L_{k+1}\ap\la\BB_{k+1}\ra$ so $FL_{(k+1)}\top [\BB_{k+1}]\ap
FL_{(k+1)}\top FL_{k+1}=L_{(k)}=V_i$. Thus $V_i\rep
[a_1,\ldots,a_{i+1}]$. This also happens when $i=n-1$, i.e. when
$k=t-1$ and $K_t$ is unary because $V_{n-1}\ap FK_{(k)}\rep FL\ap
[a_1,\ldots,a_n]$. 

(ii) We have $V_i=FL_{(k)}\pp [\AA_{k+1}]$ so $\det V_i=\AA_{k+1}\det
FL_{(k)}$, which is an approximation for $a_{1,i}$ by Lemma
3.2(ii). If $\a_{i-1}+\a_i>2e$ then  $\FF_k\sb 4\AA_{k+1}\ww_{k+1}\1$
by [B3, Lemma 3.7(ii)] so $FK_{(k)}\rep FL_{(k)}\pp [\AA_{k+1}]=V_i$
By O'Meara's 93:28(iii). But $FK_{(k)}\ap [a_1,\ldots,a_{i-1}]$ (we
have $i-1=n_k$) so $[a_1,\ldots,a_{i-1}]\rep V_i$. This also holds if
$i=1$ (see the remark following Definition 10).

(iii) We have $V_i=FL_{(k)}\top [\AA_k]$ so $\det V_i=\AA_k\det
FL_{(k)}$, which is an approximation for $a_{1,i}$ by Corollary
3.3(ii). If $\a_i+\a_{i+1}>2e$ then $\FF_k\sb 4\AA_k\ww_k\1$
by [B3, Lemma 3.7(i)] (we have $i+1=n_k>n_{k-1}+1$). By O'Meara's
93:28(iv) we have $FL_{(k)}\rep FK_{(k)}\pp [\AA_k]$
i.e. $V_i=FL_{(k)}\top [\AA_k]\rep FK_{(k)}\ap [a_1,\ldots,a_{i+1}]$. (We
have $i+1=n_k$.) The relation $V_i\rep [a_1,\ldots,a_{i+1}]$ also holds
when $i=n-1$, i.e. when $k=t$. In this case we have $V_{n-1}\ap
FL_{(t)}\top [\AA_t]\rep FL_{(t)}=FL\ap [a_1,\ldots,a_n]$. \qed

\blm (i) If $r_k\leq r\leq r_{k+1}$ then $\od\nn L^{\p^r}=\min\{
u_k+2(r-r_k),u_{k+1}\}$. If $r\leq r_1$ then $\od\nn L^{\p^r}=u_1$. If $r\geq
r_t$ then $\od\nn L^{\p^r}=u_t+2(r-r_t)$.

(ii) If $u_k<u_{k+1}$ and $u_k-2r_k<u_{k-1}-2r_{k-1}$ for some $1\leq k\leq t$
then $\nn L_k=\nn L^{\ss_k}$ regardless of the splitting. (If $k=1$ or $t$
then we ignore the condition $u_k-2r_k<u_{k-1}-2r_{k-1}$ resp. $u_k<u_{k+1}$.)
\elm
\pf (i) If $r_k\leq r\leq r_{k+1}$ then $L\cap\p^r L^\*\spq L\cap\p^{r_{k+1}}
L^\*$ and $L\cap\p^r L^\*\spq p^{r-r_k}L\cap\p^r L^\*
=\p^{r-r_k}(L\cap\p^{r_k} L^\* )$ so $L^{p^r}\spq L^{\p^{r_{k+1}}}$
and $L^{\p^r}\spq\p^{r-r_k}L^{\p^{r_k}}$. Thus \\ $\od\nn
L^{\p^r}\leq\min\{ \od\nn\p^{r-r_k}L^{\p^{r_k}},\od\nn
L^{\p^{r_{k+1}}}\}=\min\{ u_k+2(r-r_k),u_{k+1}\}$. On the other hand
$L^{\p^r}=\p^rL_{(k)}^\*\pp
L_{(k+1)}^*$ so $\od\nn L^{\p^r}=\min\{ 2r+\od\nn L_{(k)}^\*,\od\nn
L_{(k+1)}^*\}$. But $L_{(k+1)}^*\sbq L^{\p^{r_{k+1}}}$ so $\od\nn
L_{(k+1)}^*\geq\od\nn L^{\p^{r_{k+1}}}=u_{k+1}$ and
$\p^{r_k}L_{(k)}^\*\sbq L^{\p^{r_k}}$ so $\od\nn L_{(k)}^\*\geq
-2r_k+\od\nn L^{\p^{r_k}}=u_k-2r_k$. Hence $\od\nn L^{\p^r}\geq\min\{
u_k+2(r-r_k),u_{k+1}\}$.

The other two statements of (i) follow by a similar argument.

(ii) We have $L^{\ss_k}=\ss_kL_{(k-1)}^\*\pp L_k\pp L_{(k+1)}^*$. But
$\od\nn  L_{(k+1)}^*\geq u_{k+1}>u_k=\od\nn L^{\ss_k}$ and $\od\nn
L_{(k-1)}^\*\geq u_{k-1}-2r_{k-1}>u_k-2r_k$ (see the proof of (i)
above) so $\od\ss_k\nn L_{(k-1)}^\* >u_k=\od\nn L^{\ss_k}$. Therefore
$\nn L^{\ss_k}=\nn L_k$ as claimed. \qed

\blm If $i=n_{k-1}+1=n_k-1$, $R_i<R_{i+2}$ or $i=n-1$ and
$R_{i-1}<R_{i+1}$ or $i=1$ then $L_k$ is binary and $\nn L_k=\nn
L^{\ss_k}$ regardless of the Jordan splitting.
\elm
\pf We have $n_k-n_{k-1}=2$ so $L_k$ is binary. We have $i-1=n_{k-1}$
and $i+1=n_k$ so $2r_{k-1}-u_{k-1}=R_{i-1}<R_{i+1}=2r_k-u_k$,
i.e. $u_k-2r_k<u_{k-1}-2r_{k-1}$. Also $i=n_{k-1}+1$ and $i+2=n_k+1$
so $u_{k-1}=R_i<R_{i+2}=u_k$. By Lemma 3.5(ii) we have $\nn L_k=\nn
L^{\ss_k}$. (If $i=1$ then $k=1$ so the condition
$u_k-2r_k<u_{k-1}-2r_{k-1}$ from Lemma 3.5(ii) is ignored. If $i=n-1$
then $k=t$ so the condition $u_k<u_{k+1}$ is ignored.) \qed

\blm Let $1\leq i\leq n-1$. In all of the following situations $V_i$
 is an approximation for $[a_1,\ldots,a_i]$:

(i) $i=n_k$ and $V_i=FL_{(k)}$

(ii) $i=n_k+1<n_{k+1}$, $R_i=R_{i+2}$ and $V_i=FL_{(k)}\pp
[\AA_{k+1}]$.

(iii) $i=n_k-1>n_{k-1}$, $R_{i-1}=R_{i+1}$ and $V_i=FL_{(k)}\top
[\AA_k]$.

(iv) $i=n_{k-1}+1=n_k-1$, $R_i<R_{i+2}$ or $i=n-1$, $R_{i-1}<R_{i+1}$
or $i=1$, $\AA_k$ is a norm generator for $L^{\ss L_k}$
s.t. $\AA_k\in Q(L_k)$ and $V_i=L_{(k-1)}\pp [\AA_k]$ or
$V_i=L_{(k)}\top [\AA_k]$.

In all other cases $R_{i-1}=R_{i+1}$ and $R_i=R_{i+1}$ so
$\a_{i-1}+\a_i\leq 2e$ and $\a_i+\a_{i+1}\leq 2e$. Hence $V_i\sim
[a_1,\ldots,a_i]$ iff $\dim V_I=i$ and $\det V_i\sim a_{1,i}$.
\elm
\pf If we don't have both $R_{i-1}=R_{i+1}$ and $R_i=R_{i+2}$ then, by
Remark 3.1, $i$ is of the form $n_k-1,n_k$ or $n_k+1$. If $i$ is of
the form $n_k$ we are in case (i). Otherwise $i=n_k+1<n_{k+1}$ or
$i=n_k-1>n_{k-1}$ for some $k$. If $i=n_k+1<n_{k+1}$ then either
$i=n_k+1\leq n_{k+1}-2$, which implies $R_i=R_{i+2}=u_k$ so we are in
the case (ii), or $i=n_k+1=m_{k+1}-1$. If $i=n_k-1>n_{k-1}$ then
either $i=n_k-1\geq n_{k-1}+2$, which implies
$R_{i-1}=R_{i+1}=2r_k-u_k$ so we are in the case (iii), or
$i=n_k-1=n_{k-1}+1$. So we are left to the case when
$i=n_{k-1}+1=n_k-1$ for some $k$. If $R_i=R_{i+2}$ or
$R_{i-1}=R_{i+1}$ then we are in the case (ii) resp. (iii). Otherwise
we are in the case (iv). 

(i) is just Lemma 3.4(i).

(ii) We have $V_i\sim_l [a_1,\ldots,a_i]$ by Lemma 3.4(ii). Since
$\det V_i\sim a_{1,i}$ and $R_i=R_{i+2}$, which implies
$\a_i+\a_{i+1}\leq 2e$, we also have $V_i\sim_r [a_1,\ldots,a_i]$.

(iii) We have $V_i\sim_r [a_1,\ldots,a_i]$ by Lemma 3.4(iii). Since
$\det 
V_i\sim a_{1,i}$ and $R_{i-1}=R_{i+1}$, which implies
$\a_{i-1}+\a_i\leq 2e$, we also have $V_i\sim_l [a_1,\ldots,a_i]$.

(iv) We claim that there is $\AA'_k$, another norm generator for
$L^{\ss L_k}$, s.t. $FL_k\ap [\AA_k,\AA'_k]$. By Lemma 3.6 $L_k$ is
binary and $\nn L_k=\nn L^{\ss_k}$ so any norm generator of $L_k$ is
also a norm generator of $L^{\ss L_k}$. So if $L_k\ap\[\BB_k,\BB\]$
then $\BB_k$ is a norm generator for $L^{\ss L_k}$. We have
$L_k^\*\ap\[\BB\1,\BB_k\1\]$ so
$L_k=\p^{r_k}L_k^\*\ap\[\pi^{2r_k}\BB\1,\pi^{2r_k}\BB_k\1\]$. Now both
$\BB_k$ and $\BB'_k:=\pi^{2r_k}\BB\1$ are norm generators for $L_k$
and thus for $L^{\ss L_k}$ and we also have
$FL_k\ap [\BB_k,\BB ]\ap [\BB_k,\BB'_k]$. Let $\AA_k=\e\BB_k$ and let
$\AA'_k:=\e\BB'_k$. Since $\AA_k,\BB_k$ are norm generators for
$L^{\ss L_k}$ we have $\AA_k\ap\BB_k\mod\ww_k$ so $d(\e
)=d(\AA_k\BB_k)\geq\ord\AA_k\1\ww_k$. Thus $d(\AA'_k\BB'_k)=d(\e
)\geq\ord\AA_k\1\ww_k={\AA'_k}\1\ww_k$, which implies
$\AA'_k\ap\BB'_k\mod\ww_k$. Since $\BB'_k$ is a norm generator for
$L^{\ss L_k}$ so is $\AA'_k$. On the other hand $\e\BB_k=\AA_k\rep
FL_k\ap [\BB_k,\BB'_k]$ so $FL_k\ap [\e\BB_k,\e\BB'_k]\ap
[\AA_k,\AA'_k]$. 

If $V_i=FL_{(k-1)}\pp [\AA_k]$ then $V_i\sim_l[a_1,\ldots,a_i]$ by
Lemma 3.4(ii). Also since $FL_k\ap [\AA_k,\AA'_k]$ we have $V_i\ap
FL_{(k-1)}\pp FL_k\top [\AA'_k]=FL_{(k)}\top [\AA'_k]$. Therefore
$V_i\sim_r[a_1,\ldots,a_i]$ by Lemma 3.4(iii).

Similarly if $V_i=FL_{(k)}\top [\AA_k]$ then $V_i=FL_{(k-1)}\pp
FL_k\top [\AA_k]\ap FL_{(k-1)}\pp [\AA'_k]$. Thus $V_i$ approximates
$[a_1,\ldots,a_i]$ both to the right and left by Lemma 3.4(iii) and
(ii). \qed 

{\bf Remark} If $i=2=n_k-1>n_{k-1}$ and $R_1=R_3$, in order to use
Lemma 3.7(iii), we have to ask the condition $[\AA_k]\rep L_{(k)}$ so
that $V_2=FL_{(k)}\top [\AA_k]$ makes sense. Now $n_k=3$ so $L_{(k)}$ is
ternary. Thus we have to avoid the case when $FL_{(k)}$ is anisotropic
and $\AA_k=-\det FL_{(k)}$ in $\ff/\fs$. This can be easily acheieved
by simply noting that if $\AA_k$ is a norm generator for $L^{\ss L_k}$
then so is, say, $\D\AA_k$. (We have $4\AA_k\oo\sb 2\AA_k\oo\sbq\ww_k$
so $\D\AA_k\in\AA_k\oo^2+4\AA_k\oo\sbq\GG_k$.) 

\blm Definitions 9 and 10 are independent of the choice of BONGs.
\elm
\pf Suppose $L\ap [b_1,\ldots,b_n]$ relative to another good BONG. By
[B3, Theorem 3.1(iii)] we have $d(a_{1,i}b_{1,i})\geq\a_i$ so the
conditions $d(a_{1,i}X_i)\geq\a_i$ and $d(b_{1,i}X_i)\geq\a_i$ are
equivalent. So $X_i$ approximates $a_{1,i}$ iff it approximates
$b_{1,i}$.

Let $V_i$ be with $\dim V_i=i$ s.t. $\det V_i$ approximates $a_{1,i}$
and hence $b_{1,i}$. Let $V_i\ap [c_1,\ldots,c_i]$. We have $\det
V_i=c_{1,i}$ and so $d(c_{1,i}a_{1,i})\geq\a_i$ and
$d(c_{1,i}b_{1,i})\geq\a_i$.

Suppose that $\a_{i-1}+\a_i>2e$. By Lemma 1.5(i) an even number of the
following statements are true: $[b_1,\ldots,b_{i-1}]\rep
[a_1,\ldots,a_i]$, $[a_1,\ldots,a_{i-1}]\rep [c_1,\ldots,c_i]$,
$[b_1,\ldots,b_{i-1}]\rep [c_1,\ldots,c_i]$,
$(c_{1,i}a_{1,i},a_{1,i-1}b_{1,i-1})_\p =1$. The first follows from
[B3, Theorem 3.1] since $\a_{i-1}+\a_i>2e$. The fourth follows from
$d(c_{1,i}a_{1,i})+d(a_{1,i-1}b_{1,i-1})\geq\a_i+\a_{i-1}>2e$. Thus
the second and third are equivalent. So $[a_1,\ldots,a_{i-1}]\rep V_i$
iff $[b_1,\ldots,b_{i-1}]\rep V_i$. Thus $V_i\sim_l [a_1,\ldots,a_i]$
iff $V_i\sim_l [b_1,\ldots,b_i]$.

Suppose that $\a_i+\a_{i+1}>2e$. By Lemma 1.5(i) an even number of the
following sentences are true: $[b_1,\ldots,b_i]\rep
[a_1,\ldots,a_{i+1}]$, $[c_1,\ldots,c_i]\rep [b_1,\ldots,b_{i+1}]$,
$[c_1,\ldots,c_i]\rep [a_1,\ldots,a_{i+1}]$,
$(b_{1,i}c_{1,i},a_{i+1}b_{1,i+1})_\p =1$. The first follows from
[B3, Theorem 3.1] since $\a_i+\a_{i+1}>2e$ and the fourth is true
because
$d(c_{1,i}b_{1,i})+d(a_{i+1}b_{1,i+1})\geq\a_i+\a_{i+1}>2e$. So the
second and the third are equivalent, i.e. $V_i\rep
[a_1,\ldots,a_{i+1}]$ iff $V_i\rep [b_1,\ldots,b_{i+1}]$. Thus
$V_i\sim_r [a_1,\ldots,a_i]$ iff $V_i\sim_r [b_1,\ldots,b_i]$. \qed

Let now $M,N$ be like in $\S 2$. We now show that conditions
(i)-(iv) can be written in terms of Jordan decompositions of $M,N$. We
know from [B3] how to write the numbers $R_i,S_i,\a_1,\b_i$ in terms
of Jordan decompositions.

Let $X_i,Y_i,V_i,W_i$ be approximations for
$a_{1,i},b_{1,i},[a_1,\ldots,a_i]$ and $[b_1,\ldots,b_i]$.

\bff For any $\e\in\ff$ we have $d[\e a_{1,i}b_{1,j}]=d[\e
X_iY_j]:=\min\{ d(\e X_iY_j)\a_i,\b_j\}$. Indeed $d(\e
X_iY_j)\geq\min\{d(\e
a_{1,i}b_{1,j}),d(a_{1,i}X_i),d(b_{1,j}Y_j)\}\geq\min\{ d(\e
a_{1,i}b_{1,j}),\a_i,\b_j\} =d[\e a_{1,i}b_{1,j}]$. Since also
$\a_i,\b_j\geq d[\e a_{1,i}b_{1,j}]$ we have $d[\e X_iY_j]\geq d[\e
a_{1,i}b_{1,j}]$. Similarly $d[\e a_{1,i}b_{1,j}]\geq d[\e
X_iY_j]$. Note that, together with $d(\e X_iY_j)\geq d[\e
a_{1,i}b_{1,j}]$, we also have $d(\e X_ib_{1,j}),d(\e a_{1,i}Y_j)\geq
d[\e a_{1,i}b_{1,j}]$. 

Similarly, $d[\e a_{i,j}]=d[\e X_{i-1}X_j]:=\min\{ d(\e
X_{i-1}X_j),\a_{i-1},\a_j\}$.

In particular, Definition 4 becomes 
$$A_i=\min\{ (R_{i+1}-S_i)/2+e, R_{i+1}-S_i+d[-X_{i+1}Y_{i-1}],
R_{i+1}+R_{i+2}-S_{i-1}-S_i+d[X_{i+2}Y_{i-2}]\}$$
and, if $n\leq m-2$, then 
$$S_{n+1}+A_{n+1}=\min\{ R_{n+2}+d[-X_{n+2}Y_n],
R_{n+2}+R_{n+3}-S_n+d[X_{n+3}Y_{n-1}]\},$$ 
with the terms that don't make sense ignored.

Since $X_i,Y_j$ can be written in terms of
Jordan decompositions the same happens for the expressions
$d[\e a_{1,i}b_{1,j}]$, $d[\e a_{i,j}]$ and $A_i$.
\eff

\blm Assume that $M,N$ satisfy (i) of the main theorem. Then
conditions (ii)-(iv) are equivalent to:

(ii') $d[X_iY_i]\geq A_i$ for any $1\leq i\leq\min\{ m-1,n\}$.

(iii') $W_{i-1}\rep V_i$ for any $1<i\leq\min\{ m-1,n+1\}$
s.t. $R_{i+1}>S_{i-1}$ and $A_{i-1}+A_i>2e+R_i-S_i$.

(iv') If $1<i\leq\min\{ m-2,n+1\}$ and $S_i\geq R_{i+2}>S_{i-1}+2e\geq
R_{i+1}+2e$ then $W_{i-1}\rep V_{i+1}$.
\elm
\pf We have $d[X_iY_i]=d[a_{1,i}b_{1,i}]$ so (ii) and (ii') are
equivalent.

Suppose now that $R_{i+1}>S_{i-1}$ and $A_{i-1}+A_i>2e+R_i-S_i$. In
order to prove that (iii) and (iii') are equivalent we have to show
that $[b_1,\ldots,b_{i-1}]\rep [a_1,\ldots,a_i]$ iff $W_{i-1}\rep
V_i$. Denote $V_i\ap [a'_1,\ldots,a'_i]$ and $W_i\ap
[b'_1,\ldots,b'_{i-1}]$. Then $a'_{1,i}\det V_i\ap a_{1,i}$ and
$b'_{1,i-1}=\det W_{i-1}\ap b_{1,i-1}$ so
$d(a_{1,i}a'_{1,i})=\geq\a_i$ and
$d(b_{1,i-1}b'_{1,i-1})\geq\b_{i-1}$. Also by
Lemma 1.2 $d(a'_{1,i}b_{1,i})\geq d[a_{1,i}b_{1,i}]\geq A_i$ and
$d(-a'_{1,i}b_{1,i-2})\geq d[-a_{1,i}b_{1,i-2}]$.

We will prove that $[b_1,\ldots,b_{i-1}]\rep
[a_1,\ldots,a_i]$ iff $[b_1,\ldots,b_{i-1}]\rep [a'_1,\ldots,a'_i]$
and\\ $[b_1,\ldots,b_{i-1}]\rep [a'_1,\ldots,a'_i]$ iff
$[b'_1,\ldots,b'_{i-1}]\rep [a'_1,\ldots,a'_i]$.

We now prove that $[b_1,\ldots,b_{i-1}]\rep [a_1,\ldots,a_i]$ iff
$[b_1,\ldots,b_{i-1}]\rep [a'_1,\ldots,a'_i]$. By Lemma 1.5(i) and
(ii) it is enough to show that $[a_1,\ldots,a_{i-1}]\rep
[a'_1,\ldots,a'_i]$ and $(a_{1,i-1}b_{1,i-1},a'_{1,i}a_{1,i})_\p =1$
or $[a'_1,\ldots,a'_i]\rep [a_1,\ldots,a_{i+1}]$ and
$(a_{1,i}a'_{1,i},-a_{1,i+1}b_{1,i-1})_\p =1$.

By Lemma 2.18(i) we have $\a_i+A_{i-1}>2e$ or
$\a_i+d[-a_{1,i+1}b_{1,i-1}]>2e$. Suppose first that
$\a_i+A_{i-1}>2e$. Then $\a_{i-1}+\a_i\geq A_{i-1}+\a_i>2e$, which
implies $[a_1,\ldots,a_{i-1}]\rep V_i\ap [a'_1,\ldots,a'_i]$ (we have
$V_i\sim_l[a_1,\ldots,a_i]$) and also
$d(a_{1,i-1}b_{1,i-1})+d(a'_{1,i}a_{1,i})>A_{i-1}+\a_i>2e$ so
$(a_{1,i-1}b_{1,i-1},a'_{1,i}a_{1,i})_\p =1$ and we are done.

Suppose now $\a_i+d[-a_{1,i+1}b_{1,i-1}]>2e$. We have $i=m-1$ or
$\a_i+\a_{i+1}\geq\a_i+d[-a_{1,i+1}b_{1,i-1}]$, which implies
$[a'_1,\ldots,a'_i]\ap V_i\rep [a_1,\ldots,a_{i+1}]$ (we have
$V_i\sim_r[a_1,\ldots,a_i]$) and also
$d(a_{1,i}a'_{1,i})+d(-a_{1,i+1}b_{1,i-1})
\geq\a_i+d[-a_{1,i+1}b_{1,i-1}]>2e$ so
$(a'_{1,i}a_{1,i},-a_{1,i+1}b_{1,i-1})_\p =1$ and we are done.

We now prove that $[b_1,\ldots,b_{i-1}]\rep [a'_1,\ldots,a'_i]$ iff
$[b'_1,\ldots,b'_{i-1}]\rep [a'_1,\ldots,a'_i]$. By Lemma 1.5(i) and
(iii) it is enough to show that $[b'_1,\ldots,b'_{i-1}]\rep
[b_1,\ldots,b_i]$ and $(b_{1,i-1}b'_{1,i-1},a'_{1,i}b_{1,i})_\p =1$ or
$[b_1,\ldots,b_{i-2}]\rep [b'_1,\ldots,b'_{i-1}]$ and
$(b'_{1,i-1}b_{1,i-1},-a'_{1,i}b_{1,i-2})_\p =1$.

By Lemma 2.18(ii) we have $\b_{i-1}+A_i>2e$ or
$\b_{i-1}+d[-a_{1,i}b_{1,i-2}]>2e$. Suppose first that
$\b_{i-1}+A_i>2e$. Then $\b_{i-1}+\b_i\geq A_{i-1}+\b_i>2e$, which
implies $[b'_1,\ldots,b'_{i-1}]\ap W_{i-1}\rep [b_1,\ldots,b_i]$ (we
have $W_{i-1}\sim_r[b_1,\ldots,b_{i-1}]$) and also
$d(b_{1,i-1}b'_{1,i-1})+d(a'_{1,i}b_{1,i})\geq \b_{i-1}+A_i>2e$ so
$(b_{1,i-1}b'_{1,i-1},a'_{1,i}b_{1,i})_\p =1$ and we are
done. 

Suppose now $\b_{i-1}+d[-a_{1,i}b_{1,i-2}]>2e$. Then $i=2$ or
$\b_{i-2}+\b_{i-1}\geq\b_{i-1}+d[-a_{1,i}b_{1,i-2}]>2e$, which implies
$[b_1,\ldots,b_{i-2}]\rep W_{i-1}\ap [b'_1,\ldots,b'_{i-1}]$ (we have
$W_{i-1}\sim_l[b_1,\ldots,b_{i-1}]$) and
$d(b'_{1,i-1}b_{1,i-1})+d(-a'_{1,i}b_{1,i-2})
\geq\b_{i-1}+d[-a_{1,i}b_{1,i-2}]>2e$ so\\
$(b'_{1,i-1}b_{1,i-1},-a'_{1,i}b_{1,i-2})_\p =1$ and we are
done. 

Finally, for the equivalence between (iv) and (iv') assume that
$S_i\geq R_{i+2}>S_{i-1}+2e\geq R_{i+1}+2e$. Then $R_{i+2}-R_{i+1}>2e$
so $\a_{i+1}>2e$. Since $d(\det
V_{i+1}a_{1,i+1})\geq\a_{i+1}>2e$ we have $\det
V_{i+1}a_{1,i+1}\in\fs$. Also $\a_{i+1}+\a_{i+2}\geq\a_{i+1}>2e$ so
$V_{i+1}\rep [a_1,\ldots,a_{i+2}]$. Together with $\det
V_{i+1}\a_{1,i+1}\in\fs$, this implies $V_{i+1}\ap
[a_1,\ldots,a_{i+1}]$. Similarly $S_i-S_{i-1}>2e$ implies $W_{i-1}\ap
[b_1,\ldots,b_{i-1}]$. Thus $[b_1,\ldots,b_{i-1}]\rep
[a_1,\ldots,a_{i+1}]$ iff $W_{i-1}\rep V_{i+1}$. \qed 

{\bf Remark} As seen from the proof of (iii)$\zz$(iii') if
$\a_i+A_{i-1}>2e$ we can replace the condition that $V_i\sim
[a_1,\ldots,a_i]$ by $V_i\sim_l[a_1,\ldots,a_i]$ and if
$\a_i+d[-a_{1,i+1}b_{1,i-1}]>2e$ we can replace it by
$V_i\sim_r[a_1,\ldots,a_i]$. Similarly, if $\b_{i-1}+A_i>2e$ we can
replace the condition that $W_{i-1}\sim [b_1,\ldots,b_{i-1}]$ by
$W_{i-1}\sim_r[b_1,\ldots,b_{i-1}]$ and if
$\b_{i-1}+d[-a_{1,i}b_{1,i-2}]>2e$ we can replace it by
$W_{i-1}\sim_l[b_1,\ldots,b_{i-1}]$.

\bco The main theorem is independent of the BONGs.
\eco
\pf If $M\ap [a'_1,\ldots,a'_m]$ and $N\ap [b'_1,\ldots,b'_n]$
relative to other good BONGs then $X_i:=a'_{1,i}$, $Y_i:=b'_{1,i}$,
$V_i:=[a'_1,\ldots,a'_i]$ and $W_i:=[b'_1,\ldots,b'_i]$ are
approximations for $a_{1,i},b_{1,i},[a_1,\ldots,a_i]$ and
$[b_1,\ldots,b_i]$. (See [B3, Theorem 3.1(iii) and (iv)].) 

If we apply Lemma 3.10 with the approximations above we 
get the main theorem in terms of the new BONGs. \qed

\section{Proof of transitivity}

Let $M\ap\[ a_1,\ldots,a_n\]$, $N\ap\[\ b_1,\ldots,b_n\]$ and $K\ap\[
c_1,\ldots,c_n\]$. We denote by $R_i,S_i,T_i$ the $R_i$'s corresponding
to $M,N,K$ and by $\a_i,\b_i,\c_i$ their $\a_i$'s. Let $A_i=A_i(M,N)$,
$B_i=A_i(N,K)$ and $C_i=A_i(M,K)$.

In this section we will prove that $N\leq M$ and $K\leq N$ imply
$K\leq M$. We will use the techniques and results from \S2.

Although the transitivity could be proved in the general case we can
restrict ourselves to the case when the three lattices
have the same rank, which allows the use of duality. This can cut the
number of cases to be treated by half.

\subsection{Duality}

\bff For the dual lattices we have $K^\*\ap\[ a_1^\*,\ldots,a_n^\*\]$,
$N^\*\ap\[\ b_1^\*,\ldots,b_n^\*\]$ and $M^\*\ap\[ c_1^\*,\ldots,c_n^\*\]$
where $a_i^\* =c_{n+1-i}\1$, $b_i^\* =b_{n+1-i}\1$ and $c_i^\*
=a_{n+1-i}\1$. Let $R_i^\*,S_i^\*,T_i^\*$ be the $R_i$'s corresponding to
$K^\*,N^\*,M^\*$ and by $\a_i^\*,\b_i^\*,\c_i^\*$ their $\a_i$'s. Let
$A_i^\* =A_i(K^\*,N^\* )$, $B_i=A_i(N^\*,M^\* )$ and
$C_i=A_i(K^\*,M^\* )$.

We have $R_i^\* =-T_{n+1-i}$, $S_i^\* =-S_{n+1-i}$, $T_i^\*
=-R_{n+1-i}$, $\a_i^\* =\c_{n-i}$, $\b_i^\* =\b_{n-i}$, $\c_i^\*
=\a_{n-i}$, $A_i^\* = B_{n-i}$, $B_i^\* = A_{n-i}$ and $C_i^\* =
C_{n-i}$.

Like in 2.4, in many cases, given a formula or a statement for
$M,N,K$ at some index $i$, we will consider the same formula or
statement for $K^\*,N^\*,M^\*$ at index $n+1-i$. When we do so and
then consider the equivalent statements in terms of $M,N,K$ we get:

$R_{i+k},S_{i+k},T_{i+k}$ become $-T_{i-k},-S_{i-k},-R_{i-k}$.

$d[\e a_{i+k,i+l}]$, $d[\e b_{i+k,i+l}]$, $d[\e c_{i+k,i+l}]$ become
$d[\e c_{i-l,i-k}]$, $d[\e b_{i-l,i-k}]$, $d[\e a_{i-l,i-k}]$.

$\a_{i+k},\b_{i+k},\c_{i+k},A_{i+k},B_{i+k},C_{i+k}$ become
$\c_{i-k-1},\b_{i-k-1}\a_{i-k-1},B_{i-k-1},A_{i-k-1},C_{i-k-1}$.

$d[\e a_{1,i+k}b_{1,i+l}]$, $d[\e b_{1,i+k}c_{1,i+l}]$, $d[\e
a_{1,i+k}c_{1,i+l}]$,  become
$d[\e b_{1,i-l-1}c_{1,i-k-1}]$,\\ $d[\e a_{1,i-l-1}b_{1,i-k-1}]$,
$d[\e a_{1,i-l-1}c_{1,i-k-1}]$,

$[b_1,\ldots,b_{i+k}]\rep [a_1,\ldots,a_{i+l}]$, $[c_1,\ldots,c_{i+k}]\rep
[b_1,\ldots,b_{i+l}]$, $[c_1,\ldots,c_{i+k}]\rep [a_1,\ldots,a_{i+l}]$ become
$[c_1,\ldots,c_{i-k-1}]\rep [b_1,\ldots,b_{i-l-1}]$,
$[b_1,\ldots,b_{i-k-1}]\rep [a_1,\ldots,a_{i-l-1}]$,\\
$[c_1,\ldots,c_{i-k-1}]\rep [a_1,\ldots,a_{i-l-1}]$.

The property that $i+k$ is essential for $M,N$, for $N,K$ or $M,K$
becomes $i-k$ is essential for $N,K$, for $M,N$ resp. $M,K$. (See the
remark following Definition 7.)

Sometimes instead of $n+1-i$ we take $n-i$. In this case all indices
in the dual formulas and statements described above are increased by
$1$.
\eff

\vskip 0.5 cm

We assume that $N\leq M$ and $K\leq N$. We want to prove that $K\leq
M$. We prove that $M,K$ satisfy the conditions (i)-(iv) of the main
theorem.

{\bf Proof of 2.1(i)} Condition 2.1(i) for $R_i$'s and $T_i$'s follows
from the corresponding condition for $R_i$'s and $S_i$'s and that for
$S_i$'s si $T_i$'s. (We have $\rr (M)\leq\rr (N)\leq\rr (K)$ so $\rr
(M)\leq\rr (K)$.)

\subsection{A key lemma}

For the proof of (ii) and (iii) we need the following key lemma:

\blm If $i$ is an essential index for $M$ and $K$ then:

(i) If $i>1$ then $C_{i-1}\leq A_{i-1},B_{i-1}$ if
$R_{i+1}+S_i>T_{i-2}+T_{i-1}$ or $i\in\{ 2,n\}$ and $C_{i-1}\leq
R_i-R_{i+1}+A_i$ otherwise.

(ii) If $i<n$ then $C_i\leq A_i,B_i$ if $R_{i+1}+R_{i+2}>S_i+T_{i-1}$
or $i\in\{ 1,n-1\}$ and $C_i\leq T_{i-1}-T_i+B_{i-1}$ otherwise.
\elm
\pf We only have to prove (i) since (ii) will follow by duality at
index $n+1-i$.

Condition that $i$ is essential means $R_{i+1}>T_{i-1}$ (if $2\leq
i\leq n-1$) and $R_{i+1}+R_{i+2}>T_{i-2}+T_{i-1}$ (if $3\leq i\leq
n-2$).

Suppose  $R_{i+1}+S_i>T_{i-2}+T_{i-1}$ or $i\in\{ 2,n\}$. This means
$R_{i+1}+S_i>T_{i-2}+T_{i-1}$ whenever this makes sense. (If $i=2$
then $T_{i-2}$ is not defined. If $i=n$ then $R_{i+1}$ is not
defined.)

We have $S_i+S_{i+1}>T_{i-2}+T_{i-1}$, provided that this makes sense
(i.e. when $i\notin\{ 2, n\}$). Indeed, if $S_{i+1}\geq R_{i+1}$ then
$S_i+S_{i+1}\geq R_{i+1}+S_i>T_{i-2}+T_{i-1}$. Otherwise
$S_i+S_{i+1}\geq R_{i+1}+R_{i+2}>T_{i-2}+T_{i-1}$.

Also $R_{i+1}>S_{i-1}$, provided that this makes sense (i.e. when
$i\neq n$). Indeed, otherwise $S_{i-1}\geq R_{i+1}>T_{i-1}$ so
$S_{i-1}+S_i\leq T_{i-2}+T_{i-1}<R_{i+1}+S_i$ so
$R_{i+1}>S_{i-1}$. This implies $R_i\leq S_i$, which also happens if
$i=n$. Also, by Lemma 2.7.(i) and (ii), in the formula for $A_{i-1}$
we can replace $d[a_{1,i+1}b_{1,i-3}]$ by $d[-a_{1,i+1}b_{1,i-1}]$
and in the formula for $A_i$ we can replace $d[a_{1,i+2}b_{1,i-2}]$
by $d[-a_{1,i}b_{1,i-2}]$.

Suppose that $C_{i-1}>A_{i-1}=\min\{ (R_i-S_{i-1})/2+e,
R_i-S_{i-1}+d[-a_{1,i}b_{1,i-2}],
R_i+R_{i+1}-S_{i-2}-S_{i-1}+d[-a_{1,i+1}b_{1,i-1}]\}$.

Suppose first that $A_{i-1}=(R_i-S_{i-1})/2+e$. Since
$(R_i-T_{i-1})/2+e,R_i-(T_{i-2}+T_{i-1})/2+e\geq C_{i-1}>A_{i-1}$ we
get $S_{i-1}>T_{i-1}$ and $R_i+S_{i-1}>T_{i-2}+T_{i-1}$. But the first
inequality implies $T_{i-2}+T_{i-1}\geq S_i+S_{i-1}\geq
R_i+S_{i-1}$. Contradiction.

Suppose now that
$A_{i-1}=R_i+R_{i+1}-S_{i-2}-S_{i-1}+d[-a_{1,i+1}b_{1,i-1}]$. Since
$R_i+R_{i+1}-S_{i-2}-S_{i-1}+d[a_{1,i+1}c_{1,i-3}]\geq
R_i+R_{i+1}-T_{i-2}-T_{i-1}+d[a_{1,i+1}c_{1,i-3}]\geq C_{i-1}>A_{i-1}$
we get $d[a_{1,i+1}c_{1,i-3}]>d[-a_{1,i+1}b_{1,i-1}]$ so
$d[-a_{1,i+1}b_{1,i-1}]=d[-b_{1,i-1}c_{1,i-3}]$. Also
$R_i+R_{i+1}-S_{i-2}-S_{i-1}+d[-c_{i-2,i-1}]\geq
R_i+R_{i+1}-T_{i-2}-T_{i-1}+d[-c_{i-2,i-1}]\geq
R_i+R_{i+1}-2T_{i-1}+\c_{i-2}>R_i-T_{i-1}+\c_{i-2}\geq
C_{i-1}>A_{i-1}$ so
$d[-c_{i-2,i-1}]>d[-a_{1,i+1}b_{1,i-1}]=d[-b_{1,i-1}c_{1,i-3}]$. Thus
$d[-a_{1,i+1}b_{1,i-1}]=d[b_{1,i-1}c_{1,i-1}]\geq B_{i-1}$. Hence
$A_{i-1}\geq R_i+R_{i+1}-S_{i-2}-S_{i-1}+B_{i-1}$. Since
$S_i+S_{i+1}>T_{i-2}+T_{i-1}$ we have $B_{i-1}=\( B_{i-1}=\min\{
(S_i-T_{i-1})/2+e, S_i-T_{i-1}+d[-b_{1,i}c_{1,i-2}],
S_i+S_{i+1}-2T_{i-1}+\c_{i-2}, 2S_i-T_{i-2}-T_{i-1}+\b_i\}$. If
$B_{i-1}=(S_i-T_{i-1})/2+e$ then $R_i-(T_{i-2}+T_{i-1})/2+e\geq
C_{i-1}>A_{i-1}\geq R_i+R_{i+1}-S_{i-2}-S_{i-1}+B_{i-1}\geq
R_i+R_{i+1}-T_{i-2}-T_{i-1}+B_{i-1}=
R_i+R_{i+1}-T_{i-2}-T_{i-1}+(S_i-T_{i-1})/2+e$. This implies
$T_{i-2}+2T_{i-1}>2R_{i+1}+S_i>R_{i+1}+S_i+T_{i-1}$ so
$T_{i-2}+T_{i-1}>R_{i+1}+S_i$. Contradiction. If
$B_{i-1}=2S_i-T_{i-2}-T_{i-1}+\b_i\geq
S_{i-2}+S_i-T_{i-2}-T_{i-1}+\b_{i-2}$ then
$A_{i-1}\geq R_i+R_{i+1}-S_{i-2}-S_{i-1}+B_{i-1}\geq
R_i+R_{i+1}+S_i-S_{i-1}-T_{i-2}-T_{i-1}+\b_{i-2}>R_i-S_{i-1}+\b_{i-2}$.
Contradiction. If $B_{i-1}=S_i+S_{i+1}-2T_{i-1}+\c_{i-2}\geq
S_{i-2}+S_{i-1}-2T_{i-1}+\c_{i-2}$ then
$C_{i-1}>R_i+R_{i+1}-S_{i-2}-S_{i-1}+B_{i-1}\geq
R_i+R_{i+1}-2T_{i-1}+\c_{i-2}>R_i-T_{i-1}+\c_{i-2}$. Contradiction.

Hence $B_{i-1}=S_i-T_{i-1}+d[-b_{1,i}c_{1,i-2}]$ so
$C_{i-1}>A_{i-1}\geq
R_i+R_{i+1}-S_{i-2}-S_{i-1}+S_i-T_{i-1}+d[-b_{1,i}c_{1,i-2}]\geq
R_i+R_{i+1}-S_{i-1}-T_{i-1}+d[-b_{1,i}c_{1,i-2}]$. Now
$R_i+R_{i+1}-S_{i-1}-T_{i-1}+d[-a_{1,i}c_{1,i-2}]>
R_i-T_{i-1}+d[-a_{1,i}c_{1,i-2}]\geq
C_{i-1}>R_i+R_{i+1}-S_{i-1}-T_{i-1}+d[-b_{1,i}c_{1,i-2}]$. Thus
$d[-a_{1,i}c_{1,i-2}]>d[-b_{1,i}c_{1,i-2}]$ so
$d[-b_{1,i}c_{1,i-2}]=d[a_{1,i}b_{1,i}]$. Also
$R_i+R_{i+1}-S_{i-1}-T_{i-1}+d[-a_{1,i}b_{1,i-2}]>
R_i-S_{i-1}+d[-a_{1,i}b_{1,i-2}]\geq A_{i-1}\geq
R_i+R_{i+1}-S_{i-1}-T_{i-1}+d[-b_{1,i}c_{1,i-2}]$. Thus
$d[-a_{1,i}b_{1,i-2}]>d[-b_{1,i}c_{1,i-2}]=d[a_{1,i}b_{1,i}]$. It
follows that $d[-b_{1,i}c_{1,i-2}]=d[-b_{i-1,i}]\geq
S_{i-1}-S_i+\b_{i-1}\geq S_{i-2}-S_i+\b_{i-2}$. Hence $A_{i-1}\geq
R_i+R_{i+1}-S_{i-2}-S_{i-1}+S_i-T_{i-1}+d[-b_{1,i}c_{1,i-2}]\geq
R_i+R_{i+1}-S_{i-1}-T_{i-1}+\b_{i-2}>R_i-S_{i-1}+\b_{i-2}\geq
A_{i-1}$. Contradiction.

In conclusion $A_{i-1}=R_i-S_{i-1}+d[-a_ib_{i-2}]$. We have 2 cases:

a. $S_i\geq T_{i-2}$ or $i=2$. It follows that $S_{i-1}\leq T_{i-1}$
so $R_i-S_{i-1}+d[-a_{1,i}c_{1,i-2}]\geq
R_i-T_{i-1}+d[-a_{1,i}c_{1,i-2}]\geq
C_{i-1}>A_{i-1}=R_i-S_{i-1}+d[-a_{1,i}b_{1,i-2}]$. Thus
$d[-a_{1,i}c_{1,i-2}]>d[-a_{1,i}b_{1,i-2}]$ so
$d[-a_{1,i}b_{1,i-2}]=d[b_{1,i-2}c_{i-2}]\geq B_{i-2}$. So
$C_{i-1}>A_{i-1}\geq R_i-S_{i-1}+B_{i-2}=\min\{
R_i-(S_{i-1}+T_{i-2})/2+e, R_i-T_{i-2}+d[-b_{1,i-1}c_{1,i-3}],
R_i+S_i-T_{i-3}-T_{i-2}+d[b_{1,i}c_{1,i-4}]\}$. But
$R_i-(S_{i-1}+T_{i-2})/2+e\geq R_i-(T_{i-2}+T_{i-1})/2+e\geq C_{i-1}$
so it can be ignored.

If $R_i-S_{i-1}+B_{i-2}=R_i+S_i-T_{i-3}-T_{i-2}+d[b_{1,i}c_{1,i-4}]$
note that $R_i+S_i-T_{i-3}-T_{i-2}+d[-c_{i-3,i-2}]\geq
R_i+S_i-2T_{i-2}+\c_{i-3}\geq R_i-T_{i-2}+\c_{i-3}\geq
R_i-T_{i-1}+\c_{i-2}\geq C_{i-1}$. Also
$R_i+S_i-T_{i-3}-T_{i-2}+d[-a_{1,i}c_{1,i-2}]\geq
R_i-T_{i-1}+d[-a_{1,i}c_{1,i-2}]\geq C_{i-1}$. (We have $T_{i-1}\geq T_{i-3}$
and $S_i\geq T_{i-2}$.) Since
$C_{i-1}>R_i-S_{i-1}+B_{i-2}=R_i+S_i-T_{i-3}-T_{i-2}+d[b_{1,i}c_{1,i-4}]$
we have $d[-c_{i-3,i-2}],d[-a_{1,i}c_{1,i-2}]>d[b_{1,i}c_{1,i-4}]$ so
$d[b_{1,i}c_{1,i-4}]=d[a_{1,i}b_{1,i}]\geq A_i$. Thus
$R_i-S_{i-1}+B_{i-2}\geq R_i+S_i-T_{i-3}-T_{i-2}+A_i\geq
R_i+S_i-T_{i-2}-T_{i-1}+A_i$. 

If $R_i-S_{i-1}+B_{i-2}=R_i-T_{i-2}+d[-b_{1,i-1}c_{1,i-3}]$ then note
that $R_i-T_{i-2}+d[-c_{i-2,i-1}]\geq R_i-T_{i-1}+\c_{i-2}\geq
C_{i-1}>R_i-S_{i-1}+B_{i-2}$. Thus
$d[-c_{i-2,i-1}]>d[-b_{1,i-1}c_{1,i-3}]$. So
$d[-b_{1,i-1}c_{1,i-3}]=d[b_{1,i-1}c_{1,i-1}]\geq B_{i-1}$. Now
$S_i+S_{i+1}>T_{i-2}+T_{i-1}$ and $S_i\geq T_{i-2}$ so $B_{i-1}=\(
B_{i-1}=\min\{ (S_i-T_{i-1})/2+e, S_i-T_{i-1}+d[-b_{1,i}c_{1,i-2}],
S_i+S_{i+1}-2T_{i-1}+\c_{i-2}\}$ and $R_i-T_{i-2}+B_{i-1}=\min\{
R_i-T_{i-2}+(S_i-T_{i-1})/2+e, R_i+S_i-T_{i-2}-T_{i-1}+d[-b_{1,i}c_{1,i-2}],
R_i+S_i+S_{i+1}-T_{i-2}-2T_{i-1}+\c_{i-2}\}$. But $S_i\geq T_{i-2}$ so
$R_i-T_{i-2}+(S_i-T_{i-1})/2+e\geq R_i-(T_{i-2}+T_{i-1})/2+e\geq
C_{i-1}>R_i-T_{i-2}+B_{i-1}$ so it can be ignored. Also
$S_i+S_{i+1}>T_{i-2}+T_{i-1}$ so
$R_i+S_i+S_{i+1}-T_{i-2}-2T_{i-1}+\c_{i-2}>R_i-T_{i-1}+\c_{i-2}\geq
C_{i-1}$ so it can also be removed. Thus
$C_{i-1}>R_i-T_{i-2}+B_{i-1}=R_i+S_i-T_{i-2}-T_{i-1}+d[-b_{1,i}c_{1,i-2}]$.
Now $R_i+S_i-T_{i-2}-T_{i-1}+d[-a_{1,i}c_{1,i-2}]\geq
R_i-T_{i-1}+d[-a_{1,i}c_{1,i-2}]\geq C_{i-1}$ so
$d[-a_{1,i}c_{1,i-2}]>d[-b_{1,i}c_{1,i-2}]$. This implies
$d[-b_{1,i}c_{1,i-2}]=d[a_{1,i}b_{1,i}]\geq A_i$. So
$R_i-S_{i-1}+B_{i-2}\geq R_i+S_i-T_{i-2}-T_{i-1}+A_i$.

In both cases we obtained $C_{i-1}>A_{i-1}\geq R_i-S_{i-1}+B_{i-2}\geq
R_i+S_i-T_{i-2}-T_{i-1}+A_i=\min\{
R_i+(R_{i+1}+S_i)/2-T_{i-2}-T_{i-1}+e,
R_i+R_{i+1}-T_{i-2}-T_{i-1}+d[-a_{1,i+1}b_{1,i-1}],
R_i+R_{i+1}+R_{i+2}-S_{i-1}-T_{i-2}-T_{i-1}+d[-a_{1,i}b_{1,i-2}]\}$.
Now $R_i+(R_{i+1}+S_i)/2-T_{i-2}-T_{i-1}+e>
R_i+(T_{i-2}+T_{i-1})/2-T_{i-2}-T_{i-1}+e=
R_i-(T_{i-2}+T_{i-1})/2+e\geq C_{i-1}$ so it can be removed. Also
$R_i+R_{i+1}+R_{i+2}-S_{i-1}-T_{i-2}-T_{i-1}+d[-a_{1,i}b_{1,i-2}]\geq
R_{i+1}+R_{i+2}-T_{i-2}-T_{i-1}+A_{i-1}>A_{i-1}$ so this one can also
be removed. Hence
$R_i+R_{i+1}-T_{i-2}-T_{i-1}+d[a_{1,i+1}c_{1,i-3}]\geq
C_{i-1}>R_i+S_i-T_{i-2}-T_{i-1}+A_i=
R_i+R_{i+1}-T_{i-2}-T_{i-1}+d[-a_{1,i+1}b_{1,i-1}]$. Thus
$d[a_{1,i+1}c_{1,i-3}]>d[-a_{1,i+1}b_{1,i-1}]$ so
$d[-a_{1,i+1}b_{1,i-1}]=d[-b_{1,i-1}c_{1,i-3}]\geq
T_{i-2}-S_{i-1}+B_{i-2}$, which implies $R_i+S_i-T_{i-2}-T_{i-1}+A_i=
R_i+R_{i+1}-T_{i-2}-T_{i-1}+d[-b_{1,i-1}c_{1,i-3}]\geq
R_i+R_{i+1}-S_{i-1}-T_{i-1}+B_{i-2}>R_i-S_{i-1}+B_{i-2}$.
Contradiction.

b. $S_i<T_{i-2}$. It follows that $S_{i-1}+S_i\leq
T_{i-2}+T_{i-1}<R_{i+1}+R_{i+2}$. Also $R_{i+1}>S_{i-1}$ so
$A_i=\( A_i=\min\{ (R_{i+1}-S_i)/2+e,
R_{i+1}-S_i+d[-a_{1,i+1}b_{1,i-1}],
R_{i+1}+R_{i+2}-2S_i+\b_{i-1}\}$. (If $i=n-1$ then we ignore
$R_{i+1}+R_{i+2}-2S_i+\b_{i-1}$. If $i=n$ then $A_i$ is not defined.)
We have $S_i+S_{i+1}>T_{i-2}+T_{i-1}$ and $S_i<T_{i-2}$ so
$S_{i+1}>T_{i-1}$. By Lemmas 2.14 and 2.9 we have $B_{i-1}=B'_{i-1}=\(
B'_{i-1}=\min\{ S_i-T_{i-1}+d[-b_{1,i}c_{i-2}],
2S_i-T_{i-1}-T_{i-2}+\b_i\}$. (If $i=n$ we ignore
$2S_i-T_{i-1}-T_{i-2}+\b_i$.) 

We have $C_{i-1}>A_{i-1}=R_i-S_{i-1}+d[-a_{1,i}b_{1,i-2}]\geq\min\{
R_i-S_{i-1}+d[-b_{i-1,i}], R_i-S_{i-1}+d[a_{1,i}b_{1,i}]\}\geq\min\{
R_i-S_i+\b_{i-1}, R_i-S_{i-1}+A_i\}\geq\min\{ R_i-S_i+B_{i-1},
R_i-S_{i-1}+A_i\} =\min\{ R_i-T_{i-1}+d[-b_{1,i}c_{1,i-2}],
R_i+S_i-T_{i-2}-T_{i-1}+\b_i, R_i-S_{i-1}+A_i\}$. Since
$R_i-T_{i-1}+d[-a_{1,i}c_{1,i-2}]\geq C_{i-1}$ we can replace above
$R_i-T_{i-1}+d[-b_{1,i}c_{1,i-2}]$ by $R_i-T_{i-1}+d[a_{1,i}b_{1,i}]$,
which is $\geq R_i-T_{i-1}+A_i$. Also
$R_i+S_i-T_{i-2}-T_{i-1}+\b_i\geq R_i+S_i-T_{i-2}-T_{i-1}+A_i$. Thus
$\min\{ R_i-S_i+B_{i-1}, R_i-S_{i-1}+A_i\}\geq\min\{ R_i-T_{i-1}+A_i,
R_i+S_i-T_{i-2}-T_{i-1}+A_i, R_i-S_{i-1}+A_i\}
=R_i+S_i-T_{i-2}-T_{i-1}+A_i$. (We have $S_i<T_{i-2}$ and
$S_{i-1}+S_i\leq T_{i-2}+T_{i-1}$ so
$R_i-T_{i-1}+A_i,R_i-S_{i-1}+A_i\geq R_i+S_i-T_{i-2}-T_{i-1}+A_i$.) It
follows that $C_{i-1}>A_{i-1}\geq R_i+S_i-T_{i-2}-T_{i-1}+A_i$. Also
$R_i-S_i+B_{i-1}\geq R_i+S_i-T_{i-2}-T_{i-1}+A_i$. Now
$R_i+S_i-T_{i-2}-T_{i-1}+A_i=\min\{
R_i+(R_{i+1}+S_i)/2-T_{i-2}-T_{i-1}+e,
R_i+R_{i+1}-T_{i-2}-T_{i-1}+d[-a_{1,i+1}b_{1,i-1}],
R_i+R_{i+1}+R_{i+2}-S_i-T_{i-2}-T_{i-1}+\b_{i-1}\}$. We have
$R_i+(R_{i+1}+S_i)/2-T_{i-2}-T_{i-1}+e>
R_i+(T_{i-2}+T_{i-1})/2-T_{i-2}-T_{i-1}+e=R_i-(T_{i-2}+T_{i-1})/2+e\geq
C_{i-1}$ so it can be removed. Also
$R_i+R_{i+1}+R_{i+2}-S_i-T_{i-2}-T_{i-1}+\b_{i-1}>R_i-S_i+\b_{i-1}\geq
R_i-S_i+B_{i-1}$ so it also can be removed. Hence
$R_i+S_i-T_{i-2}-T_{i-1}+A_i=
R_i+R_{i+1}-T_{i-2}-T_{i-1}+d[-a_{1,i+1}b_{1,i-1}]$. But
$R_{i+1}>T_{i-1}$ so $d[a_{1,i+1}c_{1,i-3}]$ in the formula for
$C_{i-1}$ can be replaced by $d[-a_{1,i+1}c_{1,i-1}]$. Thus
$R_i+R_{i+1}-T_{i-2}-T_{i-1}+d[-a_{1,i+1}c_{1,i-1}]\geq
C_{i-1}>R_i+R_{i+1}-T_{i-2}-T_{i-1}+d[-a_{1,i+1}b_{1,i-1}]$. Thus
$d[-a_{1,i+1}b_{1,i-1}]=d[b_{1,i-1}c_{1,i-1}]\geq B_{i-1}$ so
$R_i+S_i-T_{i-2}-T_{i-1}+A_i\geq
R_i+R_{i+1}-T_{i-2}-T_{i-1}+B_{i-1}>R_i-S_i+B_{i-1}$. Contradiction.
\vskip 0.5cm

So we proved $A_{i-1}\geq C_{i-1}$. Suppose now that $A_{i-1}\geq
C_{i-1}>B_{i-1}$. Since $S_i+S_{i+1}>T_{i-2}+T_{i-1}$ or $i-1\in\{
1,n-1\}$ we have $B_{i-1}=\( B_{i-1}=\min\{ (S_i-T_{i-1})/2+e,
S_i-T_{i-1}+d[-b_{1,i}c_{1,i-2}],2S_i-T_{i-2}-T_{i-1}+\b_i,
S_i+S_{i+1}-2T_{i-1}+\c_{i-2}\}$. (If $i\in\{ 2,n\}$ we ignore
$2S_i-T_{i-2}-T_{i-1}+\b_i$ and $S_i+S_{i+1}-2T_{i-1}+\c_{i-2}$.) But
$S_i+S_{i+1}-2T_{i-1}+\c_{i-2}\geq R_i+R_{i+1}-2T_{i-1}+\c_{i-2}>
R_i-T_{i-1}+\c_{i-2}\geq C_{i-1}>B_{i-1}$ so it can be
removed. Suppose $B_{i-1}=(S_i-T_{i-1})/2+e$. Then
$(R_i-T_{i-1})/2+e\geq C_{i-1}>(S_i-T_{i-1})/2+e$ so $R_i>S_i$, which
implies $R_i+R_{i+1}\leq S_{i-1}+S_i$. Also $(R_i-S_{i-1})/2+e\geq
A_{i-1}>(S_i-T_{i-1})/2+e$ and so $R_i+T_{i-1}>S_{i-1}+S_i\geq
R_i+R_{i+1}$. But this contradicts $R_{i+1}>T_{i-1}$. We have
$R_i-T_{i-1}+d[-a_{1,i}c_{1,i-2}]\geq C_{i-1}>B_{i-1}$ and
$R_i-T_{i-1}\leq S_i-T_{i-1}$. By Lemma 1.4(b), in the formula for
$B_{i-1}$, we can replace $S_i-T_{i-1}+d[-b_{1,i}c_{1,i-2}]$ by
$S_i-T_{i-1}+d[a_{1,i}b_{1,i}]$. Thus $B_{i-1}=\min\{
S_i-T_{i-1}+d[a_{1,i}b_{1,i}],2S_i-T_{i-2}-T_{i-1}+\b_i\}\geq\min\{
S_i-T_{i-1}+A_i,2S_i-T_{i-2}-T_{i-1}+A_i\}$. If
$A_i=(R_{i+1}-S_i)/2+e$ then $B_{i-1}\geq\min\{
(R_{i+1}+S_i)/2-T_{i-1}+e,
(R_{i+1}+S_i)/2+S_i-T_{i-2}-T_{i-1}+e\}$. But we have
$(R_{i+1}+S_i)/2-T_{i-1}+e>(S_i-T_{i-1})/2+e\geq B_{i-1}$ and
$(R_{i+1}+S_i)/2+S_i-T_{i-2}-T_{i-1}+e>
(T_{i-2}+T_{i-1})/2+S_i-T_{i-2}-T_{i-1}+e=S_i-(T_{i-2}+T_{i-1})/2\geq
B_{i-1}$ so we get a contradiction. Hence $A_i=A_i'$. Also
$S_i-T_{i-1},2S_i-T_{i-2}-T_{i-1}>S_i-R_{i+1}$ so
$B_{i-1}>S_i-R_{i+1}+A'_i=\min\{ d[-a_{1,i+1}b_{1,i-1}],
R_{i+2}-S_{i-1}+d[-a_{1,i}b_{1,i-2}]\}$. But
$R_{i+2}-S_{i-1}+d[-a_{1,i}b_{1,i-2}]\geq
R_i-S_{i-1}+d[-a_{1,i}b_{1,i-2}]\geq A_{i-1}>B_{i-1}$ so it can be
removed. Now $d[a_{1,i-1}b_{1,i-1}]\geq
A_{i-1}>S_i-R_{i+1}+A'_i=d[-a_{1,i+1}b_{1,i-1}]$ so
$d[-a_{1,i+1}b_{1,i-1}]=d[-a_{i,i+1}]$. Also
$d[b_{1,i-1}c_{1,i-1}]\geq
C_{i-1}>S_i-R_{i+1}+A'_i=d[-a_{1,i+1}b_{1,i-1}]$ so
$d[-a_{1,i+1}b_{1,i-1}]=d[-a_{1,i+1}c_{1,i-1}]$. Hence
$S_i-R_{i+1}+A'_i=d[-a_{1,i+1}b_{1,i-1}]=d[-a_{i,i+1}]
=d[-a_{1,i+1}c_{1,i-1}]$,which implies
$A_i=A_i'=R_{i+1}-S_i+d[-a_{i,i+1}]=R_{i+1}-S_i+d[-a_{1,i+1}c_{1,i-1}]$.
It follows that $C_{i-1}>B_{i-1}\geq\min\{
S_i-T_{i-1}+A_i,2S_i-T_{i-2}-T_{i-1}+A_i\} =\min\{
R_{i+1}-T_{i-1}+d[-a_{i,i+1}],
R_{i+1}+S_i-T_{i-2}-T_{i-1}+d[-a_{1,i+1}c_{1,i-1}]\}$. But
$R_{i+1}-T_{i-1}+d[-a_{i,i+1}]\geq R_i-T_{i-1}+\a_i\geq
C_{i-1}$ and $R_{i+1}+S_i-T_{i-2}-T_{i-1}+d[-a_{1,i+1}c_{1,i-1}]\geq
R_i+R_{i+1}-T_{i-2}-T_{i-1}+d[-a_{1,i+1}c_{1,i-1}]\geq
C_{i-1}$. ($R_{i+1}>T_{i-1}$ so $d[a_{1,i+1}c_{1,i-3}]$ can be
replaced by $d[-a_{1,i+1}c_{1,i-1}]$ in $C_{i-1}$.) Contradiction. So
$B_{i-1}\geq C_{i-1}$.
\vskip 0.5cm

We consider now the case $R_{i+1}+S_i\leq T_{i-2}+T_{i-1}$. Since also
$R_{i+1}>T_{i-1}$ we have $S_i<T_{i-2}$. Hence $S_{i-1}+S_i\leq
T_{i-2}+T_{i-1}<R_{i+1}+R_{i+2}$. Also $R_{i+1}+S_i\leq
T_{i-2}+T_{i-1}<R_{i+1}+R_{i+2}$ so $S_i<R_{i+2}$. It follows that
$R_{i+1}\leq S_{i+1}$ and also $A_i=\( A_i=\min\{ (R_{i+1}-S_i)/2+e,
R_{i+1}-S_i+d[-a_{1,i+1}b_{1,i-1}],
2R_{i+1}-S_{i-1}-S_i+\a_{i+1}\}$. (When $i=n-1$ the inequalities above
involving $R_{i+2}$ are ignored. The inequality $R_{i+1}\leq S_{i+1}$
still holds in this case. Also $2R_{i+1}-S_{i-1}-S_i+\a_{i+1}$ i
ignored if $i=n-1$.)

Suppose that $C_{i-1}>R_i-R_{i+1}+A_i=\min\{
R_i-(R_{i+1}+S_i)/2+e, R_i-S_i+d[-a_{1,i+1}b_{1,i-1}],
R_i+R_{i+1}-S_{i-1}-S_i+\a_{i+1}\}$. But
$R_i+R_{i+1}-S_{i-1}-S_i+\a_{i+1}\geq
R_i+R_{i+1}-T_{i-2}-T_{i-1}+\a_{i+1}\geq C_{i-1}$ so it can be
removed. Also $R_i-(R_{i+1}+S_i)/2+e\geq
R_i-(T_{i-2}+T_{i-1})/2+e\geq C_{i-1}$ so it can be removed and we
have $C_{i-1}>R_i-R_{i+1}+A_i=R_i-S_i+d[-a_{1,i+1}b_{1,i-1}]$. Since
$R_i-S_i+d[-a_{1,i+1}c_{1,i-1}]\geq
R_i+R_{i+1}-T_{i-2}-T_{i-1}+d[-a_{1,i+1}c_{1,i-1}]\geq
C_{i-1}>R_i-S_i+d[-a_{1,i+1}b_{1,i-1}]$ we have
$d[-a_{1,i+1}b_{1,i-1}]=d[b_{1,i-1}c_{1,i-1}]\geq B_{i-1}$ so
$C_{i-1}>R_i-R_{i+1}+A_i\geq R_i-S_i+B_{i-1}=\min\{
R_i-(S_i+T_{i-1})/2+e, R_i-T_{i-1}+d[-b_{1,i}c_{1,i-2}],
R_i+S_{i+1}-T_{i-2}-T_{i-1}+d[b_{1,i+1}c_{1,i-3}]\}$. But
$R_i-(S_i+T_{i-1})/2+e>R_i-(T_{i-2}+T_{i-1})/2+e\geq C_{i-1}$ so
it can be removed. We have $R_i-T_{i-1}+d[-a_{1,i}c_{1,i-2}]\geq
C_{i-1}>R_i-S_i+B_{i-1}$ so $R_i-T_{i-1}+d[-b_{1,i}c_{1,i-2}]$ can
be replaced by $R_i-T_{i-1}+d[a_{1,i}b_{1,i}]$. But
$R_i-T_{i-1}+d[a_{1,i}b_{1,i}]\geq
R_i-T_{i-1}+A_i>R_i-R_{i+1}+A_i\geq R_i-S_i+B_{i-1}$ so it can be
removed. So
$R_i-S_i+B_{i-1}=R_i+S_{i+1}-T_{i-2}-T_{i-1}+d[b_{1,i+1}c_{1,i-3}]$.
Now $R_i+S_{i+1}-T_{i-2}-T_{i-1}+d[a_{1,i+1}c_{1,i-3}]\geq
R_i+R_{i+1}-T_{i-2}-T_{i-1}+d[a_{1,i+1}c_{1,i-3}]\geq
C_{i-1}>R_i-S_i+B_{i-1}$. So
$d[b_{1,i+1}c_{1,i-3}]=d[a_{1,i+1}b_{1,i+1}]\geq A_{i+1}$. Hence
$C_{i-1}>R_i-R_{i+1}+A_i\geq R_i-S_i+B_{i-1}\geq
R_i+S_{i+1}-T_{i-2}-T_{i-1}+A_{i+1}=\min\{
R_i+(R_{i+2}+S_{i+1})/2-T_{i-2}-T_{i-1}+e,
R_i+R_{i+2}-T_{i-2}-T_{i-1}+d[-a_{1,i+2}b_{1,i}],
R_i+R_{i+2}+R_{i+3}-S_i-T_{i-2}-T_{i-1}+d[a_{1,i+3}b_{1,i-1}]\}$.
Since $R_{i+2}+S_{i+1}\geq R_{i+1}+R_{i+2}>T_{i-2}+T_{i-1}$ we
have $R_i+(R_{i+2}+S_{i+1})/2-T_{i-2}-T_{i-1}+e>
R_i-(T_{i-2}+T_{i-1})/2+e\geq C_{i-1}$ so it can be removed. Now
$R_i+R_{i+2}+R_{i+3}-S_i-T_{i-2}-T_{i-1}+d[-a_{i+2,i+3}]\geq
R_i+2R_{i+2}-S_i-T_{i-2}-T_{i-1}+\a_{i+2}\geq
R_i+R_{i+1}+R_{i+2}-S_i-T_{i-2}-T_{i-1}+\a_{i+1}>
R_i+R_{i+1}-T_{i-2}-T_{i-1}+\a_{i+1}\geq C_{i-1}$. So
$R_i+R_{i+2}+R_{i+3}-S_i-T_{i-2}-T_{i-1}+d[a_{1,i+3}b_{1,i-1}]$
can be replaced by
$R_i+R_{i+2}+R_{i+3}-S_i-T_{i-2}-T_{i-1}+d[-a_{1,i+1}b_{1,i-1}]$,
which is $\geq
R_i+R_{i+2}+R_{i+3}-R_{i+1}-T_{i-2}-T_{i-1}+A_i>R_i-R_{i+1}+A_i$.
(We have $R_{i+2}+R_{i+3}\geq R_{i+1}+R_{i+2}>T_{i-2}+T_{i-1}$.)
So this term can also be removed. Now
$R_i+R_{i+2}-T_{i-2}-T_{i-1}+d[-a_{i+1,i+2}]\geq
R_i+R_{i+1}-T_{i-2}-T_{i-1}+\a_{i+1}\geq
C_{i-1}>R_i+S_{i+1}-T_{i-2}-T_{i-1}+A_{i+1}=
R_i+R_{i+2}-T_{i-2}-T_{i-1}+d[-a_{1,i+2}b_{1,i}]$. Hence
$d[-a_{1,i+2}b_{1,i}]=d[a_{1,i}b_{1,i}]\geq A_i$ so
$R_i+S_{i+1}-T_{i-2}-T_{i-1}+A_{i+1}\geq
R_i+R_{i+2}-T_{i-2}-T_{i-1}+A_i>R_i-R_{i+1}+A_i$. Contradiction.
So $R_i-R_{i+1}+A_i\geq C_{i-1}$. \qed

\subsection{Proof of 2.1(ii)}

By Lemma 2.12 the condition that $d[a_{1,i}c_{1,i}]>C_i$ is vacuous if
$i,i+1$ are both nonessential for $M,K$. WLOG we may assume that $i$
is essential. (The case when $i+1$ is essential follows by duality at
index $n-i$.) We apply Lemma 4.2(ii).

If $R_{i+1}+R_{i+2}>S_i+T_{i-1}$ or $i\in\{ 1,n-1\}$ then
$d[a_{1,i}b_{1,i}]\geq A_i\geq C_i$ and $d[b_{1,i}c_{1,i}]\geq B_i\geq
C_i$. So $d[a_{1,i}c_{1,i}]\geq C_i$.

If $R_{i+1}+R_{i+2}\leq S_i+T_{i-1}$ then $T_{i-1}-T_i+B_{i-1}\geq
C_i$. We claim that $R_{i+1}+R_{i+2}\leq T_{i-1}+T_i$. Suppose not so
$T_{i-1}+T_i<R_{i+1}+R_{i+2}\leq S_i+T_{i-1}$. Hence $S_i>T_i$, which
implies $S_i+S_{i+1}\leq T_{i-1}+T_i<R_{i+1}+R_{i+2}$. It follows that
$R_{i+2}>S_i$. Since also $R_{i+1}>T_{i-1}$ ($i$ is essential) we get
$R_{i+1}+R_{i+2}>S_i+T_{i-1}$. Contradiction. So
$R_{i+1}+R_{i+2}\leq T_{i-1}+T_i$. This, together with
$R_{i+1}>T_{i-1}$, implies by Lemma 2.11(ii) that
$d[a_{1,i}c_{1,i}]\geq C_i$ is equivalent to $T_{i-1}-T_i+\c_{i-1}\geq
C_i$. But this follows from $T_{i-1}-T_i+B_{i-1}\geq C_i$ and
$\c_{i-1}\geq B_{i-1}$. \qed

\subsection{Proof of 2.1(iii)} 

Suppose that  $R_{i+1}>T_{i-1}$ and $C_{i-1}+C_i>2e+R_i-T_i$. We have
to prove that $[c_1,\ldots,c_{i-1}]\rep [a_1,\ldots,a_i]$. By Lemma
2.13 $i$ must be an essential index for $M,K$ so
$R_{i+1}+R_{i+2}>T_{i-2}+T_{i-1}$ (if $2<i<n-1$).

\blm One of the following holds:

a) $S_i+A'_i\geq T_i+C_i$

b) $-R_{i+1}+A'_i\geq -S_i+B_{i-1}$

c) $R_{i+1}<S_{i-1}$

Moreover if $R_{i+1}+R_{i+2}\leq S_i+T_{i-1}$ then $S_i+A'_i\geq
T_i+C_i$.
\elm

Take first the case $R_{i+1}+R_{i+2}\leq S_i+T_{i-1}$. If $i>2$ then
also $R_{i+1}+R_{i+2}>T_{i-2}+T_{i-1}$ and so $S_i>T_{i-2}$, which
implies $S_{i-1}\leq T_{i-1}$. (This also happens if $i=2$ since
$S_1\leq T_1$.) By Lemma 4.2(ii) we have $B_{i-1}\geq
T_i-T_{i-1}+C_i$. Suppose
$C_i>S_i-T_i+A'_i=\min\{R_{i+1}-T_i+d[-a_{1,i+1}b_{1,i-1}],
R_{i+1}+R_{i+2}-S_{i-1}-T_i+d[a_{1,i+2}b_{1,i-2}]\}$. But
$R_{i+1}-T_i+d[-a_{1,i+1}c_{1,i-1}]\geq C_i$ and
$R_{i+1}-T_i+d[b_{1,i-1}c_{1,i-1}]\geq R_{i+1}-T_i+B_{i-1}\geq
R_{i+1}-T_{i-1}+C_i>C_i$ so $R_{i+1}-T_i+d[-a_{1,i+1}b_{1,i-1}]\geq
C_i$. Hence
$C_i>R_{i+1}+R_{i+2}-S_{i-1}-T_i+d[a_{1,i+2}b_{1,i-2}]$. But
$R_{i+1}+R_{i+2}-S_{i-1}-T_i+d[a_{1,i+2}c_{1,i-2}]\geq
R_{i+1}+R_{i+2}-T_{i-1}-T_i+d[a_{1,i+2}c_{1,i-2}]\geq C_i$ so
$d[a_{1,i+2}c_{1,i-2}]>d[a_{1,i+2}b_{1,i-2}]$. Hence
$d[a_{1,i+2}b_{1,i-2}]=d[b_{1,i-2}c_{1,i-2}]\geq B_{i-2}$, which
implies $C_i>R_{i+1}+R_{i+2}-S_{i-1}-T_i+B_{i-2}=
\min\{ R_{i+1}+R_{i+2}-(S_{i-1}+T_{i-2})/2-T_i+e,
R_{i+1}+R_{i+2}-T_{i-2}-T_i+d[-b_{1,i-1}c_{1,i-3}],
R_{i+1}+R_{i+2}+S_i-T_{i-3}-T_{i-2}-T_i+d[-b_{1,i}c_{1,i-2}]\}$. (We
have $S_i>T_{i-2}$ so $d[b_{1,i}c_{1,i-4}]$ was replaced by
$d[-b_{1,i}c_{1,i-2}]$ in the formula for $B_{i-2}$. See Lemma
2.7(ii).) But $R_{i+1}+R_{i+2}-(S_{i-1}+T_{i-2})/2-T_i+e\geq
R_{i+1}+R_{i+2}-(T_{i-1}+T_{i-2})/2-T_i+e>(R_{i+1}+R_{i+2})/2-T_i+e\geq
C_i$ so it can be removed. Now
$R_{i+1}+R_{i+2}-T_{i-2}-T_i+d[-c_{i-2,i-1}]=
R_{i+1}+R_{i+2}-T_{i-1}-T_i+\c_{i-2}\geq C_i$ so we can replace
$R_{i+1}+R_{i+2}-T_{i-2}-T_i+d[-b_{1,i-1}c_{1,i-3}]$ by
$R_{i+1}+R_{i+2}-T_{i-2}-T_i+d[b_{1,i-1}c_{1,i-1}]$. But this is
$\geq R_{i+1}+R_{i+2}-T_{i-2}-T_i+B_{i-1}\geq
R_{i+1}+R_{i+2}-T_{i-2}-T_{i-1}+C_i>C_i$ so it can be removed. Thus
$C_i>R_{i+1}+R_{i+2}+S_i-T_{i-3}-T_{i-2}-T_i+d[-b_{1,i}c_{1,i-2}]$.
But $d[b_{1,i}c_{1,i-2}]\geq T_{i-1}-S_i+B_{i-1}\geq
T_i-S_i+C_i$ so
$C_i>R_{i+1}+R_{i+2}-T_{i-3}-T_{i-2}+C_i>C_i$. Contradiction. (We have
$R_{i+1}+R_{i+2}>T_{i-2}+T_{i-1}\geq T_{i-3}+T_{i-2}$.) Thus $C_i\leq
S_i-T_i+A'_i$ and a) holds. 

Suppose now that $R_{i+1}+R_{i+2}>S_i+T_{i-1}$ or $i=n-1$. If $i\neq
n-1$ then $R_{i+1}+R_{i+2}>S_{i-1}+S_i$. (Otherwise $S_{i-1}>T_{i-1}$ so
$S_{i-1}+S_i\leq T_{i-2}+T_{i-1}<R_{i+1}+R_{i+2}$.) So $A'_i=\( A'_i$. 

Assume that neither c) nor a) holds so $R_{i+1}\geq S_{i-1}$ and
$C_i>S_i-T_i+A'_i$. Since $R_{i+1}\geq S_{i-1}$ we have $A'_i=\min\{
R_{i+1}-S_i+d[-a_{1,i+1}b_{1,i-1}],
R_{i+1}+R_{i+2}-S_{i-1}-S_i+d[-a_{1,i}b_{1,i-2}]\}$ by Lemma 2.7(i). If
$A'_i=R_{i+1}-S_i+d[-a_{1,i+1}b_{1,i-1}]$ then
$R_{i+1}-T_i+d[-a_{1,i+1}c_{1,i-1}]\geq
C_i>S_i-T_i+A'_i=R_{i+1}-T_i+d[-a_{1,i+1}b_{1,i-1}]$. Hence
$d[-a_{1,i+1}b_{1,i-1}]=d[b_{1,i-1}c_{1,i-1}]\geq B_{i-1}$ so
$A'_i\geq R_{i+1}-S_i+B_{i-1}$ and we have b). 

So we can suppose that
$A'_i=R_{i+1}+R_{i+2}-S_{i-1}-S_i+d[-a_{1,i}b_{1,i-2}]\geq
R_{i+1}+R_{i+2}-R_i-S_i+A_{i-1}$. If $R_{i+2}\geq S_i$ then, together
with $R_{i+1}\geq S_{i-1}$, this implies $A'_i=\(
A'_i=R_{i+1}-S_i+d[-a_{1,i+1}b_{1,i-1}]$, a case already
discussed. (See Remark 2.8 following Definition 6.) So we can assume
that $R_{i+2}<S_i$. It follows that $i=2$ or
$R_{i+1}+S_i>R_{i+2}+R_{i+1}>T_{i-2}+T_{i-1}$, which by Lemma 4.2(i)
implies $A_{i-1}\geq C_{i-1}$. Therefore $A'_i\geq
R_{i+1}+R_{i+2}-R_i-S_i+C_{i-1}$. We get $C_i>S_i-T_i+A_i'\geq
R_{i+1}+R_{i+2}-R_i-T_i+C_{i-1}$. This implies
$2C_i>R_{i+1}+R_{i+2}-R_i-T_i+C_{i-1}+C_i>
R_{i+1}+R_{i+2}-R_i-T_i+2e+R_i-T_i$ so
$C_i>(R_{i+1}+R_{i+2})/2-T_i+e$. Contradiction. \qed

By duality at index $n+1-i$ we get:

\bco One of the following holds:

a) $-S_i+B'_{i-1}\geq -R_i+C_{i-1}$

b) $T_{i-1}+B'_{i-1}\geq S_i+A_i$

c) $S_{i+1}<T_{i-1}$

Moreover if $R_{i+1}+S_i\leq T_{i-2}+T_{i-1}$ then $-S_i+B'_{i-1}\geq
-R_i+C_{i-1}$.
\eco

\blm (i) If $R_{i+2}>S_i$ or $i=n-1$ and $-R_{i+1}+A_i\geq
-R_i+C_{i-1}$ then $[b_1,\ldots,b_i]\rep [a_1,\ldots,a_{i+1}]$ and
$(a_{1,i}b_{1,i},-a_{1,i+1}c_{1,i-1})_\p =1$.

(ii) If $S_i>T_{i-2}$ or $i=2$ and $T_{i-1}+B_{i-1}\geq T_i+C_i$ then
$[c_1,\ldots,c_{i-2}]\rep [b_1,\ldots,b_{i-1}]$ and
$(b_{1,i-1}c_{1,i-1},-a_{1,i}c_{1,i-2})_\p =1$.
\elm

\pf (i) We have $d(a_{1,i}b_{1,i})+d(-a_{1,i+1}c_{1,i-1})\geq
A_i+T_i-R_{i+1}+C_i\geq R_{i+1}-R_i+C_{i-1}+T_i-R_{i+1}+C_i>2e$ so
$(a_{1,i}b_{1,i},-a_{1,i+1}c_{1,i-1})_\p =1$.

If $A'_{i+1}>A_{i+1}$ then $[b_1,\ldots,b_i]\rep [a_1,\ldots,a_{i+1}]$
by Lemma 2.14 so we are done. So we may suppose that
$A'_{i+1}=A_{i+1}$. 

We have $-R_{i+1}+A_i\geq -R_i+C_{i-1}$ and
$-R_i+C_{i-1}+T_i+C_i>2e$. In order to prove that
$-R_{i+1}+A_i+S_{i+1}+A_{i+1}>2e$, i.e. that
$A_i+A_{i+1}>2e+R_{i+1}-S_{i+1}$ it is enough to show that
$S_{i+1}+A_{i+1}\geq T_i+C_i$. Suppose the contrary, i.e.
$C_i>S_{i+1}-T_i+A_{i+1}=\min\{ R_{i+2}-T_i+d[-a_{1,i+2}b_{1,i}],
R_{i+2}+R_{i+3}-S_i-T_i+d[-a_{1,i+1}b_{1,i-1}]\}$. (We have
$R_{i+2}>S_i$ so in the formula for $A'_{i+1}$ we can replace
$d[a_{1,i+3}b_{1,i-1}]$ by $d[-a_{1,i+1}b_{1,i-1}]$.)

We have $R_{i+2}-T_i+[-a_{i+1,i+2}]\geq R_{i+1}-T_i+\a_{i+1}\geq C_i$
so $R_{i+2}-T_i+d[-a_{1,i+2}b_{1,i}]$ can be replaced by
$R_{i+2}-T_i+d[a_{1,i}b_{1,i}]$ which is $\geq R_{i+2}-T_i+A_i$. Also
$R_{i+2}+R_{i+3}-S_i-T_i+d[a_{1,i+1}b_{1,i-1}]\geq
R_{i+1}+R_{i+2}-S_i-T_i+d[a_{1,i+1}b_{1,i-1}]\geq R_{i+2}-T_i+A_i$.

In conclusion $C_i\geq R_{i+2}-T_i+A_i\geq
R_{i+1}+R_{i+2}-R_i-T_i+C_{i-1}$ so $2C_i\geq
R_{i+1}+R_{i+2}-R_i-T_i+C_{i-1}+C_i>2e+R_{i+1}+R_{i+2}-2T_i$. Hence
$C_i>(R_{i+2}+R_{i+1})/2-T_i+e$. Contradiction.

(ii) follows from (i) by duality at index $n+1-i$.\qed
\vskip 0.5cm

We prove now that $[c_1,\ldots,c_{i-1}]\rep [a_1,\ldots,a_i]$. There are
several cases:

a) $A'_i>A_i$ and $B'_{i-1}>B_{i-1}$. By Lemma 2.14 we have
$[b_1,\ldots,b_{i-1}]\rep [a_1,\ldots,a_i]$ and\\ $[c_1,\ldots,c_{i-1}]\rep
[b_1,\ldots,b_i]$. Also $d(a_{1,i}b_{1,i})+d(b_{1,i-1}c_{1,i-1})\geq
A_i+B_{i-1}=(R_{i+1}-S_i)/2+e+(S_i-T_{i-1})/2+e=
(R_{i+1}-T_{i-1})/2+2e>2e$ so $(a_{1,i}b_{1,i},b_{1,i-1}c_{1,i-1})_\p
=1$. Lemma 1.5(i) implies $[c_1,\ldots,c_{i-1}]\rep [a_1,\ldots,a_i]$.

b) $A_i=A'_i$ and $B'_{i-1}>B_{i-1}$. By Lemma 2.14 we have
$S_{i+1}>T_{i-1}$ and $S_i>T_{i-2}$,\\ $[c_1,\ldots,c_{i-1}]\rep
[b_1,\ldots,b_i]$ and $[c_1,\ldots,c_{i-2}]\rep
[b_1,\ldots,b_{i-1}]$. (If $i=2$ we ignore $S_i>T_{i-2}$ but
$[c_1,\ldots,c_{i-2}]\rep [b_1,\ldots,b_{i-1}]$ still holds
trivially.) In particular, $S_{i-1}\leq T_{i-1}<R_{i+1}$. Also since
$R_{i+1}>T_{i-1}$ we get $R_{i+1}+S_i>T_{i-2}+T_{i-1}$ (or $i=2$),
which implies $A_{i-1},B_{i-1}\geq C_{i-1}$ by Lemma 4.2(i). 

Suppose $-R_{i+1}+A'_i\geq -S_i+B_{i-1}$. Then $A_i=A'_i\geq
R_{i+1}-S_i+(S_i-T_{i-1})/2+e=
(R_{i+1}-T_{i-1})/2+(R_{i+1}-S_i)/2+e>(R_{i+1}-S_i)/2+e$, which is
impossible. We don't have $R_{i+1}<S_{i-1}$ either so $S_i+A_i\geq
T_i+C_i$ by Lemma 4.3. Together with $-R_i+A_{i-1}\geq -R_i+C_{i-1}$
and $-R_i+C_{i-1}+T_i+C_i>2e$, this implies $-R_i+A_{i-1}+S_i+A_i>2e$,
i.e. $A_{i-1}+A_i>2e+R_i-S_i$. Since also $R_{i+1}>S_{i-1}$ we have
$[b_1,\ldots,b_{i-1}]\rep [a_1,\ldots,a_i]$. Also
$[c_1,\ldots,c_{i-1}]\rep [b_1,\ldots,b_i]$ so if we prove that
$(a_{1,i}b_{1,i},b_{1,i-1}c_{1,i-1})_\p =1$ we get
$[c_1,\ldots,c_{i-1}]\rep [a_1,\ldots,a_i]$ by Lemma 1.5(i). Also
$[c_1,\ldots,c_{i-2}]\rep [b_1,\ldots,b_{i-1}]$ so if we prove that
$(b_{1,i-1}c_{1,i-1},-a_{1,i}c_{1,i-2})_\p =1$ we get
$[c_1,\ldots,c_{i-1}]\rep [a_1,\ldots,a_i]$ by Lemma 1.5(iii). Suppose
that none of the 2 conditions is satisfied. It follows that  $2e\geq
d(a_{1,i}b_{1,i})+d(b_{1,i-1}c_{1,i-1})\geq A_i+B_{i-1}\geq
T_i-S_i+C_i+B_{i-1}$ and $2e\geq
d(b_{1,i-1}c_{1,i-1})+d(-a_{1,i}c_{1,i-2})\geq
B_{i-1}+T_{i-1}-R_i+C_{i-1}$. When we add these inequalities we get
$4e\geq 2B_{i-1}+C_{i-1}+C_i+T_{i-1}+T_i-R_i-S_i>
2((S_i-T_{i-1}/2+e)+2e+R_i-T_i+T_{i-1}+T_i-R_i-S_i=4e$. Contradiction.

c) $A'_i>A_i$ and $B_{i-1}=B'_{i-1}$. It follows from b) by
duality at index $n-i+1$.

d) $A_i=A'_i$ and $B_{i-1}=B'_{i-1}$. We make first some remarks:

We cannot have both $R_{i+1}+S_i\leq T_{i-2}+T_{i-1}$ and
$R_{i+1}+R_{i+2}\leq S_i+T_{i-1}$. By adding these would imply
$R_{i+1}+R_{i+2}+R_{i+1}+S_i\leq T_{i-2}+T_{i-1}+S_i+T_{i-1}$,
which is impossible since $R_{i+2}+R_{i+1}>T_{i-2}+T_{i-1}$
and $R_{i+1}>T_{i-1}$.

If $R_{i+1}+S_i>T_{i-2}+T_{i-1}$ or $i=2$ then
$R_{i+1}>S_{i-1}$. Otherwise $S_{i-1}\geq R_{i+1}>T_{i-1}$, which
implies $T_{i-2}+T_{i-1}\geq S_{i-1}+S_i\geq
R_{i+1}+S_i$. Contradiction.

By duality at $n+1-i$ if $R_{i+1}+R_{i+2}>S_i+T_{i-1}$ or $i=n-1$ then
$S_{i+1}>T_{i-1}$.

We will show that one of the following happens and each of them
implies that $[c_1,\ldots,c_{i-1}]\rep [a_1,\ldots,a_i]$:

1. $R_{i+1}+S_i>T_{i-2}+T_{i-1}$ or $i=2$,
$R_{i+1}+R_{i+2}>S_i+T_{i-1}$ or $i=n-1$, $S_i+A_i\geq T_i+C_i$ and
$-S_i+B_{i-1}\geq -R_i+C_{i-1}$.

2. $R_{i+2}>S_i$ or $i=n-1$, $-S_i+B_{i-1}\geq -R_i+C_{i-1}$ and
$-R_{i+1}+A_i\geq -R_i+C_{i-1}$.

3. $S_i>T_{i-2}$ or $i=2$, $S_i+A_i\geq T_i+C_i$ and
$T_{i-1}+B_{i-1}\geq T_i+C_i$.

First we prove that any of these 3 cases implies
$[c_1,\ldots,c_{i-1}]\rep [a_1,\ldots,a_i]$.

In cases 1. and 2. we have $R_{i+1}+R_{i+2}>S_i+T_{i-1}$ or
$i=n-1$. (In the case 2. this follows from $R_{i+2}>S_i$ or $i=n-1$
and $R_{i+1}>T_{i-1}$.) By Lemma 4.2(ii) we get $B_i\geq C_i$. Since
also $-S_i+B_{i-1}\geq -R_i+C_{i-1}$ we get $-S_i+B_{i-1}+T_i+B_i\geq
-R_i+C_{i-1}+T_i+C_i>2e$, i.e. $B_{i-1}+B_i>2e+S_i-T_i$. We also have
$S_{i+1}>T_{i-1}$. (In the case 1. this follows from
$R_{i+1}+R_{i+2}>S_i+T_{i-1}$ or $i=n-1$. In the case 2. we have
$R_{i+2}>S_i$ or $i=n-1$ so $S_{i+1}\geq R_{i+1}>T_{i-1}$.) Therefore
$[c_1,\ldots,c_{i-1}]\rep [b_1,\ldots,b_i]$.

By duality in cases 1. and 3. we have $[b_1,\ldots,b_{i-1}]\rep
[a_1,\ldots,a_i]$.

In case 1. we have $[b_1,\ldots,b_{i-1}]\rep [a_1,\ldots,a_i]$ and
$[c_1,\ldots,c_{i-1}]\rep [b_1,\ldots,b_i]$. Also
$d(a_{1,i}b_{1,i})+d(b_{1,i-1}c_{1,i-1})\geq A_i+B_{i-1}\geq
T_i-S_i+C_i+S_i-R_i+C_{i-1}=C_{i-1}+C_i+T_i-R_i>2e$ so
$(a_{1,i}b_{1,i},b_{1,i-1}c_{1,i-1})_\p =1$. Thus
$[c_1,\ldots,c_{i-1}]\rep [a_1,\ldots,a_i]$ by Lemma 1.5(i).

In case 2. we have $R_{i+2}>S_i$ or $i=n-1$ and $-R_{i+1}+A_i\geq
-R_i+C_{i-1}$, which by Lemma 4.5(i) implies $[b_1,\ldots,b_i]\rep
[a_1,\ldots,a_{i+1}]$ and $(a_{1,i}b_{1,i},-a_{1,i+1}c_{1,i-1})_\p
=1$. Together with $[c_1,\ldots,c_{i-1}]\rep [b_1,\ldots,b_i]$, this
implies $[c_1,\ldots,c_{i-1}]\rep [a_1,\ldots,a_i]$ by Lemma 1.5(ii).

Case 3. follows from 2. by duality at index $n+1-i$.

We prove now that if $R_{i+1}+S_i\leq T_{i-2}+T_{i-1}$ we have case
2. Indeed, condition $R_{i+1}+S_i\leq T_{i-2}+T_{i-1}$ implies
$-R_{i+1}+A_i\geq -R_i+C_{i-1}$ by Lemma 4.2(i) and $-S_i+B_{i-1}\geq
-R_i+C_{i-1}$ by Corollary 4.4. We also have $i=n-1$ or
$R_{i+1}+R_{i+2}>T_{i-2}+T_{i-1}\geq R_{i+1}+S_i$ and so
$R_{i+2}>S_i$. Thus we have 2..

By duality at index $n+1-i$ if $R_{i+1}+R_{i+2}\leq S_i+T_{i-1}$ we
have 3.. 

Suppose now that $R_{i+1}+S_i>T_{i-2}+T_{i-1}$ or $i=2$ and
$R_{i+1}+R_{i+2}>S_i+T_{i-1}$ or $i=n-1$. This implies
$R_{i+1}>S_{i-1}$ and $S_{i+1}>T_{i-1}$. If $S_i+A_i\geq T_i+C_i$ and
$-S_i+B_{i-1}\geq -R_i+C_{i-1}$ we have 1. and we are done.

Suppose now that $S_i+A_i<T_i+C_i$. Since also $R_{i+1}>S_{i-1}$ we
have $-R_{i+1}+A_i\geq -S_i+B_{i-1}$ by Lemma 4.3. If
$-S_i+B_{i-1}<-R_i+C_{i-1}$, since also $S_{i+1}>T_{i-1}$, we have
$T_{i-1}+B_{i-1}\geq S_i+A_i$ by Corollary 4.4 and so $A_i\geq
R_{i+1}-S_i+B_{i-1}\geq R_{i+1}-T_{i-1}+A_i$, which is false. Thus
$-R_{i+1}+A_i\geq -S_i+B_{i-1}\geq -R_i+C_{i-1}$. In order to
have 2. we still need $R_{i+2}>S_i$ or $i=n-1$. Suppose the contrary,
i.e. $R_{i+2}\leq S_i$. We have $C_i>S_i-T_i+A_i\geq
R_{i+1}-T_i+B_{i-1}\geq R_{i+1}+S_i-R_i-T_i+C_{i-1}\geq
R_{i+1}+R_{i+2}-R_i-T_i+C_{i-1}$. It follows that
$2C_i>R_{i+1}+R_{i+2}-R_i-T_i+C_{i-1}+C_i>2e+R_{i+1}+R_{i+2}-2T_i$
so $C_i>(R_{i+1}+R_{i+2})/2-T_i+e$ which is false. Thus we have 2..

By duality if $-S_i+B_{i-1}<-R_i+C_{i-1}$ we have 3. \qed

\subsection{Proof of 2.1(iv)}

Suppose that $T_i\geq R_{i+2}>T_{i-1}+2e\geq
R_{i+1}+2e$. We want to prove that\\ $[c_1,\ldots,c_{i-1}]\rep
[a_1,\ldots,a_{i+1}]$. We may assume that $-a_{1,i+1}c_{1,i-1}\in\fs$.

We have $R_{i+2}-R_{i+1}>2e$ and $T_i-T_{i-1}>2e$ so
$\a_{i+1}>2e$ and $\c_{i-1}>2e$. Since also
$d(-a_{1,i+1}c_{1,i-1})=\j$ we get $d[-a_{1,i+1}c_{1,i-1}]>2e$.

Suppose that $S_i\geq R_{i+2}$. We have $S_i-T_{i-1}\geq
R_{i+2}-T_{i-1}>2e$ so $S_{i-1}\leq T_{i-1}$ (see 2.2), which implies
$R_{i+2}-S_{i-1}>2e$. By Lemma 2.19 we get $[c_1,\ldots,c_{i-1}]\ap
[b_1,\ldots,b_{i-1}]$ and $[b_1,\ldots,b_{i-1}]\rep
[a_1,\ldots,a_{i+1}]$ so $[c_1,\ldots,c_{i-1}]\rep
[a_1,\ldots,a_{i+1}]$. By duality at index $n-i$, if $T_{i-1}\geq
S_{i+1}$ we get $[b_1,\ldots,b_{i+1}]\ap [a_1,\ldots,a_{i+1}]$ and
$[c_1,\ldots,c_{i-1}]\rep [b_1,\ldots,b_{i+1}]$  so again we are done.

Hence we may assume that $R_{i+2}>S_i$ and $S_{i+1}>T_{i-1}$. We will
show that\\ $[b_1,\ldots,b_i]\rep [a_1,\ldots,a_{i+1}]$ and
$[c_1,\ldots,c_{i-1}]\rep [b_1,\ldots,b_i]$. In order to apply
2.1(iii) to $M,N$ and $N,K$ we still need
$A_i+A_{i+1}>2e+R_{i+1}-S_{i+1}$ and $B_{i-1}+B_i>2e+S_i-T_i$. To do
this we apply Corollary 2.17 so it is enough to prove that
$d[-a_{1,i+1}b_{1,i-1}]>2e+S_i-R_{i+2}$ and
$d[-b_{1,i+1}c_{1,i-1}]>2e+T_{i-1}-R_{i+1}$.

Suppose that $d[-a_{1,i+1}b_{1,i-1}]\leq 2e+S_i-R_{i+2}$. Since
$2e+S_i-R_{i+2}<2e<d[-a_{1,i+1}c_{1,i-1}]$ we have
$d[-a_{1,i+1}b_{1,i-1}]=d[b_{1,i-1}c_{1,i-1}]\geq B_{i-1}=\min\{
(S_i-T_{i-1})/2+e, S_i-T_{i-1}+d[-b_{1,i}c_{1,i-2}],
S_i+S_{i+1}-T_{i-2}-T_{i-1}+d[-b_{1,i+1}c_{1,i-1}]\}$. (We have
$S_{i+1}>T_{i-1}$ so we can replace $d[b_{1,i+1}c_{1,i-3}]$ by
$d[-b_{1,i+1}c_{1,i-1}]$.) But $(S_i-T_{i-1})/2+e>2e+S_i-R_{i+2}$
since $2R_{i+2}>2e+S_i+T_{i-1}$. (We have $R_{i+2}>S_i$ and
$R_{i+2}>T_{i-1}+2e$.) Hence it can be removed. Also
$S_i-T_{i-1}+d[-b_{1,i}c_{1,i-2}]\geq S_i-T_{i-1}>2e+S_i-R_{i+2}$ so
it can be removed. Hence $2e+S_i-R_{i+2}\geq
d[-a_{1,i+1}b_{1,i-1}]\geq
B_{i-1}=S_i+S_{i+1}-T_{i-2}-T_{i-1}+d[-b_{1,i+1}c_{1,i-1}]$. It
follows that $d[-b_{1,i+1}c_{1,i-1}]\leq
2e+T_{i-2}+T_{i-1}-R_{i+2}-S_{i+1}<2e+T_{i-1}-S_{i+1}$. (We have
$R_{i+2}>T_{i-1}+ 2e\geq T_{i-2}$.) By duality at index $n-i$ the
inequality $d[-b_{1,i+1}c_{1,i-1}]\leq 2e+T_{i-1}-S_{i+1}$ implies
$d[-b_{1,i+1}c_{1,i-1}]\geq A_{i+1}=
R_{i+2}+R_{i+3}-S_i-S_{i+1}+d[-a_{1,i+1}b_{1,i-1}]$. Together with
$d[-a_{1,i+1}b_{1,i-1}]\geq
B_{i-1}=S_i+S_{i+1}-T_{i-2}-T_{i-1}+d[-b_{1,i+1}c_{1,i-1}]$, this
implies that $R_{i+2}+R_{i+3}\leq T_{i-2}+T_{i-1}$. Since
$R_{i+2}+R_{i+3}\geq 2R_{i+2}-2e$ and $T_{i-2}+T_{i-1}\leq
2T_{i-1}+2e$ we get $2R_{i+2}-2e\leq 2T_{i-1}+2e$ so $R_{i+2}\leq
T_{i-1}+2e$. Contradiction.

So $d[-a_{1,i+1}b_{1,i-1}]>2e+S_i-R_{i+2}$. By duality at index $n-i$
$d[-b_{1,i+1}c_{1,i-1}]>2e+T_{i-1}-S_{i+1}$. Thus
$[b_1,\ldots,b_i]\rep [a_1,\ldots,a_{i+1}]$ and
$[c_1,\ldots,c_{i-1}]\rep [b_1,\ldots,b_i]$, which
implies $[c_1,\ldots,c_{i-1}]\rep [a_1,\ldots,a_{i+1}]$. \qed

\section{The case $[M:N]\spq\p$}

In this section $M\sbq N$ are two lattices on the same quadratic space
and $[M:N]\spq\p$. We want to prove that $N\leq M$.

First take case $[M:N]=\oo$, i.e. when $M=N$. By Corollary 3.11 the
conditions (i)-(iv) are independent of the choice of BONGs we can
assume that $b_i=a_i$. We also have $S_i=R_i$ and $\b_i=\a_i$. Since
$R_i=S_i$ for all $i$ condition 2.1(i) is satisfied. We have
$d[a_{1,i}b_{1,i}]=\min\{ d(a_{1,i}b_{1,i}),\a_i,\b_i\} =\min\{
d(a_{1,i}a_{1,i}),\a_i,\a_i\} =\a_i$ and
$d[-a_{1,i+1}b_{1,i-1}]=\min\{
d(-a_{1,i+1}b_{1,i-1}),\a_{i+1},\b_{i-1}\} =\min\{
d(-a_{1,i+1}a_{1,i-1}),\a_{i+1},\a_{i-1}\} =d[-a_{i,i+1}]$. Hence
$A_i\leq\min\{ (R_{i+1}-S_i)/2+e,R_{i+1}-S_i+d[-a_{1,i+1}b_{1,i-1}]\}
=\min\{ (R_{i+1}-R_i)/2+e,R_{i+1}-R_i+d[a_{i,i+1}]\}
=\a_i=d[a_{1,i}b_{1,i}]$ and we have 2.1(ii). Finally, since $b_i=a_i$
for all $i$ the relations $[b_1,\ldots,b_{i-1}]\rep [a_1,\ldots,a_i]$
and $[b_1,\ldots,b_{i-1}]\rep [a_1,\ldots,a_{i+1}]$ hold
unconditionally so we have 2.1(iii) and (iv). \qed

For the case $[M:N]=\p$ we need the following lemma:

\blm If $[M:N]=\p$ then there are some splittings $N=J\pp K$ and
$M=J'\pp K$ s.t. $J\sb J'$ are either unary or binary modular and
$[J':J]=\p$. If $J,J'$ are unary then $J=\p J'$ and if they are
binary then $\p J'\sb J\sb J'$ and  $\ss J'=\p\1\ss J$. \elm
\pf Since $[M:N]=\p$ there is a primitive element $x\in N$ s.t.
$M=N+\p\1 x$. Moreover $M=N+\p\1 x$ remains true if we replace $x$
by any $x'\ev x\pmod{\p N}$. Let $N=N_1\pp\cdots\pp N_t$ be a Jordan
splitting with $\ss N_1\sp\ldots\sp\ss N_t$ and we write
$x=x_1+\cdots+x_t$ with $x_k\in N_k$. Since $x\in N$ is primitive at
least one of the $x_k$'s is primitive in $N_k$. Let $k_0$ be the
smallest $k$ s.t. $x_k\in N_k$ is primitive. It follows that $x\ev
x_{k_0}+\cdots+x_n \pmod{\p N}$ so we can replace $x$ by
$x_{k_0}+\cdots+x_n$.

We write $N=N'\pp N''$ with $N'=N_1\pp\cdots\pp N_{k_0-1}$ and
$N''=N_{k_0}\pp\cdots\pp N_t$. Then $x=x_{k_0}+\cdots+x_t\in N''$ and
since $x_{k_0}\in N_{k_0}$ is primitive we have $B(x,N'')=\ss
N_{k_0}=\ss N''$. Let $\ord\ss N''=r$. If $\ord Q(x)=r$ then we have a
splitting $N''=\oo x\pp N^*$ so we can take $J=\oo x$, $J'=\p\1 x$ and
$K=N'\pp N^*$. If $\ord Q(x)>\ord\ss N''=r$ then there is $y\in N''$
s.t. $B(x,y)=\pi^r$. Now $J=\oo x+\oo y$ is
$\p^r$-modular so we have a splitting $N''=J\pp N^*$. If we put
$K=N'\pp N^*$ we have $N=J\pp K$ and $M=J'\pp K$ where $J'=\p\1 x+\oo
y$. Since also $J=\oo x+\oo y$, we get $\p J'\sb J\sb J'$. Now
$B(\pi\1 x,y)=\pi^{r-1}$, $\ord Q(y)\geq r$ and $\ord Q(\pi\1
x)=\ord Q(x)-2\geq r-1$ so $J'$ is $\p^{r-1}$-modular and we are
done. \qed

Let $K=K_1\pp\cdots\pp K_t$ be a Jordan decomposition and let $\ord\ss
K_k=r_k$, $\ord\nn K^{\ss K_k}=u_k$ and $\dim K_{(k)}=n_k$. Let $J,J'$
be as above. Denote $a=\dim J=\dim J'$, $r=\ord\ss J$, $r'=\ord\ss
J'$. Let $k_1=\max\{ k\mid r_k\leq r'\}$ and $k_2=\max\{ k\mid
r_k<r\}$. Hence $r_{k_1}\leq r'<r\leq r_{k_2}$, $r'<r_{k_1+1}$ and
$r_{k_2-1}<r$. 

We have the splittings $M=K_1\pp\cdots\pp K_{k_1}\pp J'\pp
K_{k_1+1}\pp\cdots\pp K_t$ and $N=K_1\pp\cdots\pp K_{k_2}\pp J\pp
K_{k_2+1}\pp\cdots\pp K_t$, which are almost Jordan decompositions in
the sense that if $r_{k_1}=r'$ then $K_{k_1}\pp J'$ is the
$\p^{r_{k_1}}$ modular component of $M$ and if $r_{k_2+1}=r$ then
$J\pp K_{k_2+1}$ is the $\p^{r_{k_2+1}}$ modular component of $N$. (If
$r_1>r'$ we just take $k_1=0$ so we have the splitting $M=J'\pp
K_1\pp\cdots\pp K_t$. Similarly for $N$ when $r_1\geq r$.)

Note that if $a=2$ then $r=r'+1$ so $k_1=k_2$. If $a=1$ then $r=r'+2$
so we have either $r_{k_1+1}=r'+1$ and $k_2=k_1+1$ or $r_{k_1+1}\neq
r'+1$ and $k_1=k_2$.

\bff By the use of duality the number of cases to be considered can be
cut by half. Recall that by 2.5 conditions 2.1(i) or (iii) for
$N^\*,M^\* $ at $i$ are equivalent to the same conditions for $M,N$ at
$n+1-i$, while 2.1(ii) and (iv) are equivalent to the same conditions
for $M,N$ at $n-i$. 

Together with $[M:N]=\p$ we have $[N^\* :M^\* ]=\p$. Similar to
$M=K_1\pp\cdots K_{k_1}\pp J'\pp K_{k_1+1}\pp\cdots\pp K_t$ and
$N=K_1\pp\cdots K_{k_2}\pp J\pp K_{k_2+1}\pp\cdots\pp K_t$ we have the
splittings $N^\* =K_t^\*\pp\cdots\pp K_{k_2+1}^\*\pp J^\*\pp
K_{k_2}^\*\pp\cdots\pp K_1^\*$ and $M^\* =K_t^\*\pp\cdots\pp
K_{k_1+1}^\*\pp J^\*\pp K_{k_1}^\*\pp\cdots\pp K_1^\*$. We have $\dim
(K_t^\*\pp\cdots\pp K_{k_2+1}^\* )=n-a-n_{k_2}$ and $\dim
(K_t^\*\pp\cdots\pp K_{k_1+1}^\* )=n-a-n_{k_1}$. So the analogues of
$n_{k_1},n_{k_2}$ and $a$ for $N^\*,M^\*$ are $n-a-n_{k_2},n-a-n_{k_1}$
and $a$. 

If we prove 2.1(i) and (iii) for $i\leq n_{k_2}+1$ then similarly we
can prove the same statements for $N^\*,M^\*$ for $i\leq
n-a-n_{k_1}+1$, which are equivalent to the same statements for $M,N$
at $i\geq n_{k_1}+a$. Since $n_{k_1}+a\leq n_{k_2}+2$ 2.1(i) and (iii)
will hold at all $i$. 

Similarly if we prove 2.1(ii) and (iv) at $i\leq n_{k_2}+a-1$ then we
can prove them for $N^\*,M^\*$ at $i\leq n-a-n_{k_1}+a-1=n-n_{k_1}-1$,
which is equivalent to proving them for $M,N$ at $i\geq
n_{k_1}+1$. Since $n_{k_1}+1\leq n_{k_2}+a$ 2.1(ii) and (iv) will hold
at all $i$. 

In conclusion we may restrict ourselves to proving 2.1(i) and (iii) at
$i\leq n_{k_2}+1$ and 2.1(ii) and (iv) at $i\leq n_{k_2}+a-1$.
\eff

\subsection{Proof of 2.1(i)}

We have $u_k=\ord\nn K^{\p^{r_k}}$ and we denote $u=\ord\nn J$ and
$u'=\ord\nn J'$. If $a=1$ then $J$ and $J'$ are unary so $u=r$ and
$u'=r'$. Since $J=\p J'$ we have $r=r'+2$. If $a=2$ then
$r=r'+1$. Since $\p J'\sb J\sb J'$ we have $\p^2\nn J'\sbq\nn J\sbq\nn
J'$ and so $u'+2\geq u\geq u'$.

The scales of Jordan components of $M$ have orders $r_1,\ldots,r_t$
and $r'$ and for $N$ they have orders $r_1,\ldots,r_t$ and $r$. We
denote $v'=\ord\nn M^{\p^{r'}}$, $v'_k=\ord\nn M^{\p^{r_k}}$,
$v=\ord\nn N^{\p^{r}}$, $v_k=\ord\nn N^{\p^{r_k}}$.

\bff For any $1\leq k\leq t$ we have $v_k=\ord\nn (K\pp
J)^{\p^{r_k}}=\min\{\ord\nn K^{\p^{r_k}},\ord\nn J^{\p^{r_k}}\}$. But
$J^{\p^{r_k}}=J$ if $k\leq k_2$ and $J^{\p^{r_k}}=\p^{r_k-r}J$ if
$k>k_2$. Thus $\ord\nn J^{\p^{r_k}}=u$ if $k\leq k_2$ and $\ord\nn
J^{\p^{r_k}}=u+2(r_k-r)$ if $k>k_2$. Also $\ord\nn
K^{\p^{r_k}}=u_k$. Therefore $v_k=\min\{ u_k,u\}$ if $k\leq k_2$ and
$v_k=\min\{ u_k,u+2(r_k-r)\}$ if $k>k_2$.

Also $r_{k_2}\leq r<r_{k_2+1}$ so $\ord\nn K^{\p^r}=\min\{
u_{k_2}+2(r-r_{k_2}),u_{k_2+1}\}$ by Lemma 3.5(i). (If $r\leq r_1$,
when $k_2=0$, we ignore $u_{k_2}+2(r-r_{k_2})$; if $r>r_t$ when
$k_2=t$, we ignore $u_{k_2+1}$.) Since $N^{\p^r}=K^{\p^r}\pp J^{\p^r}$
and $\ord\nn J^{\p^r}=u$ we get $v=\ord\nn N^{\p^r}=\min\{
u_{k_2}+2(r-r_{k_2}),u_{k_2+1},u\}$.

By considering the norms of $N^{\p^{r_k}}$'s and $N^{\p^r}$ we get
$v_1\leq\ldots\leq v_{k_2}\leq v\leq v_{k_2+1}\leq\ldots\leq v_t$.

Similarly for $M$, with $k_2,v_k,r,v$ replaced by $k_1,v'_k,r',v'$,
respectively. 
\eff

\bff We now use [B3, Lemma 2.13] to write $R_i$'s and $S_i$'s in
terms of $r_k,r,r',v,v',v_k,v'_k$.

If $k\leq k_2$ then for any $n_{k-1}<i\leq n_k$ we have $S_i=v_k$ if
$i\ev n_{k-1}+1\m2$ and $S_i=2r_k-v_k$ if $i\ev n_{k-1}\m2$. This
corresponds to the $K_1\pp\cdots\pp K_{k_2}$ part of the splitting
$N=K_1\pp\cdots\pp K_{k_2}\pp J\pp K_{k_2+1}\pp\cdots\pp
K_t$. Corresponding to $J$ we have $S_{k_2+1}=r=u$ if $a=1$ and
$S_{k_2+1}=v$, $S_{k_2+2}=2r-v$ if $a=2$. Finally, corresponding to the
$K_{k_2+1}\pp\cdots\pp K_t$ part if $k\geq k_2+1$ and $n_{k-1}+a<i\leq
n_k+a$ we have $S_i=v_k$ if $i\ev n_{k-1}+a+1\m2$ and $S_i=2r_k-v_k$
if $i\ev n_{k-1}+a\m2$. (Note that these relations hold also when
$r_{k_2+1}=r$, i.e. when the Jordan splitting of $N$ is
$K_1\pp\cdots\pp K_{k_2}\pp (J\pp K_{k_2+1})\pp K_{k_2+2}\pp\cdots\pp
K_t$.) 

Suppose that $v_k=u$ for some $k\leq k_2$. Since $v_k=\min\{ u_k,u\}$ we have
$u\leq u_k$. Then for any $k\leq l\leq k_2$ we have $u\leq u_k\leq u_l$ so
$v_l=\min\{ u_l,u\} =u$. Also $v=\min\{ u_{k_2}+2(r-r_{k_2}),u_{k_2+1},u\}$
and $u\leq u_k\leq u_{k_2}\leq u_{k_2}+2(r-r_{k_2})$
and $u\leq u_k\leq u_{k_2+1}$ so $v=u$.

Similarly for the lattice $M$ but with $S_i,k_2,v_k,r,v$ replaced by
$R_i,k_1,v'_k,r'$ and $v'$.
\eff

We start now our proof. We will show that in fact $R_i\leq S_i$ or
$R_i+R_{i+1}=S_{i-1}+S_i$ for all $i$. By 5.2 we may restrict
ourselves to the case $i\leq n_{k_2}+1$.

Suppose first that $i\leq n_{k_1}$. Then $n_{k-1}<i\leq n_k$ for some
$k\leq k_1$. We have $v_k=\min\{ u_k,u\}\geq\min\{ u_k,u'\} =v'_k$
with equality unless $u'<u$ and $u'<u_k$. If $v_k>v'_k$ then
$v'_k=u'<u_k$. This implies that $v'_l=u'$ for any $k\leq l\leq k_1$
and $v'=u'$. (See 5.4.) 

Suppose that $S_i<R_i$. If $i\ev n_{k-1}+1\m2$ then $S_i=v_k\geq
v'_k=R_i$. So we may assume $i\ev n_{k-1}\m2$. This implies that
$n_{k-1}+2\leq i\leq n_k$ so $S_{i-1}+S_i=2r_k$. If $i<n_k$ then
$R_i+R_{i+1}=2r_k=S_{i-1}+S_i$ so we are done. So we can assume that
$i=n_k$. We have $2r_k-v_k=S_i<R_i=2r_k-v'_k$ so $v'_k<v_k$. But this
implies that $v'_l=u'$ for any $k\leq l\leq k_1$ and $v'=u'$. If
$k<k_1$ then $R_{i+1}=v'_{k_1+1}$ while if $k=k_1$ then
$R_{i+1}=v'$. In both cases $R_{i+1}=u'$. Together with
$R_i=2r_k-v'_k=2r_k-u'$, this implies
$R_i+R_{i+1}=2r_k=S_{i-1}+S_i$. So we proved (i) for $i\leq n_{k_1}$.

If $k_1=k_2$ then we proved (i) for $i\leq n_{k_1}=n_{k_2}$ so we are
left to prove it for $i=n_{k_1}+1=n_{k_2}+1$. But at this index we
have $R_i=r'<r=S_i$ if $a=1$ and $R_i=v'=\min\{
u_{k_2}+2(r'-r_{k_2}),u_{k_2+1},u'\}\leq\min\{
u_{k_2}+2(r-r_{k_2}),u_{k_2+1},u\} =v=S_i$ if $a=2$ so we are done. 

Suppose now that $k_1\neq k_2$ so we have $a=1$, $k_2=k_1+1$ and
$r_{k_1+1}=r'+1=r-1$. Also $u=r$ and $u'=r'=r-2$. We have to prove (i)
for $n_{k_1}+1\leq i\leq n_{k_2}+1$. We have
$v_{k_1+1}=v_{{k_2}}=\min\{ u_{k_2},u\} =\min\{ u_{k_1+1},r\}$ and
$v'_{k_1+1}=\min\{ u_{k_1+1},u'+2(r_{k_1+1}-r')\} =\min\{
u_{k_1+1},r\}$. (We have $u'+2(r_{k_1+1}-r)=r-2+2(r'+1-r')=r$.) Thus
$v'_{k_1+1}=v_{k_1+1}=r-1$ if $u_{k_1+1}=r_{k_1+1}=r-1$ and
$v'_{k_1+1}=v_{k_1+1}=r$ if $u_{k_1+1}>r_{k_1+1}=r-1$.

Suppose that $v'_{k_1+1}=v_{k_1+1}=r-1=r_{k_1+1}$. Then the
$\p^{r-1}=\p^{r_{k_1+1}}$ modular component of both $M$ and $N$ is
proper. Thus $R_i=r_{k_1+1}=r-1$ for $n_{k_1}+2=n_{k_1}+a+1<i\leq
n_{k_1+1}+a=n_{k_2}+1$ and $S_i=r_{k_1+1}=r-1$ for
$n_{k_1}+1=n_{k_2-1}+1\leq i\leq n_{k_2}$. Also $S_{n_{k_2}+1}=v=r$
and $R_{n_{k_1}+1}=v'=r'=r-2$. In conclusions the
sequences $R_{n_{k_1}+1},\ldots,R_{n_{k_2}+1}$ and
$S_{n_{k_1}+1},\ldots,S_{n_{k_2}+1}$ are $r-2,r-1,\ldots,r-1$ and
$r-1,\ldots,r-1,r$, respectively. Therefore $R_i\leq S_i$ for
$n_{k_1}+1\leq i\leq n_{k_2}+1$ and we are done.

Suppose now that $v'_{k_1+1}=v_{k_1+1}=r>r-1=r_{k_1+1}$. The sequence
$R_{n_{k_1}+2}=R_{n_{k_1}+a+1},\ldots,R_{n_{k_1+1}+a}=R_{n_{k_2}+1}$
is made of copies the pair $v'_{k_1+1},2r_{k_1+1}-v'_{k_1+1}$,
i.e. $r,r-2$. Similarly the sequence
$S_{n_{k_1}+1},\ldots,S_{n_{k_1+1}}=S_{n_{k_2}}$ is made of copies of
the pair $v_{k_1+1},2r_{k_1+1}-v_{k_1+1}$, i.e. $r,r-2$. Also
$R_{n_{k_1}+1}=v'=r'=r-2$ and $S_{n_{k_2}+1}=v=r$. In conclusion the
sequences $R_{n_{k_1}+1},\ldots,R_{n_{k_2}+1}$ and
$S_{n_{k_1}+1},\ldots,S_{n_{k_2}+1}$ are $r-2,r,r-2,\ldots,r,r-2$
resp. $r,r-2,r,\ldots,r-2,r$. If $n_{k_1}+1<i<n_{k_2}+1$ then
$R_i+R_{i+1}=2r-2=S_{i-1}+S_i$, while if $i\in\{
n_{k_1}+1,n_{k_2}+1\}$ then $R_i=r-2<r=S_i$ so we are done. \qed

\subsection{Consequences of 2.1(i)}

Since property (i) is true when $[M:N]=\p$ and it is transitive it is
true for $M,N$ arbitrary with $m=n$. As seen from the Lemma 2.20 it is
also true when $n<m$. 

Before proving 2.1(ii) we give some consequences of 2.1(i). These will
shorten the proof of 2.1(ii) and will be useful in the future
sections. First we prove some properties of $(\bbb,\leq )$.

\blm Let $x,y\in\bbb$ with $x=(x_1,\ldots,x_m)$ and $y=(y_1,\ldots,y_n)$
s.t. $x\leq y$. We have:

(i) $x_1+\cdots +x_i\leq y_1+\cdots +y_i$ for any $1\leq i\leq n$

If $m=n$ then:

(ii) $x_i+\cdots +x_n\leq y_i+\cdots +y_n$ for any $1\leq i\leq n$.

(iii) If $x_1+\cdots +x_n=y_1+\cdots +y_n$ then $x=y$.

(iv) If $x_1+\cdots +x_n+k=y_1+\cdots +y_n$ for some $k\geq 0$ then
for any $1\leq i\leq j\leq n$ we have $x_i+\cdots +x_j+k\geq y_i+\cdots +y_j$
with equality iff  $x_1+\cdots +x_{i-1}=y_1+\cdots +y_{i-1}$ and
$x_{j+1}+\cdots +x_n=y_{j+1}+\cdots +y_n$.
\elm
\pf (i) If $i$ is even we add the inequalities $x_1+x_2\leq y_1+y_2$,
$x_3+x_4\leq y_3+y_4,\ldots,~x_{i-1}+x_i\leq y_{i-1}+y_i$. If $i$ is
odd we add $x_1\leq y_1,~x_2+x_3\leq y_2+y_3,\ldots,~x_{i-1}+x_i\leq
y_{i-1}+y_i$. 

(ii) follows from (i) by duality. We have $x^\*=(-x_n,\ldots,-x_1)$,
$y^\* =(-y_n,\ldots,-y_1)$ and $y^\*\leq x^\*$. By (i) we get 
$-y_n-\cdots -y_i\leq -x_n-\cdots -x_i$ so $x_i+\cdots +x_n\leq
y_i+\cdots +y_n$. 

(iii) For any $1\leq j\leq n$ we use (i) and (ii) and get
$x_1+\cdots +x_j\leq y_1+\cdots +y_j$ and $x_{j+1}+\cdots +x_n\leq
y_{j+1}+\cdots +y_n$. Since $x_1+\cdots +x_n=y_1+\cdots +y_n$ we must have
equalities. In particular, $x_1+\cdots +x_j=y_1+\cdots +y_j$ for all
$j$. This implies $x_i=y_i$ for all $i$ so $x=y$.

(iv) follows from $x_1+\cdots +x_{i-1}\leq y_1+\cdots +y_{i-1}$,
$x_{j+1}+\cdots +x_n\leq y_{j+1}+\cdots +y_n$ and
$(x_1+\cdots +x_{i-1})+(x_i+\cdots +x_j+k)+(x_{j+1}+\cdots +x_n)
=(y_1+\cdots +y_{i-1})+(y_i+\cdots +y_j)+(y_{j+1}+\cdots +y_n)$. \qed

\blm Let $x,y\in\bbb$ with $x=(x_1,\ldots,x_m)$ and
$y=(y_1,\ldots,y_n)$ s.t. $x\leq y$ and $x\neq y$. We have:

(i) If $x_1+\cdots +x_a=y_1+\cdots +y_a$ for some $1\leq a\leq n$ then
$x_i+x_{i+1}=y_i+y_{i+1}$ for any $1\leq i<a$ with $i\ev a+1\m2$, and
$x_i=\min\{ y_i,x_{a+1}\}$ for any $1\leq i\leq\min\{ a+1,n\}$ with
$i\ev a+1\m2$. Moreover if there is $1\leq c\leq a$ s.t. $x_c\neq y_c$
and $c$ is minimal with this property then $c\ev a+1\m2$, $x_c<y_c$
and $x_c=x_{c+2}=\ldots =x_{a+1}$. 

(ii) If $m=n$ and $x_b+\cdots +x_n=y_b+\cdots +y_n$ for some $1\leq
b\leq n$ then $x_{i-1}+x_i=y_{i-1}+y_i$ for any $b<i\leq n$ with $i\ev
b-1\m2$ and $y_i=\max\{ x_i,y_{b-1}\}$ for any $b-1\leq i\leq n$ with
$i\ev b-1\m2$. Moreover if there is $b\leq d\leq n$ s.t. $x_d\neq y_d$
and $d$ is maximal with this property then $d\ev b-1\m2$, $x_d<y_d$
and $y_d=y_{d-2}=\ldots =y_{b-1}$. 
\elm
\pf (i) If $1\leq i<a$, $i\ev a+1\m2$, then $x_1+\cdots +x_{i-1}\leq
y_1+\cdots +y_{i-1}$, $x_i+x_{i+1}\leq y_i+y_{i+1}$,
$x_{i+2}+x_{i+3}\leq y_{i+2}+y_{i+3}$,\ldots, $x_{a-1}+x_a\leq
y_{a-1}+y_a$. By adding we get $x_1+\cdots +x_a\leq y_1+\cdots
+y_a$. But $x_1+\cdots +x_a=y_1+\cdots +y_a$ so we must have
equalities. In particular, $x_i+x_{i+1}=y_i+y_{i+1}$.

Since $x\leq y$ we have $m\geq n\geq a$. If $m=n=a$ then
$x_1+\cdots +x_n=y_1+\cdots +y_n$ so $x=y$ by Lemma 5.5(iii). Hence
$m>a$. We note that if $n\geq a+1$ then $x_1+\cdots +x_{a+1}\leq
y_1+\cdots +y_{a+1}$ by Lemma 5.5(i), which, together with
$x_1+\cdots +x_a=y_1+\cdots +y_a$, implies $x_{a+1}\leq y_{a+1}$. Thus
$x_i=\min\{ y_i,x_{a+1}\}$ holds at $i=a+1$.

If $x_j=y_j$ for all $1\leq j\leq a$ then for any $1\leq i<a$
with $i\ev a+1\m2$ we have $x_i+x_{i+1}=y_i+y_{i+1}$ and $\min\{
y_i,x_{a+1}\} =\min\{ x_i,x_{a+1}\} =x_i$ so we are done. Otherwise
let $1\leq c\leq a$ be minimal s.t. $x_c\neq y_c$.

By Lemma 5.5(i) $x_1+\cdots +x_c\leq y_1+\cdots +y_c$ and by the
minimality of $c$ we have $x_i=y_i$ for $1\leq i\leq c-1$ so $x_c\leq
y_c$. But $x_c\neq y_c$ so $x_c<y_c$. Suppose now that $c\ev
a\m2$. Then $x_1+\cdots +x_{c-1}=y_1+\cdots +y_{c-1}$, $x_c<y_c$,
$x_{c+1}+x_{c+2}\leq y_{c+1}+y_{c+2}$, $x_{c+3}+x_{c+4}\leq
y_{c+3}+y_{c+4}$, \ldots, $x_{a-1}+x_a\leq y_{a-1}+y_a$. By adding we
get $x_1+\cdots +x_a<y_1+\cdots +y_a$. Contradiction. So $c\ev
a+1\m2$. 

We prove now by induction that for any $c\leq i\leq a+1$ with $i\ev
a+1\m2$ we have $x_i=x_c$. At $i=c$ this statement is trivial. Suppose
now that $i>c$. By the induction hypothesis we have $x_{i-2}=x_c$. Since
$i-2\ev a+1\ev c\m2$ and $i-2\geq c$ we have $x_{i-2}=x_c<y_c\leq y_{i-2}$. As
proved above we have $x_{i-2}+x_{i-1}=y_{i-2}+y_{i-1}$ and so
$x_{i-1}>y_{i-1}$. This implies $x_{i-1}+x_i\leq
y_{i-2}+y_{i-1}=x_{i-2}+x_{i-1}$ i.e. $x_i\leq x_{i-2}$ which implies
$x_i=x_{i-2}=x_c$. 

We prove now that $x_i=\min\{ y_i,x_{a+1}\}$ for any
$1\leq i<a$ with $i\ev a+1\m2$. If $i<c$ then $x_i=y_i$ and, since
$i<a+1$ and $i\ev a+1\m2$, we have $x_i\leq x_{a+1}$. Thus $x_i=\min\{
x_i,x_{a+1}\} =\min\{ y_i,x_{a+1}\}$. If $i\geq c$ then
$x_i=x_{a+1}$ and also $x_i=x_c<y_c\leq y_i$ ($i\ev a+1\ev c\m2$ and
$c\leq i$). Hence $x_i=x_{a+1}=\min\{ y_i,x_{a+1}\}$.

(ii) follows from (i) by duality. If $y^\* =(y_1^\*,\ldots,y_n^\* )$
and $x^\* =(x_1^\*,\ldots,x_n^\* )$ then $y_i^\* =-y_{n-i+1}$,
$x_i^\* =-x_{n-i+1}$ and $y^\*\leq x^\*$. It is easy to see that
(ii) is equivalent to (i) with $x,y,a,c$ replaced by
$y^\*,x^\*,n-b+1,n-d+1$. \qed

\blm Suppose $N$ is given and let $k\geq 0$ be maximal with the
property that $S_1=S_3=\ldots=S_{2k+1}$. Let $y\in N$ be a norm
generator and let $pr_{y^\pp}N\ap\[ a_2,\ldots,a_n\]$ relative to a
good BONG $x_2,\ldots,x_n$. If $R_i=\ord a_i$ then:

(i) $R_i=S_i$ for any $i>2k+1$

(ii) If $1<i<2k+1$ is even then $R_i+R_{i+1}=S_i+S_{i+1}=S_i+S_1$,
$R_i\geq S_i$ and if $R_i=S_i$ then $R_j=S_j$ for $j\geq i$.
\elm
\pf We have $\ord Q(y)=\ord Q(y_1)=S_1$. Let $s\gg 0$ s.t. $S_1-2s\leq
R_3$ and let $M=\p^{-s}y+N$. Let $x_1=\pi^{-s}y$ and let $a_1=\ord
Q(x_1)$. We have $\ord a_1=R_1$, where $R_1=S_1-2s$. We claim that
$M\ap\[ a_1,\ldots,a_n\]$ relative to the good BONG
$x_1,\ldots,x_n$. Since $M\sbq\p^{-s}N$ we have $\nn
M\sbq\nn\p^{-s}N=\nn \p^{-s}N=\p^{S_1-2s}=Q(x_1)\oo$. Hence $x_1\in M$
is a norm generator. We have
$pr_{x_1^\pp}M=pr_{y^\pp}(\p^{-s}y+N)=pr_{y^\pp}N=\[ x_2,\ldots,x_n\]$
so $M=\[ x_1,\ldots,x_n\]$. The BONG  $x_2,\ldots,x_n$ is good and
$R_1=S_1-2s\leq R_3$ so the BONG $x_1,\ldots,x_n$ is also good. We
have $N\sbq M$ so they satisfy 2.1(i). Hence $\rr (M)\leq \rr(N)$.

Now $y\in N$ is a norm generator so it is primitive and
$M=\p^{-s}y+N$, which implies $[M:N]=\p^s$. Hence $S_1+S_2+\cdots
+S_n=\ord\det N=2s+\ord\det M=2s+R_1+R_2+\cdots+R_n=S_1+R_2+\cdots
+R_n$ so $S_2+\cdots+S_n=R_2+\cdots+R_n$. Therefore we can apply Lemma
5.6(ii) with $x=\rr (M)$, $y=\rr (N)$ and $b=2$. If $1<i<2k+1$ is
even then $2=b<i+1\leq 2k+1\leq n$ and $i+1\ev 1=b-1\m2$ so
$R_i+R_{i+1}=S_i+S_{i+1}$ and $S_{i+1}=\max\{ R_{i+1},S_1\}\geq
R_{i+1}$ so $S_i\leq R_i$. We still need to prove that if $R_i=S_i$
then $R_j=S_j$ for $j\geq i$, i.e. the last part of (ii), and
$R_j=S_j$ for $j>2k+1$, i.e. (i).

If $R_j=S_j$ for all $j\geq 2$ then both statements are
trivial. Otherwise let $2\leq d\leq n$ be maximal s.t. $R_d\neq
S_d$. By Lemma 5.6(ii) we have $d\ev b-1=1\m2$, $R_d<S_d$ and
$S_d=S_{d-2}=\ldots=S_{b-1}$. If we write $d=2l+1$ this means
$S_1=S_3=\ldots=S_{2l+1}$. By the maximality of $k$ we get $l\leq
k$. If $j>2k+1$ then $j>2l+1=d$ so $R_j=S_j$. For the other statement
suppose that $R_i=S_i$. Assume first that $i<2l+1$. Then $i+1\leq
2l+1=d$ so $R_{i+1}\leq R_d<S_d=S_1=S_{i+1}$ ($i+1$ and $d$ are odd
and $\leq 2k+1$). Since also $R_i+R_{i+1}=S_i+S_{i+1}$ we get
$R_i>S_i$. Contradiction. Since $i$ is even we get $i>2l+1$ so for any
$j\geq i$ we have $j> 2l+1=d$ so $R_j=S_j$. \qed

\bco Suppose that $N\ap\[ b,a_2\ldots,a_n\]$ relative to a (possibly
bad) BONG $y,x_2,\ldots,x_n$ s.t. the BONG $x_2,\ldots,x_n$ is good
and $R_i=\ord a_i$. Let $l\geq 1$ be maximal with the property that
$R_{2l+1}<S_1$. (If $R_3\geq S_1$ or $n=2$ we just take $l=0$.) Then
$S_i$'s are given in terms of $S_1$ and $R_2,\ldots,R_n$ by $S_i=R_i$
for $i>2l+1$, $S_1=S_3=\ldots=S_{2l+1}$ and $S_i=R_i+R_{i+1}-S_1$ for
$1<i<2l+1$ even.
\eco
\pf $y\in N$ is a norm generator and $pr_{y^\pp}N=\[ x_2,\ldots,x_n\]$
so we are in the situation of Lemma 5.7.  Let $k\geq 0$ be maximal
s.t. $S_1=S_3=\ldots=S_{2k+1}$. If $j>k$ then $S_1<S_{2j+1}=R_{2j+1}$
by Lemma 5.7(i) so we must have $l\leq k$. In particular,
$S_1=S_3=\ldots=S_{2l+1}$ and by Lemma 5.7(ii) for any $1<i<2l+1\leq
2k+1$ even we have $R_i+R_{i+1}=S_i+S_1$ so $S_i=R_i+R_{i+1}-S_1$.

So we are left to show that $S_i=R_i$ for $i>2l+1$. If $l=k$ this
follows directly from Lemma 5.7(i). If $l<k$ then $R_{2l+2}\geq
S_{2l+2}$ and $R_{2l+2}+R_{2l+3}=S_{2l+2}+S_1$, by Lemma 5.7(ii), and
$R_{2l+3}\geq S_1$, by the maximality of $l$. It follows that
$R_{2l+3}=S_1$ and $R_{2l+2}=S_{2l+2}$. By Lemma 5.7(ii) this implies
$R_i=S_i$ for $i\geq 2l+2$, as claimed. \qed

\bco (i) If $y\in N$ is norm generator then $\nn pr_{y^\pp}N\sbq
\p^{S_2}$ with equality iff $y$ is a first element in a good BONG
of $N$.

(ii) If $S_1<S_3$ or $S_2-S_1=2e$ or $n\leq 2$ then any norm generator
of $N$ is a first element in a good BONG.
\eco
\pf: (i) We use the notations from Lemma 5.7. By Lemma 5.7(ii) we have
$\ord\nn pr_{y^\pp}N=R_2\geq S_2$ so  $\nn pr_{y^\pp}N\sbq
\p^{S_2}$. (If $k=0$ then $R_2=S_2$ by Lemma 5.7(i).) If $R_2=S_2$ then
$R_i=S_i$ for any $i\geq 2$. Hence $N\ap\[ Q(y),a_2,\ldots,a_n\]$
relative to the good BONG $y,x_2,\ldots,x_n$. (The orders of
$Q(y),a_2,\ldots,a_n$ are $S_1,\ldots,S_n$.) Conversely, if $y$ is a
first element in a good BONG $y_1'=y,y_2',\ldots,y_n'$ of $N$ then
$\ord Q(y_i')=S_i$ by [B3, Lemma 4.7]. Hence $\ord\nn
pr_{y^\pp}N=\ord\nn\[ y_2',\ldots,y_n'\]=S_2$. 

(ii) If $S_1<S_3$ then $k=0$. Since $2>2k+1$ we have $R_2=S_2$
by Lemma 5.7(i) so we can apply (i). If $S_1=S_3$ and $S_2-S_1=2e$
then $S_3-S_2=S_1-S_2=-2e$. Suppose that $y$ is not a first element in
a good BONG of $M$ so $R_2>S_2$ by (i). By Lemma 5.7(ii) we also have
$R_2+R_3=S_2+S_3$ so
$R_3-R_2=R_2+R_3-2R_2<S_2+S_3-2S_2=S_3-S_2=-2e$. But this is
impossible. If $n\leq 2$ then all BONGs are good so our statement
follows. \qed

\bco Let $N\sbq M$ be lattices with $M\ap\[ a_1,\ldots,a_m\]$ and
$N\ap\[b_1,\ldots,b_n\]$ relative to the good BONGs $x_1,\ldots,x_m$
and $y_1,\ldots,y_n$ with $\ord a_i=R_i$, $\ord b_i=S_i$. Suppose that
$R_1=S_1$,\ldots,$R_i=S_i$ and one of the following happens: $i\geq
m-1$, $R_{i+1}=S_{i+1}$, $R_i<R_{i+2}$, $R_{i+1}-R_i=2e$. Then there
is a good BONG of $M$ that begins with $y_1,\ldots,y_i$.
\eco
\pf If $i=m$ then $n=m$ and $\ord
volM=R_1+\cdots+R_m=S_1+\cdots+S_m=\ord volN$ so $M=N$ and
$y_1,\ldots,y_m$ is a good BONG of $M$. For the other cases we use
induction on $i$. 

If $i=1$ the condition $R_1=S_1$ implies that $y_1$ is a norm
generator for $M$. If $R_1<R_3$ or $R_2-R_1=2e$ or $m\leq 2$
(i.e. $1\geq m-1$) then $y_1$ is a first element in a BONG of $M$ by
Corollary 5.9(ii). If if $R_2=S_2$ then $\nn pr_{y_1^\pp}M\spq\nn
pr_{y_1^\pp}N=\p^{S_2}=\p^{R_2}$. By Corollary 5.9(i) we have $\nn
pr_{y_1^\pp}M=\p^{R_2}$ and $y_1$ is a first element in a BONG of
$M$.

For $i>1$ we use the induction step. Since $R_1=S_1$ and $R_2=S_2$ we
have that $y_1$ is a first element in a good BONG of $M$ by the case
$i=1$ proved above. Thus $M\ap\[ b_1,a'_2,\ldots,a'_m\]$ relative to
another good BONG $y_1,x'_2,\ldots,x'_m$. We have $\ord a'_i=R_i$. If
$M^*=pr_{y_1^\pp}M$ and $N^*=pr_{y_1^\pp}N$ then $N^*\sbq M^*$. Since
$\rr (M^*)=(R_2,\ldots,R_m)$ and $\rr (N^*)=(S_2,\ldots,S_n)$ the
lattices $M^*,N^*$ satisfy the hypothesis of our corollary at index
$i-1$. By the induction hypothesis the good BONG $x'_2,\ldots,x'_m$ of
$M^*=pr_{y_1^\pp}M$ can be replaced by another good BONG
$y_2,\ldots,y_i,x''_{i+1},\ldots,x''_m$. Therefore
$y_1,\ldots,y_i,x''_{i+1},\ldots,x''_m$ is a good BONG for $M$. \qed

\subsection{Proof of 2.1(ii).} 

By 5.2 we can restrict to the case when $i\leq n_{k_2}+a-1$.

\bff If $a=2$ we make a further reduction by using duality.
In this case $k_1=k_2$, $R_{n_{k_1}+1}=v'$ and $S_{n_{k_1}+1}=v$ so
$R_{n_{k_1}+1}\leq S_{n_{k_1}+1}\leq R_{n_{k_1}+1}+2$. We claim that
the statement 2.1(ii) in the case $i=n_{k_1}+1$ and $S_i=R_i+2$ is
equivalent to statement 2.1(ii) for $N^\*,M^\*$ in the case
$i=n_{k_1}+1$ and $S_i=R_i$. Indeed, the case $i=n_{k_1}+1$ and
$S_i=R_i$ for $N^\*,M^\*$ means statement 2.1(ii) at
$i=n-n_{k_1}-a+1=n-n_{k_1}-1$ when $S_{n-n_{k_1}-1}^\*
=R_{n-n_{k_1}-1}^\*$. ($n_{k_1}=n_{k_2}$ corresponding to $N^\*,M^\*$
is $n-n_{k_1}-a$. See 5.2.) But 2.1(ii) for $N^\*,M^\*$ at
$i=n-n_{k_1}-1$ is equivalent to 2.1(ii) for $M,N$ at 
$i=n_{k_1}+1$. Also $S_{n-n_{k_1}-1}^\* =-R_{n_{k_1}+2}$ and
$R_{n-n_{k_1}-1}^\* =-S_{n_{k_1}+2}$ so $S_{n-n_{k_1}-1}^\*
=R_{n-n_{k_1}-1}^\*$ means $S_{n_{k_1}+2}=R_{n_{k_1}+2}$. But
$S_{n_{k_1}+1}+S_{n_{k_1}+2}=2r=2r'+2=R_{n_{k_1}+1}+R_{n_{k_1}+2}+2$
so $S_{n_{k_1}+2}=R_{n_{k_1}+2}$ is equivalent to
$S_{n_{k_1}+1}=R_{n_{k_1}+1}+2$, as claimed.

As a consequence, when we prove (ii) at $i=n_{k_1}+1$ when $a=2$ we
can ignore the case $S_i=R_i+2$ and consider only the cases $S_i=R_i$
and $S_i=R_i+1$.
\eff

\bff Let $\AA_k$ be norm generators for $L^{\ss L_k}$. If $a=1$
let $J\ap\la\AA\ra$ and $J'\ap\la\AA'\ra$. We can take $\AA
=\pi^2\AA'$ since $J=\p J$. If $a=2$ then let $J\ap\[\AA,\(\AA\]$ and
$J\ap\[\AA',\(\AA'\]$. If $u=u'$ then $\nn J=\nn J'$ so any norm
generator of $J$ is also a norm generator of $J'$. So in this case we
can take $\AA'=\AA$. Also $\det J=\pi^2\det J$ so we may assume
$\(\AA=\pi^2\(\AA'$. If $u=u'+2$ then $\ord\nn\p J'=u'+2=u=\ord\nn J$
and $\p J'\sbq J$ so a norm generator for $\p J'$ is also a norm
generator for $J$. Since $\p J'\ap\[\pi^2\AA',\pi^2\(\AA'\]$ we can
take $\AA =\pi^2\AA'$. Since $\det J=\pi^2\det J'$ we can take $\(\AA
=\(\AA'$.

We denote by $\BB_k$, $\BB'_k$, $\BB$ and $\BB'$ some norm generators
for $N^{\p^{r_k}}$, $M^{\p^{r_k}}$, $N^{\p^r}$ and
$M^{\p^{r'}}$ respectively.

If $k\leq k_2$ then $N^{\p^{r_k}}=K^{\p^{r_k}}\pp J$ and $\AA_k$ and
$\AA$ are norm generators for $K^{\p^{r_k}}$ and $J$ with orders $u_k$
and $u$ we can take $\BB_k=\AA_k$, if $v_k=u_k$, or $\BB_k=\AA$, if
$v_k=u$. 

If $k>k_2$ then $N^{\p^{r_k}}=K^{\p^{r_k}}\pp\p^{r_k-r}J$ and $\AA_k$
and $\pi^{2(r_k-r)}\AA$ are norm generators for $K^{\p^{r_k}}$ and
$\p^{r_k-r}J$ $\AA_k$ with orders $u_k$ and $u+2(r_k-r)$ we can take
$\BB_k=\AA_k$, if $v_k=u_k$, or $\BB_k=\pi^{2(r_k-r)}\AA$, if
$v_k=u+2(r_k-r)$.

Now $\p^{r-r_{k_2}}K^{\p^{r_{k_2}}}$, $K^{\p^{r_{k_2+1}}}$ and $J$ are
sublattices of $N^{\p^r}$ with norms of orders and $v=\min\{
u_{k_2}+2(r-r_{k_2}),u_{k_2+1},u\}$. Since
$\pi^{2(r-r_{k_2})}\AA_{k_2}$, $\AA_{k_2+1}$ and $\AA$ are norm
generators for the three lattices we can take $\BB
=\pi^{2(r-r_{k_2})}\AA_{k_2}$, if $v=u_{k_2}+2(r-r_{k_2})$, or
$\BB =\AA_{k_2+1}$, if $v=u_{k_2+1}$ or $\BB =\AA$, if  or $v=u$. 

Similarly for the lattice $M$ with $k_2,u,v,v_k,\AA,\BB,\BB_k$
replaced by $k_1,u',v',v'_k,\AA',\BB',\BB'_k$.
\eff


\blm Suppose that $i\leq n_{k_1}+a-1$ and, if $i=n_{k_1}+1$, then
$R_i=S_i$. Then:

(i) If $S_i\neq R_i+1$ then the approximations $X_i,Y_i$ can be chosen
equal.

(ii) If $S_i=R_i+1$ then $a_{1,i}b_{1,i}$ has an odd order.
\elm
\pf (i) Suppose first that $i\leq n_{k_1}$. Then $n_{k-1}+1\leq i\leq
n_k$ for some $1\leq k\leq k_1$. We use Lemma 3.2. If $i\ev
n_{k-1}\m2$ then we can choose $X_i=Y_i=(-1)^{(i-n_{k-1})/2}\det
FK_{(k-1)}$. Suppose now that $i\ev n_{k-1}+1\m2$ so $X_i=\BB'_k\det
FK_{(k-1)}$ and $Y_i=\BB_k\det FK_{(k-1)}$. We have $R_i=v'_k$ and
$S_i=v_k$. $S_i\neq R_i+1$ implies $v_k\neq v_k'+1$ so
$v_k=v_k'$ or $v_k=v_k'+2$. (We have $v_k=\min\{
u_k,u\}$, $v'_k=\min\{ u_k,u'\}$ and $0\leq u-u'\leq 2$.) If
$v_k=v_k'$ then either $v_k=v'_k=u_k$ or $v_k=v'_k=u=u'$. In the first
case we take $\BB'_k=\BB_k=\AA_k$ while in the second we take
$\BB'_k=\AA'=\AA =\BB_k$. Since $\BB'_k=\BB_k$ we have $X_i=Y_i$. If
$v_k=v_k'+2$ then $v_k'=u'$ and $u=u'+2=v_k$. Thus we can take
$\BB'_k=\AA'$ and $\BB_k=\AA =\pi^2\AA'$. It follows that
$\BB'_k=\BB_k$ in $\ff/\fs$ so $X_i=Y_i$.

If $i>n_{k_1}$ then $a=2$ and $i=n_{k_1}+1$ so, by hypothesis,
$R_i=S_i$. Also $k_1=k_2$. By Lemma 3.2, in both cases when the
$\p^r$-modular component of $N$ is $J$ or $J\pp K_{k_1+1}$, we can
take $Y_i=\BB\det FK_{(k_1)}$. By Corollary 3.3, in both cases when
the $\p^{r'}$-modular component of $M$ is $J'$ or $K_{k_1}\pp J'$ we
can take $X_i=\BB'\det F(K_{(k_1)}\pp J')$ for any norm generator of
$M^{\p^{r'}}$, $\BB'$. (We have $\dim (K_{(k_1)}\pp J')=n_{k_1}+2$.)
Now $M^{\p^{r'}}=K^{\p^{r'}}\pp J'\sp K^{\p^r}\pp J=N^{\p^r}$ and,
since $v=S_i=R_i=v'$, we have $\nn M^{\p^{r'}}=\nn
N^{\p^r}$. Therefore $\BB$, the norm generator of $N^{\p^r}$, is also
a norm generator of $M^{\p^{r'}}$, and so is $-\BB$. Therefore we can
take $\BB'=-\BB$ so $X_i=-\BB\det F(K_{(k_1)}\pp J')=-Y_i\det FJ'$.
In order to prove that $Y_i$ also is an approximation for $a_{1,i}$ we
need to prove that $d(-\det FJ')=d(X_iY_i)\geq\a_i$. Since $J'\sbq
M^{\p^{r'}}$ we have $\ww J'\sbq\ww M^{\p^{r'}}$. By [B3, Lemma
2.16(i)] we have $\ord\ww J'=(R_1+\a_1)(J')\leq
(R_1+R_2-R_1+d(-a_{1,2}))(J')=R_2(J')+d(-\det FJ')=2r'-u'+d(-\det
FJ')$ and $\ord\ww M^{\p^{r'}}=R_i+\a_i=v'+\a_i$. Hence
$2r'-u'+d(-\det FJ')\geq v'+\a_i$. Since $v'\geq r'\geq 2r'-u'$ we get
$d(-\det FJ')=d(X_iY_i)\geq\a_i$, as claimed.

(ii) Now for $k\leq k_1$ the sequence $R_{n_{k-1}+1},\ldots,R_{n_k}$
is $v'_k,2r_k-v'_k,\ldots,v'_k,2r_k-v'_k$ if $v'_k>r_k$ and it is
$r_k,r_k,\ldots,r_k$ if $v'_k=r_k$. In both cases
$\sum_{j=n_{k-1}+1}^{n_k}R_j=(n_k-n_{k-1})r_k$ and for $n_{k-1}\leq
i\leq n_k$ with $k\leq k_1$ s.t. $i\ev n_{k-1}\m2$ or $i\ev n_k\m2$ we
have $\sum_{j=n_{k-1}+1}^iR_j=(i-n_{k-1})r_k$. Similarly for $N$ so
for any $n_{k-1}\leq i\leq n_k$ with $i\ev n_{k-1}\m2$ we have
$R_1+\cdots +R_i=S_1+\cdots
+S_i=\sum_{l=1}^{k-1}(n_l-n_{l-1})r_l+(i-n_{k-1})r_k$.

Now if $i\leq n_{k_1}+a-1$ and $S_i=R_i+1$ then either $a=2$ and
$i=n_{k_1}+1$ or $n_{k-1}<i\leq n_k$ for some $k\leq k_1$ and $i\ev
n_{k-1}+1\m2$. (If $i\ev n_{k-1}\m2$ then $R_i=2r_k-v'_k\geq
2r_k-v_k=S_i$.) In both cases we have
$R_1+\cdots+R_{i-1}=S_1+\cdots+S_{i-1}$ (see above). Together with
$S_i=R_i+1$, this implies $\ord b_{1,i}=S_1+\cdots+ S_i=R_1+\cdots+
R_i+1=\ord a_{1,i}+1$. Thus $\ord a_{1,i}b_{1,i}$ is odd. \qed



\blm Let $L=L_1\pp\cdots\pp L_t$ be a Jordan splitting. If $i=n_k$ for
some $1\leq k\leq t-1$ then $\a_i=\min\{ -2r_k+\od\sum\dd (\a\b
),(u_k+u_{k+1})/2-r_k+e\}$ where $\a$ runs over $\GG_k$ and $\b$ over
$\GG_{k+1}$. Also $(u_k+u_{k+1})/2-r_k+e=(R_{i+1}-R_i)/2+e$.

(Here we use the notations from [B3, \S2].)
\elm
\pf We have $R_{i+1}=u_{k+1}$ and $R_i=2r_k-u_k$ by [B3, Lemma 2.13(i)
and (ii)]. Hence $(R_{i+1}-R_i)/2+e=(u_k+u_{k+1})/2-r_k+e$. Also
$R_{i+1}-R_i$ is odd iff $u_k+u_{k+1}$ is odd.

If $R_{i+1}-R_i$ is odd then $u_k+u_{k+1}$ is odd so $\sum\dd (\a\b
)=\AA_k\AA_{k+1}\oo =\p^{u_k+u_{k+1}}$. Therefore $\min\{
-2r_k+\od\sum\dd (\a\b ),(u_k+u_{k+1})/2-r_k+e\} =\min\{
u_k+u_{k+1}-2r_k,(u_k+u_{k+1})/2-r_k+e\} =\min\{
R_{i+1}-R_i,(R_{i+1}-R_i)/2+e\}$, which is equal to $\a_i$ by [B3,
Corollary 2.9(ii)].

If $R_{i+1}-R_i$ is even, so $u_k+u_{k+1}$ is even, then
$\a_i=\od\FF_k$, by [B3, Lemma 2.16(ii)] and $\ss_k^2\FF_k=\sum\dd
(\a\b )+2\p^{(u_k+u_{k+1})/2+r_k}$. Thus
$2r_k+\od\FF_k=\min\{\od\sum\dd (\a\b ),(u_k+u_{k+1})/2+r_k+e\}$. So
$\a_i=\od\FF_k=\min\{ -2r_k+\od\sum\dd (\a\b
),(u_k+u_{k+1})/2-r_k+e\}$. \qed 

\blm Let $L$ be a binary $\p^r$-modular lattice of norm $\p^u$. Then
there is $y\in L$ with $\ord Q(z)=u+1$ iff either $L\ap\pi^r\aa$ or
$u<r+e$ and $d(-\det FL)=2u-2r+1$. 
\elm
\pf Our statement is invariant to scaling so we may assume that $L$ is
unimodular, i.e. $r=0$.

Suppose first that $r=e$ so $L\ap\aa$ or $\ab$ relative to sopme basis
$x,y$. In the first case we may take $z=x+\pi y$ so $\ord Q(z)=\ord
2\pi =e+1$. In the second case $FL\ap [2,-2\D ]$ so $Q(FL)\setminus\{
0\} =2\ooo\fs$. In particular, $FL$ doesn't represent elemnts of $F$
of order $e+1$.

Suppose now that $u<e$. If $-\det L=1+\a$ with $\dd (1+\a )=\a\oo$ and
$\ord\a =d$ then $d(-\det FL)=d$. We write as in [OM, 93:17] $\ap A(a,-\a
a\1 )$, relative to some basis $x,y$, where $\ord a=u$. We have
$d-u=\ord\a\1 a\geq\ord\nn L=u$ so $d\geq 2u$. But $d$ cannot be odd
and $<2e$ so $d\geq 2u+1$. If $d=2u+1$ then we take $z=y$ and we have
$\ord Q(z)=\ord\a\1 a=d-u=u+1$. Conversely, assume that $d>2u+1$ and
there is $z\in J$ with $\ord Q(z)=u+1$. If $z=sx+ty$ with $s,t\in\oo$
then $\ord (s^2a +2st-t^2\a a\1 )=u+1$. But $\ord t^2\a a\1\geq
d-u>u+1$ and $\ord 2st\geq e\geq u+1$ so $\ord (s^2u+2st)=u+1$ and
$\ord s^2\a\geq u+1$ so $\ord s>0$. Thus $\ord s^2\a\geq u+2$ and
$\ord 2st\geq e+1\geq u+2$ so $\ord (s^2+2st)\geq
u+2$. Contradiction. \qed

\bff Note that if $s\leq r'$ then also $s<r$ and so
$M^{\p^s}=K^{\p^s}\pp J'\sp K^{\p^s}\pp J=N^{\p^s}$.
\eff

\blm Let $i\leq n_{k_1}+a-1$ s.t. $R_i=S_i$. We have:

(i) $\a_i\leq\b_i$.

(ii) $R_j=S_j$ for all $j\leq i$.
\elm
\pf (i) If $n_{k-1}<i<n_k$ for some $k\leq k_1$ then $R_i+\a_i=\ord\ww
M^{\p^{r_k}}$ and $S_i+\b_i=\ord\ww N^{\p^{r_k}}$ ([B3, Lemma
2.16(i)]). Since $M^{\p^{r_k}}\sp N^{\p^{r_k}}$ we have $R_i+\a_i\leq 
S_i+\b_i=R_i+\b_i$ so $\a_i\leq\b_i$.

If $i=n_k$ with $k\leq k_1$ then by Lemma 5.14 we have $\b_i=\min\{
-2r_k+\ord\sum_1\dd (\a\b ),(S_{i+1}-S_i)/2+e\}$, where $\a,\b$ in
$\sum_1$ run over $\GG N^{\p^{r_k}}$ and $\GG N^{\p^s}$, respectively,
where $s=\min\{ r_{k+1},r\}$. (We have $s=r$ when $k=k_1=k_2$. In this
case, in the Jordan decomposition for $N$, the scale following
$\p^{r_k}$ is $\p^r$, not $\p^{r_{k+1}}$.) Similarly $\a_i=\min\{
-2r_k+\ord\sum_2\dd (\a\b ),(R_{i+1}-R_i)/2+e\}$, where $\a,\b$ in
$\sum_2$ run over $\GG M^{\p^{r_k}}$ and $\GG M^{\p^{s'}}$,
respectively, where $s'=\min\{ r_{k+1},r'\}$. By 5.16 $M^{\p^{r_k}}\sp
N^{\p^{r_k}}$ and, since $s'\leq s$, $M^{\p^{s'}}\sp N^{\p^{s'}}\spq
N^{\p^s}$. So $\sum_1\dd (\a\b )\sbq\sum_2\dd (\a\b )$, which implies
$-2r_k+\ord\sum_1\dd (\a\b )\geq -2r_k+\ord\sum_1\dd (\a\b )$. Also
since $S_i=R_i$ and $S_i+S_{i+1}\geq R_i+R_{i+1}$) we have
$(R_{i+1}-R_i)/2+e\leq (S_{i+1}-S_i)/2+e$. It follows that
$\a_i\leq\b_i$.

We have an exception when $k=k_1$ and $r=r_{k_1}$. In
this case $K_k\pp J'$ is the $\p^{r_k}=\p^{r'}$ modular Jordan
component of $M$. Hence $R_i+\a_i=\ord\ww M^{\p^{r_k}}\leq\ord\ww
N^{\p^{r_k}}\leq S_i+\b_i=R_i+\b_i$ so $\a_i\leq\b_i$. (If $\dim
K_k>1$ then $n_{k-1}<i-1<n_k$ so $\ord\ww
N^{\p^{r_k}}=S_{i-1}+\b_{i-1}\leq S_i+\b_i$.  If $K_k$ is unary then
$\ord\BB_k\1\ww N^{\p^{r_k}}=\min\{ \b_{i-1},\b_i,e\}$ by [B3,
Corollary 2.17(ii)] so $\ord\ww
N^{\p^{r_k}}\leq\ord\BB_k+\b_i=r_k+\b_i=S_i+\b_i$.) 

Finally, if $a=2$ and $i=n_{k_1}+1$ then $M^{\p^{r'}}\sp
N^{\p^{r'}}\spq N^{\p^{r}}$ so $R_i+\a_i=\ord\ww
M^{\p^{r'}}\leq\ord\ww N^{\p^r}=S_i+\b_i=R_i+\b_i$ so $\a_i\leq\b_i$.

(ii) We have $n_{k-1}<i\leq n_k$ for some $k\leq k_1$ or $a=2$ and
$i=n_{k_1}+1$. Let $j<i$. We have $n_{k'-1}<i\leq n_{k'}$ for some
$k'\leq k$, resp. $k'\leq k_1$ (when $i=n_{k_1}+1$). Suppose first
that $n_{k-1}<i\leq n_k$ for some $k\leq k_1$. If $i\ev n_{k-1}+1\m2$
then $R_i=v'_k$ and $S_i=v_k$ while if $i\ev n_{k-1}\m2$
then $R_i=2r_k-v'_k$ and $S_i=2r_k-v_k$. In both cases $R_i=S_i$
implies $v'_k=v_k$. Now $v_k=\min\{ u_k,u\}$ and $v_{k'}=\min\{
u_{k'},u\}$. Since $u_{k'}\leq u_k$ we have $v_{k'}=\min\{
u_{k'},v_k\}$. Similarly $v'_{k'}=\min\{ u_{k'},v'_k\}$ and since
$v'_k=v_k$ we get $v'_{k'}=v_{k'}$. If $i=n_{k_1}+1$ then
$v'=R_i=S_i=v$. We have $v_{k'}\leq v\leq u$ and $v_{k'}=\min\{
u_{k'},u\}$ so $v_{k'}=\min\{ u_{k'},v\}$ and similarly
$v'_{k'}=\min\{ u_{k'},v'\}$. But $v'=v$ so we get again
$v'_{k'}=v_{k'}$. This implies $R_j=S_j=v_{k'}$ or $2r_{k'}-v_{k'}$,
depending on the parity of $j$. \qed

We start now our proof in the case $i\leq n_{k_1}+a-1$. By the
reduction from 5.11, if $a=2$ and $i=n_{k_1}+1$ then we may assume that
$S_i=R_i$ or $S_i=R_i+1$. We have 2 cases: 

1. $S_i\neq R_i+1$. By Lemma 5.13(i) we can assume that $X_i=Y_i$ so
$d(X_iY_i)=\j$. Hence $d[a_{1,i}b_{1,i}]=d[X_iY_i]=\min\{
d(X_iY_i),\a_i,\b_i\} =\min\{\a_i,\b_i\}$. So we have to prove that
$\min\{\a_i,\b_i\}\geq A_i$. We have several subcases:

a. $R_i=S_i$. By Lemma 5.17(i) we have $\a_i\leq\b_i$ so we must prove
that $\a_i\geq A_i$.

By Lemma 5.17(ii) we have $R_1=S_1$,\ldots,$R_i=S_i$, which, by
Corollary 5.10, implies that we can change the good BONG of $M$
s.t. $a_j=b_j$ for $1\leq j\leq i-1$.

If $\a_i=(R_{i+1}-R_i)/2+e$ then $\a_i=(R_{i+1}-S_i)/2+e\geq A_i$. If
$\a_i=R_{i+1}-R_i+d(-a_{i,i+1})$ then
$\a_i=R_{i+1}-S_i+d(-a_{1,i+1}a_{1,i-1})=
R_{i+1}-S_i+d(-a_{1,i+1}b_{1,i-1})\geq A_i$. If
$\a_i=R_{i+1}-R_j+d(-a_{j,j+1})$ with $j<i$ then
$\a_i=R_{i+1}-S_j+d(-b_{j,j+1})=R_{i+1}-S_i+S_i-S_j+d(-b_{j,j+1})\geq
R_{i+1}-S_i+\b_{i-1}\geq A_i$. Finally, if
$\a_i=R_{j+1}-R_i+d(-a_{j,j+1})$ with $j>i$ then
$\a_i=R_{i+1}-S_i+R_{j+1}-R_{i+1}+d(-a_{j,j+1})\geq
R_{i+1}-S_i+\a_{i+1}\geq A_i$ so we are done.

For the following cases, b. and c., we have both $S_i\neq R_i$ and
$S_i\neq R_i+1$. This rules out the case $i=n_{k_1}+1$. Hence
$n_{k-1}+1\leq i\leq n_k$ for some $1\leq k\leq k_1$.

b. $R_i<S_i$. We cannot have $i\ev n_{k-1}\m2$ since this would imply
$R_i=2r_k-v'_k\geq 2r_k-v_k=S_i$. Thus $i\ev n_{k-1}+1\m2$ and we have
$\min\{ u_k,u'\}=v_k'=R_i<S_i=v_k=\min\{ u_k,u\}$. This implies
$r_k\leq v'_k=u'<u_k$ so $K_k$ is improper. Thus $n_k\ev
n_{k-1}\m2$ so $i\ev n_k+1\m2$. In particular, $i<n_k$. If $i\neq
n_k-1$ then $i+2\leq n_k-1$ so $R_{i+2}=v'_k=R_i$. If $i=n_k-1$ and
$k<k_1$ then $R_{i+2}=v'_{k+1}=\min\{ u_{k+1},u'\}=u'=v'_k=R_i$. (We
have $u_{k+1}\geq u_k>u'$.) Finally, if $i=n_k-1$ and $k=k_1$ we have
$R_{i+2}=v'\leq u'=R_i$ so $R_i=R_{i+2}$. So anyways we have
$R_i=R_{i+2}$, which implies $R_i+\a_i=R_{i+1}+\a_{i+1}$ by [B3,
Corollary 2.3(i)]. It follows that $A_i\leq
R_{i+1}-S_i+\a_{i+1}=R_i-S_i+\a_i<\a_i$. On the other hand by 5.16
$M^{\p^{r_k}}\sp N^{\p^{r_k}}$ so $S_i+\b_i=\ord\ww
N^{\p^{r_k}}\geq\ord\ww M^{\p^{r_k}}=R_i+\a_i$. It follows that
$\b_i\geq R_i-S_i+\a_i\geq A_i$ and we are done. 

c. $R_i>S_i$. This implies $R_i+R_{i+1}\leq S_{i-1}+S_i$ and
$R_{i+1}<S_{i-1}$. So $R_{i-1}<S_{i-1}$, which, by the same reasoning
from b., implies $R_{i-1}=R_{i+1}$. (In the proof of b. $R_i<S_i$
implies $R_i=R_{i+2}$.) We have $A_i\leq
R_{i+1}-S_i+\b_{i-1}<S_{i-1}-S_i+\b_{i-1}\leq\b_i$. We have to prove
that also $\a_i\geq A_i$. Note that we cannot have $i\ev n_{k-1}+1\m2$
since this would imply $R_i=v'_k\leq v_k=S_i$. So $i\ev n_{k-1}\m2$
and we have $2r_k-v'_k=R_i>S_i=2r_k-v_k$ so $v'_k<v_k$. We have
$\min\{ u_k,u'\} =v'_k<v_k=\min\{ u_k,u\}$ and $u-u'\leq 2$ so
$v'_k=v_k-1$ or $v'_k=v_k-2$. 

If $v'_k=v_k-1$ then $R_i=2r_k-v'_k=2r_k-v_k+1=S_i+1$ and
$S_{i-1}=v_k=v'_k+1=R_{i-1}+1$. By Lemma 5.13(ii) $\ord
a_{1,i-1}b_{1,i-1}$ is odd. Also $\ord
a_{i,i+1}=R_i+R_{i+1}=R_{i-1}+R_i=2r_k$ is even so $\ord
a_{1,i+1}b_{1,i-1}$ is odd. Hence $d[-a_{1,i+1}b_{1,i-1}]=0$, which
implies $A_i\leq
R_{i+1}-S_i+d(-a_{1,i+1}b_{1,i-1})=R_{i+1}-S_i=R_{i+1}-R_i+1$. So it
suffices to prove that $\a_i\geq R_{i+1}-R_i+1$. By [B3, Lemma
2.7(iii) and Corollary 2.8(i)] it is enough to prove that
$R_{i+1}-R_i$ is even and $\leq 2e$. We have
$R_{i+1}-R_i=R_{i-1}-R_i=v'_k-(2r_k-v'_k)=2(v'_k-r_k)$ so it is
even. Also $v'_k-r_k<v_k-r_k\leq e$ so $R_{i+1}-R_i<2e$ and we are
done. 

If $v'_k=v_k-2$ then $R_i=S_i+2$. Since $\min\{ u_k,u'\}
=v'_k=v_k-2=\min\{ u_k,u\} -2$ and $u'\geq u-2$ we have $v'_k=u'$,
$v_k=u$ and $u'=u-2$. By 5.4 we get $v'=u'$ and $v=u$. By 5.12 we can
take $\BB_k=\BB =\AA$ and $\BB'_k=\BB'=\AA'$. Also, since $u'=u-2$, we
may assume that $\AA=\pi^2\AA'$. Let $x\in J'$ with $Q(x)=\AA'$. Now
$x$ is a norm generator for $J'$ and $\nn J'\sp\nn J$ so $x\notin
J$. Since $[J':J]=\p$ we have $J'=J+\oo x$ which implies that
$M^{\p^{r_k}}=N^{\p^{r_k}}+\oo x$. (We have
$M^{\p^{r_k}}=K^{\p^{r_k}}\pp J'$ and $N^{\p^{r_k}}=K^{\p^{r_k}}\pp
J$.) Now $\AA=\pi^2\AA'$ is a norm generator for $N^{\p^{r_k}}$ while
$\AA'$ is a norm generator for both $M^{\p^{r_k}}$ and $\oo x$. By
[B3, Lemma 2.11] we have $\ww M^{\p^{r_k}}=\ww N^{\p^{r_k}}+\ww (\oo
x)+{\AA'}\1\dd (\AA'\AA )+2\ss M^{\p^{r_k}}=\ww
N^{\p^{r_k}}+2\AA'\oo+0+2\p^{r_k}$. But $2\AA'\oo\sbq 2\ss
J'=2\p^{r'}\sbq 2\p^{r_k}$ and $2\p^{r_k}=2\ss N^{\p^{r_k}}\sbq\ww
N^{\p^{r_k}}$ so $\ww M^{\p^{r_k}}=\ww N^{\p^{r_k}}$. Hence
$S_{i-1}+\b_{i-1}=\ord\ww N^{\p^{r_k}}=\ord\ww
M^{\p^{r_k}}=R_{i-1}+\a_{i-1}$. (Note that $n_{k-1}+1\leq i\leq n_k$
and $i\ev n_{k-1}\m2$ so $n_{k-1}<i-1<n_k$ and we can apply [B3, Lemma
2.16(i)].) Since $S_{i-1}+S_i\geq R_i+R_{i+1}$ we get $A_i\leq
R_{i+1}-S_i+\b_{i-1}\leq S_{i-1}-R_i+\b_{i-1}=
R_{i-1}-R_i+\a_{i-1}\leq\a_i$ and we are done. 

2. Now we consider the case $S_i=R_i+1$. By Lemma 5.13(ii), $\ord
a_{1,i}b_{1,i}$ is odd, which implies $d[a_{1,i}b_{1,i}]=0$. So we
have to prove that $0\geq A_i$. Since $R_i<S_i$ we have either
$n_{k-1}<i\leq n_k$ for some $k\leq k_1$ and $i\ev n_{k-1}+1\m2$ or
$a=2$ and $i=n_{k_1}+1$. (See case 1.b.) In the first case, by the
same reasoning from case 1.b., we get $R_i=R_{i+2}$ and so
$R_i+\a_i=R_{i+1}+\a_{i+1}$ and $A_i\leq
R_{i+1}-S_i+\a_{i+1}=R_i-S_i+\a_i=\a_i-1$. Note that $\nn
M^{\p^{r_k}}=\p^{v'_k}=\p^{R_i}$ and $\nn
N^{\p^{r_k}}=\p^{v_k}=\p^{S_i}$. Now $\GG M^{\p^{r_k}}\spq\GG
N^{\p^{r_k}}$ so $\GG M^{\p^{r_k}}$ contains elements of order
$S_i=R_i+1$. But $\nn M^{\p^{r_k}}=\p^{R_i}$ so $\ww
M^{\p^{r_k}}\spq\p^{R_i+1}$. It follows that $R_i+\a_i=\ord\ww
M^{\p^{r_k}}\leq R_i+1$. (Same as in case 1.b., we have $n_k\ev
n_{k-1}\m2$ so $i\ev n_k+1\m2$, which implies $n_{k-1}<i<n_k$. This is
why we can apply [B3, Lemma 2.16(i)].) Thus $\a_i\leq 1$ and
$A_i\leq\a_i-1\leq 0$.

If $a=2$ and $i=n_{k_1}+1$ then $R_i=v'=\min\{
u',u_{k_1+1},2(r'-r_{k_1})+u_{k_1}\}$ and $S_i=v=\min\{
u,u_{k_1+1},2(r-r_{k_1})+u_{k_1}\}$. We cannot have $R_i=u_{k_1+1}$
since this would imply $R_i\geq S_i$. Also we cannot have
$S_i=2(r-r_{k_1})+u_{k_1}=2(r'-r_{k_1})+u_{k_1}+2$ since this would
imply $S_i\geq R_i+2$. Thus $R_i=v'=\min\{ u',2(r'-r_{k_1})+u_{k_1}\}$
and $S_i=v=\min\{ u,u_{k_1+1}\}$. We have 3 cases:

If $R_i=v'=2(r'-r_{k_1})+u_{k_1}$ then $u\geq u'\geq
2(r'-r_{k_1})+u_{k_1}\geq u_{k_1}$. Thus $v_{k_1}=\min\{ u_{k_1},u\}
=u_{k_1}$ so $S_{i-1}=2r_{k_1}-u_{k_1}\ev
2(r'-r_{k_1})+u_{k_1}=R_i=S_i-1\m2$. Since $S_i-S_{i-1}$ is odd we
have $\b_{i-1}\leq S_i-S_{i-1}$ so $A_i\leq R_{i+1}-S_i+\b_{i-1}\leq
R_{i+1}-S_{i-1}$. But $R_{i+1}=2r'-v'=2r_{k_1}-u_{k_1}=S_{i-1}$ so
$A_i\leq 0$.

If $S_i=v=u_{k_1+1}$ then $u_{k_1+1}\leq u\leq u'+2\leq
2(r_{k_1+1}-r')+u'$ so $v'_{k_1+1}=\min\{
u_{k_1+1},2(r_{k_1+1}-r')+u'\} =u_{k_1+1}$. Thus
$R_{i+2}=v'_{k_1+1}=u_{k_1+1}=S_i=R_i+1\ev R_{i+1}+1\m2$. (We have
$R_i+R_{i+1}=2r'\ev 0\m2$.) Since $R_{i+2}-R_{i+1}$ is odd
we have $\a_{i+1}\leq R_{i+2}-R_{i+1}$ so $A_i\leq
R_{i+1}-S_i+\a_{i+1}\leq R_{i+2}-S_i=0$.

If $R_i=v'=u'$ and $S_i=v=u$ then $u=u'+1$ so
$R_{i+1}-S_i=2r'-u'-u=2r'-2u'-1$. Now $J\sb J'$ and $\nn
J=\p^u=\p^{u'+1}=\p\nn J'$ so $J'$ represents elements of order
$\ord\nn J'+1$. By Corollary 5.15 either $J'$ is hyperbolic or
$d(-\det FJ')=2u'-2r'+1$. In the first case $u'=r'+e$ so
$R_{i+1}-S_i=2r'-2u'-1=-2e-1<-2e$, which implies $A_i\leq
(R_{i+1}-S_i)/2+e<0$. In the second case we use Lemma 3.2(i) and
we have the approximations $X_{i+1}=\det F(K_{(k_1)}\pp J')$ and
$Y_{i-1}=\det FK_{(k_1)}$. Thus $A_i\leq
R_{i+1}-S_i+d(-X_{i+1}Y_{i-1})=R_{i+1}-S_i+d(-\det
FJ)=2r'-2u'-1+2u'-2r'+1=0$.
\vskip 3mm

Thus we have proved (ii) in the case $i\leq n_{k_1}+a-1$. If
$k_1=k_2$ then we are done. Otherwise $a=1$ and $K$ contains a
$\p^{r'+1}=\p^{r-1}$ modular component i.e. $k_2=k_1+1$ and
$r_{k_1+1}=r'+1=r-1$. We have to prove (ii) for $n_{k_1}+1\leq i\leq
n_{k_2}$. We have 2 cases:

3. $K_{k_1+1}$ is proper. As seen in the proof of (i) the sequences
$R_{n_{k_1}+1},\ldots,R_{n_{k_2}+1}$ and
$S_{n_{k_1}+1},\ldots,S_{n_{k_2}+1}$ are $r-2,r-1,\ldots,r-1$ and
$r-1,\ldots,r-1,r$. Since also
$R_1+\cdots+R_{n_{k_1}}=\sum_{l=1}^{k_1}(n_l-n_{l-1})r_l=
S_1+\cdots+S_{n_{k_1}}$ (see the proof of Lemma 5.13(ii)) we
have $R_1+\cdots+R_i=S_1+\cdots+S_i-1$. So $\ord a_{1,i}b_{1,i}$ is
odd, which implies $d[a_{1,i}b_{1,i}]=0$. Thus we have to prove that
$0\geq A_i$. But we also have $R_{i+1}=S_i=r-1$ which, together
with $R_1+\cdots+R_i=S_1+\cdots+S_i-1$, implies
$R_1+\cdots+R_{i+1}=S_1+\cdots+S_{i-1}+2r-3$. Hence $\ord
a_{1,i+1}b_{1,i-1}$ is odd. It follows that $A_i\leq
R_{i+1}-S_i+d(-a_{1,i+1}b_{1,i-1})=(r-1)-(r-1)+0=0$.

4. $K_{k_1+1}$ is improper. This time the sequences
$R_{n_{k_1}+1},\ldots,R_{n_{k_2}+1}$ and
$S_{n_{k_1}+1},\ldots,S_{n_{k_2}+1}$ are $r-2,r,\ldots,r,r-2$ and
$r,r-2,\ldots,r-2,r$. Note that since $\ss J'=\p^{r-2}$, $\ss J=\p^r$
and $r_{k_1+1}=r_{k_2}=r-1$ we have
$M^{\p^{r_{k_1+1}}}=K^{\p^{r-1}}\pp\p J'=K^{\p^{r-1}}\pp
J=N^{\p^{r_{k_2}}}$. We denote $A=\ord\ww
M^{\p^{r_{k_1+1}}}=\ww N^{\p^{r_{k_2}}}$. By [B3, Lemma 2.16(i)]
$R_j+\a_j=\ord\ww M^{\p^{r_{k_1+1}}}=A$ for $n_{k_1}+2=n_{k_1}+a+1\leq
j\leq n_{k_1+1}+a-1=n_{k_2}$. This relation holds also for
$j=n_{k_1}+1$ because at this index we have $R_j=R_{j+2}=r-2$ so, by
[B3, Corollary 2.3], $R_j+\a_j=R_{j+1}+\a_{j+1}=A$. Similarly for
$n_{k_1}+1\leq j\leq n_{k_1+1}-1=n_{k_2}-1$ we have $S_j+\b_j=\ord\ww
N^{\p^{r_{k_1+1}}}=A$. This also holds for $j=n_{k_2}$ since at this
index we have $S_{j-1}=S_{j+1}=r$ so
$S_j+\b_j=S_{j-1}+\b_{j-1}=A$. Thus $R_j+\a_j=S_j+\b_j=A$ holds for
any $n_{k_1}+1\leq j\leq n_{k_2}$.

We have $u_{k_1+1}>r_{k_1+1}=r-1$ and $u=r$ so $v_{k_1+1}=\min\{
u_{k_1+1},u\} =r$. By 5.12 $\AA$ is a norm generator for
$M^{\p^{r_{k_1+1}}}=N^{\p^{r_{k_1+1}}}$ and so is $-\AA$. If
$n_{k_1}+1\leq i\leq n_{k_2}-1$ we use Lemma 3.2 to find some
approximations $X_i,Y_i$ for $a_{1,i}$ and $b_{1,i}$. We take
$Y_i=(-1)^{(i-n_{k_1})/2}\det FK_{(k_1)}$ if $i\ev n_{k_1}\m2$ and
$Y_i=\pm\AA\det FK_{(k_1)}$ if $i\ev n_{k_1}+1\m2$. For $X_i$ we
note that the Jordan splitting of $M$ begins with $K_1\pp\cdots\pp
K_{k_1}\pp J'$, of rank $n_{k_1}+1$. Therefore we can take
$X_i=(-1)^{(i-n_{k_1}-1)/2}\det F(K_{(k_1)}\pp
J')=(-1)^{(i-n_{k_1}-1)/2}\AA\det FK_{(k_1)}$ if $i\ev
n_{k_1}+1\m2$ and $X_i=\pm\AA\det F(K_{(k_1)}\pp J')=\pm\det
FK_{(k_1)}$ if $i\ev n_{k_1}\m2$. (We have $\det FJ'=\AA'=\AA$ in
$\ff/\fs$.) With a good choice of the $\pm$ signs we have
$X_i=Y_i$. If $i=n_{k_2}$ we take $Y_i=\det FK_{k_2}$. For $X_i$
we use Corollary 3.3(ii).  Since $\BB_{k_1+1}=\AA$ is a norm generator
for $M^{\p^{k_1+1}}$ and the Jordan splitting of $M$ starts with
$K_1\pp\cdots\pp K_{k_1}\pp J'\pp K_{k_1+1}=K_{(k_2)}\pp J'$, of rank
$n_{k_2}+1$, we can take $X_i=\BB_{k_1+1}\det (FK_{(k_2)}\pp J')=\det
FK_{(k_2)}$. So again $X_i=Y_i$.

Thus $d[a_{1,i}b_{1,i}]=d[X_iY_i]=\min\{ d(X_iY_i),\a_i,\b_i\}
=\min\{\a_i,\b_i\} =\min\{ A-R_i,A-S_i\}$. But one of $R_i,S_i$ is
$r-2$ and the other $r$, depending on the parity of $i$. So
$d[a_{1,i}b_{1,i}]=A-r$ and we have to prove that $A-r\geq A_i$.

We have $R_{i+1}=S_i=r$ or $r-2$ depending on the parity of $i$. If
$n_{k_1}+1<i<n_{k_2}$ then $A_i\leq\min\{
R_{i+1}-S_i+\a_{i+1},R_{i+1}-S_i+\b_{i-1}\} =\min\{\a_{i+1},\b_{i-1}\}
=\min\{ A-R_{i+1},A-S_{i-1}\} =A-r$. (One of $R_{i+1},S_{i-1}$ is
$r-2$ and the other is $r$.) If $i=n_{k_1}+1$ then $A_i\leq
R_{i+1}-S_i+\a_{i+1}=\a_{i+1}=A-R_{i+1}=A-r$. If $i=n_{k_2}$ then
$A_i\leq R_{i+1}-S_i+\b_{i-1}=\b_{i-1}=A-S_{i-1}=A-r$. \qed

\subsection{Proof of 2.1(iii)}

By the reduction from 5.2 we can restrict to
indices $i\leq n_{k_2}+1$. Suppose now that $R_{i+1}>S_{i-1}$ and
$A_{i-1}+A_i>2e+R_i-S_i$. We want to prove that $W_{i-1}\rep V_i$. 

Note that we cannot have both $\a_{i-1}+\a_i\leq 2e$ and
$\a_i+\a_{i+1}\leq 2e$ since this would imply
$\a_i+A_{i-1}\leq\a_{i-1}+\a_i\leq 2e$ and
$\a_i+d[-a_{1,i+1}b_{1,i-1}]\leq\a_i+\a_{i+1}\leq 2e$, which
contradicts Lemma 2.18(i). This implies that the index $i$ and the
lattice $M$ are in one of the cases (i)-(iv) of the Lemma
3.7. Similarly we cannot have both $\b_{i-2}+\b_{i-1}\leq 2e$ and
$\b_{i-1}+\b_i\leq 2e$ since this would contradict Lemma 2.18(ii). So
$i-1$ and $N$ are in one of the cases (i)-(iv) of Lemma 3.7. Also we
cannot have $R_{i-1}+R_i=R_{i+1}+R_{i+2}$ since this would imply
$R_{i-1}=R_{i+1}$ and $R_i=R_{i+2}$ and so $\a_{i-1}+\a_i\leq 2e$ and
$\a_i+\a_{i+1}\leq 2e$. (See [B3, Lemma 3.3].)

Suppose first that $i\leq n_{k_1}+1$ so $i-1\leq n_{k_1}$. Then
$n_{k-1}+1\leq i-1\leq n_k$ for some $1\leq k\leq k_1$. Since $i-1$
and $N$ are in one of the cases (i)-(iv) of the Lemma 3.7 we have
$i-1=n_{k-1}+1$, $n_k-1$ or $n_k$. We consider the three cases:

If $i-1=n_k-1$, so $i=n_k$, since also $i-1\geq n_{k-1}+1$, $i-1$ and
$N$ are not in the case (i) of Lemma 3.7. Suppose they are in case
(ii) so $S_{i-1}=S_{i+1}$. We have $v'_k=\min\{ u_k,u'\}\leq \min\{
u_k,u\} =v_k$, so $R_i=2r_k-v'_k\geq 2r_k-v_k=S_i$, and
$R_{i+1}>S_{i-1}=S_{i+1}$ so $R_i+R_{i+1}>S_i+S_{i+1}$, which is
impossible. So we are in the case (iii) or (iv) and we can take
$W_{i-1}=FK_{(k)}\top [\BB_k]$. (In the case (iv) we also ask the
condition that $\BB_k\rep FK_k$.)  If $r_k<r'$ then $i=n_k$ and $M$
are in te case (i) of Lemma 3.7. So $V_i=FK_{(k)}$, which implies
$W_{i-1}\rep V_i$. If $r_k=r'$ then $k=k_1$ and $M_{(k)}=K_{(k)}\pp
J'$ so $n_{k-1}+2\leq i=n_k=\dim M_{(k)}-a$. If $a=2$ then $i$ and $M$
cannot be in any of the cases of Lemma 3.7. If $a=1$ then we can only
be in the case (iii) so $V_i=FM_{(k)}\top [\BB'_k]=FK_{(k)}\pp FJ'\top
[\BB'_k]$. But $J'\ap [\AA']$ splits the $\p^{r_k}$ modular component
of $M$ so we can take $\BB'_k=\AA'$. This implies $V_i\ap FK_{(k)}$
and again $W_{i-1}\rep V_i$.

If $i-1=n_k$, so $i=n_k+1$, we use Lemma 3.7(i) and take
$W_{i-1}=FK_{(k)}$. If $i$ and $M$ are in the case (i) of Lemma 3.7
then either $k<k_1$ and $K_{k+1}$ is unary or $k=k_1$ and $J'$ is
unary (i.e. $a=1$). In the first case we take $V_i=FK_{(k+1)}$ while
in the second $V_i=F(K_{(k)}\pp J')$. In both cases $W_{i-1}\rep
V_i$. If $i$ and $M$ are in the case (ii)  or (iv) of Lemma 3.7 then
either $k_1<k$ and the $\p^{r_{k+1}}$ modular component of $M$ is not
unary or $k_1=k$ and $a=\dim J'=2$. In the first case we take
$V_i=FK_{(k)}\pp [\BB'_{k+1}]$, while in the second $V_i=FK_{(k)}\pp
[\BB']$. In both cases $W_{i-1}\rep V_i$. So we are left with the case
when $i$ and $M$ are in the case (iii). This means $R_{i-1}=R_{i+1}$
and we are in one of the following three cases:

a) $r_{k+1}<r'$ and $K_{k+1}$ is binary.

b) $r_{k+1}=r'$ and $K_{k+1}$ and $J$ are both unary.

c) $k=k_1$ and $a=2$.

These are all the cases when $i=n_k+1$ is not of the form $\dim
M_{(l)}$ for some $l$ but $i+1=n_k+2$ is. The Jordan splitting of $M$
starts with $K_1\pp\cdots\pp K_k\pp L$ where $L$ is binary. In the case
a) we have $L=K_{k+1}$, in the case b) we have $k+1=k_1$ and
$L=K_{k+1}\pp J'=K_{k_1}\pp J'$ and in the case c) we have $L=J$. We
will show that, in fact, only case c) can occur and, moreover,
$r_k=r_{k_1}=r'$. First note that $i-1=n_k$ so
$2r_k-v'_k=R_{i-1}=R_{i+1}>S_{i-1}=2r_k-v_k$. Thus $\min\{ u_k,u'\}
=v'_k<v_k=\min\{ u_k,u\}$ which implies $v'_k=u'<u_k$. In the case a)
this implies $v'_{k+1}=u'$ by 5.4. On the other hand
$2r_k-v'_k=R_{i-1}=R_{i+1}=2r_{k+1}-v'_{k+1}$ (we have
$i-1=n_k$ and $i+1=n_k+2=n_{k+1}$). Together with $v'_k=u'=v'_{k+1}$,
this implies $r_k=r_{k+1}$, which is false. In the case b) we have
$k+1=k_1$ and $a=1$ so $v'=u'=r'$. Also $i-1=n_k$ and
$i+1=n_k+2=n_{k+1}+1=n_{k_1}+1$ so
$2r_k-v'_k=R_{i-1}=R_{i+1}=v'=r'$. But $r'>r_k\geq
2r_k-v'_k$. Contradiction. So we are left with the case c). Note that
$v'_k=u'$ implies $v'=u'=v'_k$ by 5.4. We have $i-1=n_k$ and
$i+1=n_k+2=n_{k_1}+2$ so $2r_k-v'_k=R_{i-1}=R_{i+1}=2r'-v'$, which,
together with $v'_k=v'$, implies $r_k=r'$. Thus the $\p^{r_k}=\p^{r'}$
modular component of $M$ is $K_k\pp J'$ and so $M_{(k)}=K_{(k)}\pp J'$
and $\dim M_{(k)}=n_k+2$. Since $i=n_k+1$ and $M$ are in the case
(iii) of Lemma 3.7 we take $V_i=FM_{(k)}\top [\BB'_k]$. But $v'_k=u'$
so by 5.12 we can take $\BB'_k=\AA'$. Since $FJ'\ap [\AA',\(\AA']$ we
get $V_i=FK_{(k)}\pp FJ'\top [\AA']\ap FK_{(k)}\pp [\(\AA']$, which
represents $W_{i-1}=FK_{(k)}$.

If $i-1=n_{k-1}+1$, since we have already treated cases $i-1=n_k-1$
and $i-1=n_k$ we may assume $i-1\leq n_k-2$ so $n_{k-1}+2=i\leq
n_k-1$. If $r_k=r'$ then the $\p^{r_k}$ modular component of $M$ is
$K_k\pp J'$ and so $\dim M_{(k)}=\dim (K_{(k)}\pp J')=n_k+a\geq
n_k+1$. It follows that $\dim M_{(k-1)}+2=n_{k-1}+2=i\leq
n_k-1\leq\dim M_{(k)}-2$ so $i$ and $M$ cannot be in any of the cases
of Lemma 3.7. Thus $r_k<r'$ and so $\dim M_{(k)}=\dim
K_{(k)}=n_k$. Now $n_{k-1}+2=i\leq n_k-1$ so $i$ and $M$ can only be
in the case (iii) of Lemma 3.7 and this happens only when
$i=n_k-1$. This implies $\dim K_{k}=n_k-n_{k-1}=3$. It follows that
$K_k$ is  proper and we have $K_k\ap\la c_1,c_2,c_3\ra$, where
$c_1,c_2,c_3$ are all norm generators for both $M^{\p^{r_k}}$ and
$N^{\p^{r_k}}$. Now $n_{k-1}+1=i-1=n_k-2$ so $i-1$ and $N$ must be in
the case (ii) of Lemma 3.7 and we can take $W_{i-1}=FK_{(k-1)}\pp
[c_1]$. Also $i$ and $M$ are in the case (iii) of the Lemma 3.7 and we
can take $V_i=FK_{(k)}\top [c_3]=FK_{(k-1)}\pp FK_k\top [c_3]\ap
FK_{(k-1)}\pp [c_1,c_2]$ so $W_{i-1}\rep V_i$. (We have $FK_k\ap
[c_1,c_2,c_3]$.)

So we have proved 2.1(iii) for $i\leq n_{k_1}+1$. If $k_1=k_2$ then we
proved 2.1(iii) for $i\leq n_{k_2}+1$ so we are done. Otherwise we
have $k_2=k_1+1$, $a=1$ and $r_{k_2}=r-1=r'+1$ and we have to consider
cases $n_{k_1}+2\leq i\leq n_{k_2}+1$. If $i=n_{k_2}+1$ then the
Jordan splitting of $M$ begins with $K_1\pp\cdots\pp K_{k_1}\pp J'\pp
K_{k_1+1}=K_{(k_2)}\pp J'$ of dimension $n_{k_2}+1=i$ so, by Lemma
3.7(i), we can take $V_i=F(K_{(k_2)}\pp J')$. Also the Jordan splitting
of $N$ starts with $K_1\pp\cdots\pp K_{k_2}=K_{(k_2)}$ of dimension
$n_{k_2}=i-1$ so, by Lemma 3.7(i), we can take $W_{i-1}=FK_{(k_2)}$
and we have $W_{i-1}\rep V_i$. So we may restrict to the case
$n_{k_1}+2\leq i\leq n_{k_2}$.

The sequences $R_{n_{k_1}+1},\ldots,R_{n_{k_2}+1}$ and
$S_{n_{k_1}+1},\ldots,S_{n_{k_2}+1}$ are $r-2,r-1,\ldots,r-1$ and
$r-1,\ldots,r-1,r$ when $K_{k_1+1}$ is proper and they are
$r-2,r,\ldots,r,r-2$ and $r,r-2,\ldots,r-2,r$ otherwise.

If $K_{k_1+1}$ is improper then for any $n_{k_1}+2\leq i\leq n_{k_2}$
we have $S_{j-1}=S_{j+1}$ so $\b_{j-1}+\b_j\leq 2e$. We cannot have
$n_{k_1}+3\leq i\leq n_{k_2}$ since this would imply
$\b_{j-2}+\b_{j-1}\leq 2e$ and $\b_{j-1}+\b_j\leq 2e$. Thus
$i=n_{k_1}+2$, which implies $R_{i+1}=r-2<r=S_{i-1}$ and so (iii) is
vacuous.

Finally, if $K_{k_1+1}$ is proper then for any $n_{k_1}+2\leq i\leq
n_{k_2}$ we have $R_{i+1}=r-1=S_{i-1}$ so (iii) is vacuous. \qed

\subsection{Proof of 2.1(iv)} 

By 5.2 we can restrict to indices $i\leq n_{k_2}+a-1$. Suppose that
$S_i\geq R_{i+2}>S_{i-1}+2e\geq R_{i+1}+2e$. We want to prove that
$W_{i-1}\rep V_i]$. Note that $S_i>R_{i+1}+2e\geq R_i$ and
$S_{i+1}\geq S_i-2e>R_{i+1}$. 

Suppose first that $i\leq n_{k_1}$. We have $n_{k-1}+1\leq i\leq n_k$
for some $k\leq k_1$. Note that $v'_k=\min\{ u_k,u'\}\leq\min\{
u_k,u\} =v_k$. So if $i\ev n_k\m2$ then $R_i=2r_k-v'_k\geq
2r_k-v_k=S_i$, while if $i\ev n_k+1\m2$ then $R_{i+1}=2r_k-v'_k\geq
2r_k-v_k=S_{i+1}$ so we get a contradiction.

If $a=1$ and $k_1=k_2$ then $n_{k_1}=n_{k_2}+a-1$ so we are done. So
we are left with the cases $a=2$, $i=n_{k_2}+a-1=n_{k_1}+1$ and $a=1$,
$k_2=k_1+1$, $n_{k_1}+1\leq i\leq n_{k_2}+a-1=n_{k_2}$. In the first
case $i-1=n_{k_1}$ and $N$ are in the case (i) of Lemma 3.7 and we
take $W_{i-1}=FK_{(k_1)}$. Also the Jordan splitting of $M$ starts
with $K_1\pp\cdots\pp K_{k_1}\pp J'=K_{(k_1)}\pp J'$ and $\dim
(K_{(k_1)}\pp J')=n_{k_1}+2$. Thus $i+1=n_{k_1}+2$ and $M$ are in the
case (i) of Lemma 3.7 so we take $V_{i+1}=F(K_{(k_1)}\pp J')$ and we
have $W_{i-1}\rep V_{i+1}$. In the second case the sequences
$R_{n_{k_1}+1},\ldots,R_{n_{k_2}+1}$ and
$S_{n_{k_1}+1},\ldots,S_{n_{k_2}+1}$ are $r-2,r-1,\ldots,r-1$ and
$r-1,\ldots,r-1,r$ when $K_{k_1+1}$ is proper and they are
$r-2,r,\ldots,r,r-2$ and $r,r-2,\ldots,r-2,r$ otherwise. Since
$n_{k_1}+1\leq i\leq n_{k_2}$ the only possibility to have both
$R_i<S_i$ and $R_{i+1}<S_{i+1}$ is when $K_{k_1+1}$ is unary and
$i=n_{k_1}+1=n_{k_2}$. In this case $K_{k_1+1}$ is proper and the
sequences $R_{n_{k_1}+1},\ldots,R_{n_{k_2}+1}$ and
$S_{n_{k_1}+1},\ldots,S_{n_{k_2}+1}$ are $r-2,r-1$ and $r-1,r$. Now
$i-1=n_{k_1}$ and $N$ are in the case (i) of Lemma 3.7 and we take
$W_{i-1}=FK_{(k_1)}$. Also the Jordan splitting of $M$ starts with
$K_1\pp\cdots\pp K_{k_1}\pp J'\pp K_{k_1+1}$ of dimension
$n_{k_1}+2$. Thus $i+1=n_{k_1}+2$ and $M$ are in the case (i) of Lemma
3.7 and we can take $V_{i+1}=F(K_{(k_1)}\pp J'\pp K_{k_1+1})$ so
$W_{i-1}\rep V_{i+1}$.

\section{More on the case $[M:N]=\p$}

In this section we give a more detailed description of two lattices
$M,N$ satisfying $[M:N]=\p$. Some of these results could be obtained
by methods similar to those of section $\S$5. However, since we want
to keep the use of Jordan splittings to a minimum, we will use the
necessity of 2.1(i)-(iv), which we have already proved. We also use
the relation $R_1+\cdots +R_n+2=\ord volM+2=\ord volN=S_1+\cdots
+S_n$, which follows from $[M:N]=\p$.
\vskip 3mm

First we give some properties of the invariant $\WW (L)\in\bbb$
introduced in $\S$1 (see 1.11).

\bff If $L$ is a lattice of dimension $n$ with $R_i(L)=R_i$,
$\a_i(L)=\a_i$ and $\WW (L)=(x_1,\ldots,x_{2n-2})$ then the sequence
$x_1+x_2,x_2+x_3,\ldots,x_{2n-3}+x_{2n-2}$ is
$R_1+R_2,2R_2+\a_2-\a_1,R_2+R_3,2R_3+\a_3-\a_2,\ldots,R_{n-1}+R_n$. In
particular, $x_1+\cdots+x_{2n-2}=
(x_1+x_2)+(x_3+x_4)+\cdots+(x_{2n-3}+x_{2n-2})=
(R_1+R_2)+(R_2+R_3)+\cdots+(R_{n-1}+R_n)=R_1+2R_2+\cdots+2R_{n-1}+R_n$.
\eff


\bpr If the lattices $M,N$ satisfy 2.1(i) and (ii), in particular, if
$N\rep M$, then $\WW (M)\leq\WW (N)$.
\epr
\pf Let $\WW (M)=(x_1,\ldots,x_{2m-2})$ and $\WW
(N)=(y_1,\ldots,y_{2n-2})$. The condition that $x_j\leq y_j$ or
$x_j+x_{j+1}\leq y_{j-1}+y_j$ for $j=2i-1$, resp. $j=2i$ can be
written as:

a) $R_i+\a_i\leq S_i+\b_i$ or $R_i+R_{i+1}\leq 2S_i+\b_i-\b_{i-1}$.

b) $R_{i+1}-\a_i\leq S_{i+1}-\b_i$ or $2R_{i+1}+\a_{i+1}-\a_i\leq
S_i+S_{i+1}$.

First we reduce to the case $m=n$, which will allow the use of
duality. To do this we use the same technique and notations from the
proof of Lemma 2.20. The fact that $M,N$ satisfy 2.1(i) and (ii)
implies that $M,\( N$ satisfy the same conditions. Now
$R_i(\( N)=R_i(N)=S_i$ for $1\leq i\leq n$ and $\a_i(\(
N)=\a_i(N)=\b_i$ for $1\leq i\leq n-1$ so the first $2n-2$ entries of
$\WW (\( N)\in\bbb_{2m-2}$ are the entries of $\WW (N)$. Thus $\WW
(M)\leq\WW (\( N)$ implies $\WW (M)\leq\WW (N)$.

So we can assume $m=n$. Since $M,N$ satisfy 2.1(i) and (ii) so will do
$N^\*,M^\*$. Since statement b) is equivalent to a) for the lattices
$N^\*,M^\*$ at index $n-i$ we can restrict to proving only statement
a).

Suppose a) is not true. So $R_i+\a_i>S_i+\b_i$ and, if $i>1$,
$R_i+R_{i+1}>2S_i+\b_i-\b_{i-1}\geq S_{i-1}+S_i$. In particular,
$R_i\leq S_i$ and $R_{i+1}>S_{i-1}$. Also
$R_{i+1}+R_{i+2}>S_{i-1}+S_i$ and $R_i+R_{i+1}>S_{i-2}+S_{i-1}$ so
$A_i=\( A_i=\min\{ (R_{i+1}-S_i)/2+e,
R_{i+1}-S_i+d[-a_{1,i+1}b_{1,i-1}], R_{i+1}+R_{i+2}-2S_i+\b_{i-1}\}$
and $A_{i-1}=\( A_{i-1}=\min\{ (R_i-S_{i-1})/2+e,
R_i-S_{i-1}+d[-a_{1,i}b_{1,i-2}], 2R_i-S_{i-2}-S_{i-1}+\a_i\}$. (If
$i=1$ we ignore the terms and relations that are not defined. However
$R_i\leq S_i$ still holds in this case.)

If $A_i=(R_{i+1}-S_i)/2+e$ then $S_i+\b_i\geq
S_i+A_i=(R_{i+1}+S_i)/2+e\geq
(R_i+R_{i+1})/2+e=R_i+(R_{i+1}-R_i)/2+e\geq
R_i+\a_i$. Contradiction. If $A_i=R_{i+1}+R_{i+2}-2S_i+\b_{i-1}$ then
$\b_i\geq A_i=R_i+R_{i+1}-2S_i+\b_{i-1}$ so
$2S_i+\b_i-\b_{i-1}\geq R_i+R_{i+1}$. Contradiction. Thus
$A_i=R_{i+1}-S_i+d[-a_{1,i+1}b_{1,i-1}]$. Since
$R_{i+1}-S_i+d[-a_{i,i+1}]\geq R_i-S_i+\a_i>\b_i\geq A_i$ we have
$d[-a_{i,i+1}]>d[-a_{1,i+1}b_{1,i-1}]$ so
$d[-a_{1,i+1}b_{1,i-1}]=d[a_{1,i-1}b_{1,i-1}]\geq A_{i-1}$. Thus
$A_i\geq R_{i+1}-S_i+A_{i-1}$. In particular, $S_i+\b_i\geq
S_i+A_i\geq R_{i+1}+A_{i-1}$.

If $A_{i-1}=(R_i-S_{i-1})/2+e>(R_i-R_{i+1})/2+e$ then $S_i+\b_i\geq
R_{i+1}+A_{i-1}>(R_i+R_{i+1})/2+e=R_i+(R_{i+1}-R_i)/2+e\geq
R_i+\a_i$. Contradiction. If $A_{i-1}=2R_i-S_{i-2}-S_{i-1}+\a_i$ then
$S_i+\b_i\geq R_{i+1}+A_{i-1}=2R_i+R_{i+1}-S_{i-2}-S_{i-1}+\a_i>
R_i+\a_i$. (We have $R_i+R_{i+1}>S_{i-2}+S_{i-1}$.)
Contradiction. Thus $A_{i-1}=R_i-S_{i-1}+d[-a_{1,i}b_{1,i-2}]$ and we
have $d[a_{1,i}b_{1,i}]\geq A_i\geq R_{i+1}-S_i+A_{i-1}=
R_i+R_{i+1}-S_{i-1}-S_i+d[-a_{1,i}b_{1,i-2}]>d[-a_{1,i}b_{1,i-2}]$.
Hence $d[-a_{1,i}b_{1,i-2}]=d[-b_{i-1,i}]$. So $\b_i\geq A_i\geq
R_i+R_{i+1}-S_{i-1}-S_i+d[-b_{i-1,i}]\geq
R_i+R_{i+1}-2S_i+\b_{i-1}$. Thus $2S_i+\b_i-\b_{i-1}\geq
R_i+R_{i+1}$. Contradiction. \qed

\blm If $M,N$ satisfy 2.1(i) and (ii) and $i\leq\min\{ m-1,n\}$
s.t. $R_j=S_j$ for all $1\leq j\leq i$ then $A_j=\a_j$ for all $1\leq
j\leq i$. 

In particular, $\a_i(L)=A_i(L,L)$ for any lattice $L$ and for any $i$. (We
take $M=N=L$ above.)
\elm
\pf We use induction on $j$. For $j=1$ we have $A_1=\min\{
(R_2-S_1)/2+e,R_2-S_1+d[-a_{1,2}]\} =\min\{
(R_2-R_1)/2+e,R_2-R_1+d[-a_{1,2}]\} =\a_1$. Suppose now $j>1$. We have
$A_j=\min\{ (R_{j+1}-S_j)/2+e, R_{j+1}-S_j+d[-a_{1,j+1}b_{1,j-1}],
R_{j+1}+R_{j+2}-S_{j-1}-S_j+d[a_{1.j+2}b_{1,j-2}]\} =\min\{
(R_{j+1}-R_j)/2+e, R_{j+1}-R_j+d[-a_{1,j+1}b_{1,j-1}],
R_{j+1}+R_{j+2}-R_{j-1}-R_j+d[a_{1.j+2}b_{1,j-2}]\}$. By the
induction hypothesis $d[a_{1,j-1}b_{j-1}]\geq A_{j-1}=\a_{j-1}$ and, if $i\geq
3$, $d[a_{1,j-2}b_{1,j-2}]\geq A_{j-2}=\a_{j-2}$.

We have $(R_{j+1}-R_j)/2+e\geq\a_j$.

We have $R_{j+1}-R_j+d[a_{1,j-1}b_{1,j-1}]\geq R_{j+1}-R_j+\a_{j-1}\geq\a_j$
and $R_{j+1}-R_j+d[-a_{j,j+1}]\geq\a_j$ so
$R_{j+1}-R_j+d[-a_{1,j+1}b_{1,j-1}]\geq\a_j$ by the domination
principle.

We have $R_{j+1}+R_{j+2}-R_{j-1}-R_j+d[a_{1.j-2}b_{1,j-2}]\geq
R_{j+1}+R_{j+2}-R_{j-1}-R_j+\a_{j-2}\geq
R_{j+2}-R_j+\a_j\geq\a_j$. (This inequality also holds if $j=2$,
when $d[a_{1,j-2}b_{1,j-2}]=\j$.) Also
$R_{j+1}+R_{j+2}-R_{j-1}-R_j+d[-a_{j-1,j}]\geq
R_{j+1}+R_{j+2}-2R_j+\a_{j-1}\geq R_{j+1}-R_j+\a_{j-1}\geq\a_j$ and
$R_{j+1}+R_{j+2}-R_{j-1}-R_j+d[a_{j+1,j+2}]\geq
2R_{j+1}-R_{j-1}-R_j+\a_{j+1}\geq R_{j+1}-R_j+\a_{j+1}\geq\a_j$ so
$R_{j+1}+R_{j+2}-R_{j-1}-R_j+d[a_{1.j+2}b_{1,j-2}]\geq\a_j$ by the
domination principle.

Hence $A_j\geq\a_j$. Since also $\a_j\geq d[a_{1,j}b_{1,j}]\geq A_j$
we have $A_j=\a_j$. \qed

We are now able to prove the sufficiency of 2.1(i)-(iv) in a
particular case.

\bpr The main theorem is true for lattices of same rank and volume.
\epr
\pf Suppose that $N\leq M$  and they have the same rank $n$ and same
volume. We will prove that $M\ap N$ by using [B3, Theorem 3.1].

We have $\rr (M)\leq\rr (N)$ and $R_1+\cdots+R_n=\ord volM=\ord
volN=S_1+\cdots+S_n$ so $\rr (M)=\rr (N)$ by Lemma 5.5(iii) Thus
$R_i=S_i~\forall~1\leq i\leq n$.

Denote $\WW (M)=(x_1,\ldots,x_{2n-1})$ and $\WW
(N)=(y_1,\ldots,y_{2n-1})$. We have $\WW (M)\leq\WW (N)$ and
$x_1+\cdots+x_{2n-2}=R_1+2R_2+\cdots+2R_{n-1}+R_n=
S_1+2S_2+\cdots+2S_{n-1}+S_n=y_1+\cdots+y_{2n-2}$. (See 6.1.) By Lemma
5.5(iii) we get $\WW (M)=\WW (N)$ so for any $1\leq i\leq n-1$ we
have $R_i+\a_i=x_{2i-1}=y_{2i-1}=S_i+\b_i=R_i+\b_i$ so $\a_i=\b_i$. So
we already have the conditions (i) and (ii) of [B3, Theorem 3.1].

By Lemma 6.3 we also have $A_i=\a_i$ for all $1\leq i\leq n-1$. It
follows that for any $1\leq i\leq n-1$ we have $d(a_{1,i}b_{1,i})\geq
d[a_{1,i}b_{1,i}]\geq A_i=\a_i$ so we have the (iii) part of [B3,
Theorem 3.1]. Finally, for (iv) if $\a_{i-1}+\a_i>2e$ then, by [B3,
Lemma 3.3], we have $R_{i+1}>R_{i-1}=S_{i-1}$ and also
$A_{i-1}+A_i=\a_{i-1}+\a_i>2e=2e+R_i-S_i$ so by 2.1(iii) we get
$[b_1,\ldots,b_{i-1}]\rep [a_1,\ldots,a_i]$, i.e. we have the (iv) part of
[B3, Theorem 3.1].

Note that we didn't need to use condition (iv) of the main theorem
which is void in the case when $R_i=S_i$ for $1\leq i\leq n$. Indeed,
we cannot have $R_i=S_i>R_{i+1}+2e$. \qed

We want now to study more in detail the case $[M:N]=\p$. Then
$R_i=S_i$ and $\a_i=\b_i=A_i$ will no longer hold but these numbers
will differ from each other by little.

\blm If $M,N$ satisfy 2.1(i) and (ii), in particular if $N\rep M$, and
$i\leq \min\{ m-1,n\}$ is an index s.t. $R_1+\cdots+R_i\ev
S_1+\cdots+S_i+1\m2$ then $R_{i+1}\leq S_i$ or $R_{i+1}+R_{i+2}\leq
S_{i-1}+S_i$. 
\elm
\pf Since $\ord a_{1,i}b_{1,i}=R_1+\cdots+R_i+S_1+\cdots+S_i$ is odd
we have $0=d[a_{1,i}b_{1,i}]\geq A_i=\min\{ (R_{i+1}-S_i)/2+e,
R_{i+1}-S_i+d[-a_{1,i+1}b_{1,i-1}],
R_{i+1}+R_{i+2}-S_{i-1}-S_i+d[a_{1,i+2}b_{1,i-2}]\}$. Thus one of the
following holds: $0\geq (R_{i+1}-S_i)/2+e>(R_{i+1}-S_i)/2$, which
implies $R_{i+1}<S_i$; $0\geq R_{i+1}-S_i+d[-a_{1,i+1}b_{1,i-1}]\geq
R_{i+1}-S_i$, which implies $R_{i+1}\leq S_i$; $0\geq
R_{i+1}+R_{i+2}-S_{i-1}-S_i+d[a_{1,i+2}b_{1,i-2}]\geq
R_{i+1}+R_{i+2}-S_{i-1}-S_i$, which implies $R_{i+1}+R_{i+2}\leq
S_{i-1}+S_i$. \qed 

\blm (i) Let $1\leq i\leq j\leq n$ with $i\ev j\m2$. If $R_i=R_j=R$
then $R_i\ev\ldots\ev R_j\ev R\m2$ and $R_{k+1}-R_k\leq 2e$ for $i\leq
k<j$. Also $R_i+\cdots+R_j\ev R\m2$.

(ii) If $1\leq i<j\leq n$, $j\ev i+1\m2$ and $R_i\geq R_j$ then
$R_i+\cdots+R_j$ is even
\elm
\pf (i) We have $R_i\leq R_{i+2}\leq\ldots\leq R_j$ so $R_i=R_j=R$
implies $R_i=R_{i+2}=\ldots=R_j=R$. In particular, $R_l\ev R\m2$ holds
for $i\leq l\leq j$, $l\ev i\m2$. For $i<l<j$, $l\ev i+1\m2$, at least
one of the numbers $R_{l+1}-R_l=R-R_l$ and $R_l-R_{l-1}=R_l-R$ is
$\leq 0$ and so it is even. So again $R_l\ev R\m2$. 

As a consequence, $R_i+\cdots+R_j\ev (j-i+1)R\ev R\m2$. 

Let $i\leq k<j$. If $k\ev i\m2$ then $R_k=R_{k+2}=R$ so
$R_{k+1}-R_k=-(R_{k+2}-R_{k+1})\leq 2e$. If $k\ev i+1\m2$ then
$R_{k-1}=R_{k+1}=R$ so $R_{k+1}-R_k=-(R_k-R_{k-1})\leq 2e$.

(ii) For any $i\leq l\leq j-1$, $l\ev i\m2$ we have $R_l\geq R_i\geq
R_j\geq R_{l+1}$. (We have $l+1\ev i+1\ev j\m2$.) Thus
$R_{l+1}-R_l\leq 0$ so it is even. So $R_l+R_{l+1}$ is even. By
adding, $R_i+\cdots+R_j$ is even.  \qed

\blm If $[M:N]=\p$ then there are $1\leq s\leq t\leq t'\leq u\leq n$
with $s\ev t\m2$ and $t'\ev u\m2$ and the integers $R,S$ such that:

(i) $R_i=S_i$ for $i<s$ and for $i>u$

(ii) If $s\leq i<t$ then $R_i=R$ and $S_i=R+1$ if $i\ev s\m2$ and
$R_i=S_i+1$ if $i\ev s+1\m2$. Also $R_t=R$.

(iii) If $t'<i\leq u$ then $R_i=S-1$ and $S_i=S$ if $i\ev u\m2$ and
$R_i=S_i+1$ if $i\ev u-1\m2$. Also $S_{t'}=S$.

(iv) One of the following happens:

1. $R+1=S-1=:T$, $t\ev t'\m2$, $R_i=T-1$ and $S_i=T+1$ for $t\leq i\leq
t'$, $i\ev t\m2$ and $R_i=S_i+2$ if $t<i<t'$, $i\ev t+1\m2$.

2. $R+1=S-1=:T$, $t<t'$, $R_i=T$ for $t<i\leq t'$, $S_i=T$ for $t\leq
i<t'$.

3. $t'=t+1$, $S_t=R+1$ and $R_{t'}=S-1$.
\elm

(See [OSU2, Lemma 2.1.2] for another form of this result.)

\pf Let $s$ and $u$ be the smallest resp. the largest indices $i$
s.t. $R_i\neq S_i$. We have $R_i=S_i$ for $i<s$ and for $i>u$ so
(i) holds.

Since $vol N=\p^2vol M$ we have $S_1+\cdots+S_n=R_1+\cdots+R_n+2$. By
Lemma 5.5(iv) we have $S_i\leq R_i+2$ for $1\leq i\leq n$. Assume
first that there is $i_0$ s.t. $S_{i_0}=R_{i_0}+2$. Then $s\leq
i_0\leq u$. Let $S:=S_{i_0}$, $R:=R_{i_0}$ and $T:=R+1=S-1$. By Lemma
5.5(iv) we have $R_1+\cdots+R_{i_0-1}=S_1+\cdots+S_{i_0-1}$ and
$R_{i_0+1}+\cdots+R_n=S_{i_0+1}+\cdots+S_n$. If $s<i_0$ then we can
use Lemma 5.6(i) for $x=\rr (M)$, $y=\rr (N)$, $a=i_0-1$ and $c=s$
and if $u>i_0$ then we can use Lemma 5.6(ii) with $b=i_0+1$ and
$d=u$. We have $R_s<S_s$, $s\ev i_0\m2$,
$R_s=R_{s+2}=\ldots=R_{i_0}=R$ and $R_i+R_{i+1}=S_i+S_{i+1}$ for
$1\leq i<i_0-1$, $i\ev i_0\ev s\m2$. Also $R_u<S_u$, $i_0\ev u\m2$,
$S=S_{i_0}=S_{i_0+2}=\ldots=S_u$ and $R_{i-1}+R_i=S_{i-1}+S_i$ for
$i_0+1<i\leq n$, $i\ev i_0\ev u\m2$. We have $R=R_s<S_s\leq
S_{s+2}\leq\ldots\leq S_{i_0}=S=R+2$. Let $s\leq t\leq i_0$, $t\ev
s\m2$ be minimal s.t. $S_t=R+2$. If $s\leq i<t$, $i\ev s\m2$ then
$S_i=R+1$. If $s<i<t$, $i\ev s+1\m2$ then $R_{i-1}+R_i=S_{i-1}+S_i$,
$R_{i-1}=R$ and $S_{i-1}=R+1$ so $R_i=S_i+1$. Also $R_t=R$ so (ii)
holds. If $t\leq i\leq i_0$, $i\ev t\m2$ we have $S_i=R+2=T+1$ and
$R_i=R=T-1$. If $t<i<i_0$, $i\ev t+1\m2$ then
$R_{i-1}+R_i=S_{i-1}+S_i$, $R_{i-1}=T-1$ and $S_{i-1}=T+1$ so
$R_i=S_i+2$. Similarly $S-2=R=R_{i_0}\leq R_{i_0+2}\leq\ldots\leq
R_u<S_u=S$. Let $i_0\leq t'\leq u$, $t'\ev u\m2$ be maximal
s.t. $R_{t'}=S-2$.  If $t'<i\leq u$, $i\ev u\m2$ then $R_i=S-1$. If
$t'<i<u$, $i\ev u-1\m2$ then $R_i+R_{i+1}=S_i+S_{i+1}$, $R_{i+1}=S-1$
and $S_{i+1}=S$ so $R_i=S_i+1$. Also $S_{t'}=S$ so (iii) holds. If
$i_0\leq i\leq t'$, $i\ev t'\m2$ we have $R_i=S-2=T-1$ and
$S_i=S=T+1$. If $i_0<i<t'$, $i\ev t'-1\m2$ then
$R_i+R_{i+1}=S_i+S_{i+1}$, $R_{i+1}=T-1$ and $S_{i+1}=T+1$ so
$R_i=S_i+2$. In conclusion, $t\ev i_0\ev t'\m2$ and for any $t\leq
i\leq t'$ with $i\ev t\ev t'\m2$ we have $R_i=T-1$ and $S_i=T+1$ (both
when $i\leq i_0$ and when $i\geq i_0$; see above) and for any $t<i<t'$
with $i\ev t+1\ev t'-1\m2$ we  have $R_i=S_i+2$ (both when $i<i_0$ and
when $i>i_0$; see above).Hence (iv) 1. holds.

Suppose now that $S_i\neq R_i+2$ for all $i$ so whenever $S_i>R_i$ we
must have $S_i=R_i+1$. For $0\leq i\leq n$ we define
$f(i):=S_1+\cdots+S_i-R_1-\cdots-R_i$. In particular, $f(i)=f(0)=0$
for $i<s$ and $f(i)=f(n)=2$ for $i\geq u$. Since
$S_1+\cdots+S_n=R_1+\cdots+R_n+2$ we have $0\leq f(i)\leq 2$ by Lemma
5.5(i) and (iv). Let $t-1$ and $t'$ be the largest index
s.t. $f(i)=0$ resp. the smallest index s.t. $f(i)=2$. We have $s\leq
t$ and $u\geq t'$.

We have $f(t-1)=0$ and $f(i)>0$ for $i\geq t$. In particular,
$S_t-R_t=f(t)-f(t-1)=f(t)>0$, which implies $S_t=R_t+1$ so
$f(t)=1$. We also have $t\geq s$ and
$R_1+\cdots+R_{t-1}=S_1+\cdots+S_{t-1}$. If $s<t$ we use Lemma 5.6(i)
with $a=t-1$ and $c=s$. We have $s\ev t\m2$, $R_s<S_s$ so $S_s=R_s+1$,
$R_s=R_{s+2}=\ldots=R_t=:R$ and $R_i+R_{i+1}=S_i+S_{i+1}$ for any
$1\leq i<t$, $i\ev t\ev s\m2$. Now $R+1=R_s+1=S_s\leq
S_{s+2}\leq\ldots\leq S_t=R_t+1=R+1$ so $S_i=R+1$ if $s\leq i\leq t$,
$i\ev s\m2$. For $s<i<t$, $i\ev s+1\m2$ we have
$R_{i-1}+R_i=S_{i-1}+S_i$, $R_{i-1}=R$ and $S_{i-1}=R+1$ so
$R_i=S_i+1$. Therefore (ii) holds. 

Similarly $f(t')=2$ and $f(i)<2$ for $i<t'$. In particular,
$S_{t'}-R_{t'}=f(t')-f(t'-1)=2-f(t'-1)>0$, which implies
$S_{t'}=R_{t'}+1$ so $f(t'-1)=1$. Now $t'\leq u$ and, since $f(t')=2$
we have $R_1+\cdots+R_{t'}+2=S_1+\cdots+S_{t'}$ so
$R_{t'+1}+\cdots+R_n=S_{t'+1}+\cdots+S_n$. If $t'<u$ we use Lemma
5.6(ii) with $b=t'+1$ and $d=u$. We have $t'\ev u\m2$, $R_u<S_u$ so
$S_u=R_u+1$, $S_{t'}=S_{t'+2}=\ldots=S_u=:S$ and
$R_{i-1}+R_i=S_{i-1}+S_i$ for any $t'<i\leq n$, $i\ev t'\ev u\m2$. Now
$S-1=S_{t'}-1=R_{t'}\leq R_{t'+2}\leq\ldots\leq R_u=S_u-1=S-1$ so
$R_i=S-1$ for any $t'\leq i\leq u$, $i\ev u\m2$. For $t'<i<u$, $i\ev
u-1\m2$ we have $R_i+R_{i+1}=S_i+S_{i+1}$, $R_{i+1}=S-1$ and
$S_{i+1}=S$ so $R_i=S_i+1$. Therefore (iii) holds.

We note that $f(i)=0$ for $i<s$. If $s\leq i\leq t$ it is easy to
see that $f(i)=1$ if $i\ev s\m2$ and $f(i)=0$ if $i\ev s+1\m2$. Hence
$f(i)\leq 1$ for $i\leq t$. Since $f(t')=2$ we have $t<t'$. If
$t'=t+1$, since $S_t=R+1$ and $R_{t'}=S-1$ (iv) 3. holds.

Assume now that $t'>t+1$. Since we already have $S_t=R+1$ and
$R_{t'}=S-1$, in order to prove that (iv) 2. holds, we still need to
show that $S-R=2$ and if $T=R+1=S-1$ then $R_i=S_i=T$ for $t<i<t'$. 

We note first that for $t\leq i<t'$ we have $0<f(i)<2$ so
$f(i)=1$. This implies that for $t<i<t'$ we have
$S_i-R_i=f(i)-f(i-1)=1-1=0$ so $R_i=S_i$. Also for $t\leq i<t'$ we
have $f(i)=1$ so $S_1+\cdots+S_i=R_1+\cdots+R_i+1$, which, by Lemma
6.5 implies that $R_{i+1}\leq S_i$ or $R_{i+1}+R_{i+2}\leq
S_{i-1}+S_i$. 

Take first $i=t$ and suppose that $R_{t+1}+R_{t+2}\leq
S_{t-1}+S_t$. If $s<t$ then $R_{t-1}=S_{t-1}+1$ and $S_t=R+1=R_t+1$ so
$R_{t-1}+R_t=S_{t-1}+S_t\geq R_{t+1}+R_{t+2}$, which implies
$R_t=R_{t+2}$. If $s=t$ then $R_{t-1}=S_{t-1}$, which, together with
$R_t+1=S_t$ implies $R_{t-1}+R_t+1= S_{t-1}+S_t\geq
R_{t+1}+R_{t+2}$. Hence $R_{t-1}=R_{t+1}$ or $R_t=R_{t+2}$. But we
cannot have $R_{t-1}=R_{t+1}$ since this would imply
$S_{t-1}=R_{t-1}=R_{t+1}=S_{t+1}$ so, by Lemma 6.6(i), we have both
$S_{t-1}\ev S_t\m2$ and $R_{t-1}\ev R_t\m2$. But this is impossible
since $R_{t-1}=S_{t-1}$ and $R_t+1=S_t$. Thus, in both cases,
$R_t=R_{t+2}$. If $t+2<t'$ then $S_t-1=R_t=R_{t+2}=S_{t+2}$, which is 
impossible. If $t+2=t'$ then $S_t=R_{t}+1=R_{t+2}+1=S=S_{t+2}$. Since
$R_t=R_{t+2}$ and $S_t=S_{t+2}$ we have both $R_{t+1}\ev R_t=R\m2$ and
$S_{t+1}\ev S_t=R+1\m2$ by Lemma 6.6(i). But this is impossible since
$R_{t+1}=S_{t+1}$. Thus $R_{t+1}\leq S_t=R+1$. Suppose that
$R_{t+1}=S_{t+1}\leq R$. Then $S_{t+1}-S_t=R_{t+1}-R-1$ and
$R_{t+1}-R_t=R_{t+1}-R$ are both $\leq 0$ and of opposite
parities. Thus one of them is a negative odd integer, which is
impossible. So $R_{t+1}=S_{t+1}=R+1$.

We use induction to prove that $R_i=S_i=R+1$ for $t+1\leq i<t'$. Case
$i=t+1$ was done above. Suppose now that $i>t+1$. We already have
$R_i=S_i$, as seen above. Note that $R_{i-1}=S_{i-1}=R+1$ by the
induction step and $S_{i-2}=R+1$ either by the induction step, if
$i>t+2$, or because $S_t=R+1$, if $i=t+2$. We have $t\leq i-1<t'$ so
either $R_i\leq S_{i-1}$ or $R_i+R_{i+1}\leq S_{i-2}+S_{i-1}$. In the
first case we have $S_i=R_i\leq S_{i-1}=R+1=S_{i-2}$ so
$S_i=S_{i-2}=R+1$ and we are done. In the second case
$R+1+R_i=R_{i-1}+R_i\leq R_i+R_{i+1}\leq S_{i-2}+S_{i-1}=2R+2$ so
$S_i=R_i\leq R+1=S_{i-2}$. It follows that $R_i=S_i=R+1$.

So we have proved that $R_i=S_i=R+1$ for $t<i<t'$. We are left to
prove that $S-R=2$. We have either $R_{t'}\leq S_{t'-1}$
or $R_{t'}+R_{t'+1}\leq S_{t'-2}+S_{t'-1}$. In the first case $S-1\leq
R+1$ and in the second $R+1+S-1=R_{t'-1}+R_{t'}\leq
R_{t'}+R_{t'+1}\leq S_{t'-2}+S_{t'-1}=R+1+R+1$ so again $S-1\leq
R+1$. Suppose now that $S\leq R+1$. It follows that
$R_{t'}-R_{t'-1}=S-1-(R+1)<0$ and $S_{t'}-S_{t'-1}=S-(R+1)\leq 0$. Now
$R_{t'}-R_{t'-1}$ and $S_{t'}-S_{t'-1}$ are both $\leq 0$ and of
opposite parities so one of them is a negative integer, which is
impossible. Thus $S-R=2$ as claimed. \qed 

\bdf A pair of lattices $M,N$ with $[M:N]=\p$ is said to be of type I,
II or III if it satisfies the condition 1., 2. resp. 3. of Lemma
6.7(iv).
\edf

\bff If the pair $M,N$ is of type I, II or III then so is the pair
$N^\*,M^\*$. Moreover, if the indices $s,t,t',u$ corresponding to
$N^\*,M^\*$ are $s^\*,t^\*,t'^\*,u^\*$ then $s^\* =n+1-u$, $t^\*
=n+1-t'$, $t'^\* =n+1-t$ and $u^\* =n+1-s$.

In the proof of the next lemma we will use the duality  in order to
shorten the proof. To do this we refer to the remark above and to the
results from 2.4.
\eff

{\bf Remark} With the notations from Lemma 6.7 the indices $0\leq
i\leq n$ s.t. $S_1+\cdots+S_i=R_1+\cdots+R_i+1$ are $s\leq i<t$ with
$i\ev s\m2$,  $t'<i\leq u$ with $i\ev u-1\m2$ and, if $M,N$ are of the
type II or III, $t\leq i<t'$.

\blm Suppose that $[M:N]=\p$. With the notations of Lemma 6.7 we have:

(i) If $S_1+\cdots+S_i=R_1+\cdots+R_i+1$ then $\a_i,\b_i\leq 1$ and
$A_i=0$, unless $\a_i=\b_i=0$. The case $\a_i=\b_i=0$ happens only
when $M,N$ are of type III, $i=t$ and
$S-R-1=R_{t+1}-R_t=S_{t+1}-S_t=-2e$.

(ii) If $S_1+\cdots+S_i\neq R_1+\cdots+R_i+1$ then $|\a_i-\b_i|\leq 2$
and $A_i=\min\{\a_i,\b_i\}$.

(iii) If $1\leq i<s$ then
$A_i=\a_i=\min\{\b_i,R_{s-1}-R_i+\a_{s-1}\}$.

(iv) If $u\leq i<n$ then
$A_i=\b_i=\min\{\a_i,S_{i+1}-S_{u+1}+\b_u\}$.

(v) If $M,N$ are of type I then $R_i+\a_i=S_i+\b_i$ for any $t\leq
i<t'$.
\elm
\pf Let $\WW (M)=(x_1,\ldots,x_{2n-2})$ and $\WW
(N)=(y_1,\ldots,y_{2n-2})$.

Note that, in part (ii), since $\a_i,\b_i\geq A_i$, in order to prove
that $A_i=\min\{\a_i,\b_i\}$ it is enough to show that $A_i\geq\a_i$
or $A_i\geq\b_i$. In particular, if $1\leq i<s$ we have by Lemma
6.3 $A_i=\a_i$. By duality at index $n-i$ we get $A_i=\b_i$ for $u\leq
i<n$. Thus $A_i=\min\{\a_i,\b_i\}$ holds for both $1\leq i<s$ and
$u\leq i<n$.

Also, in the part (i), since $R_1+\cdots+R_i+1=S_1+\cdots+S_i$, we
have that $\ord a_{1,i}b_{1,i}$ is odd so $0=d[a_{1,i}b_{1,i}]\geq
A_i$. Thus in order to prove $A_i=0$ we only need $A_i\geq 0$.

We prove first (iii). We note that for $1\leq i\leq s-2$ we have
$x_{2i-1}+x_{2i}=R_i+R_{i+1}=S_i+S_{i+1}=y_{2i-1}+y_{2i}$. By adding
we get $x_1+\cdots+x_{2s-4}=y_1+\cdots+y_{2s-4}$. By Lemma 5.6(i) we
have $x_j=\min\{ y_j,x_{2s-3}\}$ for any $1\leq j\leq 2s-3$, $j$
odd. If we take $j=2i-1$ with $1\leq i<s$ we have $R_i=S_i$ and
$x_{2i-1}=\min\{ y_{2i-1},x_{2s-3}\}$ i.e. $R_i+\a_i=\min\{ 
S_i+\b_i,R_{s-1}+\a_{s-1}\}=\min\{ R_i+\b_i,R_{s-1}+\a_{s-1}\}$. Thus
$\a_i=\min\{\b_i,R_{s-1}-R_i+\a_{s-1}\}$.

(iv) follows from (iii) by duality at index $n-i$.

Next we prove that the $\a_i,\b_i\leq 1$ part of (i) holds for $t\leq
i<t'$ if $M,N$ are of type II or III. It is enough to prove $\a_i\leq
1$ since $\b_i\leq 1$ will follow by duality at index $n-i$. If we
have type II then for any $t\leq i<t'$ we have $\a_i\leq
R_{i+1}-R_t+d(-a_{t,t+1})=T-(T-1)+0=1$ ($\ord
a_{t,t+1}=R_t+R_{t+1}=2T-1$ so $d(-a_{t,t+1})=0$). If we have type III
then $t'=t+1$ so we only have to prove that $\a_t\leq 1$. If $s=t$
then $S_{t-1}=R_{t-1}$ while if $s<t$ then $S_{t-1}=R_{t-1}-1$. In
both cases $S_{t-1}\leq R_{t-1}\leq R_{t+1}$. Similarly
$S_{t+2}=R_{t+2}$ if $t+1=t'=u$ and $S_{t+2}=R_{t+2}-1$ otherwise. So
$S_t\leq S_{t+2}\leq R_{t+2}$. If $S_t=S_{t+2}=R_{t+2}$ then
$R_{t+1}+1=S_{t+1}\ev S_{t+2}=R_{t+2}\m2$ so $\ord a_{t+1,t+2}$ is odd
which implies that $\a_t\leq R_{t+2}-R_t+d(-a_{t+1,t+2})=1$. (We have
$R_{t+2}=S_t=R_t+1$ and $d(-a_{t+1,t+2})=0$.) So we may assume that
$R_{t+2}>S_t$. Together with $R_{t+1}\geq S_{t-1}$, this implies
$R_{t+1}+R_{t+2}>S_{t-1}+S_t$ so $A_t=\( A_t=\min\{ (R_{t+1}-S_t)/2+e,
R_{t+1}-S_t+d[-a_{1,t+1}b_{1,t-1}]\}$. (Same happens if $t=1$ or $n-1$, when
$S_{t-1}$ resp. $R_{t+2}$ is not defined.) Since
$S_1+\cdots+S_t=R_1+\cdots+R_t+1$ we have that $\ord a_{1,t}b_{1,t}$ is odd so
$0=d[a_{1,t}b_{1,t}]\geq A_t$. If $0\geq A_t=(R_{t+1}-S_t)/2+e$ we get
$R_{t+1}-R_t=R_{t+1}-S_t+1\leq 1-2e$ so $R_{t+1}-R_t=-2e$ and we have
$\a_t=0$. Otherwise $0\geq A_t=R_{t+1}-S_t+d[-a_{1,t+1}b_{1,t-1}]$ so
$d[-a_{1,t+1}b_{1,t-1}]\leq S_t-R_{t+1}=(R+1)-(S-1)=R-S+2$. If we
prove that $d[-a_{t,t+1}]\leq R-S+2$ then $\a_t\leq
R_{t+1}-R_t+d[-a_{t,t+1}]\leq S-1-R+R-S+2=1$ and we are done. Since
$d[-a_{1,i+1}b_{1,i-1}]\leq R-S+2$ it is enough to prove that
$d[a_{1,i-1}b_{1,i-1}]>R-S+2$. To do this, by domination principle, it
is enough to prove that $d[a_{1,s-1}b_{1,s-1}]>R-S+2$ and
$d[-a_{i,i+1}],d[-b_{i,i+1}]>R-S+2$ for $s\leq i<t$, $i\ev s\m2$. If
$s=1$ then $d[a_{1,s-1}b_{1,s-1}]=\j >R-S+2$. If $s>1$ then
$d[a_{1,s-1}b_{1,s-1}]\geq A_{s-1}=\a_{s-1}$ so if
$d[a_{1,s-1}b_{1,s-1}]\leq R-S+2$ then $\a_t\leq
R_{t+1}-R_s+\a_{s-1}\leq S-1-R+R-S+2=1$ and we are done. So we may
assume $d[a_{1,s-1}b_{1,s-1}]>R-S+2$. Let now $s\leq i<t$, $i\ev
s\m2$. If $d[-a_{i,i+1}]\leq R-S+2$ then $\a_t\leq
R_{t+1}-R_i+d[-a_{i,i+1}]=S-1-R+d[-a_{i,i+1}]\leq 1$ and we are
done. So we may assume $d[-a_{i,i+1}]>R-S+2$. Also $d[-b_{i,i+1}]\geq
S_i-S_{i+1}+\b_i\geq S_i-S_{i+1}=(R+1)-(R_{i+1}-1)\geq
R-R_{t+1}+2=R-S+3$.

Now we prove that $\a_i,\b_i\leq 1$ for indices $i$
s.t. $R_1+\cdots+R_i+1=S_1+\cdots+S_i$, other than $t\leq i<t'$. Take
first $s\leq i<t$, $i\ev s\m2$. If $M,N$ are of type II or III we have
proved that $\a_t,\b_t\leq 1$. But $R_i+\a_i\leq R_t+\a_t$,
$S_i+\b_i\leq S_t+\b_t$, $R_i=R_t=R$ and $S_i=S_t=R+1$. Thus
$\a_i\leq\a_t\leq 1$ and $\b_i\leq\b_t\leq 1$. If we have type I and
$s<t$ then $R_{t-2}=R_t=T-1$ so $R_{t-1}\ev R_t=T-1\m2$. Hence
$S_{t-1}=R_{t-1}-1\ev T\m2$. Since also $S_t=T+1$ we get that $\ord
b_{t-1,t}=S_{t-1}+S_t$ is odd so $d(-b_{t-1,t})=0$. Hence
$\b_i\leq S_t-S_i+d(-b_{t-1,t})=T+1-T+0=1$.

So we are left with proving that $\a_i\leq 1$ for $s\leq i<t$, $i\ev
s\m2$, in the case of type I. We have $R_{s+1}\leq
R_{s+3}\leq\ldots\leq R_{t-1}$. Let $s\leq l< t$, $l\ev s\m2$ be
minimal s.t. $R_{l+1}=R_{t-1}$. We have
$R_{l+1}=R_{l+3}=\ldots=R_{t-1}$ and, if $s<l$, $R_{l-1}<R_{l+1}$. We
prove first that $\a_l\leq 1$. We have
$R_1+\cdots+R_l+1=S_1+\cdots+S_l$ so $\ord a_{1,l}b_{1,l}$ is odd,
which implies $0=d[a_{1,l}b_{1,l}]\geq A_l=\min\{ (R_{l+1}-S_l)/2+e,
R_{l+1}-S_l+d[-a_{1,l+1}b_{1,l-1}],
R_{l+1}+R_{l+2}-S_{l-1}-S_l+d[a_{1,l-2}b_{1,l+2}]\}$. If $0\geq
(R_{l+1}-S_l)/2+e$ then $-2e\geq R_{l+1}-S_l=S_{l+1}-S_l+1$ so
$S_{s+1}-S_s<-2e$, which is impossible. Suppose $0\geq
R_{l+1}+R_{l+2}-S_{l-1}-S_l+d[a_{1,l-2}b_{1,l+2}]\geq
R_{l+1}+R_{l+2}-S_{l-1}-S_l$. We have $S_l=T$ and $R_{l+2}=T-1$ so
$0\geq R_{l+1}-S_{l-1}-1$ i.e. $R_{l+1}\leq S_{l-1}+1$. Since also
$R_{l+1}=S_{l+1}+1\geq S_{l-1}+1$ we have $R_{l+1}=S_{l-1}+1$. If
$s<l$ then $R_{l+1}>R_{l-1}=S_{l-1}+1$. Contradiction. If $s=l$ then
$R_{l+1}=S_{l-1}+1=R_{l-1}+1$. Since $R_l=R_{l+2}=T-1$ we also have
$R_l\ev R_{l+1}=R_{l-1}+1\m2$. Thus $\ord a_{l-1,l}$ is odd and so
$\a_l\leq R_{l+1}-R_{l-1}+d(-a_{l-1,l})=1$. (We have
$R_{l+1}=R_{l-1}+1$ and $d(-a_{l-1,l})=0$.) So we are left with the
case $0\geq R_{l+1}-S_l+d[-a_{1,l+1}b_{1,l-1}]$ so
$d[-a_{1,l+1}b_{1,l-1}]\leq S_l-R_{l+1}=T-R_{l+1}$. If
$d[-a_{l,l+1}]\leq T-R_{l+1}$ then $\a_l\leq
R_{l+1}-R_l+d[-a_{l,l+1}]\leq R_{l+1}-(T-1)+T-R_{l+1}=1$. Since
$d[-a_{1,l+1}b_{1,l-1}]\leq T-R_{l+1}$ it is enough to prove that
$d[a_{1,l-1}b_{1,l-1}]>T-R_{l+1}$. To do this, by domination
principle, it is enough to prove that
$d[a_{1,s-1}b_{1,s-1}]>T-R_{l+1}$ and
$d[-a_{j,j+1}],d[-b_{j,j+1}]>T-R_{l+1}$ for any $s\leq j<l$, $j\ev
s\m2$. We have $d[a_{1,s-1}b_{1,s-1}]\geq A_{s-1}=\a_{s-1}$ so if
$d[a_{1,s-1}b_{1,s-1}]\leq T-R_{l+1}$ then $\a_l\leq
R_{l+1}-R_s+\a_{s-1}\leq R_{l+1}-(T-1)+T-R_{l+1}=1$ and we are
done. So we may assume that $d[a_{1,s-1}b_{1,s-1}]>T-R_{l+1}$. Let
$s\leq j<l$, $j\ev s\m2$. If $d[-a_{j,j+1}]\leq T-R_{l+1}$ then by
Remark 1.1 $\a_l\leq R_{l+1}-R_j+d[-a_{j,j+1}]\leq
R_{l+1}-(T-1)+T-R_{l+1}=1$ and we are done. So we may assume that
$d[-a_{j,j+1}]>T-R_{l+1}$. Finally, we have $d[-b_{j,j+1}]\geq
S_j-S_{j+1}+\b_j\geq S_j-S_{j+1}=T-(R_{j+1}-1)>T-R_{l+1}+1$. (We have
$R_{j+1}\leq R_{l-1}<R_{l+1}$.) Suppose now that $s\leq i<t$ with
$i\ev s\m2$ is arbitrary. If $i\leq l$ then $R_i+\a_i\leq R_l+\a_l$
and $R_i=R_l=T-1$ so $\a_i\leq\a_l\leq 1$. If $l\leq i$ then
$R_{l+1}=R_{i+1}$ so $-R_{l+1}+\a_l\geq -R_{i+1}+\a_i$ implies
$\a_i\leq\a_l\leq 1$. 

The statement $\a_i,\b_i\leq 1$ for indices $t'<i<u$, $i\ev u-1\m2$
follows by duality at index $n-i$ from the case $s\leq i<t'$, $i\ev
s\m2$.

Note that if $s\leq i<t$ with $i\ev s\m2$ then $R_i=T-1=S_i-1$ and
$R_{i+1}=S_{i+1}+1$ so $R_{i+1}-R_i=S_{i+1}-S_i+2\geq 2-2e$. Hence
$\a_i\neq 0$ so $\a_i=1$. Dually, if $t'<i<u$, $i\ev u-1\m2$ then
$\b_i=1$.

We now prove (v). We may assume that $t<t'$ since otherwise (v) is
vacuous. First we prove that $x_{2t-1}\leq y_{2t-1}$ i.e.
$R_t+\a_t\leq S_t+\b_t$. If $t=1$ then $x_1\leq y_1$ follows from
$\WW (M)\leq\WW (N)$ so we may assume that $t>1$. Suppose first
that $s<t$. We have $\a_{t-2}=1$ and $R_{t-2}=R_t=T-1$ so by [B3,
Corollary 2.3(i)] we have
$x_{2t-2}=R_t-\a_{t-1}=R_{t-1}-\a_{t-2}=R_{t-1}-1=S_{t-1}$. Now
$R_{t-2}=R_t$ also implies $R_{t-1}\ev R_{t-2}=T-1\m2$ so
$S_{t-1}=R_{t-1}-1\ev T=S_t-1\m2$. Thus $S_t-S_{t-1}$ is odd. Also
$S_{t-1}\geq S_{t-2}-2e=T-2e$. Hence $S_t-S_{t-1}=T+1-S_{t-1}\leq
2e+1$. If $S_t-S_{t-1}\leq 2e$, since also $S_t-S_{t-1}$ is odd, we
get $\b_t=S_t-S_{t-1}$ so $y_{2t-2}=S_t-\b_{t-1}=S_{t-1}=x_{2t-2}$. If
$S_t-S_{t-1}=2e+1$ then $\b_{t-1}=(2e+1)/2+e=2e+1/2$ so
$y_{2t-2}=S_t-\b_{t-1}=S_{t-1}+2e+1-(2e+1/2)=S_{t-1}+1/2=x_{2t-2}+1/2$.
So in both cases $y_{2t-2}\leq x_{2t-2}+1/2$. Together with
$x_{2t-2}+x_{2t-1}\leq y_{2t-2}+y_{2t-1}$, this implies
$x_{2t-1}\leq y_{2t-1}+1/2$. But $R_t=R_{t+2}=T-1$ and
$S_t=S_{t+2}=T+1$ so $R_t-R_{t-1}$ and $S_t-S_{t-1}$ are even.
This implies that $\a_t,\b_t$ are integers and so are
$x_{2t-1}=R_t+\a_t$ and $y_{2t-1}=S_t+\b_t$. Therefore
$x_{2t-1}\leq y_{2t-1}+1/2$ implies $x_{2t-1}\leq y_{2t-1}$.
Suppose now that $s=t$. We have $R_t=T-1$ and $S_t=T+1$ so
$R_t+\a_t\leq S_t+\b_t$ means $\b_t\geq\a_t-2$. But $\b_t\geq A_t$ so
it is enought to prove that $A_t\geq\a_t-2$,
i.e. $A_s\geq\a_s-2$. Note that $R_{s+2}=T-1<T+1=S_s$ so $A_s=A'_s$ by
Lemma 2.14. Also $R_{s+1}=S_{s+1}+2>S_{s-1}$ so $A'_s=\min\{
R_{s+1}-S_s+d[-a_{1,s+1}b_{1,s-1}],
R_{s+1}+R_{s+2}-S_{s-1}-S_s+d[-a_{1,s}b_{1,s-2}]\} =\min\{
R_{s+1}-R_s+d[-a_{1,s+1}b_{1,s-1}]-2,
R_{s+1}-R_{s-1}+d[-a_{1,s}b_{1,s-2}]-2\}$. (We have
$S_s=T+1=R_s+2=R_{s+2}+2$ and $S_{s-1}=R_{s-1}$.) Now
$R_{s+1}-R_s+d[a_{1,s-1}b_{1,s-1}]-2\geq
R_{s+1}-R_s+A_{s-1}-2=R_{s+1}-R_s+\a_{s-1}-2\geq\a_s-2$ and
$R_{s+1}-R_s+d[-a_{s,s+1}]-2\geq\a_s-2$ so
$R_{s+1}-R_s+d[-a_{1,s+1}b_{1,s-1}]-2\geq\a_s-2$. Also
$R_{s+1}-R_{s-1}+d[a_{1,s-2}b_{1,s-2}]-2\geq
R_{s+1}-R_{s-1}+A_{s-2}-2=R_{s+1}-R_{s-1}+\a_{s-2}-2\geq\a_s-2$
and $R_{s+1}-R_{s-1}+d[-a_{s-1,s}]-2\geq
R_{s+1}-R_s+\a_{s-1}-2\geq\a_s-2$ so
$R_{s+1}-R_{s-1}+d[-a_{1,s}b_{1,s-2}]-2\geq\a_s-2$. (Note that
$R_{s+1}-R_{s-1}+d[a_{1,s-2}b_{1,s-2}]-2\geq\a_s-2$ also holds if
$s=2$, when $d[a_{1,s-2}b_{1,s-2}]=\j$.) Hence $A_s=A'_s\geq\a_s-2$,
as desired. So we have proved that $R_t+\a_t\leq S_t+\b_t$,
i.e. $x_{2t-1}\leq y_{2t-1}$. By duality we get $-S_{t'}+\b_{t'-1}\leq
-R_{t'}+\a_{t'-1}$, i.e. $x_{2t'-2}\leq y_{2t'-2}$. Since $\WW
(M)\leq\WW (N)$, $x_{2t-1}\leq y_{2t-1}$ and $x_{2t'-2}\leq
y_{2t'-2}$ we get $(x_{2t-1},\ldots,x_{2t'-2})\leq
(y_{2t-1},\ldots,y_{2t'-2})$. But for $t\leq i\leq t'$ we have
$R_i=T-1=S_i-2$ if $i\ev t\m2$ and $R_i=S_i+2$ otherwise. So for
any $t\leq i<t'$ we have
$x_{2i-1}+x_{2i}=R_i+R_{i+1}=S_i+S_{i+1}=y_{2i-1}+y_{2i}$. By
adding we get $x_{2t-1}+\cdots+x_{2t'-2}=y_{2t-1}+\cdots+y_{2t'-2}$,
which by Lemma 5.5(iii) implies that
$(x_{2t-1},\ldots,x_{2t'-2})=(y_{2t-1},\ldots,y_{2t'-2})$. In
particular, for $t\leq i<t'$ we have
$R_i+\a_i=x_{2i-1}=y_{2i-1}=S_i+\b_i$ and we are done.

Now prove now (ii) for indices $s\leq i<u$. Take first $s<i<t$,
$i\ev s+1\m2$ (when $s<t$). As noted in the beginning of the
proof, in order to prove $A_i=\min\{\a_i,\b_i\}$, it is enough to
show that $A_i\geq\a_i$, which will imply $A_i=\a_i$. We have
$\a_{i-1}=1$ by (i) and $R_{i-1}=R_{i+1}=R$. By [B3, Corollary 2.3(i)]
$R_{i+1}-\a_i=R_i-\a_{i-1}$, which implies
$\a_i=R_{i+1}-R_i+\a_{i-1}=R_{i+1}-R_i+1$. So we have to prove
that $A_i\geq\a_i=R_{i+1}-R_i+1=R_{i+1}-S_i$. We use induction on
$i$. Now $R_{i+1}=R<R+1=S_{i-1}$ so $A_i=A'_i$ by Lemma 2.14. If
$A_i'=R_{i+1}-S_i+d[-a_{1,i+1}b_{1,i-1}]$ then $A'_i\geq
R_{i+1}-S_i$ sowe are done. So we may assume that $2\leq i\leq
n-2$ and $A'_i=R_{i+1}+R_{i+2}-S_{i-1}-S_i+d[a_{1,i-2}b_{1,i+2}]=
R_{i+2}-R_i+d[a_{1,i+2}b_{1,i-2}]$. (We have
$S_{i-1}=R+1=R_{i+1}+1$ and $S_i=R_i-1$.) We have
$R_{i+2}-R_i+d[-a_{i-1,i}]\geq d[-a_{i-1,i}]\geq
R_{i-1}-R_i+\a_{i-1}=R_{i+1}-R_i+1=\a_i$. Also
$R_{i+2}-R_i+d[-a_{i+1,i+2}]\geq R_{i+1}-R_i+\a_{i+1}\geq\a_i$. If
$i>2$ then $A_{i-2}=\a_{i-2}$. (If $i>s+1$ it follows from the
induction hypothesis; if $i=s+1$ it follows from the case $i<s$,
already discussed.) Thus $R_{i+2}-R_i+d[a_{1,i-2}b_{1,i-2}]\geq
d[a_{1,i-2}b_{1,i-2}]\geq A_{i-2}=\a_{i-2}\geq
R_{i-1}-R_{i+1}+\a_i=\a_i$. Same happens if $i=2$, when
$d[a_{1,s-2}b_{1,s-2}]=\j$. By domination principle
$A'_i=R_{i+2}-R_i+d[a_{1,i+2}b_{1,i-2}]\geq\a_i$, as desired. We
prove now that $|\a_i-\b_i|\leq 2$. We already have $\b_i\geq
A_i=\a_i$ so we still need $\b_i\leq\a_i+2$. If $i<t-1$ or $i=t-1$
and $M,N$ are of the type II or III then $S_{i-1}=S_{i+1}=R+1$ and
$\b_{i-1}\leq 1$ so $S_{i+1}-\b_i=S_i-\b_{i-1}\geq S_i-1$, which
implies $\b_i\leq S_{i+1}-S_i+1=R_{i+1}-R_i+3=\a_i+2$. (We have
$S_{i+1}=R+1=R_{i+1}+1$ and $S_i=R_i-1$.) If $i=t-1$ and we have
type I then $R_{i-1}=R_{i+1}=T-1$ so $R_i\ev R_{i+1}=T-1\m2$. Thus
$S_i=R_i+1\ev T=S_{i+1}-1\m2$. Thus $S_{i+1}-S_i$ is odd which
implies $\b_i\leq S_{i+1}-S_i=R_{i+1}-R_i+3=\a_i+2$. (We have
$S_{i+1}=T+1=R_{i+1}+2$ and $S_i=R_i-1$.)

By duality at index $n-i$ (ii) also holds for $t'<i<u$ with $i\ev
u\m2$. Moreover at these indices we have
$A_i=\b_i=S_{i+1}-S_i+1=R_{i+1}-S_i$. Also $\a_i\geq A_i=\b_i$ and
$\a_i\leq\b_i+2$.

We are left to prove (ii) for indices $t\leq i<t'$ when $M,N$ are of
type I. For such indices we have $R_i+\a_i=S_i+\b_i$. If $i\ev
t\m2$ then $S_i=T+1=R_i+2$ so $\b_i=\a_i-2$ and if $i\ev t+1\m2$ then
$S_i=R_i-2$ so $\b_i=\a_i+2$. In both cases $|\a_i-\b_i|=2$. We use
the induction on $i$ to prove that $A_i\geq\b_i=\a_i-2$ if $i\ev t\m2$
and $A_i\geq\a_i$ if $i\ev t+1\m2$. This will imply $A_i=\b_i$
resp. $A_i=\a_i$.

Take first the case $i\ev t\m2$. We have to prove that
$A_i\geq\b_i=\a_i-2$. We have $R_{i+1}=S_{i+1}+1>S_{i-1}$ so
$A_i=\min\{ (R_{i+1}-S_i)/2+e, R_{i+1}-S_i+d[-a_{1,i+1}b_{1,i-1}],
R_{i+1}+R_{i+2}-S_{i-1}-S_i+d[-a_{1,i}b_{1,i-2}]\}=\min\{
(R_{i+1}-R_i)/2+e-1, R_{i+1}-R_i+d[-a_{1,i+1}b_{1,i-1}]-2,
R_{i+1}-S_{i-1}+d[-a_{1,i}b_{1,i-2}]-2\}$. (We have
$S_i=T+1=R_i+2=R_{i+2}+2$.) If $A_i=(R_{i+1}-R_i)/2+e-1$ then
$A_i\geq\a_i-1$ so we are done. If $i>1$ then $A_{i-1}=\a_{i-1}$ both
if $i=t$, by the case $i<t$ already studied, and if $i>t$, by the
induction hypothesis, since $i-1\ev t+1\m2$. Thus
$R_{i+1}-R_i+d[a_{1,i-1}b_{1,i-1}]-2\geq
R_{i+1}-R_i+A_{i-1}-2=R_{i+1}-R_i+\a_{i-1}-2\geq\a_i-2$. (This also
holds if $i=1$, when $d[a_{1,i-1}b_{1,i-1}]=\j$.) Since also
$R_{i+1}-R_i+d[-a_{i,i+1}]-2\geq\a_i-2$ we have
$R_{i+1}-R_i+d[-a_{1,i+1}b_{1,i-1}]-2\geq\a_i-2$. Finally,
$R_{i+1}-S_{i-1}+d[-a_{1,i}b_{1,i-2}]-2\geq
R_{i+1}-R_i+A_{i-1}-2=R_{i+1}-R_i+\a_{i-1}-2\geq\a_i-2$.

Take now the case $i\ev t+1\m2$. We have to prove that
$A_i\geq\a_i$. Since $R_{i+1}=T-1<T+1=S_{i-1}$ we have by Lemma
2.14 $A_i=A'_i=\min\{ R_{i+1}-S_i+d[-a_{1,i+1}b_{1,i-1}],
R_{i+1}+R_{i+2}-S_{i-1}-S_i+d[a_{1,i+2}b_{1,i-2}]\}=\min\{
R_{i+1}-R_i+d[-a_{1,i+1}b_{1,i-1}]+2,
R_{i+2}-R_i+d[a_{1,i+2}b_{1,i-2}]\}$. (We have $S_i=R_i-2$ and
$S_{i-1}=T+1=R_{i+1}+2$.) Now $A_{i-1}=\b_{i-1}=\a_{i-1}-2$ by the
induction hypothesis ($i-1\ev t\m2$). So
$R_{i+1}-R_i+d[a_{1,i-1}b_{1,i-1}]+2\geq
R_{i+1}-R_i+A_{i-1}+2=R_{i+1}-R_i+\a_{i-1}\geq\a_i$. Also
$R_{i+1}-R_i+d[-a_{i,i+1}]+2\geq\a_i+2>\a_i$ so
$R_{i+1}-R_i+d[-a_{1,i+1}b_{1,i-1}]+2\geq\a_i$. We have
$R_{i+2}-R_i+d[-a_{1,i}b_{1,i-2}]\geq d[-a_{1,i}b_{1,i-2}]\geq
S_{i-1}-R_i+A_{i-1}=R_{i+1}-R_i+\a_{i-1}\geq\a_i$. (We have
$S_{i-1}=T+1=R_{i+1}+2$ and $A_{i-1}=\a_{i-1}-2$.) Since also
$R_{i+2}-R_i+d[-a_{i+1,i+2}]\geq\a_i$ (see Remark 1.1) we have
$R_{i+2}-R_i+d[a_{1,i+2}b_{1,i-2}]\geq\a_i$. So $A_i\geq\a_i$ and we
are done.

To complete the proof of (ii) we need to show that $|\a_i-\b_i|\leq 2$
for $1\leq i<s$ and for $u\leq i<n$. If $1\leq i<s$ then
$\a_i\leq\b_i$ by (iii) so we still have to prove that
$\b_i\leq\a_i+2$. Take first $i=s-1$. If $s<t$ or $s=t$ and we have
type II or III then $S_s=R_s+1$ so $R_1+\cdots+R_s+1=S_1+\cdots+S_s$ which
implies $\b_s\leq 1$ so $\b_{s-1}\leq S_s-S_{s-1}+\b_s\leq
S_s-S_{s-1}+1=R_s-R_{s-1}+2$. If $R_s-R_{s-1}\leq 2e$ then
$\a_{s-1}\geq R_s-R_{s-1}$ so $\b_{s-1}\leq\a_{s-1}+2$. If
$R_s-R_{s-1}>2e$ then $S_s-S_{s-1}=R_s-R_{s-1}+1>2e$ so
$\b_{s-1}=(S_s-S_{s-1})/2+e=(R_s-R_{s-1})/2+e+1/2=\a_{s-1}+1/2$. So we
are left with the case when $s=t$ and we have type I. We have
$S_s=T+1=R_s+2$. Note that if $\a_{s-1}=(R_s-R_{s-1})/2+e$ then
$\b_{s-1}\leq (S_s-S_{s-1})/2+e=(R_s-R_{s-1})/2+e+1=\a_{s-1}+1$ so we
are done. Hence we may assume that
$(R_s-R_{s-1})/2+e>\a_{s-1}=R_s-R_{s-1}+d[-a_{s-1,s}]$. Assume first
that $3\leq s\leq n-1$. We have $\a_s\geq A_s=\b_s$. This happens in
all cases $s=t<t'$, $s=t'<u$ and $s=u$; see the proof of
$A_i=\min\{\a_i,\b_i\}$ in these cases. We have $\a_{s-1}\geq
R_s-R_{s+1}+\a_s\geq R_s-R_{s+1}+\b_s\geq
R_s-R_{s+1}+S_{s-1}-S_s+\b_{s-1}=S_{s-1}-R_{s+1}+\b_{s-1}-2$ so if
$R_{s+1}=R_{s-1}=S_{s-1}$ we get $\b_{s-1}\leq\a_{s-1}+2$, as
claimed. So we may assume that $R_{s+1}>R_{s-1}=S_{s-1}$. Together
with $R_s\geq R_{s-2}=S_{s-2}$, this implies $\a_{s-1}=A_{s-1}=\(
A_{s-1}=\min\{ (R_s-S_{s-1})/2+e,R_s-S_{s-1}+d[-a_{1,s}b_{1,s-2}]\}
=\min\{ (R_s-R_{s-1})/2+e,R_s-R_{s-1}+d[-a_{1,s}b_{1,s-2}]\}$. This
also holds if $s=2$ or $n$, i.e. if $s-1=1$ or $n-1$. But
$\a_{s-1}<(R_s-R_{s-1})/2+e$ so
$R_s-R_{s-1}+d[-a_{s-1,s}]=\a_{s-1}=A_{s-1}=
R_s-R_{s-1}+d[-a_{1,s}b_{1,s-2}]$. Thus
$d[-a_{s-1,s}]=d[-a_{1,s}b_{1,s-2}]$. If $d[-a_{s-1,s}]\geq
d[-b_{s-1,s}]$ then $\b_{s-1}\leq S_s-S_{s-1}+d[-b_{s-1,s}]\leq
R_s-R_{s-1}+d[-a_{s-1,s}]+2=\a_{s-1}+2$ so we are done. So we may
assume that $d[-b_{s-1,s}]>d[-a_{s-1,s}]$. If $s<n$ then
$d[a_{1,s}b_{1,s}]\geq A_s=\b_s\geq d[-b_{s-1,s}]$ and same happens if
$s=n$, when $d[a_{1,s}b_{1,s}]=\j$. By the domination principle we get
$d[-a_{1,s}b_{1,s-2}]\geq d[-b_{s-1,s}]>d[-a_{s-1,s}]$. Contradiction.

So $\b_{s-1}\leq\a_{s-1}+2$. If $1\leq i<s$ is arbitrary then
$\b_i\leq S_{s-1}-S_i+\b_{s-1}\leq R_{s-1}-R_i+\a_{s-1}+2$. Therefore
$\a_i=\min\{\b_i,R_{s-1}-R_i+\a_{s-1}\}\geq\{ \b_i,\b_i-2\} =\b_i-2$,
as claimed.

By duality at index $n-i$ we have $\b_i\leq\a_i\leq\b_i+2$ for
$u\leq i<n$.

To complete the proof we have to prove the $A_i=0$ part of (i). If
$s\leq i<t$, $i\ev s\m2$ then $\a_i=1$. If $t'<\leq i<u$, $i\ev
u-1\m2$ then $\b_i=1$. If $M,N$ are of type II and $t\leq i<t'$
then $\a_i=\b_i=1$. So if $S_1+\cdots+S_i=R_1+\cdots+R_i+1$ then
$\a_i=\b_i=0$ can only happen if $M,N$ are of type II and $i=t$.
In this case note that $R_{t+1}-R_t=(S-1)-R=S-(R+1)=S_{t+1}-S_t$
so we have $\a_t=\b_t=0$ iff $S-R-1=R_{t+1}-R_t=S_{t+1}-S_t=-2e$.

We have to prove that $A_i=\min\{ (R_{i+1}-S_i)/2+e,
R_{i+1}-S_i+d[-a_{1,i+1}b_{1,i-1}],
R_{i+1}+R_{i+2}-S_{i-1}-S_i+d[a_{1,i+2}b_{1,i-2}]\}\geq 0$. To do this
we will show that $R_{i+1}-S_i\geq -2e$ so $(R_{i+1}-S_i)/2+e\geq 0$,
$R_{i+1}-S_i+d[-a_{1,i+1}b_{1,i-1}]\geq 0$ and $R_{i+1}+R_{i+2}\geq
S_{i-1}+S_i$ so $R_{i+1}+R_{i+2}-S_{i-1}-S_i+d[a_{1,i+2}b_{1,i-2}]\geq
R_{i+1}+R_{i+2}-S_{i-1}-S_i\geq 0$.

If $s\leq i<t$, $i\ev s\m2$ then $R_{i+1}=S_{i+1}+1$ and
$R_i=R_{i+2}=R=S_i-1$ so $R_{i+1}+R_{i+2}=S_i+S_{i+1}\geq S_{i-1}+S_i$
and $R_{i+1}-S_i=S_{i+1}-S_i+1\geq -2e+1$. Also $i+1\ev s+1\m2$ so, as
proved above, we have $A_{i+1}=\a_{i+1}=R_{i+2}-R_{i+1}+1$, which
implies $R_{i+1}-S_i+d[a_{1,i+1}b_{1,i+1}]\geq
R_{i+1}-S_i+A_{i+1}=R_{i+2}-S_i+1=0$. Since also
$R_{i+1}-S_i+d[-b_{i,i+1}]\geq R_{i+1}-S_{i+1}+\b_i=1+\b_i>0$ we have
$R_{i+1}-S_i+d[-a_{1,i+1}b_{1,i-1}]\geq 0$.

The case $t'<i<u$ follows by duality at index $n-i$.

If $M,N$ are of type II and $t\leq i<t'$ then $R_{i+1}=S_i=T$ so
$R_{i+1}-S_i=0>-2e$ and
$R_{i+1}-S_i+d[-a_{1,i+1}b_{1,i-1}]=d[-a_{1,i+1}b_{1,i-1}]\geq
0$. We have $S_{i-1}\leq S_{i+1}\leq T+1$ and we cannot have equality
since this would imply $S_i-S_{i-1}=T-(T+1)=-1$. So $S_{i-1}\leq
T$. Similarly $R_{i+2}\geq R_i\geq T-1$ and we cannot have equality
since this would imply $R_{i+2}-R_{i+1}=(T-1)-T=-1$. So $R_{i+2}\geq
T$. Hence $R_{i+1}+R_{i+2}\geq 2T\geq S_{i-1}+S_i$.

Finally, consider the case when $M,N$ are of type III, $i=t$ and
$S-R-1=R_{t+1}-R_t=S_{t+1}-S_t\geq 2-2e$ so $\a_t=\b_t=1$ and $S-R\geq
3-2e$. Hence $R_{t+1}-S_t=(S-1)-(R+1)\geq 1-2e$. We have $R_{t+2}\geq
S_{t+2}\geq S_t$ ($R_{t+2}=S_{t+2}+1$ if $t+1=t'<u$ and
$R_{t+2}=S_{t+2}$ if $t+1=u$) and $R_{t+1}\geq R_{t-1}\geq S_{t-1}$
($R_{t-1}=S_{t-1}+1$ if $s<t$ and $R_{t-1}=S_{t-1}$ if $s=t$). Thus
$R_{t+1}+R_{t+2}\geq S_{t-1}+S_t$. Now both when $s<t$ and when $s=t$
we have $A_{t-1}=\a_{t-1}$ so $R_{t+1}-S_t+d[a_{1,t-1}b_{1,t-1}]\geq
R_{t+1}-S_t+\a_{t-1}\geq R_t-S_t+\a_t=R-(R+1)+1=0$. (Same happens if
$t=1$, when $d[a_{1,t-1}b_{1,t-1}]=\j$.) Since also
$R_{t+1}-S_t+d[-a_{t,t+1}]\geq R_t-S_t+\a_t=0$ we have
$R_{t+1}-S_t+d[-a_{1,t+1}b_{1,t-1}]\geq 0$. \qed


\bff {\bf Remark} Note that it is possible that the pair $M,N$ is both of
type II and III. More precisely if $M,N$ are of type II and $t'=t+1$
they are also of type III. Also if $M,N$ are of type III and $S-R=2$
then they are also of type II.

If $M,N$ are both of type II and III we say they are of type
II/III. In general, unless otherwise specified, the case II/III will
be treated as case II. So if $M,N$ are of the type III we will assume
that $S-R\neq 2$. Thus $R_{t+1}-R_t=S_{t+1}-S_t=S-R-1\neq 1$. Since
$\a_1,\b_1\leq 1$, this implies that $S-R-1$ is $\leq 0$ and even,
i.e. $S-R$ is $\leq 1$ and odd.
\eff

We now to see what the parities of $R_i$'s and $S_i$'s are for indices
$s\leq i\leq u$ when we have each of the three types.

\blm (i) If $M,N$ are of type I then $S_s\ev\ldots\ev S_{t-1}\ev T\m2$
and $S_t\ev\ldots\ev S_u\ev T+1\m2$. Also $R_s\ev\ldots\ev R_{t'}\ev
T-1\m2$ and $R_{t'+1}\ev\ldots\ev R_u\ev T\m2$. In particular, if
$s<t$ then $S_t-S_{t-1}$ is odd and if $t'<u$ then $R_{t'+1}-R_{t'}$
is odd. 

(ii) If $M,N$ are of type II then $S_s\ev\ldots\ev S_{t'-1}\ev T\m2$
and $S_{t'}\ev\ldots\ev S_u\ev T+1\m2$. Also $R_s\ev\ldots\ev R_t\ev
T-1\m2$ and $R_{t+1}\ev\ldots\ev R_u\ev T\m2$.

(iii) If $M,N$ are of type III (but not II/III) then $S_s\ev\ldots\ev
S_u\ev S\m2$ and $R_s\ev\ldots\ev R_u\ev R\m2$.
\elm
\pf We will use Lemma 6.6(i). Recall that $s\ev t\m2$ and $t'\ev
u\m2$. Moreover, if the type is I then $t\ev t'\m2$ while if the type
is III then $t'=t+1$. 

(i) We have $S_t=S_u=T+1$ so $S_t\ev\ldots\ev S_u\ev T+1\m2$ and, if
$s<t$, $S_s=S_{t-2}=T$ so $S_s\ev\ldots\ev S_{t-2}\ev T\m2$. We are
left to prove that if $s<t$ then $S_{t-1}\ev T\m2$. We have
$R_{t-2}=R_t=T-1$ so $R_{t-1}\ev R_t=T-1\m2$. Thus
$S_{t-1}=R_{t-1}-1\ev T\m2$. 

(ii) We have $S_s=S_t=T$ so $S_s\ev\ldots\ev S_t\ev
T=S_{t+1}=\ldots=S_{t'-1}\m2$ and $S_{t'}=S_u$ so $S_{t'}\ev\ldots\ev
S_u=T+1\m2$. 

(iii) We have $S_s=S_t=R+1$ so $S_s\ev\ldots\ev S_t\ev R+1\ev S\m2$
and $S_{t+1}=S_u=S$ so $S_{t+1}\ev\ldots\ev S_u\ev S\m2$. (Recall
$S-R$ is odd so $R+1\ev S\m2$.) 

The parts of (i)-(iii) regarding the parities of $R_i$'s are
similar. They follow by duality. \qed

\bff Note that, as seen from the proof of Lemma 6.9, we have $\a_i=1$
if $s\leq i<t$, $i\ev s\m2$ and $\b_i=1$ for $t'<i<u$, $i\ev
t'+1\m2$. (Not merely $\a_i\leq 1$ resp. $\b_i\leq 1$, as stated by
Lemma 6.9(i).)

Also note that for $s<i<t$, $i\ev s+1\m2$ we have $\a_{i-1}=1>0$ so
$R_i-R_{i-1}\geq 2-2e$. Since $R_{i-1}=R_{i+1}=R$ we get $R_{i+1}-R_i\leq
2e-2$. Similarly, by duality at index $n-i$, we have $S_{i+1}-S_i\leq 2e-2$
for $t'\leq i<u$, $i\ev t'\m2$. 
\eff

\bff Let's see how we apply Lemma 6.9(ii), namely when is $A_i=\a_i$
and when $A_i=\b_i$.

If $i<s$ then $A_i=\a_i$. If $i\geq u$ then $A_i=\b_i$.

If $s<i<t$, $i\ev s+1\m2$ then $A_i=\a_i=R_{i+1}-R_i+1$. If $t'\leq
i<u$, $i\ev u\m2$ we have $A_i=\b_i=S_{i+1}-S_i+1$.

If $M,N$ are of type I and $t\leq i<t'$ then $R_i+\a_i=S_i+\b_i$. If
$i\ev t\m2$ then $S_i=T+1=R_i+2$ so $\a_i=\b_i+2$. If $i\ev t+1\m2$
then $R_i=S_i+2$ so $\b_i=\a_i+2$. Therefore $A_i=\b_i=\a_i-2$ if
$i\ev t\m2$ and $A_i=\a_i=\b_i-2$ if $i\ev t+1\m2$.
\eff

\bff We claim that for $s\leq i\leq t$ we have $R_i\leq R$. Indeed if
$i\ev s\m2$ then $R_i=R$ and if $i\ev s+1\m2$ then $\a_{i-1}=1$
so $R_i-R_{i-1}\leq 1$. But we cannot have equality because this would
imply $R_{i+i}-R_i=R_{i-1}-R_i=-1$. (We have $R_{i+1}=R_{i-1}=R$.) So
$R_i-R_{i-1}\leq 0$ i.e. $R_i\leq R_{i-1}=R$.

By duality at index $n-i+1$ we have the similar result, $S_i\geq S$
for $t'\leq i\leq u$.
\eff

\bff We prove that for $s\leq i<u$ we have $S_{i+1}-S_i\leq 2e+1$ with
equality iff $M,N$ are of type I, $s<t$, $i=t-1$ and
$S_{t-1}=T-2e$. 

We have $S_{t'}=S_u=S$ so $S_{i+1}-S_i\leq 2e$ holds for $t'\leq i<u$
by Lemma 6.6(i). Regarding indices $t\leq i<t'$: if we have type I
then $S_t=S_{t'}=T+1$ so $S_{i+1}-S_i\leq 2e$ holds for $t\leq i<t'$;
if we have type II then the sequence $S_t,\ldots,S_{t'}$ is
$T,\ldots,T,T+1$ so $S_{i+1}-S_i\leq 1<2e$ holds for $t\leq i<t'$; if
we have type III then $\b_t\leq 1$ so $S_{t+1}-S_t\leq 1<2e$. If we
have type II or III then $S_s=S_t=R+1$ so, by Lemma 6.6(i),
$S_{i+1}-S_i\leq 2e$ holds for $s\leq i<t$. If the type is I then
$S_s=S_{t-2}=R+1=T$ so $S_{i+1}-S_i\leq 2e$ only holds for $s\leq
i<t-2$. So we are left with the cases $i=t-2$ and $t-1$. If 
$s\leq i=t-2$ then $S_{t-1}-S_t\leq 2e$ follows from $\b_{t-2}\leq
1$. If $s<i=t-1$ then $S_t=T+1=S_{t-2}+1$ so
$S_t-S_{t-1}=1-(S_{t-1}-S_{t-2})\geq 2e+1$ with equality when
$S_{t-1}=S_t-2e-1=T-2e$. 

By duality at index $n-i$ we have the similarly result,
$R_{i+1}-R_i\leq 2e+1$ for $s\leq i<u$ with equality iff $M,N$ are of
type I, $t'<u$, $i=t'$ and $R_{t'+1}=T+2e$. \eff

\bff If $S_1+\cdots+S_i\neq R_1+\cdots+R_i+1$ then either
$\b_i\geq\a_i\geq\b_i-2$ and $A_i=\a_i$ or $\a_i\geq\b_i\geq\a_i-2$
and $A_i=\b_i$. In the first case $\a_i\geq d[a_{1,i}b_{1,i}]\geq
A_i=\a_i$ so $d[a_{1,i}b_{1,i}]=\a_i$. 

Moreover, if $K\ap\[ c_1,\ldots,c_k\]$ is a third lattice, $0\leq
j\leq k$ and $\e\in\ff$ then $d[\e
  a_{1,i}c_{1,j}]\leq\a_i=d[a_{1,i}b_{1,i}]$ so by the domination
principle we have $d[\e a_{1,i}c_{1,j}]=\min\{
d[\e b_{1,i}c_{1,j}],d[a_{1,i}b_{1,i}]\} =\min\{
d[\e b_{1,i}c_{1,j}],\a_i\}$. Since also $\a_i\geq\b_i-2\geq d[\e
b_{1,i}c_{1,j}]-2$ we have $d[\e b_{1,i}c_{1,j}]\geq d[\e
a_{1,i}c_{1,j}]\geq d[\e b_{1,i}c_{1,j}]-2$. In particular, if $j=0$
we get $d[\e a_{1,i}]=\min\{ d[\e b_{1,i}],\a_i\}$ and $d[\e
  b_{1,i}]\geq d[\e a_{1,i}]\geq d[\e b_{1,i}]-2$.

Similarly, when $A_i=\b_i$ and $\a_i\geq\b_i\geq\a_i-2$, we have the
dual relations $d[a_{1,i}b_{1,i}]=\b_i$, $d[\e b_{1,i}c_{1,j}]=\min\{
d[\e a_{1,i}c_{1,j}],\b_i\}$ and $d[\e a_{1,i}c_{1,j}]\geq d[\e
b_{1,i}c_{1,j}]\geq d[\e a_{1,i}c_{1,j}]-2$. Also $d[\e
b_{1,i}]=\min\{ d[\e a_{1,i}],\b_i\}$ and $d[\e a_{1,i}]\geq d[\e
b_{1,i}]\geq d[\e a_{1,i}]-2$.
\eff

\section{Reduction to the case $\nn M=\nn N$}

In this section we use the notations from $\S$4 for the lattices
$M,N,K$. Given $M,K$ two lattices on the same quadratic space
s.t. $K\leq M$ and $\nn K\sb\nn M$, we will show that there is $N\sb
M$ s.t. $K\leq N$. This reduces the proof to the case $\nn M=\nn
K$. Indeed $K\leq M$, more precisely the condition (i) of the main
theorem, implies by Lemma 5.5(i) that $\ord volM=R_1+\cdots+R_n\leq
T_1+\cdots+T_n=\ord vol K$. Suppose that our theorem is not true and
let $M,K$ be two lattices s.t. $K\leq M$ but $K\notrep M$ and,
moreover, $M,K$ with this property are s.t. $\ord volK-\ord volM$ is
minimal. Assume that $\nn K\sb\nn M$. Let $N\sb M$ s.t. $K\leq N$. We
have $\ord volK-\ord volN<\ord volK-\ord volM$. By the minimality of
$\ord volK-\ord volM$ we have $K\rep N$, which, together with $N\sb M$
implies $K\rep M$. Contradiction. Thus $\nn K=\nn M$.

\blm Let $L$ be a lattice with $R_i(L)=R_i$. Let $L'=\{ x\in L\mid x$
not a norm generator$\}$. If $R_2-R_1\neq -2e$ then $L'$ is a lattice
and $[L:L']=\p$.

Also $R_2-R_1>-2e$ if and only if $\ss L\sb\h 12\nn L$.
\elm
\pf We have $\ord\nn L=R_1$ and by [B1, Corollary 4.4(iv)] $\od\ss
L=\min\{ R_1,(R_1+R_2)/2\}$. Hence $\ss L\sb\h 12\nn L$ iff
$\min\{ R_1,(R_1+R_2)/2\} >R_1-e$, i.e. iff $(R_1+R_2)/2>R_1-e$,
which is equivalent to $R_2-R_1>-2e$ 

Since $\nn L=\p^{R_1}$ we have $x\in L'$ iff $x\in L$ and $\od
Q(x)>R_1$. Obviously if $x\in L'$ and $\a\in\oo$ then $\a x\in
L'$. So we still have to prove that if $x,y\in L'$ then $x+y\in
L'$. Suppose the contrary, i.e. that there are $x,y\in L'$
s.t. $x+y\in L\setminus L'$. This means $\od Q(x),\od Q(y)>R_1$
and $\od Q(x+y)=R_1$. By [B1, Lemma 3.19] we have $\oo x+\oo
y\ap\pi^{R_1-e}\aa$ and so $\pi^{R_1-e}\aa\rep L$. But this is
impossible since $\ss L\sb\p^{R_1-e}$.

Let $x_1\in L$ be a norm generator. We have $x_1\in L\setminus L'$ but
$\pi x_1\in L'$. In order to prove that $[L:L']=\p$ it is enough to
prove that $L=L'+\oo x_1$. Let $x\in L$. If $x\in L'$ then $x\in
L'+\oo x_1$ and we are done. If $x\in L\setminus L'$ then $x$ is a
norm generator of $L$ and so $\od Q(x)=R_1=\od Q(x_1)$. Thus
$Q(x)/Q(x_1)\in\ooo$ and so there is $\a\in\ooo$
s.t. $-Q(x)/Q(x_1)\ev\a^2(\mo\p )$. It follows that $\od
(Q(x)/Q(x_1)+\a^2)>0$ and so $\od (Q(x)+\a^2Q(x_1))>\od
Q(x_1)=R_1$. On the other hand $B(x,x_1)\in\ss L\sb\p^{R_1-e}$ so
$\od 2\a B(x,x_1)=\od B(x,x_1)+e>R_1$. Together with $\od
(Q(x)+\a^2Q(x_1))>R_1$, this implies that $\od Q(x-\a x_1)=\od
(Q(x)+\a^2Q(x_1)-2\a B(x,x_1))>R_1$. Therefore $x-\a x_1\in L'$. It
follows that $x=(x-\a x_1)+\a x_1\in L'+\oo x_1$. \qed

\subsection{The case $R_2-R_1>-2e$}

Suppose first that $\nn M\sp 2\ss M$, i.e. that $R_2-R_1>-2e$. By
Lemma 7.1 $M':=\{ x\in M|x\text{ not a norm generator}\}$ is a
lattice and we have $[M:M']=\p$. We will take $N=M'$. We have $\nn
N\sb\nn M$, so $S_1>R_1$ and $[M:N]=\p$ so $M,N$ satisfy Lemmas 6.7
and 6.9 with $s=1$. In particular, $t$ is odd. Also if we have type I
then $t'$ and $u$ are odd and if we have type III then $t'=t+1$ and
$u$ are even.

\blm If $M,N$ are as above then:

(i) If $M,N$ are of type I then $S_1+\cdots+S_i$ is $\ev iT\m2$ if
$i<t$ and it is $\ev i(T+1)\m2$ if $t\leq i\leq u$. Also
$R_1+\cdots+R_i$ is $\ev i(T-1)\m2$ if $i\leq t'$ and it is $\ev
iT-1\m2$ if $t\leq i\leq u$.

(ii) If $M,N$ are of type II then $S_1+\cdots+S_i$ is $\ev iT\m2$ if
$i<t'$ and it is $\ev i(T+1)+t'-1\m2$ if $t'\leq i\leq u$. Also
$R_1+\cdots+R_i$ is $\ev i(T-1)\m2$ if $i\leq t$ and it is $\ev
iT-1\m2$ if $t<i\leq u$. 

(iii) If $M,N$ are of type III then $S_1+\cdots+S_i\ev iS\m2$ and
$R_1+\cdots+R_i\ev iR\m2$ for any $1\leq i\leq u$.
\elm
\pf We use Lemma 6.11 for the parities of $R_i$'s and $S_i$'s.

(i) We have $S_j\ev T\m2$ for $1\leq j<t$ so $S_1+\cdots +S_i\ev
iT\m2$ for $i<t$. In particular, $S_1+\cdots+S_{t-1}$ is even ($t$ is
odd). If $t\leq j\leq u$ then $S_j\ev T+1\m2$ so for $t\leq i\leq
u$ we have $S_t+\cdots +S_i\ev (i-t+1)(T+1)\ev i(T+1)\m2$. Since also
$S_1+\cdots +S_{t-1}$ is even we get $S_1+\cdots+S_i\ev i(T+1)\m2$.

Similarly we have $R_1+\cdots +R_i\ev i(T-1)\m2$ for $i\leq t'$. In
particular, $R_1+\cdots +R_{t'}\ev T-1\m2$ ($t'$ is odd.) If $t'<i\leq
u$ then $R_{t'+1}+\cdots +R_i\ev (i-t')T\ev (i-1)T\m2$. Together with
$R_1+\cdots +R_{t'}\ev T-1\m2$, this implies $R_1+\cdots +R_i\ev
T-1+(i-1)T=iT-1\m2$.

(ii) We have $S_1+\cdots +S_i\ev iT\m2$ for $i<t'$. In particular,
$S_1+\cdots +S_{t'-1}\ev (t'-1)T\m2$. If $t'\leq i\leq u$ then
$S_{t'}+\cdots +S_i\ev (i-t'+1)(T+1)\m2$ so $S_1+\cdots +S_i\ev
(t'-1)T+(i-t'+1)(T+1)\ev i(T+1)+t'-1\m2$.

We have $R_1+\cdots +R_i\ev i(T-1)\m2$ for $i\leq t$. In particular,
$R_1+\cdots +R_t\ev T-1\m2$ ($t$ is odd). If $t<i\leq u$ then
$R_{t+1}+\cdots +R_i\ev (i-t)T\ev (i-1)T\m2$ so $R_1+\cdots +R_i\ev
T-1+(i-1)T=iT-1\m2$.

(iii) Follows trivially from Lemma 6.11(iii). \qed

\blm (i) If $i<j$ and $R_i+\a_i=R_j+\a_j$ then $R_i\ev\ldots\ev
R_j\m2$ and $\a_i\ev\ldots\ev\a_j\m2$. Also for any $i\leq k\leq j$ we
have $R_{k+1}-R_k\leq 2e$ and $\a_k\leq 2e$.

(ii) If $i<j$ and $-R_{i+1}+\a_i=-R_{j+1}+\a_j$ then
$R_{i+1}\ev\ldots\ev R_{j+1}\m2$ and $\a_i\ev\ldots\ev\a_j\m2$. Also
for any $i\leq k\leq j$ we have $R_{k+1}-R_k\leq 2e$ and $\a_k\leq
2e$. \elm

\pf (i) Note that the sequence $(R_k+\a_k)$ is increasing so
$R_i+\a_i=R_j+\a_j$ implies
$R_i+\a_i=R_{i+1}+\a_{i+1}=\ldots=R_j+\a_j$. Hence the statements
$R_i\ev\ldots\ev R_j\m2$ and $\a_i\ev\ldots\ev\a_j\m2$ are equivalent
so we only have to prove one of them. Also $R_{k+1}-R_k\leq 2e$ is
equivalent to $\a_k\leq 2e$ by [B1, Corollary 2.8(ii)].

Take first the case when $i=k$, $j=k+1$. If $\a_k<(R_{k+1}-R_k)/2+e$
and $\a_{k+1}<(R_{k+2}-R_{k+1})/2+e$ then both $\a_k$ and $\a_{k+1}$
are odd and $<2e$ so we are done. If $\a_k=(R_{k+1}-R_k)/2+e$ then
$(R_k+R_{k+1})/2+e=R_k+\a_k=R_{k+1}+\a_{k+1}$ so
$\a_{k+1}=e-(R_{k+1}-R_k)/2$. But $R_{k+1}-R_k\geq -2e$ ao
$\a_{k+1}\leq 2e$. We have $e-(R_{k+1}-R_k)/2=\a_{k+1}\geq 0$ so
$R_{k+1}-R_k\leq 2e$. Also $R_{k+1}-R_k$ cannot be odd since this
would imply that $\a_k=R_{k+1}-R_k\neq (R_{k+1}-R_k)/2+e$. Thus
$R_k\ev R_{k+1}\m2$ and we are done. Finally, if
$\a_k<(R_{k+1}-R_k)/2+e$ and $\a_{k+1}=(R_{k+2}-R_{k+1})/2+e$ then
$(R_k+R_{k+1})/2+e>R_k+\a_k=R_{k+1}+\a_{k+1}=(R_{k+1}+R_{k+2})/2+e$ so
$R_k>R_{k+2}$. Contradiction.

For the general case we just note that
$R_i+\a_i=R_{i+1}+\a_{i+1}$,\ldots,$R_{j-1}+\a_{j-1}=R_j+\a_j$ so we
can apply the case above for $k=i,i+1,\ldots,j-1$.

(ii) follows by applying (i) to the dual lattice with $i,j$ replaced
by $n-j$ and $n-i$. \qed

\blm Let $1\leq i\leq j\leq n$, $i\ev j\m2$. If $R_i=R_j$ then:

(i) If $j<n$ then $d[(-1)^{(j-i+2)/2}a_{i,j+1}]\geq
R_j-R_{j+1}+\a_j$. 

(ii) If $i>1$ then $d[(-1)^{(j-i+2)/2}a_{i-1,j}]\geq
R_{i-1}-R_i+\a_{i-1}$. 

(iii) If $i<j$ then
$d[(-1)^{(j-i)/2}a_{i,j-1}]=R_{j-2}-R_{j-1}+\a_{j-2}=\a_{j-1}$ and
$d[(-1)^{(j-i)/2}a_{i+1,j}]=R_{i+1}-R_{i+2}+\a_{i+1}=\a_i$.
\elm
\pf (i) Since $R_i=R_j$ we have $R_i=R_{i+2}=\ldots=R_j$. For
any $i\leq k\leq j$, $k\ev i\m2$ we have $R_k=R_j$ so
$d[-a_{k,k+1}]\geq R_k-R_{k+1}+\a_k=R_j-R_{k+1}+\a_k\geq
R_j-R_{j+1}+\a_j$. By the domination principle we have
$d[(-1)^{(j-i+2)/2}a_{i,j+1}]\geq R_j-R_{j+1}+\a_j$.

(ii) follows from (i) applied to the dual lattice with $i,j$ replaced
by $n-j+1,n-i+1$.

(iii) We have $R_{j-2}=R_j$ and $R_i=R_{i+2}$ so
$R_{j-2}+\a_{j-2}=R_{j-1}+\a_{j-1}$ and
$-R_{i+1}+\a_i=-R_{i+2}+\a_{i+1}$ by [B3, Corollary 2.3(i)]. Since
$R_i=R_{j-2}$ and $R_{i+2}=R_j$ we can use (i) and (ii) and we get
$\a_{j-1}\geq d[(-1)^{(j-i)/2}a_{i,j-1}]\geq
R_{j-2}-R_{j-1}+\a_{j-2}=\a_{j-1}$ and $\a_i\geq
d[(-1)^{(j-i)/2}a_{i+1,j}]\geq R_{i+1}-R_{i+2}+\a_{i+1}=\a_i$ and we
are done. \qed

\blm If $1\leq i\leq j<n$, $i\ev j\m2$, $R_i=R$ and $R_{j+1}=R-2e$
then $d[(-1)^{(j-i+2)/2}a_{i,j+1}]\geq 2e$, $\[ a_1,\ldots,a_{i+1}\]$
is an orthogonal sum of copies of $\h 12\pi^R\aa$ and $\h 12\pi^R\ab$
and $[a_1,\ldots,a_{i+1}]\ap\hh\pp\cdots\pp\hh$ or
$\hh\pp\cdots\pp\hh\pp [\pi^R,-\D\pi^R]$, where $\hh$ is the
hyperbolic plane.

(Equivalently, $[a_1,\ldots,a_{i+1}]\ap\hh\pp\cdots\pp\hh\pp
[\pi^R,-\eta\pi^R]$, where $\eta\in\{ 1,\D\}$.)
\elm
\pf We have $R=R_i\leq R_j\leq R_{j+1}+2e=R$ so $R_i=R_j=R$. By Lemma
7.4(i) we get $d[(-1)^{(j-i+2)/2}a_{i,j+1}]\geq R_j-R_{j+1}+\a_j$. But
$R_{j+1}-R_j=-2e$ so $\a_j=0$ and $R_j-R_{j+1}+\a_j=2e$. Thus
$d[(-1)^{(j-i+2)/2}a_{i,j+1}]\geq 2e$.

Since $R_i=R_j=R$ we have $R_i=R_{i+2}=\ldots=R_j=R$. Now
$R-2e=R_i-2e\leq R_{i+1}\leq R_{j+1}=R-2e$ so
$R_{i+1}=R_{i+3}=\ldots=R_{j+1}=R-2e$. Hence for any $i<l<j$, $l\ev i+1\m2$
we have $R_l=R-2e<R=R_{l+1}$ so we have the splitting $\[
a_i,\ldots,a_{j+1}\]\ap\[ a_i,a_{i+1}\]\pp\cdots\pp\[ a_j,a_{j+1}\]$.

For any $i\leq l\leq j$, $l\ev i\m2$ we have
$R_{l+1}-R_l=R-2e-R=-2e$. It follows that $a_{l+1}/a_i\in -\h 14\os$
or $-\h\D 4\os$, which implies that, up to scalling, $\[
a_l,a_{l+1}\]$ is isometric to $\aa$ or $\ab$. But $\ord\nn\[
a_l,a_{l+1}\] =R_l=R$ so $\[ a_l,a_{l+1}\]\ap\h 12\pi^R\aa$ or $\h
12\pi^R\ab$. Therefore $\[ a_i,\ldots,a_{i+1}\]\ap\[
a_i,a_{i+1}\]\pp\cdots\pp\[ a_j,a_{j+1}\]$ is an orthogonal sum of
copies of $\h 12\pi^R\aa$ and $\h 12\pi^R\ab$. 

Since $\ab\pp\ab\ap\aa\pp\aa$ we have $\[ a_i,\ldots,a_{j+1}\]\ap\h
12\pi^R\aa\pp\cdots\pp\h 12\pi^R\aa$ or $\h 12\pi^R\aa\pp\cdots\pp\h
12\pi^R\aa\pp\h 12\pi^R\ab$. Since the underlying quadratic spaces
corresponding to $\h 12\pi^R\aa$ and $\h 12\pi^R\ab$ are $\hh$ and
$[\pi^R,-\D\pi^R]$ we have $[a_i,\ldots,a_{i+1}]\ap\hh\pp\cdots\pp\hh$
or $\hh\pp\cdots\pp\hh\pp [\pi^R,-\D\pi^R]$.

Finally, note that $\hh\ap [\pi^R,-\pi^R]$ so
$[a_i,\ldots,a_{i+1}]\ap\hh\pp\cdots\pp\hh\pp [\pi^R,-\eta\pi^R]$,
where $\eta\in\{ 1,\D\}$. \qed

\blm Suppose that $M,N$ are of type I. Then for any $1<i\leq u+1$ even
we have $d[(-1)^{i/2}b_{1,i}]\geq S_{i-1}-S_i+\b_{i-1}$. If $1<i<u$ is
even then $d[(-1)^{i/2}b_{1,i}]=\b_i$ if $S_{i+1}-S_i\leq 2e$ and
$d[(-1)^{i/2}b_{1,i}]\geq 2e$ if $S_{i+1}-S_i=2e+1$.
\elm
\pf Let $1<i\leq u+1$ be even. If $i\leq t-3$ then $S_1=S_{i+1}=T$ so
we can apply Lemma 7.4(iii) and we get
$d[(-1)^{i/2}b_{1,i}]=S_{i-1}-S_i+\b_{i-1}=\b_i$ and we are done.

For $i=t-1$ note that $S_1=S_{t-2}=T$ so by Lemma 7.4(i) we have
$d[(-1)^{(t-1)/2}b_{1,t-1}]\geq S_{t-2}-S_{t-1}+\b_{t-2}$. Now
$\b_{t-2}\leq 1$ by Lemma 6.9(i). If
$T+1-S_{t-1}=S_t-S_{t-1}<2e+1$, i.e. if $S_{t-1}>T-2e$ then
$S_{t-1}-S_{t-2}=S_{t-1}-T>-2e$ so $\b_{t-2}\neq 0$ so
$\b_{t-2}=1$ and so $S_{t-2}-S_{t-1}+\b_{t-2}=T-S_{t-1}+1$. Also
$S_t-S_{t-1}$ is odd and $<2e+1$ so
$\b_{t-1}=S_t-S_{t-1}=T+1-S_{t-1}=S_{t-2}-S_{t-1}+\b_{t-2}$.
Therefore $\b_{t-1}\geq d[(-1)^{(t-1)/2}b_{1,t-1}]\geq
S_{t-2}-S_{t-1}+\b_{t-2}$ implies
$d[(-1)^{(t-1)/2}b_{1,t-1}]=\b_{t-1}$, as claimed. If
$S_t-S_{t-1}=2e+1$, i.e. if $S_{t-1}=T-2e$ then
$S_{t-1}-S_{t-2}=-2e$ so $\b_{t-2}=0$ and
$S_{t-2}-S_{t-1}+\b_{t-2}=2e$. Hence in this case
$d[(-1)^{(t-1)/2}b_{1,t-1}]\geq 2e$.

Suppose now that $t<i\leq u+1$. Then $S_t=S_{i-1}=T+1$ so by Lemma
7.4(i) we have $d[(-1)^{(i-t+1)/2}b_{t,i}]\geq
S_{i-1}-S_i+\b_{i-1}$. In order to prove $d[(-1)^{i/2}b_{1,i}]\geq
S_{i-1}-S_i+\b_{i-1}$ we still need to show that
$d[(-1)^{(t-1)/2}b_{1,t-1}]\geq S_{i-1}-S_i+\b_{i-1}$. But, as seen
above, $d[(-1)^{(t-1)/2}b_{1,t-1}]$ is $=\b_{t-1}$ or $\geq 2e$ so
it is enough to prove that $\b_{t-1},2e\geq S_{i-1}-S_i+\b_{i-1}$. We
have $S_t=S_{i-1}=T+1$ so $\b_{t-1}\geq
S_t-S_i+\b_{i-1}=S_{i-1}-S_i+\b_{i-1}$ and $\b_{i-1}\leq
(S_i-S_{i-1})/2+e$ so $S_{i-1}-S_i+\b_{i-1}\leq e-(S_i-S_{i-1})/2\leq
e-(-2e)/2=2e$ and we are done.

Finally, note that if $t<i<u$ then $S_{i-1}=S_{i+1}=T+1$ so
$S_{i-1}+\b_{i-1}=S_i+\b_i$ and $S_{i-1}-S_i+\b_{i-1}=\b_i$. Therefore
$\b_i\geq d[(-1)^{i/2}b_{1,i}]\geq S_{i-1}-S_i+\b_{i-1}$ implies
$d[(-1)^{i/2}b_{1,i}]=\b_i$. \qed

\blm If $M,N$ are of type I and $t<i\leq u+1$ is even then
$d[(-1)^{i/2}a_{1,i}]\geq R_{i-1}-R_i+2$. (If $u=n$ we ignore the case
$i=u+1$.)
\elm
\pf Let $t<i<u$ be even. If $t'<i<u$ then $S_i\geq S=T+1$ by 6.14 so
$R_i=S_i+1\geq T+2$ so $R_{i-1}-R_i+2=T-R_i+2\leq 0\leq
d[(-1)^{i/2}a_{1,i}]$ so we are done. Similarly, if $i=u+1$ and $t'<u$
then $R_{u+1}\geq R_{u-1}=S_{u-1}+1\geq T+2$ so
$R_u-R_{u+1}+2=T-R_{u+1}+2\leq 0\leq d[(-1)^{(u+1)/2}a_{1,u+1}]$ and
we are done. Thus we may assume that $t<i<t'$ or $i=u+1$ and $t'=u$.

In the first case $R_{i-1}=T-1=R_1$ so by Lemma 7.4(i) we have
$d[(-1)^{i/2}a_{1,i}]\geq R_{i-1}-R_i+\a_{i-1}$. But $i-1$ is odd so
$\a_{i-1}=\b_{i-1}+2\geq 2$ (see 6.13.) so $d[(-1)^{i/2}a_{1,i}]\geq
R_{i-1}-R_i+2$.

In the second case $S_{u+1}=R_{u+1}$ and, since $t'=u$, we have
$S_u=T+1=R_u+2$ so $R_u-R_{u+1}+2=S_u-S_{u+1}$. By Lemma 7.4(i) we
have $d[(-1)^{(u+1)/2}b_{1,u+1}]\geq S_u-S_{u+1}+\b_u\geq
S_u-S_{u+1}$. Since also $A_{u+1}=\b_{u+1}$ or $u+1=n$ (see 6.13.) we
have by 6.16 $d[(-1)^{(u+1)/2}a_{1,u+1}]\geq
d[(-1)^{(u+1)/2}b_{1,u+1}]\geq S_u-S_{u+1}=R_u-R_{u+1}+2$. \qed

\blm If $M,N$ are of type III (but not II/III) then for any $t<i\leq
u$ even we have
$d[(-1)^{i/2}a_{1,i}]=d[(-1)^{i/2}b_{1,i}]=R-S+2$. Also $\a_t=\b_t=1$
and $S-R\geq 3-2e$. (In particular, the exceptional case $\a_i=\b_i=0$
from Lemma 6.9(i) does not occur.)
\elm
\pf We note that $R_t=R_1=R$ and $R_{t+1}\geq R_2$ so
$S-R-1=S_{t+1}-S_t=R_{t+1}-R_t\geq R_2-R_1>-2e$. So $\a_t,\b_t>0$,
which, together with $\a_i,\b_i\leq 1$ (Lemma 6.9(i)), implies
$\a_t=\b_t=1$. Also $S-R-1=R_{t+1}-R_t\geq 2-2e$ so $S-R\geq 3-2e$.

Take first $i=t+1$. By Lemma 6.9(i) we have $A_t=0$. Since
$R_{t+1}-S_t=(S-1)-(R+1)>-2e$ we have
$(R_{t+1}-S_t)/2+e>0=A_t$. Suppose
$0=A_t=R_{t+1}+R_{t+2}-S_{t-1}-S_t+d[a_{1,t+2}b_{1,t-2}]\geq
R_{t+1}+R_{t+2}-S_{t-1}-S_t$. But $R_{t+2}\geq S_{t+2}\geq S_t$. (We
have $R_{t+2}=S_{t+2}+1$ if $t+1=t'<u$ and $R_{t+2}=S_{t+2}$ if
$t+1=u$.) Also $R_{t+1}\geq R_{t-1}=S_{t-1}+1>S_{t-1}$ and so
$R_{t+1}+R_{t+2}>S_{t-1}+S_t$. Contradiction. It follows that
$0=A_t=R_{t+1}-S_t+d[-a_{1,t+1}b_{1,t-1}]$ so
$d[-a_{1,t+1}b_{1,t-1}]=S_t-R_{t+1}=R-S+2$. Now $S_1=S_{t-2}=R+1$ so
by Lemma 7.4(i) we have $d[(-1)^{(t-1)/2}b_{1,t-1}]\geq
S_{t-2}-S_{t-1}+\b_{t-2}\geq S_{t-2}-S_{t-1}\geq R-S+3$. (We have
$S_{t-2}=R+1$ and $S_{t-1}=R_{t-1}-1\leq R_{t+1}-1=S-2$.) Note that
this inequality holds also if $t=1$, when
$d[(-1)^{(t-1)/2}b_{1,t-1}]=\j$. Together with
$d[-a_{1,t+1}b_{1,t-1}]=R-S+2$, this implies
$d[(-1)^{(t+1)/2}a_{1,t+1}]=R-S+2$, as claimed. 

Now if $t+1<i\leq u$ is even then $R_{t+3}=R_i=S-1$ so by Lemma
7.4(ii) we have $d[(-1)^{(i-t-1)/2}a_{t+2,i}]\geq
R_{t+2}-R_{t+3}+\a_{t+2}\geq R-S+3$. (We have $\a_{t+2}\geq 0$,
$R_{t+3}=S-1$ and $R_{t+2}=S_{t+2}+1\geq S_t+1=R+2$.) Together with
$d[(-1)^{(t+1)/2}a_{1,t+1}]=R-S+2$, this implies
$d[(-1)^{i/2}a_{1,i}]=R-S+2$, as claimed.

Finally, note that $A_i=\b_i$ or $i=n$ (both when $t<i<u$, $i$ even,
and when $i=u$; see 6.13). By 6.16 we have
$d[(-1)^{i/2}b_{1,i}]=\min\{ d[(-1)^{i/2}a_{1,i}],\b_i\}$. ($\b_i$ is
ignored if $i=n$.) So in order to prove that
$d[(-1)^{i/2}b_{1,i}]=R-S+2$ we still need $\b_i\geq R-S+2$ (when
$i<n$). If $S_{i+1}-S_i\geq 2e$ then $\b_i\geq 2e>R-S+2$. If
$S_{i+1}-S_i<2e$ then $\b_i\geq S_{i+1}-S_i$ and if $S_{i+1}-S_i$ is
even then $\b_i\geq S_{i+1}-S_i+1$. But $S_{i+1}\geq S_t=R+1$ and
$S_i=S$ so $S_{i+1}-S_i\geq R-S+1$. If $S_{i+1}-S_i=R-S+1$ then, since
$R-S+1$ is even (see Remark 6.10.), we have $\b_i\geq
S_{i+1}-S_i+1=R-S+2$. Otherwise $\b_i\geq S_{i+1}-S_i\geq R-S+2$. \qed


\blm Suppose that $\nn M\sp 2\ss M$, i.e. $R_2-R_1>-2e$, $K\leq M$ and
$\nn K\sb\nn M$. If $N:=M'$, as defined above, then $K\leq N$.
\elm
\pf The condition $\nn K\sb\nn M$ simply means $T_1>R_1=R$,
i.e. $T_1\geq R+1$. If $M,N$ are in the case I or II this means
$T_1\geq T$. We now prove 2.1(i)-(iv) for $N,K$.
\vskip 3mm

{\bf Proof of 2.1(i)} We must prove that for any $1\leq i\leq n$ we
have $S_i>T_i$ or $T_{i-1}+T_i\geq S_i+S_{i+1}$. There are several
cases:

1. $1\leq i<t$. If $i$ is odd then $S_i=R+1\leq T_1\leq T_i$ and we
are done. 

If $i$ is even and $S_i>T_i$ then also $R_i=S_i+1>T_i$ so
$T_{i-1}+T_i\geq R_i+R_{i+1}$. Unless $M,N$ have
type I and $i=t-1$, we have $R_i=S_i+1$ and $R_{i+1}=R=S_{i+1}-1$
so $S_i+S_{i+1}=R_i+R_{i+1}\leq T_{i-1}+T_i$. In the exceptional
case we have $S_i=R_i-1$ and $S_{i+1}=T+1=R_{i+1}+2$ so
$S_i+S_{i+1}=R_i+R_{i+1}+1$ so the only possibility for (i) to
fail is to have $T_{i-1}+T_i=R_i+R_{i+1}=S_i+S_{i+1}-1$. Suppose this
happens so we have type I, $T_{t-1}<S_{t-1}=R_{t-1}-1$,
i.e. $T_{t-1}\leq R_{t-1}-2$, and $T_{t-2}+T_{t-1}=R_{t-1}+R_t$. We
have $T_{t-2}>R_t$ so $C_{t-1}=C_{t-1}'$ by Lemma 2.14. Hence
$R_t-R_{t-1}+1=\a_{t-1}\geq C_{t-1}=C'_{t-1}=\min\{
R_t-T_{t-1}+d[-a_{1,t}c_{1,t-2}],
R_t+R_{t+1}-T_{t-2}-T_{t-1}+d[a_{1,t+1}c_{1,t-3}]\}$. But
$T_{t-1}\leq R_{t-1}-2$ so $R_t-T_{t-1}+d[-a_{1,t}c_{1,t-2}]\geq
R_t-T_{t-1}\geq R_t-R_{t-1}+2$. Hence $R_t-R_{t-1}+1\geq
R_t+R_{t+1}-T_{t-2}-T_{t-1}+d[a_{1,t+1}c_{1,t-3}]=
R_{t+1}-R_{t-1}+d[a_{1,t+1}c_{1,t-3}]$. (We have
$T_{t-2}+T_{t-1}=R_{t-1}+R_t$.) Thus $R_t-R_{t+1}+1\geq
d[a_{1,t+1}c_{1,t-3}]$. By Lemma 7.7 we have $d[(-1)^{(t+1)/2}
a_{1,t+1}]\geq R_t-R_{t+1}+2$ so $d[(-1)^{(t-3)/2}c_{1,t-3}]\leq
R_t-R_{t+1}+1$. But for any $1\leq j\leq t-4$ odd we have
$d[-c_{j,j+1}]\geq T_j-T_{j+1}+\c_j\geq T_j-T_{j+1}\geq
R_t-R_{t+1}+3$. (We have $T_j\geq T_1\geq T=R_t+1$ and
$T_{j+1}\leq T_{t-1}\leq R_{t-1}-2\leq R_{t+1}-2$.) By domination
principle we get $d[(-1)^{(t-3)/2}c_{1,t-3}]\geq R_t-R_{t+1}+3$.
Contradiction.

2. If $t'<i<u$, $i\ev t'+1\m2$ then $R_i=S_i+1$ and
$R_{i+1}=S-1=S_{i+1}-1$ so either $S_i<R_i\leq T_i$ or
$S_i+S_{i+1}=R_i+R_{i+1}\leq T_{i-1}+T_i$ and we are done.

3. If $i>u$ then either $S_i=R_i\leq T_i$ or
$S_i+S_{i+1}=R_i+R_{i+1}\leq T_{i-1}+T_i$. 

We are left with the cases $t\leq i<t'$ and $t'\leq i\leq u$, $i\ev
t'\m2$. We consider separately the cases when $M,N$ are of type I, II
or III. 

4. First assume we have type I. If $t<i<t'$ is even
then $R_i=S_i+2$ and $S_{i+1}=T+1=R_{i+1}+2$ so either $S_i<R_i\leq
T_i$ or $S_i+S_{i+1}=R_i+R_{i+1}\leq T_{i-1}+T_i$.

So we are left with the cases $t\leq i<t'$ odd and $t'\leq i\leq u$,
$i\ev t'\ev 1\m2$. The two cases are in fact just one, namely $t\leq
i\leq u$, $i$ odd. Let $i$ be such an index. We have $T_i\geq T_1\geq
T$ so if $T_i<S_i=T+1$ then $T_1=T_i=T$. By Lemma 6.6(i)
$T_1+\cdots +T_i\ev T\m2$. On the other hand, by Lemma 7.2(i)
$R_1+\cdots + R_i\ev i(T-1)\ev T-1\m2$. Since $T_1+\cdots +T_i\ev
R_1+\cdots +R_i+1\m2$ we have $0=d[a_{1,i}c_{1,i}]\geq C_i$. If $0\geq
(R_{i+1}-T_i)/2+e$ then $-2e\geq R_{i+1}-T_i\geq R_2-R_1-1$. (We have
$R_{i+1}\geq R_2$ and $T_i=T=R_1+1$.) So $R_2-R_1\leq 1-2e$, which
implies $R_2-R_1=-2e$. Contradiction. If $0\geq
R_{i+1}-T_i+d[-a_{1,i+1}c_{1,i-1}]$ then $d[-a_{1,i+1}c_{1,i-1}]\leq
T_i-R_{i+1}=T-R_{i+1}$. By Lemma 7.7 we have
$d[(-1)^{(i+1)/2}a_{1,i+1}]\geq R_i-R_{i+1}+2=T-R_{i+1}+1$. Together
with $d[-a_{1,i+1}c_{1,i-1}]\leq T-R_{i+1}$, this implies
$d[(-1)^{(i-1)/2}c_{1,i-1}]\leq T-R_{i+1}$. But $T_1=T_{i-2}=T$ so by
Lemma 7.4(i) we have $d[(-1)^{(i-1)/2}c_{1,i-1}]\geq
T_{i-2}-T_{i-1}+\c_{i-2}$. Thus $T-R_{i+1}\geq T-T_{i-1}+\c_{i-2}$,
which implies $T_{i-1}\geq T_{i-1}-\c_{i-2}\geq R_{i+1}\geq
S_{i+1}$. ($R_{i+1}=S_{i+1}+1$ if $i<u$ and $R_{i+1}=S_{i+1}$ if
$i=u$.) If this inequality is strict then $T_{i-1}\geq S_{i+1}+1$
which, together with $T_i=T=S_i-1$, implies $T_{i-1}+T_i\geq
S_i+S_{i+1}$ so we are done. So we may assume that $T_{i-1}=S_{i+1}$
so $T_{i-1}=R_{i+1}=S_{i+1}$ and $\c_{i-1}=0$. Thus
$-2e=T_i-T_{i-1}=T-R_{i+1}\geq d[-a_{1,i+1}c_{1,i-1}]$, which is
impossible. Finally, if $0\geq
R_{i+1}+R_{i+2}-T_{i-1}-T_i+d[a_{1,i+2}c_{1,i-2}]$ then
$T_{i-1}+T_i\geq R_{i+1}+R_{i+2}\geq S_i+S_{i+1}$ and we are
done. (In all cases $i\leq t'-2$, $t'\leq i\leq u-2$ and $i=u$ we have
$R_{i+1}+R_{i+2}=S_{i+1}+S_{i+2}\geq S_i+S_{i+1}$.)

5. Assume now that we have type II. If $i=t$ then $T_i\geq T_1\geq
T=S_i$ so we are done. If $t<i\leq t'-2$ then either $T_i\geq
R_i=T=S_i$ or $T_{i-1}+T_i\geq R_i+R_{i+1}=2T=S_i+S_{i+1}$. 

Take now $i=t'-1$. If and $i$ is odd then $T_i\geq T_1\geq T=S_i$ so
we are done. So we can assume that $i=t'-1$ is even, i.e. $t'$ is
odd. Suppose that 2.1(i) fails so $T=S_{t'-1}>T_{t'-1}$ and
$2T+1=S_{t'-1}+S_{t'}>T_{t'-2}+T_{t'-1}$. Since $T_{t'-1}<T\leq T_1$
we have by Lemma 6.6(ii) that $T_1+\cdots +T_{t'-1}$ is even. On the
other hand $R_1+\cdots +R_{t'-1}\ev (t'-1)T-1\ev 1\m2$ by Lemma
7.2(ii). Since $T_1+\cdots + T_{t'-1}\ev R_1+\cdots + R_{t'-1}+1\m2$ we
have $T_{t'-1}\geq R_{t'}$ or $T_{t'-2}+T_{t'-1}\geq R_{t'}+R_{t'+1}$
by Lemma 6.5. In the first case $T_{t'-1}\geq
R_{t'}=T$. Contradiction. In the second case note that
$R_{t'+1}\geq S_{t'+1}\geq S_{t'-1}=T$. (We have $R_{t'+1}=S_{t'+1}+1$
or $S_{t'+1}$, corresponding to $t'<u$ resp. $t'=u$.) But we cannot
have equality, because this would imply
$S_{t'+1}-S_{t'}=T-(T+1)=-1$. Thus $R_{t'+1}\geq T+1$
and we get $T_{t'-2}+T_{t'-1}\geq R_{t'}+R_{t'+1}=T+R_{t'+1}\geq
2T+1$. Contradiction. 

So we are left with the case when $t'\leq i\leq u$, $i\ev
t'\m2$. Suppose that $T+1=S_i>T_i$. We prove that $S_i+S_{i+1}\leq
T_{i-1}+T_i$. 

If $t'$ is even so $i$ is even then $T_i\leq T\leq T_1$ so
$T_1+\cdots +T_i$ is even by Lemma 6.6(ii). If $i$ is odd then $T_i\leq
T\leq T_1$ implies that $T_1=T_i=T$ so $T_1+\cdots +T_i\ev T\m2$ by
Lemma 6.6(i). In both cases $T_1+\cdots +T_i\ev iT\m2$. On the other
hand by Lemma 7.2(ii) we have $R_1+\cdots +R_i\ev iT-1\m2$. Thus
$T_1+\cdots +T_i\ev R_1+\cdots +R_i+1\m2$, which implies $T_i\geq
R_{i+1}$ or $T_{i-1}+T_i\geq R_{i+1}+R_{i+2}$. In the first case note
that $R_{i+1}\geq S_{i+1}$ ($R_{i+1}=S_{i+1}+1$ if $i\leq u-2$ and
$R_{i+1}=S_{i+1}$ if $i=u$) so $T_i\geq S_{i+1}$. Now $S_{i+1}\geq
S_{t'-1}=T$. But we cannot have $S_{i+1}=T$ since this would imply
$S_{i+1}-S_i=T-(T+1)=-1$. Thus $T_i\geq S_{i+1}\geq
T+1=S_i$. Contradiction. In the second case we have
$R_{i+1}+R_{i+2}=S_{i+1}+S_{i+2}$ (both when $i\leq u-2$ and when
$i=u$). So $T_{i-1}+T_i\geq R_{i+1}+R_{i+2}=S_{i+1}+S_{i+2}\geq
S_i+S_{i+1}$ and we are done.

6. Now suppose that we have type III. Then $i=t$ or $t+1=t'\leq i\leq
u$, $i\ev t'\ev 0\m2$. In the first case $T_i\geq
T_1\geq R+1=S_i$ so we are done. Suppose now that $t+1\leq i\leq u$ is
even. Suppose that $T_i<S_i=S$ so $T_i\leq S-1$. We have by Lemma 7.8
$d[(-)^{i/2}a_{1,i}]=R-S+2$. Also for any $1\leq j\leq i-1$ odd we
have $T_{j+1}\leq T_i\leq S-1$ and $T_j\geq T_1\geq R_1+1=R+1$. Thus
$d[-c_{j,j+1}]\geq T_j-T_{j+1}+\c_j\geq T_j-T_{j+1}\geq R-S+2$. If we
have equality then $T_j-T_{j+1}=R-S+2$ so $T_{j+1}-T_j=S-R-2$. But
this is impossible since $S-R-2$ is odd and $\leq  -1$. ($S-R\leq 1$
and it is odd by 6.10.) So $d[-c_{j,j+1}]>R-S+2$. By domination
principle we get $d[(-1)^{i/2}c_{1,i}]>R-S+2$, which, together with
$d[(-1)^{i/2}a_{1,i}]=R-S+2$, implies $R-S+2=d[a_{1,i}c_{1,i}]\geq
C_i$. We may assume that $R_{i+1}+R_{i+2}>T_{i-1}+T_i$ since otherwise
$T_{i-1}+T_i\geq R_{i+1}+R_{i+2}\geq S_i+S_{i+1}$. (In both cases when
$i\leq u-2$ or $i=u$ we have $R_{i+1}+R_{i+2}=S_{i+1}+S_{i+2}\geq
S_i+S_{i+1}$.) Thus $C_i=\( C_i$ and, since $T_i\leq S-1=R_i\leq
R_{i+2}$, we can ignore $R_{i+1}+R_{i+2}-2T_i+\c_{i-1}$ from the
definition of $\( C_i$. Note that $R_{i+1}\geq S_{i+1}$
($R_{i+1}=S_{i+1}+1$ if $i\leq u-2$ and $R_{i+1}=S_{i+1}$ if
$i=u$). We have $R_{i+1}\geq S_{i+1}\geq S_t=R+1$ and $T_i\leq S-1$ so
$R_{i+1}-T_i\geq R-S+2$. So if $R-S+2\geq (R_{i+1}-T_i)/2+e$ then
$R-S+2\geq (R-S+2)/2+e$ so $R-S+2\geq 2e$, i.e. $2-2e\geq S-R$, which
contradicts Lemma 7.8. Suppose now that $R-S+2\geq
R_{i+1}-T_i+d[-a_{1,i+1}c_{1,i-1}]\geq R_{i+1}-T_i$. Since also
$R_{i+1}-T_i\geq R-S+2$ we must have equality, which implies that
$R_{i+1}=S_{i+1}=R+1$ and $T_i=S-1$ (see above). Thus
$S_i+S_{i+1}=S+R+1$. On the other hand $T_{i-1}\geq T_1\geq R+1$. If
$T_{i-1}=R+1$ then $T_i-T_{i-1}=(S-1)-(R+1)=S-R-2$, which is
impossible because $S-R-2$ is odd and $\leq -1$. So $T_{i-1}\geq R+2$,
which, together with $T_i=S-1$, implies $T_{i-1}+T_i\geq
S+R+1=S_i+S_{i+1}$ and we are done. So we are left with the case
$R-S+2\geq\( C_i=2R_{i+1}-T_{i-1}-T_i+\a_{i+1}$, which implies that
$T_{i-1}+T_i\geq 2R_{i+1}+S-R-2+\a_{i+1}$. But $S=S_i$ and
$R_{i+1}\geq S_{i+1}\geq S_1=R+1$. So $2R_{i+1}+S-R-2+\a_{i+1}\geq
S_{i+1}+R+1+S_i-R-2+\a_{i+1}=S_i+S_{i+1}-1+\a_{i+1}\geq S_i+S_{i+1}$
with the exception of the case when $\a_{i+1}=0$,
i.e. $R_{i+2}-R_{i+1}=-2e$, and $R_{i+1}=S_{i+1}=R+1$. Suppose that
this happens. Since $R_{i+2}\geq R_2$ and $R_{i+1}=R+1=R_1+1$ we get
$-2e=R_{i+2}-R_{i+1}\geq R_2-R_1-1$, i.e. $R_2-R_1\leq 1-2e$, so
$R_2-R_1=-2e$. This contradicts the hypothesis. Thus $T_{i-1}+T_i\geq
S_i+S_{i+1}$ and we are done.
\vskip 3mm

{\bf Proof of 2.1(ii)} We have again several cases:

1. $M,N$ are of type I, $i=t-1$, $S_t-S_{t-1}=2e+1$ and
$T_{t-1}=S_{t-1}$. Since $S_t-T_{t-1}>2e$ the condition
$d[b_{1,t-1}c_{1,t-1}]\geq B_{t-1}$ is equivalent by Corollary 2.10
to $b_{1,t-1}c_{1,t-1}\in\fs$.

We have $T_{t-1}=S_{t-1}=S_t-2e-1=T-2e$. Now $T\leq T_1\leq
T_2\leq T_{t-1}+2e=T$ so $T_1=T$. Since $S_1=T_1=T$ and
$S_{t-1}=T_{t-1}=T-2e$ we have by Lemma 7.5 that both
$[b_1,\ldots,b_{t-1}]$ and $[c_1,\ldots,c_{t-1}]$ are
isometric to $\hh\pp\cdots\pp\hh$ or to $\hh\pp\cdots\pp\hh\pp
[\pi^T,-\D\pi^T]$. If $[b_1,\ldots,b_{t-1}]\ap [c_1,\ldots,c_{t-1}]$
then $b_{1,t-1}c_{1,t-1}\in\fs$ and we are done. So we may assume that
$[b_1,\ldots,b_{t-1}]$ and $[c_1,\ldots,c_{t-1}]$ are one
$\ap\hh\pp\cdots\pp\hh$ and the other $\ap\hh\pp\cdots\pp\hh\pp
[\pi^T,-\D\pi^T]$. We also have
$d[(-1)^{(t-1)/2}b_{1,t-1}],d[(-1)^{(t-1)/2}b_{1,t-1}]\geq 2e$ so
$d(b_{1,t-1}c_{1,t-1})\geq 2e$.

We want to prove that $[b_1,\ldots,b_{t-1}]\rep [a_1,\ldots,a_t]$ and
$[c_1,\ldots,c_{t-1}]\rep [a_1,\ldots,a_t]$. In all cases $t<t'$,
$t=t'<u$ and $t=u$ we have $R_{t+1}\geq S_{t+1}\geq
S_t-2e=T-2e+1>T-2e=S_{t-1}=T_{t-1}$ so in order to use 2.1(iii) we
still need to prove that $A_{t-1}+A_t\geq 2e+R_t-S_t$ and
$C_{t-1}+C_t\geq 2e+R_t-T_t$. By 6.13
$A_{t-1}=\a_{t-1}=R_t-R_{t-1}+1=2e-1$. ($R_t=T-1$ and
$R_{t-1}=S_{t-1}+1=T-2e+1$.) Also, by Lemma 6.9(ii)
$A_t=\b_t\geq 0$ (in all cases $t<t'$, $t=t'<u$ and $t=u$) so
$A_{t-1}+A_t\geq 2e-1>2e-2=2e+R_t-S_t$. We prove now that
$C_{t-1}+C_t>2e+R_t-T_t$. To do this, by Corollary 2.17, it is enough
to prove that $d[-a_{1,t+1}c_{1,t-1}]>2e+T_{t-1}-R_{t+1}$. Since
$d[b_{1,t-1}c_{1,t-1}]\geq 2e>2e+T_{t-1}-R_{t+1}$ it is enough to
prove that $d[-a_{1,t+1}b_{1,t-1}]>2e+T_{t-1}-R_{t+1}=T-R_{t+1}$. But
by Lemma 2.12 $A_{t-1}+A_t>2e+R_t-S_t$ implies
$d[-a_{1,t+1}c_{1,t-1}]>e+(S_{t-2}+S_{t-1})/2-R_{t+1}=
e+(T+T-2e)/2-R_{t+1}=T-R_{t+1}$ so we are done.

Hence $[b_1,\ldots,b_{t-1}]\rep [a_1,\ldots,a_t]$ and
$[c_1,\ldots,c_{t-1}]\rep [a_1,\ldots,a_t]$, which implies that
$[a_1,\ldots,a_t]$ represents both $\hh\pp\cdots\pp\hh$ and
$\hh\pp\cdots\pp\hh\pp [\pi^T,-\D\pi^T]$. By Lemma 7.2(i) $\ord
a_{1,t}=R_1+\cdots +R_t\ev tT-1\ev T-1\m2$. Since also $\det
(\hh\pp\cdots\pp\hh\pp [\pi^T,-\D\pi^T])$ has an even order we have
$[a_1,\ldots,a_t]\ap\hh\pp\cdots\pp\hh\pp [\pi^T,-\D\pi^T]\pp
[\pi^{T-1}\e ]$ for some $\e\in\ooo$. Hence
$\hh\pp\cdots\pp\hh\rep [a_1,\ldots,a_t]$ implies that $\hh\rep
[\pi^T,-\D\pi^T,\pi^{T-1}\e ]$, which is impossible since
$[\pi^T,-\D\pi^T,\pi^{T-1}\e ]$ is anisotropic.

2. $1\leq i<t$ is odd. We note that $\b_{i+1}\leq
S_i-S_{i+1}+1$. Indeed, if we have type II or III or if we have type I
and $i<t-2$ then $S_i=S_{i+2}=R+1$ so $S_{i+1}+\b_{i+1}=S_i+\b_i\leq
S_i+1$ so $\b_{i+1}\leq S_i-S_{i+1}+1$. If we have type I and $i=t-2$
then $S_t-S_{t-1}$ is odd by Lemma 6.11(i) so $\b_{t-1}\leq
S_t-S_{t-1}=S_{t-2}-S_{t-1}+1$. (We have $S_t=T+1=S_{t-2}+1$.) Thus
$B_i\leq S_{i+1}-T_i+\b_{i+1}\leq S_i-T_i+1$. But $T_i\geq T_1\geq
R+1=S_i$. If $T_i\geq S_i+1$ then $C_i\leq S_i-T_i+1\leq 0$
which makes $d[b_{1,i}c_{1,i}]\geq C_i$ trivial. So we may assume that
$T_i=S_i=R+1$ so we only have $C_i\leq S_i-T_i+1=1$. Since also
$T_i+T_{i+1}\geq S_i+S_{i+1}$ ($N,K$ satisfy (i) as proved above) we
get $T_{i+1}\geq S_{i+1}$. Now $R+1=T_i\geq T_1\geq R+1$ so
$T_1=R+1$. By Lemma 6.6(i) $T_1+\cdots +T_i\ev R+1\m2$. Also
$S_1=S_i=R+1$ so $S_1+\cdots +S_i\ev R+1\m2$. Thus $\ord
b_{1,i}c_{1,i}=S_1+\cdots +S_i+T_1+\cdots +T_i$ is even, which implies
that $d(b_{1,i}c_{1,i})\geq 1$. If also $\b_i,\c_i\geq 1$ we get
$d[b_{1,i}c_{1,i}]\geq 1\geq C_i$ and we are done. Otherwise $\b_i$ or
$\c_i$ is $0$. If $\b_i=0$ then $-2e=S_{i+1}-S_i=S_{i+1}-T_i$, while
if $\c_i=0$ then $-2e=T_{i+1}-T_i\geq S_{i+1}-T_i$. In both cases
$S_{i+1}-T_i\leq -2e$, which implies $C_i\leq (S_{i+1}-T_i)/2+e\leq 0$
and again we are done.

3. $i$ is even and either $1<i<t$ or we have type I and
$1<i<u$. Suppose that $S_{i+1}-T_i>2e$. By Lemma 2.2 we have
$S_i\leq T_i$ so $S_{i+1}-S_i>2e$. This implies by 6.15 that we have
type I, $i=t-1$ and and $S_{i+1}-S_i=2e+1$. Since
$2e+1=S_{i+1}-S_i\geq S_{i+1}-T_i>2e$ we get $S_i=T_i$ so we are in
the case 1., already discussed. So we may assume that $S_{i+1}-T_i\leq
2e$, which implies that $B_i\leq (S_{i+1}-T_i)/2+e\leq 2e$.

In order to prove $d[b_{1,i}c_{1,i}]\geq B_i$ we will show that
$d[(-1)^{i/2}b_{1,i}],d[(-1)^{i/2}c_{1,i}]\geq B_i$.

We prove that $d[(-1)^{i/2}c_{1,i}]\geq B_i$. By domination principle
there is $1\leq j\leq i-1$ odd s.t. $d[(-1)^{i/2}c_{1,i}]\geq
d[-c_{j,j+1}]\geq T_j-T_{j+1}+\c_j\geq T_j-T_i+\c_{i-1}$. If $T_j\geq
S_{i+1}$ then $T_j-T_i+\c_{i-1}\geq S_{i+1}-T_i+\c_{i-1}\geq B_i$ and
we are done. So we may assume that $S_{i+1}>T_j$. We have $T_j\geq
T_1\geq R+1$. If $M,N$ are of type II or III (so $1<i<t$) or of type I
and $1\leq i\leq t-3$ then $S_{i+1}=R+1\leq T_j$. So we may assume
that we have type I, $t-1\leq i<u$ and we have $S_{i+1}=T+1>T_j\geq
T_1\geq T$, i.e. $T_j=T=S_{i+1}-1$. So we have
$d[(-1)^{i/2}c_{1,i}]\geq T_j-T_{j+1}+\c_j\geq
T_j-T_i+\c_{i-1}=S_{i+1}-T_i+\c_{i-1}-1$. Suppose that we have
equality. It follows that $-T_{j+1}+\c_j=-T_i+\c_{j-1}$ so $T_{j+1}\ev
T_{j+2}\ev\ldots\ev T_i\m2$ by Lemma 7.3(ii). This implies that
$T_{j+1}+\cdots +T_{i-1}$ is even. Also $T\leq T_1\leq T_j=T$ so
$T_1=T_j=T$, which implies $T_1+\cdots +T_j\ev T\m2$ so $T_1+\cdots
+T_{i-1}\ev T\m2$, as well. On the other hand, since $t\leq i+1\leq u$
we have $S_1+\cdots +S_{i+1}\ev (i+1)(T+1)\ev T+1\m2$ by Lemma
7.2(i). Thus $\ord b_{1,i+1}c_{1,i-1}$ is odd so $B_i\leq
S_{i+1}-T_i+d[-b_{1,i+1}c_{1,i-1}]=S_{i+1}-T_i$. If $T_{j+1}-T_j=-2e$
then $T_j-T_{j+1}+\c_j=2e\geq B_i$ and we are done. Otherwise
$T_j-T_{j+1}+\c_j\geq T_j-T_{j+1}+1\geq T_j-T_i+1=S_{i+1}-T_i\geq B_i$
and we are done. If $d[(-1)^{i/2}c_{1,i}]>S_{i+1}-T_i+\c_{i-1}-1$ and
$d[(-1)^{i/2}c_{1,i}],\c_{i-1}\in\ZZ$ then $d[(-1)^{i/2}c_{1,i}]\geq
S_{i+1}-T_i+\c_{i-1}\geq B_i$. If $d[(-1)^{i/2}c_{1,i}]\notin\ZZ$ then
$\c_i\notin\ZZ$, which implies $\c_i>2e$, and we have
$d[(-1)^{i/2}c_{1,i}]=\c_i>2e\geq B_i$. If $\c_{i-1}\notin\ZZ$ then
$T_i-T_{i-1}$ is odd and $>2e$ so $B_i\leq S_{i+1}-(T_{i-1}+T_i)/2+e<
S_{i+1}-(T_{i-1}+T_{i-1}+2e)/2+e=S_{i+1}-T_{i-1}$. We have
$T_{i-1}\geq T_1\geq T$. If $T_{i-1}\geq T+1=S_{i+1}$ then $B_i<0$ and
we are done. If $T_{i-1}=T_1=T$ then $T_1+\cdots +T_{i-1}\ev T\m2$ by
Lemma 6.6(i). We also have $S_1+\cdots +S_{i+1}\ev T+1\m2$ so again
$B_i\leq S_{i+1}-T_i+d[-b_{1,i+1}c_{1,i-1}]=S_{i+1}-T_i<0$.
($T_i>T_{i-1}+2e=T+2e>T+1=S_{i+1}$.)

We prove now that $d[(-1)^{i/2}b_i]\geq B_i$. By Lemma 7.6 if
$S_{i+1}-S_i=2e+1$ then $d[(-1)^{i/2}b_i]\geq 2e\geq B_i$ so we are
done. If $S_{i+1}-S_i\leq 2e$ then $d[(-1)^{i/2}b_i]=\b_i$ so we have
to prove that $\b_i\geq B_i$. If $S_{i+1}-S_i=2e$ then $\b_i=2e\geq
B_i$ so we are done. If $S_{i+1}-S_i<2e$ we claim that
$\b_i=\a_i+2$. If we have type II or III (so $i<t$) or type I and
$i\leq t-3$ then $S_{i-1}=S_{i+1}=R+1$ so
$-S_i+\b_{i-1}=-S_{i+1}+\b_i$. Now $\b_{i-1}\leq 1$ by Lemma
6.9(i). But $S_i-S_{i-1}=S_i-S_{i+1}>-2e$ so $\b_{i-1}=1$ and so
$\b_i=S_{i+1}-S_i+1=R_{i+1}-R_i+3=\a_i+2$. (We have
$S_{i+1}=R+1=R_{i+1}+1$ and $S_i=R_i-1$.) If $t<i<t'$ (when we have
type I) then $R_i+\a_i=S_i+\b_i$ so
$\b_i=R_i-S_i+\a_i=\a_i+2$. Finally, if we have type I and $i=t-1$ then
$S_t-S_{t-1}<2e$ and it is odd so
$\b_{t-1}=S_t-S_{t-1}=R_t-R_{t-1}+3=\a_{t-1}+2$. (We have $S_t=T+1=R_t+2$ and
$R_{t-1}=S_{t-1}+1$.) Since $\b_i=\a_i+2\geq C_i+2$ it is enough to prove that
$C_i+2\geq B_i$. 

If $C_i=(R_{i+1}-T_i)/2+e$ then $C_i+2>(R_{i+1}+2-T_i)/2+e\geq
(S_{i+1}-T_i)/2+e\geq B_i$ so we are done. (In all cases $S_{i+1}\leq
R_{i+1}+2$.) 

Suppose now that $C_i=R_{i+1}-T_i+d[-a_{1,i+1}c_{1,i-1}]$. If we have
type II or III or if we have type I and $i+1\leq t-2$ then
$S_{i+1}=R+1=R_{i+1}+1$ and $\b_{i+1}\leq 1$ so $C_i+2\geq
R_{i+1}-T_i+2=S_{i+1}-T_i+1\geq S_{i+1}-T_i+\b_{i+1}\geq B_i$. If we
have type I and $t\leq i+1\leq t'$ then $S_{i+1}=T+1=R_{i+1}+2$ and
$A_{i+1}=\b_{i+1}$ so $d[-a_{1,i+1}c_{1,i-1}]\geq
d[-b_{1,i+1}c_{1,i-1}]$ by 6.16. It follows that $C_i+2\geq
R_{i+1}-T_i+d[-b_{1,i+1}c_{1,i-1}]+2=
S_{i+1}-T_i+d[-b_{1,i+1}c_{1,i-1}]\geq B_i$.

Suppose now that
$C_i=R_{i+1}+R_{i+2}-T_{i-1}-T_i+d[a_{1,i+2}c_{1,i-2}]$. Assume first
that $R_{i+1}+R_{i+2}=S_{i+1}+S_{i+2}$. Since also
$d[a_{1,i+2}c_{1,i-2}]+2\geq d[b_{1,i+2}c_{1,i-2}]$, by 6.16, we get
$C_i+2\geq S_{i+1}+S_{i+2}-T_{i-1}-T_i+d[b_{1,i+2}c_{1,i-2}]\geq B_i$
and we are done. If $R_{i+1}+R_{i+2}<S_{i+1}+S_{i+2}$ we have either
type II or III and $i=t-1$ or we have type I and $i=t'-1$ or
$u-1$. Suppose that we have type III, including the case of type
II/III, and $i=t-1$. We have $R_{i+1}+R_{i+2}+2=S_{i+1}+S_{i+2}$ and,
both when $i+2=t+1=t'<u$ and when $i+2=u$, we have $A_{i+2}=\b_{i+2}$
or $i+2=n$ so $d[a_{1,i+2}c_{1,i-2}]\geq d[b_{1,i+2}c_{1,i-2}]$ by
6.16. Thus
$C_i+2=R_{i+1}+R_{i+2}-T_{i-1}-T_i+d[a_{1,i+2}c_{1,i-2}]+2\geq
S_{i+1}+S_{i+2}-T_{i-1}-T_i+d[b_{1,i+2}c_{1,i-2}]\geq B_i$. If
we have type II and $i=t-1$ note first that we can assume $t'>t+1$
since otherwise we have the type II/III, treated above. It follows that
$S_{i+1}+S_{i+2}=2T=R_{i+1}+R_{i+2}+1$ and also $1=\b_{i+2}\geq
d[b_{1,i+2}c_{1,i-2}]$. Thus
$C_i+2=R_{i+1}+R_{i+2}-T_{i-1}-T_i+d[a_{1,i+2}c_{1,i-2}]+2\geq
R_{i+1}+R_{i+2}-T_{i-1}-T_i+2=S_{i+1}+S_{i+2}-T_{i-1}-T_i+1\geq
S_{i+1}+S_{i+2}-T_{i-1}-T_i+d[b_{1,i+2}c_{1,i-2}]\geq B_i$. If we have
type I suppose first that $i=u-1$. Then $R_{i+2}=S_{i+2}$,
$S_{i+1}=T+1$ and $R_{i+1}=T-1$ or $T$ corresponding to $t'=u$
resp. $t'<u$. Hence $S_{i+1}+S_{i+2}\leq R_{i+1}+R_{i+2}+2$. Also
$i+2=n$ or $A_{i+2}=\b_{i+2}$ so $d[b_{1,i+2}c_{1,i-2}]\leq
d[a_{1,i+2}c_{1,i-2}]$ by 6.16. Therefore
$C_i+2=R_{i+1}+R_{i+2}-T_{i-1}-T_i+d[a_{1,i+2}c_{1,i-2}]+2\geq 
S_{i+1}+S_{i+2}-T_{i-1}-T_i+d[b_{1,i+2}c_{1,i-2}]\geq B_i$. So we are
left with the case $i=t'-1$ and $t'<u$. (If $t'=u$ then $i=u-1$.) Then
$S_{i+1}=T+1=R_{i+1}+2$ and $R_{i+2}=S_{i+2}+1$ so
$S_{i+1}+S_{i+2}=R_{i+1}+R_{i+2}+1$ and $1=\b_{i+2}\geq
d[b_{1,i+2}c_{1,i-2}]$. Thus
$C_i+2=R_{i+1}+R_{i+2}-T_{i-1}-T_i+d[a_{1,i+2}c_{1,i-2}]+2\geq
R_{i+1}+R_{i+2}-T_{i-1}-T_i+2=S_{i+1}+S_{i+2}-T_{i-1}-T_i+1\geq
S_{i+1}+S_{i+2}-T_{i-1}-T_i+d[b_{1,i+2}c_{1,i-2}]\geq B_i$. 

4. $M,N$ are of type I and $t\leq i<u$, $i$ odd. We have $A_i=\b_i$ so
$d[b_{1,i}c_{1,i}]=\b_i$ or $d[a_{1,i}c_{1,i}]$ by 6.16. In the first case
note that $S_i=S_{i+2}=T+1$ so $S_i+\b_i=S_{i+1}+\b_{i+1}$ and we have
$\b_i=S_{i+1}-S_i+\b_{i+1}$. If $T_i\geq S_i$ then $\b_i\geq
S_{i+1}-T_i+\b_{i+1}\geq B_i$ so we are done. Otherwise we have
$T_i<S_i$ so $i>1$ and $T_{i-1}+T_i\geq
S_i+S_{i+1}=S_{i+1}+S_{i+2}$. Also $T_{i-1}>S_{i+1}$. We have $T\leq
T_1\leq T_{i-2}\leq T_i<S_i=T+1$ which implies that $T_1=T_{i-2}=T$ so
$T_1+\cdots +T_{i-2}\ev T\m2$ by Lemma 6.6(i). We also have
$S_1+\cdots +S_{i+2}\ev (i+2)(T+1)\ev T+1\m2$ by Lemma 7.2(i) and so
$\ord b_{1,i+2}c_{1,i-2}$ is odd. Hence $B_i\leq
S_{i+1}+S_{i+2}-T_{i-1}-T_i+d[b_{1,i+2}c_{1,i-2}]=
S_{i+1}+S_{i+2}-T_{i-1}-T_i\leq 0$ and we are done. If
$d[b_{1,i}c_{1,i}]=d[a_{1,i}c_{1,i}]$ then $d[b_{1,i}c_{1,i}]\geq C_i$
so it is enough to prove that $C_i\geq B_i$. If
$C_i=(R_{i+1}-T_i)/2+e$ then $C_i>(S_{i+1}-T_i)/2+e\geq B_i$ (we have
$R_{i+1}=S_{i+1}+2$ or $S_{i+1}+1$ if $i+1<t'$
resp. $i+1>t'$). Suppose now that
$C_i=R_{i+1}-T_i+d[a_{1,i+1}c_{1,i-1}]$. If $t<i+1<t'$ then
$R_{i+1}=S_{i+1}+2$ and $d[a_{1,i+1}c_{1,i-1}]+2\geq
d[b_{1,i+1}c_{1,i-1}]$ by 6.16 so
$C_i=S_{i+1}-T_i+d[a_{1,i+1}c_{1,i-1}]+2\geq
S_{i+1}-T_i+d[b_{1,i+1}c_{1,i-1}]\geq B_i$. If $t'<i+1<u$ then
$R_{i+1}=S_{i+1}+1$ and $1=\b_{i+1}\geq d[-b_{1,i+1}c_{1,i-1}]$ so
$C_i\geq R_{i+1}-T_i=S_{i+1}-T_i+1\geq
S_{i+1}-T_i+d[b_{1,i+1}c_{1,i-1}]\geq B_i$. We are left with the case
$C_i=R_{i+1}+R_{i+2}-T_{i-1}-T_i+d[a_{1,i+2}c_{1,i-2}]$. In both cases
$i+2\leq t'$ and $t'<i+2\leq u$ we have
$R_{i+1}+R_{i+2}=S_{i+1}+S_{i+2}$ and $A_{i+2}=\b_{i+2}$ or $i+2=n$
so, by 6.16, $d[a_{1,i+2}c_{1,i-2}]\geq d[b_{1,i+2}c_{1,i-2}]$. This
implies that $C_i\geq
S_{i+1}+S_{i+2}-T_{i-1}-T_i+d[b_{1,i+2}c_{1,i-2}]\geq B_i$.

5. We have type II and $t\leq i<t'-1$. We have $B_i\leq
S_{i+1}-T_i+\b_{i+1}=T-T_i+1$. If $T_i>T$ then 
$B_i\leq 0\leq d[b_{1,i}c_{1,i}]$. So we may assume that $T_i\leq
T$. f $i$ is odd then $T\leq T_1\leq T_i\leq T$ implies $T_1=T_i=T$ so
$T_1+\cdots +T_i\ev T\m2$ by Lemma 6.6(i). If $i$ is even then
$T_1\geq T\geq T_i$ so $T_1+\cdots +T_i$ is even by Lemma 6.6(ii). In
both cases $T_1+\cdots +T_i\ev iT\m2$. But also $S_1+\cdots +S_i\ev
iT\m2$ by Lemma 7.2(ii). So $\ord b_{1,i}c_{1,i}$ is even, which
implies $d(b_{1,i}c_{1,i})\geq 1$. We also have $\b_i=1$ and $\c_i\geq
1$ since otherwise $T_i+T_{i+1}=2T_i-2e\leq
2T-2e<2T=S_i+S_{i+1}$. Therefore $d[b_{1,i}c_{1,i}]=1$ so we have to
prove that $1\geq B_i$. 

If $T_i=T$ then $B_i\leq S_{i+1}-T_i+\b_{i+1}=T-T+1=1$ and we are
done. Suppose now that $T_i<T=S_i$ so $T_{i-1}+T_i\geq
S_i+S_{i+1}=2T$. We have $T_i<S_i\leq S_{i+2}$ so $B_i\leq
S_{i+1}+S_{i+2}-T_{i-1}-T_i+d[-b_{1,i+2}c_{1,i}]$. If $i<t'-2$ then
$S_{i+1}+S_{i+2}=2T\leq T_{i-1}+T_i$ so $B_i\leq
d[-b_{1,i+2}c_{1,i}]\leq\b_{i+2}=1$. If $i=t'-2$ then
$S_{i+1}+S_{i+2}=2T+1\leq T_{i-1}+T_i+1$ so $B_i\leq
1+d[-b_{1,i+2}c_{1,i}]$. But $\ord b_{1,i}c_{1,i}$ is even so $\ord
b_{1,i+2}c_{1,i}=\ord b_{1,i}c_{1,i}+S_{i+1}+S_{i+2}=\ord
b_{1,i}c_{1,i}+2T+1$ is odd so $d[-b_{1,i+2}c_{1,i}]=0$. Thus $B_i\leq
1$. 

6. We have 
type II or III and $t'-1\leq i<u$, $i\ev t'-1\m2$. We have
$\b_i=1$, $R_{i+1}=S-1=S_{i+1}-1$ and $S_1+\cdots +S_i=R_1+\cdots
+R_i+1$ so $\ord a_{1,i}c_{1,i}$ and $\ord b_{1,i}c_{1,i}$ have
opposite parities. Thus one of $d(a_{1,i}c_{1,i})$ and
$d(b_{1,i}c_{1,i})$ is $0$ and the other one is $\geq 1$.

Assume first that $T_1+\cdots +T_i\ev S_1+\cdots +S_i\m2$. It follows
that $d[a_{1,i}c_{1,i}]=d(a_{1,i}c_{1,i})=0$ and
$d(b_{1,i}c_{1,i})\geq 1$. If $\c_i=0$ then $T_{i+1}-T_i=-2e$ so
$d[b_{1,i}c_{1,i}]=0=(T_{i+1}-T_i)/2+e\geq B_i$ by 2.6 and we are
done. So we may assume that $\c_i\geq 1$. Together with
$d(b_{1,i}c_{1,i})\geq 1$ and $\b_i=1$, this implies
$d[b_{1,i}c_{1,i}]=1$ so we have to prove that $1\geq B_i$. Since
also $0=d[a_{1,i}c_{1,i}]\geq C_i$ it is enough to prove that
$B_i\leq C_i+1$. If $C_i=(R_{i+1}-T_i)/2+e$ then
$C_i+1>(R_{i+1}+1-T_i)/2+e=(S_{i+1}-T_i)/2+e\geq B_i$. If
$C_i=R_{i+1}-T_i+d[-a_{1,i+1}c_{1,i-1}]$ then note that
$A_{i+1}=b_{i+1}$ or $i+1=n$ (both when $i+1<u$ and when $i+1=u$)
so $d[-a_{1,i+1}c_{1,i-1}]\geq d[-b_{1,i+1}c_{1,i-1}]$ by 6.16. It
follows that $C_i+1\geq R_{i+1}-T_i+d[-b_{1,i+1}c_{1,i-1}]+1=
S_{i+1}-T_i+d[-b_{1,i+1}c_{1,i-1}]\geq B_i$. If
$C_i=R_{i+1}+R_{i+2}-T_{i-1}-T_i+d[a_{1,i+2}c_{1,i-2}]$ then we
have the cases $i+1<u$ and $i+1=u$. If $i+1<u$ then
$R_{i+1}+R_{i+2}=S_{i+1}+S_{i+2}$ and
$d[b_{1,i+2}c_{1,i-2}]\leq\b_{i+2}=1$ so $C_i+1\geq
R_{i+1}+R_{i+2}-T_{i-1}-T_i+1\geq
S_{i+1}+S_{i+2}-T_{i-1}-T_i+d[b_{1,i+2}c_{1,i-2}]\geq B_i$. If
$i+1=u$ then $R_{i+1}+R_{i+2}+1=S_{i+1}+S_{i+2}$ and
$A_{i+2}=\b_{i+2}$ or $i+2=n$ so $d[a_{1,i+2}c_{1,i-2}]\geq
d[b_{1,i+2}c_{1,i-2}]$. Thus
$C_i+1=S_{i+1}+S_{i+2}-T_{i-1}-T_i+d[a_{1,i+2}c_{1,i-2}]\geq
S_{i+1}+S_{i+2}-T_{i-1}-T_i+d[b_{1,i+2}c_{1,i-2}]\geq B_i$.

Assume now that $T_1+\cdots +T_i\ev S_i+\cdots +S_i+1\m2$. Note that
if $i$ is even then $t'$ is odd so we have type II. (If the type is
III then $t'=t+1$ is even.) Also $S_1+\cdots +S_i$ is even by Lemma
7.2(ii) so $T_1+\cdots +T_i$ is odd. (We have $S_1+\cdots +S_i\ev
iT\ev 0\m2$ if $i=t'-1$ and $S_1+\cdots +S_i\ev i(T+1)+t'-1\ev 0\m2$
if $t'<i<u$. Recall, $i$ is even and $t'$ odd.) 

Suppose first that $T_i\ev S=S_{i+1}\m2$. Since $T_1+\cdots +T_i\ev
S_i+\cdots +S_i+1\m2$ we get $T_1+\cdots +T_{i-1}\ev S_i+\cdots
+S_{i+1}+1\m2$. It follows that $\ord b_{1,i+1}c_{1,i-1}$ is odd so
$B_i\leq S_{i+1}-T_i+d[-b_{1,i+1}c_{1,i-1}]=S_{i+1}-T_i$. If
$S_{i+1}\leq T_i$ then $B_i\leq S_{i+1}-T_i\leq 0$ and we are
done. Suppose $T_i<S_{i+1}$. Since $S-R\leq 2$ we have $T_i\leq
S_{i+1}-2=S-2\leq R<T_1$, which implies that $i$ is even and, by Lemma
6.6(ii), $T_1+\cdots +T_i$ is even. Contradiction. 

So we may assume that $T_i\ev S_{i+1}+1\m2$ so $T_1+\cdots +T_{i-1}\ev
S_i+\cdots +S_{i+1}\m2$.

Suppose first that $i$ is even, i.e. $t'$ is odd. This implies that we
have type II and $T_1+\cdots +T_i$ is odd. Note that $S_{i+1}=T+1$ so
$T_i\ev S_{i+1}+1\ev T\m2$. We want to show that there is $k<i$
s.t. $T_{k+1}-T_k$ is odd and $T_k>T$. Since $T_1+\cdots +T_i$ is odd
there is $1\leq j<i$ odd s.t. $T_j+T_{j+1}$ is odd. We have $T_j\geq
T_1\geq T$. If $T_j>T$ we can take $k=j$. If $T_j=T$ note that
$T_{j+1}\ev T_j+1=T+1\ev T_i+1\m2$ so we cannot have $j+1=i$. Thus
$j\leq i-2$. Also since $T_{j+1}-T_j$ is odd we have $T_{j+1}>T_j=T$
and $T_{j+2}>T_j=T$. (If $T_j=T_{j+2}$ then $T_j\ev T_{j+1}\m2$.) We
have $T_{j+1}\ev T_i+1\m2$. Let $j+1\leq k<i$ be maximal s.t. $T_k\ev
T_i+1\m2$. It follows that $T_{k+1}\ev T_i\m2$ so $T_{k+1}-T_k$ is
odd. Also, depending on the parity of $k$, we have $T_k\geq T_{j+1}$
or $T_k\geq T_{j+2}$ so $T_k>T$, as required. Now $\ord c_{k,k+1}$ is
odd so $\c_{i-1}\leq T_i-T_k+d(-c_{k,k+1})=T_i-T_k\leq T_i-T-1$. It
follows that $B_i\leq S_{i+1}-T_i+\c_{i-1}=T+1-T_i+\c_{i-1}\leq 0$ and
we are done. 

Assume now that $i$ is odd, i.e. $t'$ is even. Suppose first that we
have type II. Since $i+1$ and $t'$ are even we have by Lemma 7.2(ii)
that $S_1+\cdots +S_{i+1}\ev (i+1)(T+1)+t'-1\ev 1\m2$, i.e. it is odd
and so is $T_1+\cdots +T_{i-1}$. Also $T_i\ev S_{i+1}+1\ev T\m2$. By
the same reasoning as in the previous case there is an index $1\leq
k<i$ s.t. $T_{k+1}-T_k$ is odd and $T_k>T$ and the proof follows like
before. (This time we use the fact that $T_1+\cdots +T_{i-1}$ is odd
so there is $1\leq j<i-1$ odd s.t. $T_j+T_{j+1}$ is odd.) Suppose now
that we have type III. We have $T_i\geq T_1\geq R+1$ and $T_i\ev
S_{i+1}+1=S+1\ev R\m2$ so $T_i\geq R+2$. We have
$d[(-1)^{(i+1)/2}b_{1,i+1}]=R-S+2$ by Lemma 7.8. If
$d[-b_{1,i+1}c_{1,i-1}]\leq R-S+2$ then $B_i\leq
S_{i+1}-T_i+d[-b_{1,i+1}c_{1,i-1}]\leq S-(R+2)+R-S+2=0$ so we are
done. Otherwise
$d[(-1)^{(i+1)/2}b_{1,i+1}]=R-S+2<d[-b_{1,i+1}c_{1,i-1}]$, which
implies that $d[(-1)^{(i-1)/2}c_{1,i-1}]=R-S+2$. By the domination
principle there is $1\leq j<i-1$ odd s.t. $R-S+2\geq d[-c_{j,j+1}]\geq
T_j-T_{j+1}+\c_j\geq T_j-T_i+\c_{i-1}$ so $\c_{i-1}\leq
T_i-T_j+R-S+2$, which implies that $B_i\leq S_{i+1}-T_i+\c_{i-1}\leq
S_{i+1}-T_i+T_i-T_j+R-S+2=-T_j+R+2$. We have $T_j\geq T_1\geq R+1$. If
the equality is strict then $B_i\leq 0$ and we are done. Otherwise
$T_j=R+1\ev T_i+1\m2$, which implies, by Lemma 7.3(ii), that the
inequality $-T_{j+1}+\c_j\geq -T_i+\c_{i-1}$ is strict. So this time
$B_i\leq S_{i+1}-T_i+\c_{i-1}<T_j-R+2=1$. If $\c_{i-1}\in\ZZ$ then
$B_i\leq 0$ so we are done. Otherwise $T_i-T_{i-1}>2e$ and
$\c_{i-1}=(T_i-T_{i-1})/2+e$ and we only have $B_i\leq
S_{i+1}-T_i+\c_{i-1}=S_{i+1}-(T_{i-1}+T_i)/2+e\leq 1/2$. But
$T_i-T_{i-1}\geq 2e+1$ so $S_{i+1}-T_i\leq
S_{i+1}-(T_{i-1}+T_i)/2-e-1/2\leq -2e$, which implies $B_i\leq
(S_{i+1}-T_i)/2+e\leq 0$ so we are done. 

7. We have type II or III, $t'\leq i<u$ and $i\ev t'\m2$. We have
$A_i=\b_i$ so, by 6.16, either
$d[b_{1,i}c_{1,i}]=d[a_{1,i}c_{1,i}]\geq C_i$ or
$d[b_{1,i}c_{1,i}]=\b_i=S_{i+1}-S_i+1$. In the first case it is enough
to prove that $C_i\geq B_i$. If
$C_i=(R_{i+1}-T_i)/2+e=(S_{i+1}-T_i+1)/2+e$ then
$C_i>(S_{i+1}-T_i)/2+e\geq B_i$. If
$C_i=R_{i+1}-T_i+d[-a_{1,i+1}c_{1,i-1}]$ then $C_i\geq
R_{i+1}-T_i=S_{i+1}-T_i+1=S_{i+1}-T_i+\b_{i+1}\geq B_i$. If
$C_i=R_{i+1}+R_{i+2}-T_{i-1}-T_i+d[a_{1,i+2}c_{1,i-2}]$ then
$R_{i+1}+R_{i+2}=S_{i+1}+S_{i+2}$. We also have $i=n$ or
$A_{i+2}=\b_{i+2}$ (both when $i+2<u$ and when $i+2=u$) so
$d[a_{1,i+2}c_{1,i-2}]\geq d[b_{1,i+2}c_{1,i-2}]$. It follows that
$C_i\geq S_{i+1}+S_{i+2}-T_{i-1}-T_i+d[b_{1,i+2}c_{1,i-2}]\geq
B_i$. Suppose now that $d[b_{1,i}c_{1,i}]=\b_i=S_{i+1}-S_i+1$ so we
have to prove that $S_{i+1}-S_i+1\geq B_i$. If $S_i\leq T_i$ then
$B_i\leq S_{i+1}-T_i+\b_{i+1}=S_{i+1}-T_i+1\leq S_{i+1}-S_i+1$ so we
are done. Thus we may assume that $T_i<S_i$. We will use this to prove
that $d[b_{1,i}c_{1,i}]=d[a_{1,i}c_{1,i}]$, a case already
discussed. Suppose first we have type II. Since $i\ev t'\m2$ we
have by Lemma 7.2(ii) $R_1+\cdots +R_i\ev iT-1\ev i(T+1)+t'-1\ev
S_1+\cdots +S_i\m2$. On the other hand $T_i\leq S_i-1=T$, which
implies that $T_1+\cdots +T_i\ev iT\m2$. (If $i$ is even then $T_i\leq
T\leq T_1$ implies that $T_1+\cdots +T_i$ is even by Lemma 6.6(ii). If
$i$ is odd it implies that $T_1=T_i=T$ so $T_1+\cdots +T_i\ev T\m2$ by
Lemma 6.6(i).) Thus $\ord a_{1,i}c_{1,i}$ and $\ord b_{1,i}c_{1,i}$
are both odd so $d[a_{1,i}c_{1,i}]=d[b_{1,i}c_{1,i}]=0$. If we have
type III then $i\ev t'\m2$ so it is even. By Lemma 7.8 we have
$d[(-1)^{i/2}a_{1,i}]=d[(-1)^{i/2}b_{1,i}]=R-S+2$ so if we prove
that $d[(-1)^{i/2}c_{1,i}]>R-S+2$ then 
$d[a_{1,i}c_{1,i}]=d[b_{1,i}c_{1,i}]=R-S+2$ and we are done. For
any $1\leq j<i$ odd we have $T_{j+1}\leq T_i\leq S_i-1=S-1$ and
$T_j\geq T_1\geq R+1$ so $T_j-T_{j+1}\geq R-S+2$. If
$T_{j+1}-T_j=-2e$ then $\c_j=0$ so $d[-c_{j,j+1}]\geq
T_j-T_{j+1}+\c_j=2e\geq R-S+3$. If $T_{j+1}-T_j>-2e$ then $\c_j\geq
1$ so $d[-c_{j,j+1}]\geq T_j-T_{j+1}+\c_j\geq T_j-T_{j+1}+1\geq
R-S+3$. By domination principle we get $d[(-1)^{i/2}c_{1,i}]\geq
R-S+3$.

8. $i\geq u$. Note that $S_1+\cdots +S_n=R_1+\cdots +R_n+2$ and
$R_k=S_k$ for $k>u$ so $S_1+\cdots +S_j=R_1+\cdots +R_j+2$ for $j\geq
u$.

We have $A_i=\b_i$ so $d[b_{1,i}c_{1,i}]=\min\{
d[a_{1,i}c_{1,i}],\b_i\}$ by 6.16. Suppose first that
$d[b_{1,i}c_{1,i}]=d[a_{1,i}c_{1,i}]\geq C_i$. So it is enough to
prove that $C_i\geq B_i$. But $A_{i+1}=\b_{i+1}$ and
$A_{i+2}=\b_{i+2}$ (or $i+1=n$ resp. $i+2=n$) so, by 6.16,
$d[-a_{1,i+1}c_{1,i-1}]\geq d[-b_{1,i+1}c_{1,i-1}]$ and
$d[a_{1,i+2}c_{1,i-2}]\geq d[b_{1,i+2}c_{1,i-2}]$. Also
$R_{i+1}=S_{i+1}$ and (if $i+2\leq n$) $R_{i+2}=S_{i+2}$. Hence
$C_i=\min\{(R_{i+1}-T_i)/2+e, R_{i+1}-T_i+d[-a_{1,i+1}c_{1,i-1}],
R_{i+1}+R_{i+2}-T_{i-1}-T_i+d[a_{1,i+2}c_{1,i-2}]\}
\geq\min\{(S_{i+1}-T_i)/2+e, S_{i+1}-T_i+d[-b_{1,i+1}c_{1,i-1}],
S_{i+1}+S_{i+2}-T_{i-1}-T_i+d[b_{1,i+2}c_{1,i-2}]\} =B_i$ and we are
done.

So we may assume that
$d[b_{1,i}c_{1,i}]=\b_i<d[a_{1,i}c_{1,i}]\leq\a_i$ and we have to
prove that $\b_i\geq B_i$.

Now $S_1+\cdots +S_i=R_1+\cdots +R_i+2$. So we may assume that
$T_1+\cdots +T_i\ev S_1+\cdots +S_i\ev R_1+\cdots +R_i\m2$ since
otherwise $\ord a_{1,i}c_{1,i}$ and $\ord b_{1,i}c_{1,i}$ are both odd
so $d[a_{1,i}c_{1,i}]=d[b_{1,i}c_{1,i}]=0$.

Since $\b_i=\min\{\a_i, S_{i+1}-S_{u+1}+\b_u\}$ and $\b_i<\a_i$ we
have $\b_i=S_{i+1}-S_{u+1}+\b_u$ so
$-S_{u+1}+\b_u=-S_{i+1}+\b_i$. Hence for any $u\leq j\leq i$ we have
$-S_{j+1}+\b_j=-S_{u+1}+\b_u$, i.e. $\b_j=S_{j+1}-S_{u+1}+\b_u$. Also,
since $\b_i<\a_i$, we have
$-S_{j+1}+\b_j=-S_{i+1}+\b_i<-R_{i+1}+\a_i\leq -R_{j+1}+\a_j$ and so
$\b_j<\a_j$. (We have $R_{i+1}=S_{i+1}$ and $R_{j+1}=S_{j+1}$.) In
particular, $\b_j<\a_j\leq (R_{j+1}-R_j)/2+e=(S_{j+1}-S_j)/2+e$ for
any $u<j\leq i$. Suppose that $\b_u=(S_{u+1}-S_u)/2+e$. If $u<i$ then
$(S_{i+1}-S_i)/2+e>\b_i=S_{i+1}-S_{u+1}+\b_u=S_{i+1}-(S_u+S_{u+1})/2+e$,
which implies that $S_u+S_{u+1}>S_i+S_{i+1}$. Contradiction
($u<i$). If $i=u$ then $\b_u=(S_{u+1}-S_u)/2+e\geq B_u$ by 2.6 and we
are done. So we may assume that $\b_j<(S_{j+1}-S_j)/2+e$ holds for
$j=u$ as well. By [B3, Lemma 2.7(iv)] this implies that $\b_j$ is
odd and $0<\b_j<2e$ for $u\leq j\leq i$. Also
$\b_j=S_{j+1}-S_j+d[-b_{j,j+1}]$.

Also $-S_{u+1}+\b_u=-S_{i+1}+\b_i$ implies $S_{u+1}\ev\ldots\ev
S_{i+1}\m2$ by Lemma 7.3(ii).

Let now $u\leq j<k\leq i+1$, $k\ev j+1\m2$. For any $j\leq l<k$, $l\ev
j\m2$ we have $d[-b_{l,l+1}]=S_l-S_{l+1}+\b_l=S_l-S_{u+1}+\b_u\geq
S_j-S_{u+1}+\b_u$. By domination principle
$d[(-1)^{(k-j+1)/2}b_{j,k}]\geq S_j-S_{u+1}+\b_u$.

Similarly if $u<j<k\leq i+1$, $k\ev j+1\m2$ then
$d[-a_{l,l+1}]=R_l-R_{l+1}+\a_l>S_l-S_{l+1}+\b_l\geq
S_j-S_{u+1}+\b_u$ for any $j\leq l<k$, $l\ev j\m2$. (We have $u<j\leq
i$ so $\a_l>\b_l$, $R_l=S_l$ and $R_{l+1}=S_{l+1}$.) So
$d[(-1)^{(k-j+1)/2}a_{j,k}]>S_j-S_{u+1}+\b_u$ by domination
principle. 

Suppose first that $S_u=R_u+1$. Since $S_{u+1}-S_u<2e$ we have
$R_{u+1}-R_u=S_{u+1}-S_u+1\leq 2e$. If $S_{u+1}-S_u$ is even so
$R_{u+1}-R_u$ is odd then
$\a_u=R_{u+1}-R_u=S_{u+1}-S_u+1\leq\b_u$. Contradiction. So
$S_{u+1}-S_u$ is odd. Thus $S_{u+1}\ev S_u+1=S+1\m2$. It follows that
$S_{u+1}\ev\ldots\ev S_{i+1}\ev S+1\m2$. Since also $S_{u+1}-S_u<2e$
we have $\b_u=S_{u+1}-S_u=S_{u+1}-S$ so for any $u\leq j\leq i$ we
have $\b_j=S_{j+1}-S_{u+1}+\b_u=S_{j+1}-S$. So we must prove that
$B_i\leq S_{i+1}-S$.

Take first the case when the type is I or II. We have $S=T+1$ so
we must prove that $B_i\leq S_{i+1}-T-1$. Also $S_{u+1}\ev\ldots\ev
S_{i+1}\ev S+1\ev T\m2$. By Lemma 7.2(ii) if we have type II then
$S_1+\cdots +S_u\ev u(T+1)+t'-1\ev uT+1\m2$. (We have $t'\ev u\m2$.)
The same happens if the type is I, when $u$ is odd and we have by
Lemma 7.2(i) $S_1+\cdots +S_u\ev u(T+1)\ev uT+1\m2$. If $u\leq j\leq
i+1$ then $S_{u+1}\ev\ldots\ev S_j\ev T\m2$ so $S_{u+1}+\cdots +S_j\ev
(j-u)T\m2$, which, together with $S_1+\cdots +S_u\ev uT+1\m2$,
implies $S_1+\cdots +S_j\ev jT+1\m2$. Also note that $S_1+\cdots
+S_j=R_1+\cdots +R_j+2$ so $R_1+\cdots +R_j\ev jT+1\m2$ as well. In
particular, if $j$ is even then $S_1+\cdots +S_j$ and $R_1+\cdots
+R_j$ are odd. 

Suppose first that $i$ is even. Then $R_1+\cdots +R_i$ and $S_1+\cdots
+S_i$ are odd and so is $T_1+\cdots +T_i$. Thus there is $1\leq j<i$
odd s.t. $T_j+T_{j+1}$ is odd. It follows that $T_i\geq
T_{j+1}>T_j\geq T_1\geq T$ so $T_i\geq T+1$. If $T_i\ev T+1\m2$ then,
since $T_1+\cdots +T_i$ is odd, we have $T_1+\cdots +T_{i-1}\ev
T\m2$. On the other hand $S_1+\cdots +S_{i+1}\ev (i+1)T+1\ev T+1\m2$
so $\ord b_{1,i+1}c_{1,i-1}$ is odd. It follows that $S_i\leq
S_{i+1}-T_i+d[-b_{1,i+1}c_{1,i-1}]=S_{i+1}-T_i\leq S_{i+1}-T-1$ and we
are done. So we may assume that $T_i\ev T\m2$. We claim that
$T_{k+1}-T_k$ is odd and $T_k>T$ for some $1\leq k<i$. Let
$1\leq j<i$ be odd s.t. $T_j+T_{j+1}$ is odd. We have $T_j\geq T_1\geq
T$. If $T_j>T$ we can take $k=j$. If $T_j=T$ Then $T_{j+1}\ev
T_j+1=T+1\ev T_i+1\m2$ so $j+1<i$. Let $j+1\leq k\leq i$ be the
largest index s.t. $T_k\ev T+1\m2$. Since $T_i\ev T\m2$ we have $k<i$
and $T_{k+1}\ev T\m2$ so $T_{k+1}-T_k$ is odd. Also $T_k\geq T_{j+1}$
or $T_k\geq T_{j+2}$, depending on the parity of $k$. Since
$T_{j+1}-T_j$ is odd we have $T_{j+1}>T_j$ and $T_{j+2}>T_j$ so
$T_k>T_j=T$. Since $\ord c_{k,k+1}=T_k+T_{k+1}$ is odd we have
$\c_{i-1}\leq T_i-T_k+d(-c_{k,k+1})=T_i-T_k\leq T_i-T-1$, which implies
that $B_i\leq S_{i+1}-T_i+\c_{i-1}\leq S_{i+1}-T-1$ so we are done.

If $i$ is odd then $i+1$ is even so $S_1+\cdots +S_{i+1}$ is odd. Also
$R_1+\cdots +R_i\ev S_1+\cdots +S_i\ev iT+1\ev T+1\m2$ and so
$T_1+\cdots +T_i\ev T+1\m2$ as well. We have $T_i\geq T_1\geq T$. If
$T_1=T_i=T$ then $T_1+\cdots +T_i\ev T\m2$ by Lemma
6.6(i). Contradiction. Thus $T_i\geq T+1$. If $T_1+\cdots +T_{i-1}$ is
even, since also $S_1+\cdots +S_{i+1}$ is odd, we have that $\ord
b_{1,i+1}c_{1,i-1}$ is odd so $B_i\leq
S_{i+1}-T_i+d[-b_{1,i+1}c_{1,i-1}]=S_{i+1}-T_i\leq S_{i+1}-T-1$ and we
are done. So we may assume that $T_1+\cdots +T_{i-1}$ is odd. Since
also $T_1+\cdots +T_i\ev T+1\m2$ we have $T_i\ev T\m2$. By a similar
reasoning as in the previous case there is $1\leq k<i$
s.t. $T_{k+1}-T_k$ is odd and $T_k>T$ and the proof follows the
same. (This time we use the fact that $T_1+\cdots +T_{i-1}$ is odd to
find a $1\leq j<i-1$ odd s.t. $T_{j+1}-T_j$ is odd.)

Suppose now that we have type III. We have $S_{u+1}\ev\ldots\ev
S_{i+1}\ev S+1\ev R\m2$ so for $u\leq j\leq i+1$ we have
$S_{u+1}+\cdots +S_j\ev (j-u)R\ev jR\m2$ ($u$ is even). Since also
$S_1+\cdots +S_u\ev uS\ev 0\m2$ by Lemma 7.2(iii) we have $S_1+\cdots
+S_j\ev jR\m2$. Since $S_1+\cdots +S_j=R_1+\cdots +R_j+2$ we also have
$R_1+\cdots +R_j\ev jR\m2$. In particular, if $j$ is even then
$R_1+\cdots +R_j$ and $S_1+\cdots +S_j$ are even. 

By Lemma 7.8 we have $d[(-1)^{u/2}a_{1,u}]=R-S+2$. Let $u<j\leq
i+1$ be even. We have
$d[(-1)^{(j-u)/2}a_{u+1,j}]>S_{u+1}-S_{u+1}+\b_u=\b_u=S_{u+1}-S$.
But $S_{u+1}\geq S_1=R+1$ and $S_{u+1}\ev S+1\ev R\m2$
so $S_{u+1}\geq R+2$. So
$d[(-1)^{(j-u)/2}a_{j+1,u}]>S_{u+1}-S\geq R-S+2$. It follows
that $d[(-1)^{j/2}a_{1,j}]=R-S+2$. Also $\b_j=S_{j+1}-S\geq
S_{u+1}-S\geq R-S+2$. Since $A_j=\b_j$ we have
$d[(-1)^{j/2}b_{1,j}]=\min\{ d[(-1)^{j/2}a_{1,j}],\b_j\}=R-S+2$. 

Suppose first that $i$ is even. We may assume that
$d[(-1)^{i/2}c_{1,i}]=R-S+2$ since otherwise
$d[(-1)^{i/2}a_{1,i}]=d[(-1)^{i/2}b_{1,i}]=R-S+2$ implies
$d[a_{1,i}c_{1,i}]=d[b_{1,i}c_{1,i}]$, a case already discussed. By
domination principle there is $1\leq j<i$ odd s.t. $R-S+2\geq
d[-c_{j,j+1}]\geq T_j-T_{j+1}+\c_j\geq T_j-T_i+\c_{i-1}$. Hence
$B_i\leq S_{i+1}-T_i+\c_{i-1}\leq S_{i+1}-T_j+R-S+2$. Suppose
that $B_i>S_{i+1}-S$. It follows that $T_j<R+2$. Since $T_j\geq T_1\geq
R+1$ we have $T_1=T_j=R+1$ so $T_1+\cdots +T_j\ev R+1\m2$ by Lemma
6.6(i). Suppose that the inequality $R-S+2\geq T_j-T_i+\c_{i-1}$ is
strict. It follows that $B_i\leq
S_{i+1}-T_i+\c_{i-1}<S_{i+1}-T_j+R-S+2=S_{i+1}-S+1$. If
$\c_{i-1}\in\ZZ$ then $B_i\leq S_{i+1}-S$ and we are done. Otherwise
$T_i-T_{i-1}$ is odd and $>2e$ so we have $B_i\leq
S_{i+1}-(T_{i-1}+T_i)/2+e<S_{i+1}-(2T_{i-1}+2e)/2+e=S_{i+1}-T_{i-1}$.
But $T_{i-1}\geq R+1\geq S$ so again $B_i\leq S_{i+1}-S$. So we may
assume that $R-S+2=T_j-T_i+\c_{i-1}$, which implies that
$-T_{j+1}+\c_j=-T_i+\c_{i-1}$ so $T_{j+1}\ev\ldots\ev T_i\m2$ by Lemma
7.3(ii). It follows that $T_{j+1}+\cdots +T_{i-1}$ is even. Since also
$T_1+\cdots +T_j\ev R+1\m2$ we have $T_1+\cdots +T_{i-1}\ev
R+1\m2$. On the other hand $S_1+\cdots +S_{i+1}\ev (i+1)R\ev
R\m2$. Thus $\ord b_{1,i+1}c_{1,i-1}$ is odd, which implies that
$B_i\leq S_{i+1}-T_i+d[-b_{1,i+1}c_{1,i-1}]=S_{i+1}-T_i$. But
$R-S+2=T_j-T_i+\c_{i-1}=R+1-T_i+\c_{i-1}$ so $T_i=S-1+\c_{i-1}\geq
S-1$. But $T_i\ev R+1\ev S\m2$. (We have $T_1+\cdots +T_{i-1}\ev
R+1\m2$ and $T_1+\cdots +T_i\ev S_1+\cdots +S_i\ev iR\ev 0\m2$.) Thus
$T_i\geq S$ and $B_i\leq S_{i+1}-T_i\leq S_{i+1}-S$ so we are done.

Suppose now that $i$ is odd. Since $i+1$ is even $S_1+\cdots +S_{i+1}$
is even. Also $R_1+\cdots +R_i\ev S_1+\cdots +S_i\ev iR\ev R\m2$ so
$T_1+\cdots +T_i\ev R\m2$ as well. We have $T_i\geq T_1\geq R+1$. If
$T_1=T_i=R+1$ then $T_1+\cdots +T_i\ev R+1\m2$ by Lemma
6.6(i). Contradiction. So we may assume that $T_i\geq R+2$. We have
$d[(-1)^{(i+1)/2}b_{1,i+1}]=R-S+2$. If $d[-b_{1,i+1}c_{1,i-1}]\leq
R-S+2$ then $B_i\leq S_{i+1}-T_i+R-S+2\leq S_{i+1}-S$ so we are
done. So we may assume that
$d[-b_{1,i+1}c_{1,i-1}]>R-S+2=d[(-1)^{(i+1)/2}b_{1,i+1}]$, which
implies $d[(-1)^{(i-1)/2}c_{1,i-1}]=R-S+2$. Same as in the previous
case $R-S+2\geq d[-c_{j,j+1}]\geq T_j-T_{j+1}+\c_j\geq
T_j-T_i+\c_{i-1}$ for some $1\leq j<i-1$ odd. Like before we get
$B_i\leq S_{i+1}-T_i+\c_{i-1}\leq S_{i+1}-T_j+R-S+2$. Suppose that
$B_i>S_{i+1}-S$ so $T_j<R+2$. Since also $T_j\geq T_1\geq R+1$ we get
$T_1=T_j=R+1$ so $T_1+\cdots +T_j\ev R+1\m2$ by Lemma 6.6(i). Suppose
that the inequality $R-S+2\geq T_j-T_i+\c_{i-1}$ is strict. If
$\c_{i-1}\in\ZZ$ then $B_i\leq S_{i+1}-S$, same as above. If
$\c_{i-1}\notin\ZZ$ then $T_i-T_{i-1}>2e$ and
$\c_{i-1}=(T_i-T_{i-1})/2+e$ so
$R-S+2>T_j-T_i+\c_{j-1}=R+1+e-(T_{i-1}+T_i)/2>R+1+e-(T_i-2e+T_i)/2
=R+1+2e-T_i$ so $T_i>S-1+2e$. We have $T_i\geq S+2e$ and
$S_{i+1}\geq S_{t'}=S$ (both $t'$ and $i+1$ are even) so
$S_{i+1}+T_i\geq 2S+2e$, which implies $B_i\leq
(S_{i+1}-T_i)/2+e\leq S_{i+1}-S$. Suppose now that
$R-S+2=T_j-T_i+\c_{i-1}$. Like in the previous case we get
$-T_{j+1}+\c_j=-T_i+\c_{i-1}$ so $T_{j+1}\ev\ldots\ev T_i\m2$. This
implies that $T_{j+1}+\cdots +T_i$ is even ($i,j$ are both odd). Since
also $T_1+\cdots +T_j\ev R+1\m2$ we have $T_1+\cdots +T_i\ev R+1\m2$.
Contradiction.

Suppose now that $S_u=R_u+2$, i.e. $M,N$ are of type I and $t'=u$. We
have $A_u=\b_u$. If $2\leq u\leq n-2$ then $S_u\leq S_{u+2}=R_{u+2}$
and $S_{u-1}<R_{u-1}\leq R_{u+1}$. ($R_{u-1}=S_{u-1}+2$ if $t<t'=u$
and $R_{u-1}=S_{u-1}+1$ if $t=u$.) Hence $S_{u-1}+S_u<R_{u+1}+R_{u+2}$
and $\b_u=A_u=\( A_u=\min\{ (R_{u+1}-S_u)/2+e,
R_{u+1}-S_u+d[-a_{1,u+1}b_{1,u-1}]\}$. Same happens if $u=1$ or
$n-1$. We have $\b_u<(S_{u+1}-S_u)/2+e=(R_{u+1}-S_u)/2+e$ so
$\b_u=R_{u+1}-S_u+d[-a_{1,u+1}b_{1,u-1}]$,
i.e. $d[-a_{1,u+1}b_{1,u-1}]=S_u-R_{u+1}+\b_u=S_u-S_{u+1}+\b_u$. We
claim that $d[(-1)^{(u-1)/2}b_{1,u-1}]>S_u-S_{u+1}+\b_u$, which, by
domination principle, will imply
$d[(-1)^{(u+1)/2}a_{1,u+1}]=S_u-S_{u+1}+\b_u$. If $u=1$ then
$d[(-1)^{(u-1)/2}b_{1,u-1}]=\j$ and we are done. Otherwise we use
Lemma 7.6. If $S_u-S_{u-1}=2e+1$ then $d[(-1)^{(u-1)/2}b_{1,u-1}]\geq
2e>S_u-S_{u+1}+\b_u$. (We have $S_{u+1}-S_u\geq -2e$ so
$S_{u+1}-S_u+2e\geq (S_{u+1}-S_u)/2+e>\b_u$.) If $S_u-S_{u-1}\leq 2e$
then $d[(-1)^{(u-1)/2}b_{1,u-1}]=\b_{u-1}$. We also have
$\b_{u-1}=\a_{u-1}+2$. (If $t<t'=u$ then this follows from Lemma
6.9(v). If $u=t$ we just note that $S_t-S_{t-1}$ is odd and $<2e$ so
$\b_{u-1}=S_u-S_{u-1}=(R_u+2)-(R_{u-1}-1)=R_u-R_{u-1}+3=\a_{u-1}+2$.)
Therefore $d[(-1)^{(u-1)/2}b_{1,u-1}]=\a_{u-1}+2\geq
R_u-R_{u+1}+\a_u+2=S_u-2-S_{u+1}+\a_u+2=S_u-S_{u+1}+\a_u>
S_u-S_{u+1}+\b_u$.

If $u<j\leq i+1$ is even, since $u+1$ is also even, we have
$d[(-1)^{(j-u-1)/2}a_{u+2,j}]>S_{u+2}-S_{u+1}+\b_u\geq
S_u-S_{u+1}+\b_u=d[(-1)^{(u+1)/2}a_{1,u+1}]$. It follows that
$d[(-1)^{j/2}a_{1,j}]=S_u-S_{u+1}+\b_u$. Since also
$\b_j=S_{j+1}-S_{u+1}+\b_u\geq S_u-S_{u+1}+\b_u$ ($j+1$ and $u$ are
both odd) we have $d[(-1)^{j/2}b_{1,j}]=\min\{
d[(-1)^{j/2}a_{1,j}],\b_j\}=S_u-S_{u+1}+\b_u$.

Suppose first that $i$ is even. We have
$d[(-1)^{i/2}a_{1,i}]=d[(-1)^{i/2}b_{1,i}]=S_u-S_{u+1}+\b_u$.
We may assume that $d[(-1)^{i/2}c_{1,i}]=S_u-S_{u+1}+\b_u$ since
otherwise $d[a_{1,i}c_{1,i}]=d[b_{1,i}c_{1,i}]$, a case already
discussed. It follows that for some $1\leq j<i$ odd we have
$S_u-S_{i+1}+\b_i=S_u-S_{u+1}+\b_u\geq d[-c_{j,j+1}]\geq
T_j-T_{j+1}+\c_j\geq T_j-T_i+\c_{j-1}$ so $B_i\leq
S_{i+1}-T_i+\c_{i-1}\leq S_u-T_j+\b_i=T+1-T_j+\b_i$. We have $T_j\geq
T_1\geq T$. If $T_j\geq T+1$ then $B_i\leq T+1-T_j+\b_i\leq\b_i$
so we are done. Otherwise $T_1=T_j=T$ so $T_1+\cdots +T_j\ev T\m2$ by
Lemma 6.6(i). If the inequality $S_u-S_{i+1}+\b_i\geq
T_j-T_i+\c_{i-1}$ is strict then $B_i\leq
S_{i+1}-T_i+\c_{i-1}<T+1-T_j+b_i=\b_i+1$. If $\c_{i-1}\in\ZZ$, since
also $\b_i\in\ZZ$, we have $B_i\leq S_{i+1}-T_i+\c_{i-1}\leq\b_i$ and
we are done. If $\c_{i-1}\notin\ZZ$ then $T_i-T_{i-1}>2e$. We have
$B_i\leq
S_{i+1}-(T_{i-1}+T_i)/2+e<S_{i+1}-(2T_{i-1}+2e)/2+e=S_{i+1}-T_{i-1}$. Suppose
that $B_i>\b_i=S_{i+1}-S_{u+1}+\b_u$. Then $T_{i-1}<S_{u+1}-\b_u$. We
have $T_{i-1}\geq T_1=T$. Now $S_{u+1}-S_u<2e$ so if $S_{u+1}-S_u$
is even then $\b_u\geq S_{u+1}-S_u+1$ so $T\leq
T_{i-1}<S_{u+1}-\b_u\leq S_u-1=T$. Contradiction. Thus $S_{u+1}-S_u$ is
odd and we have $S_{u+1}-\b_u=S_u=T+1$ so $T\leq T_{i-1}<T+1$,
i.e. $T_{i-1}=T$. We have $S_1+\cdots S_i\ev T_1+\cdots T_i\m2$,
$T_i\ev T_{i-1}+1=T+1\m2$ and $S_{i+1}\ev S_{u+1}\ev S_u+1=T+2\m2$ so
$S_1+\cdots S_{i+1}\ev T_1+\cdots T_{i-1}+1\m2$. So $\ord
b_{1,i+1}c_{1,i-1}$ is odd, which implies that $B_i\leq
S_{i+1}-T_i+d[-b_{1,i+1}c_{1,i-1}]=S_{i+1}-T_i$. But $T_i\geq
T_{i-1}+2e=T+2e>T+1=S_{u+1}-\b_u$ so $B_i\leq
S_{i+1}-S_{u+1}+\b_u=\b_i$. Suppose now that
$S_u-S_{i+1}+\b_i=T_j-T_i+\c_{i-1}$. It follows that
$-T_{j+1}+\c_j=-T_i+\c_{i-1}$ so $T_{j+1}\ev\ldots\ev T_i\m2$, by
Lemma 7.3(ii), so $T_{j+1}+\cdots +T_{i-1}$ is even. Together with
$T_1+\cdots +T_j\ev T\m2$, this implies $T_1+\cdots +T_{i-1}\ev
T\m2$. On the other hand $S_1+\cdots +S_u\ev u(T+1)\ev T+1\m2$ by
Lemma 7.2(i) and $S_{u+1}\ev\ldots\ev S_{i+1}\m2$ so $S_{u+1}+\cdots
+S_{i+1}$ is even ($u$ is odd and $i$ even). Hence $S_1+\cdots
+S_{i+1}\ev T+1\m2$ so $\ord b_{1,i+1}c_{1,i-1}$ is odd. It follows
that $B_i\leq S_{i+1}-T_i+d[-b_{1,i+1}c_{1,i-1}]=S_{i+1}-T_i$. Since
$T+1-S_{i+1}+\b_i=S_u-S_{i+1}+\b_i=T_j-T_i+\c_{i-1}=T-T_i+\c_{i-1}$
we have $B_i\leq S_{i+1}-T_i=1+\b_i-\c_{i-1}$. If $\c_{i-1}\geq 1$
then $B_i\leq\b_i$ and we are done. If $\c_{i-1}=0$ then
$S_{i+1}-T_i=1+\b_i$. But $\b_i$ is odd so $S_{i+1}\ev
T_i\m2$. Together with $S_1+\cdots +S_{i+1}\ev T_1+\cdots
+T_{i-1}+1\m2$, this implies $S_1+\cdots +S_i\ev T_1+\cdots
+T_i+1\m2$. Contradiction.

Suppose now that $i$ is odd.  We have $S_1+\cdots +S_u\ev T+1\m2$ and
$S_{u+1}\ev\ldots\ev S_i\m2$ so $S_{u+1}+\cdots +S_i$ is even ($u,i$
are both odd). Thus $T_1+\cdots +T_i\ev S_1+\cdots +S_i\ev T+1\m2$. We
have $T_i\geq T_1\geq T$. If $T_i=T_1=T$ then $T_1+\cdots +T_i\ev
T\m2$ by Lemma 6.6(i). Contradiction. So $T_i\geq T+1$. We have 
$d[(-1)^{(i+1)/2}b_{1,i+1}]=S_u-S_{u+1}+\b_u=S_u-S_{i+1}+\b_i$. If
$d[(-1)^{(i-1)/2}c_{1,i-1}]\neq S_u-S_{i+1}+\b_i$ then
$d[-b_{1,i+1}c_{1,i-1}]\leq S_u-S_{i+1}+\b_i$ so $B_i\leq
S_{i+1}-T_i+d[-b_{1,i+1}c_{1,i-1}]\leq
S_u-T_i+\b_i=T+1-T_i+\b_i\leq\b_i$. So we may assume that
$d[(-1)^{(i-1)/2}c_{1,i-1}]=S_u-S_{i+1}+\b_i$. Similarly as in the
previous case we have some $1\leq j<i-1$ odd
s.t. $S_u-S_{i+1}+\b_i=S_u-S_{u+1}+\b_u\geq d[-c_{j,j+1}]\geq
T_j-T_{j+1}+\c_j\geq T_j-T_i+\c_{i-1}$. As before we get $B_i\leq\b_i$
if $T_j>T$. Otherwise $T_1=T_j=T$ so $T_1+\cdots +T_j\ev T\m2$ by
Lemma 6.6(i). We also get $B_i\leq\b_i$ if the inequality
$S_u-S_{i+1}+\b_i\geq T_j-T_i+\c_{j-1}$ is strict and
$\c_{i-1}\in\ZZ$. (See the previous case.) If
$S_u-S_{i+1}+\b_i>T_j-T_i+\c_{j-1}$ and $\c_{i-1}\notin\ZZ$ then
$\c_{i-1}>2e$ so
$T+1-S_{i+1}+\b_i=S_u-S_{i+1}+\b_i>T_j-T_i+\c_{i-1}=T-T_i+2e$, which 
implies $S_{i+1}-T_i+2e<1+\b_i$. Since $\b_i\in\ZZ$ we have
$S_{i+1}-T_i+2e\leq\b_i$ so $B_i\leq (S_{i+1}-T_i)/2+e\leq\b_i/2<\b_i$
and we are done. If $S_u-S_{i+1}+\b_i=T_j-T_i+\c_{j-1}$ then we have
again $-T_{j+1}+\c_j=-T_i+\c_{i-1}$ so $T_{j+1}\ev\ldots\ev
T_i\m2$. Hence $T_{j+1}+\cdots +T_i$ is even ($j$ and $i$ are both
odd)m which, together with $T_1+\cdots +T_j\ev T\m2$, implies
$T_1+\cdots +T_i\ev T\m2$. Contradiction. 

{\bf Proof of 2.1(iii)} Given an index $1<i<n$ s.t. $S_{i+1}>T_{i-1}$
and $B_{i-1}+B_i>2e+S_i-T_i$ we want to prove that
$[c_1,\ldots,c_{i-1}]\rep [b_1,\ldots,b_i]$. By Lemma 1.5(i) it is
enough to prove that $[b_1,\ldots,b_{i-1}]\rep [a_1,\ldots,a_i]$,
$[c_1,\ldots,c_{i-1}]\rep [a_1,\ldots,a_i]$ and
$(a_{1,i}b_{1,i},b_{1,i-1}c_{1,i-1})_\p =1$. By Lemma 1.5(ii) it is
enough to prove that $[b_1,\ldots,b_i]\rep [a_1,\ldots,a_{i+1}]$,\\
$[c_1,\ldots,c_{i-1}]\rep [a_1,\ldots,a_i]$ and
$(a_{1,i}b_{1,i},-a_{1,i+1}c_{1,i-1})_\p =1$.

By Lemma 2.18(i) $S_{i+1}>T_{i-1}$ and $B_{i-1}+B_i>2e+S_i-T_i$
imply that $\b_i+B_{i-1}>2e$ or $\b_i+d[-b_{1,i+1}c_{1,i-1}]>2e$. We
consider the two cases.

Assume first that $\b_i+B_{i-1}>2e$. In particular,
$\b_{i-1}+\b_i>2e$.

1. $1<i<t$ or $M,N$ are of type II or III and $1<i\leq t$. Suppose
first that $i$ is odd. Then $\b_i\leq 1$ so $\b_{i-1}>2e-\b_i\geq
2e-1$ so $S_i-S_{i-1}\geq 2e$. Since $i<t$ or $i=t$ and the type is
not I we cannot have $S_i-S_{i-1}=2e+1$ (see 6.15) so we must have
$S_i-S_{i-1}=2e$, i.e. $S_{i-1}=S_i-2e=R+1-2e$. Also
$S_i-(T_{i-2}+T_{i-1})/2+e\geq B_{i-1}>2e-\b_i\geq 2e-1$ so
$T_{i-2}+T_{i-1}<2S_i-2e+2=2R-2e+4$. Thus $2T_{i-2}-2e\leq
T_{i-2}+T_{i-1}<2R-2e+4$, which implies $T_{i-2}<R+2$. Since also
$T_{i-2}\geq T_1\geq R+1$ we have $T_1=T_{i-2}=R+1$. It follows that
$R+1+T_{i-1}=T_{i-2}+T_{i-1}<2R-2e+4$ so $T_{i-1}<R-2e+3$. Hence
$T_{i-1}-T_{i-2}<(R-2e+3)-(R+1)=2-2e$, which implies that
$T_{i-1}-T_{i-2}=-2e$ so $T_{i-1}=R+1-2e$. Now $S_1=T_1=R+1$ and
$S_{i-1}=T_{i-1}=R+1-2e$ so by Lemma 7.5 we have
$[b_1,\ldots,b_{i-1}]\ap\hh\pp\cdots\pp\hh\pp
[\pi^{R+1},-\eta_1\pi^{R+1}]$ and
$[c_1,\ldots,c_{i-1}]\ap\hh\pp\cdots\pp\hh\pp
[\pi^{R+1},-\eta_2\pi^{R+1}]$ with $\eta_1,\eta_2\in\{ 1,\D\}$. In
order to prove that $[b_1,\ldots,b_i]\ap\hh\pp\cdots\pp\hh\pp
[\pi^{R+1},-\eta_1\pi^{R+1},b_i]$ represents
$[c_1,\ldots,c_{i-1}]\ap\hh\pp\cdots\pp\hh\pp
[\pi^{R+1},-\eta_2\pi^{R+1}]$ it is enough to prove that
$[-\eta_1\pi^{R+1},b_i]$ represents $-\eta_2\pi^{R+1}$, which is
equivalent to $(\eta_1\eta_2,\eta_1\pi^{R+1}b_i)_\p =1$. But this
follows from the fact that $\eta_1\eta_2$ is a power of $\D$ and $\ord
b_i=S_i=R+1$ so $\ord \eta_1\pi^{R+1}b_i=2R+2$ is even.

If $i$ is even then, with the exception of the case when $M,N$ are of
type I and $i=t-1$, we have $R+1\leq T_1\leq
T_{i-1}<S_{i+1}=R+1$. Contradiction. So we may assume that we have
type I and $i=t-1$. We have $T\leq T_1\leq T_{t-2}<S_t=T+1$ so
$T_1=T_{t-2}=T$. We have $1+\b_{t-1}\geq\b_{t-2}+\b_{t-1}>2e$ so
$\b_{t-1}>2e-1$, which implies $S_t-S_{t-1}>2e-1$. So we must have
$S_t-S_{t-1}=2e+1$ and so $S_{t-1}=T-2e$. ($S_t-S_{t-1}$ is odd and
$\leq 2e+1$. See Lemma 6.11(i) and 6.15.) Thus
$\b_{t-1}=(S_t-S_{t-1})/2+e=2e+1/2$ which, if $t\geq 5$, implies
$T-2e-(T_{t-3}+T)/2+e=S_{t-1}-(T_{t-3}+T_{t-2})/2+e\geq
B_{t-2}>2e-\b_{t-1}=-1/2$ and so $T-2e+1>T_{t-3}\geq T_2\geq
T_1-2e=T-2e$. So $T_{t-3}=T-2e$. We have $S_1=T_1=T$ and
$S_{t-1}=T_{t-3}=T-2e$ so Lemma 7.5 both $\[ b_1,\ldots,b_{t-1}\]$ and
$\[ c_1,\ldots,c_{t-3}\]$ are orthogonal sums of copies of $\h
12\pi^T\aa$ and $\h 12\pi^T\ab$. (Same happens if $t=3$, when $\[
c_1,\ldots,c_{t-3}\]$ is just the zero lattice.) It follows that $\[
b_1,\ldots,b_{t-1}\]\ap\[ c_1,\ldots,c_{t-3}\]\pp\h 12\pi^T\aa$ or$\[
c_1,\ldots,c_{t-3}\]\pp\h 12\pi^T\ab$ and so $[b_1,\ldots,b_{t-1}]\ap
[c_1,\ldots,c_{t-3}]\pp\hh$ or $[c_1,\ldots,c_{t-3}]\pp
[\pi^T,-\D\pi^T]$. So in order to prove that $[c_1,\ldots,c_{t-2}]\rep
[b_1,\ldots,b_{t-1}]$ we still need $c_{t-2}\rep\hh$,
resp. $c_{t-2}\rep [\pi^T,-\D\pi^T]$. But this follows from $\ord
c_{t-2}=T_{t-2}=T$. 

2. $t\leq i<t'$. If the type is II or III then we may assume that
$t<i<t'$. ($i=t$ was treated during the case 1.) If the type is III
then $t'=t+1$ so this case is vacuous. If the type is II then we note
that $\b_{i-1}+\b_i=1+1=2\leq 2e$. Contradiction.

So the type is I and we may assume that $i$ is odd since otherwise
$S_{i-1}=S_{i+1}=T+1$ so $\b_{i-1}+\b_i\leq 2e$. By 6.13
$A_{i-1}=\a_{i-1}\geq\b_{i-1}-2$ and $A_i=\b_i$ so
$A_{i-1}+A_i\geq\b_{i-1}+\b_i-2>2e-2=2e+R_i-S_i$. Since also
$R_{i+1}=S_{i+1}+2>S_{i-1}$ we have $[b_1,\ldots,b_{i-1}]\rep
[a_1,\ldots,a_i]$. Also $d(a_{1,i}b_{1,i})+d(b_{1,i-1}c_{1,i-1})\geq
A_i+B_{i-1}=\b_i+B_{i-1}>2e$, which implies that
$(a_{1,i}b_{1,i},b_{i-1}c_{i-1})_\p =1$. In order to prove that
$[c_1,\ldots,c_{i-1}]\rep [b_1,\ldots,b_i]$ we still have to prove
that $[c_1,\ldots,c_{i-1}]\rep [a_1,\ldots,a_i]$, i.e. that
$R_{i+1}>T_{i-1}$ and $C_{i-1}+C_i>2e+R_i-T_i$. We have
$R_{i+1}>S_{i+1}>T_{i-1}$ and, by Lemma 2.12, $C_{i-1}+C_i>2e+R_i-T_i$
is equivalent to
$d[-a_{1,i+1}c_{1,i-1}]+d[-a_{1,i}c_{1,i-2}]>2e+T_{i-1}-R_{i+1}$,
$d[-a_{1,i+1}c_{1,i-1}]>e+(T_{i-2}+T_{i-1})/2-R_{i+1}$ and
$d[-a_{1,i}c_{1,i-2}]>e+T_{i-1}-(R_{i+1}+R_{i+2})/2$. These will
follows from the corresponding conditions for $N$ and $K$. By 6.16
$A_i=\b_i$ implies $d[-a_{1,i}c_{1,i-2}]\geq d[-b_{1,i}c_{1,i-2}]$ and
we also have $d[-a_{1,i+1}c_{1,i-1}]\geq d[-b_{1,i+1}c_{1,i-1}]-2$. We
have $d[-a_{1,i+1}c_{1,i-1}]+d[-a_{1,i}c_{1,i-2}]\geq
d[-b_{1,i+1}c_{1,i-1}]+d[-b_{1,i}c_{1,i-2}]-2>
2e+T_{i-1}-S_{i+1}-2=2e+T_{i-1}-R_{i+1}$. Also
$d[-a_{1,i+1}c_{1,i-1}]\geq d[-b_{1,i+1}c_{1,i-1}]-2>
e+(T_{i-2}+T_{i-1})/2-S_{i+1}-2=e+(T_{i-2}+T_{i-1})/2-R_{i+1}$ and
$d[-a_{1,i}c_{1,i-2}]\geq d[-b_{1,i}c_{1,i-2}]>
e+T_{i-1}-(S_{i+1}+S_{i+2})/2=e+T_{i-1}-(R_{i+1}+R_{i+2})/2$. 

3. $t'\leq i<u$. Note that if $t'\leq j<u$ then
$\b_j=S_{j+1}-S_j+1\leq 2e-1$ if $j\ev t'\m2$ and $\b_j=1$ if $j\ev
t'+1\m2$ (see 6.12). Thus if $t'<i<u$, regardless of the parity of
$i$, we have $\b_{i-1}+\b_i\leq 2e-1+1=2e$. So we are left with the
case $i=t'$. If $M,N$ are of type II or III then $\b_{t'-1}=1$ and
$\b_{t'}\leq 2e-1$ so $\b_{t'-1}+\b_{t'}\leq 2e$. Hence we may assume
that we have type I. We have $A_{t'}=\b_{t'}$ and, in both cases
$t<t'$ and $1<t=t'$, we have $A_{t'-1}=\a_{t'-1}\geq\b_{t'-1}-2$ so
$A_{t'-1}+A_{t'}\geq\b_{t'-1}+\b_{t'}-2>2e-2=2e+R_{t'}-S_{t'}$. Since
also $R_{t'+1}=S_{t'+1}+1>S_{t'-1}$ we have $[b_1,\ldots,b_{t'-1}]\rep
[a_1,\ldots,a_{t'}]$. Also $d[a_{1,t'}b_{1,t'}]=\b_{t'}$, by 6.16, and
$d[b_{1,t'-1}c_{1,t'-1}]\geq B_{t'-1}$ so
$d[a_{1,t'}b_{1,t'}]+d[b_{1,t'-1}c_{1,t'-1}]\geq\b_{t'}+B_{t'-1}>2e$.
Hence $(a_{1,t'}b_{1,t'},b_{1,t'-1}c_{1,t'-1})_\p =1$. So we still
need $[c_1,\ldots,c_{t'-1}]\rep [a_1,\ldots,a_{t'}]$. Now
$R_{t'+1}=S_{t'+1}+1>T_{t'-1}$ so we still have to prove that
$C_{t'-1}+C_{t'}>2e+S_{t'}-R_{t'}$. This is equivalent by Lemma
2.12 to
$d[-a_{1,t'+1}c_{1,t'-1}]+d[-a_{1,t'}c_{t'-2}]>2e+T_{t'-1}-R_{t'+1}$,
$d[-a_{1,t'+1}c_{1,t'-1}]>e+(T_{t'-2}+T_{t'-1})/2-R_{t'+1}$ and
$d[-a_{1,t'}c_{t'-2}]>e+T_{t'-1}-(R_{t'+1}+R_{t'+2})/2$. These will
follow from the similar conditions corresponding to
$B_{t'-1}+B_{t'}>2e+S_{t'}-T_{t'}$. We have $R_{t'+1}=S_{t'+1}+1$ and
$R_{t'+2}=T=S_{t'+2}-1$, $d[-b_{1,t'+1}c_{1,t'-1}]\leq\b_{t'+1}=1\leq
d[-a_{1,t'+1}c_{1,t'-1}]+1$ and $A_{t'}=\b_{t'}$ so
$d[-a_{1,t'}c_{1,t'-2}]\geq d[-b_{1,t'}c_{1,t'-2}]$, by 6.16. Hence
$d[-a_{1,t'+1}c_{1,t'-1}]+d[-a_{1,t'}c_{1,t'-2}]\geq
d[-b_{1,t'+1}c_{1,t'-1}]-1+d[-b_{1,t'}c_{1,t'-2}]>
2e+T_{t'-1}-S_{t'+1}-1=2e+T_{t'-1}-R_{t'+1}$,
$d[-a_{1,t'+1}c_{1,t'-1}]\geq d[-b_{1,t'+1}c_{1,t'-1}]-1>
e+(T_{t'-2}+T_{t'-1})/2-S_{t'+1}-1=e+(T_{t'-2}+T_{t'-1})/2-R_{t'+1}$
and $d[-a_{1,t'}c_{1,t'-2}]\geq d[-b_{1,t'}c_{1,t'-2}]>
e+T_{t'-1}-(S_{t'+1}+S_{t'+2})/2=e+T_{t'-1}-(R_{t'+1}+R_{t'+2})/2$.

4. $i\geq u$. We have $R_{i+1}=S_{i+1}>T_{i-1}$. To prove
$C_{i-1}+C_i>2e+R_i-T_i$ we will prove the equivalent conditions from
Lemma 2.12. They will follow from the similar conditions corresponding
to $B_{i-1}+B_i>2e+S_i-T_i$. We have $R_{i+1}=S_{i+1}$,
$R_{i+2}=S_{i+2}$ and by 6.16 $A_i=\b_i$ and $A_{i+1}=\b_{i+1}$ imply
$d[-a_{1,i}c_{1,i-2}]\geq d[-b_{1,i}c_{1,i-2}]$ and
$d[-a_{1,i+1}c_{1,i-1}]\geq d[-b_{1,i+1}c_{1,i-1}]$. Hence
$d[-a_{1,i+1}c_{1,i-1}]+d[-a_{1,i}c_{1,i-2}]\geq
d[-b_{1,i+1}c_{1,i-1}]+d[-b_{1,i}c_{1,i-2}]>
2e+S_{i+1}-T_{i-1}=2e+R_{i+1}-T_{i-1}$, $d[-a_{1,i+1}c_{1,i-1}]\geq
d[-b_{1,i+1}c_{1,i-1}]>e+(T_{i-2}+T_{i-1})/2-S_{i+1}=
e+(T_{i-2}+T_{i-1})/2-R_{i+1}$ and $d[-a_{1,i}c_{1,i-2}]\geq
d[-b_{1,i}c_{1,i-2}]>e+T_{i-1}-(S_{i+1}+S_{i+2})/2=
e+T_{i-1}-(R_{i+1}+R_{i+2})/2$. Thus $C_{i-1}+C_i>2e+R_i-T_i$ and so
$[c_1,\ldots,c_{i-1}]\rep [a_1,\ldots,a_i]$. Also $A_i=\b_i$ implies
$d[a_{1,i}b_{1,i}]=\b_i$ and so
$d[a_{1,i}b_{1,i}]+d[b_{1,i-1}c_{1,i-1}]\geq\b_i+B_{i-1}>2e$. It
follows that $(a_{1,i}b_{1,i},b_{1,i-1}c_{1,i-1})_\p =1$. So we still
need $[b_1,\ldots,b_{i-1}]\rep [a_1,\ldots,a_i]$ i.e. we have to prove
that $R_{i+1}>S_{i-1}$ and $A_{i-1}+A_i>2e+R_i-S_i$. We have
$\b_{i-1}+\b_i>2e$ so $S_{i-1}<S_{i+1}=R_{i+1}$. If $i>u$ then
$A_{i-1}+A_i=\b_{i-1}+\b_i>2e=2e+R_i-S_i$ so we are left with the case
$i=u$. If $S_u=R_u+2$ then we have type I and $t'=u$. We have
$A_u=\b_u$ and, in both cases $t<t'=u$ and $t=u$, we have
$A_{u-1}=\a_{u-1}\geq\b_{u-1}-2$. Thus
$A_{u-1}+A_u\geq\b_{u-1}+\b_u-2>2e-2=2e+R_u-S_u$. If $S_u=R_u+1$ then
$S_1+\cdots +S_u=R_1+\cdots +R_u+1$. (We have $S_1+\cdots
+S_n=R_1+\cdots +R_n+2$, $S_u=R_u+1$ and $S_j=R_j$ for $j>u$.) So
$\b_{u-1}\leq 1$ and $A_{u-1}=0\geq\b_{u-1}-1$. Since also $A_u=\b_u$
we have $A_{u-1}+A_u\geq\b_{u-1}+\b_u-1>2e-1=2e+R_u-S_u$.

Suppose now that $\b_i+d[-b_{1,i+1}c_{1,i-1}]>2e$. In particular, if
$i\leq n-2$, $\b_i+\b_{i+1}>2e$. We consider several cases:

5. $i<t$. Assume frist that $i$ is odd. This implies that we have type
I and $i=t-2$ since otherwise $S_i=S_{i+2}=R+1$ so $\b_i+\b_{i+1}\leq
2e$. We have $S_t=T+1$, $S_{t-2}=T$ and $S_{t-1}\ev T\m2$. Now
$S_t-S_{t-1}$ is odd. If $S_{t-1}>T-2e$ so $S_t-S_{t-1}<2e+1$ then
$\b_{t-1}=S_t-S_{t-1}\leq 2e-1$, which, together with $\b_{t-2}\leq
1$, implies $\b_{t-2}+\b_{t-1}\leq 2e$. Contradiction. So we can
assume that $S_{t-1}=T-2e$. But $T_{t-3}\geq T_2\geq T_1-2e\geq
T-2e=S_{t-1}$, which contradicts the hypothesis. (We have
$S_{i+1}>T_{i-1}$.) 

Suppose now that $i$ is even. We can assume that we have type I and
$i=t-1$ since otherwise $S_{i+1}=R+1\leq T_1\leq T_{i-1}$, contrary to
the hypothesis. Then $T\leq T_1\leq T_{t-2}<S_t=T+1$ so
$T_1=T_{t-2}=T$. which implies by Lemma 6.6(i) that $T_1+\cdots
+T_{t-2}\ev T\m2$. Also if $t\geq 5$ then $T_{t-4}=T$. On the other
hand $S_1+\cdots +S_t\ev (T+1)t\ev T+1\m2$ by Lemma 7.2(i) so $\ord
b_{1,t}c_{1,t-2}$ is odd and so $d[-b_{1,t}c_{t-2}]=0$. Thus
$\b_{t-1}=\b_{t-1}+d[-b_{1,t}c_{1,t-2}]>2e$, which implies
$S_t-S_{t-1}>2e$ so $S_t-S_{t-1}=2e+1$ and $S_{t-1}=T-2e$. If $t=3$
then $S_1=T$ and $S_2=T-2e$ so $[b_1,b_2]\ap\hh$ or
$[\pi^T,-\D\pi^T]$. But $\ord c_1=T_1=T$ so, in both cases, $c_1\rep
[b_1,b_2]$. If $t>3$ then, by Lemma 2.12,
$0=d[-b_{1,t}c_{1,t-2}]>e+(T_{t-3}+T_{t-2})/2-S_t=
e+(T_{t-3}+T)/2-T-1$ so $T_{t-3}<T+2-2e$. It follows that
$T_{t-3}-T_{t-4}=T_{t-3}-T<2-2e$, which implies that
$T_{t-3}-T_{t-4}=-2e$ so $T_{t-3}=T-2e$. Now $S_1=T_1=T$ and
$S_{t-1}=T_{t-3}=T-2e$. But we have encountered this situation during
the proof of case 1. for $i$ even. Again we get that
$[b_1,\ldots,b_{t-1}]\ap [c_1,\ldots,c_{t-3}]\pp\hh$ or
$[c_1,\ldots,c_{t-3}]\pp [\pi^T,-\D\pi^T]$ and, since $\ord c_{t-2}=T$
we have both $c_{t-2}\rep\hh$ and $c_{t-2}\rep
[\pi^T,-\D\pi^T]$. Therefore $[c_1,\ldots,c_{t-2}]\rep
[b_1,\ldots,b_{t-1}]$. 

Suppose now that $t\leq i<t'$. We consider the three cases when $M,N$
are of type I, II or III:

6. If we have type I note first that $i$ cannot be odd since this
would imply $S_i=S_{i+2}=T+1$ so $\b_i+\b_{i+1}\leq 2e$. So $i$ is
even. We have $T\leq T_1\leq T_{i-1}<S_{i+1}=T+1$ so $T_1=T_{i-1}=T$,
which implies by Lemma 6.6(i) that $T_1+\cdots +T_{i-1}\ev T\m2$. On
the other hand by Lemma 7.2(i) we have $S_1+\cdots +S_{i+1}\ev
(i+1)(T+1)\ev T+1\m2$ so $\ord b_{1,i+1}c_{1,i-1}$ is odd and
$d[-b_{1,i+1}c_{1,i-1}]=0$. Hence $\b_i=\b_i+d[-b_{i+1}c_{1,i-1}]>2e$,
which implies that $S_{i+1}-S_i>2e$. But this is impossible. (We have
$i\neq t-1$. See 6.15.)

7. If we have type II then we may assume that $i=t'-1$ since otherwise
$\b_i+\b_{i+1}=1+1=2\leq 2e$. We have $T_{t'-2}<S_{t'}=T+1$. If $t'$
is odd then $T\leq T_1\leq T_{t'-2}$ so $T_1=T_{t'-2}=T$ and
$T_1+\cdots +T_{t'-2}\ev T\m2$. If $t'$ is even then $T_1\geq T\geq 
T_{t'-2}$ so $T_1+\cdots +T_{t'-2}$ is even. (See Lemma 6.6(i) and
(ii)) In both cases $T_1+\cdots +T_{t'-2}\ev t'T\m2$. On the other
hand $S_1+\cdots +S_{t'}\ev t'(T+1)+t'-1\ev t'T+1\m2$ by Lemma
7.2(ii). Hence $\ord b_{1,t'}c_{1,t'-2}$ is odd and
$d[-b_{1,t'}c_{1,t'-2}]=0$. It follows that
$2e<\b_{t'-1}+d[-b_{1,t'}c_{1,t'-2}]=1+0=1$. Contradiction.

8. If we have type III then $t\leq i<t'=t+1$ so $i=t$. We have
$2e<\b_t+d[-b_{1,t+1}c_{1,t-1}]=1+d[-b_{1,t+1}c_{1,t-1}]$ so
$d[-b_{1,t+1}c_{1,t-1}]>2e-1$. Since also
$d[(-1)^{(t+1)/2}b_{1,t+1}]=R-S+2\leq 2e-1$, by Lemma 7.8, we have
$d[(-1)^{(t-1)/2}c_{1,t-1}]=R-S+2$. It follows that for some $1\leq
j<t-1$ odd we have $R-S+2\geq d[-c_{j,j+1}]\geq T_j-T_{j+1}+\c_j\geq
T_j-T_{j+1}$. So $T_{j+1}-T_j\geq S-R-2$. But $T_j\geq T_1\geq R+1$
and $T_{j+1}\leq T_{t-1}<S_{t+1}=S$  so $T_{j+1}-T_j<S-R-1$. Thus
$T_{j+1}-T_j=S-R-2$, which is impossible since $S-R-2$ is odd and $\leq
-1$.

9. $t'\leq i<u$. Same as in case 3., for any $t'\leq j<u$ we have
$\b_j\leq 2e-1$ if $j\ev t'\m2$ and $\b_j\leq 1$ if $j\ev
t'+1\m2$. So if $t'\leq i<u-1$, regardless of the parity of $i$, we
have $\b_i+\b_{i+1}\leq 2e-1+1=2e$. So we are left with the case
$i=u-1$. We have
$1+d[-b_{1,u}c_{1,u-2}]=\b_{u-1}+d[-b_{1,u}c_{1,u-2}]>2e$ so
$d[-b_{1,u}c_{1,u-2}]>2e-1$.

Suppose first we have type I or II. Then $T_{u-2}<S_u=T+1$. If $u$
is odd then $T_{u-2}\geq T_1\geq T$ so $T_1=T_{u-2}=T$, which
implies that $T_1+\cdots +T_{u-2}\ev T\m2$. If $u$ is even then $T_1\geq
T\geq T_{u-2}$ implies that $T_1+\cdots +T_{u-2}$ is even. (See Lemma
6.6(i) and (ii).) So $T_1+\cdots +T_{u-2}\ev uT\m2$, regardless the
parity of $u$. On the other hand if we have type II then $S_1+\cdots
+S_u\ev u(T+1)+t'-1\ev uT+1\m2$ (we have $t'\ev u\m2$). Same happens if
we have type I, when $u$ is odd and $S_1+\cdots +S_u\ev u(T+1)\ev
uT+1\m2$. (See Lemma 7.2(i) and (ii).) Thus $S_1+\cdots +S_u\ev
T_1+\cdots +T_{u-2}+1\m2$. So $\ord b_{1,u}c_{1,u-2}$ is odd and we
have $0=d[-b_{1,u}c_{1,u-2}]>2e-1$. Contradiction.

Suppose now that the type is III. Then $u$ is even. We have
$d[(-1)^{u/2}b_{1,u}]=R-S+2\leq 2e-1$, by Lemma 7.8 and
$d[-b_{1,u}c_{1,u-2}]>2e-1$ so $d[(-1)^{(u-2)/2}c_{1,u-2}]=R-S+2$. By
domination principle there is $1\leq j<u-2$ odd s.t. $R-S+2\geq
d[-c_{j,j+1}]\geq T_j-T_{j+1}+\c_j\geq T_j-T_{j+1}$ so $T_{j+1}-T_j\geq
S-R-2$. But $T_j\geq T_1\geq R+1$ and $T_{j+1}\leq T_{u-2}<S_u=S$ so
$T_{j+1}-T_j<S-R-1$. Thus $T_{j+1}-T_j=S-R-2$. But this is
impossible since $S-R-2$ is odd and $\leq -1$.

10. $i\geq u$. We have $[c_1,\ldots,c_{i-1}]\rep [a_1,\ldots,a_i]$ by
the same reasoning as in case 4.. We have
$A_i+A_{i+1}=\b_i+\b_{i+1}>2e=2e+R_{i+1}-S_{i+1}$. Also
$\b_i+\b_{i+1}>2e$ implies $S_i<S_{i+2}=R_{i+2}$. Hence
$[b_1,\ldots,b_i]\rep [a_1,\ldots,a_{i+1}]$. Finally, $A_i=\b_i$ so
$d[a_{1,i}b_{1,i}]=\b_i$ and $A_{i+1}=\b_{i+1}$ so
$d[-a_{1,i+1}c_{1,i-1}]\geq d[-b_{1,i+1}c_{1,i-1}]$. (See 6.16.) Hence
$d[a_{1,i}b_{1,i}]+d[-a_{1,i+1}c_{1,i-1}]\geq
\b_i+d[-b_{1,i+1}c_{1,i-1}]>2e$, which implies
$(a_{1,i}b_{1,i},-a_{1,i+1}c_{1,i-1})_\p =1$.

{\bf Proof of 2.1(iv)} Let $1<i\leq n-2$ s.t. $T_i\geq
S_{i+2}>T_{i-1}+2e\geq S_{i+1}+2e$. We want to ptove that
$[c_1,\ldots,c_{i-1}]\rep [b_1,\ldots,b_{i+1}]$. 

Since $S_{i+2}>S_{i+1}+2e$ we have either $i+1\geq u$ or we have type
I, $i+2=t$ and $S_t-S_{t-1}=2e+1$, i.e. $S_{t-1}=T-2e$.

If $i+1\geq u$ then $R_{i+2}=S_{i+2}>S_{i+1}+2e$ so
$[b_1,\ldots,b_{i+1}]\ap [a_1,\ldots,a_{i+1}]$ by Lemma 2.19. Also by
Lemma 2.19 $R_{i+2}=S_{i+2}>T_{i-1}+2e$ implies
$[c_1,\ldots,c_{i-1}]\rep [a_1,\ldots,a_{i+1}]\ap
[b_1,\ldots,b_{i+1}]$. 

In the other case we have $i=t-2$ so $T_{t-2}\geq S_t>T_{t-3}+2e\geq
S_{t-1}+2e$ and we want to prove that $[c_1,\ldots,c_{t-3}]\rep
[b_1,\ldots,b_{t-1}]$. Now $T+1=S_t>T_{t-3}+2e\geq S_{t-1}+2e=T$ so
$T_{t-3}=T-2e$. We have $T\leq T_1\leq T_2+2e\leq T_{t-3}+2e=T$ so
$T_1=T_{t-4}=T$. Since $S_1=T_1=T$ and $S_{t-1}=T_{t-3}=T-2e$, by
Lemma 7.5, both $\[ b_1,\ldots,b_{t-1}\]$ and $\[
c_1,\ldots,c_{t-3}\]$ are orthogonal sums of copies of $\h 12\pi^T\aa$
and $\h 12\pi^T\ab$. Since $\aa\pp\aa\ap\ab\pp\ab$ this implies $\[
c_1,\ldots,c_{t-3}\]\rep\[ b_1,\ldots,b_{t-1}\]$ so
$[c_1,\ldots,c_{t-3}]\rep [b_1,\ldots,b_{t-1}]$. \qed

\subsection{The case when $R_2-R_1=-2e$}

We now assume that $R_2-R_1=-2e$, i.e. that $\nn M=2\ss M$. If
$R_1<R_3$ and $a_2/a_1\in -\h 14\D\os$ then $M':=\{ x\in M\mid x\text{
not a norm generator}\}$ is still a lattice but this time
$[M:M']=\p^2$. In all the other cases $M'$ is no longer a lattice so
we cannot take $N=M'$. 

\blm Let $x_1,\ldots,x_t$ and $x_{t+1},\ldots,x_n$ be two good BONGs
with $\ord Q(x_i)=R_i$. Let $1\leq s\leq t<u\leq n$ and let $\[
x_s,\ldots,x_t\]\pp\[ x_{t+1},\ldots,x_u\] =\[ y_s,\ldots,y_u\]$. If
$R_{s-1}\leq R_{t+1}$, $R_t\leq R_{u+1}$ and
$x_1,\ldots,x_{s-1},y_s,\ldots,y_u,x_{u+1},\ldots,x_n$ is a good BONG
then $\[ x_1,\ldots,x_t\]\pp\[ x_{t+1},\ldots,x_n\] =\[
x_1,\ldots,x_{s-1},y_s,\ldots,y_u,x_{u+1},\ldots,x_n\]$.

(We ignore the condition $R_{s-1}\leq R_{t+1}$ if $s=1$ and we ignore
$R_t\leq R_{u+1}$ if $u=n$.)
\elm
\pf Let $S_i=\ord Q(y_i)$.  We consider first the case $u=n$. We may
assume that $s>1$ since otherwise the statement is trivial. We have
$\[ x_{t+1},\ldots,x_n\]\sbq\[ y_s,\ldots,y_n\]$ so $R_{t+1}=\ord\nn\[
x_{t+1},\ldots,x_n\]\geq\ord\nn\[ y_s,\ldots,y_n\] =S_s$. Since the
BONG $x_1,\ldots,x_{s-1},y_s,\ldots,y_n$ is good we have $R_{s-2}\leq
S_s\leq R_{t+1}$ (if $s\geq 3$). By hypothesis we also have
$R_{s-1}\leq R_{t+1}$ so $R_i\leq R_{t+1}$ for any $i\leq s-1$. (We
have $R_i\leq R_{s-2}$ or $R_i\leq R_{s-1}$, depending on the parity
of $i$.) We prove by decreasing induction that $\[
x_i,\ldots,x_t\]\pp\[ x_{t+1},\ldots,x_n\] =\[
x_i,\ldots,x_{s-1},y_s,\ldots,y_n\]$ for $1\leq i\leq s$. For $i=s$ we
have $\[ x_s,\ldots,x_t\]\pp\[ x_{t+1},\ldots,x_n\] =\[
y_s,\ldots,y_n\]$ by hypothesis. Let now $i<s$. We have $R_i\leq
R_{t+1}$ so $\nn (\[ x_i,\ldots,x_t\]\pp\[ x_{t+1},\ldots,x_n\]
)=\p^{R_i}=Q(x_i)\oo$ so $x_i$ is a norm generator. By the induction
hypothesis $pr_{x_i^\pp}(\[ x_i,\ldots,x_t\]\pp\[ x_{t+1},\ldots,x_n\]
)=\[ x_{i+1},\ldots,x_t\]\pp\[ x_{t+1},\ldots,x_n\]=\[
x_{i+1},\ldots,x_{t-1},y_t,\ldots,y_n\]$ so $\[ x_i,\ldots,x_t\]\pp\[
x_{t+1},\ldots,x_n\] =\[ x_i,\ldots,x_{t-1}y_t,\ldots,y_n\]$. When
$i=1$ we get $\[ x_1,\ldots,x_t\]\pp\[ x_{t+1},\ldots,x_n\] =\[
x_1,\ldots,x_{t-1},y_t,\ldots,y_n\]$.

For the general case we use duality. Since
$x_1,\ldots,x_{s-1},y_s,\ldots,y_u,x_{u+1},\ldots,x_n$ is a good BONG,
so is $y_s,\ldots,y_u,x_{u+1},\ldots,x_u$ and thus
$x_n^\*,\ldots,x_{u+1}^\*,y_u^\*,\ldots,y_s^\*$ also. Since $\[
x_s,\ldots,x_t\]\pp\[ x_{t+1},\ldots,x_u\] =\[ y_s,\ldots,y_u\]$ we
have $\[ x_u^\*,\ldots,x_{t+1}^\*\]\pp\[ x_t^\*,\ldots,x_s^\*\] =\[
y_u^\*,\ldots,y_s^\*\]$. Since $\ord Q(x_{u+1}^\* )=-R_{u+1}\leq
-R_t=\ord Q(x_t^\* )$ we have, by the case $u=n$ proved above, that
$\[ x_n^\*,\ldots,x_{t+1}^\*\]\pp\[ x_t^\*,\ldots,x_s^\*\] =\[
x_n^\*,\ldots,x_{u+1}^\*,y_u^\*,\ldots,y_s^\*\]$. By duality $\[
x_s,\ldots,x_t\]\pp\[ x_{t+1},\ldots,x_n\] =\[
y_s,\ldots,y_u,x_{u+1},\ldots,x_n\]$. Since also $R_{s-1}\leq R_{t+1}$
and $x_1,\ldots,x_{s-1}y_s,\ldots,y_u,x_{u+1},\ldots,x_n$ is a good
BONG we have, again by the case $u=n$, $\[ x_1,\ldots,x_t\]\pp\[
x_{t+1},\ldots,x_n\] =\[
x_1,\ldots,x_{s-1},y_s,\ldots,y_u,x_{u+1},\ldots,x_n\]$. \qed

\blm Suppose that $M,N$ are two ternary lattices over the same
quadratic space, $R_1=R_3$ and $R_i=S_i$ for $1\leq i\leq 3$. Then
$M\ap N$ iff $\a_1=\b_1$.
\elm
\pf We use [B3, Theorem 3.1]. The necessity of $\a_1=\b_1$ follows
from (ii). Conversely, suppose that $\a_1=\b_1$. Condition (i) follows
from the hypothesis. Since $R_1=R_3$ we have $R_1+\a_1=R_2+\a_2$ so
$\a_2=R_1-R_2+\a_1$ (see [B3, Corollary 2.3(i)]). Similarly
$\b_2=S_1-S_2+\b_1=R_1-R_2+\a_1$ so $\a_2=\b_2$. Thus (ii) holds.

We prove now (iii). By Lemma 7.4(iii) $R_1=R_3$ implies
$d[-a_{1,2}]=\a_2$ and $d[-a_{2,3}]=\a_1$. Similarly
$d[-b_{1,2}]=\b_2=\a_2$ and $d[-b_{2,3}]=\b_1=\a_1$. Now
$d[a_{1,3}b_{1,3}]=\j$ (we have $a_{1,3}b_{1,3}\in\fs$) and
$d[-a_{2,3}],d[-b_{2,3}]\geq\a_1$. By domination principle
$d(a_1b_1)\geq d[a_1b_1]\geq\a_1$. Also
$d[-a_{1,2}]=\a_2$ and $d[-b_{1,2}]=\b_2=\a_2$ so
$d(a_{1,2}b_{1,2})\geq d[a_{1,2}b_{1,2}]\geq\a_2$. So (iii) holds.

Finally, $R_1=R_3$ implies $\a_1+\a_2\leq 2e$ so (iv) is vacuous. \qed

\blm If $a\in\ff$ we have:

(i) $\h 12\pi a\aa\pp\la a\ra\ap\[ a\e,-\pi^{2-2e}a\e,a\]$ relative to
some good BONG for any $\e\in\ooo$.

(ii) $\h 12\pi a\ab\pp\la a\ra\ap\[ a\e,-\h\D 4\pi^2a\e\eta,a\eta\]$
relative to some good BONG for any $\e,\eta\in\ooo$ with $d(\e
)=1$, $d(\eta )=2e-1$ and $(\e,\eta )_\p=-1$. 
\elm
\pf By scaling we may assume that $a=1$. 

Let $M$ be the lattice from (i) or (ii). We have $M=M_1\pp M_2$, where
$M_1\ap\h 12\pi\aa$ or $\h 12\pi\aa$ and $M_2\ap\la 1\ra$. We want to
prove that $M\ap N$, where $N\ap\[ b_1,b_2,b_3\]$, with $b_1,b_1,b_3$
being $\e,-\pi^{2-2e}\e,1$ and $\e,\h\D 4\pi^2\e\eta,\eta$,
respectively. 

First we prove that a lattice $N\ap\[ b_1,b_2,b_3\]$ exists. We use
[B1, Lemmas 4.3(ii) and 3.6]. If $S_i=\ord b_i$ then
$(S_1,S_2,S_3)=(0,2-2e,0)$. Since $S_1\leq S_3$ we still have to prove
that $b_2/b_1,b_3/b_2\in\aaa$. Since $S_3-S_2=2e-2\geq 0$ we have
$b_3/b_2\in\aaa$. For $b_2/b_1$ we have $S_2-S_1+2e=2\geq 0$ and
$S_2-S_1+d(-b_1b_2)=2-2e+d(-b_1b_2)\geq 0$. (In the case (i) we have
$d(-b_1b_2)=\j$ and in the case (ii) we have $d(\D)=2e$ and $d(\eta
)=2e-1$ so $d(-b_1b_2)=d(\D\eta )=2e-1$.)

Now we prove that $FM=FN$. We have $\det FM=\det FN=-1$ or $-\D$
corresponding to (i) and (ii), respectively. We still have to prove
that $FM$ are either both isotropic or both anisotropic. In the case
(i) it is obvious that both $FM$ and $FN\ap [\e,-\e,1]$ are
isotropic. In the case (ii) $FM\ap [\pi,-\pi\D,1]$ is anisotropic and
$FN=[b_1.b_2,b_3]$ is aisotropic because $(-b_1b_2,-b_2,b_3)_\p
=(\D\eta,\D\e )_\p =-1$. (We have $(\D,\D )_\p =(\D,\e )_\p =(\D,\eta
)_\p =1$ and $(\e,\eta )_\p =-1$.)

We now prove that $R_i:=R_i(M)=S_i$. If $e>1$ then $M=M_1\pp M_2$ is a
Jordan splitting with $\ss_k:=\ss M_k=\p^{r_k}$, where $r_1=1-e$ and
$r_2=0$. Let $\nn M^{\ss_k}=u_k$. We have $\nn M_1=\p$ and $\nn
M_2=\oo$ so $\nn M=\oo$ so $r_1=0$. Thus $R_1=u_1=0$ and
$R_2=2u_1-r_1=2-2e$. Also $M_2$ is unary so $R_3=r_2=0$. If $e=1$ then
$1-e=0$ so $FM$ is unimodular and of odd rank. Thus
$R_1=R_2=R_3=0$. In both cases
$(R_1,R_2,R_3)=(0,2-2e,0)=(S_1,S_2,S_3).$.

Since $R_1=R_3$ and $R_i=S_i$ for all $i$ we may apply Lemma 7.11. But
$R_2-R_1=S_2-S_1=2-2e$ which uniquely defines $\a_1,\b_1$, namely
$\a_1=\b_1=1$.(See [B3, Corollary 2.9(i)].) \qed

\bco If $\ord b=S$ then $\[
\pi^S,-\pi^{S-2e+2},\ldots,\pi^S,-\pi^{S-2e+2},b\]\ap\h
12\pi^{S+1}\aa\pp\cdots\pp\h 12\pi^{S+1}\aa\pp\la b\ra$. 
\eco 
\pf Let $b_1,\ldots,b_{2s+1}$ be the sequence
$\pi^S,-\pi^{S-2e+2},\ldots,\pi^S,-\pi^{S-2e+2},b$ and let $S_i=\ord
b_i$. We have $S_i=S_{i+2}$ for $1\leq i\leq 2s-1$. Also for
$1\leq i\leq 2s-1$ we have $b_{i+1}/b_i=-\pi^{2-2e}$ or
$-\pi^{2e-2}$ so $b_{i+1}/b_i\in\aaa$. Same happens when $i=2s$
since $\ord b_{2s+1}/b_{2s}=2e-2\geq 0$. So there is a lattice
$L\ap\[ \pi^S,-\pi^{S-2e+2},\ldots,\pi^S,-\pi^{S-2e+2},b\]$ relative
to a good BONG. Since $S_{2i}=S-2e+2\leq S=S_{2i+1}$ for $i\leq
i\leq s$ we have
$L\ap\[\pi^S,-\pi^{S-2e+2}\]\pp\cdots\pp\[\pi^S,-\pi^{S-2e+2}\]\pp\la
b\ra$ by [B1, Corollary 4.4(i)]. To prove that $L\ap \h
12\pi^{S+1}\aa\pp\cdots\pp\h 12\pi^{S+1}\aa\pp\la b\ra$ it is enough
to show that $\[\pi^S,-\pi^{S-2e+2}\]\pp <b>\ap\h 12\pi^{S+1}\aa\pp\la
b\ra$. (The general case follows easily by induction.) But $\ord
b=S$ so by Lemma 7.12(i) we have $\h 12\pi^{S+1}\aa\pp\la b\ra\ap\h
12\pi b\aa\pp\la b\ra\ap\[ b\e,-\pi^{2-2e}b\e,b\]$ for any
$\e\in\ooo$. But $\ord b=S$ so $b\e =\pi^S$ for some $\e\in\ooo$.
It follows that $\h 12\pi^{S+1}\aa\pp\la
b\ra\ap\[\pi^S,-\pi^{S-2e+2},b\]$. But $\ord\pi^{S-2e+2}=S-2e+2\leq
S=\ord b$ so
$\[\pi^S,-\pi^{S-2e+2},b\]\ap\[\pi^s,-\pi^{S-2e+2}\]\pp\la b\ra$ and
we are done. \qed

\blm Suppose that $a_2/a_1\in -\h\D 4\ooo^2$ and $R_3>R_1$ or
$n=2$. We denote $R_1=R$. Let $s\geq 2$ be the smallest even integer
such that $R_{s+2}>R-2e+1$ or $s+2>n$. Let $N=M':=\{ x\in M\mid x$ not
a norm generator$\}$. We have:

(i) If $s\geq 4$ then $R_3=R_5=\ldots=R_{s-1}=R+1$ and
$R_4=R_6=\ldots=R_s=R-2e+1$.

(ii) $N$ is a lattice and we have:

1. If $R_{s+1}\geq R+2$ or $s=n$ then $N\ap\[
a_3,\ldots,a_s,\pi^2a_1,\pi^2a_2,a_{s+1},\ldots,a_n\]$ relative to the
good BONG $x_3,\ldots,x_s,\pi x_1,\pi x_2,x_{s+1},\ldots,x_n$.

2. If $R_{s+1}=R+1$ then for any $\e,\eta\in\ooo$ with $d(\e )=1$,
$d(\eta )=2e-1$ s.t. $(\e,\eta )_\p =-1$ we have $N\ap\[
a_3,\ldots,a_s,a_{s+1}\e,-\pi^{2-2e}a_{s+1}\e\eta\D,
a_{s+1}\eta,a_{s+2}\ldots,a_n\]$ relative to some good BONG
$x_3,\ldots,x_s,y_{s-1},y_s,y_{s+1},x_{s+2},\ldots,x_n$.
\elm
\pf (i) We have $R_3\geq R_1+1=R+1$ so $R_4\geq R_3-2e\geq R-2e+1$. So
$R-2e+1\leq R_4\leq R_6\leq\ldots\leq R_s$. But by the minimality of
$s$ we have $R_s\leq R-2e+1$ so $R_4=R_6=\ldots=R_s=R-2e+1$. For
$3\leq i<s$ odd we have $R+1\leq R_3\leq R_i\leq R_{i+1}+2e=R+1$. Thus
$R_3=R_5=\ldots=R_{s-1}=R+1$.

(ii) Since $R_3>R-2e=R_2$ we have the splitting $L=J\pp K$, with $J=\[
x_1,x_2\]$ and $K=\[ x_3,\ldots,x_n\]$. Since $\nn J=R_1<R_3=\nn K$ we
have $N=M'=J'\pp K$, where $J':=\{ x\in J\mid x\text{ not a norm
generator}\}$. But $a_2/a_1\in -\h\D 4\ooo^2$ so $J\ap\h 12a_1\ab\ap\h
12\pi^R\ab$, which implies that $J'=\p J=\[\pi x_1,\pi x_2\]\ap\h
12\pi^{R+2}\ab$. (This happens because every primitive element of
$\ab$ is a norm generator.) Since $R_s=R-2e+1<R+1\leq R_{s+1}$ (if
$2<s<n$) we have $K=\[ x_3,\ldots,x_s\]\pp\[ x_{s+1},\ldots,x_n\]$. If
$s>2$ then $R_s=R-2e+1<R+2=\ord Q(\pi x_1)$, $R_{s-1}=R+1<R+2=\ord
Q(\pi x_1)$ and $R_s=R-2e+1<R-2e+2=\ord Q(\pi x_2)$. By [B1, Corollary
4.4(v)] we get that $x_3,\ldots,x_s,\pi x_1,\pi x_2$ is a good BONG
for $\[ x_3,\ldots,x_s\]\pp\[\pi x_1,\pi x_2\]$. Same happens if
$s=2$, when $x_3,\ldots,x_s$ are ignored. So $M=\p J\pp K=\[\pi
x_1,\pi x_2\]\pp\[ x_3,\ldots,x_s\]\pp\[ x_{s+1},\ldots,x_n\]
=\[x_3,\ldots,x_s,\pi x_1,\pi x_2\]\pp\[ x_{s+1},\ldots,x_n\]$.

If $R_{s+1}\geq R+2$ then $\ord Q(\pi x_2)=R-2e+2\leq R_{s+2}$ (if
$s+2\leq n$), $\ord Q(\pi x_1)=R+2\leq R_{s+1}$ and $\ord Q(\pi
x_2)=R-2e+2<R+2\leq R_{s+1}$ (if $s+1\leq n$) so
$N=\[x_3,\ldots,x_s,\pi x_1,\pi x_2\]\pp\[ x_{s+1},\ldots,x_n\]
=\[x_3,\ldots,x_s,\pi x_1,\pi x_2,x_{s+1},\ldots,x_n\]$ by [B1,
Corollary 4.4(v)].

If $R_{s+1}=R+1$ then, by Lemma 7.12(ii), for any $\e,\eta\in\ooo$
with $d(\e )=1$, $d(\eta )=2e-1$ and $(\e,\eta )_\p =-1$ we have
$\[\pi x_1,\pi x_2\]\pp\oo x_{s+1}\ap\h 12\pi^{R+2}\ab\pp\la
a_{s+1}\ra\ap\h 12\pi a_{s+1}\ab\pp\la a_{s+1}\ra\ap\[
a_{s+1}\e,-\pi^{2-2e}a_{s+1}\e\eta\D,a_{s+1}\eta\]$ relative to some
good BONG $y_{s-1},y_s,y_{s+1}$. (We have $\ord a_{s+1}=R_{s+1}=R+1$.)
Let $S_i=\ord Q(y_i)$ so $S_{s-1}=S_{s+1}=R+1$ and $S_s=R-2e+3$. To
prove that $\[ x_3,\ldots,x_s,\pi x_1,\pi x_2\]\pp\[
x_{s+1},\ldots,x_n\] =\[
x_3,\ldots,x_s,y_{s-1},y_s,y_{s+1},x_{s+2},\ldots,x_n\]$ we use Lemma
7.10. We have to show that
$x_3,\ldots,x_s,y_{s-1},y_s,y_{s+1},x_{s+2},\ldots,x_n$ is a good
BONG, $R_s\leq R_{s+1}$ (if $2<s<n$) and $\ord Q(\pi x_2)\leq R_{s+2}$
(if $s+2\leq n$).

If $s>2$ then $R_{s-1}=R+1=S_{s-1}$ and $R_s=R-2e+1<R-2e+3=S_s$.
If $s+2\leq n$ then $R_{s+2}\geq R-2e+2$ and we cannot have
equality since this would imply
$R_{s+2}-R_{s+1}=(R-2e+2)-(R+1)=1-2e$, which is odd and negative.
So $S_s=R-2e+3\leq R_{s+2}$ and we also have
$S_{s+1}=R+1=R_{s+1}\leq R_{s+3}$ (if $s+2\leq n$ resp. $s+3\leq
n$). These inequalities, together with the fact that
$x_3,\ldots,x_s$, $~~y_{s-1},y_s,y_{s+1}$ and $x_{s+2},\ldots,x_n$ are
good BONGs and with the fact that
$Q(y_{s-1})/Q(x_s),Q(x_{s+2})/Q(y_{s+1})\in\aaa$, will imply that
$x_3,\ldots,x_s,y_{s-1},y_s,y_{s+1},x_{s+2},\ldots,x_n$ is a good
BONG. Now $\ord Q(y_{s-1})/Q(x_s)=S_{s-1}-R_s=(R+1)-(R-2e+1)=2e\geq 0$
so $Q(y_{s-1})/Q(x_s)\in\aaa$ (if $s>2$). We also have
$Q(x_{s+2})/Q(y_{s+1})=\eta\1 a_{s+2}/a_{s+1}$. Now
$a_{s+2}/a_{s+1}\in\aaa$ so
$R_{s+2}-R_{s+1}+d(-a_{s+1}a_{s+2})\geq 0$. But $d(\eta )=2e-1$ so
$R_{s+2}-R_{s+1}+d(\eta )\geq (R-2e+3)-(R+1)+2e-1=1>0$. By
domination principle we get $R_{s+2}-R_{s+1}+d(-\eta\1
a_{s+1}a_{s+2})\geq 0$. Since also $R_{s+2}-R_{s+1}+2e\geq 0$ we have
$\eta\1 a_{s+2}/a_{s+1}\in\aaa$.

Finally, $R_s=R-2e+1<R+1=R_{s+1}$ (if $2<s$) and $\ord
Q(\pi x_2)=R-2e+2<R-2e+3\leq R_{s+2}$ (if $s+2\leq n$) so we are
done. \qed 

In the context of Lemma 7.14 we say that $M$ is of type I if
$R_{s+1}\geq R+2$ or $s=n$ and it is of type II if
$R_{s+1}=R+1$. (I.e. if $M$ is in case 1. resp. 2. of Lemma 7.14(ii).)

We have $N\ap\[ b_1,\ldots,b_n\]$ relative to a good BONG where the
sequence $b_1,\ldots,b_n$ is
$a_3,\ldots,a_s,\pi^2a_1,\pi^2a_2,a_{s+1},\ldots,a_n$, if we have
case I, resp. \\$a_3,\ldots,a_s,a_{s+1}\e,
-\pi^{2-2e}a_{s+1}\e\eta\D,a_{s+1}\eta,a_{s+2},\ldots,a_n$ if the case
is II. 

\blm With the notations of Lemma 7.14 we have $R_i=S_i$, $\a_i=\b_i$
and $[a_1,\ldots,a_i]\ap [b_1,\ldots,b_i]$ for any $i\geq s+1$. If $M$
is of type I then $[a_1,\ldots,a_i]\ap [b_1,\ldots,b_i]$ also holds
for $i=s$. 
\elm
\pf We have $b_{s+1}=a_{s+1}$ in the case 1. and $b_{s+1}=\eta
a_{s+1}$ in the case 2. so $S_{s+1}=R_{s+1}$ is both cases. Also
$a_i=b_i$ so $R_i=S_i$ for $i\geq s+2$. If $i\geq s+1$ (or $i=s$ in
the case 1.) we have $a_j=b_j$ for $j\geq i+1$ so
$[a_{i+1},\ldots,a_n]\ap [b_{i+1},\ldots,b_n]$, which, together with
$[a_1,\ldots,a_n]\ap FM=FN\ap [b_1,\ldots,b_n]$, implies
$[a_1,\ldots,a_i]\ap [b_1,\ldots,b_i]$.

Let $K=\[ x_{s+2},\ldots,x_n\]\ap\[ a_{s+2},\ldots,a_n\]\ap\[
b_{s+2},\ldots,b_n\]$. For $i>s+1$ we have by [B3, Lemmas 2.1 and
2.4] $\a_i=\min\{\a_{i-s-1}(K),R_{i+1}-R_{s+2}+a_{s+1}\}$.
($\a_{i-s-1}(K)$ replaces $(R_{i+1}-R_i)/2+e$ and the terms
corresponding to indices $s+2\leq j\leq n-1$.
$R_{i+1}-R_{s+2}+a_{s+1}$ replaces the terms with $1\leq j\leq s+1$.)
Similarly $\b_i=\min\{\a_{i-s-1}(K),S_{i+1}-S_{s+2}+b_{s+1}\}$.
Since $R_{i+1}=S_{i+1}$ and $R_{s+2}=S_{s+2}$ it is enough to
prove that $\a_{s+1}=\b_{s+1}$. 

We have $\a_{s+1}=\min\{ (R_{s+2}-R_{s+1})/2+e,
R_{s+2}-R_{s+1}+d(-a_{s+1,s+2}), R_{s+2}-R_{s+1}+\a_1(K),
R_{s+2}-R_{s+1}+\a_s\}$. (By [B3, Lemmas 2.1 and 2.4]
$R_{s+2}-R_{s+1}+\a_1(K)$ replaces the terms with $j\geq s+2$ and
$R_{s+2}-R_{s+1}+\a_s$ replaces the terms with $j\leq s$.) Similarly
$\b_{s+1}=\min\{ (S_{s+2}-S_{s+1})/2+e,
S_{s+2}-S_{s+1}+d(-b_{s+1,s+2}), S_{s+2}-S_{s+1}+\a_1(K),
S_{s+2}-S_{s+1}+\b_s\} =\min\{ (R_{s+2}-R_{s+1})/2+e,
R_{s+2}-R_{s+1}+d(-b_{s+1,s+2}), R_{s+2}-R_{s+1}+\a_1(K),
R_{s+2}-R_{s+1}+\b_s\}$. Note that $(R_{s+2}-R_{s+1})/2+e$ and
$R_{s+2}-R_{s+1}+\a_1(K)$ appear in the sets defining both $\a_{s+1}$
and $\b_{s+1}$. 

If $M$ is of type I then $R_{s+1}-R_s\geq (R+2)-(R-2e+1)>2e$
($R_s=R-2e$ if $s=2$ and $R_s=R-2e+1$ if $s>2$) and
$S_{s+1}-S_s\geq (R+2)-(R-2e+2)=2e$. It follows that
$R_{s+2}-R_{s+1}+\a_s, R_{s+2}-R_{s+1}+\b_s\geq R_{s+2}-R_{s+1}+2e\geq
(R_{s+2}-R_{s+1})/2+e$ (we have $R_{s+2}-R_{s+1}\geq -2e$). Thus
$R_{s+2}-R_{s+1}+\a_s$ and $R_{s+2}-R_{s+1}+\b_s$ can be ignored. Also
in this case $a_{s+1}=b_{s+1}$ so
$d(-a_{s+1,s+2})=d(-b_{s+1,s+2})$. Hence $\a_{s+1}=\b_{s+1}$.

If $M$ is of type II then $R_{s+1}-R_s=R+1-R_s\geq 2e$ ($R_s=R-2e$ or
$R-2e+1$) and $S_{s+1}-S_s=(R+1)-(R-2e+3)=2e-2$. Thus
$R_{s+2}-R_{s+1}+\a_s, R_{s+2}-R_{s+1}+\b_s\geq
R_{s+2}-R_{s+1}+2e-1\geq (R_{s+2}-R_{s+1})/2+e$. (We have
$R_{s+2}-R_{s+1}\geq (R-2e+2)-(R+1)=1-2e$ so $R_{s+2}-R_{s+1}\geq
2-2e$.) Thus $R_{s+2}-R_{s+1}+\a_s$ and $R_{s+2}-R_{s+1}+\b_s$ can be
ignored. We also have $d(-b_{s+1,s+2})=d(-\eta a_{s+1,s+2})$. If
$d(-a_{s+1,s+2})<2e-1=d(\eta )$ then
$d(-a_{s+1,s+2})=d(-b_{s+1,s+2})$ so we are done. Otherwise we have
$d(-a_{s+1,s+2}),d(-b_{s+1,s+2})\geq 2e-1$ so
$R_{s+2}-R_{s+1}+d(-a_{s+1,s+2}), R_{s+2}-R_{s+1}+d(-b_{s+1,s+2})\geq
R_{s+2}-R_{s+1}+2e-1\geq (R_{s+2}-R_{s+1})/2+e$ so
$R_{s+2}-R_{s+1}+d(-a_{s+1,s+2})$ and
$R_{s+2}-R_{s+1}+d(-b_{s+1,s+2})$ can be ignored and we have again
$\a_{s+1}=\b_{s+1}$. \qed

\blm Let $M,N$ be like in Lemma 7.14. If $K\leq M$ and $\nn K\sb\nn M$
then $K\leq N$.
\elm
\pf The condition $\nn K\sb\nn M$ means $T_1\geq R_1+1=R+1$. By
hypothesis $a_2/a_1\in -\h\D 4\ooo^2$ and $\ord a_1=R$ so $\[
a_1,a_2\]\ap\h 12\pi^R\ab$ and $[a_1,a_2]\ap [\pi^R,-\D\pi^R]$. Also
since $R_2-R_1=-2e$ we have $d[-a_{1,2}]\geq R_1-R_2+\a_1=2e$. Note
also that if $s>2$ then $R_3=R+1$ and $R_s=R+1-2e$ so by Lemma 7.5 we
have $d[(-1)^{(s-2)/2}a_{3,s}]\geq 2e$ and
$[a_3,\ldots,a_s]\ap\hh\pp\ldots\pp\hh\pp [\pi^{R+1},-\xi\pi^{R+1}]$
with $\xi\in\{ 1,\D\}$. Since $d[-a_{1,2}]\geq 2e$ and
$d[(-1)^{(s-2)/2}a_{3,s}]\geq 2e$ we have $d[(-1)^{s/2}a_{1,s}]\geq
2e$. Also if $s>2$ then $[a_1,\ldots,a_s]\ap [a_1,a_2]\pp
[a_3,\ldots,a_n]\ap\hh\pp\cdots\pp\hh\pp
[\pi^R,-\D\pi^R,\pi^{R+1},-\xi\pi^{R+1}]$. 

We now prove the conditions 2.1(i)-(iv) for $N,K$.
\vskip 3mm

{\bf Proof of 2.1(i).} If $i\geq s+1$ then we have either
$T_i\geq R_i=S_i$ or $T_{i-1}+T_i\geq R_i+R_{i+1}=S_i+S_{i+1}$ so we
are done. If is odd and $1\leq i\leq s-3$ or $1\leq i\leq s-1$ when
the type is II then $T_i\geq T_1\geq R+1=S_i$ so we are done. If
$1<i\leq s-2$ is even then $T_i\geq T_2\geq T_1-2e\geq R+1-2e=S_i$ so
we are done. Thus we are left with the cases $i=s-1$ or $s$, if the
type is I, and $i=s$, if the type is II. 

Suppose first that $M$ is of type I. We claim that if (i) fails at
$i=s-1$ or $s$ then $T_1=T_{s-1}=R+1$ and, if $s>2$,
$T_{s-2}=R+1-2e$. Take first $i=s-1$. Suppose that
$T_{s-1}<S_{s-1}=R+2$. But $T_{s-1}\geq T_1\geq R+1$ so
$T_1=T_{s-1}=R+1$. If $s=2$ then we have our claim. If $s>2$ then also
$T_{s-3}=R+1$. Suppose that
$T_{s-2}+T_{s-1}<S_{s-1}+S_s=2R-2e+4$. Since $T_{s-1}=R+1$ we have
$T_{s-2}<R-2e+3$. Since also $T_{s-3}=R+1$ we get
$T_{s-2}-T_{s-3}<2-2e$ so $T_{s-2}-T_{s-3}=-2e$,
i.e. $T_{s-2}=R+1-2e$. If $i=s$ and $T_s<S_s=R-2e+2$ then $R+1\leq
T_1\leq T_{s-1}\leq T_s+2e<R+2$, which implies $R+1=T_{s-1}=T_s+2e$
and we are done. 

Since $T_1=R+1$ and $T_{s-2}=R+1-2e$ we have by Lemma
7.5 that $[c_1,\ldots,c_{s-2}]\ap\hh\pp\cdots\pp\hh\pp
[\pi^{R+1},-\xi'\pi^{R+1}]$ where $\xi'\in\{ 1,\D\}$ and
$d[(-1)^{(s-2)/2}c_{1,s-2}]\geq 2e$. (This also holds if $s=2$, when
$d[(-1)^{(s-2)/2}c_{1,s-2}]=\j$.) We want to prove that
$[c_1,\ldots,c_{s-1}]\rep [a_1,\ldots,a_s]$. We have
$R_{s+1}>R+1=T_{s-1}$. Also $d[(-1)^{(s-2)/2}c_{1,s-2}]\geq 2e$ and
$d[(-1)^{s/2}a_{1,s}]\geq 2e$ so $d[-a_{1,s}c_{1,s-2}]\geq
2e>2e+T_{s-1}-R_{s+1}$, which, by Corollary 2.17, implies
$C_{s-1}+C_s>2e+R_s-T_s$. So $[c_1,\ldots,c_{s-1}]\rep
[a_1,\ldots,a_s]$.

If $s=2$ then $c_1\rep [a_1,a_2]\ap [\pi^R,-\D\pi^R]$, which is
impossible since $\ord c_1=T_1=R+1$. If $s>2$ then
$[a_1,\ldots,a_s]\ap\hh\pp\cdots\pp\hh\pp
[\pi^R,-\D\pi^R,\pi^{R+1},-\xi\pi^{R+1}]$. Since $\ord
c_{s-1}=T_{s-1}=R+1$ we have $[-\xi'\pi^{R+1},c_{s-1}]\ap
[-\D\xi'\pi^{R+1},\D c_{s-1}]$. But $\xi,\xi'\in\{ 1,\D\}$ so
$\xi=\xi'$ or $\D\xi'$ in $\ff/\fs$ so
$[-\xi'\pi^{R+1},c_{s-1}]\ap [-\xi\pi^{R+1},c'_{s-1}]$, where
$c'_{s-1}=c_{s-1}$ or $\D c_{s-1}$. Thus
$[c_1,\ldots,c_{s-1}]\ap\hh\pp\cdots\pp\hh\pp
[\pi^{R+1},-\xi'\pi^{R+1},c_{s-1}]\ap\hh\pp\cdots\pp\hh\pp
[\pi^{R+1},-\xi\pi^{R+1},c'_{s-1}]$. Since $[c_1,\ldots,c_{s-1}]\rep
[a_1,\ldots,a_s]$ we get $c'_{s-1}\rep [\pi^R,-\D\pi^R]$. But this
is impossible since $\ord c'_{s-1}=\ord c_{s-1}=R+1$.

Suppose now that $M$ is of type II and (i) fails at $i=s$. Hence
$T_s<S_s=R-2e+3$ and $T_{s-1}+T_s<S_s+S_{s+1}=2R-2e+4$. We have
$2T_{s-1}-2e\leq T_{s-1}+T_s<2R-2e+4$ so $T_{s-1}<R+2$. But
$T_{s-1}\geq T_1\geq R+1$ so $T_1=T_{s-1}=R+1$. Since also
$T_s<R-2e+3$ we get $T_s-T_{s-1}<(R-2e+3)-(R+1)=2-2e$ so
$T_s-T_{s-1}=-2e$ so $T_s=R-2e+1$. Since $T_1=R+1$ and $T_s=R+1-2e$ we
have by Lemma 7.5 that $d[(-1)^{s/2}c_{1,s}]\geq 2e$ and
$[c_1,\ldots,c_s]\ap\hh\pp\cdots\pp\hh\pp
[\pi^{R+1},-\xi'\pi^{R+1}]$, where $\xi'\in\{ 1,\D\}$. We want to
prove that $[c_1,\ldots,c_s]\rep [a_1,\ldots,a_{s+1}]$. This statement
is trivial if $n=s+1$, when $[a_1,\ldots,a_{s+1}]\ap FM=FN$. If
$s+2\leq n$ then $R_{s+2}>R-2e+1=T_s$. In order to prove that
$C_s+C_{s+1}>2e+R_{s+1}-T_{s+1}$ it is enough, by Corollary 2.17, to
show that $d[-a_{1,s+2}c_{1,s}]>2e+T_s-R_{s+2}$. We have
$d[(-1)^{s/2}a_{1,s}],d[(-1)^{s/2}c_{1,s}]\geq 2e$ so
$d[a_{1,s}c_{1,s}]\geq 2e>2e+T_s-R_{s+2}$. Hence it is enough to prove
that $d[-a_{s+1,s+2}]>2e+T_s-R_{s+2}$. But $d[-a_{s+1,s+2}]\geq
R_{s+1}-R_{s+2}+\a_{s+1}=2e+T_s-R_{s+2}+\a_{s+1}$. (We have
$R_{s+1}=R+1=T_s+2e$.) So we have to prove that $\a_{s+1}>0$. But this
follows from $R_{s+2}-R_{s+1}>(R-2e+1)-(R+1)=-2e$. 

We have $[a_1,\ldots,a_{s+1}]\ap\hh\pp\cdots\pp\hh\pp
[\pi^R,-\D\pi^R,\pi^{R+1},-\xi\pi^{R+1},a_{s+1}]$. Since $\ord
a_{s+1}=R_{s+1}=R+1$ we have $[-\xi\pi^{R+1},a_{s+1}]\ap
[-\D\xi\pi^{R+1},\D a_{s+1}]$. But $\xi,\xi'\in\{ 1,\D\}$ so
$\xi'=\xi$ or $\D\xi$ in $\ff/\fs$. Hence
$[-\xi\pi^{R+1},a_{s+1}]\ap [-\xi'\pi^{R+1},a'_{s+1}]$,
where $a'_{s+1}=a_{s+1}$ or $\D a_{s+1}$. It follows that
$[a_1,\ldots,a_{s+1}]\ap\hh\pp\cdots\pp\hh\pp
[\pi^R,-\D\pi^R,\pi^{R+1},-\xi'\pi^{R+1},a'_{s+1}]$. We also have
$[c_1,\ldots,c_s]\ap\hh\pp\cdots\pp\hh\pp [\pi^{R+1},-\xi'\pi^{R+1}]$
so $[c_1,\ldots,c_s]\rep [a_1,\ldots,a_{s+1}]$ implies $\hh\rep
[\pi^R,-\D\pi^R,a'_{s+1}]$. But this is impossible since $\ord
a'_{s+1}=\ord a_{s+1}=R+1$ so $[\pi^R,-\D\pi^R,a'_{s+1}]$ is
anisotropic. 
\vskip 3mm

Before proving 2.1(ii) and 2.1(iii) we note that if $s>2$ then for
$1<i\leq s-3$ odd we have $T_{i-1}\geq T_2\geq T_1-2e\geq
R-2e+1=S_{i+1}$ so $i$ is not essential for $N,K$. Also for $2\leq
i\leq s-4$ even resp. for $2\leq i\leq s$ even when $M$ is of type II
we have $T_{i-1}\geq T_1\geq R+1=S_{i+1}$ so again $i$ is not
essential. Hence $i$ is not essential for $N,K$ for $1<i\leq s-3$, if
$M$ is of type I, and for $1<i\leq s-2$ or $i=s$, if $M$ is of type
II.  By Lemma 2.13, for such indices 2.1(iii) needs not be
verified. Also, by Lemma 2.12, 2.1(ii) needs not be verified for
$1<i\leq s-4$, if $M$ is of type I, and for $1<i\leq s-3$, if $M$ is
of type II. 

Also note that $[a_1,\ldots,a_i]=[b_1,\ldots,b_i]$ so
$a_{1,i}=b_{1,i}$ in $\ff/\fs$ for $i\geq s+1$. For such indices
$[c_1,\ldots,c_{i-1}]\rep [b_1,\ldots,b_i]$ is equivalent to
$[c_1,\ldots,c_{i-1}]\rep [a_1,\ldots,a_i]$. Also if $i\geq s+1$, $j$
is arbitrary and $\xi\in\ff$ we have $d(\xi a_{1,i}c_{1,j})=d(\xi
b_{1,i}c_{1,j})$ and also $\a_i=\b_j$. It follows that $d(\xi
a_{1,i}c_{1,j})=\min\{ d(\xi a_{1,i}c_{1,j}),\a_i,\c_j\} =\min\{ d(\xi
b_{1,i}c_{1,j}),\b_i,\c_j\} =d[\xi b_{1,i}c_{1,j}]$. 

In particular, for $i\geq s$ we have $C_i=\min\{ (R_{i+1}-T_i)/2+e,
R_{i+1}-T_i+d[-a_{1,i+1}c_{1,i-1}],
R_{i+1}+R_{i+2}-T_{i-1}-T_i+d[a_{1,i+2}c_{1,i-2}]\} =\min\{
(S_{i+1}-T_i)/2+e, S_{i+1}-T_i+d[-b_{1,i+1}c_{1,i-1}],
S_{i+1}+S_{i+2}-T_{i-1}-T_i+d[b_{1,i+2}c_{1,i-2}]\} =B_i$.
\vskip 3mm

{\bf Proof of 2.1(ii).} Take first $i=1$. If $s>2$ then $B_1\leq
(S_2-T_1)/2+e\leq ((R-2e+1)-(R+1))/2+e=0$ so 2.1(ii) holds
trivially. 


Suppose now that $i\geq s+1$. We have
$d[b_{1,i}c_{1,i}]=d[a_{1,i}c_{1,i}]$ and $B_i=C_i$ so
$d[b_{1,i}c_{1,i}]\geq B_i$ follows from $d[a_{1,i}c_{1,i}]\geq
C_i$. So we are left with the cases $s-3\leq i\leq s$, if $M$ is of
type I, and $s-2\leq i\leq s$, if $M$ is of type II.

Take first $i=s$. Suppose that 2.1(ii) does not hold so we have
$d[b_{1,s}c_{1,s}]<B_{s}=C_{s}\leq d[a_{1,s}c_{1,s}]$. Thus
$d[a_{1,s}b_{1,s}]=d[b_{1,s}c_{1,s}]<B_{s}\leq d[a_{1,s}c_{1,s}]$. We
cannot have $d[a_{1,s}b_{1,s}]=\a_s$ since $\a_s\geq
d[a_{1,s}c_{1,s}]$. If $d[a_{1,s}b_{1,s}]=\b_s$ then note that if $M$
is of type II then $S_{s+1}-S_s=(R+1)-(R-2e+3)=2e-2$ and if it is of
type I then $S_{s+1}-S_s\geq (R+2)-(R-2e+2)=2e$. In both cases
$\b_s=(S_{s+1}-S_s)/2+e\geq B_s$ so we get a contradiction. So we are
left with the case $d[a_{1,s}b_{1,s}]=d(a_{1,s}b_{1,s})$. If we have
type I then $a_{1,s}=b_{1,s}$ in $\ff/\fs$ so
$d(a_{1,s}b_{1,s})=\j>B_s$. If we have type II then
$a_{1,s+1}=b_{1,s+1}$ and $b_{s+1}=\eta a_{s+1}$. Thus
$a_{1,s}b_{1,s}\in\eta\fs$ and we have $d(a_{1,s}b_{1,s})=d(\eta
)=2e-1=(S_{s+1}-S_s)/2+e\geq B_s$.

Suppose now that $i=s-1$. We have $T_{s-1}\geq T_1\geq R+1$. Suppose
first that $T_{s-1}\geq R+2$. If $M$ is of type I then
$B_{s-1}\leq (S_s-T_{s-1})/2+e\leq ((R-2e+2)-(R+2))/2+e=0$ and we are
done. If $M$ is of type II then $B_{s-1}\leq
(S_s+S_{s+1})/2-T_{s-1}+e\leq ((R-2e+3)+(R+1))/2-(R+2)+e=0$ so we are
done. So we can assume that $T_1=T_3=\ldots=T_{s-1}=R+1$. If $M$ is of
type I then we cannot have $s=2$ because this would imply $S_1=R+2\geq
T_1$, which contradicts 2.1(i). So $s>2$ and we have
$T_{s-2}+R+1=T_{s-2}+T_{s-1}\geq S_{s-2}+S_{s-1}=(R-2e+1)+(R+2)$ so
$T_{s-2}\geq R-2e+2$. But we cannot have $T_{s-2}=R-2e+2$ since this
would imply $T_{s-2}-T_{s-3}=1-2e$. So $T_{s-2}\geq R-2e+3$. It
follows that $B_{s-1}\leq S_s-(T_{s-3}+T_{s-2})/2+e\leq
R-2e+2-((R+1)+(R-2e+3))/2+e=0$ and we are done. Suppose now $M$ is of
type II. We have $B_{s-1}\leq
(S_s+S_{s+1})/2-T_{s-1}+e=((R-2e+3)+(R+1))/2-(R+1)+e=1$. Now
$S_1=S_{s-1}=R+1$ and $T_1=T_{s-1}=R+1$. By Lemma 6.6(i) we have
$S_1+\cdots +S_{s-1}\ev T_1+\cdots +T_{s-1}\ev R+1\m2$ so $\ord
b_{1,s-1}c_{1,s-1}$ is even. It follows that
$d(b_{1,s-1}c_{1,s-1})\geq 1$. Now $T_{s-1}=S_{s-1}=R+1$ and by 2.1(i)
$T_{s-1}+T_s\geq S_{s-1}+S_s$ so $T_s-T_{s-1}\geq
S_s-S_{s-1}=(R-2e+3)-(R+1)=2-2e$. It follows that $\b_{s-1}=1$ and
$\c_{s-1}\geq 1$. Together with $d(b_{1,s-1}c_{1,s-1})\geq 1$, these
imply $d[b_{1,s-1}c_{1,s-1}]=1\geq B_{s-1}$. 

Let now $i=s-2$. Suppose first that $M$ is of type I. We have
$T_{s-3}\geq T_1\geq R+1$. If $T_{s-3}\geq R+2=S_{s-1}$ then $s-2$ is
not an essential index for $N,K$. We also have $T_{s-2}\geq
T_{s-3}-2e\geq R-2e+2=S_s$ so $s-1$ is not essential either and thus
the statement 2.1(ii) becomes vacuous at $i=s-2$. Thus we may assume
that $T_1=T_{s-3}=R+1$. If $T_{s-2}-T_{s-3}\geq 2-2e$ then
$\c_{s-3}\geq 1$ and $T_{t-2}\geq R-2e+3$. We have $B_{s-2}\leq
(S_{s-1}+S_s)/2-T_{s-2}+e=((R+2)+(R-2e+2))/2-T_{s-2}+e=R+2-T_{s-2}$. Since
$S_1=R+1$ and $S_{s-2}=R+1-2e$ we have by Lemma 7.5
$d[(-1)^{(s-2)/2}b_{1,s-2}]\geq 2e>R+2-T_{s-2}$. Also
$T_1=T_{s-3}=R+1$ which, by Lemma 7.4(i), implies
$d[(-1)^{(s-2)/2}c_{1,s-2}]\geq
T_{s-3}-T_{s-2}+\c_{s-3}=R+1-T_{s-2}+\c_{s-3}\geq R+2-T_{s-2}$. Hence
$d[b_{1,s-2}c_{1,s-2}]\geq R+2-T_{s-2}\geq B_{s-2}$. So we are left
with the case when $T_{s-2}-T_{s-3}=-2e$, i.e. $T_{s-2}=R-2e+1$. Since
$S_{s-1}-T_{s-2}=2e+1>2e$, by Corollary 2.10, we have to prove that
$b_{1,s-2}c_{1,s-2}\in\fs$. Since $T_1=S_1=R+1$ and
$T_{s-2}=S_{s-2}=R-2e+1$ by Lemma 7.5 both $\[ c_1,\ldots,c_{s-2}\]$
and $\[ b_1,\ldots,b_{s-2}\] =\[ a_3,\ldots,a_s\]$ are orthogonal sums
of copies of $\h 12\pi^{R+1}\aa$ and $\h 12\pi^{R+1}\ab$. If $\[
b_1,\ldots,b_{s-2}\]\ap\[ c_1,\ldots,c_{s-2}\]$ then
$b_{1,s-2}c_{1,s-2}\in\fs$ and we are done. Otherwise, since
$\aa\pp\aa\ap\ab\pp\ab$, we have $\[ c_1,\ldots,c_{s-2}\]\pp\h
12\pi^{R+1}\aa\ap\[ b_1,\ldots,b_{s-2}\]\pp\h 12\pi^{R+1}\ab\ap\[
a_3,\ldots,a_s\]\pp\h 12\pi^{R+1}\ab$. It follows that
$[c_1,\ldots,c_{s-2}]\pp\hh\ap [a_3,\ldots,a_s]\pp
[\pi^{R+1},-\D\pi^{R+1}]$ and $a_{3,s}c_{1,s-2}\in\D\fs$. Since also
$-a_{1,2}\in\D\fs$ we have $-a_{1,s}c_{1,s-2}\in\fs$. If $s<n$ then
$R_{s+1}\geq R+2>R+1=T_{s-2}+2e$ so by Lemma 2.19 we have
$[c_1,\ldots,c_{s-2}]\rep [a_1,\ldots,a_s]$. Same happens obviously if
$s=n$. Since also $-a_{1,s}c_{1,s-2}\in\fs$ we have
$[a_1,\ldots,a_s]\ap [c_1,\ldots,c_{s-2}]\pp\hh\ap [a_3,\ldots,a_s]\pp
[\pi^{R+1},-\D\pi^{R+1}]$ so $[\pi^R,-\D\pi^R]\ap
[a_1,a_2]\ap[\pi^{R+1},-\D\pi^{R+1}]$, which is false.

Suppose now that $M$ is of type II. In this case $s-2$ is not
essential for $N,K$. So we may assume that $s-1$ is essential since
otherwise 2.1(ii) becomes vacuous at $i=s-2$ by Lemma 2.12. It follows
that $2T_{s-3}-2e\leq T_{s-3}+T_{s-2}<S_s+S_{s+1}=R-2e+3+R+1=2R-2e+4$
so $R+1\leq T_1\leq T_{s-3}<R+2$ and so $T_1=T_{s-3}=R+1$. Thus
$R+1+T_{s-2}=T_{s-3}+T_{s-2}<2R-2e+4$, i.e. $T_{s-2}<R-2e+3$ so
$T_{s-2}-T_{s-3}<2-2e$. Hence $T_{s-2}-T_{s-3}=-2e$ and so
$T_{s-2}=R+1-2e$. Since $S_1=T_1=R+1$, $S_{s-2}=T_{s-2}=R+1-2e$ we
have $d[(-1)^{(s-2)/2}b_{1,s-2}]\geq 2e$ and
$d[(-1)^{(s-2)/2}c_{1,s-2}]\geq 2e$ by Lemma 7.5. It follows that
$d[b_{1,s-2}c_{1,s-2}]\geq 2e=(S_{s-1}-T_{s-2})/2+e\geq B_{s-2}$. (We
have $S_{s-1}=R+1$ and $T_{s-2}=R-2e+1$.)

We are left with $i=s-3$ and only if $M$ is of type I. We have
$T_{s-3}\geq T_1\geq R+1$ and $S_{s-2}=R-2e+1$ so $B_{s-3}\leq
(S_{s-2}-T_{s-3})/2+e\leq 0$ and 2.1(ii) becomes trivial.
\vskip 3mm

{\bf Proof of 2.1(iii).} Let $1<i<n$ s.t. $T_{i-1}<S_{i+1}$ and
$B_{i-1}+B_i>2e+S_i-T_i$. We want to prove that
$[c_1,\ldots,c_{i-1}]\rep [b_1,\ldots,b_i]$. 

If $i\geq s+1$ then $B_{i-1}=C_{i-1}$ and $B_i=C_i$ so
$C_{i-1}+C_i=B_{i-1}+B_i>2e+S_i-T_i=2e+R_i-T_i$. We also
have $T_{i-1}<S_{i+1}=R_{i+1}$ so $[c_1,\ldots,c_{i-1}]\rep
[a_1,\ldots,a_i]\ap [b_1,\ldots.b_i]$. So we are left with the cases
$s-2\leq i\leq s$ if $M$ is of type I and $i=s-1$ if it is
of type II. (All the other indices are not essential.)

Take first $i=s$, when $M$ is of type I. We have $[b_1,\ldots,b_s]\ap
[a_1,\ldots,a_s]$ so it is enough to prove that
$[c_1,\ldots,c_{s-1}]\rep [a_1,\ldots,a_s]$, i.e. that
$T_{s-1}<R_{s+1}$ and $C_{s-1}+C_s>2e+R_s-T_s$. We have
$T_{s-1}<S_{s+1}=R_{s+1}$ and $C_{s-1}+C_s>2e+R_s-T_s$ is equivalent
by Lemma 2.12 to
$d[-a_{1,s+1}c_{1,s-1}]+d[-a_{1,s}c_{1,s-2}]>2e+T_{s-1}-R_{s+1}$,
$d[-a_{1,s+1}c_{1,s-1}]>e+(T_{s-2}+T_{s-1})/2-R_{s+1}$ and
$d[-a_{1,s}c_{1,s-2}]>e+T_{s-1}-(R_{s+1}+R_{s+2})/2$. But we have
$R_{s+1}=S_{s+1}$, $R_{s+2}=S_{s+2}$ and
$d[-a_{1,s+1}c_{1,s-1}]=d[-b_{1,s+1}c_{1,s-1}]$ so the conditions
above will follow from the corresponding conditions for $N,K$
provided that $d[-a_{1,s}c_{1,s-2}]\geq d[-b_{1,s}c_{1,s-2}]$.
Suppose that $d[-a_{1,s}c_{1,s-2}]<d[-b_{1,s}c_{1,s-2}]$. We have
$d[-a_{1,s}c_{1,s-2}]=\min\{
d(-a_{1,s}c_{1,s-2}),\a_s,\c_{s-2}\}$. But $\c_{s-2}\geq
d[-b_{1,s}c_{1,s-2}]$ and
$d(-a_{1,s}c_{1,s-2})=d(-b_{1,s}c_{1,s-2})\geq
d[-b_{1,s}c_{1,s-2}]$ so we most have $d[-a_{1,s}c_{1,s-2}]=\a_s$.
But $R_{s+1}=S_{s+1}\geq R+2$ and $R_s<R-2e+2=S_s$. ($R_s=R-2e$ or
$R-2e+1$ if $s=2$ resp. $s>2$.) Thus $R_{s+1}-R_s>S_{s+1}-S_s\geq 2e$
and so $\a_s=(R_{s+1}-R_s)/2+e>(S_{s+1}-S_s)/2+e=\b_s\geq
d[-b_{1,s}c_{1,s-2}]$. Contradiction.

Take now $i=s-1$. Suppose first that $M$ is of type I. We have
$T_{s-2}<S_s=R-2e+2$. But $T_{s-3}\geq T_1\geq R+1$ so
$T_{s-2}\geq T_{s-3}-2e\geq R-2e+1$. It follows that
$T_{s-3}=R-2e+1$ and $T_1=T_{s-3}=R+1$. But we have
encountered this situation during the proof of 2.1(ii) at $i=s-2$,
when we proved that $\[ c_1,\ldots,c_{s-2}\]\ap\[
b_1,\ldots,b_{s-2}\]$ and so $[c_1,\ldots,c_{s-2}]\rep
[b_1,\ldots,b_{s-1}]$. Suppose now $M$ is of type II. Since $s-1$ is
essential we have $2T_{s-3}-2e\leq
T_{s-3}+T_{s-2}<S_s+S_{s+1}=R-2e+3+R+1=2R-2e+4$ so $T_{s-3}<R+2$.
But $T_{s-3}\geq T_1\geq R+1$ so $T_1=T_{s-3}=R+1$. We
also have $T_{s-2}<S_s=R-2e+3$ so $T_{s-2}-T_{s-3}\leq 2-2e$ and so
$T_{s-2}-T_{s-3}=-2e$, i.e. $T_{s-2}=R-2e+1$. Since $S_1=T_1=R+1$ and
$S_{s-2}=T_{s-2}=R+1-2e$ we have, by Lemma 7.5,
$[b_1,\ldots,b_{s-2}]\ap\hh\pp\cdots\pp\hh\pp
[\pi^{R+1},-\xi\pi^{R+1}]$ and
$[c_1,\ldots,c_{s-2}]\ap\hh\pp\cdots\pp\hh\pp
[\pi^{R+1},-\xi'\pi^{R+1}]$ for some $\xi,\xi'\in\{ 1,\D\}$. It
follows that $[b_1,\ldots,b_{s-1}]\ap\hh\pp\cdots\pp\hh\pp
[\pi^{R+1},-\xi\pi^{R+1},b_{s-1}]$. But $\ord b_{s-1}=S_{s-1}=R+1$
and $\xi\xi'$ is a power of $\D$ so $-\pi^{R+1}\xi'\rep
[-\pi^{R+1}\xi,b_{s-1}]$, which implies that
$[c_1,\ldots,c_{s-2}]\ap\hh\pp\cdots\pp\hh\pp
[\pi^{R+1},-\pi^{R+1}\xi']\rep\hh\pp\cdots\pp\hh\pp
[\pi^{R+1},-\pi^{R+1}\xi,b_{s-1}]\ap [b_1,\ldots,b_{s-1}]$.

Finally, we take $i=s-2$, when $M$ is of type I. We have $R+1=T_1\leq
T_{s-3}<S_{s-1}=R+2$ so $T_1=T_3=\ldots=T_{s-3}=R+1$. Now $i=s-2$ is
essential so if $s>4$ then
$T_{s-4}+R+1=T_{s-4}+T_{s-3}<S_{s-1}+S_s=R+2+R-2e+2=2R-2e+4$
so $T_{s-4}<R-2e+3$. It follows that $T_{s-4}-T_{s-5}<2-2e$, which
implies that $T_{s-4}-T_{s-5}=-2e$ so $T_{s-4}=R-2e+1$. Now
$S_1=T_1=R+1$ and $S_{s-2}=T_{s-4}=R+1-2e$ so, by Lemma 7.5, $\[
b_1,\ldots,b_{s-2}\]$ and $\[
c_1,\ldots,c_{s-4}\]$ are orthogonal sums of copies of $\h
12\pi^{R+1}\aa$ and $\h 12\pi^{R+1}\ab$. (Same happens if $s=4$,
when $\[ c_1,\ldots,c_{s-4}\]$ is just the zero lattice.) It follows
that $\[ b_1,\ldots,b_{s-2}\]\ap\[ c_1,\ldots,c_{s-4}\]\pp\h
12\pi^{R+1}\aa$ or $\[ c_1,\ldots,c_{s-4}\]\pp\h 12\pi^{R+1}\ab$,
which implies that $[b_1,\ldots,b_{s-2}]\ap
[c_1,\ldots,c_{s-4}]\pp\hh$ or $[c_1,\ldots,c_{s-4}]\pp
[\pi^{R+1},-\pi^{R+1}\D ]$. But $\ord c_{s-3}=R+1$ so $c_{s-3}$ is
represented by both $\hh$ and $[\pi^{R+1},-\pi^{R+1}\D ]$. It follows
that $[c_1,\ldots,c_{s-3}]$ is represented by both
$[c_1,\ldots,c_{s-4}]\pp\hh$ and $[c_1,\ldots,c_{s-4}]\pp
[\pi^{R+1},-\pi^{R+1}\D ]$ and thus by $[b_1,\ldots,b_{s-2}]$.
\vskip 3mm

{\bf Proof of 2.1(iv).} Let $1<i\leq n-2$ s.t. $S_{i+1}+2e\leq
T_{i-1}+2e<S_{i+2}\leq T_i$. We want to prove that
$[c_1,\ldots,c_{i-1}]\rep [b_1,\ldots,b_{i+1}]$. 

If $i\geq s-1$, when $M$ is of type I, or if $i\geq s$, when $M$ is of
type II, we have $R_{i+2}=S_{i+2}>T_{i-1}+2e$, which implies by Lemma
2.19 that $[c_1,\ldots,c_{i-1}]\rep [a_1,\ldots,a_{i+1}]\ap
[b_1,\ldots,b_{i+1}]$. 

If $M$ is of type I and $i\leq s-2$ the condition $S_{i+1}+2e<S_{i+2}$
implies that $i=s-3$. We have
$T_{s-4}+2e<S_{s-1}=R+2$ so $T_{s-4}\leq R+1-2e$. But $T_{s-5}\geq
T_1\geq R+1$ so $T_{s-4}\geq T_{s-5}-2e\geq R+1-2e$. Thus
$T_1=R+1$ and $T_{s-4}=R+1-2e$. Also $S_1=R+1$ and $S_{s-2}=R-2e+1$ so
by Lemma 7.5 $\[ b_1,\ldots,b_{s-2}\] $ and $\[ c_1,\ldots,c_{s-4}\]$
are orthogonal sums of copies of $\h 12\pi^{R+1}\aa$ and $\h
12\pi^{R+1}\ab$. It follows that $\[ c_1,\ldots,c_{s-4}\]\rep\[
b_1,\ldots,b_{s-2}\]$ so $[c_1,\ldots,c_{s-4}]\rep [
b_1,\ldots,b_{s-2}]$. 

Finally, if $M$ is of type II then the condition $S_{i+1}+2e<S_{i+2}$
is not satisfied for any value of $i\leq s-1$. \qed

We now consider the remaining case when $R_2-R_1=-2e$ and either
$a_2/a_1\in -\h 14\ooo^2$ or $R_1=R_3$.

\blm Suppose $M$ is as above. Let $R:=R_1$ and let $s\geq 2$ be even
and maximal with the property that $R_s=R-2e$. Let $k$ be maximal with
the property that $M$ splits $k$ copies of $\h 12\pi^R\aa$. Then
$k\geq 1$
and we have the following 3 cases:

1. If $R_{s+1}>R$ or $s=n$ and $(-1)^{s/2}a_{1,s}\in\fs$ then
$k=s/2$.

2. If $R_{s+1}>R$ or $s=n$ and $(-1)^{s/2}a_{1,s}\in\D\fs$ then
$k=s/2-1$.

3. If $R_{s+1}=R$ then $k=s/2$
\elm
\pf First we show that $k\leq s/2$. Suppose the contrary so $M$
represents $N$, an orthogonal sum of $s/2+1$ copies of $\h
12\pi^R\aa$. We have $N\ap\[ b_1,\ldots,b_{s+2}\]$ where
$b_1=b_3=\ldots=b_{s+1}=\pi^R$ and $b_2=b_4=\ldots=b_{s+2}=-\pi^{R-2e}$. We
have $N\rep M$ so $N\leq M$. It follows that $R_{s+1}+R_{s+2}\leq
S_{s+1}+S_{s+2}=R+R-2e=2R-2e$. But $R_{s+1}\geq R_1=R$ and
$R_{s+2}\geq R_2=R-2e$ so we must have $R_{s+1}=R$ and
$R_{s+2}=R-2e$, which contradicts the maximality of $s$.

We have $R=R_1\leq R_{s-1}\leq R_s+2e=R$ so $R_1=R_{s-1}=R$. Since also
$R_s=R-2e$ we have, by Lemma 7.5, that $\[ a_1,\ldots,a_s\]$ is an
orthogonal sum of copies of $\h 12\pi^R\aa$ and $\h 12\pi^R\ab$ and so
it is isomorphic to $\h 12\pi^R\aa\pp\cdots\pp\h 12\pi^R\aa$ or $\h
12\pi^R\aa\pp\cdots\pp\h 12\pi^R\aa\pp\h 12\pi^R\ab$. Considering
determinants we have $\[ a_1,\ldots,a_s\]\ap\h
12\pi^R\aa\pp\cdots\pp\h 12\pi^R\aa$ if $a_{1,s}\in (-1)^{s/2}\fs$ and
$\[ a_1,\ldots,a_s\]\ap\h 12\pi^R\aa\pp\cdots\pp\h 12\pi^R\aa\pp\h
12\pi^R\ab$ if $a_{1,s}\in (-1)^{s/2}\D\fs$. In both cases if $s+1\leq
n$ we have $R_{s+1}\geq R_1=R>R-2e=R_s$ so $M\ap\[
a_1,\ldots,a_s\]\pp\[ a_{s+1},\ldots,a_n\]$. Hence $M$ splits at least
$s/2-1$ copies of $\h 12\pi^R\aa$. So $s/2-1\leq k\leq s/2$.

If $a_{1,s}\in (-1)^{s/2}\fs$, in particular in the case 1., then $M$
splits $\[ a_1,\ldots,a_s\]\ap\h 12\pi^R\aa\pp\cdots\pp\h 12\pi^R\aa$ so
$k=s/2$. Same happens in the case 3. if $a_{1,s}\in
(-1)^{s/2}\fs$. If $a_{1,s}\in (-1)^{s/2}\D\fs$ then $\ord
a_{s+1}/a_s=R_{s+1}-R_s=R-(R-2e)=2e$ so $\D\in g(a_{s+1}/a_s)$. (See
[B1, Definition 6].) By [B1, 3.12] we have $\[ a_1,\ldots,a_n\]\ap\[
a_1,\ldots,\D a_s,\D a_{s+1},a_{s+2},\ldots,a_n\]$. By this change of
BONGs the product $a_{1,s}$ is replaced by $\D a_{1,s}$ which belongs
to $(-1)^{s/2}\fs$ so we are back to the previous case. So again
$k=s/2$. 

Finally, in the case 2. suppose that $k=s/2$ so $M$ represents $N$, an
orthogonal sum of $s/2$ copies of $\h 12\pi^R\aa$. We have $N\ap\[
b_1,\ldots,b_s\]$, where $b_1=b_3=\ldots=b_{s-1}=\pi^R$ and
$b_2=b_4=\ldots=b_s=-\pi^{R-2e}$. If $s>n$ then $R_{s+1}>R$ and
$S_s=R-2e$ so $R_{s+1}-S_s>2e$. By Corollary 2.10 we have
$a_{1,s}b_{1,s}\in\fs$. Same happens when $s=n$ since
$[a_1,\ldots,a_n]\ap [b_1,\ldots,b_n]$. But
$(-1)^{s/2}a_{1,s}\in\D\fs$ and $(-1)^{s/2}b_{1,s}\in\fs$ so
$a_{1,s}b_{1,s}\in\D\fs$. Contradiction. So $k=s/2-1$. \qed

In the context of Lemma 7.17 we say that $M$ is of type I, II or III
if it is in the case 1., 2., resp 3. of Lemma 7.17.

With the notations of Lemma 7.17 let $M=\hh_1\pp\cdots\pp\hh_k\pp P$
with $\hh_i\ap\h 12\pi^R\aa$ relative to some basis $t_i,u_i$.
Note that $Q(t_i)=Q(u_i)=0$ and $Q(t_i+u_i)=\pi^R$. So if
$M'=\{x\in M\mid x$ not a norm generator$\}$ then $t_i,u_i\in M'$
but $t_i+u_i\notin M'$. Therefore this time $M'$ is not a lattice
and we cannot take $N=M'$.

In this case we take $N=\hh'_1\pp\cdots\pp\hh'_k\pp P$, where
$\hh'_i\ap\h 12\pi^{R+1}\aa$ relative to the basis $t_i,\pi
u_i$. Since $\hh'_i\sb\hh_i$ we have $N\sb M$.

\blm With the notations above we have:

(i) If $M$ is of type I we can choose its BONG
s.t. $a_1=a_3=\ldots=a_{s-1}=\pi^R$,
$a_2=a_4=\ldots=a_s=-\pi^{R-2e}$. We have
$$N\ap\[\pi^{R+1},-\pi^{R-2e+1},\ldots,\pi^{R+1},-\pi^{R-2e+1},
a_{s+1},\ldots,a_n\]$$
relative to a good BONG $y_1,\ldots,y_s,x_{s+1},\ldots,x_n$.

(ii) If $M$ is of type II we can choose its BONG
s.t. $a_1=a_3=\ldots=a_{s-1}=\pi^R$, $a_2=-\D\pi^{R-2e}$,
$a_4=a_6=\ldots=a_s=-\pi^{R-2e}$. We have
$$N\ap\[\pi^R,-\D\pi^{R-2e},
\pi^{R+1},-\pi^{R-2e+1},\ldots,\pi^{R+1},-\pi^{R-2e+1},
a_{s+1},\ldots,a_n\]$$
relative to a good BONG $x_1,x_2,y_3,\ldots,y_s,x_{s+1},\ldots,x_n$.

(iii) If $M$ is of type III we can choose its BONG
s.t. $a_1=a_3=\ldots=a_{s-1}=\pi^R$,
$a_2=a_4=\ldots=a_s=-\pi^{R-2e}$. We have
$$N\ap\[
\pi^R,-\pi^{R-2e+2},\ldots,\pi^R,-\pi^{R-2e+2},a_{s+1},\ldots,a_n\]$$
relative to some good BONG $y_1,\ldots,y_{s+1},x_{s+2},\ldots,x_n$.
\elm
\pf As seen from the proof of Lemma 7.17, if $M$ is of type I then $\[
a_1,\ldots,a_s\]\ap\h 12\pi^R\aa\pp\cdots\pp\h 12\pi^R\aa$. Same can
be assumed in the case III by making, if necessary, the change of BONG
$\[ a_1,\ldots,a_n\]\ap\[ a_1,\ldots,a_{s-1},\D a_s,\D
a_{s+1},a_{s+2}\ldots,a_n\]$. Since $\h 12\pi^R\aa\pp\cdots\pp\h
12\pi^R\aa\ap\[\pi^R,-\pi^{R-2e},\ldots,\pi^R,-\pi^{R-2e}\]$ we can
change the BONG of $M$ s.t. the sequence $a_1,\ldots,a_s$ is
$\pi^R,-\pi^{R-2e},\ldots,\pi^R,-\pi^{R-2e}$, as claimed. If $M$ is of
type II then $\[ a_1,\ldots,a_s\]\ap\h 12\pi^R\ab\pp\h
12\pi^R\aa\pp\cdots\pp\h 12\pi^R\aa\ap \\ \[\pi^R,-\D\pi^{R-2e}\]\pp
\[\pi^R,-\pi^{R-2e}\]\pp\cdots\pp\[\pi^R,-\pi^{R-2e}\]\ap \\
\[\pi^R,-\D\pi^{R-2e},\pi^R,-\pi^{R-2e},\ldots,\pi^R,-\pi^{R-2e}\]$.
Therefore we can change the BONG of $M$ s.t. the sequence $a_1,..,a_s$
is $\pi^R,-\D\pi^{R-2e},\pi^R,-\pi^{R-2e},\ldots,\pi^R,-\pi^{R-2e}$,
as claimed. So we have proved the first part of (i)-(iii).

Before proving the second part we note that in the cases I and III we
have \\ $\[ x_1,\ldots,x_s\]=H_1\pp\cdots\pp H_k$, with $H_i\ap\h
12\pi^R\aa$. Since both $\[ x_1,\ldots,x_s\]$ and
$\hh_1\pp\cdots\pp\hh_k$ split $M$ there is by [OM, Theorem 93:14]
some $\s\in O(M)$ s.t. $\hh_1\pp\cdots\pp\hh_k=\s (\[ x_1,\ldots,x_s\]
)=\[\s x_1,\ldots,\s x_s\]$. We have $M=\s M=\[\s x_1,\ldots,\s x_n\]$
so if we replace the BONG of $M$ by $\s x_1,\ldots,\s x_n$ we may
assume that $\hh_1\pp\cdots\pp\hh_k=\[ x_1,\ldots,x_s\]$. We have
$M=\[ x_1,\ldots,x_s\]\pp\[ x_{s+1},\ldots,x_n\]$ so in these cases
$P=\[ x_{s+1},\ldots,x_n\]$. If $M$ is of type II then $R_2<R_3$ and
$R_s<R_{s+1}$ so by [B1, Corollary 4.4 (i)] we have $M=\[
x_1,x_2\]\pp\[ x_3,\ldots,x_s\]\pp\[ x_{s+1},\ldots,x_n\]$. Now both
$\[ x_3,\ldots,x_s\]\ap\[\pi^R,-\pi^{R-2e},\ldots,\pi^R,-\pi^{R-2e}\]$
and $\hh_1\pp\cdots\pp\hh_k$ split $M$ and are $\ap\h
12\pi^R\aa\pp\cdots\pp\h 12\pi^R\aa$ so there is $\s\in O(M)$
s.t. $\hh_1\pp\cdots\pp\hh_k=\s (\[ x_3,\ldots,x_s\] )=\[\s
x_3,\ldots,\s x_s\]$. When we replace the BONG of $M$ by $\s
x_1,\ldots,\s x_n$ we may assume that $\hh_1\pp\cdots\pp\hh_k=\[
x_3,\ldots,x_s\]$. So in this case $P=\[ x_1,x_2\]\pp\[
x_{s+1},\ldots,x_n\]$. We now prove the second part of (i)-(iii).

(i) We have $\hh'_1\pp\cdots\pp\hh'_k\ap\h
12\pi^{R+1}\aa\pp\cdots\pp\h 12\pi^{R+1}\aa\ap \\
\[\pi^{R+1},-\pi^{R-2e+1},\ldots,\pi^{R+1},-\pi^{R-2e+1}\]$ relative
to a good BONG $y_1,\ldots,y_s$. If we denote $S_i:=\ord Q(y_i)$ then
$S_1=S_3=\ldots=S_{s-1}=R+1$ and $S_2=S_4=\ldots=S_s=R-2e+1$. We have
$N=\hh'_1\pp\cdots\pp\hh'_k\pp P=\[ y_1,\ldots,y_s\]\pp\[
x_{s+1},\ldots,x_n\]$. We want to prove that
$y_1,\ldots,y_s,x_{s+1},\ldots,x_n$ is a good BONG for $N$. If $s=n$
this statement is trivial. Otherwise we use [B1, Corollary 4.4(v)]. We
have $R_{s+1}\geq R+1=S_{s-1}$ and $R_{s+1}>R-2e+1=S_s$. Also if
$s\leq n-2$ then $R_{s+2}\geq R_{s+1}-2e\geq R-2e+1=S_s$ so we are
done.

(ii) We have $\hh'_1\pp\cdots\pp\hh'_k\ap
\[\pi^{R+1},-\pi^{R-2e+1},\ldots,\pi^R,-\pi^{R-2e+1}\]$ relative to a
good BONG $y_3,\ldots,y_s$. Denote $S_i=\ord Q(y_i)$. We have
$M=\hh_1\pp\cdots\pp\hh_k\pp P=\[ y_1,\ldots,y_s\]\pp\[ x_1,x_2\]\pp\[
x_{s+1},\ldots,x_n\]$. By the same reasoning above
$y_3,\ldots,y_s,x_{s+1},\ldots,x_n$ is a good BONG for $\[
y_3,\ldots,y_s\]\pp\[ x_{s+1},\ldots,x_n\]$. In order to prove that
$x_1,x_2,y_3,\ldots,y_s,x_{s+1},\ldots,x_n$ is a good BONG for $M=\[
x_1,x_2\]\pp\[ y_3,\ldots,y_s,x_{s+1},\ldots,x_n\]$ we use [B1, Corollary
4.(v)]. We have $R_1=R<R+1=S_3$, $R_2=R-2e+1<R+1=S_3$ and
$R_2=R-2e<R-2e+1=S_4$ so we are done.

(iii) Note that, by the maximality of $s$, if $s\leq n-2$ then
$R_{s+2}\neq R-2e$ so $R_{s+2}-R_{s+1}=R_{s+2}-R\neq -2e$. It follows
that $R_{s+2}-R_{s+1}\geq 2-2e$ so $R_{s+2}\geq R-2e+2$. We have
$\hh'_1\pp\cdots\pp\hh'_k\ap
\[\pi^{R+1},-\pi^{R-2e+1},\ldots,\pi^{R+1},-\pi^{R-2e+1}\]$ relative
to some BONG $x'_1,\ldots,x'_s$ and $P=\[ x_{s+1},\ldots,x_n\]$. If
we denote $R'_i:=\ord Q(x'_i)$ then $R'_1=R'_3=\ldots=R'_{s-1}=R+1$
and $R'_2=R'_4=\ldots=R'_s=R-2e+1$. By Corollary 7.13 we have $\[
x'_1,\ldots,x'_s\]\pp\[ x_{s+1}\]\ap\h
12\pi^{R+1}\aa\pp\cdots\pp\h 12\pi^{R+1}\aa\pp\la
a_{s+1}\ra\ap\[\pi^R,-\pi^{R-2e+2},\ldots,
\pi^R,-\pi^{R-2e+2},a_{s+1}\]$
relative to some good BONG $y_1,\ldots,y_{s+1}$. (We have $\ord
a_{s+1}=R_{s+1}=R+1$.) If $S_i:=\ord Q(y_i)$ then
$S_1=S_3=\ldots=S_{s+1}=R$ and $S_2=S_4=\ldots=S_s=R-2e+2$. Now
$y_1,\ldots,y_{s+1},x_{x+2},\ldots,x_n$ is a good
BONG. ($y_1,\ldots,y_{s+1}$ is a good BONG and so is
$y_{s+1},x_{s+2},\ldots,x_n$ since
$Q(y_{s+1})=Q(x_{s+1})=a_{s+1}$. Also if $s+2\leq n$ then
$S_s=R-2e+2\leq R_{s+2}$.) In order to prove that
$y_1,\ldots,y_{s+1},x_{x+2},\ldots,x_n$ is a BONG for $N=\[
x'_1,\ldots,x'_s\]\pp\[ x_{s+1},\ldots,x_n\]$ we use Lemma 7.10. If
$n\geq s+2$ we have $R'_s=R-2e+2\leq R_{s+2}$ so we are done. \qed

\blm We have $[a_1,\ldots,a_i]\ap [b_1,\ldots,b_i]$ if $M$ is of type
III or if $i\geq s$ or if $i$ is even. We also have $\a_i=\b_i$ for
$i\geq s+1$.
\elm
\pf If $M$ is of type III then $a_j=b_j$ in $\ff/\fs$ at any
$j$ so $[a_1,\ldots,a_i]\ap [b_1,\ldots,b_i]$ for any $i$. If the type
is I or II then $[a_j,a_{j+1}]\ap [b_j,b_{j+1}]$ for any $i\leq j<s$
odd and $a_j=b_j$ for any $j\geq s+1$. It follows that
$[a_1,\ldots,a_i]\ap [b_1,\ldots,b_i]$ for $1<i\leq s$ even and for
$i\geq s+1$.

Let $K\ap\[ a_{s+1},\ldots,a_n\]\ap\[ b_{s+1},\ldots,b_n\]$. We prove
now that $\a_i=\b_i$ for $i\geq s+1$. By [B3, Lemmas 2.1 and 2.4(i)]
we have
$\a_i=\min\{\a_{i-s}(K),R_{i+1}-R_{s+1}+\a_s\}$. ($\a_{i-s}(K)$
replaces $(R_{i+1}-R_i)/2+e$ and the terms with $j\geq s+1$ in the
definition of $\a_i$, while $R_{i+1}-R_{s+1}+\a_s$ replaces the terms
with $j\leq s$.) Similarly
$\b_i=\min\{\a_{i-s}(K),S_{i+1}-S_{s+2}+\b_{s+1}\}$. If $M$ is of type
I or II then $R_{s+1}-R_s\geq (R+1)-(R-2e)=2e+1$ and $S_{s+1}-S_s\geq
(R+1)-(R-2e+1)=2e$. If $M$ is of type III then
$R_{s+1}-R_s=R-(R-2e)=2e$ and $S_{s+1}-S_s=R-(R-2e+2)=2e-2$. In all
cases we have $\a_s=(R_{s+1}-R_s)/2+e$ and
$\b_s=(S_{s+1}-S_s)/2+e$. Thus
$R_{i+1}-R_{s+1}+\a_s=R_{i+1}-(R_s+R_{s+1})/2+e\geq
R_{i+1}-(R_i+R_{i+1})/2+e=(R_{i+1}-R_i)/2+e\geq\a_{i-s}(K)$.
Similarly  $S_{i+1}-S_{s+1}+\b_s\geq
(S_{i+1}-S_i)/2+e=(R_{i+1}-R_i)/2+e\geq\a_{i-s}(K)$ Thus
$\a_i=\b_i=\a_{i-s}(K)$. \qed

\blm With the notations above if  $K\leq M$ and $\nn K\sb\nn M$ then
$K\leq N$.
\elm
\pf Note that $\nn K\sb\nn M$ means $T_1>R_1=R$ so $T_1\geq R+1$.
\vskip 3mm

{\bf Proof of 2.1(i).} If $i\geq s+1$ then either $T_i\geq R_i=S_i$
or $T_{i-1}+T_i\geq R_i+R_{i+1}=S_i+S_{i+1}$ so we are done. Suppose
now $i\leq s$. If $i$ is odd then, regardless of the type of $M$,
$T_i\geq T_1\geq R+1\geq S_i$ so we are done. If $i$ is even and $M$
is of type I or II then $T_i\geq T_2\geq T_1-2e\geq R-2e+1\geq S_i$ so
we are done. If $M$ is of type III then $T_{i-1}+T_i\geq T_1+T_2\geq
2T_1-2e\geq 2(R+1)-2e=(R-2e+2)+R=S_i+S_{i+1}$ so we are done.

Before proving 2.1(ii) and 2.1(iii) we note that for any $i\geq s$ we
have $[a_1,\ldots,a_i]\ap [b_1,\ldots,b_i]$ so $a_{1,i}=b_{1,i}$ in
$\ff/\fs$. Thus if $i\geq s$, $j$ is arbitrary and $\e\in\ff$ then
$d(\e a_{1,i}c_{1,j})=d(\e b_{1,i}c_{1,j})$. If $i\geq s+1$ then
also $\a_{i+1}=\b_{i+1}$ and we get $d[\e a_{1,i}c_{1,j}]=d[\e
b_{1,i}c_{1,j}]$. It follows that for any $i\geq s$ we have
$C_i=\min\{ (R_{i+1}-T_i)/2+e, R_{i+1}-T_i+d[-a_{1,i+1}c_{1,i-1}],
R_{i+1}+R_{i+2}-T_{i-1}-T_i+d[a_{1,i+2}c_{1,i-2}]\} =\min\{
(S_{i+1}-T_i)/2+e, S_{i+1}-T_i+d[-b_{1,i+1}c_{1,i-1}],
S_{i+1}+S_{i+2}-T_{i-1}-T_i+d[b_{1,i+2}c_{1,i-2}]\} =B_i$.

Also note that if $M$ is of type I or II then for $1<i<s$ even we have
$T_{i-1}\geq T_1\geq R+1=S_{i+1}$ and for $1<i<s$ odd we have
$T_{i-1}\geq T_2\geq T_1-2e\geq R-2e+1=S_{i+1}$. Thus $i$ is not
essential for $1<i\leq s-1$. By Lemmas 2.12 and 2.13 it follows that
2.1(ii) is vacuous for $1<i\leq s-2$ and 2.1(iii) is vacuous for
$1<i\leq s-1$. If $M$ is of type III then for $1<i\leq s$ even we have
$T_{i-1}\geq T_1>R=S_{i+1}$ and if $1<i<s$ is odd then
$T_{i-2}+T_{i-1}\geq T_1+T_2\geq 2T_i-2e\geq
2(R+1)-2e=(R-2e+2)+R=S_{i+1}+S_{i+2}$. So $i$ is not essential for
$1<i\leq s$. So this time 2.1(ii) is vacuous for $1<i\leq s-1$ and
2.1(iii) for $1<i\leq s$.
\vskip 3mm

{\bf Proof of 2.1(ii).} If $M$ is of type II or III then $S_2=R-2e$
resp. $R-2e+1$ so $B_1\leq (S_2-T_1)/2+e\leq
((R-2e+1)-(R+1))/2+e=0$. If $M$ is of type III then $B_i\leq
(S_2+S_3)/2-T_1+e\leq (R-2e+2+R)/2-(R+1)+e=0$. Hence 2.1(ii) is
trivial at $i=1$. Take now $i\geq s+1$. We have
$d[b_{1,i}c_{1,i}]=d[a_{1,i}c_{1,i}]\geq C_i=B_i$ so we are done. So
we are left with the cases $i=s$ and, if $M$ is of type I or II,
$i=s-1$. 

Take $i=s$. We have $d[b_{1,s}c_{1,s}]=\min\{
d(b_{1,s}c_{1,s}),\b_s,\c_s\}$. Now
$d(b_{1,s}c_{1,s})=d(a_{1,s}b_{1,s})\geq d[a_{1,s}c_{1,s}]$ and
$\c_s\geq d[a_{1,s}c_{1,s}]$. But $d[a_{1,s}c_{1,s}]\geq C_s=B_s$ so
we may assume that $d[b_{1,s}c_{1,s}]=\b_s$. If $M$ is of type I or II
then $S_{s+1}-S_s\geq (R+1)-(R-2e+1)=2e$ and if it is of type III then
$S_{s+1}-S_s=R-(R-2e+2)=2e-2$. It follows that
$\b_s=(S_{s+1}-S_s)/2+e$, which is $\geq B_s$ by 2.6.

Take now $i=s-1$, when $M$ is of type I or II. We have $T_{s-1}\geq
T_1\geq R+1$ so $B_{s-1}\leq (S_s-T_{s-1})/2+e\leq
((R-2e+1)-(R+1))/2+e=0$, which makes 2.1(ii) trivial. 
\vskip 3mm

{\bf Proof of 2.1(iii).} Let $1<i<n$ s.t. $S_{i+1}>T_{i-1}$ and
$B_{i-1}+B_i>2e+S_i-T_i$. It follows that $i$ is essential, which
rules out the cases $1<i\leq s-1$, if $M,N$ are of type I or II and
$1<i\leq s$ if $M$ is of type III. We prove that
$[c_1,\ldots,c_{i-1}]\rep [b_1,\ldots,b_i]$. 

If $i\geq s+1$ then $R_{i+1}=S_{i+1}>T_{i-1}$ and
$C_{i-1}+C_i=B_{i-1}+B_i>2e+S_i-T_i=2e+R_i-T_i$. It follows that
$[c_1,\ldots,c_{i-1}]\rep [a_1,\ldots,a_i]\ap [b_1,\ldots,b_i]$ and we
are done. So we are left with the case $i=s$ and only when $M$ is of
type I or II. Since $[b_1,\ldots,b_s]\ap [a_1,\ldots,a_s]$ it is
enough to prove that $[c_1,\ldots,c_{s-1}]\rep [a_1,\ldots,a_s]$
i.e. that $R_{s+1}>T_{s-1}$ and $C_{s-1}+C_s>2e+R_s-T_s$. We have
$R_{s+1}=S_{s+1}>T_{s-1}$ and $C_{s-1}+C_s>2e+R_s-T_s$ is equivalent
by Lemma 2.12 to
$d[-a_{1,s+1}c_{1,s-1}]+d[-a_{1,s}c_{1,s-2}]>2e+T_{s-1}-R_{s+1}$,
$d[-a_{1,s+1}c_{1,s-1}]>e+(T_{s-2}+T_{s-1})/2-R_{s+1}$ (if $s>2$) and
$d[-a_{1,s}c_{1,s-2}]>e+T_{s-1}-(R_{s+1}+R_{s+2})/2$ (if $s\leq
n-2$). But $R_{s+1}=S_{s+1}$, $R_{s+2}=S_{s+2}$ and
$d[-a_{1,s+1}c_{1,s-1}]=d[-b_{1,s+1}c_{1,s-1}]$ so the inequalities
above will follow from the corresponding inequalities for $N,K$
provided that $d[-a_{1,s}c_{1,s-2}]\geq d[-b_{1,s}c_{1,s-2}]$. We have
$d[-a_{1,s}c_{1,s-2}]=\min\{
d(-a_{1,s}c_{1,s-2}),\a_s,\c_{s-2}\}$. But
$d(-a_{1,s}c_{1,s-2})=d(-b_{1,s}c_{1,s-2})\geq d[-b_{1,s}c_{1,s-2}]$
and $\c_{s-2}\geq d[-b_{1,s}c_{1,i-s}]$ so it is enough to prove that
$\a_s\geq\b_s\geq d[-b_{1,i}c_{1,i-2}]$. We have $R_{s+1}=S_{s+1}\geq
R+1$ and $R_s=R-2e<R-2e+1=S_s$ so $R_{s+1}-R_s>S_{s+1}-S_s\geq 2e$,
which implies that $\a_s=(R_{s+1}-R_s)/2+e>(S_{s+1}-S_s)/2+e=\b_s$. 
\vskip 3mm

{\bf Proof of 2.1(iv).} Let $1<i\leq n-2$ s.t. $S_{i+1}+2e\leq
T_{i-1}+2e<S_{i+2}\leq T_i$. We want to prove that
$[c_1,\ldots,c_{i-1}]\rep [b_1,\ldots,b_{i+1}]$. If $i\geq s-1$ then
$T_{i-1}+2e<S_{i+2}=R_{i+2}$, which implies by Lemma 2.19 that
$[c_1,\ldots,c_{i-1}]\rep [a_1,\ldots,a_{i+1}]\ap
[b_1,\ldots,b_{i+1}]$. Suppose now that $i\leq s-2$. We have
$S_{i+1}+2e<S_{i+2}$ and $i+1<s$. But the only index $j<s$
s.t. $S_j+2e<S_{j+1}$ is $j=2$ and only if $M$ is of type II. So we
must have $i+1=2$, which contradicts our assumption that $1<i\leq
n-2$. \qed

\section{First element of a good BONG}

Throughout the sections $\S$8 and $\S$9 we will make extensive use of
[B3, Lemmas 2.1 and 2.4]. For convenience we state these results
here. 

In the formulas for $\a_i$ some of the terms $(R_{i+1}-R_i)/2+2+e$,
$R_{i+1}-R_j+d(-a_{j,j+1})$ for $j\leq i$ and
$R_{j+1}-R_i+d(-a_{j,j+1})$ for $j\geq i$ may be replaced by
expressions involving other $\a_h$ invariants. These expressions are
$\geq\a_i$ but $\leq$ than the terms they replace.

If $k\leq i<l$ then $\a_{i-k+1}(\[ a_k,\ldots,a_l\] )$ is the
minimum of the set
$$\{ (R_{i+1}-R_i)/2+e\}\cup\{ R_{i+1}-R_j+d(-a_{j,j+1})\mid k\leq
j\leq i\}\cup\{ R_{j+1}-R_i+d(-a_{j,j+1})\mid i\leq j<l\}$$
so it can replace $(R_{i+1}-R_i)/2+e$ and the terms corresponding to
$k\leq j<l$. 

More generally, if $k\leq h<l$ then:

If $h\leq i$ then $R_{i+1}-R_{h+1}+\a_h$ replaces all terms with
$j\leq h$ and $R_{i+1}-R_{h+1}+\a_{h-k+1}(\[ a_k,\ldots,a_l\] )$ those
with $k\leq h\leq i$. 

If $h\geq i$ then $R_i-R_h+\a_h$ replaces all terms with
$j\geq h$ and $R_i-R_h+\a_{h-k+1}(\[ a_k,\ldots,a_l\] )$ those with
$i\leq h<l$.  

Unlike we did in the proofs of Lemmas 7.15 and 7.19, we will leave to
the reader to determine what terms in the definitions of various
$\a_i$'s are replaced by every expression.
\vskip 2mm

In this section we will find out what numbers can be the length of a
first element in a good BONG of a lattice. Namely, if $M$ is a lattice
and $b_1\in\ff$ we want to find necessary and sufficient conditions
s.t. $M\ap\[ a'_1,\ldots,a'_n\]$ relative to some good BONG with
$a'_1=b_1$. It turns out that this problem is more difficult in the
case when $F/\QQ_2$ is totally ramified, i.e. when the residual field
$\oo/\p$ has only 2 elements. We encountered a similar situation in
[B1]. 

The differences in the proofs in the case when the extension
$F/\QQ_2$ is totally ramified (i.e. when $|\oo/\p |=2$) come from the
following two lemmas.

\blm (i) If $|\oo/\p |>2$ then for any $a\in\ff$ with $d(a)\neq 0,2e$
there is $b\in\ff$ s.t. $d(a)=d(b)=d(ab)$.

(ii) If $|\oo/\p |=2$ and $a,b\in\ff$ then $d(a)=d(b)=d(ab)$ iff
$a,b\in\fs$. Equivalently, if $d(a)=d(b)<\j$ then $d(ab)>d(a)$.
\elm
\pf If $d(a)=d(b)=\j$ then also $d(ab)=\j$. If $d(a)=d(b)=0$ then
$\ord a,\ord b$ are odd so $\ord ab$ is even and we have $d(ab)>0$. If
$d(a)=d(b)=2e$ then $a,b\in\D\fs$ so $ab\in\fs$ which implies
$d(ab)=\j$. These are in agreement with both (i) and (ii). From now on
we will assume that $d:=d(a)$ is odd and $1\leq d\leq 2e-1$. By
multiplying $a$ by a square, we may assume that $a=1+\pi^d\e$ for some
$\e\in\ooo$. 

For any $\a\in\oo$ we denote by $\widehat\a$ its class in $\oo/\p$.

(i) Since $|\oo/\p |>2$ there is $\eta\in\oo$
s.t. $\widehat\eta\notin\{\widehat 0,\widehat\e\}$. Hence
$\widehat\eta,\widehat{\e+\eta}\neq 0$, which implies that $\eta,\e
+\eta\in\ooo$. Take $b=1+\pi^d\eta$. We have $ab=1+\pi^d(\e +\eta
)+\pi^{2d}\e\eta$. Since $\eta,\e+ \eta\in\ooo$ we have $\ord\pi^d\eta
=\ord\pi^d(\e+ \eta )=d$. Since $\ord\pi^{2d}\e\eta=2d>d$ we also have
$\ord (\pi^d(\e +\eta )+\pi^{2d}\e\eta)=d$. But $d$ is odd and $<2e$
so $d(b)=d(ab)=d$.

(ii) Take $b\in\ooo$ s.t. $d(b)=d(a)=d$. By multiplying $b$ by a
square we may assume that $b=1+\pi^d\eta$ with $\eta\in\ooo$. Since
$\oo/\p =\{\widehat 0,\widehat 1\}$ and $\e,\eta\in\ooo$ we have
$\widehat\e =\widehat\eta =\widehat 1$ so $\widehat{\e+ \eta}=\widehat
0$, i.e. $\e +\eta\in\p$. It follows that $\ord\pi^d(\e +\eta
)>d$. Since also $2d>d$ we have $\ord (\pi^d(\e +\eta )+\pi^{2d}\e\eta
)>d$ and so $d(ab)=d(1+\pi^d(\e +\eta )+\pi^{2d}\e\eta )>d$. \qed

\blm Let $d,d'\in d(\ff )=\{ 0,1,3,\ldots,2e-1,2e,\j \}$ and let
$a\in\ff$ with $d(a)=d$.

(i) There is $b\in\ff$ with $d(b)=d'$ s.t. $(a,b)_\p =-1$ iff
$d+d'\leq 2e$.

(ii) If $|\oo/\p |>2$ then there is $b\in\ff$ with $d(b)=d'$
s.t. $(a,b)_\p =1$ iff $\{ d,d'\}\neq\{ 0,2e\}$.

(iii)If $|\oo/\p |=2$ then there is $b\in\ff$ with $d(b)=d'$
s.t. $(a,b)_\p =1$ iff $d+d'\neq 2e$.
\elm
\pf Let $d=a(a)$ and $d'=d(b)$. If $d+d'>2e$ then $(a,b)_\p =1$. If
$d=2e$ and $d'=0$ then $a\in\D\fs$ and $\ord b$ is odd so $(a,b)_\p
=-1$. Similarly when $d=0,d'=2e$. These are consistent with (i),(ii)
and (iii). From now on we assume that $d+d'\leq 2e$ and $\{
d,d'\}\neq\{ 0,2e\}$.  

We assume now that $d+d'=2e$ and $\{ d,d'\}\neq\{ 0,2e\}$ so
$0<d,d'<2e$.  Then (i) follows from [H, Lemma 3]. For (ii) we use
Lemma 8.1(i). There are $c,c'\in\ff$ with $d(c)=d(c')=d(cc')=d$. Then
at least one of $(a,c)_\p,(a,c')_\p,(a,cc')_\p$ is $1$. So we can
choose $b\in\{ c,c',cc'\}$ with $d(b)=d'$ such that $(a,b)_\p =1$. For
(iii) let $c\in\ff$ with $d(c)=d'$ s.t. $(a,c)_\p =-1$. Then for any
$b\in\ff$ with $d(b)=d(c)=d'$ we have by Lemma 8.1(ii) $d(bc)>d'$ so
$d(a)+d(bc)>d+d'=2e$, which implies $(a,bc)_\p =1$. Hence $(a,b)_\p
=(a,c)_\p =-1$. 

Assume now that $d+d'<2e$. Then there is $c\in\ff$ with $d(c)=2e-d$
s.t. $(a,c)_\p=-1$. (See the case $d+d'=2e$ of (i).) Take now
$c'\in\ff$ with $d(c')=d'$. Since $d(c)=2e-d>d'$ we have
$d(cc')=d'$. Also $(a,c)_\p=-1$ implies $(a,cc')_\p =-(a,c')_\p$. So
for any $\pm$ sign there is a good choice of $b\in\{ c',cc'\}$
s.t. $d(b)=d'$ and $(a,b)_\p =\pm 1$. This proves (i), (ii) and (iii)
in the case $d+d'<2e$. \qed

\blm If $M$ is quaternary, $R_1=R_3$, $R_2=R_4$ and $\e\in\ooo$ with
$d(\e )\geq\a_1$ then $M\ap\[ a'_1,a'_2,a'_3,a'_4\]$ relative to some
other good BONG s.t. $a'_1=\e a_i$.
\elm
\pf First we prove that if $N$ is another lattice with $FM\ap FN$ then
$M\ap N$ iff $R_i=S_i$ for all $i$'s, $\a_1=\b_1$ and
$d(a_1b_1)\geq\a_1$. The necessity follows from [B3, Theorem
3.1]. Conversely, condition (i) of [B3, Theorem 3.1] is just $R_i=S_i$
for all $i$'s. For condition (ii) we have $\a_1=\b_1$. But
$R_1+R_2=R_3+R_4$. By [B3, Corollary 2.3] we have
$R_1+\a_1=R_2+\a_2=R_3+\a_3$. Similarly
$S_1+\b_1=S_2+\b_2=S_3+\b_3$. Since $R_i=S_i$ for all $i$'s and
$\a_1=\b_1$ we have $\a_2=\b_2$ and $\a_3=\b_3$ as well. Thus (ii)
holds. We have $d(a_1b_1)\geq\a_1$ so (iii) is satisfied at
$i=1$. Since $R_1=R_3$ we have by Lemma 7.4(iii) $d(-a_{1,2})\geq
d[-a_{1,2}]=\a_2$ and $d(-a_{2,3})\geq d[-a_{2,3}]=\a_1$. Similarly
$d(-b_{1,2})\geq\b_2=\a_2$ and $d(-b_{2,3})\geq\b_1=\a_1$. It follows
that $d(a_{1,2}b_{1,2})\geq\a_2$ and
$d(a_{2,3}b_{2,3})\geq\a_1$. Since also $d(a_1b_1)\geq\a_1$ we have
$d(a_{1,3}b_{1,3})\geq\a_1=\a_3$. (We have $R_1+\a_1=R_3+\a_3$ and
$R_1=R_3$.) Thus (iii) holds. Finally, $R_1=R_3$ and $R_2=R_4$ imply
that $\a_1+\a_2,\a_2+\a_3\leq 2e$ so (iv) is vacuous. 

Next we prove that, assuming that $b_1,b_2,b_3,b_4\in\ff$ with
$S_i:=\ord b_i=R_i$ for all $i$'s, then there is a lattice $N\ap\[
b_1,b_2,b_3,b_4\]$ with $\a_1=\b_1$ iff $d(-b_{1,2})\geq\a_2$,
$d(-b_{2,3})\geq\a_1$ and $d(-b_{3,4})\geq\a_2$ and, if $\a_1\neq
(R_2-R_1)/2+e$, at least one of the three inequalities becomes
equality. First we prove that the three inequalities imply the the
existence of $N$ with $N\ap\[ b_1,b_2,b_3,b_4\]$. For $i=1,2$ we have
$S_i=R_i=R_{i+2}=S_{i+2}$ so we still need $b_{i+1}/b_i\in\aaa$
for $i=1,2,3$. We have $S_{i+1}-S_i=R_{i+1}-R_i\geq -2e$ so we
still need
$R_{i+1}-R_i+d(-b_{i,i+1})=S_{i+1}-S_i+d(-b_{i,i+1})\geq 0$. We
have $R_2-R_1+d(-b_{1,2})\geq R_2-R_1+\a_2=\a_1\geq 0$,
$R_3-R_2+d(-b_{2,3})\geq R_1-R_2+\a_1=\a_2\geq 0$ and
$R_4-R_3+d(-b_{3,4})\geq R_2-R_1+\a_2=\a_1\geq 0$. (We have
$R_1+\a_1=R_2+\a_2$.) For the equality $\a_1=\b_1$ note that
$\b_1=\min\{ (S_2-S_1)/2+e, S_2-S_1+d(-b_{1,2}),
S_3-S_1+d(-b_{2,3}), S_4-S_1+d(-b_{3,4})\} =\min\{ (R_2-R_1)/2+e,
R_2-R_1+d(-b_{1,2}), d(-b_{2,3}), R_2-R_1+d(-b_{3,4})\}$. So
$\a_1=\b_1$ iff $R_2-R_1+d(-b_{1,2}), d(-b_{2,3}),
R_2-R_1+d(-b_{3,4})$ are all $\geq\a_1$ and, if
$\a_1<(R_2-R_1)/2+e$, at least one of them is $=\a_1$. But
$R_2-R_1+d(-b_{1,2})\geq\a_1$ is equivalent to $d(-b_{1,2})\geq
R_1-R_2+\a_1=\a_2$ and, similarly, $R_2-R_1+d(-b_{3,4})\geq\a_1$ iff
$d(-b_{3,4})\geq\a_2$. 

We now start our proof. If we take $a'_i=\e a_i$ then $M^\e\ap\[
a'_1,a'_2,a'_3,a'_4\]$. We have $R_i(M^\e )=\ord a'_i=R_i$ and
$\a_i(M^\e )=\a_i$. ($\a_i$'s are invariant to scalling.) Also
$d(a_1a'_1)=d(\e )\geq\a_1$. So we still need $FM\ap FM^\e$. Now $\det
FM^\e =a'_{1,4}=a_{1,4}=\det FM$ (in $\ff/\fs$) so the only case when
$FM\ap FM^\e$ fails is when $S(FM^\e )=-S(FM)$. Suppose that this
happens. By [OM, 58:3] we have $S(FM^\e )=(\e,\det FM)_\p S(FM)$ so
$(\e,a_{1,4})_\p =-1$. Now $d(-a_{1,2}),d(-a_{3,4})\geq\a_2$. (We take
$b_i=a_i$ in the reasoning above.) Hence $d(a_{1,4})\geq\a_2$. It
follows that $2e\geq d(\e )+d(a_{1,4})\geq d(\e )+\a_2$. 

We claim that $\a_1$ is an odd integer. Otherwise $\a_1=(R_2-R_1)/2+e$
so $\a_1+\a_2=R_1-R_2+2\a_1=2e$. We have $2e\geq d(\e
)+d(a_{1,4})\geq\a_1+\a_2=2e$ so we must have equalities. Thus $d(\e
)=\a_1$, $d(a_{1,4})=\a_2$ and $d(\e )+d(a_{1,4})=2e$. Since $d(\e
)=\a_1$ is not odd one of $d(\e )+d(a_{1,4})=2e$ is $0$ an the other
is $2e$. But $\e\in\ooo$ and $\ord a_{1,4}=R_1+R_2+R_3+R_4=2R_1+2R_2$
is even so $d(\e ),d(a_{1,4})>0$. Contradiction. Thus $\a_1$ is odd
and so is $\a_2=R_1-R_2+\a_1$. (We have $R_1=R_3$ so $R_1\ev R_2\m2$.)
Also $R_1=R_3$ implies $\a_1+\a_2\leq 2e$ so $\a_1,\a_2<2e$. It
follows that $\a_1,\a_2\in d(\ooo )=\{ 1,3,\ldots,2e-1,2e,\j\}$.

We will show that there is a lattice $N$ s.t. $R_i=S_i$ at all $i$'s,
$\a_1=\b_1$, $a_1=b_1$, $\det FM=\det FN$ and $FM\not\ap FN$
(i.e. $S(FN)=-S(FM)$). Now $d(-a_{3,4})\geq\a_2$. We have 2 cases:

1. $d(-a_{3,4})=\a_2$. We take $b_1,b_2,b_3,b_4$ to be the sequence
$a_1,a_2,\eta a_3,\eta a_4$, where $\eta\in\ooo$ with $d(\eta
)=2e-\a_2=2e-d(-a_{3,4})$ s.t. $(\eta,-a_{3,4})_\p =-1$. We have
$S_i:=\ord b_i=R_i$ for all $i$'s. To prove that $N\ap\[
b_1,b_2,b_3,b_4\]$ exists and $\a_1=\b_1$ we note that
$d(-b_{1,2})=d(-a_{1,2})\geq\a_2$, $d(-b_{2,3})=d(-\eta
a_{2,3})\geq\a_1$ (we have $d(\eta )\geq 2e-\a_2\geq\a_1$ and
$d(-a_{2,3})\geq\a_1$) and $d(-b_{3,4})=d(-a_{3,4})=\a_2$. We have
$\det FN=a_{1,4}=\det FN$. Since $(\eta,-a_{3,4})_\p =-1$ we have
$[b_3,b_4]=[\eta a_3,\eta a_4]\not\ap [a_3,a_4]$. It follows that
$FN\ap [a_1,a_2,b_3,b_4]\not\ap [a_1,a_2,a_3,a_4]\ap FM$. So $N$
satisfies all the required conditions.

2. $d(-a_{3,4})>\a_2$. Take $\eta\in\ooo$ with $d(\eta )=\a_2$.
Let $b_1,b_2,b_3,b_4$ be the sequence $a_1,\eta a_2,\eta a_3,a_4$.
We have $S_i:=\ord b_i=R_i$ for all $i$'s. To prove that $N\ap\[
b_1,b_2,b_3,b_4\]$ exists and $\a_1=\b_1$ we note that
$d(-b_{1,2})=d(-\eta a_{1,2})\geq\a_2$,
$d(-b_{2,3})=d(-a_{2,3})\geq\a_1$ and $d(-b_{3,4})=d(-\eta
a_{3,4})=\a_2$. (We have $d(-a_{1,2})\geq\a_2$, $d(-a_{3,4})>\a_2$ and
$d(\eta )=\a_2$.) We have $\det FM=\det FN$. If also $FM\not\ap FN$
then $N$ satisfies all the required conditions. If $FM\ap FN$ then
$M\ap N$. So $M\ap\[ b_1,b_2,b_3,b_4\]$ relative to some other
BONG. By this change of BONGs $a_1$ is not changed (we have $a_1=b_1$)
but $d(-a_{3,4})$ becomes $d(-b_{3,4})=\a_2$. So after this change of
BONGs we are in the case 1., already studied.

We take $a'_i=\e b_i$. We have $N^\e\ap\[ a'_1,a'_2,a'_3,a'_4\]$. In
particular, $a'_1=\e b_1=\e a_1$, as required. We have $\ord a'_i=R_i$
and $\a_i(N^\e )=\a_i(N)=\a_i$. In particular, $\a_1(N^\e )=\a_1$. So
we still need $FM\ap FN^\e$. Now $\det FN^\e =\det FN=\det FM$. Also
by [OM, 58:3] we have $S(FN^\e )=(\e,\det FN)_\p S(FN)=(\e,\det FM)_\p
S(FN)=-S(FN)=S(FM)$. Thus $FN^\e\ap FM$. \qed

\blm Let $L$ be as usual and let $1\leq i\leq j<n$.

(i) If $R_i+\a_i=R_j+\a_j$ then $R_i+\a_i=\ldots=R_j+\a_j$. If
$-R_{i+1}+\a_i=-R_{j+1}+\a_j$ then
$-R_{i+1}+\a_i=\ldots=-R_{j+1}+\a_j$. 

(ii) If $R_i+\a_i=R_j+\a_j$ and $\a_k=(R_{k+1}-R_k)/2+e$ for some
$i\leq k\leq j$ then $R_i+R_{i+1}=R_k+R_{k+1}$ and
$\a_l=(R_{l+1}-R_l)/2+e$ holds for any $i\leq l\leq k$.

(iii) If $-R_{i+1}+\a_i=-R_{j+1}+\a_j$ and $\a_k=(R_{k+1}-R_k)/2+e$
for some $i\leq k\leq j$ then $R_k+R_{k+1}=R_j+R_{j+1}$ and
$\a_l=(R_{l+1}-R_l)/2+e$ holds for any $k\leq l\leq j$.
\elm
\pf (i) It follows from the fact that the sequences $(R_i+\a_i)$ and
$(-R_{i+1}+\a_i)$ are increasing resp. decreasing. 

(ii) We have $\a_i\leq (R_{i+1}-R_i)/2+e$ so $(R_i+R_{i+1})/2+e\geq
R_i+\a_i=R_k+\a_k=(R_k+R_{k+1})/2+e$ so $R_i+R_{i+1}\geq
R_k+R_{k+1}$. But $i\leq k$ so we must have
$R_i+R_{i+1}=R_k+R_{k+1}$. Together with $\a_k=(R_{k+1}-R_k)/2+e$,
this implies that $\a_l=(R_{l+1}-R_l)/2+e$ holds for any $i\leq l\leq
k$ by [B3, Corollary 2.3(iii)].

(iii) We have $\a_j\leq (R_{j+1}-R_j)/2+e$ so
$-(R_k+R_{k+1})/2+e=-R_{k+1}+\a_k=-R_{j+1}+\a_j\leq
-(R_j+R_{j+1})/2+e$ so $R_k+R_{k+1}\geq R_j+R_{j+1}$. But $k\leq j$ so 
we must have $R_k+R_{k+1}=R_j+R_{j+1}$. Together with
$\a_k=(R_{k+1}-R_k)/2+e$, this implies that $\a_l=(R_{l+1}-R_l)/2+e$
holds for any $k\leq l\leq j$ by [B3, Corollary 2.3(iii)]. \qed

\blm Let $A:=\{ 1\leq i<n\mid\a_i<(R_{i+1}-R_i)/2+e\}$ and 
$$\begin{aligned} B:=\{ 1\leq i<n\mid R_i+\a_i<R_j+\a_j\text{ whenever }
R_i+R_{i+1}<R_j+R_{j+1},\\ 
-R_{i+1}+\a_i<-R_{j+1}+\a_j\text{ whenever
}R_i+R_{i+1}>R_j+R_{j+1}\}.
\end{aligned}$$
Let $C=\{ i\in A\cap B\mid\a_i=R_{i+1}-R_i+d(-a_{i,i+1})\}$. 

Then for any $i\in A$ there is $j\in C$ s.t. either $1\leq j\leq i$ and
$\a_i=R_{i+1}-R_{j+1}+\a_j=R_{i+1}-R_j+d(-a_{j,j+1})$ or $i\leq j<n$
and $\a_i=R_j-R_i+\a_j=R_{j+1}-R_i+d(-a_{j,j+1})$.

If moreover $i\in B$ then $j$ satisfies $R_i+R_{i+1}=R_j+R_{j+1}$.
\elm
\pf Note that if $1\leq i,k<n$ and $R_i+R_{i+1}=R_k+R_{k+1}$ then
$R_i+\a_i=R_k+\a_k$ and $-R_{i+1}+\a_i=-R_{k+1}+\a_k$ so $i\in B$ iff
$k\in B$. Also $\a_i=(R_{i+1}-R_i)/2+e$ iff $\a_k=(R_{k+1}-R_k)/2+e$
so $i\in A$ iff $k\in A$. (See [B3, Corollary 2.3].)

Recall that the sequences $(R_i+R_{i+1})$ and $(R_i+\a_i)$  are
increasing, while $(-R_{i+1}+\a_i)$ is decreasing. We have two cases:

1. $i\in A\cap B$. Since $\a_i<(R_{i+1}-R_i)/2+e$
we have $\a_i=R_{i+1}-R_j+d(-a_{j,j+1})$ for some $1\leq j\leq i$ or
$\a_i=R_{j+1}-R_i+d(-a_{j,j+1})$ for some $i\leq j<n$. We show that
in both cases $R_i+R_{i+1}=R_j+R_{j+1}$ and
$\a_j=R_{j+1}-R_j+d(-a_{j,j+1})$. Indeed, in the first case we have
$\a_i=R_{i+1}-R_j+d(-a_{j,j+1})=R_{i+1}-R_{j+1}+R_{j+1}-R_j+d(-a_{j,j+1})\geq
R_{i+1}-R_{j+1}+\a_j\geq\a_i$ so we have equalities. In particular,
$\a_j=R_{j+1}-R_j+d(-a_{j,j+1})$. Also
$-R_{i+1}+\a_i=-R_{j+1}+\a_j$. Since $i\in B$ we cannot have
$R_i+R_{i+1}>R_j+R_{j+1}$ so $R_i+R_{i+1}=R_j+R_{j+1}$. In the second
case we have
$\a_i=R_{j+1}-R_i+d(-a_{j,j+1})=R_j-R_i+R_{j+1}-R_j+d(-a_{j,j+1})\geq
R_j-R_i+\a_j\geq\a_i$ so we have equalities. In particular,
$\a_j=R_{j+1}-R_j+d(-a_{j,j+1})$. Also $R_i+\a_i=R_j+\a_j$. Since
$i\in B$ we cannot have $R_i+R_{i+1}<R_j+R_{j+1}$ so
$R_i+R_{i+1}=R_j+R_{j+1}$. Now since $i\in A\cap B$ and
$R_i+R_{i+1}=R_j+R_{j+1}$ we also have $j\in A\cap B$. Since also
$\a_j=R_{j+1}-R_j+d(-a_{j,j+1})$ we have $j\in C$, as claimed.

2. $i\in A\setminus B$. Since $i\in A$ we have
$\a_i<(R_{i+1}-R_i)/2+e$. Since $i\notin B$ we have
$R_i+R_{i+1}>R_k+R_{k+1}$ and $-R_{i+1}+\a_i=-R_{k+1}+\a_k$ or
$R_i+R_{i+1}<R_k+R_{k+1}$ and $R_i+\a_i=R_k+\a_k$ for some $k$.

In the first case we have $i>k$. WLOG we may assume that $k$ is
minimal with the property that $-R_{k+1}+\a_k=-R_{i+1}+\a_i$ so
$-R_{l+1}+\a_l>-R_{k+1}+\a_k$ for $l<k$. We cannot have
$\a_k=(R_{k+1}-R_k)/2+e$ since this would imply
$\a_i=(R_{i+1}-R_i)/2+e$ by Lemma 8.4(iii). Thus $k\in A$. Suppose
that $k\notin B$. We have $-R_{l+1}+\a_l>-R_{k+1}+\a_k$ whenever
$l<k$, in particular, whenever $R_l+R_{l+1}<R_k+R_{k+1}$. Thus
$k\notin B$ implies that $R_k+\a_k=R_l+\a_l$ and
$R_k+R_{k+1}<R_l+\a_l$ for some $l$. Let $h:=\min\{ l,i\}$. Since
$R_k+R_{k+1}<R_i+R_{i+1}$ and $R_k+R_{k+1}<R_l+R_{l+1}$ we have
$R_k+R_{k+1}<R_h+R_{h+1}$ so $k<h$. By Lemma 8.4(i), since
$R_k+\a_k=R_l+\a_l$ and $k<h\leq l$ we have $R_k+\a_k=R_h+\a_h$. Also
since $-R_{k+1}+\a_k=-R_{i+1}+\a_i$ and $k<h\leq i$ we have
$-R_{k+1}+\a_k=-R_{h+1}+\a_h$. By subtracting we get
$R_k+R_{k+1}=R_h+R_{h+1}$. Contradiction. Thus $k\in A\cap B$. By the
case 1. there is $j\in C$ with $R_k+R_{k+1}=R_j+R_{j+1}$. By [B3,
Corollary 2.3(i)] $-R_{j+1}+\a_j=-R_{k+1}+\a_k=-R_{i+1}+\a_i$ so
$\a_i=R_{i+1}-R_{j+1}+\a_j=R_{i+1}-R_j+d(-a_{j,j+1})$. (We have $j\in
C$ so $\a_j=R_{j+1}-R_j+d(-a_{j,j+1})$.)

In the second case we have $i<k$ and WLOG we may assume that $k$ is
maximal with the property that $R_i+\a_i=R_k+\a_k$. So
$R_l+\a_l>R_k+\a_k$ for $l>k$. We cannot have $\a_k=(R_{k+1}-R_k)/2+e$
since this would imply $\a_i=(R_{i+1}-R_i)/2+e$ by Lemma 8.4(ii). Thus
$k\in A$. Suppose that $k\notin B$. We have $R_l+\a_l>R_k+\a_k$
whenever $l>k$, in particular, whenever
$R_l+R_{l+1}>R_k+R_{k+1}$. Thus $k\notin B$ implies that
$-R_{k+1}+\a_k=-R_{l+1}+\a_l$ and $R_k+R_{k+1}>R_l+R_{l+1}$ for some
$l$.  Let $h:=\max\{ l,i\}$. Since $R_k+R_{k+1}>R_i+R_{i+1}$ and
$R_k+R_{k+1}>R_l+R_{l+1}$ we have $R_k+R_{k+1}>R_h+R_{h+1}$ so
$k>h$. By Lemma 8.4(i), since $-R_{k+1}+\a_k=-R_{l+1}+\a_l$ and
$k>h\geq l$ we have $-R_{k+1}+\a_k=-R_{h+1}+\a_h$. Also since
$R_k+\a_k=R_i+\a_i$ and $k>h\geq i$ we have $R_k+\a_k=R_h+\a_h$. By
subtracting we get $R_k+R_{k+1}=R_h+R_{h+1}$. Contradiction. Thus
$k\in A\cap B$. By the case 1. there is $j\in C$ with
$R_k+R_{k+1}=R_j+R_{j+1}$. By [B3, Corollary 2.3(i)]
$R_j+\a_j=R_k+\a_k=R_i+\a_i$ so
$\a_i=R_j-R_i+\a_j=R_{j+1}-R_i+d(-a_{j,j+1})$. (We have $j\in C$ so
$\a_j=R_{j+1}-R_j+d(-a_{j,j+1})$.) \qed 

\blm Let $M\ap\[ a_1,\ldots,a_n\]$ and let $b_1,\ldots,b_n\in\ff$
s.t. $S_i:=\ord b_i=R_i$ and $a_{1,n}b_{1,n}\in\fs$. If
$d(a_{1,i}b_{1,i})\geq\a_i$ then:

(i) There is a lattice $N\ap\[ b_1,\ldots,b_n\]$.

(ii) $\a_i\leq\b_i$ for $1\leq i\leq n-1$.

(iii) If $M$ has property A then $\a_i=\b_i$ for $1\leq i\leq n-1$. 

(Recall that property A is equivalent to $R_i<R_{i+2}$ for $1\leq
i\leq n-2$. See [B1, Definition 7, Corllary 4.2(i) and Lemma 4.3(i)].)
\elm
\pf (i) We have $S_i=R_i\leq R_{i+2}=S_{i+2}$ for $1\leq i\leq
n-2$. We still need $b_{i+1}/b_i\in\aaa$ for $1\leq i\leq n-1$. We
have $S_{i+1}-S_i=R_{i+1}-R_i\geq -2e$ so we still have to prove that
$S_{i+1}-S_i+d(-b_{i,i+1})=R_{i+1}-R_i+d(-b_{i,i+1})\geq 0$. In fact
we will show that $R_{i+1}-R_i+d(-b_{i,i+1})\geq\a_i$. We have
$R_{i+1}-R_i+d(-a_{i,i+1})\geq\a_i$ so, by the domination principle,
it is enough to prove that $R_{i+1}-R_i+d(a_{1,i-1}b_{1,i-1}),
R_{i+1}-R_i+d(a_{1,i+1}b_{1,i+1})\geq\a_i$. We have
$R_{i+1}-R_i+d(a_{1,i-1}b_{1,i-1})=\j$ if $i=1$ and
$R_{i+1}-R_i+d(a_{1,i-1}b_{1,i-1})\geq R_{i+1}-R_i+\a_{i-1}\geq\a_i$
otherwise. We have $R_{i+1}-R_i+d(a_{1,i+1}b_{1,i+1})=\j$ if $i=n-1$
and $R_{i+1}-R_i+d(a_{1,i+1}b_{1,i+1})\geq
R_{i+1}-R_i+\a_{i+1}\geq\a_i$ otherwise.

As a consequence, $d(-b_{i,i+1})\geq R_i-R_{i+1}+\a_i$ for $1\leq
i\leq n-1$.

(ii) Let $1\leq i\leq n-1$. We have
$(S_{i+1}-S_i)/2+e=(R_{i+1}-R_i)/2+e\geq\a_i$. For $j\leq i$ we have
$S_{i+1}-S_j+d(-b_{j,j+1})\geq
R_{i+1}-R_j+R_j-R_{j+1}+\a_j=R_{i+1}-R_{j+1}+\a_j\geq\a_i$. For $j\geq
i$ we have $S_{j+1}-S_i+d(-b_{j,j+1})\geq
R_{j+1}-R_i+R_j-R_{j+1}+\a_j=R_j-R_i+\a_j\geq\a_i$. By the definition
of $\b_i$ we have $\b_i\geq\a_i$.

(iii) We use the notations of Lemma 8.5. Let $1\leq i\leq n-1$. If
$\a_i=(R_{i+1}-R_i)/2+e$ then, since $\a_i\leq\b_i\leq
(S_{i+1}-S_i)/2+e=(R_{i+1}-R_i)/2+e$, we have $\a_i=\b_i$. If
$\a_i<(R_{i+1}-R_i)/2+e$ then by Lemma 8.5 there is $j\in C$ s.t. either
$1\leq j\leq i$ and
$\a_i=R_{i+1}-R_{j+1}+\a_j=R_{i+1}-R_j+d(-a_{j,j+1})$ or $i\leq j\leq
n-1$ and $\a_i=R_j-R_i+\a_j=R_{j+1}-R_i+d(-a_{j,j+1})$. 

We have $j\in C$ so $\a_j=R_{j+1}-R_j+d(-a_{j,j+1})$. We claim
that $d(-a_{j,j+1})=d(-b_{j,j+1})$, i.e. that
$\a_j=R_{j+1}-R_j+d(-b_{j,j+1})$. By the domination principle it is
enough to prove that $\a_j<R_{j+1}-R_j+d(a_{1,j-1}b_{1,j-1})$ and
$\a_j<R_{j+1}-R_j+d(a_{1,j+1}b_{1,j+1})$. The two conditions are
trivial if $j=1$, resp. $j=n-1$, since
$d(a_{1,0}b_{1,0})=d(a_{1,n}b_{1,n})=\j$. We have $-R_{j+1}+\a_j\leq
-R_j+\a_{j-1}$ (if $j>1$) and $R_j+\a_j\leq R_{j+1}+\a_{j+1}$ (if
$j<n-1$). But since $j\in C$ the two inequalities are strict. (From
property A we have $R_{j-1}+R_j<R_j+R_{j+1}$ and
$R_j+R_{j+1}<R_{j+1}+R_{j+2}$.) So if $i>1$ then
$R_{j+1}-R_j+d(a_{1,j-1}b_{1,j-1})\geq R_{j+1}-R_j+\a_{j-1}>\a_j$ and
if $j<n-1$ then Thus $R_{j+1}-R_j+d(a_{1,j+1}b_{1,j+1})\geq
R_{j+1}-R_j+\a_{j+1}>\a_j$. 

So $d(-a_{j,j+1})=d(-b_{j,j+1})$. It follows that either $j\leq i$ and
$\a_i\leq\b_i\leq
S_{i+1}-S_j+d(-b_{j,j+1})=R_{i+1}-R_j+d(-a_{j,j+1})=\a_i$ or $j\geq i$
and $\a_i\leq\b_i\leq
S_{j+1}-S_i+d(-b_{j,j+1})=R_{j+1}-R_i+d(-a_{j,j+1})=\a_i$. Hence
$\a_i=\b_i$. \qed 

\bff Recall some consequences that follow whenever
$R_{i-1}=R_{i+1}$. 

We have $R_{i-1}\ev R_i\ev R_{i+1}\m2$ by Lemma 6.6(i).

We have $R_{i-1}+\a_{i-1}=R_i+\a_i$ by [B3, Corollary 2.3(i)] so
$\a_i=R_{i-1}-R_i+\a_{i-1}=R_{i+1}-R_i+\a_{i-1}$ and
$\a_{i-1}=R_i-R_{i-1}+\a_i=R_i-R_{i+1}+\a_i$. In particular,
$\a_{i-1}\ev\a_i\m2$. (We have $R_i\ev R_{i+1}\m2$.)

By [B3, Corollary 2.3(iii)] we have $\a_{i-1}=(R_i-R_{i-1})/2+e$ iff
$\a_i=(R_{i+1}-R_i)/2+e$. Also note that $\a_{i-1}+\a_i\leq
(R_i-R_{i-1})/2+e+(R_{i+1}-R_i)/2+e=2e$. Moreover $\a_{i-1}+\a_i=2e$
iff $\a_{i-1}=(R_i-R_{i-1})/2+e$ and $\a_i=(R_{i+1}-R_i)/2+e$. In
conclusion, $\a_{i-1}=(R_i-R_{i-1})/2+e$ iff $\a_i=(R_{i+1}-R_i)/2+e$
iff $\a_{i-1}+\a_i=2e$.

By Lemma 7.4(iii) we have $d[-a_{i-1,i}]=\a_i$ and
$d[-a_{i,i+1}]=\a_{i-1}$. In particular, $d(-a_{i-1,i})\geq\a_i$ and
$d(-a_{i,i+1})\geq\a_{i-1}$.
\eff

\blm (i) There is $\e\in\ooo$ with $d(\e )=\a_1$ s.t. $M\ap\[
a'_1,\ldots,a'_n\]$ relative to some good BONG with $a'_1=\e a_1$
unless $\a_1=(R_2-R_1)/2+e$ and one of the following happens:

(a) $\a_1\notin d(\ooo )$ i.e. $\a_i\notin\{ 1,3,\ldots,2e-1,2e\}$.

(b) $|\oo/\p |=2$, $d[-a_{1,2}]=e-(R_2-R_1)/2$ and
$\a_2>e-(R_2-R_1)/2$ (if $n\geq 3$).

(c) $|\oo/\p |=2$, $R_1=R_3$, $\a_3>(R_2-R_1)/2+e$ (if $n\geq 4$)
and $[a_1,a_2,a_3]$ is not isotropic. \elm
\pf Note that $d[-a_{1,2}]=\min\{ d(-a_{1,2}),\a_2\}$ (with $\a_2$
ignored if $n=2$). So the condition $d[-a_{1,2}]=e-(R_2-R_1)/2$ and
$\a_2>e-(R_2-R_1)/2$ (if $n\geq 3$) from (b) is equivalent to
$d(-a_{1,2})=e-(R_2-R_1)/2$ and $\a_2>e-(R_2-R_1)/2$ (if $n\geq
3$).

Suppose that $M\ap\[ a'_1,\ldots,a'_n\]$ relative to some good BONG
with $a'_1=\e a_1$ and $d(\e )=\a_1$. Suppose that
$\a_1=(R_2-R_1)/2+e$. We have to prove that none of (a)-(c) happens.

We have $\a_1=d(\e )\in d(\ooo )$ so (a) doesn't happen.

Suppose that (b) happens. This implies
$d(-a_{1,2})=e-(R_2-R_1)/2$. (See
above.) We also have $\a_1+\a_2>(R_2-R_1)/2+e+e-(R_2-R_1)/2=2e$ (if
$n\geq 3$) so $a'_1\rep [a_1,a_2]$. (The same happens obviously if
$n=2$, when $FM\ap [a_1,a_2]$.) It follows that $(a'_1a_1,-a_{1,2})_\p
=1$. But $d(a'_1a_1)=d(\e )=\a_1=(R_2-R_1)/2+e$ and
$d(-a_{1,2})=e-(R_2-R_1)/2$ so $d(a'_1a_1)+d(-a_{1,2})=2e$. Since
$|\oo/\p |=2$ we get by Lemma 8.2(iii) $(a'_1a_1,-a_{1,2})_\p
=-1$. Contradiction.

Suppose that (c) happens. For consequences of $R_1=R_3$ we refer to
8.7. Since $\a_1=(R_2-R_1)/2+e$ we have
$\a_2=(R_3-R_2)/2+e=e-(R_2-R_1)/2$ and $\a_1+\a_2=2e$. It follows that
$d(-a_{1,2})\geq\a_2=e-(R_2-R_1)/2$ and
$d(-a_{2,3})\geq\a_1=(R_2-R_1)/2+e$. If any of these inequalities is
strict then $d(-a_{1,2})+d(-a_{2,3})>2e$ so $(-a_{1,2},-a_{2,3})_\p
=1$ so $[a_1,a_2,a_3]$ is isotropic. Contradiction. Thus
$d(-a_{1,2})=e-(R_2-R_1)/2$ and $d(-a_{2,3})=(R_2-R_1)/2+e$. If $n\geq
4$ then $\a_3>(R_2-R_1)/2+e=\a_1$ so $\a_2+\a_3>\a_1+\a_2=2e$. It
follows that $[a'_1,a'_2]\rep [a_1,a_2,a_3]$. (Same happens if $n=3$,
when $FM\ap [a_1,a_2,a_3]$.) It follows that
$[a'_1,a'_2,a_{1,3}a'_{1,2}]\ap [a_1,a_2,a_3]$ so it is
anisotropic. Hence $(-a'_{1,2},-a_{1,3}a'_1)_\p =-1$. But
$d(-a'_{1,2})\geq\a_2$. Also $d(a'_1a_1)=d(\e )=\a_1=d(-a_{2,3})$ so
by Lemma 8.1(ii) we have $d(-a_{1,3}a'_1)>\a_1$. It follows that
$d(-a'_{1,2})+d(-a_{1,3}a'_1)>\a_1+\a_2=2e$ so
$(-a'_{1,2},-a_{1,3}a'_1)_\p =1$. Contradiction.

For the reverse implication we use induction on $n$. We have to prove
that if either $\a_1<(R_2-R_1)/2+e$ or none of (a)-(c) holds then
$M\ap\[ a'_1,\ldots,a'_n\]$ relative to some good BONG and $a'_1=\e
a_1$ for some $\e\in\ooo$ with $d(\e )=\a_1$.

Take first $n=2$. We want to prove that there is $\e\in g(-a_2/a_1)$
with $d(\e )=\a_1$. If we do so then $M\ap\[\e a_1,\e a_2\]$ so we can
take $a'_1=\e a_1$ and $a'_2=\e a_2$ and we are done. By [B3, Lemma
  1.1] we have $g(a_2/a_1)=\upon{\a_1}{-a_{1,2}}$. So we have to find
some $\e\in\ooo$ with $d(\e )=\a_1$ and $(\e,-a_{1,2})_\p =1$. 

If $\a_1<(R_2-R_1)/2+e$ then $\a_1<2e$ is odd so $\a_1\in d(\ooo
)$. Also $\a_1=R_2-R_1+d(-a_{1,2})$, which implies
$d(-a_{1,2})=R_1-R_2+\a_1$. Since $\a_1\in d(\ooo )$ and
$\a_1+d(-a_{1,2})=R_1-R_2+2\a_1<2e$, by Lemma 8.2(ii) and (iii) there
is $\e\in\ooo$ with $d(\e )=\a_1$ and $(\e,-a_{1,2})_\p =1$ so we are
done. 

Suppose now that $\a_1=(R_2-R_1)/2+e$ and neither (a) nor (b)
holds ((c) cannot hold if $n=2$). It follows that $\a_1\in
d(\ooo )$ and either $|\oo/\p |>2$ or $d(-a_{1,2})=d[-a_{1,2}]\neq
e-(R_2-R_1)/2$. Note that $d(-a_{1,2})\geq
R_1-R_2+\a_1=e-(R_2-R_1)/2$. If $d(-a_{1,2})>e-(R_2-R_1)/2$ then
$g(a_2/a_1)=\upo{(R_2-R_1)/2+e}=\upo{\a_1}$. (See [B1, Definition 6,
case (iii)].) Take $\e\in\ooo$ with $d(\e )=\a_1$. We have
$\e\in\upo{\a_1}=g(a_2/a_1)$ so we are done. If
$d(-a_{1,2})=e-(R_2-R_1)/2$ we must have $|\oo/\p |>2$. Since
$|\oo/\p |>2$ and $\a_1\in d(\ooo )$ there is by Lemma 8.2(iii) some
$\e\in\ooo$ s.t. $d(\e )=\a_1$ and $(\e,-a_{1,2})_\p =1$ and we are done. (We
cannot have $\a_1=0$ and $d(-a_{1,2})=2e$ since $\a_1\in d(\ooo )$. We cannot
have $\a_1=2e$ and $d(-a_{1,2})=0$ since these would imply that both
$R_2-R_1=2e$ and $R_1+R_2=\ord a_{1,2}$ is odd.) 

For $n\geq 3$ we use the induction. If we denote $M^*\ap\[
a_2,\ldots,a_n\]$ then $\a_1=\min\{\a_1(\[ a_1,a_2\]
),R_2-R_1+\a_1(M^*)\} =\min\{
(R_2-R_1)/2+e,R_2-R_1+d(-a_{1,2}),R_2-R_1+\a_1(M^*)\}$.

Take first the case when $\a_1<(R_2-R_1)/2+e$. If
$d(-a_{1,2})\leq\a_1 (M^*)$ then $\a_1=\min\{
(R_2-R_1)/2+e,R_2-R_1+d(-a_{1,2})\} =\a_1(\[ a_1,a_2\] )$. By the
binary case there is $\e\in\ooo$ with $d(\e )=\a_1(\[ a_1,a_2\]
)=\a_1$ s.t. $\[ a_1,a_2\]\ap\[\e a_1,\e a_2\]$. It follows that
$M\ap\[ \e a_1,\e a_2,a_3,\ldots,a_n\]$ relative to some good BONG
and we are done. If $d(-a_{1,2})>\a_1(M^*)$ then
$\a_1=R_2-R_1+\a_1(M^*)$. Note that we cannot have
$\a_1(M^*)=(R_3-R_2)/2+e$ since this would imply that
$\a_1(M^*)\geq (R_1-R_2)/2+e$ so $\a_1=R_2-R_1+\a_1(M^*)\geq
(R_2-R_1)/2+e$. Since $\a_1(M^*)<(R_3-R_2)/2+e$ we have by the
induction hypothesis that $M^*\ap\[ a'_2,\ldots,a'_n\]$ relative to a
good BONG with $a'_2=\eta a_2$ and $d(\eta )=\a_1(M^*)$. It
follows that $M\ap\[ a_1,a'_2,\ldots,a'_n\]$ relative to a good BONG.
By this change of BONGs $d(-a_{1,2})$ is replaced by
$d(-a_1a'_2)=d(-\eta a_{1,2})$. But $d(\eta
)=\a_1(M^*)>d(-a_{1,2})$ so $d(-a_1a'_2)=\a_1(M^*)$. Hence we are
back to the previous case ($d(-a_{1,2})\leq\a_1 (M^*)$) and the
proof follows the same.

Suppose now that $\a_1=(R_2-R_1)/2+e$ and none of (a)-(c) happens. In
particular, $\a_1\in d(\ooo )$. We have $(R_2-R_1)/2+e\geq\a_1(\[
a_1,a_2\] )\geq\a_1$ so $\a_1(\[ a_1,a_2\] )=\a_1=(R_2-R_1)/2+e$. If
$|\oo/\p |>2$ or $d(-a_{1,2})\neq e-(R_2-R_1)/2$ then by the binary
case we have $\[ a_1,a_2\]\ap\[\e a_1,\e a_2\]$ for some $\e\in\ooo$
with $d(\e )=\a_1(\[ a_1,a_2\] )=\a_1$. It follows that
$M\ap\[ \e a_1,\e a_2,a_3,\ldots,a_n\]$ relative to some good BONG and
we are done. So we may assume that $|\oo/\p |=2$ and
$d(-a_{1,2})=e-(R_2-R_1)/2$. Since also $\a_2\geq
R_1-R_2+\a_1=e-(R_2-R_1)/2$ we have $d[-a_{1,2}]=e-(R_2-R_1)/2$. Since
(b) doesn't happen we cannot have $\a_2>e-(R_2-R_1)/2$. Hence
$\a_2=e-(R_2-R_1)/2$.

For any $2\leq i\leq n-1$ we have $\a_i=\min\{
R_{i+1}-R_1+d(-a_{1,2}),\a_{i-1}(M^*)\}$. But
$R_{i+1}-R_1+d(-a_{1,2})=R_{i+1}-(R_1+R_2)/2+e\geq
R_{i+1}-(R_i+R_{i+1})/2+e=(R_{i+1}-R_i)/2+e\geq\a_{i-1}(M^*)$. So
$\a_i=\a_{i-1}(M^*)$. In particular, $\a_1(M^*)=\a_2=e-(R_2-R_1)/2$.

If $M^*\ap\[ a'_2,\ldots,a'_n\]$ relative to a good BONG with
$a'_2=\eta a_2$ and $\eta\in\ooo$, $d(\eta
)=\a_1(M^*)=e-(R_2-R_1)/2$ then $M\ap\[ a_1,a'_2,\ldots,a'_n\]$
relative to some good BONG. Note that $a_1$ is preserved by this
change of BONGs but $d(-a_{1,2})$ is replaced by $d(-a_1a'_2)=d(\eta
a_{1,2})$. Since $d(-a_{1,2})=d(\eta )=e-(R_2-R_1)/2$ we have by
Lemma 8.1(ii) $d(\eta a_{1,2})>e-(R_2-R_1)/2$. Since after this change
of BONGs we no longer have $d(-a_{1,2})=e-(R-2-R_1)/2$ we are done by
the previous reasoning. Thus we may assume that such writing $M^*\ap\[
a'_2,\ldots,a'_n\]$ doesn't exist. By  the induction step we have
$e-(R_2-R_1)/2=\a_2=(R_3-R_2)/2+e$, i.e. $R_1=R_3$, and also one of
the conditions (a)-(c) holds for $M^*$.

If (a) holds for $M^*$ note that $\a_1\in d(\ooo )\setminus\{\j\}$ and
$\a_1(M^*)=e-(R_2-R_1)/2=2e-\a_1$. Then $\a_1(M^*)\notin d(\ooo )$
implies $\a_1=2e$ so that $\a_1(M^*)=0$. It follows that
$R_2-R_1=2e$. Thus $d(-a_{1,2})>0=e-(R_2-R_1)/2$ and we have
$g(a_2/a_1)=\upo{(R_2-R_1)/2+e}=\upo{2e}$. Thus $\D\in g(a_2/a_1)$ and
we have $\[ a_1,a_2\]\ap\[\D a_1,\D a_2\]$ so $M\ap\[\D a_1,\D
a_2,a_3,\ldots,a_n\]$. But $d(\D )=2e=\a_1$ so we are done.

Suppose now that (b) happens for $M^*$. Then $|\oo/\p |=2$,
$d(-a_{2,3})=e-(R_3-R_2)/2=(R_2-R_1)/2+e$ and, if $n\geq 4$, then
$\a_3=\a_2(M^*)>e-(R_3-R_2)/2=(R_2-R_1)/2+e$. If we prove
that $\[ a_1,a_2,a_3\]$ is not isotropic then (c) happens (for $M$) so
we get a contradiction. But this is equivalent to
$(-a_{1,2},-a_{2,3})_\p =-1$, which follows, by Lemma 8.2(iii), from
$|\oo/\p |=2$ and
$d(-a_{1,2})+d(-a_{2,3})=e-(R_2-R_1)/2+(R_2-R_1)/2+e=2e$. 

Finally, if (c) happens for $M^*$ then $R_2=R_4$. Since also $R_1=R_3$,
$\[ a_1,a_2,a_3,a_4\]$ is in the situation from Lemma 8.3. Since
$(R_2-R_1)/2+e\geq\a_1(\[ a_1,a_2,a_3,a_4\] )\geq\a_1$ we have $\a_1(\[
a_1,a_2,a_3,a_4\] )=\a_1=(R_2-R_1)/2+e$. Take $\e\in\ooo$ arbitrary
with $d(\e )=\a_1=\a_1(\[ a_1,a_2,a_3,a_4\] )$. By Lemma 8.3 we have
$\[ a_1,a_2,a_3,a_4\]\ap\[ a'_1,a'_2,a'_3,a'_4\]$ with $a'_1=\e a_1$
relative to some good BONG. It follows that $M\ap\[
a'_1,a'_2,a'_3,a'_4,a_5,\ldots,a_n\]$ relative to some good BONG. \qed

By duality we get:

\bco There is $\e\in\ooo$ with $d(\e )=\a_{n-1}$ s.t. $M\ap\[
a'_1,\ldots,a'_n\]$ relative to some other good BONG with $a'_n=\e a_n$
unless $\a_{n-1}=(R_n-R_{n-1})/2+e$ and one of the following happens:

(a) $\a_{n-1}\notin d(\ooo )$.

(b) $|\oo/\p |=2$, $d[-a_{n-1,n}]=e-(R_n-R_{n-1})/2$ and
$\a_{n-2}>e-(R_n-R_{n-1})/2$ (if $n\geq 3$).

(c) $|\oo/\p |=2$, $R_{n-2}=R_n$, $\a_{n-3}>(R_n-R_{n-1})/2+e$ (if
$n\geq 4$) and $[a_{n-2},a_{n-1},a_n]$ is not isotropic.
\eco

\bco There is a change of good BONGs $\[ a_2,\ldots,a_n\]\ap\[
a'_2,\ldots,a'_n\]$ s.t. $\a_1=\a_1(\[ a_1,a'_2\] )$.
\eco
\pf We have $\a_1=\min\{\a_1(\[ a_1,a_2\] ),R_2-R_1+\a_1(\[
a_2,\ldots,a_n\] )\}$. If $\a_1=\a_1(\[ a_1,a_2\] )$ then we just take
$a'_i=a_i$ for all $i\geq 2$ and we are done. So we may assume that
$R_2-R_1+\a_1(\[ a_2,\ldots,a_n\] )=\a_1<\a_1(\[ a_1,a_2\] )=\min\{
(R_2-R_1)/2+e,R_2-R_1+d(-a_{1,2})\}$. Hence $\a_1(\[ a_2,\ldots,a_n\] )
<d(-a_{1,2})$ and $\a_1(\[ a_2,\ldots,a_n\] )<e-(R_2-R_1)/2\leq
(R_3-R_2)/2+e$. So by Lemma 8.8 we have $\[ a_2,\ldots,a_n\]\ap\[
a'_2,\ldots,a'_n\]$ relative to some good BONG s.t. $a'_2=\e a_2$ for
some $\e\in\ooo$ with $d(\e )=\a_1(\[ a_2,\ldots,a_n\] )$. Since also
$d(-a_{1,2})>\a_1(\[ a_2,\ldots,a_n\] )$ we have $d(-a_1a'_2)=d(-\e
a_{1,2})=\a_1(\[ a_2,\ldots,a_n\] )$. It follows that
$R_2-R_1+d(-a_1a'_2)=R_2-R_1+\a_1(\[ a_2,\ldots,a_n\] )=\a_1\leq
(R_2-R_1)/2+e$. Thus $\a_1=\min\{ (R_2-R_1)/2+e,R_2-R_1+d(-a_1a'_2)\}
=\a_1(\[ a_1,a'_2\] )$. \qed

\bco If $1\leq i\leq n-1$ then for a good choice of the good BONG of
$M$ we have $\a_i=\a_1(\[ a_i,a_{i+1}\] )$.
\eco
\pf For the case $i=1$ we use Corollary 8.10 We  have $\[
a_2,\ldots,a_n\]\ap\[ a'_2,\ldots,a'_n\]$ relative to a good BONG and
$\a_1=\a_1(\[ a_1,a'_2\] )$. But $\[ a_2,\ldots,a_n\]\ap\[
a'_2,\ldots,a'_n\]$ implies $M\ap\[ a_1,a'_2,\ldots,a'_n\]$. So if we
replace $a_1,\ldots,a_n$ by $a_1,a'_2,\ldots,a'_n$ we get
$\a_1=\a_1(\[ a_1,a_2\] )$. 

The case $i=n-1$ is done by duality. We have $M^\*\ap\[
a_n\1,\ldots,a_1\1\]$. Denote $\a_i^\*:=\a_i(M^\* )$. By the case
$i=1$ above, for a good choice of a good BONG for $M^\*$, and so for
$M$, we have $\a_1^\* =\a_1(\[ a_n\1,\a_{n-1}\1\] )=\a_1((\[
a_{n-1},a_n\] )^\* )$. But $\a_1^\* =\a_{n-1}$ and $\a_1((\[
a_{n-1},a_n\] )^\* )=\a_1(\[ a_{n-1},a_n\] )$ so $\a_{n-1}=\a_1(\[
a_{n-1},a_n\] )$. 

For the general case note that $\a_i=\min\{\a_i(\[
a_1,\ldots,a_{i+1}\]),\a_1(\[ a_i,\ldots,a_n\] )\}$. If $\a_i=\a_i(\[
a_1,\ldots,a_{i+1}\])$ then, by applying the case $i=n-1$, for a good
choice of the BONG of $\[ a_1,\ldots,a_{i+1}\]$ we have $\a_1(\[
a_i,a_{i+1}\] )=\a_i(\[ a_1,\ldots,a_{i+1}\])=\a_i$. If $\a_i=\a_1(\[
a_i,\ldots,a_n\] )$ then, by applying the case $i=1$, for a good
choice of the BONG of $\[ a_i,\ldots,a_n\]$ we have $\a_1(\[
a_i,a_{i+1}\] )=\a_1(\[ a_i,\ldots,a_n\])=\a_i$. \qed

\blm Suppose that $M,N$ have arbitrary ranks $m$ and $n$ and
$R_1=S_1$. We have:

(i) $A_1=\a_1$ and $A'_1=\a'_1$.

(ii) If $n\geq 2$ then $A_2=\min\{
(R_3-S_2)/2+e,R_3-S_2+d[-a_{1,3}b_1]\}$ and
$A'_2=R_3-S_2+d[-a_{1,3}b_1]$. If $n=1$ then
$S_2+A_2=R_3+d[-a_{1,3}b_1]$.
\elm
\pf (i) Since $R_1=S_1$ we have $A_1=\min\{
(R_2-R_1)/2+e,R_2-R_1+d[-a_{1,2}]\} =\a_1$ and
$A'_1=R_2-R_1+d[-a_{1,2}]=\a'_1$.

(ii) We have to prove that the term containing $d[a_{1,4}]$ in the
definition of $A_2,A'_2$ or $S_2+A_2$ can be removed. We have
$S_1=R_1\leq R_3$ so by Lemma 2.7(i) we can replace $d[a_{1,4}]$ by
$d[-a_{1,2}]$. Now $d[-a_{1,2}]\geq R_1-R_2+\a_1\geq S_1-R_4+\a_3\geq
S_1-R_4+d[-a_{1,3}b_1]$.

Thus if $n\geq 2$ then $R_3+R_4-S_1-S_2+d[-a_{1,2}]\geq
R_3-S_2+d[-a_{1,3}b_1]$ and so it can be removed from $A_2$ and
$A'_2$. If $n=1$ then $R_3+R_4-S_1+d[-a_{1,2}]\geq R_3+d[-a_{1,3}b_1]$
so $S_2+A_2=\min\{ R_3+d[-a_{1,3}b_1], R_3+R_4-S_1+d[-a_{1,2}]\}
=R_3+d[-a_{1,3}b_1]$. \qed

\bff We want now to write explicitely what $\[b_1\]\leq M$ means when
$S_1:=\ord b_1=R_1$. To do this we use Lemma 8.12 with $N=\[
b_1\]$. Condition (i) of the main theorem follows from $R_1=S_1$.

For (ii), since $\a_1\geq d[a_1b_1]$, the condition $d[a_1b_1]\geq
A_1=\a_1$ is equivalent to $d[a_1b_1]=\a_1$. Since
$d[a_1,b_1]=\min\{ d(a_1b_1),\a_1\}$ it is also equivalent to
$d(a_1b_1)\geq\a_1$.

The condition $R_3>S_1$ and $A_1+(S_2+A_2)>2e+R_2$ from (iii) means
$R_1<R_3$ and $\a_1+R_3+d[-a_{1,3}b_1]>2e+R_2$,
i.e. $\a_1+d[-a_{1,3}b_1]>2e+R_2-R_3$. So (iii) means $b_1\rep
[a_1,a_2]$ whenever $R_1<R_3$ and
$\a_1+d[-a_{1,3}b_1]>2e+R_2-R_3$. (Note that even if $n>3$ we may
assume that $S_2+A_2=S_2+A'_2=R_3+d[-a_{1,3}b_1]$. If $A_2\neq A'_2$
then, by Lemma 2.14 we have $A_1+A_2'>A_1+A_2>2e+R_2-S_2$. So, in both
cases, $A_1+(S_2+A_2)>2e+R_2$ iff $A_1+(S_2+A'_2)>2e+R_2$.) 

Since $S_1=R_1\leq R_3$ the condition $R_3+2e\leq S_1+2e<R_4$ from
(iv) means $R_1=R_3$ and $R_4-R_3>2e$. So (iv) means $b_1\rep
[a_1,a_2,a_3]$ whenever $R_1=R_3$ and $R_4-R_3>2e$.

Also if $n=3$ then the condition $b_1\rep FM$ needs only be verified
if $R_1=R_3$. If $R_1<R_3$ then either
$\a_1+d[-a_{1,3}b_1]>2e+R_2-R_3$ so from (iii) we have $b_1\rep
[a_1,a_2]\rep FM$ or $\a_1+d[-a_{1,3}b_1]\leq 2e+R_2-R_3$ so
$d(-a_{1,3}b_1)=d[-a_{1,3}b_1]\neq\j$ so $-a_{1,3}b_1\notin\fs$, which
implies that $b_1\rep [a_1,a_2,a_3]\ap FM$.

In conclusion, $\[ b_1\]\leq M$ iff the following hold:

a. $d[a_1b_1]=\a_1$, i.e. $d(a_1b_1)\geq\a_1$. 

b. $b_1\rep [a_1,a_2]$ if $R_1<R_3$ and
$\a_1+d[-a_{1,3}b_1]>2e+R_2-R_3$ or if $n=2$. 

c. $b_1\rep [a_1,a_2,a_3]$ if $R_1=R_3$ and either $R_4-R_3>2e$ or
$n=3$. 
\eff 

\blm Suppose that $S_1:=\ord b_1=R_1$ and $\[ b_1\]\leq M$. Then
$M\ap\[ a'_1,\ldots,a'_n\]$ relative to some good BONG with $a'_1=b_1$
unless one of the following happens:

(a) $R_1=R_3$, $\a_2+d[-a_{1,3}b_1]>2e$ and $[a_1,a_2,a_3]$ is not
isotropic.

(b) $R_1=R_3$, $|\oo/\p |=2$, $\a_1<(R_2-R_1)/2+e$,
$\a_2+d[-a_{1,3}b_1]=2e$, $[a_1,a_2,a_3]$ is isotropic and, if $n\geq
4$, $\a_2+\a_3>2e$. 

(c) $R_1=R_3$, $|\oo/\p |=2$, $R_2<R_4$,
$d[-a_{1,3}b_1]=\a_3=(R_4-R_3)/2+e$,
$d[a_{1,4}]=\a_2=e-(R_4-R_3)/2$, $[a_1,a_2,a_3,a_4]\top [b_1]$ is
not isotropic and, if $n\geq 5$, $\a_4>e-(R_4-R_3)/2$. \elm 

\pf First we prove that conditions (a)-(c) are independent of the
choice of BONGs. Now since $d[-a_{1,3}b_1]$ and $d[a_{1,4}]$ are
independent of BONGs all conditions from (a)-(c) are independent
on BONGs with the exception of the conditions involving
$[a_1,a_2,a_3]$ (for (a) and (b)) or $[a_1,a_2,a_3,a_4]\top [b_1]$
(for (c)). So we still have to prove that if $M\ap\[
a'_1,\ldots,a'_n\]$ relative to some other good BONG and all the
conditions not involving $[a_1,a_2,a_3]$ or $[a_1,a_2,a_3,a_4]\top
[b_1]$ form (a), (b) or (c) are satisfied, then $[a_1,a_2,a_3]$ is
isotropic iff $[a'_1,a'_2,a'_3]$ is isotropic (in the case of (a) and
(b)) resp. $[a_1,a_2,a_3,a_4]\top [b_1]$ is isotropic iff
$[a'_1,a'_2,a'_3,a'_4]\top [b_1]$ is isotropic (in the case of (c)).

For consequences of $R_1=R_3$ see 8.7.

First we consider the case of (a) and (b). If $n=3$ then
$[a_1,a_2,a_3]\ap [a'_1,a'_2,a'_3]\ap FM$ so we may assume that
$n\geq 4$. In the case of (b) we have $\a_2+\a_3>2e$ and same happens
in the case of (a), when $\a_2+\a_3\geq\a_2+d[-a_{1,3}b_1]>2e$. By
[B3, Theorem 3.1(iv)] we have $[a'_1,a'_2]\rep [a_1,a_2,a_3]$ so
$[a_1,a_2,a_3]\ap [a'_1,a'_2,a_{1,3}a'_{1,2}]$. Thus we have to
prove that $[a'_1,a'_2,a_{1,3}a'_{1,2}]$ is isotropic iff
$[a'_1,a'_2,a'_3]$ is. But this is equivalent to
$(-a'_{1,2},-a_{1,3}a'_1)_\p =(-a'_{1,2},-a'_{2,3})_\p$, i.e. to
$(-a'_{1,2},a_{1,3}a'_{1,3})_\p =1$. Now $R_1=R_3$ so
$d(-a'_{1,2})\geq R_1-R_2+\a_1=\a_2$. Hence
$d(-a'_{1,2})+d(a_{1,3}a'_{1,3})\geq\a_2+\a_3>2e$, which implies
$(-a'_{1,2},a_{1,3}a'_{1,3})_\p =1$.

Take now the case of (c). Let $V\ap [a_1,a_2,a_3,a_4]\top [b_1]$
and $V'\ap [a'_1,a'_2,a'_3,a'_4]\top [b_1]$. We have to prove that
$V,V'$ are either both isotropic or both anisotropic, i.e. that
$SV=SV'$. If $n=4$ then $[a_1,a_2,a_3,a_4]\ap [a'_1,a'_2,a'_3,a'_4]\ap
FM$ so we may assume that $n\geq 5$. We have $\a_3+\a_4>2e$ so
$[a'_1,a'_2,a'_3]\rep [a_1,a_2,a_3,a_4]$. It follows that $V\ap
[a'_1,a'_2,a'_3,a_{1,4}a'_{1,3}]\top [b_1]$. Together with $V'\ap
[a'_1,a'_2,a'_3,a'_4]\top [b_1]$, this implies that $V\pp
[a'_4]\ap V'\pp [a_{1,4}a'_{1,3}]$. Considering Hasse symbols we
get $SV (a'_4,a'_4\det V)_\p
=SV'(a_{1,4}a'_{1,3},a_{1,4}a'_{1,3}\det V')$. (By [OM, 58:3] we
have $S(V\pp [\a ])=SV(\det V,\a )_\p S([\a ])=SV(\det V,\a )_\p
(\a,\a )_\p=SV(\a,\a\det V)_\p$.) But $\det V=a_{1,4}b_1$ and
$\det V'=a'_{1,4}b_1$ so $SV(a'_4,a_{1,4}a'_4b_1)_\p
=SV'(a_{1,4}a'_{1,3},a_{1,4}a'_4b_1)_\p$. It follows that $SV=SV'$ is
equivalent to $(a'_4,a_{1,4}a'_4b_1)_\p
=(a_{1,4}a'_{1,3},a_{1,4}a'_4b_1)_\p$, i.e. to
$1=(a_{1,4}a'_{1,4},a_{1,4}a'_4b_1)_\p
=(a_{1,4}a'_{1,4},-a'_{1,3}b_1)_\p$. But $d(a_{1,4}a'_{1,4})\geq\a_4$
and $d(-a'_{1,3}b_1)\geq d[-a'_{1,3}b_1]=d[-a_{1,3}b_1]=\a_3$ so
$d(a_{1,4}a'_{1,4})+d(-a'_{1,3}b_1)\geq\a_3+\a_4>2e$, which implies
$(a_{1,4}a'_{1,4},-a'_{1,3}b_1)_\p =1$.

We now prove that if (a), (b) or (c) from the hypothesis is satisfied
then we cannot have $M\ap\[ a'_1,\ldots,a'_n\]$ relative to some good
BONG s.t. $a'_1=b_1$. Since (a) - (c) are independent of the choice
of BONGs we may assume that $a_i=a'_i$. In particular, $a_1=b_1$,
which implies that $d(-a_{2,3})=d(-a_{1,3}b_1)$.

If (a) holds then $d(-a_{1,2})\geq R_1-R_2+\a_1=\a_2$ and
$d(-a_{2,3})=d(-a_{1,3}b_1)\geq d[-a_{1,3}b_1]$. It follows that
$d(-a_{1,2})+d(-a_{2,3})\geq\a_2+d[-a_{1,3}b_1]>2e$, which
implies that $(-a_{1,2},-a_{2,3})_\p =1$ so $[a_1,a_2,a_3]$ is
isometric. Contradiction. 

If (b) holds then $\a_2+d[-a_{1,3}b_1]=2e<\a_2+\a_3$ so
$d[-a_{1,3}b_1]<\a_3$. It follows that  $d[-a_{1,3}b_1]=\min\{
d(-a_{1,3}b_1),\a_3\} =d(-a_{1,3}b_1)=d(-a_{2,3})$. (Same happens if
$n=3$, when $d[-a_{1,3}b_1]=d(-a_{1,3}b_1)$ by definition.) Now
$\a_1=\min\{ (R_2-R_1)/2+e, R_2-R_1+d(-a_{1,2}), R_3-R_1+d(-a_{2,3}),
R_3-R_1+\a_3\} =\min\{ (R_2-R_1)/2+e, R_2-R_1+d(-a_{1,2}),
d(-a_{2,3}), \a_3\}$. ($\a_3$ is ignored if $n=3$.) But
$\a_1<(R_2-R_1)/2+e$. Since $R_1=R_3$ this implies
$\a_1+\a_2<2e=\a_2+d[-a_{1,3}b_1]$. It follows that
$\a_1<d[-a_{1,3}b_1]$ so we have both
$\a_1<d(-a_{1,3}b_1)=d(-a_{2,3})$ and $\a_1<\a_3$ (if $n\geq 4$). Thus
$\a_1=R_2-R_1+d(-a_{1,2})$, i.e. $d(-a_{1,2})=R_1-R_2+\a_1=\a_2$. It
follows that $d(-a_{1,2})+d(-a_{2,3})=\a_2+d[-a_{1,3}b_1]=2e$. Since
$|\oo/\p |=2$ we get by Lemma 8.2(iii) $(-a_{1,2},-a_{2,3})_\p =-1$ so
$[a_1,a_2,a_3]$ is not isometric. Contradiction.

Suppose now that (c) holds. Since $a_1=b_1$ we have that
$[a_2,a_3,a_4]\ap [a_1,a_2,a_3,a_4]\top [b_1]$ is not
isotropic. Therefore $(-a_{2,3},-a_{3,4})_\p =-1$, which implies that
$d(-a_{2,3})+d(-a_{3,4})\leq 2e$. Since $R_1=R_3$ and
$\a_2=e-(R_4-R_3)/2<(R_3-R_2)/2+e$ we have by 8.7 that
$\a_1<(R_2-R_1)/2+e$ and also $\a_1+\a_2<2e=\a_2+\a_3$,
i.e. $\a_1<\a_3$. We also have $d(-a_{2,3})=d(-a_{1,3}b_1)\geq
d[-a_{1,3}b_1]=\a_3>\a_1$. Now, same as for case (b), we have
$\a_1=\min\{ (R_2-R_1)/2+e, R_2-R_1+d(-a_{1,2}),
d(-a_{2,3}),\a_3\}$. But $(R_2-R_1)/2+e,d(-a_{2,3}),\a_3>\a_1$ so
$\a_1=R_2-R_1+d(-a_{1,2})$, i.e. $d(-a_{1,2})=R_1-R_2+\a_1=\a_2$. If
$n\geq 5$ then $\a_4>e-(R_4-R_3)/2=\a_2=d[a_{1,4}]$ so
$d(a_{1,4})=d[a_{1,4}]=\a_2$. Same happens if $n=4$, when
$d[a_{1,4}]=d(a_{1,4})$ by definition. Since also $d(-a_{1,2})=\a_2$
and $|\oo/\p |=2$ we have $d(-a_{3,4})>\a_2$ by Lemma 8.1(ii). Since
also $d(-a_{2,3})\geq\a_3$ we get
$d(-a_{2,3})+d(-a_{3,4})>\a_2+\a_3=2e$. Contradiction.

We prove now that, conversely, if $\[ b_1\]\leq M$ and none of (a)-(c)
holds then $M\ap\[ a'_1,\ldots,a'_n\]$ relative to some good BONG
s.t. $a'_1=b_1$. For conditions equivalent to $\[ b_1\]\leq M$ we
refer to 8.13.

We use the invariance of (a)-(c) under change of BONGs to make some
reductions. 

(I) By Corollary 8.10 we may assume that $\a_1=\a_1(\[ a_1,a_2\] )$.

(II) By changing, if necessary, the BONG of $\[ a_3,\ldots,a_n\]$ we
may assume by Corollary 8.10 that $\a_1(\[ a_2,a_3\] )=\a_1(\[
a_2,\ldots,a_n\] )$. Note that $a_1,a_2$ are not changed by this
reduction step so the previous reduction step is preserved. 

Note that the reduction (I) implies that $\a_1(\[ a_1,\ldots,a_i\]
)=\a_1$ for any $i\geq 2$. (We have $\a_1(\[ a_1,a_2\] )\geq\a_1(
a_1,\ldots,a_i\] )\geq\a_1$.) Similarly the reduction (II) implies
$\a_1(\[ a_2,\ldots,a_i\] )=\a_1(\[ a_2,\ldots,a_n\] )$ for any $i\geq
3$ and so $\a_2(\[ a_1,\ldots,a_i\])=\a_2$. (We have $\a_2(\[
a_1,\ldots,a_i\] )=\min\{ R_3-R_1+d(-a_{1,2}),\a_1(\[ a_2,\ldots,a_i\]
)\}$ and $\a_2=\min\{ R_3-R_1+d(-a_{1,2}),\a_1(\[ a_2,\ldots,a_n\]
)\}$.) 

Let $b_1=\e a_1$ so $d(\e )\geq\a_1$. Suppose first that $\e
a_1=b_1\rep [a_1,a_2]$ so $\e\in\N (-a_{1,2})$. Since $d(\e
)\geq\a_1=\a_1(\[ a_1,a_2\])$ we have $\e\in\upon{\a_1(\[
a_1,a_2\])}{-a_{1,2}}=g(a_2/a_1)$. (See [B3, 3.8].) Thus $\[
a_1,a_2\]\ap\[\e a_1,\e a_2\]$ so $M\ap\[\e a_1,\e
a_2,a_3,\ldots,a_n\]$. Since $b_1=\e a_1$ we are done.

So we may assume that $b_1\notrep [a_1,a_2]$. By 8.13 it follows
that $n\geq 3$ and either $R_1=R_3$ or $R_1<R_3$ and
$\a_1+d[-a_{1,3}b_1]\leq 2e+R_2-R_3$. We consider the two cases:

1. $R_1<R_3$ and $\a_1+d[-a_{1,3}b_1]\leq 2e+R_2-R_3$. We have
$d[a_1b_1]=\a_1\geq R_2-R_3+\a_2$ and $d[-a_{2,3}]\geq R_2-R_3+\a_2$
so $d[-a_{1,3}b_1]\geq R_2-R_3+\a_2$. It follows that
$\a_1+R_2-R_3+\a_2\leq\a_1+d[-a_{1,3}b_1]\leq 2e-R_2-R_3$ so
$\a_1+\a_2\leq 2e$.

Suppose first that $n=3$. Then $d[-a_{1,3}b_1]=d(-a_{1,3}b_1)$. If
$\a_2=0$ then $R_3-R_2=-2e$ so $R_2-R_1>R_2-R_3=2e$. Thus
$d(a_1b_1)\geq\a_1>2e$, i.e. $b_1\in a_1\ooo^2$, which implies $M\ap\[
b_1,a_2,\ldots,a_n\]$. So we may assume that $\a_2>0$. Since
$\a_2+d(-a_{1,3}b_1)\leq 2e$ we get $d(-a_{1,3}b_1)<2e$. It follows
that there is $\eta\in\ooo$ with $d(\eta )=2e-d(-a_{1,3}b_1)$
s.t. $(\eta,-a_{1,3}b_1)_\p =-1$. We claim that $M\ap\[\e a_1,\e\eta
a_2,\eta a_3\]$ relative to some good BONG. To do this we use [B3,
Theorem 3.1] as well as Lemma 8.6. Let $\( M\ap\[ a'_1,a'_2,a'_3\]$
with $a'_1=b_1=\e a_1$, $a'_2=\e\eta a_2$ and $a'_3=\eta a_3$. We have
$R_i(\( M)=\ord a'_i=R_i$ so (i) of the classifiction theorem
holds. We have $d(a_1a'_1)=d(\e )\geq\a_1$ and
$d(a_{1,2}a'_{1,2})=d(\eta )=2e-d(-a_{1,3}b_1)\geq\a_2$ so (iii) also
holds. Together with $a_{1,3}a'_{1,3}\in\fs$, this implies by Lemma
8.6(i) that the lattice $\( M\ap\[ a',a'_2,a'_3\]$ exists. Also since
$R_1<R_3$, i.e. $M$ has property A, we have by Lemma 8.6(iii) that
$\a_i(\( M)=\a_i$ for $i=1,2$ so (ii) also holds. Finally,
$\a_1+\a_2\leq 2e$ so (iv) is vacuous. We are left to prove that $F\(
M\ap FM$. Now $[a_1,a_2,a_3]$, $[\e a_1,\e a_2,a_3]$ and $[\e
  a_1,\e\eta a_2,\eta a_3]$ have the same determinant $a_{1,3}$. Since
$\e a_1=b_1\notrep [a_1,a_2]$ we have $[a_1,a_2]\not\ap [\e a_1,\e
  a_2]$ so $[a_1,a_2,a_3]\not\ap [\e a_1,\e a_2,a_3]$. Also $b_1=\e
a_1$ so $\e a_{2,3}=a_{1,3}b_1$ in $\ff/\fs$, which implies $(\eta,-\e
a_{2,3})_\p =-1$. Thus $[\e a_2,a_3]\not\ap [\e\eta a_2,\eta a_3]$ so
$[\e a_1,\e a_2,a_3]\not\ap [\e a_1,\e\eta a_2,\eta a_3]$. It follows
that $[a_1,a_2,a_3]\ap [\e a_1,\e\eta a_2,\eta a_3]$, i.e. $FM\ap F\(
M$. 

Suppose now that $n\geq 4$. We make the reduction step (I). Let $\(
M\ap\[ a_1,a_2,a_3\]$. Let $\(\a_i:=\a_i(\( M)$ for $i=1,2$. The
reduction step (I) implies $\(\a_1=\a_1$. 

We assume first that $\a_1+d(-a_{1,3}b_1)\leq 2e+R_2-R_3$. We prove
that $\[ b_1\]\leq \(M$. This follows from 8.13. Indeed,
$d(a_1b_1)\geq\a_1=\(\a_1$ so a. holds. Also b. follows from
$(\a_1+d[-a_{1,3}b_1])(b_1,\(
M)=\(\a_1+d(-a_{1,3}b_1)=\a_1+d(-a_{1,3}b_1)\leq
2e+R_2-R_3$. (Here $(\a_1+d[-a_{1,3}b_1])(b_1,\( M)$ is the
$\a_1+d[-a_{1,3}b_1]$ corresponding to $b_1$ and $\( M$.) Finally,
c. is vacuous since $R_1<R_3$. Since $\[ b_1\]\leq\( M$, by
applying the case $n=3$ to $b_1$ and $\( M$, we get $\[
a_1,a_2,a_3\]\ap\[ a'_1,a'_2,a'_3\]$ relative to a good BONG, with
$a'_1=b_1$, which implies that $M\ap\[
a'_1,a'_2,a'_3,a_4,\ldots,a_n\]$ so we are done.

Therefore we may assume that $\a_1+d(-a_{1,3}b_1)>2e+R_2-R_3\geq
\a_1+d[-a_{1,3}b_1]$ whenever the reduction step (I) holds. Since
$d[-a_{1,3}b_1]=\min\{ d(-a_{1,3}b_1),\a_3\}$ we have
$d[-a_{1,3}b_1]=\a_3<d(-a_{1,3}b_1)$. In particular,
$\a_1+\a_3=\a_1+d[-a_{1,3}b_1]\leq 2e+R_2-R_3$.

Suppose that $R_1+\a_1=R_2+\a_2$. We make the reductions steps (I) and
(II). From (II) we have $\a_2=\(\a_2=\min\{ (R_3-R_2)/2+e,
R_3-R_2+d(-a_{2,3}), R_3-R_1+d(-a_{1,2})\}$. We have
$R_3-R_1+d(-a_{1,2})\geq R_3-R_2+\a_1=R_3-R_1+\a_2>\a_2$. Also
$\a_2=R_1-R_2+\a_1\leq
R_1-R_2+(R_2-R_1)/2+e=(R_1-R_2)/2+e<(R_3-R_2)/2+e$. Thus
$\a_2=R_3-R_2+d(-a_{2,3})$, i.e. $d(-a_{2,3})=R_2-R_3+\a_2$. Since
also $d(a_1b_1)\geq\a_1=R_2-R_1+\a_2>R_2-R_3+\a_2=d(-a_{2,3})$ we get
$d(-a_{1,3}b_1)=R_2-R_3+\a_2\leq\a_3$. Contradiction.

Hence $R_1+\a_1<R_2+\a_2$ so $\a_1<R_2-R_1+\a_2$. It follows that,
regardless of the choice of the BONG, we have $\a_1=\min\{\a_1(\[
a_1,a_2\] ),R_2-R_1+\a_2\} =\a_1(\[ a_1,a_2\] )$. Thus the reduction
step (I) and its consequences, $\a_1=\a_1(\[ a_1,\ldots,a_i\] )$ for
any $i\geq 3$ and $\a_1+d(-a_{1,3}b_1)>2e+R_2-R_3$, hold regardless of
the choice of BONG. 

By Corollary 8.11, after a change of BONGs for $M$, we may assume that
$\a_3=\a_1(\[ a_3,a_4\] )$. (Note that this change of BONGs preserves
the reduction step (I) (see above) but not necessarily the step (II),
which is not needed at this time.) Since $\a_3\leq\a_1 (\[
a_3,\ldots,a_n\] )\leq\a_1(\[ a_3,a_4\] )$ we have $\a_3=\a_1 (\[
a_3,\ldots,a_n\] )$. 

If $\a_1 (\[ a_3,\ldots,a_n\] )=\a_3<(R_4-R_3)/2+e$ then by Lemma 8.8
applied to $\[ a_3,\ldots,a_n\]$ we get that there is $\eta\in\ooo$
with $d(\eta )=\a_1 (\[ a_3,\ldots,a_n\] )=\a_3$ s.t. $\[
a_3,\ldots,a_n\]\ap\[ a'_3,\ldots,a'_n\]$ relative to some good BONG
and $a'_3=\eta a_3$. Hence $M\ap\[ a_1,a_2,a'_3,\ldots,a'_n\]$. By
this change of BONGs $d(-a_{1,3}b_1)$ is replaced by
$d(-a_{1,2}a'_3b_1)=d(-\eta a_{1,3}b_1)$. Since
$d(-a_{1,3}b_1)>\a_3=d(\eta )$ we have $d(-\eta a_{1,3}b_1)=\a_3$. So
after this change of BONGs the inequality $d(-a_{1,3}b_1)>\a_3$
doesn't hold anymore and we are done. 

So we may assume that $\a_3=(R_4-R_3)/2+e$. Note that, regardless of
the BONG, we have $(R_4-R_3)/2+e\geq\a_1(\[ a_3,a_4\]
)\geq\a_3$. Hence $\a_1(\[ a_3,a_4\] )=\a_3=(R_4-R_3)/2+e$ holds
regardless of the BONG. So, at this time, the BONG of $M$ is
considered to be arbitrary. Now $-R_2+\a_1\geq -R_4+\a_3$ so
$2e+R_2-R_3\geq\a_1+\a_3\geq R_2-R_4+2\a_3=2e+R_2-R_3$. Thus we must
have equality so $-R_2+\a_1=-R_3+\a_2=-R_4+\a_3=-(R_3+R_4)/2+e$. Hence
$\a_2=e-(R_4-R_3)/2$.

Now $d(-a_{3,4})\geq R_3-R_4+\a_3=e-(R_4-R_3)/2$. If $\a_1(\[
a_3,a_4\] )=\a_3\in d(\ooo )$ and either $|\oo/\p |>2$ or
$d(-a_{3,4})>e-(R_4-R_3)/2$ then by Lemma 8.8 applied to $\[
a_3,a_4\]$ we have $\[ a_3,a_4\]\ap\[\eta a_3,\eta a_4\]$ for some
$\eta\in\ooo$ with $d(\eta )=\a_1(\[ a_3,a_4\] )=\a_3$. Thus $M\ap\[
a_1,a_2,\eta a_3,\eta a_4,a_5,\ldots,a_n\]$ relative to some good
BONG. After this change of BONGs $d(-a_{1,3}b_1)$ is replaced by
$d(-\eta a_{1,3}b_1)$. But $d(-a_{1,3}b_1)>\a_3=d(\eta )$ so $d(-\eta
a_{1,3}b_1)=\a_3$. So after this change of BONGs we no longer have
$d(-a_{1,3}b_1)>\a_3$ and we are done.

So we may assume that either $a_3\notin d(\ooo )$ or $|\oo/\p |=2$ and
$d(-a_{3,4})=e-(R_4-R_3)/2$. (I.e. that (a) or (b) of Lemma 8.8 holds
for $\[ a_3,a_4\]$.) Now
$-R_2+\a_1=-R_4+\a_3=e-(R_3+R_4)/2<e-(R_1+R_2)/2$. (We have $R_1<R_3$
and $R_2\leq R_4$.) Thus $\a_1<(R_2-R_1)/2+e$. It follows that $\a_1$
is odd and $<2e$. Since $-R_2+\a_1=-R_4+\a_3$ we have by Lemma 7.3(ii)
that $\a_1\ev\a_2\ev\a_3\m2$ so $\a_2,\a_3$ are odd as well and also
$\a_2,\a_3\leq 2e$. Thus $\a_1,\a_2,\a_3\in d(\ooo )$. Since $a_3\in
d(\ooo )$ we must have $|\oo/\p |=2$ and
$d(-a_{3,4})=e-(R_4-R_3)/2=\a_2$. Also this holds regarless of the
BONG since at this time we don't use any reduction. 

We now use the reduction step (II). This implies that $\a_2(\[
a_1,a_2,a_3\] )=\a_2$. If $\[ a_1,a_2,a_3\]$ satisfies the
conditions of Corollary 8.9 then there is $\eta\in\ooo$ with
$d(\eta )=\a_2(\[ a_1,a_2,a_3\] )=\a_2$ s.t. $\[
a_1,a_2,a_3\]\ap\[ a'_1,a'_2,a'_3\]$ relative to some other BONG 
with $a'_3=\eta a_3$. After this chage of BONGs $d(-a_{3,4})$ is
replaced by $d(-a'_3a_4)=d(-\eta a_{3,4})$. But $d(\eta
)=\a_2=d(-a_{3,4})$ and $|\oo/\p |=2$ so by Lemma 8.1 we have $d(-\eta
a_{3,4})>\a_2$. So after this change of BONGs we no longer have
$d(-a_{3,4})=\a_2$ and we are done. So we  may assume that
$\a_2=(R_3-R_2)/2+e$ and $\[ a_1,a_2,a_3\]$ is in one of the cases
(a)-(c) of Corollary 8.9. Since $e-(R_4-R_3)/2=\a_2=(R_3-R_2)/2+e$ we
have $R_2=R_4$. Now $\a_2\in d(\ooo )$ so (a) doesn't hold. We have
$-R_2+\a_1=-R_4+\a_3=-R_2+\a_3$ so
$\a_1=\a_3=(R_4-R_3)/2+e=e-(R_3-R_2)/2$ and so (b) doesn't hold.
Finally, $R_1<R_3$ so (c) doesn't hold and we are done.

2. $R_1=R_3$. See 8.7 for consequences. In particular, since $R_1\ev
R_2\ev R_3\m2$ we have that $\ord a_{1,2},\ord a_{2,3}$ are even. Also
$\ord a_{1,3}b_1=2R_1+R_2+R_3$ is even. Thus
$a_{1,2},a_{2,3},a_{1,3}b_1\in\os$. 

Take first the case $n=3$. Suppose that neither (a) nor (b)
holds. ((c) is vacuous when $n=3$.) We are looking for an
$\eta\in\ooo$ s.t. $M\ap\[\e a_1,\e\eta a_2,\eta a_3\]$. Let $\(
M\ap\[ a'_1,a'_2,a'_3\]$ with $a'_1=\e a_1$, $a'_2=\e\eta a_2$ and
$a'_3=\eta a_3$. We claim that a lattice $\( M\ap\[ a'_1,a'_2,a'_3\]$
exists and $M\ap\( M$ iff the following happen:

1) $d(\eta )\geq\a_2$.

2) $\a_1=(R_2-R_1)/2+e$, $d(-\e a_{2,3})=\a_1$ or $d(-\eta
a_{1,2})=\a_2$.

3) $(-a_{1,2},-a_{2,3})_\p =(-\eta a_{1,2},-\e a_{2,3})_\p$.

Indeed $\ord a'_i=\ord a_i=R_i$ so the condition (i) of the
classification theorem in [B3] holds. For (iii) note that
$d(a_1a'_1)=d(\e )\geq\a_1$, by hypothesis, and that
$d(a_{1,2}a'_{1,2})\geq\a_2$ is equivalent to $d(\eta )\geq\a_2$
i.e. to 1). Since also $a_{1,3}a'_{1,3}\in\fs$, by Lemma 8.6
these imply that $\( M$ exists and also $\a_i\leq\(\a_i:=\a_i(\(
M)$ for $i=1,2$. We have $\det FM=\det F\( M$ so $FM\ap F\( M$ is
equivalent to $(-a_{1,2},-a_{2,3})_\p =(-\eta a_{1,2},-\e
a_{2,3})_\p$, i.e. to 3). So assuming that 1) and 3) hold, by Lemma
7.11 we have $M\ap\( M$ iff $\a_1=\(\a_1:=\a_1(\( M)$. But this
equivalent to 2). Indeed, we have $\a_1\leq\(\a_1=\min\{
(R_2-R_1)/2+e, R_2-R_1+d(-a'_{1,2}), R_3-R_1+d(-a'_{2,3})\} =\min\{
(R_2-R_1)/2+e,R_2-R_1+d(-\e a_{2,3}), d(-\eta a_{1,2})\}$ so
$\a_1=\(\a_1$ iff $\a_1=(R_2-R_1)/2+e$, $d(-\e a_{2,3})=\a_1$ or
$\a_1=R_2-R_1+d(-\eta a_{1,2})$, i.e. $d(-\eta
a_{1,2})=R_1-R_2+\a_1=\a_2$.

Note that $d(-\e a_{2,3})=d(-a_{1,3}b_1)=d[-a_{1,3}b_1]$ so $(-\eta
a_{1,2},-\e a_{2,3})_\p =(-\eta a_{1,2},-a_{1,3}b_1)_\p$. Also note
that $d(-a_{2,3})\geq\a_1$ and $d(a_1b_1)\geq\a_1$ so
$d(-a_{1,3}b_1)\geq\a_1$. 

Since $d(-a_{1,2})\geq\a_2$ the condition 1) is equivalent to
$d(-\eta a_{1,2})\geq\a_2$.

Suppose first that $[a_1,a_2,a_3]$ is not isotropic so
$(-a_{1,2},-a_{2,3})=-1$ and the condition 3) means $(-\eta
a_{1,2},-a_{1,3}b_1)_\p =-1$. We have $d(-a_{1,2})+d(-a_{2,3})\leq
2e$. Since (a) doesn't hold we also have $\a_2+d(-a_{1,3}b_1)\leq
2e$. In particular, $\a_2\leq 2e$. If $\a_2$ is odd then $\a_2\in
d(\ooo )$. Since $\a_2+d(-a_{1,3}b_1)\leq 2e$ there is, by Lemma
8.2(i), some $\eta\in\ff$ s.t. $d(-\eta a_{1,2})=\a_2$ and $(-\eta
a_{1,2},-a_{1,3}b_1)_\p =-1$. But $a_{1,2}\in\os$ and $d(-\eta
a_{1,2})=\a_2\in d(\ooo )$ so we may assume $\eta\in\ooo$. For this
choice of $\eta$ all conditions 1) - 3) are satisfied. If $\a_2$ is
even then $\a_2=(R_3-R_2)/2+e$ so $\a_1=(R_2-R_1)/2+e$ and
$\a_1+\a_2=2e$. Now $d(-a_{1,2})\geq\a_2$ and $d(-a_{2,3})\geq\a_1$ so
$2e\geq d(-a_{1,2})+d(-a_{2,3})\geq\a_1+\a_2=2e$. It follows that
$d(-a_{1,2})=\a_2$ and $d(-a_{2,3})=\a_1=2e-\a_2$ are even. But
$-a_{1,2},-a_{2,3}\in\os$ so we must have
$d(-a_{1,2})=d(-a_{2,3})=2e$. But this contradicts
$d(-a_{1,2})+d(-a_{2,3})=2e$. 

Suppose now that $[a_1,a_2,a_3]$ is isotropic so
$(-a_{1,2},-a_{2,3})=1$ and the condition 3) means $(-\eta
a_{1,2},-a_{1,3}b_1)_\p =1$. We have to find an $\eta\in\ooo$ that
satisfies conditions 1) - 3) above. If $\a_1=(R_2-R_1)/2+e$ or
$d(-a_{1,3}b_1)$ then the condition 2) is satisfied and conditions 1)
and 3) are satisfied if we choose $\eta\in\ooo$ s.t. $-\eta
a_{1,2}\in\fs$. Thus we may assume that $\a_1<(R_2-R_1)/2+e$ and
$\a_1<d(-a_{1,3}b_1)$. This implies that also $\a_2<(R_3-R_2)/2+e$ so
$\a_2$ is odd and $<2e$ and so $\a_2\in d(\ooo )$ and also $\a_2\neq
0,2e$. If $|\oo/\p |>2$ or $\a_2+d(-a_{1,3}b_1)\neq 2e$ there is by
Lemma 8.2(ii) and (iii) some $\eta\in\ff$ s.t. $-\eta a_{1,2}$
satisfies both $d(-\eta a_{1,2})=\a_2$ and $(-\eta
a_{1,2},-a_{1,3}b_1)_\p =1$. Since $a_{1,2}\in\os$ and $d(-\eta
a_{1,2})=\a_2\in d(\ooo )$ we may assume that $\eta\in\ooo$. So $\eta$
satisfies all required conditions. Hence we can assume that $|\oo/\p
|=2$ and $\a_2+d(-a_{1,3}b_1)=2e$. Also $\a_1<(R_2-R_1)/2+e$ so (b)
holds. Contradiction. 

Suppose now that $n=4$. We will prove that we can change the BONG
of $M$ s.t. $\[ b_1\]\leq \( M:=\[ a_1,a_2,a_3\]$ and neither of the
conditions (a) and (b) holds for $\( M$. This would imply by the
ternary case that $\( M\ap\[ a'_1,a'_2,a'_3\]$ relative to some good
BONG with $a'_1=b_1$ and so $M\ap\[ a'_1,a'_2,a'_3,a_4\]$. (Recall
that (c) is vacuous when the dimension is 3.) Since $R_1=R_3$ the
condition $\[ b_1\]\rep\( M$ is equivalent by 8.13 to
$d(a_1b_1)\geq\(\a_1$, where $\(\a_i:=\a_i(\( M)$, and $b_1\rep F\(
M\ap [a_1,a_2,a_3]$. The conditions (a) and (b) for $\( M$ are the
same as for $M$ but with $\a_1,\a_2$ and $d[-a_{1,3}b_1]$ replaced by
$\(\a_1,\(\a_2$ and $d(-a_{1,3}b_1)$ and with the condition
$\a_2+\a_3>2e$ from (b) ignored. (Here $\(\a_i:=\a_i(\( M)$.)

Since the case $R_2=R_4$ is treated by Lemma 8.3 we may assume that
$R_2<R_4$. 

We will reduce to the case when $\(\a_i=\a_i$ for $i=1,2$. Note
that $R_1=R_3$ implies both $R_1+\a_1=R_2+\a_2$ and
$R_1+\(\a_1=R_2+\(\a_2$ so $\a_1=\(\a_1$ iff $\a_2=\(\a_2$. So it is
enough to have $\a_1=\(\a_1$. This can be obtained by using the
reduction step (I). If $\a_i=\(\a_i$ for $i=1,2$ then
$d(a_1b_1)\geq\a_1=\(\a_1$. In order to have $\[ b_1\]\leq\( M$ we
still need $b_1\rep[a_1,a_2,a_3]$. Suppose the contrary. It follows
that $[a_1,a_2,a_3]$ is not isotropic and $-a_{1,3}b_1\in\fs$ so
$\a_2+d(-a_{1,3})=\j >2e$. Hence (a) holds for $b_1$ and $\( M$. Thus
the condition that $\[ b_1\]\leq\( M$ is included in the condition
that (a) doesn't hold for $b_1$ and $\( M$. Therefore it needs not be
verified. 

We have three cases:

a. $\a_2+\a_3>2e$. We have $d[-a_{1,3}b_1]=\min\{
d(-a_{1,3}b_1),\a_3\}$. If the condition $\a_2+d(-a_{1,3}b_1)>2e$
resp. $\a_2+d(-a_{1,3}b_1)=2e$ from (a) resp. (b) corresponding to
$\( M$ holds then, since also $\a_2+\a_3>2e$, the condition
$\a_2+d[-a_{1,3}b_1]>2e$ resp. $\a_2+d[-a_{1,3}b_1]=2e$ from (a)
resp. (b) corresponding to $M$ will also hold. The other
conditions are the same for both $\( M$ and $M$ with the exception
of $\a_2+\a_3>2e$ from (b), which holds anyways. Thus if $\( M$
satisfies the condition (a) or (b) then so does $M$, contrary to
our assumption. So neither (a) nor (b) hold for $b_1$ and $\( M$
and we are done.

b. $\a_2+\a_3<2e$. We may assume that $\a_2+d(-a_{1,3}b_1)\geq 2e$
since otherwise neither (a) nor (b) holds for $b_1$ and $\( M$. It
follows that $d(-a_{1,3}b_1)>\a_3$. Now $2\a_3\leq
R_4-R_3+\a_2+\a_3<R_4-R_3+2e$ so $\a_3<(R_4-R_3)/2+e$.

Suppose first that $\a_1\geq\a_3$. From the reduction step (I) we have
$\a_i=\(\a_i$ for $i=1,2$. We have $\a_3=\min\{\a_1(\[ a_3,a_4\] ),
R_4-R_3+\a_2\}$. But $R_4-R_3+\a_2=R_4-R_2+\a_1>\a_1\geq\a_3$ so
$\a_3=\a_1(\[ a_3,a_4\] )$. Since $\a_1(\[ a_3,a_4\]
)=\a_3<(R_4-R_3)/2+e$ there is by Lemma 8.8 some $\eta\in\ooo$ with
$d(\eta )=\a_1(\[ a_3,a_4\] )=\a_3$ s.t. $\[ a_3,a_4\]\ap\[\eta
a_3,\eta a_4\]$ so $M\ap\[ a_1,a_2,\eta a_3,\eta a_4\]$. By this
change of BONGs $a_1,a_2$ are preserved and so is the reduction step
(I). On the other hand $d(-a_{1,3}b_1)$ is replaced by $d(-\eta
a_{1,3}b_1)$. Since $d(\eta )=\a_3<d(-a_{1,3}b_1)$ we have $d(-\eta
a_{1,3}b_1)=\a_3$. So after this change of BONGs $\a_2+d(-a_{1,3}b_1)$
becomes $\a_2+\a_3<2e$. This implies that neither (a) nor (b) holds
and we are done.

If $\a_1<\a_3$ then, regardless of the choice of BONG, we have
$\a_1=\min\{\(\a_1,R_3-R_1+\a_3\} =\min\{\(\a_1,\a_3\} =\(\a_1$. It
follows that $\a_i=\(\a_i$  for $i=1,2$ holds regarless of the
BONG. If $\a_3=\a_1(\[ a_3,a_4\] )$ then the proof follows the same as
in the case $\a_1\geq\a_3$ above. Otherwise we use Corollary 8.11 and,
after a change of BONGs, we have $\a_3=\a_1(\[ a_3,a_4\] )$, as
desired. As seen above, the property that $\a_i=\(\a_i$  for $i=1,2$
is not altered by this change of BONGs and the proof follows as
before. 

c. $\a_2+\a_3=2e$. First we claim that $\a_1,\a_2,\a_3$ are all odd
and $<2e$ so they all belong to $d(\ooo )$. We have $\a_2+\a_3=2e$ so
$\a_2,\a_3$ have the same parity. If they are both even then
$\a_2=(R_3-R_2)/2+e$ and $\a_3=(R_4-R_3)/2+e$ so
$\a_2+\a_3=(R_4-R_2)/2+2e>2e$. Contradiction. So they are both odd
and $<2e$ (because $\a_2+\a_3=2e$). Since $R_1=R_3$ we have
$\a_1\ev\a_2\m2$ and $\a_1+\a_2\leq 2e$ by 8.7 so $\a_1$ is odd and
$<2e$ as well. 

Next we claim that $\a_i=\(\a_i$ for $i=1,2$ and $\a_3=\a_1(\[
a_3,a_4\] )$ hold regardless of the choice of BONG. To prove that
$\a_i=\(\a_i$ for $i=1,2$ it is enough to show that $\a_2=\(\a_2$. We
have $\a_2=\min\{\(\a_2,R_3-R_2+\a_3\}$. If $\a_2=(R_3-R_2)/2+e$ then
$\a_2=\(\a_2$ follows from $\a_2\leq\(\a_2\leq (R_3-R_2)/2+e$. If
$\a_2<(R_3-R_2)/2+e$ then $2\a_2<R_3-R_2+2e=R_3-R_2+\a_2+\a_3$ so
$\a_2<R_3-R_2+\a_3$, which again implies $\a_2=\(\a_2$. We have
$\a_3=\min\{\a_1(\[ a_3,a_4\] ), R_4-R_3+\a_2\}$. If
$\a_3=(R_4-R_3)/2+e$ then $\a_3=\a_1(\[ a_3,a_4\] )$ follows from
$\a_3\leq \a_1(\[ a_3,a_4\] )\leq (R_4-R_3)/2+e$. If
$\a_3<(R_4-R_3)/2+e$ then $2\a_3<R_4-R_3+2e=R_4-R_3+\a_2+\a_3$ so
$\a_3<R_4-R_3+\a_2$, which implies $\a_3=\a_1(\[ a_3,a_4\] )$. 

Next we reduce to the case when $d(-a_{1,2})>\a_2$. In general,
$d(-a_{1,2})\geq d[-a_{1,2}]=\a_2$. Suppose that
$d(-a_{1,2})=\a_2$. Take first the case when $\a_3<(R_4-R_3)/2+e$. By
Corollary 8.11 we can change the BONG of $\( M$ s.t. $\a_1(\[
a_2,a_3\] )=\a_2$. Let $K\ap\[ a_2,a_3,a_4\]$. Since
$\a_2\leq\a_1(K)\leq\a_1(\[ a_2,a_3\] )$ we have $\a_1(K)=\a_2$. Since
$\a_3=\a_1(\[ a_3,a_4\] )<(R_4-R_3)/2+e$ we have
$\a_3=R_4-R_3+d(-a_{3,4})$. Let $\eta\in\ooo$ s.t. $d(-\eta
a_{1,2})>\a_2$ (e.g. $\eta\in -a_{1,2}\fs$).  Since $d(-a_{1,2})=\a_2$
we have $d(\eta )=\a_2$. We apply the ternary case of our lemma,
already proved, for $b_2=\eta a_2$ and $K$. We have $\rr
(K)=(R_2,R_3,R_4)$ and $R_2<R_4$ so none of the conditions (a)-(c) of
our lemma holds. So we only have to prove
that $\[ b_2\]\leq K$. By Lemma 8.13 this means $d(a_2b_2)\geq\a_1(K)$
and if $\a_1(K)+d(-a_{2,4}b_2)>2e+R_3-R_4$ then $b_2\rep
[a_2,a_3]$. Now $d(a_2b_2)=d(\eta )=\a_2=\a_1(K)$. For the other
condition note that $b_2=\eta a_2$ implies that
$d(-a_{2,4}b_2)=d(-\eta a_{3,4})$. Since $\a_3<(R_4-R_3)/2+e$ we have
$2\a_3<R_4-R_3+2e=R_4-R_3+\a_2+\a_3$ so
$d(-a_{3,4})=R_3-R_4+\a_3<\a_2=d(\eta )$. It follows that $d(-\eta
a_{3,4})=R_3-R_4+\a_3$ so
$\a_1(K)+d(-a_{2,4}b_2)=\a_2+R_3-R_4+\a_3=2e+R_3-R_4$. Thus the
condition $b_2\rep [a_2,a_3]$ is not necessary. So by the ternary case
we have $\[ a_2,a_3,a_4\]\ap K\ap\[ a'_2,a'_3,a'_3\]$ with
$a'_2=b_2=\eta a_2$. Thus $M\ap\[ a_1,a'_2,a'_3,a'_4\]$. By this
change of BONGs $d(-a_{1,2})$ is replaced by $d(-a_1a'_2)=d(-\eta
a_{1,2})>\a_2$ so we are done. Suppose now that
$\a_3=(R_4-R_3)/2+e$. We have $\a_1(\[ a_2,a_3\]
)=\a_2=2e-\a_3=(R_3-R_4)/2+e<(R_3-R_2)/2+e$. It follows that
$\a_2=\a_1(\[ a_2,a_3\] )=R_3-R_2+d(-a_{2,3})$ so
$d(-a_{2,3})=R_2-R_3+\a_2=\a_1$. Also $\a_2<(R_3-R_2)/2+e$ implies
that $\a_1+\a_2<2e$ and $d(-a_{2,3})=\a_1$ is odd and $<2e$. We are
looking for some $\eta\in\ooo$ s.t. $-\eta a_{1,2}$ satisfies both
$d(-\eta a_{1,2})>\a_2$ and $(-\eta a_{1,2},-a_{2,3})_\p
=(-a_{1,2},-a_{2,3})_\p$. Recall that $a_{1,2}\in\os$. If
$(-a_{1,2},-a_{2,3})_\p =1$ we take $\eta\in\ooo$ s.t. $-\eta
a_{1,2}\in\fs$, which implies both $d(-\eta a_{1,2})>\a_2$ and $(-\eta
a_{1,2},-a_{2,3})_\p =1$. If $(-a_{1,2},-a_{2,3})_\p =-1$ we take
$\eta\in\ooo$ s.t. $d(-\eta a_{1,2})=2e-d(-a_{2,3})=2e-\a_1>\a_2$ and
$(-\eta a_{1,2},-a_{2,3})_\p =-1$. Now $d(-a_{1,2})=\a_2<d(-\eta
a_{1,2})$ so $d(\eta )=\a_2=\a_1(\[ a_2,a_3\] )$ and $(-\eta
a_{1,2},-a_{2,3})_\p =(-a_{1,2},-a_{2,3})_\p$ so $(\eta,-a_{2,3})_\p
=1$. Thus $\eta\in\upon{\a_1(\[ a_2,a_3\] )}{-a_{2,3}}=g(a_3/a_2)$
(see [B3, 3.8]), which implies that $M\ap\[ a_1,\eta a_2,\eta
a_3,a_4\]$. By this change of BONGs $d(-a_{1,2})$ is replaced by
$d(-\eta a_{1,2})$, which is $>\a_2$ so we are done. 

Suppose first that $\[ a_3,a_4\]$ satisfies the conditions of Lemma
8.8. We have $\[ a_3,a_4\]\ap\[\eta a_3,\eta a_4\]$ for some
$\eta\in\ooo$ with $d(\eta )=\a_1(\[ a_3,a_4\] )=\a_3$. Hence $M\ap\[
a_1,a_2,\eta a_3,\eta a_4\]$. We claim that after this change of BONGs
$b_1$ and $\( M$ no longer satisfy (a) or (b). First note that after
this change of BONGs $(-a_{1,2},-a_{2,3})_\p$ is replaced by
$(-a_{1,2},-\eta a_{2,3})_\p$. But $d(-a_{1,2})+d(\eta )>\a_2+\a_3=2e$
so $(-a_{1,2},\eta )_\p =1$. Hence $(-a_{1,2},-\eta a_{2,3})_\p
=(-a_{1,2},-\eta a_{2,3})_\p$, which implies that the property that
$[a_1,a_2,a_3]$ is isotropic or anisotropic is preserved.

If $b_1$ and $\( M$ satisfy the condition (a) then $[a_1,a_2,a_3]$ is
anisotropic. After the change of BONGs $[a_1,a_2,a_3]$ remains
anisotropic so the new $\( M$ cannot satisfy the condition (b). On
the other hand $\a_2+d(-a_{1,3}b_1)$ is replaced by $\a_2+d(-\eta
a_{1,3}b_1)$. But $\a_2+d(-a_{1,3}b_1)>2e$ and $\a_2+d(\eta
)=\a_2+\a_3=2e$. By the domination principle $\a_2+d(-\eta
a_{1,3}b_1)=2e$ so the new $\( M$ doesn't satisfy (a) either.

If $b_1$ and $\( M$ satisfy the condition (b) then $[a_1,a_2,a_3]$ is
anisotropic. After the change of BONGs $[a_1,a_2,a_3]$ remains
isotropic so the new $\( M$ can not satisfy the condition (a). Again
$\a_2+d(-a_{1,3}b_1)$ is replaced by $\a_2+d(-\eta a_{1,3}b_1)$. Now
$\a_2+d(-a_{1,3}b_1)=2e=\a_2+\a_3$ so $d(-a_{1,3}b_1)=\a_3=d(\eta
)$. Since $|\oo/\p |=2$ we have by Lemma 8.1(ii) $d(-\eta
a_{1,3}b_1)>\a_3$. Hence $\a_2+d(-\eta a_{1,3}b_1)>\a_2+\a_3=2e$ so
the new $\( M$ doesn't satisfy (b) either.

Suppose now that $\[ a_3,\a_4\]$ doesn't satisfy the conditions of
Lemma 8.8. Hence $\a_3=\a_1(\[ a_3,a_4\] )=(R_4-R_3)/2+e$ and we have
one of the cases (a)-(c) of the Lemma 8.8. Now (a) doesn't hold
because $\a_3\in d(\ooo )$ and (c) doesn't hold because the dimension
is 2. Hence we have (b), i.e. $|\oo/\p |=2$ and
$d(-a_{3,4})=e-(R_4-R_3)/2=2e-\a_3=\a_2$. Since also
$d(-a_{1,2})>\a_2$ we have
$d[a_{1,4}]=d(a_{1,4})=\a_2=e-(R_4-R_3)/2$. We also have
$\a_2+d(-a_{1,3}b_1)\geq 2e=\a_2+\a_3$ (both when $b_1$ and $\( M$
satisfy (a) or (b)). It follows that $d(-a_{1,3}b_1)\geq\a_3$. Since
$d[-a_{1,3}b_1]=\min\{ d(-a_{1,3}b_1),\a_3\}$ we get
$d[-a_{1,3}b_1]=\a_3=(R_4-R_3)/2+e$. We will show that if $b_1$ and
$\( M$ satisfy (a) or (b) then $[a_1,a_2,a_3,a_4]\top [b_1]$ is not
isotropic so $b_1$ and $M$ satisfy (c), contrary to our assumption. 

Let $V\ap [a_1,a_2,a_3]$ and $V'\ap [a_1,a_2,a_3,a_4]\top [b_1]$.
We have $\det V=a_{1,3}$ and $\det V'=a_{1,4}b_1$. Now $V\pp
[a_4]\ap V'\pp [b_1]\ap [a_1,a_2,a_3,a_4]$ so $a_4\det V=b_1\det
V'=a_{1,4}$. Considering Hasse's symbols we have $SV(a_4,a_4\det
V)_\p =SV'(b_1,b_1\det V')_\p$, i.e. $SV(a_4,a_{1,4})_\p
=SV'(b_1,a_{1,4})_\p$. (Recall, in general, $S(V\pp [\a
])=SV(\a,\a\det V)_\p$.) It follows that $SV'=(a_4b_1,a_{1,4})_\p
SV=(-a_{1,3}b_1,a_{1,4})_\p SV$. If (a) holds for $b_1$ and $\( M$
then $V$ is anisotropic and
$d(a_{1,4})+d(-a_{1,3}b_1)=\a_2+d(-a_{1,3}b_1)>2e$, which implies that
$(-a_{1,3}b_1,a_{1,4})_\p =1$. It follows that $SV'=SV$. Since $V$ is
anisotropic so is $V'$. If (b) holds then
$d(a_{1,4})+d(-a_{1,3}b_1)=\a_2+d(-a_{1,3}b_1)=2e$. Since also
$|\oo/\p |=2$ we have $(-a_{1,3}b_1,a_{1,4})_\p =-1$ by Lemma
8.2(ii). So $SV'=-SV$. But $V$ is isotropic so again $V'$ is
anisotropic, as claimed. 

Suppose now that $n\geq 5$. We denote $\( M\ap\[ a_1,a_2,a_3,a_4\]$
and $\(\a_i:=\a_i(\( M)$ for $i=1,2,3$. We will show that for a good
choice of the BONG for $M$ we have $\[ b_1\]\leq\( M$ and none of the
conditios (a)-(c) hold for $b_1$ and $\( M$. By the quaternary case we
have $\( M\ap\[ a'_1,a'_2,a'_3,a'_4\]$ relative to a good BONG with
$a'_1=b_1$. So $M\ap\[ a'_1,a'_2,a'_3,a'_4,a_5,\ldots,a_n\]$ and we
are done. 

By Corollary 8.11 we may assume that $\a_1(\[ a_3,a_4\] )=\a_3$. It
follows that $\a_i=\(\a_i$ for $i\leq 3$. (For $i\leq 3$ we have
$\a_i=\min\{\(\a_i,R_3-R_i+\a_3\}$ and $\(\a_i\leq R_3-R_i+\a_1(\[
a_3,a_4\] )=R_3-R_1+\a_3$.) Denote by $\( d[-a_{1,3}b_1]$ and $\(
d[a_{1,4}]$ the $d[-a_{1,3}b_1]$ and $d[a_{1,4}]$ corresponding to $\(
M$. Now $\( d[-a_{1,3}b_1]=\min\{ d(-a_{1,3}b_1),\(\a_3\}$ and
$d[-a_{1,3}b_1]=\min\{ d(-a_{1,3}b_1),\a_3\}$ so $\(\a_3=\a_3$ implies
$\( d[-a_{1,3}b_1]=d[-a_{1,3}b_1]$. Also $\( d[a_{1,4}]=d(a_{1,4})$. 

Since $\(\a_i=\a_i$ for $i\leq 3$ and $\(
d[-a_{1,3}b_1]=d[-a_{1,3}b_1]$ the condition $\[ b_1\]\leq\( M$ is
equivalent to $\[ b_1\]\leq M$ and the conditions (a) and (b) for
$b_1$ and $\( M$ are equivalent to the corresponding conditions
for $b_1$ and $M$. Thus $\[ b_1\]\leq\( M$ and $b_1$ and $\( M$
don't satisfy (a) or (b). So we may assume that (c) holds for
$b_1$ and $\( M$. It follows that $|\oo/\p |=2$, $R_2<R_4$,
$d[-a_{1,3}b_1]=\a_3=(R_4-R_3)/2+e$,
$d(a_{1,4})=\a_2=e-(R_4-R_3)/2$ and $[a_1,a_2,a_3,a_4]\top [b_1]$
is anisotropic. If $\a_4>\a_2=e-(R_4-R_3)/2$ then we also have
$d[a_{1,4}]=\min\{ d(a_{1,4}),\a_4\} =e-(R_4-R_3)/2$ so (c) also
holds for $b_1$ and $M$. Contradiction. So $\a_4\leq
e-(R_4-R_3)/2$. Since also $\a_4\geq R_3-R_4+\a_3=e-(R_4-R_3)/2$
we have $\a_4=e-(R_4-R_3)/2$.

Now for an arbitrary BONG we have $\a_3\leq\a_1(\[ a_3,a_4\] )\leq
(R_4-R_3)/2+e$ so $\a_1(\[ a_3,a_4\] )=\a_3=(R_4-R_3)/2+e$ holds
regardless of the choice of the BONG.

We have $\a_2=(R_3-R_4)/2+e<(R_3-R_2)/2+e$. By 8.7 we also
have $\a_1<(R_2-R_1)/2+e$. Hence $\a_1,a_2$ are odd and $<2e$. Since
$\a_3=2e-\a_2$, $\a_3$, too, is odd and $<2e$.

By Corollary 8.11 we may assume that $\a_1(\[ a_1,a_2\] )=\a_1$.
Since $\a_1<(R_2-R_1)/2+e$ we have $\a_1=\a_1(\[ a_1,a_2\]
)=R_2-R_1+d(-a_{1,2})$. It follows that
$d(-a_{1,2})=R_1-R_2+\a_1=\a_2$. Since also $d(a_{1,4})=\a_2$ we
have by Lemma 8.1(ii) $d(-a_{3,4})>\a_2$.

If $R_3<R_5$ we use Lemma 8.8 for $\[ a_4,\ldots,a_n\]$. We have
$\a_4=\min\{\a_1(\[ a_4,\ldots,a_n\] ),R_5-R_4+\a_3\}$. But
$R_5-R_4+\a_3>R_3-R_4+\a_3=\a_4$ (we have $\a_3=(R_4-R_3)/2+e$ and
$\a_4=e-(R_4-R_3)/2$) so $\a_1(\[ a_4,\ldots,a_n\]
)=\a_4=\a_2=e-(R_4-R_3)/2$. Since $e-(R_4-R_3)/2<(R_5-R_4)/2+e$ we can
apply Lemma 8.8. So there is $\eta\in\ooo$ with $d(\eta
)=e-(R_4-R_3)/2$ s.t. $\[ a_4,\ldots,a_n\]\ap\[ a'_4,\ldots,a'_n\]$
relative to some other good BONG, with $a'_4=\eta a_4$. It follows
that $M\ap\[ a_1,a_2,a_3,a'_4,\ldots,a'_n\]$. After this change of
BONGs $d(a_{1,4})$ is replaced by $d(a_{1,3}a'_4)=d(\eta
a_{1,4})$. But $d(\eta )=d(a_{1,4})=e-(R_4-R_3)/2$. By Lemma 8.1(ii)
we have $d(\eta a_{1,4})>e-(R_4-R_3)/2$. Therefore $b_1$ and $\[
a_1,a_2,a_3,a'_4\]$ (the new $\( M$) no longer satisfy condition (c).

If $R_3=R_5$ we apply Corollary 8.9 to $K\ap\[
a_1,a_2,a_3,a_4,a_5\]$. Since $\a_3\leq\a_3(K)\leq (R_4-R_3)/2+e$ we
have $\a_3(K)=\a_3=(R_4-R_3)/2+e$. Since $R_3=R_5$ we also have
$\a_4(K)=(R_5-R_4)/2+e=e-(R_4-R_3)/2=\a_2$ (see 8.7). We want to prove
that $K$ doesn't satisfy any of the conditions (a)-(c) of Corollary
8.9. Now $\a_4(K)=\a_2$ is odd and $<2e$ so it belongs to $d(\ooo
)$. Thus (a) doesn't hold. We have
$\a_3(K)=(R_4-R_3)/2+e=e-(R_5-R_4)/2$ so (b) doesn't hold
either. Finally, $d(-a_{3,4})>\a_2=e-(R_4-R_3)/2$ and since $R_3=R_5$
we have by 8.7 $d(-a_{4,5})\geq\a_3=(R_4-R_3)/2+e$. It follows that
$d(-a_{3,4})+d(-a_{4,5})>2e$ so $(-a_{3,4},-a_{4,5})_\p =1$. Hence
$[a_3,a_4,a_5]$ is isotropic and (c) doesn't hold. By Corollary 8.9
there is $\eta\in\ooo$ with $d(\eta )=\a_4(K)=\a_2$ s.t. $K\ap\[
a'_1,a'_2,a'_3,a'_4,a'_5\]$ relative to some good BONG with $a'_5=\eta
a_5$. It follows that $M\ap\[
a'_1,a'_2,a'_3,a'_4,a'_5,a_6,\ldots,a_n\]$. By this change of BONGs
$a_{1,4}$ is replaced by $a'_{1,4}$. But in $\ff/\fs$ we have
$a'_{1,5}=a_{1,5}=\det K$ and $a'_5=\eta a_5$ so $a'_{1,4}=\eta
a_{1,4}$. Since $d(\eta )=d(a_{1,4})=\a_2$ we have by Lemma 8.1(ii)
$d(a'_{1,4})>\a_2$. Therefore $b_1$ and $\[ a'_1,a'_2,a'_3,a'_4\]$,
the new $\( M$, no longer satisfy (c) and we ar done. \qed

\section{The final step}

We now use induction on $n$ to prove our theorem. Suppose that it is
true for dimension $<n$. Suppose that $M,N$ are 2 lattices over the 
same quadratic space $V$ with $\dim V=n$ s.t. $N\leq M$ but $N\notrep
M$. Suppose moreover that $M,N$ are chosen with the  property that
$\ord volN-\ord volM$ is minimal. From $\S$7 we now that this implies
that $\nn M=\nn N$, i.e. $R_1=S_1$. We will show that there are $x,y$,
norm generators for $M$ and $N$, respectively s.t. $Q(x)=Q(y)$ and
$pr_{y^\pp}N\leq pr_{x^\pp}M$. By the induction step $pr_{y^\pp}N\leq
pr_{x^\pp}M$ implies $pr_{y^\pp}N\rep pr_{x^\pp}M$. We have $V=Fy\pp
pr_{y^\pp}V=Fx\pp pr_{x^\pp}V$. Since $Q(x)=Q(y)$ and $pr_{y^\pp}N\rep
pr_{x^\pp}M$ there is $\phi\in O(V)$ s.t. $\phi (y)=x$ and $\phi
(pr_{y^\pp}N)\sbq pr_{x^\pp}M$. We have $\nn\phi (N)=\nn N=\nn M$, $x$
is a norm generator of $M$ and $pr_{x^\pp}\phi (N)=pr_{\phi
(y)^\pp}\phi (N)=\phi (pr_{y^\pp}N)\sbq pr_{x^\pp}M$. By [B1, Lemma
2.2] we have $\phi (N)\sbq M$ so $N\rep M$. 

In the most typical case from Lemma 9.3 we will show that $M\ap\[
a_1,\ldots,a_n\]$ and $N\ap\[ b_1,\ldots,b_n\]$ relative to some good
BONGs $x_1,\ldots,x_n$ and $y_1,\ldots,y_n$
s.t. $Q(x_1)=a_1=b_1=Q(y_1)$ and $pr_{y_1^\pp}N\ap\[
b_2,\ldots,b_n\]\leq\[ a_2,\ldots,a_n\]\ap pr_{x_1^\pp}M$. So in the
reasoning above we can take $x=x_1$ and $y=y_1$. Similarly, in the
case of Lemma 9.6 we write $N\ap\[ b_1,\ldots,b_n\]$ relative to some
good BONG $y_1,\ldots,y_n$ but $M\ap\[ a'_1,\ldots,a'_n\]$ relative to
some bad BONG $x'_1,\ldots,x'_n$ s.t. $a'_1=b_1$ and
$\[b_2,\ldots,b_n\]\leq\[ a'_2,\ldots,a'_n\]$ so we can take $x=x'_1$
and $y=y_1$. Finally, in all other cases we show that the property of
minimality for  $\ord volN-\ord volM$ is contradicted by the same
reason from $\S$7. Namely in Lemma 9.12 we show that if $\nn M=\nn K$,
$vol K\sb vol M$ and $M,K$ don't satisfy the conditions from the
hypothesis of Lemma 9.3 or Lemma 9.6 then there is $N\sb M$
s.t. $K\leq N$. Note that the very particular case from Lemma 9.6
could have been included in Lemma 9.12 but the methods used there
don't work in this case. Instead a longer proof, similar to those from
$\S$7, would have been necessary. Therefore we chose the atypical
solution of using a bad BONG for $M$.

\blm Suppose that $N\leq M$ and $R_1=S_1$. If one of the following
happens: $R_1<R_3$, $R_2=S_2$, $R_2-R_1=2e$, $R_2=R_4$ or
$d[-a_{1,3}b_1]=\a_1<\b_1$ then $M\ap\[ a'_1,\ldots,a'_n\]$ relative
to some good BONG with $a'_1=b_1$.
\elm
\pf Since $N\leq M$ and $\[ b_1\]\leq N$ we have $\[ b_1\]\leq
M$. (Note that the transitivity of $\leq$ was proved only for lattices
with the same dimension. However we can use an argument similar to
that from Lemma 2.21. Let $K$ be a lattice s.t. $FK\ap
[b_2,\ldots,b_n]$. We have $[b_1]\pp FK\ap FN\ap FM$. By Lemma 2.20
for $s\gg 0$ we have $\[ b_1\] \leq N$ iff $\[ b_1\]\pp\p^sK\leq N$
and $\[ b_1\] \leq M$ iff $\[ b_1\]\pp\p^sK\leq M$. We have $\[
b_1\]\leq N$ so $\[ b_1\]\pp\p^sK\leq N$, which, together with $N\leq
M$, implies that $\[ b_1\]\pp\p^sK\leq M$ so $\[ b_1\]\leq M$.)

Thus we have to prove that, given the conditions in the hypothesis,
none of the conditions (a) - (c) of Lemma 8.14 happens. This is
obvious if $R_1<R_3$. From now on we assume that $R_1=R_3$. For
consequences of $R_1=R_3$ see 8.7. 

If $R_2-R_1=2e$ then $g(a_2/a_1)=\upo{2e}$ and $\a_1=2e$. We have
$d(a_1,b_1)\geq\a_1=2e$ so $b_1=\e a_1$ with $d(\e )\geq 2e$. It
follows that $\e\in\upo{2e}=g(a_2/a_1)$ so $\[ a_1,a_2\]\ap\[\e a_1,\e
a_2\] =\[ b_1,\e a_1\]$ so $M\ap\[ b_1,\e a_2,a_3,\ldots,a_n\]$
relative to some good BONG.

If $R_2=R_4$ then (c) cannot hold. Suppose that (a) or (b) holds. We
have $\a_2+\a_3\geq\a_2+d[-a_{1,3}b_1]$, which is $>2e$ if (a) holds
or just $\geq 2e$ if (b) holds. But $R_2=R_4$ so $\a_2+\a_3\leq
2e$. It follows that (b) holds and $\a_2+\a_3=2e$. By 8.7 this implies
$\a_2=(R_3-R_2)/2+e$. Since $R_1=R_3$, again by 8.7, we get
$\a_1=(R_2-R_1)/2+e$ and so (b) doesn't hold. 

For the remaining two cases, $R_2=S_2$ and $d[-a_{1,3}b_1]=\a_1<\b_1$,
we denote by $\( d[-a_{1,3}b_1]:=\min\{ d(-a_{1,3}b_1),\a_3\}$, which
is $d[-a_{1,3}b_1]$ corresponding to $M$ and $\[ b_1\]$. Therefore
$d[-a_{1,3}b_1]=\min\{\( d[-a_{1,3}b_1],\b_1\}$.

If $d[-a_{1,3}b_1]=\a_1<\b_1$ then $\(
d[-a_{1,3}b_1]=d[-a_{1,3}b_1]=\a_1$. We have $\a_2+\(
d[-a_{1,3}b_1]=\a_1+\a_2\leq 2e$, with equality iff
$\a_1=(R_2-R_1)/2+e$. If $\a_1<(R_2-R_1)/2+e$ then $\a_2+\(
d[-a_{1,3}b_1]<2e$ so none of (a) - (c) happens. If
$\a_1=(R_2-R_1)/2+e$ then $\a_2+\( d[-a_{1,3}b_1]=2e$. Then we cannot
have (a) and since $\a_1=(R_2-R_1)/2+e$ we don't have (b) or
(c). (Note that if (c) holds then $R_1=R_3$ and
$\a_2=e-(R_4-R_3)/2<(R_3-R_2)/2+e$ so $\a_1<(R_2-R_1)/2+e$ by 8.7.) 

Suppose now that $R_2=S_2$. We may assume that $R_2<R_4$ or $n=3$
since the case $R_2=R_4$ was treated above. Note that if
$d[-a_{1,3}b_1]=\b_1$ then by Lemma 8.12(ii) we have $A_2=\min\{
(R_3-S_2)/2+e,R_3-S_2+d[-a_{1,3}b_1]\} =\min\{
(R_3-R_2)/2+e,R_3-R_2+\b_1\}$. But $\b_1\leq
(S_2-S_1)/2+e=(R_2-R_3)/2+e$ so $R_3-R_2+\b_1\leq (R_3-R_2)/2+e$ so
$A_2=R_3-R_2+\b_1$. Also by Lemma 6.2 we have $R_1+\a_1\leq
S_1+\b_1=R_1+\b_1$ so $\a_1\leq\b_1$, which implies $R_3-R_2+\b_1\geq
R_3-R_2+\a_1=\a_2\geq A_2=R_3-R_2+\b_1$. So if $d[-a_{1,3}b_1]=\b_1$
then $A_2=\a_2$ and $\b_1=\a_1$.

Suppose that (a) of Lemma 8.14 holds. We claim that $[b_1,b_2]\rep
[a_1,a_2,a_3]$. This obviously holds if $n=3$, when $FM\ap
[a_1,a_2,a_3]$, so we may assume that $n\geq 4$. We have $S_2=R_2<R_4$
so we still need $A_2+A_3>2e+R_3-S_3$. By Lemma 2.14 we may assume
that $A_2=A'_2$ and $A_3=A'_3$. By Lemma 8.12(ii)
$A'_2=R_3-S_2+d[-a_{1,3}b_1]$ and since $R_4>S_2$ we have
$A_3=A'_3=\min\{ R_4-S_3+d[-a_{1,4}b_{1,2}],
R_4+R_5-S_2-S_3+d[-a_{1,3}b_1]\}$. (The second term of $A'_3$ is
ignored if $n=4$.) So we have to prove that $2e+R_3-S_3<A_2+A_3=\min\{
R_3+R_4-S_2-S_3+d[-a_{1,3}b_1]+d[-a_{1,4}b_{1,2}],
R_3+R_4+R_5-2S_2-S_3+2d[-a_{1,3}b_1]\}$, i.e. that
$d[-a_{1,3}b_1]+d[-a_{1,4}b_{1,2}]>2e+S_2-R_4=2e+R_2-R_4$ and, if
$n>5$, $2d[-a_{1,3}b_1]>2e+2S_2-R_4-R_5=2e+2R_2-R_4-R_5$. Assume
first that $d[-a_{1,3}b_1]=\( d[-a_{1,3}b_1]$. Since (a) holds we have
$\a_2+d[-a_{1,3}b_1]=\a_2+\( d[-a_{1,3}b_1]>2e$ so
$d[-a_{1,3}b_1]>2e-\a_2\geq 2e-((R_3-R_2)/2+e)=e-(R_3-R_2)/2$. Thus
$2d[-a_{1,3}b_1]>2e+R_2-R_3>2e+2R_2-R_4-R_5$. ($R_4>R_2$ and $R_5\geq
R_3$.) Now $d[a_{1,2}b_{1,2}]\geq A_2=R_3-S_2+d[-a_{1,3}b_1]$ and so
$d[-a_{1,3}b_1]+d[a_{1,2}b_{1,2}]\geq
R_3-S_2+2d[-a_{1,3}b_1]>R_3-R_2+2(e-(R_3-R_2)/2)=2e>2e+R_2-R_4$. Also
$d[-a_{1,3}b_1]+d[-a_{3,4}]\geq R_3-R_4+\a_3+d[-a_{1,3}b_1]\geq
R_2-R_4+\a_2+d[-a_{1,3}b_1]>2e+R_2-R_4$. Therefore
$d[-a_{1,3}b_1]+d[-a_{1,4}b_{1,2}]>2e+R_2-R_4$. Suppose now that 
$d[-a_{1,3}b_1]=\b_1$ so $A_2=\a_2$ and $d[-a_{1,3}b_1]=\b_1=\a_1$. We
have $R_3-R_2+\a_1+R_4-R_2+\a_1\geq\a_2+\a_3\geq\a_2+\(
d[-a_{1,3}b_1]>2e$ and so
$2d[-a_{1,3}b_1]=2\a_1>2e+2R_2-R_3-R_4\geq 2R_2-R_4-R_5$. Similarly as
before, $d[-a_{1,3}b_1]+d[a_{1,2}b_{1,2}]\geq
R_3-S_2+2d[-a_{1,3}b_1]>R_3-R_2+2e+2R_2-R_3-R_4=2e+R_2-R_4$ and
$d[-a_{1,3}b_1]+d[-a_{3,4}]\geq\a_1+R_3-R_4+\a_3=R_2-R_4+\a_2+\a_3\geq
R_2-R_4+\a_2+\( d[-a_{1,3}b_1]>2e+R_2-R_4$. (By 8.7
$\a_1=R_2-R_3+\a_2$.) So
$d[-a_{1,3}b_1]+d[-a_{1,4}b_{1,2}]>2e+R_2-R_4$. Hence $[b_1,b_2]\rep
[a_1,a_2,a_3]$ so $[a_1,a_2,a_3]\ap [b_1,b_2,a_{1,3}b_{1,2}]$. But
$d(-b_{1,2})\geq S_1-S_2+\b_1=R_1-R_2+\a_1=\a_2$. It follows that
$d(-b_{1,2})+d(-a_{1,3}b_1)\geq\a_2+\( d(-a_{1,3}b_1)>2e$ and so
$(-b_{1,2},-a_{1,3}b_1)_\p =1$, which implies that
$[b_1,b_2,a_{1,3}b_{1,2}]\ap [a_1,a_2,a_3]$ is isotropic. But this
cotradicts (a). 

If we have (b) or (c) of Lemma 8.14 then $\a_2+\(
d[-a_{1,3}b_1]=2e$. (If (c) holds then $\a_2=e-(R_4-R_3)/2$ and $\(
d[-a_{1,3}b_1]=(R_4-R_3)/2+e$.) Also $\a_1<(R_2-R_1)/2+e$,
$\a_2<(R_3-R_2)/2+e$ and $\a_1+\a_2<2e$. (The three statements are
equivalent by 8.7. In the case of (b) we have $\a_1<(R_2-R_1)/2+e$ and
in the case of (c) $\a_2=e-(R_4-R_3)/2<(R_3-R_2)/2+e$.)

Suppose that $d[-a_{1,3}b_1]=\( d[-a_{1,3}b_1]$ so
$\a_2+d[-a_{1,3}b_1]=2e$. It follows that $\a_2\geq A_2=\min\{
(R_3-S_2)/2+e,R_3-S_2+d[-a_{1,3}b_1]\} =\min\{
(R_3-R_2)/2+e,R_3-R_2+2e-\a_2\}$. But $\a_2<(R_3-R_2)/2+e$, which also
implies $\a_2<R_3-R_2+2e-\a_2$. Contradiction. Thus $\(
d[-a_{1,3}b_1]>d[-a_{1,3}b_1]=\b_1$, which implies
$\a_1=\b_1=d[-a_{1,3}b_1]$ and $A_2=\a_2$.

Suppose first that (b) holds. It follows that $\a_2+\a_3>2e=\a_2+\(
d[-a_{1,3}b_1]$. Thus $\( d[-a_{1,3}b_1]<\a_3$ so $\(
d[-a_{1,3}b_1]=\min\{ d(-a_{1,3}b_1),\a_3\} =d(-a_{1,3}b_1)$. (Same
happens if $n=3$, when $\a_3$ is ignored.) It follows that
$\a_2+d(-a_{1,3}b_1)=2e$. Since $\a_1+\a_2<2e$ we also have
$d(-a_{1,3}b_1)>\a_1$. Like in the case (a) we want to prove that
$[b_1,b_2]\rep [a_1,a_2,a_3]$. This is obvious if $n=3$ so we may
assume that $n\geq 4$. Again we have $S_2=R_2<R_4$ and
$A_2+A_3>2e+R_3-S_3$ is equivalent to
$\a_1+d[-a_{1,4}b_{1,2}]=d[-a_{1,3}b_1]+d[-a_{1,4}b_{1,2}]>2e+R_2-R_4$
and, if $n\geq 5$, $2\a_1=2d[-a_{1,3}b_1]>2e+2R_2-R_4-R_5$. Now
$d[a_{1,2}b_{1,2}]\geq A_2=\a_2\geq R_3-R_4+\a_3$ and
$d[-a_{3,4}]\geq R_3-R_4+\a_3$ so $d[-a_{1,4}b_{1,2}]\geq
R_3-R_4+\a_3$. Since also $\a_1=R_2-R_3+\a_2$ we have
$\a_1+d[-a_{1,4}b_{1,2}]\geq R_2-R_4+\a_2+\a_3>2e+R_2-R_4$. Also
$-R_2+\a_1=-R_3+\a_2\geq -R_4+\a_3$ so $2\a_1\geq
R_2-R_3+\a_2+R_2-R_4+\a_3>2e+2R_2-R_3-R_4\geq 2e+2R_2-R_4-R_5$. So
$[b_1,b_2]\rep [a_1,a_2,a_3]$, which implies $[a_1,a_2,a_3]\ap
[b_1,b_2,a_{1,3}b_{1,2}]$. We have $d(-a_{1,3}b_1)>\a_1\geq
S_1-S_3+\a_1$. Also $\b_3,d(a_{1,3}b_{1,3})\geq
A_3>2e+R_3-S_3-A_2=2e+S_1-S_3-\a_2>S_1-S_3+\a_1$. Thus
$d(-b_{2,3})>S_1-S_3+\a_1$ and also $\b_3>S_1-S_3+\a_1$. (If $n=3$
we ignore $\b_3$ but we still have
$d(a_{1,3}b_{1,3})>S_1-S_3+\a_1$ because $a_{1,3}b_{1,3}\in\fs$.)
Now $\a_1=\b_1=\min\{
(S_2-S_1)/2+e,S_2-S_1+d(-b_{1,2}),S_3-S_1+d(-b_{2,3}),S_3-S_1+\b_3\}$.
(Ignore $S_3-S_1+\b_3$ if $n=3$.) But
$S_3-S_1+d(-b_{2,3}),S_3-S_1+\b_3>\a_1$ and also
$(S_2-S_1)/2+e=(R_2-R_1)/2+e>\a_1$. Thus
$\a_1=S_2-S_1+d(-b_{1,2})$ so
$d(-b_{1,2})=S_1-S_2+\a_1=R_1-R_2+\a_1=\a_2$. It follows that
$d(-b_{1,2})+d(-a_{1,3}b_1)=\a_2+d(-a_{1,3}b_1)=2e$. Since
$|\oo/\p |=2$ this implies by Lemma 8.2(iii) that
$(-b_{1,2},-a_{1,3}b_1)_\p =-1$ and so $[a_1,a_2,a_3]\ap
[b_1,b_2,a_{1,3}b_{1,2}]$ is not isotropic. But this contradicts
condition (b).

Suppose now that (c) holds. We have $\a_3=(R_4-R_3)/2+e$,
$\a_2=e-(R_4-R_3)/2<\a_4$ and
$\a_1=R_2-R_3+\a_2=e+R_2-(R_3+R_4)/2=R_2-R_4+\a_3<\a_3$. Now
$d[a_{1,4}]=\a_2<\a_4$ so $\a_2=d[a_{1,4}]=\min\{ d(a_{1,4}),\a_4\}
=d(a_{1,4})$. (This also holds if $n=4$, when $\a_4$ is ignored.)
Since also $d[-b_{1,2}]\geq S_1-S_2+\b_1=R_1-R_2+\a_1=\a_2$ we have
$d[-a_{1,4}b_{1,2}]\geq\a_2$. If $n\geq 5$ then
$\a_3+\a_4>(R_4-R_3)/2+e+e-(R_4-R_3)/2=2e$ so $R_3<R_5$.

We claim that $A_3\geq R_3-S_3+\a_3$, with equality iff $R_3=S_3$ or
$d[-a_{1,4}b_{1,2}]=\a_2$. Indeed, we have $S_3\geq S_1=R_3$ so
$(R_4-S_3)/2+e\geq R_3-S_3+(R_4-R_3)/2+e=R_3-S_3+\a_3$, with equality
iff $R_3=S_3$. Also $R_4-S_3+d[-a_{1,4}b_{1,2}]\geq
R_4-S_3+\a_2=(R_3+R_4)/2-S_3+e=R_3-S_3+\a_3$, with equality iff
$d[-a_{1,4}b_{1,2}]=\a_2$. Finally, if $n\geq 5$ then
$d[-a_{1,3}b_1]=\a_1=R_2-R_4+\a_3=S_2-R_4+\a_3$ so
$R_4+R_5-S_2-S_3+d[-a_{1,3}b_1]=R_5-S_3+\a_3>R_3-S_3+\a_3$ and
$d[-a_{4,5}]\geq R_4-R_5+\a_4>R_4-R_5+e-(R_4-R_3)/2=R_3-R_5+\a_3$ so
$R_4+R_5-S_2-S_3+d[-a_{4,5}]>R_3+R_4-S_2-S_3+\a_3>R_3-S_3+\a_3$ (we
have $S_2=R_2<R_4$). Hence
$R_4+R_5-S_2-S_3+d[a_{1,5}b_1]>R_3-S_3+\a_3$. In conclusion,
$A_3=\min\{ (R_4-S_3)/2+e, R_4-S_3+d[-a_{1,4}b_{1,2}],
R_4+R_5-S_2-S_3+d[a_{1,5}b_1]\}\geq R_3-S_3+\a_3$, with equality iff
$R_3=S_3$ or $d[-a_{1,4}b_{1,2}]=\a_2$.

We have $d(-a_{1,3}b_1)\geq\(
d[-a_{1,3}b_1]=\a_3>\a_1=\b_1\geq S_1-S_3+\b_1$. Also $A_3\geq
R_3-S_3+\a_3>S_1-S_3+\b_1$. (We have $R_3=S_1$ and
$\b_1=\a_1=R_2-R_4+\a_3<\a_3$.) It follows that
$d(a_{1,3}b_{1,3}),\b_3\geq d[a_{1,3}b_{1,3}]\geq
A_3>S_1-S_3+\b_1$. Together with $d(-a_{1,3}b_1)>S_1-S_3+\b_1$, this
implies $d(-b_{2,3})>S_1-S_3+\b_1$. Now $\a_1=\b_1=\min\{
(S_2-S_1)/2+e, S_2-S_1+d(-b_{1,2}), S_3-S_1+d(-b_{2,3}),
S_3-S_1+\b_3\}$. But $(S_2-S_1)/2+e=(R_2-R_1)/2+e>\a_1=\b_1$ and
$\b_3,d(-b_{2,3})>S_1-S_3+\b_1$ so
$S_3-S_1+d(-b_{2,3}),S_3-S_1+\b_3>\b_1$. It follows that
$\b_1=S_2-S_1+d(-b_{1,2})$ so
$d(-b_{1,2})=S_1-S_2+\b_1=R_1-R_2+\a_1=\a_2$. We also have
$d(a_{1,4})=\a_2$ and $|\oo/\p |=2$ so $d(-a_{1,4}b_{1,2})>\a_2$ by
Lemma 8.1(ii). 

There are two cases:

1. $S_1<S_3$. We have $\b_2=\min\{ (S_3-S_2)/2+e, S_3-S_2+d(-b_{2,3}),
S_3-S_2+\b_1, S_3-S_2+\b_3\}$. Now $(S_3-S_2)/2+e>(R_3-R_4)/2+e=\a_2$
and $S_3-S_2+\b_1>R_3-R_2+\a_1=\a_2$. (we have $R_3=S_1<S_3$ and
$S_2=R_2<R_4$.) Also $\b_3,d(-b_{2,3})>S_1-S_3+\b_1$ (see above) so
$S_3-S_2+d(-b_{2,3}),S_3-S_2+\b_3>S_1-S_2+\b_1=R_1-R_2+\a_1=\a_2$.
Thus $\b_2>\a_2$. Since also $\a_4>\a_2$ and $d(-a_{1,4}b_{1,2})>\a_2$
we have $d[-a_{1,4}b_{1,2}]>\a_2$. Since also $R_3=S_1<S_3$ the
inequality $A_3\geq R_3-S_3+\a_3$ becomes strict. It follows that
$A_2+A_3>\a_2+R_3-S_3+\a_3=2e+R_3-S_3$. Since also $S_2=R_2<R_4$ we
have $[b_1,b_2]\rep [a_1,a_2,a_3]$, which implies that
$[a_1,a_2,a_3]\ap [b_1,b_2,a_{1,3}b_{1,2}]$ so $[a_1,a_2,a_3,a_4]\ap
[b_1,b_2,a_{1,3}b_{1,2},a_4]$ and $[a_1,a_2,a_3,a_4]\top [b_1]\ap
[b_2,a_{1,3}b_{1,2},a_4]$. Since
$d(-a_{1,3}b_1)+d(-a_{1,4}b_{1,2})>\a_3+\a_2=2e$ we have
$(-a_{1,3}b_1,-a_{1,4}b_{1,2})_\p =1$ so $[b_2,a_{1,3}b_{1,2},a_4]\ap
[a_1,a_2,a_3,a_4]\top [b_1]$ is isotropic, which contradicts (c).

2. $S_1=S_3$. We have $\a_3\geq d[a_{1,3}b_{1,3}]\geq A_3\geq
R_3-S_3+\a_3=\a_3$ so $d(a_{1,3}b_{1,3})\geq
d[a_{1,3}b_{1,3}]=A_3=\a_3$. We claim that $[b_1,b_2,b_3]\rep
[a_1,a_2,a_3,a_4]$.  This is obvious if $n=4$, when $FM\ap
[a_1,a_2,a_3,a_4]$, so we may assume that $n\geq 5$. We have to prove
that $S_3<R_5$ and $A_3+A_4>2e+R_4-S_4$. We have $R_5>R_3=S_3$. If we
prove that also $A_4>R_4-S_4+\a_2$ then
$A_3+A_4>\a_3+R_4-S_4+\a_2=2e+R_4-S_4$ and we are done. By Lemma 2.14
may assume that $A_4=A_4'=\min\{ R_5-S_4+d[-a_{1,5}b_{1,3}],
R_5+R_6-S_3-S_4+d[-a_{1,4}b_{1,2}]\}$. Now 
$d[a_{1,3}b_{1,3}]=\a_3\geq R_4-R_5+\a_4$ and $d[-a_{4,5}]\geq
R_4-R_5+\a_4$ so $d[-a_{1,5}b_{1,3}]\geq
R_4-R_5+\a_4>R_4-R_5+\a_2$. Thus
$R_5-S_4+d[-a_{1,5}b_{1,3}]>R_4-S_4+\a_2$. Also
$d[-a_{1,4}b_{1,2}]\geq\a_2$ (see above), $R_6\geq R_4$ and
$R_5>R_3=S_3$ so
$R_5+R_6-S_3-S_4+d[-a_{1,4}b_{1,2}]>R_4-S_4+\a_2$. Thus
$A_4>R_4-S_4+\a_2$ and so  $[b_1,b_2,b_3]\rep [a_1,a_2,a_3,a_4]$. It
follows that $[a_1,a_2,a_3,a_4]\ap [b_1,b_2,b_3,a_{1,4}b_{1,3}]$ so
$[a_1,a_2,a_3,a_4]\top [b_1]\ap [b_2,b_3,a_{1,4}b_{1,3}]$. Now
$d(-a_{1,3}b_1),d(a_{1,3}b_{1,3})\geq\a_3$ so
$d(-b_{2,3})\geq\a_3$. Since also $d(-a_{1,4}b_{1,2})>\a_2$ we get
$d(-a_{1,3}b_1)+d(-a_{1,4}b_{1,2})>\a_3+\a_2=2e$. It follows that
$(-b_{2,3},-a_{1,4}b_{1,2})_\p =1$ and so
$[b_2,b_3,a_{1,4}b_{1,3}]\ap [a_1,a_2,a_3,a_4]\top [b_1]$ is
isotropic, which contradicts (c). \qed

\blm We have $M\ap\[ a'_1,a'_2,\ldots,a'_n\]$ relative to another good
BONG s.t. $a_1=a'_1$ and $\a_i=\a_{i-1}(\[ a'_2,\ldots,a'_n\] )$ for
$4\leq i<n$. Moreover if $R_1<R_3$ or $R_2=R_4$ or $R_2-R_1=2e$ then
the equality above also holds for $i=3$. 
\elm
\pf Suppose first that $R_1<R_3$ or $R_2=R_4$ or $R_2-R_1=2e$. Suppose
that our statement is true for $n=4$. If $n>4$ then by the quaternary
case we have $\[ a_1,a_2,a_3,a_4\]\ap\[ a'_1,a'_2,a'_3,a'_4\]$
relative to some good BONG s.t. $a_1=a'_1$ and $\a_3(\[
a_1,a_2,a_3,a_4\] )=\a_2(\[ a'_2,a'_3,a'_4\] )$. It follows that
$M\ap\[ a_1,\ldots,a_n\]\ap\[ a'_1,\ldots,a'_n\]$ with $a'_i:=a_i$ for
$i>4$. Now for any $i\geq 3$ our statement will follow from
$\a_i=\min\{\a_{i-1}(\[ a'_2,\ldots,a'_n\] ), R_{i+1}-R_4+\a_3(\[
a'_1,a'_2,a'_3,a'_4\] )\}$ and $\a_{i-1}(\[ a'_2,\ldots,a'_n\]
)\leq R_{i+1}-R_4+\a_2(\[ a'_2,a'_3,a'_4\] )=R_{i+1}-R_4+\a_3(\[
a_1',a'_2,a'_3,a'_4\] )$. 

Thus we have reduced to the case $n=4$. We have to prove that $M\ap\[
a_1,a'_2,a'_3,a'_4\]$ relative to
some good BONG s.t. $\a_3=\a_2(\[ a'_2,a'_3,a'_4\] )$. We have
$\a_3=\min\{\a_2(\[ a_2,a_3,a_4\] ),R_4-R_2+\a_1\}$. We can assume
that $\a_3=R_4-R_2+\a_1<\a_2(\[ a_2,a_3,a_4\] )$ since otherwise we
can take $a'_i=a_i$. We have $-R_2+\a_1=-R_4+\a_3$. This implies
$-R_2+\a_1=-R_3+\a_2=-R_4+\a_3$ so $\a_2=R_3-R_2+\a_1$. Now
$R_4-R_3+d(-a_{3,4}),R_4-R_2+d(-a_{2,3})\geq\a_2(\[ a_2,a_3,a_4\]
)>\a_3$. Thus $d(-a_{3,4})>R_3-R_4+\a_3=\a_2$ and
$d(-a_{2,3})>R_2-R_4+\a_3=\a_1$. We have $\a_3<\a_2(\[ a_2,a_3,a_4\]
)\leq (R_4-R_3)/2+e$. Since $\a_3<(R_4-R_3)/2+e$ and
$-R_2+\a_1=-R_4+\a_3$ we also have $\a_1<(R_2-R_1)/2+e$ and
$\a_2<(R_3-R_2)/2+e$ by Lemma 8.4(iii). It follows that
$\a_1,\a_2,\a_3$ are odd and $<2e$ so they belong to $d(\ooo )$. In
particular, this rules out the case $R_2-R_1=2e$, when $\a_1=2e$. Also
note that $\a_2+\a_3=2\a_3+R_3-R_4<2e$. 

Take $\e\in\ooo$ with $d(\e )=\a_2$. We
have $d(-a_{3,4})>\a_2=d(\e )$ so $d(-\e a_{3,4})=\a_2$. Since
$\a_2+\a_3<2e$ and $\a_3\in d(\ooo )$, by Lemma 8.2, for any
choice of the $\pm$ sign there is $\eta\in\ooo$ with $d(\eta )=\a_3$
s.t. $(-\e a_{3,4},\eta )_\p =\pm 1$. We choose such $\eta$ with $(-\e
a_{3,4},\eta )_\p =(\e,-a_{2,3})_\p$. Take $a'_1=a_1$, $a'_2=\e a_2$,
$a'_3=\e\eta a_3$ and $a'_4=\eta a_4$. Let $\( M:=\[
a'_1,a'_2,a'_3,a'_4\]$. We claim that $M\ap\( M$. We have $\det
[a_2,a_3,a_4]=\det [a'_2,a'_3,a'_4]$. Also $(-a'_{2,3},-a'_{3,4})_\p
=(-\eta a_{2,3},-\e a_{3,4})_\p =(-a_{2,3},-a_{3,4})_\p (\eta,-\e
a_{3,4})_\p (-a_{2,3},\e )_\p =(-a_{2,3},-a_{3,4})_\p$. Thus
$[a'_2,a'_3,a'_4]\ap [a_2,a_3,a_4]$, which implies $F\( M\ap
FM$. Condition (i) of [B3, Theorem 3.1] is obvious. For (iii) note
that $d(a_1a'_1)=\j >\a_1$, $d(a_{1,2}a'_{1,2})=d(\e )=\a_2$ and
$d(a_{1,3}a'_{1,3})=d(\eta )=\a_3$. Together with
$a_{1,4}a'_{1,4}\in\fs$, these also imply by Lemma 8.6 that $\( M$
exists and $\(\a_i:=\a_i(\( M)\geq\a_i$. If $R_1<R_3$ we prove first
that $\(\a_1=\a_1$. We have $\a_1=\min\{ (R_2-R_1)/2+e,
R_2-R_1+d(-a_{1,2}), R_2-R_1+\a_2\}$. But $\a_1<(R_2-R_1)/2+e$ and
$R_2-R_1+\a_2=R_3-R_1+\a_1>\a_1$ so $\a_1=R_2-R_1+d(-a_{1,2})$. Since
also $R_2-R_1+d(\e )=R_2-R_1+\a_2>\a_1$ we get $\a_1=R_2-R_1+d(-\e
a_{1,2})=R_2-R_1+d(-a'_{1,2})\geq\(\a_1$ so $\(\a_1=\a_1$. It follows
that for $i=2,3$ we have $\(\a_i\leq
R_{i+1}-R_2+\(\a_1=R_{i+1}-R_2+\a_1=\a_i$ so $\(\a_i=\a_i$. Thus we
have (ii). If $R_1=R_3$ and $R_2=R_4$ then $R_1+R_2=R_3+R_4$ so
$R_1+\a_1=R_2+\a_2=R_3+\a_3$, by [B3, Corollary 2.3(i)], and similarly
for the $\(\a_i$'s. Thus in order to prove that $\(\a_i=\a_i$ holds
for all $i$'s it is enough to prove that it holds for one value of
$i$. We have $\(\a_3\leq R_4-R_3+d(-a'_{3,4})=R_4-R_3+d(-\e
a_{3,4})=R_4-R_3+\a_2=\a_3$. Hence $\(\a_3=\a_3$ and (ii) is true in
this case too. For (iv) at $i=2$ note that $a_1\rep
[a_1,a'_2]=[a'_1,a'_2]$ and for $i=3$ we have $\a_2+\a_3<2e$. Thus
$M\ap\( M$. We have $\a_3=\(\a_3=\a_2(\[ a'_2,a'_3,a'_4\] )\leq
R_4-R_3+d(-a'_{3,4})=\a_3$ (see above) so $\a_3=\a_2(\[
a'_2,a'_3,a'_4\] )$ and we are done. 

We consider now the remaining cases. We have $R_1=R_3$ and
$R_2<R_4$. Similarly as for the cases when $R_1<R_3$ or $R_2=R_4$ or
$R_2-R_1=2e$ it is enough to prove our statement when $n=5$. We have
to prove that $M\ap\[ a'_1,a'_2,a'_3,a'_4,a'_5\]$ relative to some
good BONG with $a'_1=a_1$ and $\a_4=\a_3(\[ a'_2,a'_3,a'_4,a'_5\]
)$. By Corollary 8.10 we may assume that $\a_1=\a_1(\[ a_1,a_2\]
)$. We have $\a_4=\min\{\a_3(\[ a_2,a_3,a_4,a_5\]
),R_5-R_2+\a_1\}$. We may assume that $\a_4=R_5-R_2+\a_1<\a_3(\[
a_2,a_3,a_4,a_5\] )$ since otherwise we can take $a'_i=a_i$. Since
$-R_2+\a_1=-R_5+\a_4$ we have
$-R_2+\a_1=-R_3+\a_2=-R_4+\a_3=-R_5+\a_4$. Also
$R_5-R_4+d(-a_{4,5})\geq\a_3(\[a_2,a_3,a_4,a_5\] )>\a_4$ so
$d(-a_{4,5})>R_4-R_5+\a_4=\a_3$. We have $\a_4<\a_3(\[
a_2,a_3,a_4,a_5\] )\leq (R_5-R_4)/2+e$. Since $\a_4<(R_5-R_4)/2+e$ and
$-R_2+\a_1=-R_5+\a_4$ we also have by Lemma 8.4(iii) that
$\a_i<(R_{i+1}-R_i)/2+e$ for $i=1,2,3$. Hence $\a_1,\a_2,\a_3,\a_4$
are odd and $<2e$ so they belong to $d(\ooo )$. Also note that
$\a_3+\a_4=R_4-R_5+2\a_4<2e$. 

Let $\e\in\ooo$ with $d(\e )=\a_3$. We have $d(-a_{4,5})>\a_3=d(\e )$
so $d(-\e a_{4,5})=\a_3$. Since $\a_4\in d(\ooo )$ and $\a_3+\a_4<2e$,
by Lemma 8.2, for any choice of the $\pm$ sign there is $\eta\in\ooo$
with $d(\eta )=\a_4$ and $(-\e a_{4,5},\eta )_\p =\pm 1$. We choose
such $\eta$ with $(-\e a_{4,5},\eta )_\p =(\e,-a_{3,4})_\p$. Let
$a'_1=a_1$, $a'_2=a_2$, $a'_3=\e a_3$, $a'_4=\e\eta a_4$ and
$a'_5=\eta a_5$. Let $\( M=\[a'_1,\ldots,a'_5\]$. We claim that $\(
M\ap M$. We use the classification theorem. We have $\det
[a_3,a_4,a_5]=\det [a'_3,a'_4,a'_5]=a_{3,5}$. Also
$(-a'_{3,4},-a'_{4,5})_\p =(-\eta a_{3,4},-\e a_{4,5})_\p
=(-a_{3,4},-a_{4,5})_\p (\eta,-\e a_{4,5})_\p (-a_{3,4},\e )_\p
=(-a_{3,4},-a_{4,5})_\p$. Thus $[a'_3,a'_4,a'_5]\ap
[a_3,a_4,a_5]$. Since also $a'_1=a_1$ and $a'_2=a_2$ we have $F\( M\ap
FM$. Condition (i) of the classification theorem is obvious. Note that
$d(a_1a'_1)=d(a_{1,2}a'_{1,2})=\j$, $d(a_{1,3}a'_{1,3})=d(\e )=\a_3$
and $d(a_{1,4}a'_{1,4})=d(\eta )=\a_4$ so (iii) holds. Since also
$a_{1,5}a'_{1,5}\in\fs$ by Lemma 8.6 we have that $\( M$ exists and
$\(\a_i:=\a_i(\( M)\geq\a_i$. Now $\(\a_1\leq\a_1(\[ a_1,a_2\] )=\a_1$
so $\(\a_1=\a_1$. For $i=2,3,4$ we have $\(\a_i\leq
R_{i+1}-R_2+\(\a_1=R_{i+1}-R_2+\a_1=\a_i$ so $\(\a_i=\a_i$. Thus we
also have (ii). For (iv) at $i=2,3$ we note that $[a_1]\rep
[a_1,a_2]$, $[a_1,a_2]\rep [a_1,a_2,a'_3]$ and for $i=4$ we have
$\a_3+\a_4<2e$ so (iv) is vacuous. So $M\ap\( M$. Also
$\a_4=\(\a_4\leq \a_3(\[ a'_2,a'_3,a'_4,a'_5\] )\leq
R_5-R_4+d(-a'_{4,5})=R_5-R_4+d(-\e a_{4,5})=R_5-R_4+\a_3=\a_4$ so we
are done. \qed 

\blm Suppose that $N\leq M$ and $R_1=S_1$. Suppose that either
$M,N$ satisfy the conditions of Lemma 9.1 or $S_2=R_2+1$ and
$R_1=R_5$. Then there is a choice of the good BONGs of $M,N$ s.t.
$a_1=b_1$ and if $M^*=\[ a_2,\ldots,a_n\]$ and $N^*=\[
b_2,\ldots,b_n\]$ then $A_i=A_{i-1}(M^*,N^*)$ for all indices $2\leq
i<n$ s.t. $i-1$ or $i$ is essential for $M^*,N^*$. For such a choice
of BONGs the lattices we have $N^*\leq M^*$. 
\elm 
\pf We have $R_1+R_2\leq S_1+S_2=R_1+S_2$ so $R_2\leq S_2$. By Lemma
8.12(i) we have $A_1=\a_1$. In particular, $\b_1\geq A_1=\a_1$. 

With the exception of the case $S_2=R_2+1$ and $R_1=R_5$ we can use
Lemma 9.1 and we have $a_1=b_1$ for a good choice of the BONG of $M$
regardless of the choice of the BONG of $N$.

If $a_1=b_1$ then for $2\leq i<n$ we define $\a^*_i:=\a_{i-1}(M^*)$,
$\b^*_i:=\b_{i-1}(N^*)$ and $A^*_i:=A_{i-1}(M^*,N^*)$. If $1\leq
i,j\leq n$ and $\e\in\ff$ we define $d[\e a_{2,i}b_{2,j}]:=\min\{ d(\e
a_{2,i}b_{2,j}),\a^*_i,\b^*_j\}$. (We ignore $\a^*_i$ if $i=1$ or $n$
and similarly for $\b^*_j$.) Also define $d^*[\e a_{i,j}]:=\min\{ d(\e
a_{i,j}),\a^*_{i-1},\a^*_j\}$ and similarly for $d^*[\e
b_{i,j}]$. Note that $d[\e a_{2,i}b_{2,j}]$, $d^*[\e a_{i,j}]$ and
$d^*[\e b_{i,j}]$ are the $d[\e a_{1,i-1}b_{1,j-1}]$, $d[\e
a_{i-1,j-1}]$ and $d[\e b_{i-1,j-1}]$ corresponding to $M^*,N^*$. Also
$R_i=R_{i-1}(M^*)$ and $S_i=R_{i-1}(N)$. 

For $i\geq 2$ we have $\a_i=\min\{\a^*_i,R_{i+1}-R_1+d(-a_{1,2})\}
=\min\{\a^*_i,R_{i+1}-R_2+\a_1\}\leq\a^*_i$. Similarly for $\b_i$.

Suppose that $a_1=b_1$. We have $a_{2,i}b_{2,j}=a_{1,i}b_{1,j}$ in
$\ff/\fs$. Thus $d[\e a_{2,i}b_{2,j}]=\min\{ d(\e
a_{1,i}b_{1,j}),\a^*_i,\b^*_j\}$. It follows that $d[\e
  a_{1,i}b_{1,j}]\leq d[\e a_{2,i}b_{2,j}]$ with equality unless $d[\e
a_{1,i}b_{1,j}]=\a_i<\a^*_i$ or $d[\e
  a_{1,i}b_{1,j}]=\b_j<\b^*_j$. But $\a_i<\a^*_i$ implies
$\a_i=R_{i+1}-R_2+\a_1=R_{i+1}-R_1+d(-a_{1,2})$. Similarly for
$\b_j<\b^*_j$. 

Note that for any $2\leq i<n$ we have $A^*_i=A_{i-1}(M^*,N^*)=\min\{
(R_{i+1}-S_i)/2+e, R_{i+1}-S_i+d[-a_{2,i+1}b_{2,i-1}],
R_{i+1}+R_{i+2}-S_{i-1}-S_i+d[-a_{2,i+2}b_{2,i-2}]\}$ so $A^*_i\geq
A_i$.

Suppose that the BONGs of $M,N$ satisfy the first part of the
conclusion, i.e. $A_i=A^*_i$ whenever $i-1$ or $i$ is essential for
$M^*,N^*$. Conditions (i)-(iv) of Theorem 2.1 for $M^*,N^*$ are:

(i*) For $2\leq i\leq n$ we have either $R_i\leq S_i$ or $2<i<n$ and
$R_i+R_{i+1}\leq S_{i-1}+S_i$.

(ii*) $d[a_{2,i}b_{2,i}]\geq A^*_i$ for any $2\leq i<n$.

(iii*) $[b_2,\ldots,b_{i-1}]\rep [a_2,\ldots,a_i]$ for any $2<i<n$
s.t. $R_{i+1}>S_{i-1}$ and $A^*_{i-1}+A^*_i>2e+R_i-S_i$.

(iv*) $[b_2,\ldots,b_{i-1}]\rep [a_2,\ldots,a_{i+1}]$ for any
$2<i<n-1$ s.t. $S_i\geq R_{i+2}>S_{i-1}+2e\geq R_{i+1}+2e$.

We note that (i*) for $i>2$ follows from the corresponding condition
for $M,N$. At $i=2$ we have $R_2\leq S_2$.

For (ii*) we have $d[a_{2,i}b_{2,i}]\geq d[a_{1,i}b_{1,i}]\geq
A_i=A^*_i$ whenever $i-1$ or $i$ is essential for $M^*$ and $N^*$. If
neither $i-1$ nor $i$ is essential then the condition
$d[a_{2,i}b_{2,i}]\geq A^*_i=A_{i-1}(M^*,N^*)$ is vacuous by Lemma
2.13. 

For (iii*) suppose that $R_{i+1}>S_{i-1}$ and
$A^*_{i-1}+A^*_i>2e+R_i-S_i$. If $i-1$ is essential for $M^*,N^*$ then
$A^*_{i-1}=A_{i-1}$ and $A^*_i=A_i$. Thus $A_{i-1}+A_i>2e+R_i-S_i$ so
$[b_1,\ldots,b_{i-1}]\rep [a_1,\ldots,a_i]$. Since $a_1=b_1$ this implies
$[b_2,\ldots,b_{i-1}]\rep [a_2,\ldots,a_i]$. If $i-1$ is not essential for
$M^*,N^*$ then (iii*) is vacuous by Lemma 2.13.

Finally, since $a_1=b_1$ we have that $[b_2,\ldots,b_{i-1}]\rep
[a_2,\ldots,a_{i+1}]$ is equivalent to\\ $[b_1,\ldots,b_{i-1}]\rep
[a_1,\ldots,a_{i+1}]$. Thus (iv*) is an obvious consequence of (iv).

We prove now the existence of good BONGs for $M,N$ satisfying the
conditions of our lemma. We say that an index $2\leq i\leq n-1$ is
important if $i-1$ or $i$ is essential for $M^*,N^*$. So we have to
find good BONGs for $M,N$ s.t. $A_i=A^*_i$ at all important indices
$2\leq i\leq n-1$. 

Note that the sequences $R_1,\ldots,R_n$ and $S_1,\ldots,S_n$
corresponding to $M^*$ and $N^*$ are $R_2,\ldots,R_n$ and
$S_2,\ldots,S_n$ so if $i\geq 4$ then $i$ is essential for $M,N$ iff
$i-1$ is essential for $M^*,N^*$. If $i=3$ then $3$ is essential for
$M,N$ if $S_2<R_4$ and $S_1+S_2<R_4+R_5$, while $2$ is essential for
$M^*,N^*$ if $S_2<R_4$. But if $S_2<R_4$, since also $S_1=R_1\leq
R_5$, we have $S_1+S_2<R_4+R_5$. So $3$ is essential for $M,N$ iff $2$
is essential for $M^*,N^*$. If $i=2$ then $2$ is essential for $M,N$
if $S_1<R_3$, while $1$ is essential for $M^*,N^*$ unconditionally.

Suppose first that $M,N$ satisfy the conditions from the hypothesis of
Lemma 9.1. By Lemma 9.1, regardless of the choice of the BONG of $N$,
we can choose a BONG of $M$ s.t. $a_1=b_1$. For $i=2$ we have by Lemma
8.12 $A_2=\min\{ (R_3-S_2)/2+e, R_3-S_2+d[-a_{1,3}b_1]\}$ and we also
have $A^*_2=\min\{ (R_3-S_2)/2+e, R_3-S_2+d^*[-a_{2,3}]\}$. If $n\geq
3$ then in the term
$R_{i+1}+R_{i+2}-S_{i-1}-S_i+d[a_{1,i+2}b_{1,i-2}]$ we can replace
$d[a_{1,i+2}b_{1,i-2}]$ by $d[-a_{1,i}b_{1,i-2}]$ if $S_{i-1}<R_{i+1}$
and by $d[-a_{1,i+2}b_{1,i}]$ if $S_i<R_{i+2}$ (see Lemma
2.7). Similarly for $d[a_{2,i+2}b_{2,i-2}]$ in $A^*_i$. If both
$S_{i-1}\geq R_{i+1}$ and $S_i\geq R_{i+2}$ then $i$ and $i+1$ are not
essential for $M,N$ and $i-1$ and $i$ are not essential for
$M^*,N^*$. Thus it is enough to have
$d[-a_{1,i+1}b_{1,i-1}]=d[-a_{2,i+1}b_{2,i-1}]$ for all $2\leq i\leq
n-1$. (At $i=2$ this means
$d[-a_{1,3}b_1]=d[-a_{2,3}b_{2,1}]=d^*[-a_{2,3}]$.) 

We have two cases:

1. $\b_1>d[-a_{1,3}b_1]$, $\b_1\geq\a_3$, or $d[-a_{1,3}b_1]\geq
(S_2-R_3)/2+e$. By Corollary 8.11 we can change the BONG of $\[
b_1,b_2,b_3\]$ s.t. $\a_2(\[ b_1,b_2,b_3\] )=\a_1(\[ b_2,b_3\] )$. For
any $i\geq 2$ we have $\b_i=\min\{\a_{i-1}(\[ b_2,\ldots,b_n\] ),
R_{i+1}-R_3+\a_2(\[ b_1,b_2,b_3\] )\}$. But $R_{i+1}-R_3+\a_2(\[
b_1,b_2,b_3\] )=R_{i+1}-R_3+\a_1(\[ b_2,b_3\] )\geq\a_{i-1}(\[
b_2,\ldots,b_n\] )$ and so $\b_i=\a_{i-1}(\[ b_2,\ldots,b_n\]
)=\b^*_i$. By Lemma 9.1 we can choose some good BONG for $M$
s.t. $a_1=b_1$ and, by Lemma 9.2, s.t. also $\a_i=\a_{i-1}(\[
a_2,\ldots,a_n\] )=\a^*_i$ for $i\geq 4$. If moreover $R_1<R_3$ or
$R_2=R_4$ or $R_2-R_1=2e$ then $\a_i=\a^*_i$ holds for $i\geq 3$. For
$i\geq 3$ we have $\b_{i-1}=\b^*_{i-1}$ and $\a_{i+1}=\a^*_{i+1}$ so
$d[-a_{1,i+1}b_{1,i-1}]=d[-a_{2,i+1}b_{2,i-1}]$. Thus $A_i=A^*_i$ for
all $A_i$'s with $i\geq 2$ important where $d[-a_{1,3}b_1]$ is not
involved. Now $d[-a_{1,3}b_1]$ is involved in $A_2$ and, if $R_4>S_2$,
also in $A_3$. If $d[-a_{1,3}b_1]\geq (S_2-R_3)/2+e$ then
$R_3-S_2+d[-a_{1,3}b_1]\geq (R_3-S_2)/2+e$ so $A_2=(R_3-S_2)/2+e\geq
A^*_2$ so $A_2=A^*_2$. Also if $R_4>S_2$ then
$R_4+R_5-S_2-S_3+d[-a_{1,3}b_1]\geq
R_4+R_5-S_3-(R_3+S_2)/2+e>R_4+R_5-S_3-(R_5+R_4)/2+e=(R_4+R_5)/2-S_3+e\geq
A^*_3$. (We have $R_3\leq R_5$ and $S_2<R_4$.) Since also
$(R_4-S_3)/2+e\geq A^*_3$ and
$R_4-S_3+d[-a_{1,4}b_{1,2}]=R_4-S_3+d[-a_{2,4}b_2]\geq A^*_3$ we have
$A_3\geq A^*_3$ i.e. $A_3=A^*_3$. So $A_i=A^*_i$ for all important
$i$'s. Thus we may assume that $d[-a_{1,3}b_1]<(S_2-R_3)/2+e$ so
$R_3-S_2+d[-a_{1,3}b_1]<(R_3-S_2)/2+e$, which implies
$A_2=R_3-S_2+d[-a_{1,3}b_1]$. We are left with the cases
$\b_1>d[-a_{1,3}b_1]$ and  $\b_1\geq\a_3$. We have
$d[-a_{1,3}b_1]=\min\{ d^*[-a_{2,3}], \a_3,\b_1\}$. If
$d[-a_{1,3}b_1]=d^*[-a_{2,3}]$ then $A_i=A^*_i$ for all important
$i$'s so we are done. Thus we may assume that
$d[-a_{1,3}b_1]<d^*[-a_{2,3}]$. In both cases $\b_1>d[-a_{1,3}b_1]$
and $\b_1\geq\a_3$ this implies
$d[-a_{1,3}b_1]=\a_3<d^*[-a_{2,3}]\leq\a^*_3$. It follows that
$R_1=R_3$, $R_2<R_4$ and $R_2-R_1\neq 2e$ (otherwise
$\a_3=\a^*_3$). So we are left with the cases $R_2=S_2$ and
$d[-a_{1,3}b_1]=\a_1<\b_1$. Since $\a_3<\a^*_3$ we have
$\a_3=R_4-R_2+\a_1>\a_1$. Since  $d[-a_{1,3}b_1]=\a_3$ this rules out
the case $d[-a_{1,3}b_1]=\a_1<\b_1$. So we are left with the case
$R_2=S_2$. We have $\a_2\geq
A_2=R_3-S_2+d[-a_{1,3}b_1]=R_3-R_2+\a_3\geq\a_2$ so
$\a_2=R_3-R_2+\a_3$. On the other hand $-R_2+\a_1=-R_4+\a_3$ so
$-R_2+\a_1=-R_3+\a_2=-R_4+\a_3$. Thus
$R_3-R_2+\a_3=\a_2=R_3-R_4+\a_3$, i.e. $R_2=R_4$. Contradiction. 

2. $\b_1=d[-a_{1,3}b_1]<(S_2-R_3)/2+e$ and, if $n\geq 4$,
$\b_1<\a_3$. Since $S_1=R_1\leq R_3$ we have $\b_1<(S_2-R_3)/2+e\leq
(S_2-S_1)/2+e$. Note that if $R_1=R_3$ and $R_2=R_4$ then
$R_1+R_2=R_3+R_4$ so by [B3, Corollary 2.3(i)]
$R_1+\a_1=R_3+\a_3=R_1+\a_3$ and so
$\b_1\geq\a_1=\a_3$. Contradiction. We cannot have
$d[-a_{1,3}b_1]=\a_1<\b_1$ either. If $R_2-R_1=2e$ since also
$S_1=R_1$ and $S_2\geq R_2$ we have $S_2-S_1\geq 2e$ so
$\b_1=(S_2-S_1)/2+e$. Contradiction. So we are left with the case when
either $R_1<R_3$ or $R_2=S_2$ and, if $n\geq 3$, $R_2<R_4$. 

Now $d(-a_{1,3}b_1)\geq d[-a_{1,3}b_1]=\b_1$. If $d(-a_{1,3}b_1)>\b_1$
then, since $\b_1<(S_2-S_1)/2+e$, we have by Lemma 8.8 $N\ap\[
b'_1,\ldots,b'_n\]$ relative to some good BONG s.t. $b'_1=\e b_1$ with
$d(\e )=\b_1$. We have $d(\e )=\b_1<d(-a_{1,3}b_1)$ so
$d(-a_{1,3}b'_1)=d(-\e a_{1,3}b_1)=\b_1$. Thus for a good choice of
the BONG of $N$ we have $d(-a_{1,3}b_1)=\b_1$. Also by Lemma 9.2 we
can choose the BONG s.t. $\b_i=\a_{i-1}(\[ b_2,\ldots,b_n\] )=\b^*_i$
for $i\geq 4$. (Note that this change of BONGs preserves $b_1$ and so
the property that $d(-a_{1,3}b_1)=\b_1$.) If $S_1<S_3$ or $S_2=S_4$
then $\b_i=\b^*_i$ holds for $i\geq 3$. Note that
$d(-a_{1,3}b_1)=\b_1$ holds  regardless of the BONG of $M$. Indeed, if
$M\ap\[ a'_1,\ldots,a'_n\]$ then
$d(a_{1,3}a'_{1,3})\geq\a_3>\b_1=d(-a_{1,3}b_1)$ and so
$d(-a'_{1,3}b_1)=\b_1$. Like in the case 1. we can take the BONG of
$M$ s.t. $a_1=b_1$ and $\a_i=\a^*_i$ for $i\geq 4$. We have
$d^*[-a_{2,3}]\leq d(-a_{2,3})=d(-a_{1,3}b_1)=\b_1=d[-a_{1,3}b_1]$ so
$d[-a_{1,3}b_1]=d^*[-a_{2,3}]$. For $i\geq 3$ we have
$\a_{i+1}=\a^*_{i+1}$. Since $d[-a_{1,i+1}b_{1,i-1}]=\min\{
d(-a_{1,i+1},b_{1,i-1}),\b_{i-1},\a_{i+1}\}$ and
$d[-a_{2,i+1}b_{2,i-1}]=\min\{
d(-a_{1,i+1},b_{1,i-1}),\b^*_{i-1},\a^*_{i+1}\}$ we have
$d[-a_{1,i+1}b_{1,i-1}]=d[-a_{2,i+1}b_{2,i-1}]$ unless
$d[-a_{1,i+1}b_{1,i-1}]=\b_{i-1}<\b^*_{i-1}$. Also for $i\geq 5$
(resp. for $i\geq 4$ if $S_1<S_3$ or $S_2=S_4$) we have
$\b_{i-1}=\b^*_{i-1}$ so
$d[-a_{1,i+1}b_{1,i-1}]=d[-a_{2,i+1}b_{2,i-1}]$. The same relation
holds for $i=3$ or $i=4$ if $\b_2=\b^*_2$ resp. $\b_3=\b^*_3$. So we
may assume that either $\b_2<\b^*_2$ or $\b_3<\b^*_3$. If
$\b_2<\b^*_2$ then $\b_2=S_3-S_2+\b_1$. If $\b_3<\b^*_3$ then
$\b_3=S_4-S_2+\b_1$, i.e. $-S_2+\b_1=-S_4+\b_3$, which implies
$-S_2+\b_1=-S_3+\b_2=-S_4+\b_3$. So again $\b_2=S_3-S_2+\b_1$.

We have $d(-b_{2,3})\geq S_2-S_3+\b_2=\b_1$. Suppose that
$d(-b_{2,3})=\b_1$. Then $\b^*_2\leq
S_3-S_2+d(-b_{2,3})=S_3-S_2+\b_1=\b_2$ so $\b^*_2=\b_2$. Also if
$\b^*_3>\b_3$ then $\b_3=S_4-S_2+\b_1$ and so 
$\b^*_3\leq S_4-S_2+d(-b_{2,3})=S_4-S_2+\b_1=\b_3$. So $\b^*_2=\b_2$
and $\b_3^*=\b_3$, contrary to our assumption. Thus
$d(-b_{2,3})>\b_1=d(-a_{1,3}b_1)$, which implies that
$\b_1=d(a_{1,3}b_{1,3})\geq A_3$.

We have $d[-a_{1,3}b_1]<(S_2-R_3)/2+e$ so
$R_3-S_2+d[-a_{1,3}b_1]<(R_3-S_2)/2+e$. By Lemma 8.12(ii)
$A_2=R_3-S_2+d[-a_{1,3}b_1]=R_3-S_2+\b_1$. We have
$S_3-S_2+\b_1=\b_2\geq A_2=R_3-S_2+\b_1$. So $R_3\leq S_3$. 


We claim that if $A_3<R_4-S_3+d[-a_{1,4}b_{1,2}]$, $R_5>S_3$ and
$A_4=R_5+R_6-S_2-S_3+d[-a_{1,4}b_{1,2}]$ then $2,3,4$ are non
essential indices for $M^*,N^*$. We have $R_5>S_3$ so $A_3=\min\{
(R_4-S_3)/2+e, R_4-S_3+d[-a_{1,4}b_{1,2}],
R_4+R_5-S_2-S_3+d[-a_{1,5}b_{1,3}]\}$. If
$A_3=(R_4-S_3)/2+e<R_4-S_3+d[-a_{1,4}b_{1,2}]$ then
$d[-a_{1,4}b_{1,2}]>(S_3-R_4)/2+e$ so
$R_5+R_6-S_3-S_4+d[-a_{1,4}b_{1,2}]>R_5+R_6-(R_4+S_3)/2-S_4+e>
R_5+R_6-(R_6+R_5)/2-S_4+e=(R_5+R_6)/2-S_4+e\geq
A_4$. Contradiction. If
$A_3=R_4+R_5-S_2-S_3+d[-a_{1,5}b_{1,3}]<R_4-S_3+d[-a_{1,4}b_{1,2}]$
then $R_5-S_2+d[-a_{1,5}b_{1,3}]<d[-a_{1,4}b_{1,2}]$. On the other
hand $R_5+R_6-S_3-S_4+d[-a_{1,4}b_{1,2}]=A_4\leq
R_5-S_4+d[-a_{1,5}b_{1,3}]$ and so $d[-a_{1,5}b_{1,3}]\geq
R_6-S_3+d[-a_{1,4}b_{1,2}]>R_5+R_6-S_2-S_3+d[-a_{1,5}b_{1,3}]$. Thus
$R_5+R_6<S_2+S_3$ and so $3$ is not essential for $M^*,N^*$. Since
also $R_5>S_3$ we get $R_6<S_2$, which implies both $R_4<S_2$ and
$R_6<S_4$. Thus $2$ and $4$ are also not essential for $M^*,N^*$.

Similarly if $A_4<R_5-S_4+d[-a_{1,5}b_{1,3}]$, $R_6>S_4$ and
$A_5=R_6+R_7-S_3-S_4+d[-a_{1,5}b_{1,3}]$ then $3,4,5$ are non
essential indices for $M^*,N^*$.

We have 2 subcases: 

a. $d[-a_{1,5}b_{1,3}]=d[-a_{2,5}b_{2,3}]$. Then we may assume
$d[-a_{1,4}b_{1,2}]<d[-a_{2,4}b_2]$ since otherwise
$d[-a_{1,i+1}b_{1,i-1}]=d[-a_{2,i+1}b_{2,i-1}]$ holds for all $i\geq
2$ so $A_i=A^*_i$ for all important  $i\geq 2$. Thus
$d[-a_{1,4}b_{1,2}]=\b_2<d[-a_{2,4}b_2]$. So $\b_2<\b^*_2$, which
implies $\b_2=S_3-S_2+\b_1$. Now $d[-a_{1,4}b_{1,2}]$ is involved only
in $A_3$ and, if $S_3<R_5$, in $A_4$. Thus $A_i=A^*_i$ for all
important $i$, $i\neq 3,4$. We also have $A_3=A^*_3$ unless
$A_3=R_4-S_3+d[-a_{1,4}b_{1,2}]$ and $A_4=A^*_4$ unless $R_5>S_3$ and
$A_4=R_5+R_6-S_3-S_4+d[-a_{1,4}b_{1,2}]$. If
$A_3<R_4-S_3+d[-a_{1,4}b_{1,2}]$ then $A_3=A^*_3$ so we may assume
$A_4<A^*_4$ since otherwise $A_i=A^*_i$ for all important $i\geq
2$. Thus $R_5>S_3$ and $A_4=R_5+R_6-S_3-S_4+d[-a_{1,4}b_{1,2}]$. As
seen above, this implies that $2,3,4$ are not essential for
$M^*,N^*$. Hence the equality $A_4=A^*_4$ is not required and we are
done. So we may assume that
$A_3=R_4-S_3+d[-a_{1,4}b_{1,2}]=R_4-S_3+\b_2=R_4-S_2+\b_1$. Now
$\b_1\geq A_3=R_4-S_2+\b_1$ so $R_4>S_2$. So we cannot have $R_2=S_2$
and $R_2<R_4$. So we are left with the case when $S_1=R_1<R_3$. Now
$\b_1<(S_2-S_1)/2+e$ so $\b_1=S_2-S_1+d[-b_{1,2}]$. Hence
$d[-b_{1,2}]=S_1-S_2+\b_1<R_3-S_2+\b_1=A_2\leq d[a_{1,2}b_{1,2}]$. It
follows that $R_1-R_2+\a_1\leq
d[-a_{1,2}]=d[-b_{1,2}]=S_1-S_2+\b_1=R_1-S_2+\b_1$ and so
$-R_2+\a_1\leq -S_2+\b_1$. It follows that $A_3=R_4-S_2+\b_1\geq
R_4-R_2+\a_1\geq\a_3$. But $\a_3>\b_1\geq A_3$. Contradiction. 

b. $d[-a_{1,5}b_{1,3}]<d[-a_{2,5}b_{2,3}]$. Hence
$d[-a_{1,5}b_{1,3}]=\b_3<\b^*_3$. This implies that $S_1=S_3$ and
$S_2<S_4$ and $-S_2+\b_1=-S_3+\b_2=-S_4+\b_3$. Since $R_3\leq
S_3=R_1$ we have $R_1=R_3$ so we are left with the case
$S_2=R_2<R_4$. We have $R_4>S_2$ and $R_5\geq R_1=S_1=S_3$ so
$R_4+R_5>S_2+S_3$. By Lemma 2.9 $A_3=\( A_3=\min\{
(R_4-S_3)/2+e, R_4-S_3+d[-a_{1,4}b_{1,2}]\}$. But 
$(R_4-S_3)/2+e>(S_2-R_3)/2+e>\b_1\geq A_3$ so $\b_1\geq
A_3=R_4-S_3+d[-a_{1,4}b_{1,2}]$ and so $d[-a_{1,4}b_{1,2}]\leq
S_3-R_4+\b_1<S_3-S_2+\b_1=\b_2$. It follows that
$d[-a_{1,4}b_{1,2}]=d[-a_{2,4}b_2]$. Thus
$d[-a_{1,i+1}b_{1,i-1}]=d[-a_{2,i+1}b_{2,i-1}]$ for all $i\neq 4$.
As seen above $d[-a_{1,5}b_{1,3}]$ is not involved in $A_3$. It is
only involved in $A_4$ and, if $R_6>S_4$, in $A_5$. Thus
$A_i=A^*_i$ for all important $i$, $i\neq 4,5$. Also $A_4=A^*_4$
unless $A_4=R_5-S_4+d[-a_{1,5}b_{1,3}]$ and $A_5=A^*_5$ unless
$R_6>S_4$ and $A_5=R_6+R_7-S_4-S_5+d[-a_{1,5}b_{1,3}]$. If
$A_4<R_5-S_4+d[-a_{1,5}b_{1,3}]$ then $A_4=A^*_4$ so we may
suppose that $A_5<A^*_5$ and so $R_6>S_4$ and
$A_5=R_6+R_7-S_4-S_5+d[-a_{1,5}b_{1,3}]$. (Otherwise $A_i=A^*_i$
for all important $i\geq 2$.) It follows that $3,4,5$ are not
essential for $M^*,N^*$. But then condition $A_5=A^*_5$ is not
necessary and we are done. Thus we may assume that
$A_4=R_5-S_4+d[-a_{1,5}b_{1,3}]=R_5-S_4+\b_3=R_5-S_2+\b_1$. Hence
$d[a_{1,4}b_{1,4}]\geq R_5-S_2+\b_1>S_3-R_4+\b_1\geq
d[-a_{1,4}b_{1,2}]$. (We have $R_5\geq R_1=S_3$ and
$R_4>S_2$.) Thus $S_3-R_4+\b_1\geq
d[-a_{1,4}b_{1,2}]=d[-b_{3,4}]\geq S_3-S_4+\b_3=S_3-S_2+\b_1$ and
so $S_2\geq R_4$. Contradiction.

Finally, consider the case $S_2=R_2+1$ and $R_1=R_5$. We have
$R_1=R_3=R_5$ and, by Lemma 6.6(i), $R_1\ev\ldots\ev R_5\m2$. Also
$R_2-R_1$ and $R_4-R_3$ are $\leq 2e$. Since the case $R_2=R_4$ was
already treated we may assume that $R_2<R_4$. But $R_2\ev R_4\m2$ so
$R_4\geq R_3+2$. Now $R_2-R_1$ is even and $\leq 2e$. Since the case
$R_2-R_1=2e$ was already treated we may assume that $R_2-R_1\leq
2e-2$. Thus $S_2-S_1=R_2-R_1+1\leq 2e-1$ and it is odd and so
$\b_1=S_2-S_1=R_2-R_1+1\neq (S_2-S_1)/2+e$. Since $R_2-R_1$ is even
and $<2e$ we have $\a_1\geq R_2-R_1+1=\b_1$. But also $\a_1\leq\b_1$
and so $\a_1=\b_1=R_2-R_1+1$. Since $R_4-R_3\leq 2e$ we have $\a_3\geq
R_4-R_3\geq R_2+2-R_1>\a_1=\b_1$. Since $R_1=R_3$ we have
$\a_2=R_1-R_2+\a_1=1$ and $d[-a_{2,3}]=\a_1$. (See 8.7.) Since also
$d[a_1b_1]\geq A_1=\a_1$ we have $\b_1\geq
d[-a_{1,3}b_1]\geq\a_1=\b_1$. Thus
$d[-a_{1,3}b_1]=\a_1=\b_1=R_2-R_1+1$. We have $d(-a_{1,3}b_1)\geq
d[-a_{1,3}b_1]=\b_1$. If the inequality is strict, since
$\b_1<(S_2-S_1)/2+e$ we have by Lemma 8.8  $N\ap\[ b'_1,\ldots,b'_n\]$
relative to some other good BONG s.t. $b'_1=\e b_1$ with $d(\e
)=\b_1>d(-a_{1,3}b_1)$. It follows that $d(-a_{1,3}b'_1)=d(-\e
a_{1,3}b_1)=d(\e )=\b_1$. Thus by a change of the BONG of $N$ we
may assume that $d(-a_{1,3}b_1)=d[-a_{1,3}b_1]=\b_1$. Note that this
equality holds regarless of the BONG of $M$. Indeed, if $M\ap\[
a'_1,\ldots,a'_n\]$ relative to another BONG then
$d(a_{1,3}a'_{1,3})\geq\a_3>\b_1=d(a_{1,3}b_1)$ so
$d(a'_{1,3}b_1)=\b_1$. 

We cannot have $S_1=S_3$ since this would imply that $S_2-S_1$ is
even. Since $S_1<S_3$ we can use Lemma 9.2 and we can change
$b_2,\ldots,b_n$ s.t. $\b^*_i=\b_i$ for $i\geq 3$. We use the
notation in the proof of Lemma 9.1. Let $\( d[-a_{1,3},b_1]:=\min\{
d(-a_{1,3}b_1),\a_3\}$, i.e. the $d[-a_{1,3},b_1]$ corresponding to
$M$ and $\[ b_1\]$. Then $d[-a_{1,3},b_1]\leq\( d[-a_{1,3},b_1]\leq
d(-a_{1,3},b_1)$. But $d(-a_{1,3},b_1)=d[-a_{1,3},b_1]$ so $\(
d[-a_{1,3},b_1]=d[-a_{1,3},b_1]=\a_1$. We have $R_1=R_3$ so $\a_2+\(
d[-a_{1,3},b_1]=\a_1+\a_2\leq 2e$, with equality iff
$\a_1=(R_2-R_1)/2+e$ (see 8.7). If $\a_1=(R_2-R_1)/2+e$ then we don't
have (b) of Lemma 8.14 and also $\a_2+\( d[-a_{1,3},b_1]=2e$ so we
don't have (a) either. If $\a_1<(R_2-R_1)/2+e$ then $\a_2+\(
d[-a_{1,3},b_1]<2e$ so neither of (a) and (b) from Lemma 8.14
holds. Since $R_3=R_5$ we have $\a_3+\a_4\leq 2e$ so (c) doesn't hold
either. Since none of (a) - (c) holds we can change the BONG of $M$
s.t. $a_1=b_1$. By Lemma 9.2 we can change $a_2,\ldots,a_n$ s.t.
$\a^*_i=\a_i$ for $i\geq 4$. As seen above, these changes of BONG
preserve the property that $d(a_{1,3}b_1)=\b_1$. If $i\geq 4$ then
$\a^*_{i+1}=\a_{i+1}$ and $\b^*_{i-1}=\b_{i-1}$ so 
$d[-a_{2,i+1}b_{2,i-1}]=d[-a_{1,i+1}b_{1,i-1}]$. If $i=3$ note that
$S_1=R_1\ev R_2\ev R_3\ev R_4\m2$ and $S_2=R_2+1$ so both $\ord
a_{1,4}b_{1,2}=R_1+R_2+R_3+R_4+S_1+S_2$ and $\ord
a_{2,4}b_2=R_2+R_3+R_4+S_2$ are odd. Thus
$d[-a_{2,4}b_2]=d[-a_{1,4}b_{1,2}]=0$. Finally, at $i=2$ we have
$\b_1=d[-a_{1,3}b_1]\leq d^*[-a_{2,3}]\leq
d(-a_{2,3})=d(-a_{1,3}b_1)=\b_1$ so
$d[-a_{1,3}b_1]=d^*[-a_{2,3}]=\b_1$.  Hence
$d[-a_{1,i+1}b_{1,i-1}]=d[-a_{2,i+1}b_{2,i-1}]$ for all $i\geq 2$ and
we are done. \qed

\blm If $M,N$ are ternary over the same quadratic space and
$R_1=R_3=S_1=S_3$ then $N\leq M$ iff $\a_1\leq\b_1$ and
$\a_2\geq\b_2$. 
\elm 
\pf For consequences of $R_1=R_3$ and $S_1=S_3$ see 8.7.

The necessity follows from Proposition 6.2: We have $R_1+\a_1\leq
S_1+\b_1=R_1+\b_1$ so $\a_1\leq\b_1$ and $R_3-\a_2\leq
S_3-\b_2=R_3-\b_2$ so $\a_2\geq\b_2$. 

Conversely, we have $\a_1\leq\b_1$ and
$R_1-R_2+\a_1=\a_2\geq\b_2=S_1-S_2+\b_1$. This implies $R_1-R_2\geq
S_1-S_2=R_1-S_2$ so $R_2\geq S_2$ so Theorem 2.1(i) holds at $i=2$. At
$i=1$ and $3$ we have $R_1=S_1$ and $R_3=S_3$ by hypothesis. 

We have $A_1\leq R_2-S_1+\a_2=R_2-R_1+\a_2=\a_1$ and $A_2\leq
R_3-S_2+\b_1=S_3-S_2+\b_1=\b_2$. Now $d[-a_{1,2}]=\a_2\geq\b_2$ and
$d[-b_{1,2}]=\b_2$ so $d[a_{1,2}b_{1,2}]\geq\b_2\geq A_2$. Also
$d[-a_{2,3}]=\a_1$ and $d[-b_{2,3}]=\b_1\geq\a_1$. Together with
$d[a_{1,3}b_{1,3}]=\j$, these imply $d[a_1b_1]\geq\a_1\geq A_1$. So
(ii) holds. 

Condition (iii) is vacuous (at $i=2$) because $R_3=S_1$ and (iv)
is vacuous because $n=3$. \qed 

\blm If $M$ is ternary over the quadratic space $V$, $R_1=R_3$ and
$a_i=\pi^{R_i}\e_i$ then:

(i) If $V$ is not isotropic then $\a_1$ is odd.

(ii) $M\ap\la -\pi^{R_1}\e_{1,3}\ra\pp\[\pi^{R_1+\a_1},-\pi^{R_1-\a_2}\]$
if $V$ is isotropic and $M\ap\la
-\pi^{R_1}\e_{1,3}\D\ra\pp\[\pi^{R_1+\a_1},-\pi^{R_1-\a_2}\D\]$  if
$V$ is not isotropic. 

(iii) If $\a_i=(R_{i+1}-R_i)/2+e$ then
$\[\pi^{R_1+\a_1},-\pi^{R_1-\a_2}\]\ap\pi^{(R_1+R_2)/2}\aa$. 
\elm
\pf We refer to 8.7 for consequences of $R_1=R_3$. 

(i) If $\a_1$ is even then $\a_1=(R_2-R_1)/2+e$ and so
$\a_1+\a_2=2e$. Now $d(-a_{1,2})\geq\a_2$ and $d(-a_{2,3})\geq\a_1$ so
$d(-a_{1,2})+d(-a_{2,3})\geq\a_1+\a_2=2e$. If the inequality is strict
then $(-a_{1,2},-a_{2,3})_\p =1$ so $V\ap [a_1,a_2,a_3]$ is
isotropic. Contradiction. So $d(-a_{1,2})=\a_2$ and $d(-a_{2,3})=\a_1$
are even with $d(-a_{1,2})+d(-a_{2,3})=2e$, which implies $\{
d(-a_{1,2}),d(-a_{2,3})\} =\{ 0,2e\}$. But this is impossible since
$\ord a_{1,2}=R_1+R_2$ and $\ord a_{2,3}=R_2+R_3$ are even so
$d(-a_{1,2}),d(-a_{2,3})\in d(\ooo )$. So $\a_1$ is odd.

We have $R_1+\a_1=R_2+\a_2$ so $R_1-\a_2=R_2-\a_1$. So in (ii) and
(iii) we can replace $\pi^{R_1-\a_2}$ by $\pi^{R_2-\a_1}$. First we prove that
there is a lattice $J\ap\[\pi^{R_1+\a_1},-\pi^{R_2-\a_1}\]$. To do this we
have to show that
$a:=-\pi^{R_2-\a_1}/\pi^{R_1+\a_1}=-\pi^{R_2-R_1-2\a_1}\in\aaa$. Note
that $\a_1\leq (R_2-R_1)/2+e$ so $R:=\ord a=R_2-R_1-2\a_1\geq
-2e$. Also $R$ is even so $-a\in\fs$ so $d(-a)=\j$. So both $R+2e\geq
0$ and $R+d(-a)\geq 0$ hold and so $a\in\aaa$. Similarly in order that
a lattice $J\ap\[\pi^{R_1+\a_1},-\pi^{R_2-\a_1}\D\]$ exists we need
$-\pi^{R_2-R_1-2\a_1}\D\in\aaa$. Again $R:=\ord a$ is $\geq -2e$
and even but this time $-a\in\D\fs$ and so $d(-a)=2e$. Thus both
$R+2e\geq 0$ and $R+d(-a)\geq 0$ hold and so $a\in\aaa$.

By scaling with $\pi^{-(R_1+R_2)/2}$ we may assume that $R_1+R_2=0$. If
$R:=R_1=R_3$ then $R_2=-R$. Also $R_1+\a_1=R+\a_1$ and
$R_1-\a_2=R_2-\a_1=-R-\a_1$. 

We now prove (iii). Suppose that $\a_1=(R_2-R_1)/2+e=-R+e$. If
$J\ap\[\pi^{R+\a_1},-\pi^{-R-\a_1}\]\ap\[\pi^e,\pi^{-e}\]$ then
$a(J)=-\pi^{-2e}\in -\h 14\ooo^2$ so $J\ap\pi^r\aa$ for some
$r\in\ZZ$. But $\det J=-1$ so $J\ap\aa$. 

For (ii) we have to prove that $M\ap N$, where $N=N_1\pp
N_2$, with $N_1\ap\[\pi^{R+\a_1},-\pi^{-R-\a_1}\]$ and
$N_2\ap\la-\pi^R\e_{1,3}\ra$, if $V$ is isotropic, or
$N_1\ap\[\pi^{R+\a_1},-\pi^{-R-\a_1}\D\]$ and $N_2\ap\la
-\pi^R\e_{1,3}\D\ra$, if $V$ is anisotropic. First we show that
$FM=V\ap FN$. In both 
cases when $V$ is isotropic or not we have $\det
M=a_{1,3}=\pi^{R_1+R_2+R_3}\e_{1,3}=\pi^R\e_{1,3}=\det N$ so we still
have to prove that $V$ is isotropic iff $FN$ is so. If $V$ is
isotropic, in particular if $\a_1$ is even, then
$N_1\ap\[\pi^{R+\a_1},-\pi^{-R-\a_1}\]$ so $FN_1$ is hyperbolic and so
$FN=FN_1\pp FN_2$ is isotropic. If $V$ is anisotropic then $\a_1$ is
odd by (i). We have $FN\ap
[-\pi^R\e_{1,3}\D,\pi^{R+\a_1},-\pi^{-R-\a_1}\D ]\ap
[-\pi^R\e_{1,3}\D,\pi^{R+1},-\pi^{R+1}\D ]$ ($\a_1$ is odd.) Hence
$FN$ is anisotropic. 

We have $R_1=R_3$ so by Lemma 7.11 we still have to prove that
$R_i=S_i$ for every $i$ and $\a_1=\b_1$. We assume first that
$R\geq 0$. 

Since $R+\a_1=R_1+\a_1\geq R_1-\a_2=-R-\a_1$ and
$((R+\a_1)+(-R-\a_1))/2=0$ the lattice $N_1$ is unimodular by [B1,
Corollary 3.4]. If $R=0$ then $N_2$ is also unimodular so $N$ is
unimodular. We have $\nn N=\ss N=\oo$ so $S_1=S_2=S_3=0$ by [B3, Lemma
2.13]. Hence $R_i=S_i$ for all $i$. If $R>0$ then $N=N_1\pp N_2$ is a
Jordan splitting. Let $r_i=\ord\ss N_i$ and $u_i=\ord\ss N^{\ss
L_i}$. We have $r_1=0$ and $r_2=R>0$. Since $N_2$ is unary we have
$u_2=r_2=R$. We have $\nn N_1=\p^{R+\a_1}\sbq\p^R=\nn N_2$ so $\nn
N=\p^R$. Hence $u_1=\ord\nn N=R$. By [B3, Lemma 2.13] we have
$S_1=u_1=R$, $S_2=2r_1-u_1=-R$ and $R_3=u_2=R$. So $R_i=S_i$ for all
$i$.

Now we prove that $\a_1=\b_1$. By [B3, Lemma 2.16(i)], in both cases
when $R=0$ and $R>0$, we have $R+\b_1=S_1+\b_1=\ord\ww N$ so we have
to prove that $\ord\ww N=R+\a_1$. Let $\xi=-\e_{1,3}$ or $-\e_{1,3}\D$
corresponding to $V$ isotropic or not. Then $Q(N_2)=\pi^R\xi\oo^2$. We
have $\nn N_1=\p^{R+\a_1}$ so
$Q(N)=\pi^R\xi\oo^2+Q(N_1)\sbq\pi^R\xi\oo^2+\p^{R+\a_1}$. If
$\a_1=(R_2-R_1)/2+e=-R+e$ then $\p^{R+\a_1}=2\oo=2\ss N$ and so
$\ww N=2\ss N=2\oo$. Thus $\ord\ww N=e=R+\a_1$ and we are done. If
$\a_1<(R_2-R_1)/2+e=-R+e$ then $\a_1$ is odd and we have
$\p^{R+\a_1}\sp 2\oo=2\ss N$. Since also $R+\a_1\ev R+1\m2$ the
relation $Q(N)\sbq\pi^R\xi\oo^2+\p^{R+\a_1}$ implies $\ww
N\sbq\p^{R+\a_1}$. On the other hand $\pi^{R+\a_1}\in Q(N_1)$ so
$\pi^R\xi+\pi^{R+\a_1}\in Q(N)$. Now
$(\pi^R\xi+\pi^{R+\a_1})/(\pi^R\xi )=1+\pi^{\a_1}\xi\1$. Since $\a_1$
is odd and $\a_1<-R+e\leq e<2e$ we have $d( 1+\pi^{\a_1}\xi\1 )=\a_1$,
which implies that $1+\pi^{\a_1}\xi\1\notin\oo^2+\p^{\a_1+1}$. It
follows that $\pi^R\xi+\pi^{R+\a_1}\notin\pi^R\xi\oo^2+\p^{R+\a_1+1}$
so $Q(N)\not\subseteq\pi^R\xi\oo^2+\p^{R+\a_1+1}$. Hence $\ww
N\not\sbq\p^{R+\a_1+1}$ and so $\ww N=\p^{R+\a_1}$, as claimed.

If $R<0$ we use duality. We have
$M^\*\ap\[\pi^{-R}\e_3,\pi^R\e_2,\pi^{-R}\e_1\]$ and $-R>0$. Therefore
we may apply the previous case. Since $\a_1(M^\* )=\a_2$ we get
$M^\*\ap\la -\pi^{-R}\e_{1,3}\ra\pp\[\pi^{-R+\a_2}-\pi^{R-\a_2}\]$ or
$M^\*\ap\la
-\pi^{-R}\e_{1,3}\D\ra\pp\[\pi^{-R+\a_2}-\pi^{R-\a_2}\D\]$. So 
$M^\*\ap\la -\pi^R\e_{1,3}\ra\pp\[\pi^{-R+\a_2}-\pi^{R-\a_2}\]$ or
$M^\*\ap\la
-\pi^R\e_{1,3}\D\ra\pp\[\pi^{-R+\a_2}-\pi^{R-\a_2}\D\]$. (The binary
component of $M$ is unimodular so selfdual.) But
$R+\a_1=R_1+\a_1=R_2+\a_2=-R+\a_2$. Hence the conclusion. \qed 

\blm Suppose that $N\leq M$, $R_1=R_3=S_1$, $R_2-R_1=2e-2$,
$d[-a_{1,3}b_1]\geq 2e$ and $[a_1,a_2,a_3]$ is not isotropic. Then
$M\ap\[ a'_1,\ldots,a'_n\]$ relative to some bad BONG s.t. $a_1=b_1$
and $N^*\leq N^*$, where $N^*\ap\[ b_2,\ldots,b_n\]$ and $M^*\ap\[
a'_2,\ldots,a'_n\]$.
\elm
\pf We have $2e\leq d[-a_{1,3}b_1]=\min\{
d(-a_{1,3}b_1),\a_3,\b_1\}$. Thus $-a_{1,3}b_1\in\fs$ or $\D\fs$ and
$\a_3,\b_1\geq 2e$ so $R_4-R_3,S_2-S_1\geq 2e$, i.e. $R_4\geq
R_3+2e=R_1+2e$ and $S_2\geq S_1+2e=R_1+2e$. We want to reduce to the
case when $-a_{1,3}b_1\in\D\fs$. Suppose that $-a_{1,3}b_1\in\fs$. If
$R_4>R_1+2e=S_1+2e$ then $[b_1]\rep [a_1,a_2,a_3]$ by Lemma 2.19. Same
happens if $n=3$.  But this is impossible since $[a_1,a_2,a_3]$ is not
isotropic and $-a_{1,3}b_1\in\fs$. If $R_4=R_1+2e$ then $R_4-R_3=2e$
so $g(a_4/a_3)=\upo{2e}\ni\D$. Thus $M\ap\[ a_1,a_2,\D a_3,\D
a_4,a_5,\ldots,a_n\]$. By this change of BONG $-a_{1,3}b_1$ is
replaced by $-\D a_{1,3}b_1\in\D\fs$. Also $[a_1,a_2,a_3]$ is replaced
by $[a_1,a_2,\D a_3]$ but $(-a_{1,2},-\D a_{2,3})_\p
=(-a_{1,2},-a_{2,3})_\p (-a_{1,2},\D )_\p =(-a_{1,2},-a_{2,3})_\p$ so
the property that $[a_1,a_2,a_3]$ is not isotropic is preserved by
this change of BONG. ($\ord a_{1,2}$ is even so $(-a_{1,2},\D )_\p
=1$.) 

If $n\geq 4$ then $R_4\geq R_3+2e>R_3$ so $M=M_1\pp M_2$ where $M_1\ap\[
a_1,a_2,a_3\]$ and $M_2\ap\[ a_4,\ldots,a_n\]$. Since $R_2-R_1=2e-2$ and
$R_3-R_2=2-2e$ we have $\a_1=\a_1(M_1)=2e-1$ and $\a_2=\a_2(M_1)=1$. Let
$a_i=\pi^{R_i}\e_i$. Since $FM_1\ap [a_1,a_2,a_3]$ is not isotropic we
have by Lemma 9.5 that $M_1=\oo x\pp J$ with $Q(x)=a'_1$ and $J\ap\[
a'_2,a'_3\]$, where $a'_1=-\pi^{R_1}\e_{1,3}\D$, $a'_2=\pi^{R_1+\a_1}$ and
$a'_3=-\pi^{R_1-\a_2}\D$. Since $\a_1=2e-1$ and $\a_2=1$ we have
$a'_2=\pi^{R_1+2e-1}=\pi^{R_2+1}$ and $a'_3=-\pi^{R_1-1}\D$. For $i\geq 4$ we
define $a'_i=a_i$. If $R'_i=\ord a'_i$ then $R'_2=R_2+1$, $R'_3=R_1-1=R_3-1$
and $R'_i=R_i$ for $i\neq 2,3$. Now $x\in M_1\sbq M$ and $\ord
Q(x_1)=R_1=\ord\nn M$ so $x$ is a norm generator for $M$. We have
$M^*:=pr_{x^\pp}M=J\pp M_2$, where $J=pr_{x^\pp}M_1$. Then $J\ap\[
a'_2,a'_3\]$ and $M_2\ap\[ a'_4,\ldots,a'_n\]$ relative to good
BONGs. Since $R'_2=R_2+1=R_3+2e-1<R_3+2e\leq R_4=R'_4$,
$R_3'=R_3-1<R_5$ and $R_3'=R_3-1<R_3+2e\leq R'_4$ we have $M^*\ap\[
a'_2,\ldots,a'_n\]$ relative to some good BONG. Hence $M\ap\[
a'_1,\ldots,a'_n\]$ relative to a bad BONG. (The BONG is bad because
$R'_3=R_1-1=R'_1-1$.) 

Note that $\ord a'_1=R_1=\ord b_1$ and $a'_1=-\pi^{R_1}\e_{1,3}\D\in\
-\D a_{1,3}\fs$, same as $b_1$. (We have $-a_{1,3}b_1\in\D\fs$.)
Therefore, by multiplying $a'_1$ by the square of a unit, we may
assume that $a'_1=b_1$. Also note that $[a'_1,a'_2,a'_3]\ap
[a_1,a_2,a_3]$ and $a'_i=a_i$ for $i\geq 4$ so $[a'_1,\ldots,a'_i]\ap
[a_1,\ldots,a_i]$ for any $i\geq 3$. In particular $a_{1,i}=a'_{1,i}$
in $\ff/\fs$ for $i\geq 3$. 

We prove now that $N^*\leq M^*$. $FM^*\ap FN^*$ follows from $FM\ap
FN$ and $a'_1=b_1$. (We have $FM\ap [a'_1]\pp FM^*$ and $FN\ap
[b_1]\pp FN^*$.) We prove now the conditions (i)-(iv) of Theorem 2.1.

For 2.1(i) note that $R'_2=R_2+1=R_1+2e-1<S_2$ and if $i\geq 3$ then
either $S_i\geq R_i\geq R'_i$ or $i\leq n-1$ and $S_{i-1}+S_i\geq
R_i+R_{i+1}\geq R'_i+R'_{i+1}$. ($R_i\geq R'_i$ for $i\geq 3$.)

In order to prove 2.1(ii) and (iii) we use notations and techniques
similar those from the proof Lemma 9.3 but with $a_i$ replaced
by $a'_i$. We denote $\a^*_i=\a_{i-1}(M^*)$, $\b^*_i=\a_{i-1}(N^*)$,
$A^*_i=A_{i-1}(M^*,N^*)$ and $d[\e a'_{2,i}b_{2,j}]=\min\{ d(\e
a'_{2,i}b_{2,j}),\a^*_i,\b^*_j\}$. 

We claim that $\a^*_i=\a_i$ for $i\geq 4$ and $\b^*_i=\b_i$ for $i\geq
2$. Indeed, for $i\geq 2$ we have $\b_i=\min\{\a_{i-1}(\[
b_2,\ldots,b_n\] ),S_{i+1}-S_i+\b_1\}
=\min\{\b^*_i,S_{i+1}-S_i+\b_1\}$. But $S_2-S_1\geq R_1+2e-R_1=2e$ so
$\b_1=(S_2-S_1)/2+e$, which implies
$S_{i+1}-S_2+\b_1=S_{i+1}-(S_1+S_2)/2+e\geq 
S_{i+1}-(S_i+S_{i+1})/2+e=(S_{i+1}-S_i)/2+e\geq\b^*_i$ and so
$\b_i=\b^*_i$. If $i\geq 4$ then
$\a_i=\min\{\a_{i-3}(\[\a_4,\ldots,\a_n\] ),R_{i+1}-R_4+\a_3\}$. But
$R_4-R_3\geq 2e$ so $\a_3=(R_4-R_3)/2+e$, which implies
$R_{i+1}-R_4+\a_3=R_{i+1}-(R_3+R_4)/2+e\geq
R_{i+1}-(R_i+R_{i+1})/2+e=(R_{i+1}-R_i)/2+e\geq\a_{i-3}(\[
a_4,\ldots,a_n\] )$. On the other hand $\a^*_i=\min\{\a_{i-3}(\[
a'_4,\ldots,a'_n\] ),R'_{i+1}-R'_4+\a^*_3\}$. But
$R'_4-R'_3=R_4-(R_3-1)>2e$ so $\a^*_3=(R'_4-R'_3)/2+e$ so again we get
$R'_{i+1}-R'_4+\a^*_3\geq (R'_{i+1}-R'_i)/2+e\geq\a_{i-3}(\[
a'_4,\ldots,a'_n\] )$. It follows that $\a^*_i=\a_{i-3}(\[
a'_4,\ldots,a'_n\] )=\a_i$. ($a'_i=\a_i$ for $i\geq 4$.) Also note
that $\a^*_3=(R'_4-R'_3)/2+e>(R_4-R_3)/2+e=\a_3\geq 2e$. Hence
$\a^*_i\geq\a_i$ for $i\geq 3$ with equality if $i\geq 4$.

We have $a'_1=b_1$ so $d(\e a'_{2,i}b_{2,j})=d(\e a'_{1,i}b_{1,j})$
for any $i,j\geq 1$. If moreover $i\geq 3$ then $a'_{1,i}=a_{1,i}$ in
$\ff/\fs$ so $d(\e a'_{2,i}b_{2,j})=d(\e a_{1,i}b_{1,j})$. If $i\geq
3$ and $j\geq 2$ then also $\a^*_i\geq\a_i$ with equality if $i\geq 4$
and $\b^*_j=\b_j$ so $d[\e a_{1,i}b_{1,j}]\geq d[\e a'_{2,i}b_{2,j}]$
with equality if $i\geq 4$. In particular,
$d[-a'_{2,i+1}b_{2,i-1}]=d[-a_{1,i+1}b_{1,i-1}]$ for $i\geq 3$. On the
other hand we have $-a_{1,3}b_1\in\D\fs$ so
$d(-a_{1,3}b_1)=d(-a'_{2,3})=2e$. Hence $2e=d(-a_{1,3}b_1)\geq
d[-a_{1,3}b_1]\geq 2e$, i.e. $d[-a_{1,3}b_1]=2e$ and also
$d[-a'_{2,3}]=\min\{ d(-a'_{2,3}),\a^*_3\} =2e$ ($\a^*_3>2e$). So
$d[-a'_{2,i+1}b_{2,i-1}]=d[-a_{1,i+1}b_{1,i-1}]$ holds for any $i\geq
3$. It follows, by the same reasoning as for Lemma 9.3, that
$A^*_i=A_i$ for any $i\geq 3$ s.t. $i-1$ or $i$ is essential for
$M^*,N^*$. The exception at $i=2$ comes from the fact that $R_3$ and
$R'_3$ are involved in $A_2$ resp. $A^*_2$ and $R'_3\neq R_3$. At
$i=2$ we have $d[-a_{1,3}b_1]=d[-a'_{2,3}]=2e$ so $A_2=\min\{
(R_3-S_2)/2+e,R_3-S_2+2e\}$, by Lemma 8.12, and $A^*_2=\min\{
(R'_3-S_2)/2+e,R'_3-S_2+2e\}$. But $R_3-S_2+2e\leq R_1-(R_1+2e)+2e=0$
so $R_3-S_2+2e\leq\h 12(R_3-S_2+2e)=(R_3-S_2)/2+e$ so
$A_2=R_3-S_2+2e$. Similarly since $R'_3-S_2+2e=R_3-S_2+2e-1<0$ we have
$A^*_2=R'_3-S_2+2e=A_2-1<0$.

We prove now 2.1(ii). At $i=2$ then $d[-a'_2b_2]\geq 0>A^*_2$. If
$i\geq 3$ and at least one of $i-1,i$ is essential for $M^*,N^*$ then
$A_i^*=A_i$ so $d[a'_{2,i}b_{2,i}]\geq d[a_{1,i}b_{1,i}]\geq
A_i=A^*_i$. If both $i-1,i$ are not essential then
$d[a'_{2,i}b_{2,i}]\geq A^*_i$ by Lemma 2.12.

We prove now 2.1(iii). Let $3\leq i\leq n-1$ s.t. $R'_{i+1}>S_{i-1}$
and $A^*_{i-1}+A^*_i>2e+R'_i-S_i$. By Lemma 2.13 this implies that
$i-1$ is essential for $M^*,N^*$. If $i\geq 4$ this implies
$A_{i-1}=A^*_{i-1}$ and $A_i=A^*_i$ so
$A_{i-1}+A_i=A^*_{i-1}+A^*_i>2e+R'_i-S_i=2e+R_i-S_i$. We also have
$R_{i+1}=R'_{i+1}>S_{i-1}$ so $[b_1,\ldots,b_{i-1}]\rep [a_1,\ldots,a_i]\ap
[a'_1,\ldots,a'_i]$. But $a'_1=b_1$ so $[b_2,\ldots,b_{i-1}]\rep
[a'_2,\ldots,a'_i]$. If $i=3$ then the fact that $2$ is essential only
implies $A_3=A^*_3$. Since also $A_2=A^*_2+1$ we have
$A_2+A_3=A^*_2+A^*_3+1>2e+R'_3-S_3+1=2e+R_3-S_3$. Since also
$R_4=R'_4>S_2$ we get $[b_1,b_2]\rep [a_1,a_2,a_3]\ap
[a'_1,a'_2,a'_3]$. But $b_1=a'_1$ so $[b_2]\rep [a'_2,a'_3]$.

For 2.1(iv) suppose that $3\leq i\leq n-2$ s.t. $S_i\geq
R'_{i+2}>S_{i-1}+2e\geq R'_{i+1}+2e$. We have $R_{i+1}=R'_{i+1}$ and
$R_{i+2}=R'_{i+2}$ so $S_i\geq R_{i+2}>S_{i-1}+2e\geq
R_{i+1}+2e$. Therefore $[b_1,\ldots,b_{i-1}]\rep [a_1,\ldots,a_{i+1}]\ap
[a'_1,\ldots,a'_{i+1}]$. But $a'_1=b_1$ so $[b_2,\ldots,b_{i-1}]\rep
[a'_2,\ldots,a'_{i+1}]$. \qed  

\blm Suppose that $M,N$ are binary and $R_i\leq S_i$ for $i=1,2$.

(i) If $M\ap\[\pi^{R_1}\e_1,\pi^{R_2}\e_2\]$,
$N\ap\[\pi^{S_1}\e_1,\pi^{S_2}\e_2\]$ and $R_i\ev S_i\m2$ for $i=1,2$
then $N\rep M$. 

(ii) If $FM\ap FN$ and one of $M,N$ is $\ap\pi^r\aa$ for some
$r\in\ZZ$ then $N\rep M$. 
\elm
\pf (i) If $M\ap\[\pi^{R_1}\e_1,\pi^{R_2}\e_2\]$ relative to a BONG
$x_1,x_2$ and $k:=(S_2-R_2)/2$ then the lattice $K:=\[ x_1,\pi^kx_2\]$
exists and $K\sbq M$ by [B1, Lemma 2.7(ii)]. (We have $k\geq 0$ so
$\[\pi^kx_2\]\sbq\[ x_2\]$.) Since
$K\ap\[\pi^{R_1}\e_1,\pi^{R_2+2k}\e_2\]
=\[\pi^{R_1}\e_1,\pi^{S_2}\e_2\]$ we have
$\[\pi^{R_1}\e_1,\pi^{S_2}\e_2\]\rep
M\ap\[\pi^{R_1}\e_1,\pi^{R_2}\e_2\]$. Similarly, since $-S_1\leq -R_1$
and $-S_1\ev -R_1\m2$ we have
$\[\pi^{-S_2}\e_2,\pi^{-R_1}\e_1\]\rep \[\pi^{-S_2}\e_2,\pi^{-S_1}\e_1\]$,
i.e. $K^\*\rep N^\*$. By duality we get $N\rep K\rep M$. 

(ii) We have $\nn M=\p^{R_1}\spq\p^{S_1}=\nn N$ and $\nn
M^\*=\p^{-R_2}\sbq\p^{-S_2}=\nn N^\*$. If $M\ap\pi^r\aa$ then $M$ is
$\nn M$-maximal. Since $\nn N\sbq\nn M$ we have $N\rep M$. (See [OM,
82:18, 82:22 and 91:2].) If $N\ap\pi^r\aa$ then $N^\*\ap\pi^{-r}\aa$ so
$N^\*$ is $\nn N^\*$-maximal. Since also $\nn M^\*\sbq\nn N^\*$ we get
$M^\*\rep N^\*$ so $N\rep M$. \qed 

\blm If $M,N$ are ternary over the same quadratic space $V$ and
$R_1=R_3=S_1=S_3$ then $N\rep M$ iff $\a_1\leq\b_1$ and
$\a_2\geq\b_2$.
\elm
\pf The necessity follows from Lemma 9.4 and the fact that $N\rep M$
implies $N\leq M$.

For sufficiency we use Lemma 9.5. We have $R_1\ev R_2\ev R_3\m2$ so
$\ord\det V=R_1+R_2+R_3\ev R_1\m2$. So $\det V=\pi^{R_1}\xi
=\pi^{S_1}\xi$ for some $\xi\in\ooo$. If $M\ap\[ a_1,a_2,a_3\]$ with
$a_i=\pi^{R_i}\e_i$ then $\xi=\e_{1,3}$. Similarly for $N$. By Lemma
9.5 if $V$ is isotropic then $M\ap M_1\pp M_2$ and $N\ap N_1\pp N_2$,
where $M_1\ap\[\pi^{R_1+\a_1},-\pi^{R_1-\a_2}\]$ and
$M_2\ap\la -\pi^{R_1}\xi\ra$, $N_1\ap\[\pi^{S_1+\b_1},-\pi^{S_1-\b_2}\]
=\[\pi^{R_1+\b_1},-\pi^{R_1-\b_2}\]$ and $N_2\ap\la -\pi^{S_1}\xi\ra
=\la -\pi^{R_1}\xi\ra\ap M_2$. If $V$ is not isotropic then
$M_1\ap\[\pi^{R_1+\a_1},-\pi^{R_1-\a_2}\D\]$,
$N_1\ap\[\pi^{R_1+\b_1},-\pi^{R_1-\b_2}\D\]$ and $M_2\ap
N_2\ap\la -\pi^{R_1}\xi\D\ra$.

Since $M_2\ap N_2$ it is enough to prove that $N_1\rep M_1$. Since
$FM\ap FN$ and $FM_2\ap FN_2$ we have $FM_1\ap FN_1$. Since
$R_1+\a_1\leq R_1+\b_1$ and $R_1-\a_2\leq R_1-\b_2$ we can use Lemma
9.7. We have $\a_2=R_2-R_1+\a_1\ev\a_1\m2$ and similarly
$\b_1\ev\b_2\m2$. Thus if both $\a_1$ and $\b_1$ are odd then
$R_1+\a_1\ev R_1+\b_1\m2$ and $R_1-\a_2\ev R_1-\b_2\m2$ so $N_1\rep
M_1$ by Lemma 9.7(i). If $\a_1$ is even then $\a_1=(R_2-R_1)/2+e$ and,
by Lemma 9.5(i), $V$ is isotropic. By Lemma 9.5(iii) we have
$M\ap\pi^{(R_1+R_2)/2}\aa$ and so $N_1\rep M_1$ by Lemma
9.7(ii). Similarly if $\b_1$ is even we get $N\ap\pi^{(S_1+S_2)/2}\aa$
and again $N_1\rep M_1$ by Lemma 9.7(ii). \qed

\blm If $R_1,R_2,R_3,\a_1$ are integers with $R_1=R_3$ and $V$ is a
ternary quadratic field then there is a lattice $L$ over $V$ with
$R_i(L)=R_i$ for $1\leq i\leq 3$ and $\a_1(L)=\a_1$ iff $R_1\ev
R_2\m2$, $\det V=\pi^{R_1}\xi$ for some $\xi\in\ooo$, $\max\{
0,R_2-R_1\}\leq\a_1\leq (R_2-R_1)/2+e$ and if $\a_1$ is even then
$\a_1=(R_2-R_1)/2+e$ and $V$ is isotropic. Moreover if $\a_1$ is odd
then the BONG of $L$ can be chosen s.t. $d(-a_{2,3})=\a_1$.

(Note that by Lemma 7.11 the lattice $L$ with the properties above is
unique up to an isomorphism.)
\elm
\pf We prove first the necessity. By 8.7  $R_1=R_3$ implies $R_1\ev
R_2\m2$. If $\a_i=\pi^{R_i}\e_i$ then in $\ff/\fs$ we have $\det
V=a_{1,3}=\pi^{R_1+R_2+R_3}\e_{1,3}=\pi^{R_1}\xi$, where
$\xi=\e_{1,3}$. (We have $R_1+R_2+R_3\ev R_1\m2$.) Now $\a_1\leq 
(R_2-R_1)/2+e$ by the definition of $\a_1$ and $\a_1\geq\max\{
0,R_2-R_1\}$ follows by [B3, Lemma 2.7(i) and (iii)] from
$R_2-R_1=R_2-R_3\leq 2e$. If $\a_1$ is even then $\a_1=(R_2-R_1)/2+e$
follows from [B3, Lemma 2.7(iv)] and $V$ is isotropic by Lemma
9.5(i). 

We prove now the sufficiency. Note that $\max\{
0,R_2-R_1\}\leq\a_1\leq (R_2-R_1)/2+e$ implies $0\leq (R_2-R_1)/2+e$
and $R_2-R_1\leq (R_2-R_1)/2+e$ so $-2e\leq R_2-R_1\leq 2e$.

Suppose first that $\a_1$ is even so $\a_1=(R_2-R_1)/2+e$ and $V$
is isotropic. We take $a_1=a_3=-\pi^{R_1}\xi$ and
$a_2=\pi^{R_2}\xi$. Then $-a_2/a_1=\pi^{R_2-R_1}$ and
$-a_3/a_2=\pi^{R_3-R_2}$ are squares. Since $R_2-R_1$ and $R_3-R_2=R_1-R_2$
are $\geq -2e$ and $R_2-R_1+d(-a_2/a_1)=R_3-R_2+d(-a_3/a_2)=\j >0$ we
have $a_2/a_1,a_3/a_2\in\aaa$. Together with $R_1=R_3$, this implies
that there is a lattice $L\ap\[ a_1,a_2,a_3\]$. In $\ff/\fs$ we have
$\det FL=\a_{1,3}=\pi^{R_1+R_2+R_3}\xi=\pi^{R_1}\xi=\det V$. Also
$-a_2/a_1\in\fs$ so $[a_1,a_2]$ is hyperbolic and $[a_1,a_2,a_3]$ is
isotropic, same as $V$. Thus $FL\ap V$. We have $R_i(L)=\ord a_i=R_i$
so we still have to prove that $\a_1(L)=\a_1$. Since
$d(-a_{1,2})=d(-a_{2,3})=\j$ we have $\a_1(L)=\min\{ (R_2-R_1)/2+e,
R_2-R_1+d(-a_{1,2}), R_3-R_2+d(-a_{2,3})\} =(R_2-R_1)/2+e=\a_1$.

Suppose now that $\a_1$ is odd. Since $0\leq\a_1\leq (R_2-R_1)/2+e\leq
(2e)/2+e=2e$ we have $\a_1,2e-\a_1\in d(\ooo )$. Let $\e,\eta\in\ooo$
s.t. $d(\e )=\a_1$, $d(\eta )=2e-\a_1$ and $(\e,\eta )_\p =-1$. If $V$
is isotropic let $a_1=-\pi^{R_1}\xi\e$, $a_1=\pi^{R_2}\xi\e$ and
$a_3=-\pi^{R_3}\xi$. If $V$ is not isotropic let $a_1=-\pi^{R_1}\xi\e$,
$a_2=\pi^{R_2}\xi\e\eta$ and $a_3=-\pi^{R_3}\xi\eta$. First we prove that
there is a lattice $L\ap\[ a_1,a_2,a_3\]$. We have $\ord a_i=R_i$ and
$R_1=R_3$ so we still need $a_2/a_1,a_3/a_2\in\aaa$. Now $R_2-R_1\geq
-2e$ and $R_3-R_2=R_1-R_2\geq -2e$ so we need $R_2-R_1+d(-a_{1,2}),
R_3-R_2+d(-a_{2,3})\geq 0$. If $V$ is isotropic $d(-a_{1,2})=\j$ and
if $V$ is not isotropic $d(-a_{1,2})=d(\eta )=2e-\a_1$. In both cases
$d(-a_{1,2})\geq 2e-\a_1$ so $R_2-R_1+d(-a_{1,2})\geq
R_2-R_1+2e-\a_1\geq R_2-R_1+2e-2\a_1\geq 0$. (We have $\a_1\geq 0$ and
$(R_2-R_1)/2+e\geq\a_1$.) In both cases when $V$ is isotropic or not
we have $d(-a_{2,3})=d(\e )=\a_1$ so
$R_3-R_2+d(-a_{2,3})=R_1-R_2+\a_1\geq R_1-R_2+R_2-R_1=0$. Thus $L\ap\[
a_1,a_2,a_3\]$ exists. In $\ff/\fs$ we have $\det
FL=a_{1,3}=\pi^{R_1+R_2+R_3}\xi=\pi^{R_1}\xi=\det V$. If $V$ is
isotropic then $-a_2/a_1\in\fs$ so $[a_1,a_2]$ is hyperbolic and
$FL\ap [a_1,a_2,a_3]$ is isotropic, same as $V$. If $V$ is not
isotropic then $(-a_{1,2},-a_{2,3})_\p =(\eta,\e )_\p=-1$ so $FL\ap
[a_1,a_2,a_3]$ is not isotropic, same as $V$. So $FL\ap V$. Hence we
still have to prove that $\a_1(L)=\a_1$. Now $\a_1(L)=\min\{
(R_2-R_1)/2+e, R_2-R_1+d(-a_{1,2}), R_3-R_1+d(-a_{2,3})\}$ But
$(R_2-R_1)/2+e\geq\a_1$, $R_2-R_1+d(-a_{1,2})\geq
R_2-R_1+2e-\a_1\geq\a_1$ (because $(R_2-R_1)/2+e\geq\a_1$) and
$R_3-R_1+d(-a_{2,3})=d(-a_{2,3})=\a_1$. Thus $\a_1(L)=\a_1$. \qed

\blm Let $M$ be a lattice with $R_1=R_3$ and, if $n\geq 4$, $R_4\geq
R_2+2$ and let $V\ap [a_1,a_2,a_3]$. Let $S_2=R_2+2$ and $S_i=R_i$ for
$i\neq 2$. If $\b_1$ is an  integer s.t. $S_1,S_2,S_3,\b_1$ and $V$
satisfy the conditions of Lemma 9.9 and
$\a_1\leq\b_1\leq\min\{\a_1+2,\a_3\}$ then there is a lattice
$N\sb M$, $N\ap\[ b_1,\ldots,b_n\]$, such that $[M:N]=\p$,
$R_i(N)=S_i$ for all $i$ and $\a_1(N)=\b_1$. 

Moreover $N$ and the BONG of $M$ can be chosen s.t. $b_i=a_i$ for
$i\geq 4$. 
\elm 
\pf Regardless of the BONG of $M$ we have $R_1+R_2+R_3\ev R_1=S_1\m2$
so in $\ff/\fs$ we have $\det
V=a_{1,3}=\pi^{R_1+R_2+R_3}\xi=\pi^{S_1}\xi$. So the condition that
$\det V=\pi^{S_1}\xi$ for some $\xi\in\ooo$ from Lemma 9.9 holds
regardless of the BONG.

By Corollary 8.11 we may assume, after a change of BONG, that
$\a_1=\a_1(\[ a_1,a_2\] )$. Since $\a_1\leq\a_1(\[ a_1,a_2,a_3\]
)\leq \a_1(\[ a_1,a_2\] )$ this implies $\a_1=\a_1(\[
a_1,a_2,a_3\] )$. If $\b_1$ is even then
$\b_1=(S_2-S_1)/2+e=(R_2-R_1)/2+e+1\geq\a_1+1$. If $\b_1=\a_1+2$
then $\a_1$ is even so $\a_1=(R_2-R_1)/2+e$, which implies
$\b_1=(R_2-R_1)/2+e+2$. Contradiction. So $\b_1=\a_1+1$ and
$\a_1=(R_2-R_1)/2+e$. Thus regardless of the BONG of $M$ we have
$\a_1=(R_2-R_1)/2+e\geq\a_1(\[ a_1,a_2\] )\geq\a_1$ so $\a_1(\[
a_1,a_2\] )$. So in this case we don't have to change the BONG
of $M$ in order to have $\a_1=\a_1(\[ a_1,a_2\] )$, a change of BONGs
which may alter the condition that $V\ap [a_1,a_2,a_3]$ is isotropic,
required by Lemma 9.9. 

By Lemma 9.9 there is a lattice $\( N\ap\[ b_1,b_2,b_3\]$ over $V$
with $R_i(\( N)=S_i$ and $\a_1(\( N)=\b_1$. We want to prove that $\(
N\rep\( M$, where $\( M\ap\[ a_1,a_2,a_3\]$. We have $F\( M=F\( N=V$
and $R_1=R_3=S_1=S_3$ so, by Lemma 9.8, it is enough to prove that
$\a_1(\( M)\leq\a_1(\( N)$ and $\a_2(\( M)\geq\a_2(\( N)$. We have
$\a_1(\( M)=\a_1\leq\b_1=\a_1(\( N)$ and $\a_2(\(
M)=R_1-R_2+\a_1=S_1-S_2+\a_1+2\geq S_1-S_2+\b_1=\a_2(\( N)$ (see
8.7). So $\( N\rep\( M$. 

We define $N\ap\[ b_1,\ldots,b_n\]$ with $b_i=a_i$ for $i\geq 4$. We
have $\ord b_i=S_i$ for all $i$. We prove first that a lattice $N\ap\[
b_1,\ldots,b_n\]$ does exist. If $i\neq 2$ we have $S_i=R_i\leq
R_{i+2}=S_{i+2}$ and if $i=2$ then $S_2=R_2+2\leq R_4=S_4$. For
$i=1,2$ we have $b_{i+1}/b_i\in\aaa$ because $\( N\ap\[ b_1,b_2,b_3\]$
exists. For $i\geq 4$ we have $b_{i+1}/b_i=a_{i+1}/a_i\in\aaa$ so we
still need $b_4/b_3\in\aaa$. We have $S_4-S_3=R_4-R_3\geq -2e$ so we
still need $S_4-S_3+d(-b_{3,4})\geq 0$, i.e. $R_4-R_3+d(-a_4b_3)\geq
0$. Since $a_4/a_3\in\aaa$ we have $R_4-R_3+d(-a_{3,4})\geq 0$ so we
still need $R_4-R_3+d(a_3b_3)\geq 0$. Now $a_{1,3}=b_{1,3}=\det V$ in
$\ff/\fs$ so $a_{1,3}b_{1,3}\in\fs$, which implies
$d(a_3b_3)=d(a_{1,2}b_{1,2})$. But $d(-b_{1,2})\geq S_1-S_2+\b_1$ and
$d(-a_{1,2})\geq R_1-R_2+\a_1=S_1-S_2+\a_1+2\geq S_1-S_2+\b_1$ so
$d(a_3b_3)=d(a_{1,2}b_{1,2})\geq S_1-S_2+\b_1$. This implies
$R_4-R_3+d(a_3b_3)\geq R_4-R_3+S_1-S_2+\b_1\geq\b_1\geq 0$. (We have
$R_3=S_1$ and $R_4\geq R_2+2=S_2$.) Thus $N$ exists. Obviously
$R_i(N)=S_i$. 

Since $\( N\rep\( M$ we have $\[ a_3\1,a_2\1,a_1\1\]\ap\( M^\*\rep\(
N^\*\ap\[ b_3\1,b_2\1,b_1\1\]$. By [B1, Lemma 2.7(ii)] it follows that
$M^\*\ap\[ a_n\1,\ldots,a_1\1\]\rep\[
a_n\1,\ldots,a_4\1,b_3\1,b_2\1,b_1\1\]\ap N^\*$ so $N\rep M$. We may
assume that $N\sb M$. Since $\ord volN=S_1+\cdots +S_n=R_1+\cdots
+R_n+2=\ord volN+2$ we have $[M:N]=\p$.

We still need $\a_1(N)=\b_1$. We have $\a_1(N)=\min (\{ \a_1(\[
b_1,b_2,b_3\] )\}\cup\{ S_{i+1}-S_1+d(-b_{i,i+1})|~3\leq i\leq n-1\}
)$. But $\a_1(\[ b_1,b_2,b_3\] )=\b_1$ and for $i\geq 4$ we have
$S_{i+1}-S_1+d(-b_{i,i+1})=R_{i+1}-R_3+d(-a_{i,i+1})\geq\a_3\geq\b_1$.
So we still need $\b_1\leq S_4-S_1+d(-b_{3,4})=S_4-S_1+d(-a_4b_3)$. We
have $S_4-S_1+d(-a_{3,4})=R_4-R_3+d(-a_{3,4})\geq\a_3\geq\b_1$ and
$d(a_3b_3)\geq S_1-S_2+\b_1$ (see above) so $S_4-S_1+d(a_3b_3)\geq
S_4-S_2+\b_1\geq\b_1$. Thus $S_4-S_1+d(-a_4b_3)\geq\b_1$. \qed

\blm If $L\ap\[ a_1,a_2\]$ and $R_2-R_1+d(-a_{1,2})=1$ then $L'\ap
L^{\pi\e}$ for some $\e\in\ooo$.

(As usual $L':=\{ x\in L\mid x\text{ not a norm generator}\}$.)
\elm
\pf We cannot have $R_2-R_1=-2e$ since this would imply $a_2/a_1\in
-\h 14\ooo^2$ or $-\h\D 4\ooo^2$ so $R_2-R_1+d(-a_{1,2})=\j$ or $0$,
respectively. By Lemma 7.1 $L'$ is a lattice with $[L:L']=\p$. 

Let $L'\ap\[ b_1,b_2\]$ with $\ord b_i=S_i$. By the definition of $L'$
we have $\p^{R_1}=\nn L'\sb\nn L=\p^{R_1}$, i.e. $S_1>R_1$. We have
$[L:L']=\p$ so $b_{1,2}=\det L'=\pi^2\det L=\pi^2a_{1,2}$, which
implies $S_1+S_2=R_1+R_2+2$.

We have $b_2/b_1\in\aaa$ so $S_2-S_1\geq
-d(-b_{1,2})=-d(-a_{1,2})=R_2-R_1-1$. By subtracting from
$S_1+S_2=R_1+R_2+2$ we get $2S_1\leq 2R_1+3$ so $R_1<S_1\leq
R_1+1$. So $S_1=R_1+1$, which implies $b_1=\pi\e a_1$ for some
$\e\in\ooo$. But $b_{1,2}=\pi^2a_{1,2}$ so $b_2=\pi\e a_2$ (in
$\ff/\ooo^2$). Thus $L'\ap\[\pi\e a_1,\pi\e a_2\]\ap L^{\pi\e}$. \qed

\blm Suppose that $K\leq M$, $R_1=T_1$ and $M,K$ don't satisfy the
hypothesis of Lemma 9.3 or Lemma 9.6. Then there is $N\sb M$ with
$[M:N]=\p$ s.t. $K\leq N$.
\elm
\pf Since $M,K$ don't satisfy the hypothesis of Lemma 9.3 we have
$R_1=R_3$, $R_2\neq T_2$, $R_2-R_1\neq 2e$ and $R_2<R_4$. Also if
$T_2=R_2+1$ then $R_1<R_5$. We have $R_1+R_2\leq T_1+T_2=R_1+T_2$ so
$R_2\neq T_2$ implies $R_2<T_2$. For consequences of $R_1=R_3$ see
8.7. We have that $R_2-R_1$ is even and $\leq 2e$. But $R_2-R_1\neq
2e$ so $R_2-R_1\leq 2e-2$. 

By Lemma 8.12 we have $C_1=\a_1$. Hence $\a_1\geq d[a_1c_1]\geq
C_1=\a_1$ so $d[a_1c_1]=\a_1$. By 8.7 also $d[-a_{2,3}]=\a_1$
and so $\a_1\leq d[-a_{1,3}c_1]\leq\c_1$. If
$d[-a_{1,3}c_1]=\a_1<\c_1$ then $M,K$ would satisfy the hypothesis of
Lemma 9.3. So either $d[-a_{1,3}c_1]>\a_1$ or
$d[-a_{1,3}c_1]=\a_1=\c_1$. 
\vskip 3mm

CLAIM If $R_2+2\leq R_4$, $R_2+2\leq T_2$ and $N$ is another lattice
s.t. $M,N$ are like in Lemma 9.10 then $N,K$ satisfy the condition
2.1(i) and $B_i\leq C_i$ for $i\geq 3$. Moreover $K\leq N$ iff
$\b_1\leq d[-a_{1,3}c_1]$ and $S_{i+1}-S_2+\b_1\geq B_i$ for $i\geq 2$.
\vskip 2mm

\pf We show that in fact the two conditions are equivalent to the
(ii) part of the Theorem 2.1, while (i), (iii) and (iv) hold
unconditionally.

We have $S_2=R_2+2$ and $S_i=R_i$ for $i\neq 2$ so $M,N$ have type I
(see Definition 11, following Lemma 6.7) with $s=t=t'=u=2$. It follows
that for $i\geq 2$ we have
$d[a_{1,i}b_{1,i}]=A_i=\b_i=\min\{\a_i,S_{i+1}-S_3+\b_2\}$. Also for
any $i\geq 2$, $j\geq 0$ and $\e\in\ff$ we have $d[\e
b_{1,i}c_{1,j}]=\min\{ d[\e a_{1,i}c_{1,j}],\b_i\}$. (See Lemma
6.9(iv) and 6.16.) Thus for any $i\geq 2$ we have $B_i=\min\{
(S_{i+1}-T_t)/2+e, S_{i+1}-T_i+d[-b_{1,i+1}c_{1,i-1}],
S_{i+1}+S_{i+2}-T_{i-1}-T_i+d[b_{1,i+2}c_{1,i-2}]\}\leq\min\{
(R_{i+1}-T_t)/2+e, R_{i+1}-T_i+d[-a_{1,i+1}c_{1,i-1}],
R_{i+1}+R_{i+2}-T_{i-1}-T_i+d[a_{1,i+2}c_{1,i-2}]\} =C_i$. (We have
$S_{i+1}=R_{i+1}$, $S_{i+2}=R_{i+2}$, $d[-b_{1,i+1}c_{1,i-1}]\leq
d[-a_{1,i+1}c_{1,i-1}]$ and $d[b_{1,i+2}c_{1,i-2}]\leq
d[a_{1,i+2}c_{1,i-2}]$.)

We prove first 2.1(i). If $i=1$ then $S_1=R_1=T_1$. If $i=2$ then
$S_2=R_2+2\leq T_2$. If $i\geq 3$ then either $S_i=R_i\leq T_i$ or
$S_i+S_{i+1}=R_i+R_{i+1}\leq T_{i-1}+T_i$. Next we prove that 2.1(ii)
is equivalent to $\b_1\leq d[-a_{1,3}c_1]$ and
$S_{i+1}-S_2+\b_1\geq B_i$ for $i\geq 2$. Since $S_1=T_1$ we have
$B_1=\b_1$ by Lemma 8.12. Since $S_1=S_3$ we have by 8.7 
$d[-b_{2,3}]=\b_1$ so $d[b_1c_1]\geq B_1=\b_1$ is equivalent to
$d[-b_{1,3}c_1]\geq\b_1$, which in turn, since
$d[a_{1,3}b_{1,3}]=\b_3\geq S_1-S_3+\b_1=\b_1$, is equivalent to
$d[-a_{1,3}c_1]\geq\b_1$. If $i\geq 2$, since
$d[a_{1,i}c_{1,i}]\geq C_i\geq B_i$ and $d[b_{1,i}c_{1,i}]=min\{
d[a_{1,i}c_{1,i}],\b_i\}$ the relation $d[b_{1,i}c_{1,i}]\geq B_i$
is equivalent to $\b_i\geq B_i$. But $\b_i=\min\{\a_i,
S_{i+1}-S_3+\b_2\}$ and $\a_i\geq C_i\geq B_i$ so $\b_i\geq B_i$
is equivalent to $S_{i+1}-S_2+\b_1=S_{i+1}-S_3+\b_2\geq B_i$. For
2.1(iii) and (iv) we note that for any $i\geq 3$ we have
$a_{i+1}=b_{i+1}$,\ldots,$a_n=b_n$, which, together with
$[a_1,\ldots,a_n]\ap [b_1,\ldots,b_n]$, implies $[b_1,\ldots,b_i]\ap
[a_1,\ldots,a_i]$. We prove now 2.1(iii). If $i=2$ then $S_3=T_1$ so
we have nothing to prove. If $i\geq 3$, $S_{i+1}>T_{i-1}$ and
$B_{i-1}+B_i>2e+S_i-T_i$ then $R_{i+1}=S_{i+1}>T_{i-1}$ and
$C_{i-1}+C_i\geq B_{i-1}+B_i>2e+S_i-T_i=2e+R_i-T_i$ so
$[c_1,\ldots,c_{i-1}]\rep [b_1,\ldots,b_i]\ap [a_1,\ldots,a_i]$. For
2.1(iv) if $i\geq 2$ and $T_i\geq S_{i+2}>T_{i-1}+2e\geq S_{i+1}+2e$
then also $T_i\geq R_{i+2}>T_{i-1}+2e\geq R_{i+1}+2e$ (we have
$R_{i+1}=S_{i+1}$ and $R_{i+2}=S_{i+2}$) and so
$[c_1,\ldots,c_{i-1}]\rep [a_1,\ldots,a_{i+1}]\ap [b_1,\ldots,b_{i+1}]$.
\vskip 3mm

We consider separately the cases $d[-a_{1,3}c_1]<\a_1$ and
$d[-a_{1,3}c_1]=\a_1=\c_1$.\vskip 2mm

1. $d[-a_{1,3}c_1]=\a_1=\c_1$. We have $R_3-R_2=R_1-R_2\geq 2-2e$
so $\a_2\geq 1$. 

Suppose first that $\a_2>1$ so $\a_1=R_2-R_1+\a_2\geq R_2-R_1+2$. Note
that $\c_1=\a_1\leq (R_2-R_1)/2+e<(T_2-T_1)/2+e$ so $\a_1=\c_1$ is
odd. We take $N$ to be a lattice like in Lemma 9.10 with
$\b_1=\a_1$. First we have to prove that this is possible. We have
$\a_1\geq R_2-R_1+2=S_2-S_1$ and $\a_1\geq 0$. Also $\a_1\leq
(R_2-R_1)/2+e<(S_2-S_1)/2+e$. Hence $\max\{ 0,S_2-S_1\}\leq\b_1\leq
(S_2-S_1)/2+e$ from Lemma 9.9 is fulfilled. Also $\b_1=\a_1$ is odd so
the conditions that $\b_1=(S_2-S_1)/2+e$ and $V:=[a_1,a_2,a_3]$ is
isotropic from Lemma 9.9 are not necessary. Also $\a_1\leq
R_3-R_1+\a_3=\a_3$ so the condition
$\a_1\leq\b_1\leq\min\{\a_1+2,\a_3\}$ is also fulfilled. So $N$
exists. In order to use our claim above we still have to prove that
$R_2+2\leq R_4$, $R_2+2\leq T_2$, $\b_1\leq d[-a_{1,3}c_1]$ and 
$S_{i+1}-S_2+\b_1\geq B_i$ for $i\geq 2$. If $R_4=R_2+1$ then
$R_3+R_4$ is odd so so $\a_2\leq R_4-R_2+d(-a_{3,4})=R_4-R_2=1$.
Contradiction. So $R_2+2\leq R_4$. If $T_2=R_2+1$ then
$T_2-T_1=R_2-R_1+1$ is odd and $\leq 2e-1$ so
$\a_1=\c_1=T_2-T_1=R_2-R_1+1$. Contradiction. So $R_2+2\leq T_2$.
We have $\b_1=\a_1=[-a_{1,3}c_1]$ and for $i\geq 2$ we have
$S_{i+1}-S_2+\b_1=S_{i+1}-R_2-2+\c_1\geq S_{i+1}-T_2+\c_1\geq
S_{i+1}-T_i+\c_{i-1}\geq B_i$. So $K\leq N$.

Suppose now that $\a_2=1$ so $\c_1=\a_1=R_2-R_1+\a_2=R_2-R_1+1$.
Since $\c_1=R_2-R_1+1<2e$ we have $T_2-T_1\leq\c_1\leq R_2-R_1+1$
so $T_2\leq R_2+1$. But $T_2>R_2$ so $T_2=R_2+1$. This implies
$R_1<R_5$. By Corollary 8.11 we may assume that $\a_1=\a_1(\[
a_1,a_2\] )$. Since $\a_1\leq\a_1(\[ a_1,a_2,a_3\] )\leq\a_1(\[
a_1,a_2\] )$ we have $\a_1=\a_1(\[ a_1,a_2,a_3\] )$. Since $\a_1$
is odd by Lemma 9.9, after a change of BONG for $\[
a_1,a_2,a_3\]$, we may assume that
$d(-a_{2,3})=\a_1=R_2-R_1+1=R_2-R_3+1$. If $J\ap\[ a_2,a_3\]$ then
by Lemma 9.11 we have $J\sp J'$ with $J'\ap J^{\pi\e}\ap\[\pi\e
a_2,\pi\e a_3\]$ for some $\e\in\ooo$. We take $N\ap\[
b_1,\ldots,b_n\]$ relative to a good BONG with $b_2=\pi\e a_2$,
$b_3=\pi\e a_3$ and $b_i=a_i$ for $i\neq 2,3$. First we show that such
a lattice exists. If $S_i=\ord b_i$ we have $S_2=R_2+1$, $S_3=R_3+1$
and $S_i=R_i$ for $i\neq 2,3$. We have $S_1=R_1=R_3<S_3$,
$S_2=R_2+1\leq R_4=S_4$, $S_3=R_3+1=R_1+1\leq R_5=S_5$ and
$S_i=R_i\leq R_{i+2}=S_{i+2}$ for $i\geq 4$. We now prove that
$b_{i+1}/b_i\in\aaa$. We have $\a_2=1$ so $R_3-R_2\leq 1$. But
$R_3-R_2$ is even so $R_1=R_3\leq R_2$. Thus $\ord
b_2/b_1=S_2-S_1=R_2-R_1+1>0$ and $\ord b_4/b_3=S_4-S_3=R_4-R_3-1\geq
R_2-R_3\geq 0$, which implies $b_2/b_1,b_4/b_3\in\aaa$. If $i\geq 4$
as well as if $i=2$ we have $b_{i+1}/b_i=a_{i+1}/a_i\in\aaa$. Thus $N$
exists. Since $\[ b_2,b_3\]\rep\[ a_2,a_3\]$ we get by [B1, Lemma
2.7(ii)] that $\[ a_1,b_2,b_3\]\rep\[ a_1,a_2,a_3\]$. By duality $\[
a_3\1,a_2\1,a_3\1\]\rep\[ b_3\1,b_2\1,a_1\1\]$. Again by [B1, Lemma
2.7(ii)] we have $M^\*\ap\[ a_n\1,\ldots,a_1\1\]\rep\[
a_n\1,\ldots,a_4\1,b_3\1,b_2\1,a_1\1\]\ap N^\*$ so $N\rep M$. We may
assume that $N\sb M$. Since $\ord vol
N=S_1+\cdots +S_n=R_1+\cdots +R_n+2=\ord vol M+2$ we have $[M:N]=\p$.

We prove now that $K\leq N$. First we prove 2.1(i). At $i=1$ we have
$S_1=T_1$. At $i=2$ we have $S_2=R_2+1=T_2$. We cannot have
$T_3=T_1$ since $T_2-T_1=R_2-R_1+1$ is odd. Thus $T_3\geq
T_1+1=R_3+1=S_3$. For $i\geq 4$ we have either $S_i=R_i\leq T_i$
or $S_i+S_{i+1}=R_i+R_{i+1}\leq T_{i-1}+T_i$. 

Before proving 2.1(ii) - (iv) note that $S_2=R_2+1$, $S_3=R_3+1$ and
$S_i=R_i$ for $i\neq 2,3$ so $M,N$ have type III (see Definition 11,
following Lemma 6.7) with $s=t=2$ and $t'=u=3$. Thus for $i\geq 3$ we
have $A_i=\b_i=\min\{\a_i,S_{i+1}-S_4+\b_3\}$. So
$d[a_{1,i}b_{1,i}]=\b_i$ if $i\geq 3$ and $d[\e
a_{1,i}c_{1,j}]\geq d[\e b_{1,i}c_{1,j}]$ if $i\geq 3$, $j\geq 0$
and $\e\in\ff$. (See Lemmas 6.9 and 6.16.) Also $S_i=R_i$ for
$i\geq 4$. Thus for $i\geq 3$ we have $B_i=\min\{
(S_{i+1}-T_i)/2+e, S_{i+1}-T_i+d[-b_{1,i+1}c_{1,i-1}],
S_{i+1}+S_{i+2}-T_{i-1}-T_i+d[b_{1,i+2}c_{1,i-2}]\}\leq\min\{
(R_{i+1}-T_i)/2+e, R_{i+1}-T_i+d[-a_{1,i+1}c_{1,i-1}],
R_{i+1}+R_{i+2}-T_{i-1}-T_i+d[a_{1,i+2}c_{1,i-2}]\}  =C_i$. 

We prove now 2.1(ii). We have $R_1=S_1=T_1$ so by Lemma 8.12
$A_1=C_1=\a_1=R_2-R_1+1$ and $B_1=\b_1$. But $S_2-S_1=R_2-R_1+1$,
which is odd and $<2e$ and so
$B_1=\b_1=S_2-S_1=R_2-R_1+1=A_1=C_1$. Since $d[a_1b_1]\geq A_1=B_1$
and $d[a_1c_1]\geq C_1=B_1$ we have $d[b_1c_1]\geq B_1$. At $i=2$ we
have $B_2\leq S_3-T_2+\c_1=(R_1+1)-(R_2+1)+R_2-R_1+1=1$. Now
$S_1+S_2=T_1+T_2=R_1+R_2+1$ so $\ord b_{1,2}c_{1,2}$ is even and so
$d(b_{1,2}c_{1,2})\geq 1$. Also $T_3-T_2\geq T_1-T_2=R_1-R_2-1\geq
1-2e$ and $S_3-S_2=(R_1+1)-(R_2+1)=R_1-R_2\geq 2-2e$ so
$\b_2,\c_2\geq 1$. Therefore $d[b_{1,2}c_{1,2}]\geq 1\geq B_2$.
For $i\geq 3$ we have $d[a_{1,i}c_{1,i}]\geq C_1\geq B_i$ so
$d[b_{1,i}c_{1,i}]\geq B_i$ is equivalent to
$\b_i=d[a_{1,i}b_{1,i}]\geq B_i$. We have
$\b_i=\min\{\a_i,S_{i+1}-S_4+\b_3\}$. If $\b_i=\a_i$ then
$\b_i=\a_i\geq C_i\geq B_i$. So we may assume that
$\a_i>\b_i=S_{i+1}-S_4+\b_3$ and we have to prove that
$S_{i+1}-S_4+\b_3\geq B_i$. If $S_4-S_3\geq 2e$ then
$\b_3=(S_4-S_3)/2+e$ so
$\b_i=S_{i+1}-S_4+\b_3=S_{i+1}-(S_3+S_4)/2+e\geq
S_{i+1}-(S_i+S_{i+1})/2+e=(S_{i+1}-S_i)/2+e\geq B_i$ by 2.6. So we
may assume that that $S_4-S_3<2e$. We have $S_4-S_3=R_4-R_3-1$. If
$S_4-S_3$ is even then $\b_3\geq S_4-S_3+1$ and since
$R_4-R_3=S_4-S_3+1$ is odd and $\leq 2e$ we have
$\a_3=R_4-R_3\leq\b_3$. It follows that $\a_i\leq
R_{i+1}-R_4+\a_3\leq S_{i+1}-S_4+\b_3=\b_i$. Contradiction. So we may
assume that $S_4\ev S_3-1=R_1\m2$. Since $S_4-S_3$ is odd and
$<2e$ we have $\b_3=S_4-S_3$ and so
$\b_i=S_{i+1}-S_4+\b_3=S_{i+1}-S_3=S_{i+1}-R_1-1$. Note that
$T_2=R_2+1\geq R_1+1$ and $T_3\geq R_1+1$ (see above). Thus for
any $j\geq 2$ we have $T_j\geq T_2$ or $T_3$, depending on the
parity of $j$, and so $T_j\geq R_1+1$. Since
$\b_i=S_{i+1}-S_4+\b_3$, i.e. $-S_4+\b_3=-S_{i+1}+\b_i$ we have
$R_1\ev S_4\ev\ldots\ev S_{i+1}\m2$ by Lemma 7.3(ii). Since also
$S_1=R_1$ and $S_2\ev S_3\ev R_1+1\m2$ we have $S_1+\cdots +S_{i+1}\ev
(i+1)R_1\m2$. If $T_1+\cdots +T_{i-1}\ev S_1+\cdots +S_{i+1}+1\m2$ then
$\ord b_{1,i+1}c_{1,i-1}$ is odd so $B_i\leq
S_{i+1}-T_i+d[-b_{1,i+1}c_{1,i-1}]=S_{i+1}-T_i\leq
S_{i+1}-(R_1+1)=\b_i$ and we are done. So we may assume that
$T_1+\cdots +T_{i-1}\ev S_1+\cdots +S_{i+1}\ev (i+1)R_1\m2$. We have
$S_1+\cdots +S_i=R_1+\cdots +R_i+2\ev R_1+\cdots +R_i\m2$. If
$T_1+\cdots +T_i\ev S_1+\cdots +S_i+1\m2$ then both $\ord
b_{1,i}c_{1,i}$ and $\ord a_{1,i}c_{1,i}$ are odd and so
$d[b_{1,i}c_{1,i}]=0=d[a_{1,i}c_{1,i}]\geq C_i\geq B_i$. Thus we
may assume that $T_1+\cdots +T_i\ev S_1+\cdots +S_i\m2$. Together with
$T_1+\cdots +T_{i-1}\ev S_1+\cdots +S_ {i+1}\m2$, this implies $T_i\ev
S_{i+1}\ev R_1\m2$. Since $T_2=R_2+1\ev R_1+1\ev T_i+1\m2$ there
is some $2\leq j<i$ s.t. $ T_j,T_{j+1}$ have opposite parities. It
follows that $\ord c_{j,j+1}$ is odd and so  $\c_{i-1}\leq
T_i-T_j+d(-c_{j,j+1})=T_i-T_j$. Hence $B_i\leq
S_{i+1}-T_i+\c_{i-1}\leq S_{i+1}-T_j\leq S_{i+1}-(R_1+1)=\b_i$ and we
are done.

Before proving 2.1(iii) and (iv) note that $[b_2,b_3]\ap
[a_2,a_3]$ and $b_ i=a_i$ for $i\neq 2,3$ so $[b_1,\ldots,b_i]\ev
[a_1,\ldots,a_i]$ for $i\geq 3$. We prove 2.1(iii). At $i=2$ suppose
that $B_1+B_2>2e+S_2-T_2=2e$. (We have $S_2=T_2=R_2+1$.) Since $\ord
b_{1,2}=S_1+S_2=R_1+R_2+1$ is odd we have $\b_1\leq
S_2-S_1+d(-b_{1,2})=R_2+1-R_1+0=R_2-R_1+1$ and $\b_2\leq
S_3-S_1+d(-b_{1,2})=R_1+1-R_1+0=1$ so $R_2-R_1+2\geq\b_1+\b_2\geq
B_1+B_2>2e$ and so $R_2-R_1>2e-2$. Contradiction. If $i=3$ suppose
that $T_2<S_4$ and $B_2+B_3>2e+S_3-T_3$. We have $T_2<S_4=R_4$ and
we claim that $C_2+C_3>2e+R_3-T_3$. By Lemma 2.14 we may assume
that $C_2=C'_2$. We have $C_3\geq B_3$ and $B_2\leq B'_2=\min\{
S_3-T_2+d[-b_{1,3}c_1], S_3+S_4-T_1-T_2+d[b_{1,4}]\}\leq\min\{
R_3-T_2+d[-a_{1,3}c_1]+1, R_3+R_4-T_1-T_2+d[a_{1,4}]+1\}
=C'_2+1=C_2+1$. (We have $S_3=R_3+1$, $S_4=R_4$,
$d[-b_{1,3}c_1]\leq d[-a_{1,3}c_1]$ and $d[b_{1,4}]\leq
d[a_{1,4}]$.) Thus $C_2+1+C_3\geq B_2+B_3>2e+S_3-T_3=2e+R_3-T_3+1$,
i.e. $C_2+C_3>2e+R_3-T_3$. Hence $[c_1,c_2]\rep [a_1,a_2,a_3]\ap
[b_1,b_2,b_3]$. If $i\geq 4$, $S_{i+1}>T_{i-1}$ and
$B_{i-1}+B_i>2e+S_i-T_i$ then $R_{i+1}=S_{i+1}>T_{i-1}$ and
$C_{i-1}+C_i\geq B_{i-1}+B_i>2e+S_i-T_i=2e+R_i-T_i$ and so
$[c_1,\ldots,c_{i-1}]\rep [a_1,\ldots,a_i]\ap [b_1,\ldots,b_i]$. For
2.1(iv) if $i\geq 2$ and $T_i\geq S_{i+2}>T_{i-1}+2e\geq S_{i-1}+2e$
then $R_{i+2}=S_{i+2}>T_{i-1}+2e$ so by Lemma 2.19 we have
$[c_1,\ldots,c_{i-1}]\rep [a_1,\ldots,a_{i+1}]\ap [b_1,\ldots,b_{i+1}]$.

2. $d[-a_{1,3}c_1]>\a_1$. We prove that $R_2+2\leq R_4,T_2$ so we can
apply our claim. If $R_4=R_2+1$ then $R_4-R_3=R_2-R_1+1$, which is odd
and $\leq 2e-1$ so $\a_1<d[-a_{1,3}c_1]\leq\a_3=R_2-R_1+1$. But this
is impossible since $R_2-R_1$ is even and $<2e$. So $R_4\geq
R_2+2$. Similarly, if $T_2=R_2+1$ then $T_2-T_1=R_2-R_1+1$ so this
time $\a_1<d[-a_{1,3}c_1]\leq\c_1=R_2-R_1+1$ and again we get a
contradiction. So $T_2\geq R_2+2$.

Suppose first that $\a_1<(R_2-R_1)/2+e$, which implies that $\a_1$ is
odd and also $R_2-R_1\neq 2e-2$. Thus $R_2-R_1\leq 2e-4$. We claim
that $d[-a_{1,3}c_1]\geq\a_1+2$. Suppose the contrary. We have
$\a_1<d[-a_{1,3}c_1]<\a_1+2$ so $d[-a_{1,3}c_1]$ is not an odd
integer. If $d[-a_{1,3}c_1]=d(-a_{1,3}c_1)$ this implies that
$d(-a_{1,3}c_1)=2e=(2e-4)/2+e+2\geq
(R_2-R_1)/2+e+2>\a_1+2$. Contradiction. If $d[-a_{1,3}c_1]=\a_3$
then $\a_3=(R_4-R_3)/2+e\geq
(R_2+2-R_1)/2+e=(R_2-R_1)/2+e+1\geq\a_1+2$. If $d[-a_{1,3}c_1]=\c_1$
then $\c_1=(T_2-T_1)/2+e\geq
(R_2+2-R_1)/2+e=(R_2-R_1)/2+e+1\geq\a_1+2$. Contradiction. (Note that
$(R_2-R_1)/2+e+1\geq\a_1+2$ follows from $\a_1<(R_2-R_1)/2+e$ and
$\a_1,(R_2-R_1)/2+e\in\ZZ$.) So $d[-a_{1,3}c_1]\geq\a_1+2$. In
particular, $\a_3\geq\a_1+2$.

We will take $N$ like in Lemma 9.10 s.t. $\b_1=\a_1+2$. First we
show that such $N$ exists. We have $\b_1=\a_1+2>0$ and
$\b_1=\a_1+2\geq R_2-R_1+2=S_2-S_1$. Also $\a_1\leq
(R_2-R_1)/2+e-1$ and so $\b_1=\a_1+2\leq
(R_2-R_1)/2+e+1=(R_2+2-R_1)/2+e=(S_2-S_1)/2+e$. Thus the condition
$\max\{ 0,S_2-S_1\}\leq\b_1\leq (S_2-S_1)/2+e$ from Lemma 9.9 is
fulfilled. Also $\a_1$ is odd and so is $\b_1=\a_1+2$. Thus the
conditions that $\b_1=(S_2-S_1)/2+e$ and $V:=[a_1,a_2,a_3]$ is
isotropic from Lemma 9.9 are not necessary. We have
$\a_3\geq\a_1+2=\b_1>\a_1$ so the condition
$\a_1\leq\b_1\leq\min\{\a_1+2,\a_3\}$ from Lemma 9.10 is
fulfilled. Thus $N$ exists. In order to apply our claim and prove
that $K\leq N$ we still need $d[-a_{1,3}c_1]\geq\b_1$ and
$S_{i+1}-S_2+\b_1\geq B_i$ for $i\geq 2$. But
$d[-a_{1,3}c_1]\geq\a_1+2=\b_1$ and for $i\geq 2$ we have
$S_{i+1}-S_2+\b_1=R_{i+1}-(R_2+2)+\a_1+2=R_{i+1}-R_2+\a_1\geq\a_i\geq
C_i\geq B_i$ so we are done.

Suppose now that $\a_1=(R_2-R_1)/2+e$. Since $d[-a_{1,3}c_1]>\a_1$ and
$\a_1<2e$ we have $d[-a_{1,3}c_1]\geq\a_1+1$. In particular,
$\c_1\geq\a_1+1$ and $\a_3\geq\a_1+1$. 

Since $N,K$ satisfy 2.1(i) we have $(S_{i+1}-S_i)/2+e\geq B_i$ by
2.6.

Assume first that $[a_1,a_2,a_3]$ is isotropic. We will take $N$ s.t.
$\b_1=\a_1+1=(R_2-R_1)/2+e+1=(R_2+2-R_1)/2+e=(S_2-S_1)/2+e$. We
prove first that such $N$ exists. Since $R_2-R_1$ is even and
$<2e$ we have $\a_1\geq R_2-R_1+1$ so $\b_1=\a_1+1\geq
R_2-R_1+2=S_2-S_1$ and also $\b_1=\a_1+1>0$ so the condition
$\max\{ 0,S_2-S_1\}\leq\b_1\leq (S_2-S_1)/2+e$ from Lemma 9.9 is
fulfilled. Also both conditions that $\b_1=(S_{i+1}-S_i)/2+e$ and
$V:=[a_1,a_2,a_3]$ is isotropic required by Lemma 9.9 if $\b_1$
is even are satisfied. We have $\a_3\geq\a_1+1=\b_1$ so the condition
$\a_1\leq\b_1\leq\min\{\a_1+2,\a_3\}$ from Lemma 9.10 is
satisfied. Thus $N$ with $\b_1=\a_1+1$ exists. We have
$d[-a_{1,3}c_1]\geq\a_1+1=\b_1$ and for $i\geq 2$ we have
$S_{i+1}-S_2+\b_1=S_{i+1}-(S_1+S_2)/2+e\geq
S_{i+1}-(S_i+S_{i+1})/2+e=(S_{i+1}-S_i)/2+e\geq B_i$. So $K\leq N$ by
our claim. 

Suppose now that $[a_1,a_2,a_3]$ is not isotropic. Note that
$(R_2-R_1)/2+e\geq\a_1(\[ a_1,a_2,a_3\] )\geq\a_1$ so $\a_1(\[
a_1,a_2,a_3\] )=\a_1$. Since $[a_1,a_2,a_3]$ is not isotropic we
have by Lemma 9.9 that $\a_1=\a_1(\[ a_1,a_2,a_3\] )$ is odd. If
$R_2-R_1=2e-2$ then also $d[-a_{1,3}c_1]\geq\a_1+1=2e$ and so
$M,K$ are in the case of Lemma 9.6. Thus we may assume that
$R_2-R_1\leq 2e-4$. We will take $N$ s.t.
$\b_1=\a_1=(R_2-R_1)/2+e$. We prove first that this is possible. Now
$\b_1=\a_1>0$ and, since $R_2-R_1\leq 2e-4$, we also have
$\b_1=(R_2-R_1)/2+e\geq R_2-R_1+2=S_2-S_1$. Also
$\b_1=(R_2-R_1)/2+e<(S_2-S_1)/2+e$ so the condition $\max\{
0,S_2-S_1\}\leq\b_1\leq (S_2-S_1)/2+e$ from Lemma 9.9 is
satisfied. Also $\b_1=\a_1$ is odd so the conditions that
$\b_1=(S_2-S_1)/2+e$ and $V:=[a_1,a_2,a_3]$ is isotropic from
Lemma 9.9 are not necessary. Since $\a_3\geq R_1-R_3+\a_1=\a_1=\b_1$
the condition $\a_1\leq\b_1\leq\min\{\a_1+2,\a_3\}$ from Lemma 9.10 is
satisfied. So $N$ with $\b_1=\a_1$ exists. We have
$d[-a_{1,3}c_1]>\a_1=\b_1$ so we still need $S_{i+1}-S_2+\b_1\geq
B_i$ for $i\geq 2$.

Suppose that
$B_i>S_{i+1}-S_2+\b_1=S_{i+1}-R_2+\a_1-2=S_{i+1}-(R_1+R_2)/2+e-2$ for
some $i\geq 2$. We have $(S_{i+1}-S_i)/2+e\geq
B_i>S_{i+1}-(R_1+R_2)/2+e-2$ so $R_1+R_2+4>S_i+S_{i+1}$. Also
$S_{i+1}-(T_{i-1}+T_i)/2+e\geq 
B_i>S_{i+1}-(R_1+R_2)/2+e-2$ so $R_1+R_2+4>T_{i-1}+T_i$. For any
$3\leq j\leq i$ we have $R_{i+1}-R_{j+1}+\a_j\geq\a_i\geq C_i\geq
B_i>S_{i+1}-S_2+\b_1\geq S_{i+1}-S_{j+1}+\b_j=R_{i+1}-R_{j+1}+\b_j$ so
$\a_j>\b_j$. Therefore $R_{j+1}-R_j=S_{j+1}-S_j$ must be
even. (Otherwise $\a_j=\b_j$ by [B3, Corollary 2.9].) Thus all
$R_j=S_j$ for $3\leq j\leq i+1$ have the same parity. Namely they are
$\ev R_3=R_1\m2$. The congruence $S_j\ev R_1\m2$ also holds for 
$i=1,2$ so it holds for $1\leq j\leq i+1$. (We have $S_1=R_1$ and
$S_2=R_2+2\ev R_1\m2$.) If $T_{j+1}-T_j$ is odd for some $1\leq j\leq
i-1$ then $\ord c_{j.j+1}$ is odd and so $\c_{i-1}\leq
T_i-T_j+d(-c_{j,j+1})=T_i-T_j$. It follows that $S_{i+1}-T_j\geq
S_{i+1}-T_i+\c_{i-1}\geq B_i>S_{i+1}-(R_1+R_2)/2+e-2$ so
$R_1+R_2+4-2e>2T_j$. But $T_j\geq T_1=R_1$ or $T_j\geq T_2\geq R_2+2$,
depending  on the parity of $j$. In the first case we get
$R_1+R_2+4-2e>2R_1$, which implies $R_2-R_1>2e-4$. Contradiction. In
the second case $R_1+R_2+4-2e>2(R_2+2)$ so $-2e>R_2-R_1$, which is
impossible. Thus $T_{j+1}-T_j$ is even. It follows that all $T_j$ for
$1\leq j\leq i$ have the same parity. Namely they are $\ev
T_1=R_1\m2$. Now $R_1+R_2+2=S_1+S_2\leq S_i+S_{i+1}<R_1+R_2+4$. Since
both $R_1+R_2$ and $S_i+S_{i+1}$ are even we have
$S_1+S_2=S_i+S_{i+1}=R_1+R_2+2$. By [B3, Corollary 2.3] this implies
that for any $1\leq j\leq i+1$ we have $S_j=S_1=R_1$ if $j$ is odd and
$S_j=S_2=R_2+2$ if $j$ is even. Similarly $R_2+2\leq T_2$ so
$R_1+R_2+2\leq T_1+T_2\leq T_{i-1}+T_i<R_1+R_2+4$. Since $R_1+R_2$,
$T_1+T_2$ and $T_{i-1}+T_i$ are all even we get
$T_1+T_2=T_{i-1}+T_i=R_1+R_2+2$. Since $T_1=R_1$ we have
$T_2=R_2+2$. Again by [B3, Corollary 2.3], for any $1\leq j\leq i$ we
have $T_j=T_1=R_1$ if $j$ is odd and $T_j=T_2=R_2+2$ if $j$ is
even. In particular, $T_j=S_j$ for $1\leq j\leq i$. 

Now $\c_1\geq\a_1+1=(R_2-R_1)/2+e+1=(T_2-T_1)/2+e$ so
$\c_1=\a_1+1=(T_2-T_1)/2+e$, which is even ($\a_1$ is odd). But
$T_1+T_2=T_{i-1}+T_i$ so by [B3, Corollary 2.3] we have
$\c_j=(T_{j+1}-T_j)/2+e=(S_{j+1}-S_j)/2+e$ for any $1\leq j\leq
i-1$. Also note that $\c_1=(R_2-R_1)/2+e+1<(2e-2)/2+e+1=2e$ and
$\c_1=\a_1+1>0$. Since $0<\c_1<2e$ and $\c_1$ is even, so is
$2e-\c_1=(S_1-S_2)/2+e$. If $i\geq 3$ then also
$(R_4-R_3)/2+e\geq\a_3\geq\a_1+1=(R_2-R_1)/2+e+1=(R_4-R_3)/2+e$. Thus
$\a_3=(R_4-R_3)/2+e$. Since $R_3+R_4=S_3+S_4=S_i+S_{i+1}=R_i+R_{i+1}$
this implies that $\a_j=(R_{j+1}-R_j)/2+e=(S_{j+1}-S_j)/2+e$ for any
$3\leq j\leq i$. Note that $\a_j=(S_2-S_1)/2+e=\c_1$ if $j$ is odd and
$\a_j=(S_1-S_2)/2+e=2e-\c_1$ if $j$ is even. 

We want to prove that $[c_1,\ldots,c_i]\rep
[a_1,\ldots,a_{i+1}]$. This is obvious if $i+1=n$ so we may assume
that $n\geq i+2$. Since $d[b_{1,i}c_{1,i}]\leq\b_i\leq
S_{1+1}-S_2+\b_1<B_i$ we have by Lemma 2.12 that at least one of
$i,i+1$ is an essential index for $N,K$. But $i$ cannot be
essential because $S_{i+1}=S_{i-1}=T_{i-1}$ and so $i+1$ is
essential. Thus $S_i=T_i<S_{i+2}=R_{i+2}$ and, if $n\geq i+3$,
$T_{i-1}+T_i<S_{i+2}+S_{i+3}=R_{i+2}+R_{i+3}$. We have $R_{i+2}>T_i$ and we
want to prove that $C_i+C_{i+1}>2e+R_{i+1}-T_{i+1}$. We have
$d[a_{1,i}c_{1,i}]\geq C_i\geq B_i>S_{i+1}-(R_1+R_2)/2+e-2=
S_{i+1}-(S_i+S_{i+1})/2+e-1=(S_{i+1}-S_i)/2+e-1$. But
$(S_{i+1}-S_i)/2+e-1$ is equal to $(S_2-S_1)/2+e-1=\c_1-1$ or
$(S_1-S_2)/2+e-1=2e-\c_1-1$ so it is integer and $<2e$ so the
inequality above implies $d[a_{1,i}c_{1,i}]\geq
(S_{i+1}-S_i)/2+e$. Also 
$S_{i+1}-S_i+d[-a_{1,i+1}c_{1,i-1}]=R_{i+1}-T_i+d[-a_{1,i+1}c_{1,i-1}]
\geq C_i\geq B_i>(S_{i+1}-S_i)/2+e-1$ so
$d[-a_{1,i+1}c_{1,i-1}]>(S_i-S_{i+1})/2+e-1$. Again
$(S_i-S_{i+1})/2+e-1\in\{\c_1-1,2e-\c_1-1\}$ so it is integer and $<2e$ so
$d[-a_{1,i+1}c_{1,i-1}]\geq (S_i-S_{i+1})/2+e$. In particular,
$\a_{i+1}\geq (S_i-S_{i+1})/2+e$. Now
$(S_i-S_{i+1})/2+e\in\{\c_1,2e-\c_1\}$ so it is even so if we have
equality then $\a_{i+1}$ is even and we have
$(S_i-S_{i+1})/2+e=\a_{i+1}=(S_{i+2}-S_{i+1})/2+e$ so
$S_i=S_{i+2}$. Contradiction. Thus $\a_{i+1}>(S_i-S_{i+1})/2+e$, which
implies that $d[-a_{i+1,i+2}]\geq
R_{i+1}-R_{i+2}+\a_{i+1}=S_{i+1}-S_{i+2}+\a_{i+1}>(S_i+S_{i+1})/2-S_{i+2}+e$. Together
with $d[a_{1,i}c_{1,i}]\geq
(S_{i+1}-S_i)/2+e=(S_i+S_{i+1})/2-S_i+e>(S_i+S_{i+1})/2-S_{i+2}+e$, this
implies $d[-a_{1,i+2}c_{1,i}]>(S_i+S_{i+1})/2-S_{i+2}+e$. Now
$d[-a_{1,i+1}c_{1,i-1}]+d[-a_{1,i+2}c_{1,i}]>
(S_i-S_{i+1})/2+e+(S_i+S_{i+1})/2-S_{i+2}+e=
2e+S_i-S_{i+2}=2e+T_i-R_{i+2}$. Also
$d[-a_{1,i+2}c_{1,i}]>(S_i+S_{i+1})/2-S_{i+2}+e=e+(T_{i-1}+T_i)/2-R_{i+2}$.
Finally, if $n\geq i+3$ then $d[-a_{1,i+1}c_{i-1}]\geq
(S_i-S_{i+1})/2+e=(T_i-T_{i-1})/2+e>e+T_i-(R_{i+2}+R_{i+3})/2$. (We
have $S_{i+1}=S_{i-1}=T_{i-1}$ and $R_{i+2}+R_{i+3}>T_{i-1}+T_i$.)
These imply by Lemma 2.12 that $C_i+C_{i+1}>2e+R_{i+1}-T_{i+1}$. Hence
$[c_i,\ldots,c_i]\rep [a_1,\ldots,a_{i+1}]$.

We now prove by descending induction that $[c_1,\ldots,c_j]\rep
[a_1,\ldots,a_{j+1}]$ for $2\leq j\leq i$. This statement was proved
for $j=i$. Suppose that $[c_1,\ldots,c_j]\rep [a_1,\ldots,a_{j+1}]$
for some $2<j\leq i$. Since also $[c_1,\ldots,c_{j-1}]\rep
[c_1,\ldots,c_j]$ we have by Lemma 1.5(ii) that
$[c_1,\ldots,c_{j-1}]\rep [a_1,\ldots,a_j]$ iff
$(a_{1,j}c_{1,j},-a_{1,j+1}c_{1,j-1})_\p =1$. To prove this we will
show that $d(a_{1,j}c_{1,j})+d(-a_{1,j+1}c_{1,j-1})>2e$. We claim that
$d[a_{1,j}c_{1,j}]\geq (S_{j+1}-S_j)/2+e$ and
$d[-a_{1,j+1}c_{1,j-1}]\geq (S_j-S_{j+1})/2+e$. For $j=i$ these inequalities
were already proved. Before proving them for $j<i$ note that for $3\leq
h<l\leq i+1$ with $h\ev l\m2$ we have $R_h=R_l$ so, by Lemma 7.4(iii),
$d[(-1)^{(l-h)/2}a_{h+1,l}]=\a_h=(S_{h+1}-S_h)/2+e$. Similarly for $1\leq
h<l\leq i$ with $h\ev l\m2$ we have $T_h=T_l$ so
$d[(-1)^{(l-h)/2}c_{h+1,l}]=\c_h=(S_{h+1}-S_h)/2+e$. Let $3\leq j<i$. If $j\ev
i\m2$ then $d[a_{1,i}c_{1,i}]\geq (S_{i+1}-S_i)/2+e=(S_{j+1}-S_j)/2+e$,
$d[(-1)^{(j-i)/2}a_{j+1,i}]=(S_{j+1}-S_j)/2+e$ and
$d[(-1)^{(j-i)/2}c_{j+1,i}]=(S_{j+1}-S_j)/2+e$ so $d[a_{1,j}c_{1,j}]\geq
(S_{j+1}-S_j)/2+e$. Also
$d[-a_{1,i+1}c_{1,i-1}]\geq (S_i-S_{i+1})/2+e=(S_j-S_{j+1})/2+e$,
$d[(-1)^{(i-j)/2}a_{j+2,i+1}]=(S_{j+2}-S_{j+1})/2+e=(S_j-S_{j+1})/2+e$ and
$d[(-1)^{(i-j)/2}c_{j,i-1}]=(S_j-S_{j-1})/2+e=(S_j-S_{j+1})/2+e$ and so
$d[-a_{1,j+1}c_{1,j-1}]\geq (S_j-S_{j+1})/2+e$. If $j\ev i+1\m2$ then
$d[-a_{1,i+1}c_{1,i-1}]\geq (S_i-S_{i+1})/2+e=(S_{j+1}-S_j)/2+e$,
$d[(-1)^{(j-i+1)/2}a_{j+1,i+1}]=(S_{j+1}-S_j)/2+e$ and
$d[(-1)^{(j-i-1)/2}c_{j+1,i-1}]=(S_{j+1}-S_j)/2+e$ so $d[a_{1,j}c_{1,j}]\geq
(S_{j+1}-S_j)/2+e$. Also $d[a_{1,i}c_{1,i}]\geq
(S_{i+1}-S_i)/2+e=(S_j-S_{j+1})/2+e$,
$d[(-1)^{(i-j-1)/2}a_{j+2,i}]=(S_{j+2}-S_{j+1})/2+e=(S_j-S_{j+1})/2+e$ and
$d[(-1)^{(i-j+1)/2}c_{j,i}]=(S_j-S_{j-1})/2+e=(S_j-S_{j+1})/2+e$ and so
$d[-a_{1,j+1}c_{1,j-1}]\geq (S_j-S_{j+1})/2+e$. It follows that
$d(a_{1,j}c_{1,j})\geq (S_{j+1}-S_j)/2+e$ and $d(-a_{1,j+1}c_{1,j-1})\geq
(S_j-S_{j+1})/2+e$. But, depending on the parity of $j$, $(S_{j+1}-S_j)/2+e$
and $(S_j-S_{j+1})/2+e$ are equal, in some order, with $(S_2-S_1)/2+e=\c_1$
and $(S_1-S_2)/2+e=2e-\c_1$. Since $\c_1,2e-\c_1$ are even and
$0<\c_1,2e-\c_1<2e$ the inequalities above are strict. Thus
$d(a_{1,j}c_{1,j})+d(-a_{1,j+1}c_{1,j-1})>\c_1+2e-\c_1=2e$ and we are done.

When $j=2$ we have $[c_1,c_2]\rep [a_1,a_2,a_3]$ so $[a_1,a_2,a_3]\ap
[c_1,c_2,a_{1,3}c_{1,2}]$. Now $d(-c_{1,2})\geq
T_1-T_2+\c_1=S_1-S_2+\c_1=2e-\c_1$ and $d(-a_{1,3}c_1)\geq
d[-a_{1,3}c_1]\geq\a_1+1=\c_1$. Since $\c_1,2e-\c_1$ are even and
$0<\c_1,2e-\c_1<2e$ the two inequalities are strict and so
$d(-c_{1,2})+d(-a_{1,3}c_1)>\c_1+2e-\c_1=2e$, which implies that
$(-c_{1,2},-a_{1,3}c_1)_\p =1$ and so $[c_1,c_2,-a_{1,3}c_{1,2}]\ap
[a_1,a_2,a_3]$ is isotropic. Contradiction. \qed

\section*{References}

\hskip 6mm [B1] C. N. Beli, Integral spinor norm groups over dyadic
local fields, J. Number Theory 102, No. 1, 125-182 (2003). 

[B2] C. N. Beli, Representations of integral quadratic forms over
dyadic local fields, Electron. Res. Announc. Am. Math. Soc. 12,
100-112, electronic only (2006). 

[B3] C. N. Beli, A new approach to classification of integral
quadratic forms over dyadic local fields. Trans. Am. Math. Soc. 362,
No. 3, 1599-1617 (2010).

[EH1] A. G. Earnest, J. S. Hsia, Spinor norms of local integral
rotations. II, Pacific J. Math. Volume 61, Number 1 (1975), 71-86.

[EH2] A. G. Earnest, J. S. Hsia, Spinor genera under field
extensions. II: 2 unramified in the bottom field, Am. J. Math. 100,
523-538 (1978). 

[H] J. S. Hsia, Spinor norms of local integral rotations I, Pacific
J. Math., Vol. 57 (1975), 199 - 206.

[HSX] J. S. Hsia, Y. Y. Shao, Fei Xu, Representations of indefinite
quadratic forms, J. Reine Angew. Math. 494, 129-140 (1998).

[J] D. G. James, Primitive representations by unimodular quadratic
forms, J. Number Theory 44, No.3, 356-366 (1993).

[LX] JianRui Lv, Fei Xu, Integral spinor norms in dyadic local fields
III, Sci. China, Math. 53, No. 9, 2425-2446 (2010). 

[OSU1] Y. Y. Shao, Representation theory of quadratic forms, (Ohio
State University thesis, 1994).

[OSU2] C. N. Beli, Integral spinor norm groups over dyadic local fields
and representations of quadratic lattices, (Ohio State University
thesis, 2001).

[OM] O. T. O'Meara, Introduction to Quadratic Forms, Springer-Verlag,
Berlin (1952).

[OM1] O.T. O'Meara, The Integral Representation of Quadratic Forms
over Local Fields, Amer. J. Math. 80 (1958), 843-878.

[R] C. Riehm, On the Integral Representations of Quadratic Forms over
Local Fields, Amer. J. Math. 86 (1964), 25-62.

[X] Fei Xu, Arithmetic Springer theorem on quadratic forms under field
extensions of odd degree, American Mathematical
Society. Contemp. Math. 249, 175-197 (1999).

\end{document}